\documentclass[oneside,a4paper,12pt]{book}

\pdfoutput=1

\usepackage{amsmath,amssymb,wasysym}
\usepackage[T1]{fontenc}
\usepackage[utf8x]{inputenc} 
\usepackage[french]{babel}
\usepackage{arabtex}
\usepackage{enumerate}
\usepackage{url}
\usepackage{pdfpages}
\usepackage{graphicx}
\usepackage{color}
\usepackage[normalem]{ulem}
\usepackage{hyperref}
\novocalize
\usepackage{makeidx}
\makeindex

\begin{document}
\newcommand{\shatir}{Ibn al-\v{S}\=a\d{t}ir}
\newcommand{\shirazi}{Sh\={\i}r\=az\={\i}}
\newcommand{\shariaa}{Sadr al-Shar\={\i}`a}
\newcommand{\tusi}{\d{T}\=us{\=\i}}

\begin{titlepage}
\begin{center}
\Large La \textit{Nih\=aya al-s\=ul f{\=\i} ta\d{s}\d{h}{\=\i}\d{h} al-'u\d{s}\=ul} d'{\shatir}

\bigskip\large \'Edition, traduction et commentaire

\bigskip\large Erwan Penchèvre
\end{center}
\end{titlepage}

\renewcommand{\contentsname}{Table des matières}
\tableofcontents

\pagestyle{plain}

\part*{Introduction générale}
\addcontentsline{toc}{part}{Introduction générale}
\noindent En 2009, j'avais lu l'\textit{Almageste} de Ptolémée et la \textit{Configuration des mouvements} d'Ibn al-Haytham, et j'avais suivi un cours de Régis Morelon sur l'histoire de l'astronomie ancienne. Fier de mes balbutiements en arabe, c'est avec ce bagage que je suis allé voir Roshdi Rashed. Je souhaitais étudier un texte médiéval arabe d'astronomie, il m'a suggéré de me pencher sur la \textit{Nih\=aya al-s\=ul f{\=\i} ta\d{s}\d{h}{\=\i}\d{h} al-'U\d{s}\=ul} d'{\shatir} et il a mis à ma disposition un premier manuscrit.

Le présent ouvrage est le fruit de cet exercice. L'édition du texte arabe et la traduction en français (p.~\pageref{txt_debut}-\pageref{txt_fin}) sont l'{\oe}uvre d'un débutant surtout intéressé par les questions mathématiques. Le commentaire mathématique (p.~\pageref{comm_debut}-\pageref{comm_fin}) tient donc une place qui m'est chère~; j'y ai seulement analysé la ``théorie planétaire'' contenue dans la première partie de la \textit{Nih\=aya}\footnote{On rencontrera deux autres titres commençant par le même mot, mais le titre abrégé \textit{Nih\=aya} désignera toujours la \textit{Nih\=aya al-s\=ul f{\=\i} ta\d{s}\d{h}{\=\i}\d{h} al-'u\d{s}\=ul} d'{\shatir}.} .

En guise d'introduction, j'ai réuni ici quelques éléments concernant la tradition matérielle qu'il est impossible de séparer du contexte et de l'Histoire~; j'indiquerai ensuite les conclusions générales qui se dégagent de mon étude.

\paragraph{\'Eléments biographiques}
Comme souvent ches les savants anciens, il est difficile de rassembler suffisamment d'éléments pour tracer un portrait bien net d'{\shatir}\footnote{Pour le présent paragraphe, mes principales sources sont \cite{wiedemann1928} et \cite{king1975}.}. Né à Damas en 1306, il est recueilli après la mort de son père par un oncle qui lui enseigne l'art de ciseler l'ivoire ou la nacre. \`A l'âge de dix ans, il se rend au Caire et à Alexandrie où il étudie l'astronomie. Il construit des instruments mathématiques (quadrants\footnote{\textit{Cf.} \cite{schmalzl1929} p.~100-108.}), des cadrans solaires, et des astrolabes\footnote{Outre l'astrolabe de l'Observatoire de Paris mentionné ci-après, il y a eu un astrolabe d'{\shatir} dans les collections de la Bibliothèque Nationale de France. Il est décrit par {L.~A. Sédillot} qui le date de 1337, \textit{cf.} \cite{sedillot1841} p.~191-194~; mais cet astrolabe, d'ailleurs incomplet, aurait rejoint la collection Hariri du Musée d'Art Islamique du Caire dans les années 1920. Je remercie Catherine Hofmann, Conservatrice du Département des Cartes et Plans de la B. N. F., pour cette dernière information qu'elle tient d'un ouvrage d'A. Turner à paraître en 2017. Une photographie figure dans \cite{king2005} t.~1 p.~729 fig.~8.2.}. Il y a à l'Observatoire de Paris un astrolabe qu'il aurait construit en 1326~: il avait alors seulement vingt ans\footnote{Observatoire de Paris, numéro d'inventaire 1, ancienne cote 14-1. Il en existe une bonne photographie dans le catalogue informatisé. Est gravé au dos de l'astrolabe~: \textit{\d{s}ana`ahu `al\=a b. ibr\=ah{\=\i}m b. mu\d{h}ammad b. ab{\=\i} mu\d{h}ammad b. ibr\=ah{\=\i}m sana sitta wa-`ishr\=un wa-sab`a-mi'a}. D. A. King fait l'analyse détaillée de cet instrument dans \cite{king2005} t.~2 p.~692-693.}. Un contemporain a décrit un astrolabe automatique qu'{\shatir} exposait dans sa propre demeure\footnote{Dans la traduction de Sauvaire \cite{sauvaire1896} t. 2 p.~263, l'historien al-\d{S}afad{\=\i} contemporain d'{\shatir} rapporte~: ``je l'ai vu plus d'une fois et suis entré dans son logis au mois de \textit{rama\d{d}\=an} de l'année 743 [=1345], pour examiner l'astrolabe qu'il avait inventé. Je trouvai qu'il l'avait placé dans la verticale d'un mur, dans sa demeure [...]. Cet astrolabe, qui avait la forme d'une arcade [...] mesurait trois quarts de coudée environ~; il tournait toujours, continuellement, le jour et la nuit, sans sable, ni eau, suivant les mouvements de la sphère céleste, mais il l'avait réglé sur des dispositions particulières. Cet instrument faisait connaître les heures égales et les heures de temps.''}. \textit{Muwaqqit} à la Mosquée de Damas, chargé de calculer les heures des prières, il y a aussi construit un grand cadran solaire horizontal à style incliné, daté de 1371, endommagé puis reconstruit presque à l'identique au XIXème siècle, ``le cadran solaire le plus sophistiqué de la période médiévale'' selon David A. King\footnote{\textit{Cf.} \cite{king2005} t.~2 p.~85. Louis Janin en a donné une description en 1972, \textit{cf.} \cite{ghanem1976} 108-121~; une photographie de la pierre dédicatoire est à la p.~72, un dessin à l'échelle de la copie d'al-\d{T}an\d{t}\=aw{\=\i} (1873) à la p.~112, et une photographie du style p.~111.}. {\shatir} meurt vers 1375 à Damas.

Outre des ouvrages sur les instruments mathématiques et les astrolabes, il est l'auteur de quatre traités d'astronomie théorique dont l'ordre de rédaction devait être le suivant~:

-- \textit{Nih\=aya al-gh\=ay\=a f{\=\i} a`m\=al al-falak{\=\i}y\=at}, ouvrage aujourd'hui perdu qui aurait été un \textit{z{\=\i}j} (recueil de tables astronomiques) strictement ptoléméen\footnote{d'après la description qu'en donne {\shatir} lui-même dans son \textit{Z{\=\i}j al-jad{\=\i}d}.}.

-- \textit{Ta`l{\=\i}q al-ar\d{s}\=ad}, le \textit{Commentaire des observations}, ouvrage aujourd'hui perdu dans lequel {\shatir} aurait consigné des observations astronomiques et en aurait déduit de nouveaux modèles.

-- \textit{Nih\=aya al-s\=ul f{\=\i} ta\d{s}\d{h}{\=\i}\d{h} al-'u\d{s}\=ul}, objet de la présente édition. Y sont décrits de nouveaux modèles astronomiques~; c'est un ouvrage d'exposition dans la tradition des \textit{kit\=ab al-hay'a}. Les démonstrations en sont absentes~; {\shatir} renvoie pour cela au \textit{Ta`l{\=\i}q al-ar\d{s}\=ad}.

-- \textit{Al-z{\=\i}j al-jad{\=\i}d}~: un \textit{z{\=\i}j} construit au moyen des modèles décrits dans la \textit{Nih\=aya}. De nombreuses copies sont conservées\footnote{Une description précise de ce \textit{z{\=\i}j} est donnée par Kennedy \cite{kennedy1956}.}.

Que retenir de ces quelques repères biobibliographiques~? Tôt astronome de métier, il ne faut pas chercher dans {\shatir} la figure du savant universel qu'on trouverait chez un B{\=\i}r\=un{\=\i} ou un Na\d{s}{\=\i}r al-D{\=\i}n al-\d{T}\=us{\=\i}. Son {\oe}uvre théorique se résume peut-être aux quatre ouvrages mentionnés ci-dessus. Plusieurs traces témoignent de son intérêt pour ces modèles mécaniques des cieux que sont les astrolabes~; ce sont aussi des instruments d'observation, et il est assez remarquable que l'observation soit le sujet d'un de ces titres, hélas perdu.

\paragraph{Les manuscrits de la \textit{Nihaya}}
Voici les manuscrits dont j'ai connaissance~:

- Oxford, Bodleian Library, Marsh 139, que je désignerai par la lettre A [A,\RL{'a}]. Il compte 64 \textit{folii}. Le colophon indique qu'il a été copié en 768 [=1366]. La page de titre désigne l'auteur comme suit~: \textit{al-\v{s}ay\b{h} al-im\=am ab\=u al-\d{hasan} `al{\=\i} b. ibr\=ah{\=\i}m b. mu\d{h}ammad b. al-hum\=am ab{\=\i} mu\d{h}ammad b. ibr\=ah{\=\i}m b. `abd al-ra\d{h}man al-an\d{s}\=ar{\=\i} al-muwaqqit bi-al-j\=ami` al-ummawiy al-ma`r\=uf bi-ibn al-\v{s}\=a\d{t}ir}.

- Oxford, Bodleian Library, Marsh 290 [B,\RL{b}], 63 \textit{folii}. \uwave{Une date, 941} [=1535], est indiquée f.~1r.

- Oxford, Bodleian Library, Marsh 501 [C,\RL{_h}], 71 \textit{folii}. La page de titre indique \uwave{qu'il a été copié en 984 [=1577]}.

- Oxford, Bodleian Library, Hunt 547 [D,\RL{d}]. Ce sont les f. 21 à 65 d'un codex datant de 979 [=1571].
Les portions inférieures des f.~26r-30v manquent car elles ont dû être endommagées lors d'un incendie.

- Bibliothèque de Leyde, Or 194 [E,\RL{h}], 101 \textit{folii}, sans date.

- Bibliothèque Nationale d'Egypte, 1.2.15 [F,\RL{f}], 74 \textit{folii}, abîmé, date d'environ 1325 [=1907], et ne contient que les onze premiers chapitres de la première partie\footnote{Ces informations proviennent du catalogue de D. A. King \cite{king1986}.}.

- Bibliothèque Suleymaniyye, 339 [S,\RL{s}], 70 \textit{folii}, date de 751 [=1349], selon le catalogue Internet de la bibliothèque. Saliba affirme qu'il s'agit d'une première version de la \textit{Nih\=aya}, différente de celle des autres manuscrits existant\footnote{\textit{Cf.} \cite{saliba1990} p.~36. Dans \cite{rashed1997} p.~115, Saliba indique en effet quelques différences notables concernant le chapitre sur la Lune.}.

- Jérusalem, Bibliothèque Kh\=alidiyya, n°66/5 [Q,\RL{q}]. D. A. King indique dans \cite{king1975} que ce manuscrit contient une copie de la \textit{Nih\=aya} sous le titre de \textit{Ris\=ala f{\=\i} al-hay'a al-jad{\=\i}da}.




Je n'ai eu accès qu'à A, B, C, D et E. J'ai établi le texte des onze premiers chapitres de la première partie sur la base de ces cinq manuscrits~; pour la suite, bien que m'aidant parfois de B et E, je n'ai finalement retenu que A, C et D dans l'apparat critique. Dès le début, j'ai seulement indiqué les numéros de \textit{folio} pour A, C et D en marge de mon édition.

\begin{figure}
  \begin{center}
    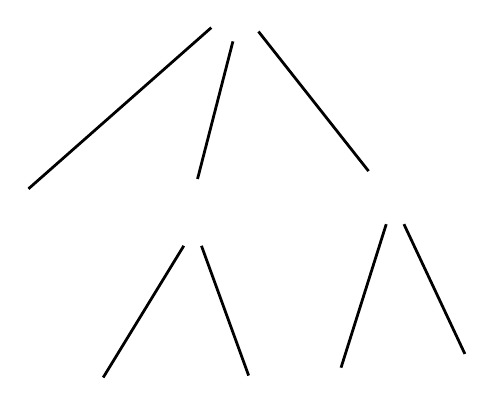
  \end{center}
  \caption{\label{stemma}\textit{stemma}}
\end{figure}

Ces cinq manuscrits ne présentent aucune différence cruciale pour l'interprétation. Pour mieux les caractériser, je m'attacherai surtout aux omissions, qui sont toujours courtes. E en commet de nombreuses\footnote{Ainsi p.~\pageref{var4} l.~14-15, p.~\pageref{var6} l.~2-3, p.~\pageref{var7} l.~16, p.~\pageref{var2} l.~7, p.~\pageref{var3} l.~12, \textit{etc.}}, parfois par saut du même au même. B présente des signes de ressemblance avec E~: certaines omissions de B sont aussi commises par E mais corrigées par des ajouts en marge de E\footnote{Ainsi p.~\pageref{var10} l.~10-11 et p.~\pageref{var11} l.~2~; mais \textit{toutes} les omissions \textit{n'ont pas} été corrigées, ainsi p.~\pageref{var1} l.~2 (un seul mot).}. Il est difficile d'en tirer une conclusion concernant la relation entre B et E car il arrive qu'une omission dans B ne soit pas commise dans E\footnote{Rarement, ainsi p.~\pageref{var12} l.~9.}. Bien que B commette beaucoup moins d'omissions que E, il présente certains signes de gaucherie qui m'ont finalement conduit à le négliger\footnote{Ainsi p.~\pageref{var27} l.~9-10 où B déplace le pronom possessif, trompant le lecteur quant à l'attribution d'une doctrine à l'auteur~; p.~\pageref{var13} l.~1-3 où B et E commettent un contre-sens (corrigé en marge de E) en inversant deux négations~; p.~\pageref{var28} l.~8 où B et E ajoutent deux mots qui nuisent au sens (dans E, ils ont ensuite été barrés)~; p.~\pageref{var25} l.~14-15, encore un contre-sens par inversion qui cette fois-ci est pire dans E, car mal corrigé. Comparées à A, C et D, les figures dans B et E sont aussi d'une qualité assez pauvre~: ainsi aux figures p.~\pageref{soleil_trajectoires} et \pageref{soleil_orbes_solides}, B et E commettent plusieurs erreurs concernant la position des orbes du Soleil dans les quadratures et les octants, et même à l'apogée dans E.}. Plusieurs omissions sont propres à C et D\footnote{Ainsi p.~\pageref{var13} l.~14-15 et p.~\pageref{var14} l.~7-8. La figure des trajectoires des centres des orbes de Mars p.~\pageref{mars_trajectoires} est aussi omise par C et D seulement.}~; certaines ont toutefois été corrigées dans C, en marge\footnote{Ainsi p.~\pageref{var17} l.~5-6. \`A la p.~\pageref{var16} l.~1-2, une longue omission par saut du même a été corrigée en marge de C~; dans D, la portion correspondante du texte a été détruite par l'incendie, mais, à en juger par l'espace disponible, il est presque certain que le copiste avait omis le même passage.}. Je n'ai pas trouvé d'omission propre à D, mais l'écriture y est parfois difficilement lisible~; c'est pourquoi j'ai aussi retenu C qui présente pourtant d'autres omissions en sus de celles partagées avec D\footnote{Ainsi p.~\pageref{var18} l.~2-4, p.~\pageref{var19} l.~15, p.~\pageref{var20} l.~5.}. De rares omissions\footnote{Ainsi p.~\pageref{var21} l. 12-13 et p.~\pageref{var22} l.~19.} propres à B, C, D et E, et non commises par A, n'ont pas suffi à me convaincre que B, C, D et E auraient un ancêtre commun non partagé avec A. Sauf une exception\footnote{Page \pageref{var23} l.~22, un mot seulement.}, je n'ai repéré aucune omission propre à A. Tous ces indices renvoient au \textit{stemma} de la figure \ref{stemma}.

D. A. King mentionne que V. Roberts avait préparé une édition et une traduction anglaise de la \textit{Nih\=aya} qui n'ont pas été publiées\footnote{\textit{Cf.} \cite{king1975}, p.~362, note 1.}. G. Saliba a indiqué, il y a trente ans, que sa propre édition critique serait bientôt prête\footnote{\textit{Cf.} \cite{saliba1987}.}, mais elle ne l'est pas encore~; quand elle verra le jour, étant l'{\oe}uvre d'un spécialiste y ayant consacré tant d'années, cette édition révèlera sans doute bien des choses qui m'ont échappé.

{\shatir} mentionne à deux reprises une ``troisième partie'' contenant des tables\footnote{\textit{Cf.} p.~\pageref{troisieme1} et \pageref{troisieme2} \textit{infra}.}, mais ceci ne figure pas dans le plan donné en introduction. Ce qui n'était peut-être qu'un projet lors de la rédaction des deux premières parties aurait pu ensuite évoluer en l'ouvrage intitulé \textit{Al-z{\=\i}j al-jad{\=\i}d}. Fuad Abbud a montré en 1962 que les tables du \textit{Z{\=\i}j al-jad{\=\i}d} peuvent être calculées au moyen des modèles de la \textit{Nih\=aya}\footnote{\textit{Cf.} \cite{ghanem1976} p.~73-80.}. Certains de ses calculs me semblent erronés, mais la conclusion semble juste. J'indiquerai dans mon commentaire le moyen de calculer les tables comme l'explique {\shatir} dans la \textit{Nih\=aya} et j'en donnerai des extraits qui se révèlent comparables aux résultats de Fuad Abbud et aux valeurs échantillonnées par lui dans le \textit{Z{\=\i}j al-jad{\=\i}d}.

\paragraph{Les concepts d'orbe et de mouvement}
En 1959, Kennedy et Roberts, citant Voltaire, enjoignaient les historiens de l'astronomie à ``cultiver leur jardin'' et ils esquissaient les linéaments d'un programme de recherche \textit{on the general topic of late medieval planetary theory} centré autour de l'{\oe}uvre d'{\shatir} et des auteurs qu'il avait lus.

Parmi ces auteurs, les collaborateurs de Na\d{s}{\=\i}r al-D{\=\i}n al-\d{T}\=us{\=\i} à l'observatoire de Maragha (XIIIème siècle) jouent un rôle prépondérant. Ces auteurs sont mieux connus aujourd'hui qu'ils ne l'étaient en 1959, spécialement grâce aux travaux de Ragep sur \d{T}\=us{\=\i} \cite{altusi1993} et de Saliba sur `Ur\d{d}{\=\i} \cite{saliba1990}, et c'est surtout par rapport à eux que j'ai essayé de situer {\shatir}.

Les modèles planétaires de ces auteurs reposant tous exclusivement sur des mouvements de rotation uniforme, je les ai tous décrits et comparés dans mon commentaire mathématique, dans un langage géométrique moderne, au moyen de rotations affines de l'espace.

Chez ces auteurs, il n'y a, à proprement parler, aucune théorie mathématique du mouvement. On ne connaît pas encore de suite aux débuts d'une théorie cinématique que constituait la \textit{Configuration des mouvements} d'Ibn al-Haytham au XIème siècle. Comparée à la \textit{Ta\b{d}kira} de Na\d{s}{\=\i}r al-D{\=\i}n al-\d{T}\=us{\=\i}, la \textit{Nih\=aya} d'{\shatir} témoigne pourtant d'une légère inflexion.

Dans l'introduction philosophique de la \textit{Nih\=aya}, {\shatir} mentionne deux doctrines.\label{mouvement} L'une classe les mouvements ainsi~:
\begin{enumerate}
\item chose mue par soi
\begin{enumerate}
\item volitif (animaux)
\item naturel
\begin{enumerate}
\item qui ne procède pas d'une manière unique (plantes)
\item qui procède d'une manière unique
\begin{enumerate}
\item rectiligne (éléments et composés)
\item circulaire (orbes)
\end{enumerate}
\end{enumerate}
\end{enumerate}
\item chose mue par autre chose
\begin{enumerate}
\item par accident
\item violent
\end{enumerate}
\end{enumerate}
L'autre doctrine conduirait à la classification suivante~:
\begin{enumerate}
\item chose mue par soi
\begin{enumerate}
\item volitif
\begin{enumerate}
\item qui ne procède pas d'une manière unique (animaux)
\item qui procède d'une manière unique (orbes)
\end{enumerate}
\item naturel
\begin{enumerate}
\item qui ne procède pas d'une manière unique (plantes)
\item qui procède d'une manière unique (éléments et composés)
\end{enumerate}
\end{enumerate}
\item chose mue par autre chose
\begin{enumerate}
\item par accident
\item violent
\end{enumerate}
\end{enumerate}
Al-\d{T}\=us\={\i} a suivi cette seconde doctrine 
quand il classe les mouvements ainsi\footnote{\textit{Cf.}
  \cite{altusi1993} p.~100.}~:
\begin{enumerate}
\item chose mue par soi
\begin{enumerate}
\item qui ne procède pas d'une manière unique, et qui a une âme
\begin{enumerate}
\item animaux
\item plantes
\end{enumerate}
\item qui procède d'une manière unique, et qui a une nature
\begin{enumerate}
\item naturel (éléments)
\item volitif (orbes)
\end{enumerate}
\end{enumerate}
\item chose mue par autre chose
\begin{enumerate}
\item par accident
\item violent
\end{enumerate}
\end{enumerate}
Le mouvement d'un orbe est-il naturel ou volitif~? C'est seulement sur
ce point que les deux doctrines diffèrent. Ibn al-\v{S}\=a\d{t}ir ne
tranche pas entièrement cette question~; il semble opter pour le
mouvement naturel, mais s'écarte probablement des autres tenants de
cette doctrine en affirmant que le mouvement des astres est
<<~simple-composé~>> et que les orbes eux-mêmes sont des corps
<<~composés~>>.

\label{simple_compo} Le mouvement individuel d'un orbe est un mouvement
circulaire uniforme, donc simple en soi, mais composé avec le
mouvement des orbes qui l'entraînent. Ainsi ce mouvement présente, du
point de vue de l'observateur (qui ne voit jamais qu'un seul
mouvement), un caractère <<~variable~>> et non uniforme, comme le dit
Ibn al-\v{S}\=a\d{t}ir. C'est, en tant que tel, un mouvement composé~;
mais il faut se garder de croire que, dans cette physique, la
composition des \emph{mouvements} ressemblerait à l'immixtion des
\emph{éléments simples} dans les corps composés\footnote{Selon
  Aristote, dans un mixte composé de deux corps, ``le résultat du
  mélange étant en acte autre que [les corps qui ont été soumis au
    mélange], mais étant encore, en puissance, l'un et l'autre [...],
  les mélanges proviennent manifestement d'éléments antérieurement
  séparés et pouvant se séparer de nouveau'' (\textit{De la génération
    et de la corruption} I 10, \textit{cf.}  \cite{aristote2005}
  p.~47).}. Il est impossible à l'observateur de distinguer ces deux
phénomènes dans le mouvement~: le mouvement simple (naturel ou
volitif) propre à l'orbe, et les mouvements (par accident) que lui
communiquent les orbes qui l'entraînent. On peut donc bien parler d'un
mouvement à la fois simple et composé\footnote{Ibn al-\v{S}\=a\d{t}ir
  insiste sur l'impuissance phénoménologique de l'observateur dans le
  sixième chapitre, p.~\pageref{impuissance} \textit{infra}.}. Chaque
mouvement simple, en tant que composé avec d'autres mouvements
communiqués au mobile <<~par accident~>>, devient ici une entité
nouvelle, un mouvement encore simple, mais <<~simple-composé~>> (et
non plus circulaire uniforme), et c'est ce mouvement dont
l'observateur est témoin. Seule la description mathématique permettra
de distinguer les composantes des mouvements des astres.

Si le concept de mouvement ``simple-composé'' gagne en réalité, c'est
d'ailleurs au détriment du concept d'``orbe'' qui semble perdre en
substance. L'astronome se garde bien de trancher sur l'essence des
orbes, et le texte semble même souligner l'ambiguïté du concept. Le
mot \textit{falak} est employé dans trois sens distincts~: c'est
(rarement) \textit{les cieux} c'est-à-dire l'ensemble des corps
célestes vus d'ici-bas\footnote{Ainsi p.~\pageref{cieux}
  \textit{infra}.}, c'est bien plus souvent \textit{un corps solide} à
symétrie sphérique, et parfois simplement \textit{un cercle
  géométrique} trajectoire du centre d'un tel corps solide au sein
d'un autre corps solide. Chaque chapitre contenant un modèle
planétaire utilise les deux dernières acceptions avec en général deux
figures, l'une représentant les sections des corps solides sphériques,
l'autre représentant les trajectoires des centres des
sphères\footnote{\textit{Cf.} par exemple les figures
  p.~\pageref{soleil_trajectoires} et \pageref{soleil_orbes_solides}
  \textit{infra} pour le Soleil.}. Le caractère impénétrable des orbes
solides impose certaines contraintes cosmologiques qu'{\shatir} ne se
prive pas de mettre à profit\footnote{Dans la conclusion de la
  première partie, p.~\pageref{conclusion} et suivantes.} pour
déterminer la taille minimale du Monde et l'ordre des planètes~; mais
si le concept d'orbe solide semble avoir sa préférence, c'est surtout
parce qu'il joue le rôle d'un \emph{référentiel solide en mouvement},
rôle essentiel dans les chapitres sur les mouvements en latitude par
exemple\footnote{\textit{Cf.} les chapitres 24 et 25,
  p.~\pageref{lat_debut}-\pageref{lat_fin}, et mon commentaire
  p.~\pageref{preliminaires}.}.

\paragraph{Les excentriques}
Si les savants de Maragha semblent avoir ignoré la naissance d'une cinématique célestre entre les mains d'Ibn al-Haytham, ils ont pourtant poursuivi la remise en cause de l'astronomie de l'\textit{Almageste} que celui-ci avait proposé dans un autre livre, ses \textit{Doutes sur Ptolémée}. La question des excentriques, du point équant et du point de prosneuse était le principal chef d'accusation contre l'astronomie de Ptolémée.

Concernant les excentriques, la position d'{\shatir} semble être encore plus sévère que celle de ses prédécesseurs~: il refuse catégoriquement les excentriques, mais ce refus se révèle plus subtil qu'il n'y paraît. Tentons de mieux cerner le problème en raisonnant d'abord par élimination.

Pour simplifier, je me limiterai ici aux mouvements en longitude. Un orbe est alors un cercle portant un point, par exemple une planète $P$, qui suit ainsi une trajectoire circulaire autour du centre de l'orbe. Soit $O$ le centre du Monde et $O_1$ le centre d'un orbe ``excentrique'' portant un point $P$. Alors $O_1\neq O$. On notera aussi, sur ce cercle, l'apogée $A$ dans la direction du point $O$. Je distinguerai quatre configurations possibles.

\begin{figure}
  \begin{center}
    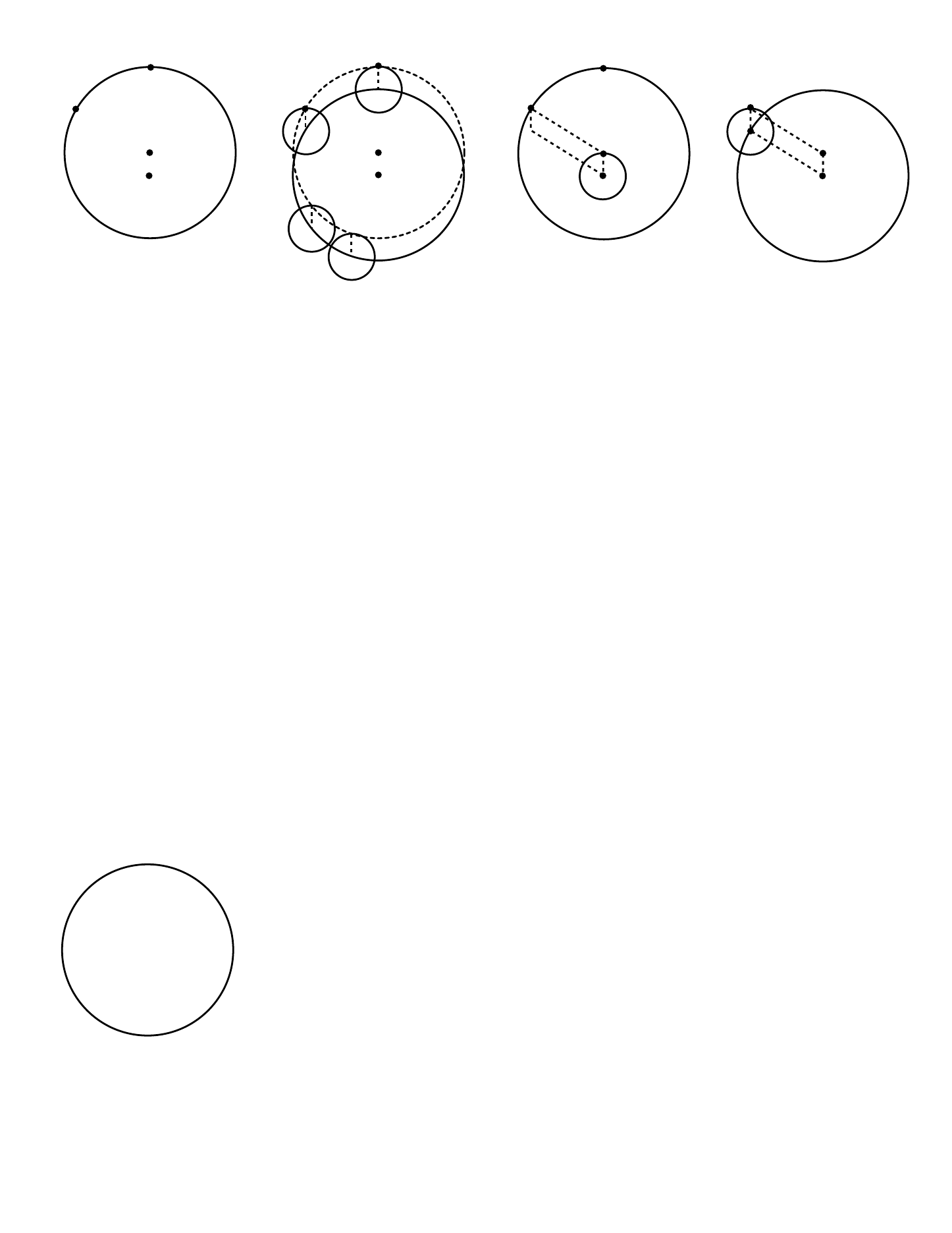
    \caption{\label{fig002}Les apories de l'orbe excentrique}
  \end{center}
\end{figure}

\emph{Cas 1}. $O_1$ et $A$ sont immobiles, et le mouvement de $P$ le long du cercle de centre $O_1$ est à vitesse angulaire constante par rapport à $O$ (\textit{cf.} figure \ref{fig002}(i)). 

\emph{Cas 2}. $O_1$ et $A$ sont immobiles, et le mouvement de $P$ le long du cercle de centre $O_1$ est à vitesse angulaire constante par rapport à $O_1$.

\emph{Cas 3}. $O_1$ est mu autour de $O$, et ce mouvement entraîne tout l'orbe excentrique et l'apogée $A$. La planète $P$ est mue le long du cercle de centre $O_1$, mais sa vitesse angulaire par rapport à $O_1$ n'est pas constante. \textit{Cf.} fig.~\ref{fig002}(iii).

\emph{Cas 4}. $O_1$ est mu autour de $O$, et ce mouvement entraîne tout l'orbe excentrique et l'apogée $A$. La planète $P$ est mue le long du cercle de centre $O_1$, à vitesse angulaire constante par rapport à $O_1$.

Les astronomes de Maragha puis {\shatir} refusent d'emblée les cas 1 et 3 qui contiennent un mouvement circulaire \textit{non uniforme} par rapport au centre du cercle~: et {\shatir} de nous dire qu'à ce compte, autant imaginer pour chaque planète un unique orbe centré en $O$ et l'emportant dans son mouvement irrégulier, avec ses accélérations, rétrogradations, \textit{etc.}\footnote{\textit{cf.} p.~\pageref{doutes_excentrique} \textit{supra}.}. Si l'on admettait des mouvements circulaires non uniformes, il ne serait plus besoin d'orbes multiples, d'excentriques, ni d'épicycles : un orbe par planète suffirait ! C'est pour cette raison purement méthodologique qu'{\shatir} refusera systématiquement tout mouvement circulaire non uniforme.

Qu'en est-il du cas 2~? Si l'on fait abstraction du mouvement diurne\footnote{qui entraîne le point $O_1$ et l'apogée $A$ en un tour toutes les vingt-quatre heures autour de $O$.}, Ptolémée concevait ainsi l'excentrique portant le Soleil qu'il jugeait exempt du mouvement de précession. Il n'est pas question ici de mouvement circulaire non uniforme. Comme les astronomes arabes, depuis le IXème siècle, imposent aussi au Soleil le mouvement de précession et le mouvement de l'Apogée, ce cas 2 n'est plus présent dans les textes d'astronomie, et il est donc difficile de deviner les objections qu'{\shatir} aurait adressées à une telle configuration. L'exercice en vaut pourtant la peine. Il me semble que l'astronome d'alors aurait pu répondre en invoquant le premier principe physique mentionné p.~\pageref{pas_de_repos} dans l'introduction de la \textit{Nih\=aya}~:
\begin{quote}
  ``\emph{Premier principe.} Le vide est impossible, et dans les orbes, le
  repos est impossible.''
\end{quote}
Il aurait fait objection que l'orbe portant l'orbe excentrique ne peut être immobile. L'existence même d'un orbe \textit{fixé} dans une position dissymétrique par rapport au centre du Monde est une brèche à la ``perfection de la circularité''. Cette brèche ne peut être compensée qu'en rompant cette fixité, en faisant l'Apogée et tout l'excentrique se mouvoir d'un mouvement de révolution autour du centre du Monde: on rétablit ainsi une sorte de symétrie dynamique. On élimine ainsi le cas 2. D'ailleurs, le cas 1 déjà éliminé tombe aussi sous cette objection.

Toutefois, une telle configuration avec un excentrique où l'apogée est immobile, et le point $P$ en révolution uniforme par rapport à $O_1$, est cinématiquement équivalente, quant à la trajectoire de $P$, à un modèle à deux orbes dont l'un, centré en $O$, entraîne un autre petit orbe portant $P$ : il suffit que le petit soit animé d'un mouvement de rotation sur lui-même dont la vitesse angulaire est l'opposé de la vitesse angulaire du grand orbe\footnote{\textit{Cf.} figure \ref{fig002}(ii) p.~\pageref{fig002}, et prop. 1 p.~\pageref{apollonius} \textit{infra}.}. {\shatir} et tous les astronomes depuis Ptolémée en étaient conscients. Ne faisant intervenir ni mouvement circulaire non uniforme, ni excentrique, ni repos, un tel modèle serait parfaitement légitime au yeux de l'astronome d'alors.

Qu'en est-il du cas 4 (fig.~\ref{fig002}(iii))~? Cette configuration ressemble beaucoup à la configuration à deux orbes, déférent et épicycle, que l'on vient de décrire~: un orbe meut le point $O_1$ circulairement autour de $O$, et un autre orbe meut le point $P$ circulairement autour de $O_1$. On parle d'''épicycle'' quand le rayon de l'orbe de $P$ est suffisamment petit pour que cet orbe ne contienne pas le point $O$, et on parle d'``excentrique'' dans le cas contraire. 

Voici une objection possible au cas 4. Le point $O$ tombe \textit{à l'intérieur} de l'orbe de $P$, mais l'orbe de $P$ doit être conçu comme un corps solide dès lors qu'un autre orbe l'entraîne dans un mouvement circulaire uniforme autour de $O$. Ce corps solide, fait d'éther, doit cependant enclore une cavité puisqu'au centre du Monde se trouve peut-être d'autres orbes -- au moins la Terre. Cette cavité doit être centrée en $O$ -- si elle était centrée en $O_1$, il y aurait du vide en son sein, puisque les orbes inférieurs et la Terre sont à symétrie sphérique de centre $O$. N'étant donc pas centrée en $O_1$ mais en $O$, cette cavité semble empêcher tout mouvement de rotation de l'orbe de $P$ sur elle-même au sein du référentiel constitué par l'orbe portant l'orbe de $P$ (\textit{cf.} figure \ref{fig002}(v))~; car la substance de l'éther était certainement conçue comme impénétrable. Chez les prédecesseurs d'{\shatir}, la solution était d'imaginer un cavité centrée en $O_1$, et de remplir le vide intermédiaire au moyen d'un autre orbe, ou d'une autre composante connexe de l'orbe portant l'orbe excentrique (\textit{cf.} fig.~\ref{fig002}(vi)), compliquant ainsi davantage la hiérarchie des corps célestes\footnote{D'ailleurs, si le repos est impossible dans les orbes, comment expliquer que ces deux composantes connexes sont immobiles l'une par rapport à l'autre~?}. Mais {\shatir} ne semble jamais s'inquiéter d'un tel problème. 

Si l'orbe de $P$ est un ``excentrique'', on peut échanger les rayons des deux orbes et obtenir ainsi une configuration cinématiquement équivalente quant à la trajectoire du point $P$, mais dans laquelle l'orbe de $P$ est un petit orbe de rayon inférieur au rayon de l'orbe le portant: il ne contient plus le point $O$ (\textit{cf.} figure \ref{fig002}(iii) et (iv)). En général, c'est la situation de la fig. \ref{fig002}(iv) qu'on rencontrera dans les modèles d'{\shatir}.

Pourquoi {\shatir} refuse-t-il donc l'excentrique du cas 4 ? Aussi étrange que cela puisse paraître, dans la \textit{Nih\=aya}, on ne trouve aucune raison générale justifiant un tel refus par un recours aux principes de la philosophie naturelle.

Résumons le problème. On rencontrait le cas 3 dans tous les modèles planétaires de Ptolémée\footnote{Sauf le Soleil mentionné dans le cas 1.}. Les astronomes de Maragha étaient parvenus à remplacer, dans chaque modèle ptoléméen, l'excentrique de la fig. \ref{fig002}(iii), dont le mouvement n'est pas uniforme par rapport à $O_1$, par un déférent en rotation uniforme par rapport à son centre. Pour ce faire, il faut changer la position de $O_1$, adjoindre un épicycle, \textit{etc.}~: je décrirai ces modèles dans mon commentaire, ils relèvent en général du cas 4. {\shatir} quant à lui préfère toujours une solution usant de la figure \ref{fig002}(iv) \textit{sans excentrique}, même quand elle est cinématiquement équivalente à la figure \ref{fig002}(iii). Pourquoi~?

G. Saliba a qualifié une telle attitude de ``purisme''\footnote{\textit{Cf.} \cite{saliba1990} p.~66.}. Ragep a récemment remarqué que cette attitude conduit, dans le cas de Mercure et Vénus, à des modèles où le centre de l'épicycle est, vu de la Terre, dans la direction du Soleil moyen, et que ce ``biais héliocentriste'' des modèles d'{\shatir} rendait davantage plausible l'hypothèse d'une transmission entre l'astronomie arabe et la science de Copernic\footnote{\textit{Cf.} \cite{ragep2016}. Sur cette thématique du rapport à Copernic qui a peu occupé mon attention une longue littérature existe déjà, et cet article en constitue certainement le meilleur point d'entrée. Les faits importants sont dans Swerdlow-Neugebauer \cite{swerdlow1984}~:

  - le modèle des planètes supérieures du \textit{Commentariolus} de Copernic est identique au modèle d'{\shatir} après changement d'origine

  - le modèle des planètes supérieures du \textit{De revolutionibus} est identique à celui de `Ur\d{d}{\=\i} (lui-même cinématiquement équivalent au modèle d'{\shatir}

  - les modèles des planètes inférieures dans le \textit{De revolutionibus} sont une adaptation facile de ceux d'{\shatir}

  - le modèle de la Lune du \textit{Commentariolus} et du \textit{De revolutionibus} est identique à celui d'{\shatir} (pas besoin de changer d'origine puisque la Lune tourne autour de la Terre)

Ragep montre que les modèles des planètes inférieures du \textit{Commentariolus} sont aussi une adaptation facile et cinématiquement équivalente des modèles d'{\shatir} (\textit{cf.} \cite{ragep2016} p.~402 pour Vénus et p.~404-406 pour Mercure)~; enfin, il remarque que le ``biais héliocentriste'' des modèles d'{\shatir} a pu faciliter le passage à l'héliocentrisme de Copernic.}; mais il concède lui-même que ce jugement \textit{a posteriori} ne peut valoir comme explication d'une telle préférence dans l'esprit du savant du quatorzième siècle.

Qu'il me soit permis de formuler une conjecture~: \emph{quand {\shatir} reproche aux savants de Maragha leurs modèles usant de la fig. \ref{fig002}(iii), c'est en fait la \emph{démonstration} de ces modèles qui est l'objet du reproche.}

Quand les auteurs de Maragha démontrent, ils le font en prouvant que le mouvement décrit par l'astre dans leur modèle est \emph{proche} du mouvement décrit par l'astre \emph{dans le modèle de Ptolémée}. Ils le font au moyen de figures géométriques. Cette technique a d'ailleurs été reprise par les commentateurs modernes à partir de Kennedy\footnote{\textit{Cf.} la figure dans \cite{ghanem1976} p.~88~; je l'ai reproduite en fig.~\ref{fig051} p.~\pageref{fig051} \textit{infra}. Je ne me priverai pas de tels outils que j'ai résumés dans les propositions 1, 2 et 3 de mon commentaire p.~\pageref{apollonius}~; mais tandis que Kennedy et Roberts y voyait la commutativité de l'addition vectorielle, parfaitement comprise par Copernic et peut-être en gestation chez des savants comme {\shatir} (\cite{ghanem1976} p.~60), il me semble que les travaux de ce dernier se situent à un autre niveau. Il s'agit de \emph{composer des rotations affines spatiales}, geste autrement plus avancé qu'une simple addition vectorielle. Des recherches récentes ont montré qu'il n'est pas incongru de parler de transformations géométriques dans les mathématiques médiévales arabes dès le IXème siècle (\textit{cf.} l'usage des projections cylindrique et conique dans le \textit{Traité sur l'astrolabe} d'al-Q\=uh{\=\i} étudié dans \cite{rashed2005} et \cite{abgrall2004}~; plus généralement, voir l'introduction de \cite{rashed2004}).} dans leurs efforts pour comprendre la genèse des modèles planétaires.

Le point de départ heuristique et le garant final de vérité étaient le modèle géométrique ptoléméen dont les paramètres seulement étaient adaptés aux observations récentes. Ainsi Na\d{s}{\=\i}r al-D{\=\i}n al-\d{T}\=us{\=\i}, pour démontrer son nouveau modèle pour la Lune\footnote{C'est un modèle qui n'utilise pourtant aucun orbe excentrique au sens ci-dessus, \textit{cf.} fig.~\ref{fig021} p.~\pageref{fig021} \textit{infra}.} prouve que la trajectoire du centre de l'épicycle y est proche du cercle excentrique de Ptolémée~: et {\shatir} lui reproche alors, avec raison, d'avoir conservé le concept d'excentrique !

Si cette conjecture est exacte, le projet d'{\shatir} était celui d'une reconstruction globale de l'astronomie, \emph{en partant des données de l'observation}, et sans plus avoir recours aux modèles de l'\textit{Almageste}. S'explique alors qu'il renvoie à un ouvrage intitulé \textit{Commentaire des observations} pour la démonstrations de ses modèles. Le scénario imaginé par Saliba \cite{saliba1987} pour déduire le modèle du Soleil d'{\shatir} me semble appuyer cette hypothèse. Le sens qu'{\shatir} donne au mot ``observation'' n'est d'ailleurs pas ambigu~: ses usages dans la \textit{Nih\=aya} confirment qu'il s'agit toujours d'observer les positions des astres à une date donnée, ou bien leur diamètre apparent, ou encore des durées (d'une éclipse). C'est à peu près le concept moderne d'observation astronomique, sans toutefois la théorie des erreurs ni aucune réflexion approfondie sur la précision des mesures~: il faudra attendre bien des siècles pour voir éclore de tels outils.

Cette conjecture que certains jugeront peut-être audacieuse s'est graduellement imposée à moi~; elle devrait gagner en plausibilité à la lecture du texte et du commentaire~; hélas, seul l'ouvrage perdu \textit{Commentaire des observations} pourrait l'infirmer ou l'étayer pleinement.

\begin{figure}
  \begin{center}
    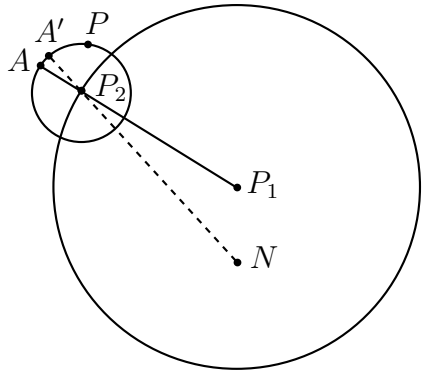
    \caption{\label{fig013}Quand l'apogée de l'épicycle suit un point autre que le centre de l'orbe qui le porte}
  \end{center}
\end{figure}

\paragraph{Origine du mouvement de l'épicycle} Si l'excentrique est le principal chef d'accusation contre l'astronomie ptoléméenne, il y a un autre thème qui revient souvent dans les doutes exprimés par {\shatir} : le fait que ``l'apogée de l'épicycle suive un point autre que le centre de l'orbe qui le porte''. Considérons un modèle contenant un orbe de centre $P_1$ portant un épicycle de centre $P_2$, et l'épicycle portant à son tour un astre en $P$ (\textit{cf.} fig.~\ref{fig013} p.~\pageref{fig013}). Les principes adoptés par {\shatir} imposent au point $P_2$ d'être immobile au sein de l'orbe centré en $P_1$. Bien sûr l'orbe centré en $P_1$ est lui-même en rotation uniforme autour d'un axe passant par $P_1$, et ce mouvement entraîne $P_2$ et l'épicycle. De même, l'épicycle est animé d'un mouvement de rotation autour d'un axe passant par $P_2$, et ce mouvement est uniforme \emph{par rapport au référentiel constitué par l'orbe centré en $P_1$}. Notons $A$ l'apogée de l'épicycle sur la droite $(P_1P_2)$ ; comme la direction $(P_1P_2)$ est immobile au sein de l'orbe centré en $P_1$, alors l'angle $AP_2P$ doit croître uniformément. Mais il en est rarement ainsi dans les dispositifs proposés depuis Ptolémée pour sauver les phénomènes ! En général, l'origine $A'$ du mouvement uniforme de $P$ autour de $P_2$ oscille autour de $A$ de sorte à rester alignée avec $P_2$ et un autre point $N$ distinct de $P_1$. Dans le deuxième modèle de la Lune dans l'\textit{Almageste}, $N$ est le centre du Monde, distinct du centre $P_1$ de l'excentrique portant l'épicycle\footnote{\textit{Cf.} \cite{pedersen1974} p.~186-187.}. Dans le troisième modèle de la Lune, $N$ est le point de \textit{prosneuse}, distinct du centre de l'excentrique et du centre du Monde. Dans les modèles planétaires de l'\textit{Almageste}, $N$ est le point \emph{équant}, distinct du centre de l'excentrique et du centre du Monde.

On retrouve ici à peu près la même attitude que vis-à-vis des excentriques chez notre auteur. Par exemple, quand al-\d{T}\=us{\=\i} conçoit un modèle à base de ``couple de \d{T}\=us{\=\i} curviligne'' pour éliminer le problème du point de prosneuse de la Lune, {\shatir} lui reproche précisément d'avoir recours à l'excentrique \emph{et} au point de prosneuse~! Mais là encore, la solution d'al-\d{T}\=us{\=\i} consiste à offrir un modèle sans expliquer sa genèse, puis à montrer que son comportement est cinématiquement équivalent à celui de Ptolémée, au moins approximativement. Ce qui gêne {\shatir}, n'est-ce-pas cette méthode démonstrative, en l'absence d'un exposé déduisant le modèle à partir de l'observation, dans un traité d'astronomie~?

\paragraph{La possibilité mathématique d'une nouvelle astronomie}
Si tel était le projet d'{\shatir}, il ne s'agissait plus d'imaginer un modèle reproduisant les effets d'une théorie géométrique transmise depuis douze siècles pour résoudre enfin le problème de son insertion dans un édifice philosophique contraignant. Il s'agissait d'imaginer une méthode universelle conduisant directement des données de l'observation continue, toujours susceptibles d'être révisées, à un modèle dont les éléments étaient choisis dans un ensemble restreint (rotations affines spatiales) mais avec une liberté de combinaison que les mathématiques permettaient d'embrasser depuis seulement quelques siècles, en composant des rotations non homocentriques et d'axes pas nécessairement parallèles. On ne peut éluder la question de la possibilité mathématique de la réussite d'un tel projet~; {\shatir} lui-même a dû se la poser. Faute de source, je suis réduit à formuler cette question dans un langage moderne. Tout mouvement curviligne, suffisamment régulier, dans l'espace affine $\mathbb{R}^3$, peut-il être décrit par une composée de rotations uniformes~?

Il ne faut certes pas se limiter aux mouvements plans, car les modèles les plus intéressants conçus par les savants de Maragha et par {\shatir} décrivent aussi les latitudes des astres\footnote{Et certes il ne faut pas non plus oublier la variation des distances des astres à la Terre. \`A la rigueur, un modèle d'Univers formé de sphères toutes homocentriques, en rotation uniforme les unes par rapport aux autres, pourrait suffire à décrire les longitudes et les latitudes des astres. Duhem s'était enthousiasmé de ce phénomène dans son analyse de l'astronomie d'Eudoxe et de la Physique péripatéticienne, \textit{cf.} \cite{duhem1913} vol.~I p.~126-129 et vol.~II p.~42-43~; mais dans les modèles purement homocentriques, la distance de chaque astre à la Terre est inévitablement constante.}. J'ai montré dans \cite{penchevre2016} comment le modèle d'{\shatir} pour Vénus rend très bien compte des mouvements en longitude \emph{et} en latitude de l'astre. Il en est de même pour les autres planètes\footnote{\textit{Cf.} commentaire mathématique \textit{infra}, p.~\pageref{begin_saturne}--\pageref{end_saturne} pour les planètes supérieures, p.~\pageref{begin_mercure}--\pageref{end_mercure} pour Mercure.}. Les prédictions en latitude ne sont pas meilleures, d'un point de vue moderne, que celles de Ptolémée, car {\shatir} s'appuyait sur les données phénoménologiques recueillies par Ptolémée et non sur de nouvelles observations~; mais il s'agit bien d'un modèle composé exclusivement de mouvements de rotation uniforme\footnote{Dans le commentaire mathématique, je montre ce qu'{\shatir} doit à ses prédécesseurs Ibn al-Haytham et al-\d{T}\=us{\=\i} concernant les latitudes.}.

Soit donc un repère spatial d'origine $O$ et un point matériel dont la position $P'$ au temps $t$ sera décrite par le vecteur $\overrightarrow{OP'}$. Les méthodes modernes de l'analyse fonctionnelle et la transformée de Fourier donnent une description du mouvement coordonnée par coordonnée comme somme de fonctions harmoniques, d'où une expression de la forme~:
$$\overrightarrow{OP'}=\int(\cos\omega t.\mathbf{u}_\omega+\sin\omega t.\mathbf{v}_\omega)\,\text{d}\omega.$$
Sous l'hypothèse que le mouvement est suffisamment est régulier, on doit pouvoir approcher une telle expression au moyen d'une somme discrète, voire même finie. Posons donc~:
$$\overrightarrow{OP'}=\sum_{k=1}^N\cos\alpha_k.\mathbf{u}_k+\sin\alpha_k.\mathbf{v}_k\qquad (\star)$$
où les $\mathbf{u}_k$, $\mathbf{v}_k$ sont des vecteurs constants et les $\alpha_k$ sont des fonctions affines de $t$, de la forme $\omega t+\phi$.
Peu importe que les fréquences soient toutes, ou non, des multiples entiers d'une même fréquence fondamentale~: c'est le cas si le mouvement est périodique, mais il n'est guère évident que les mouvements célestes le soient.
On notera
$$\overrightarrow{OP}=\sum_{k=1}^N\mathbf{u}_k$$
On va à présent essayer de décrire $P'$ comme étant l'image du point $P$ par une composée de rotations affines.

Remarquons d'abord que l'extrémité de chaque vecteur de la forme ${\cos\alpha.\mathbf{u}+\sin\alpha.\mathbf{v}}$ décrit une ellipse dans le plan $(\mathbf{u},\mathbf{v})$, et que cette ellipse est un cercle si et seulement si $\Vert\mathbf{u}\Vert=\Vert\mathbf{v}\Vert$ et $\mathbf{u}\perp\mathbf{v}$. Soit $(\mathbf{i},\mathbf{j})$ une base orthonormée du plan vectoriel $(\mathbf{u},\mathbf{v})$. Notons~:
$$\mathbf{u}=u_x\mathbf{i}+u_y\mathbf{j},\qquad
\mathbf{v}=v_x\mathbf{i}+v_y\mathbf{j}.$$
Posons~:
\begin{align*}
  \mathbf{u}_1&=\dfrac{u_x+v_y}{2}\mathbf{i}+\dfrac{u_y-v_x}{2}\mathbf{j},\\
  \mathbf{v}_1&=\dfrac{-u_y+v_x}{2}\mathbf{i}+\dfrac{u_x+v_y}{2}\mathbf{j},\\
  \mathbf{u}_2&=\dfrac{u_x-v_y}{2}\mathbf{i}+\dfrac{u_y+v_x}{2}\mathbf{j},\\
  \mathbf{v}_2&=\dfrac{u_y+v_x}{2}\mathbf{i}+\dfrac{-u_x+v_y}{2}\mathbf{j}.
\end{align*}
Il est alors facile de vérifier que~:
$$\cos\alpha.\mathbf{u}+\sin\alpha.\mathbf{v}
=(\cos\alpha.\mathbf{u}_1+\sin\alpha.\mathbf{v}_1)
+(\cos\alpha.\mathbf{u}_2+\sin\alpha.\mathbf{v}_2),$$
et $\Vert\mathbf{u}_1\Vert=\Vert\mathbf{v}_1\Vert$, $\mathbf{u}_1\perp\mathbf{v}_1$, $\Vert\mathbf{u}_2\Vert=\Vert\mathbf{v}_2\Vert$, $\mathbf{u}_2\perp\mathbf{v}_2$~; les extrémités des deux vecteurs $(\cos\alpha.\mathbf{u}_1+\sin\alpha.\mathbf{v}_1)$ et $(\cos\alpha.\mathbf{u}_2+\sin\alpha.\mathbf{v}_2)$ décrivent donc des cercles à vitesse angulaire constante dans le plan $(\mathbf{u},\mathbf{v})$. Grâce à cette récriture, on peut donc supposer que, dans $(\star)$, pour tout $k$, $\Vert\mathbf{u}_k\Vert=\Vert\mathbf{v}_k\Vert$ et $\mathbf{u}_k\perp\mathbf{v}_k$.

Notons $\mathbf{w}_k=\mathbf{u}_k\wedge\mathbf{v}_k$. Je désignerai par $R_{\alpha_k,\mathbf{w}_k}$ la rotation \emph{vectorielle} d'angle $\alpha_k$ autour du vecteur $\mathbf{w}_k$~; l'équation $(\star)$ devient alors~:
$$\overrightarrow{OP'}=\sum_{k=1}^NR_{\alpha_k,\mathbf{w}_k}(\mathbf{u}_k),$$
où seuls les $\alpha_k$ dépendent de $t$. Notons $P_1$, $P_2$,..., $P_{N+1}$ les points définis par~:
$$O=P_1,\quad \overrightarrow{P_kP_{k+1}}=\mathbf{u}_k\text{ pour }1\leq k\leq N.$$
En particulier $P=P_{N+1}$. Je désignerai par $R_{\alpha_k,\mathbf{w}_k,P_k}$ la rotation \emph{affine} d'angle $\alpha_k$ autour de la droite de vecteur directeur $\mathbf{w}_k$ passant par $P_k$. On va démontrer que $(\star)$ entraîne~:
$$P'=R_{\alpha_1,\mathbf{w}_1,P_1}\left(R_{-\alpha_1,\mathbf{w}_1,P_2}\,R_{\alpha_2,\mathbf{w}_2,P_2}\right)\left(R_{-\alpha_2,\mathbf{w}_2,P_3}\,R_{\alpha_3,\mathbf{w}_3,P_3}\right)...\left(R_{-\alpha_{N-1},\mathbf{w}_{N-1},P_N}\,R_{\alpha_N,\mathbf{w}_N,P_N}\right)(P).$$
Bien que cette représentation ne soit pas unique, on tient là une description du mouvement comme composée de rotations uniformes.

\emph{Démonstration.} On posera $\alpha_0=0$ de sorte que $R_{\alpha_0,\mathbf{w}_0,P_1}=\text{id}$. Soit $P^{(N+1)}=P$, et on définiera des points $P^{(1)}$, $P^{(2)}$,..., $P^{(N)}$ par la relation suivante pour $1\leq k\leq N$~:
$$P^{(k)}=R_{-\alpha_{k-1},\mathbf{w}_{k-1},P_k}\,R_{\alpha_k,\mathbf{w}_k,P_k}(P^{(k+1)}).$$
En termes de rotations vectorielles, ceci revient à poser~:
$$\overrightarrow{P_kP^{(k)}}=R_{-\alpha_{k-1},\mathbf{w}_{k-1}}R_{\alpha_k,\mathbf{w}_k}(\overrightarrow{P_kP^{(k+1)}}).$$
Reste à démontrer que~:
$$\overrightarrow{P_1P^{(1)}}=\sum_{k=1}^NR_{\alpha_k,\mathbf{w}_k}(\mathbf{u}_k).$$
En fait, on va démontrer, par récurrence descendante sur $k$, que pour tout $k\leq N$~:
$$\overrightarrow{P_kP^{(k)}}=R_{-\alpha_{k-1},\mathbf{w}_{k-1}}\left(\sum_{m=k}^NR_{\alpha_m,\mathbf{w}_m}(\mathbf{u}_m)\right).$$
Pour $k=N$, c'est vrai par définition de $P^{(N)}$. Supposons que l'égalité est vérifiée au rang $(k+1)$, alors~:
$$\overrightarrow{P_{k+1}P^{(k+1)}}=R_{-\alpha_k,\mathbf{w}_k}\left(\sum_{m=k+1}^NR_{\alpha_m,\mathbf{w}_m}(\mathbf{u}_m)\right).$$
Au rang $k$, on aura~:
\begin{align*}
  \overrightarrow{P_kP^{(k)}}
  &=R_{-\alpha_{k-1},\mathbf{w}_{k-1}}\,R_{\alpha_k,\mathbf{w}_k}\left(\overrightarrow{P_kP_{k+1}}+\overrightarrow{P_{k+1}P^{(k+1)}}\right)\\
  &=R_{-\alpha_{k-1},\mathbf{w}_{k-1}}\,R_{\alpha_k,\mathbf{w}_k}\left(\mathbf{u}_k+R_{-\alpha_k,\mathbf{w}_k}\left(\sum_{m=k+1}^NR_{\alpha_m,\mathbf{w}_m}(\mathbf{u}_m)\right)\right)\\
  &=R_{-\alpha_{k-1},\mathbf{w}_{k-1}}\left(R_{\alpha_k,\mathbf{w}_k}\mathbf{u}_k+\sum_{m=k+1}^NR_{\alpha_m,\mathbf{w}_m}(\mathbf{u}_m)\right)\\
  &=R_{-\alpha_{k-1},\mathbf{w}_{k-1}}\left(\sum_{m=k}^NR_{\alpha_m,\mathbf{w}_m}(\mathbf{u}_m)\right),
\end{align*}
\textit{q. e. d.}

\paragraph{Le projet d'{\shatir} a-t-il réussi~?}
Il est certes impossible qu'{\shatir} ait réussi à formuler une méthode universelle pour calculer les fréquences propres et les coefficients intervenant dans une quelconque combinaison linéaire de fonctions harmoniques pour en déduire ensuite les rotations spatiales souhaitées. Nous frôlons là l'anachronisme. Pour jauger la réussite de ses recherches d'une manière plus respectueuse de l'histoire, il faut d'abord se poser la question de l'adéquation entre les prédictions de son modèle achevé et les données de l'observation, non seulement sur le plan qualitatif (rétrogradation, occurrence d'une éclipse, variation des distances à la Terre, \textit{etc.}), mais jusque sur le plan quantitatif de la précision des prédictions. Nos moyens sont maigres puisque nous ne disposons pas des observations de l'époque~: aucune observation datée n'est recueillie dans la \textit{Nih\=aya}. Si l'on retrouvait un jour le \textit{Commentaire des observations}, il ne faut pas s'attendre non plus à un recueil consignant la précision des mesures puisqu'aucune théorie des erreurs n'existait alors~: de telles données nous seraient donc de peu d'usage, à moins d'en avoir un nombre suffisant pour faire nous-mêmes des statistiques.

Pour estimer l'adéquation du modèle à l'observation, il ne nous reste guère d'autre solution que de comparer ses prédictions aux prédictions de modèles modernes dont l'erreur, même pour une époque éloignée de la nôtre d'environ sept siècles, est certainement d'un ordre de grandeur bien moindre que celle entachant les modèles d'{\shatir}. C'est ce que j'ai fait pour chaque astre~; dans le commentaire mathématique, je donne les courbes de l'équation en longitude\footnote{L'\textit{équation en longitude} désigne l'écart entre la longitude de l'astre et sa longitude \textit{moyenne}.} et de la latitude, calculées au moyen de différents modèles à comparer (les modèles de Ptolémée, ceux des savants de Maragha, ceux d'{\shatir}), sur un même système d'axes pour chaque astre, et j'y adjoins en général le tracé des données calculées par le serveur d'éphémérides de l'IMCCE\footnote{Observatoire de Paris. Ce serveur d'éphémérides utilise la théorie planétaire moderne INPOP13c qui tient compte des perturbations entre les différents corps du système solaire et des effets relativistes. L'incertitude d'INPOP13c est de l'ordre de $15$ mas sur une période contemporaine d'un siècle, selon \cite{inpop13c} p.~11. En gros, $15\text{ mas / siècle}\times 7\text{ siècles}=0,105''$ serait l'ordre de grandeur de l'incertitude d'INPOP13c pour l'époque d'{\shatir} il y a environ sept siècles. Je souhaiterais qu'un astronome de métier confirme ce raisonnement. Le cas échéant, cet ordre de grandeur est bien inférieur à l'erreur entachant les prédictions d'{\shatir}.}. J'ai tracé ces courbes en général pour une demi-révolution, une entière, ou bien quelques révolutions de l'astre autour du Soleil~: elles permettent de se faire une première idée des défauts propres à chaque modèle.

Pour estimer plus systématiquement l'incertitude des modèles d'{\shatir}, j'ai calculé l'erreur absolue, en longitude et en latitude, \emph{par rapport aux longitudes et latitude calculées par le serveur de l'IMCCE}, pour un échantillon d'environ 2000 dates, sur une période de 60 ans à partir de l'\'Epoque de référence choisi par {\shatir} (1331 ap. J.-C.)~; j'ai calculé pour chaque astre les fréquences cumulées de l'erreur absolue en degrés. Les résultats sont consignés dans le tableau \ref{tab1}.

\begin{table}
  \begin{center}
    
    \begin{tabular}{ccccccc}
      \hline
      & Soleil & Lune & \multicolumn{2}{c}{Mercure} & \multicolumn{2}{c}{Vénus}\\
      fréquence & lon. & lon. & lon. & lat. & lon. & lat. \\
      \hline
      50 \% & $<0°6'$ & $<0°35'$ & $<1°34'$ & $<0°18'$ & $<0°22'$ & $<0°15'$ \\
      70 \% & $<0°7'$ & $<0°53'$ & $<3°6'$ & $<0°44'$ & $<0°31'$ & $<0°23'$ \\
      90 \% & $<0°10'$ & $<1°18'$ & $<5°51'$ & $<1°15'$ & $<0°46'$ & $<1°6'$ \\
      95 \% & $<0°10'$ & $<1°28'$ & $<8°53'$ & $<1°32'$ & $<1°22'$ & $<1°41'$ \\
      98 \% & $<0°10'$ & $<1°37'$ & $<10°36'$ & $<2°3'$ & $<2°10'$ & $<2°15'$ \\
      \hline
    \end{tabular}

    \medskip
    
    \begin{tabular}{ccccccc}
      \hline
      & \multicolumn{2}{c}{Mars} & \multicolumn{2}{c}{Jupiter} & \multicolumn{2}{c}{Saturne}\\
      fréquence & lon. & lat. & lon. & lat. & lon. & lat. \\
      \hline
      50 \% & $<0°28'$ & $<0°26$ & $<0°10'$ & $<0°26'$ & $<0°18'$ & $<0°15'$ \\
      70 \% & $<0°45'$ & $<0°37$ & $<0°15'$ & $<0°32'$ & $<0°25'$ & $<0°18'$ \\
      90 \% & $<1°28'$ & $<0°51$ & $<0°23'$ & $<0°45'$ & $<0°36'$ & $<0°23'$ \\
      95 \% & $<1°59'$ & $<0°56$ & $<0°27'$ & $<0°50'$ & $<0°40'$ & $<0°25'$ \\
      98 \% & $<2°36'$ & $<1°2$ & $<0°31'$ & $<0°53'$ & $<0°43'$ & $<0°30'$ \\
      \hline
    \end{tabular}

  \end{center}
\caption{\label{tab1}Fréquences cumulées de l'erreur absolue en degrés de longitude et latitude des modèles d'{\shatir}.}
\end{table}

\textit{Grosso modo}\footnote{Voici un exemple pour expliquer comment lire le tableau \ref{tab1}~: sur l'échantillon considéré, 70 \% des prédictions du modèle d'{\shatir} concernant la latitude de Saturne présentent une erreur dont la valeur absolue est inférieure à $0°18'$.}, on retiendra que la précision des prédictions délivrées par les modèles d'{\shatir} était de l'ordre \emph{du degré}. C'est peu satisfaisant, mais cet ordre de grandeur nous renseigne aussi sur la précision des observations avec les instruments d'alors.

\paragraph{La deuxième partie de la \textit{Nih\=aya}} Il ne me semble pas utile de faire un commentaire détaillé de la deuxième partie de la \textit{Nih\=aya}~: celle-ci est moins originale, elle ressemble de très près à la partie III de la \textit{Ta\b{d}kira} de \d{T}\=us{\=\i} qu'{\shatir} cite abondamment. Quelques ajouts notables dont des emprunts à Théodose, en particulier des observations empiriques concernant la durée des parties du jour sous différentes latitudes. Un chapitre fait exception~: le chapitre de géographie où {\shatir} nomme les principales villes et contrées situées sous chaque climat. Je me suis beaucoup aidé de l'atlas \cite{kennedy1987} des Kennedy pour identifier les toponymes. Le \textit{Z{\=\i}j al-jad{\=\i}d} d'{\shatir} compte parmi les sources utilisées par les Kennedy~; hélas les toponymes de la \textit{Nih\=aya} ne sont pas un sous-ensemble de ceux du {Z{\=\i}j}, ni l'inverse. Il est fait mention, dans ce chapitre de la \textit{Nih\=aya}, d'une carte des fleuves et des montagnes hélas absente des manuscrits que j'ai consultés\footnote{\textit{Cf.} p.~\pageref{carte_fleuves_montagnes} \textit{infra}. Un vague schéma fait peut-être office de carte dans les manuscrits C~f.~52v et D~f.~54v ; je ne l'ai pas reproduit.}.

\paragraph{Astronomie et cosmologie} La première partie de la \textit{Nih\=aya}, si originale, se conclut sur des considérations cosmologiques, dans la lignée du \textit{Livre des Hypothèses} de Ptolémée, sur les distances des astres à la Terre et la taille du Monde. Quelle valeur donner à un tel choix~?

Ces considérations cosmologiques montrent d'abord un {\shatir} conscient de la \emph{sensibilité} de certains calculs à la précision des données de l'observation. Premier talon d'Achille de la cosmologie antique~: la méthode attribuée à Aristarque pour calculer la distance Terre-Soleil est \emph{très sensible} à la précision d'une observation pourtant \emph{si peu précise} que pouvait l'être celle du rayon de l'ombre lors des éclipses de Lune\footnote{\textit{Cf.} note~\ref{sensibilite} p.~\pageref{sensibilite} \textit{infra}}. Cela n'empêche pas {\shatir} d'en tirer un maigre indice sur l'ordre des planètes, Lune-Mercure-Vénus-Soleil-Mars-Jupiter-Saturne selon lui contrairement à ce qu'affirmait al-`Ur\d{d}{\=\i} à Maragha ; mais autre chose semble dominer son intérêt pour la question cosmologique. Profusion des distances, surfaces, volumes, vitesses calculées dans plusieurs unités, tantôt en parasanges, tantôt en rayons terrestres, tantôt par quinzième d'heure, tantôt toutes les quatre secondes\footnote{Quatre secondes~: la durée qu'il faut à un homme ``pour compter rapidement jusqu'à six'', \textit{cf.} \pageref{quatre_secondes} \textit{infra}.}~: établissement d'ordres de grandeurs. Car ce n'est que cela. L'édifice ne donne aucune certitude sur des valeurs exactes, non seulement à cause de la sensibilité des calculs, mais aussi parce que les modèles planétaires ne prescrivent qu'une borne inférieure à l'épaisseur relative de chaque système d'orbes\footnote{Voir fig.~\ref{fig002}(viii) p.~\pageref{fig002} \textit{supra}~: rien n'empêche la présence de matière superflue en deça et au delà des épicycles dans les orbes déférents.}~! {\shatir} le répète sans cesse~: ces distances ne peuvent être moindres, mais elles pourraient être plus grandes. Il en est ainsi du rayon du Monde (supérieur ou égal à $79088$ rayons terrestres) comme de la vitesse des astres les plus éloignés, les étoiles entraînées par le mouvement diurne à une vitesse d'au moins $23$ rayons terrestres toutes les quatre secondes. Un intérêt tout physique pour la mesure du Monde.

Ces considérations ne peuvent être justifiées que si elles sont compatibles avec la théorie planétaire~: or on ne peut faire justice à {\shatir} sans insister sur le fait que, pour la première fois dans l'histoire de l'astronomie, la théorie planétaire et la cosmologie sont unies en un tout cohérent. Jusqu'alors en effet, le modèle de la Lune impliquait des variations de la distance Terre-Lune aberrantes et en contradiction flagrante avec l'observation~; mais la distance Terre-Lune déduite de la parallaxe étant un autre ingrédient essentiel de la cosmologie ancienne, cette contradiction était un second talon d'Achille pour l'édifice entier. Le modèle de la Lune d'{\shatir} donne, pour la première fois, une estimation raisonnable des variations relatives de la distance Terre-Lune\footnote{\textit{Cf.} fig. \ref{fig016} p.~\pageref{fig016} \textit{infra}.}.

\paragraph{Remerciements} Chers professeurs, collègues, amis et proches qui m'avez inspiré, conseillé, rectifié, écouté et répondu~: Roshdi Rashed, Christian Houzel, Régis Morelon, Pascal Crozet, Ahmad Hasnawi, Marouane Ben Miled, Philippe Abgrall, George Saliba, Guillaume Loizelet, Valeria Candeli, Aline Auger, Lila Lamrani, Houda Ayoub, les lecteurs anonymes de mon article \cite{penchevre2016}, et Françoise mon épouse~; ce travail n'aurait pu voir le jour sans votre aide. Merci.

\part*{\'Edition et traduction}
\addcontentsline{toc}{part}{\'Edition et traduction}
\label{txt_debut}

\newpage\phantomsection 
\index{AHAWBDBEBJBHAS@\RL{b.talimayUs}, Ptolémée}
\index{ACAUBD@\RL{'a.sl}!ACAUBD@\RL{'a.sl}, 1) fondement, 2) origine}
\index{ASACBD@\RL{sa'ala}!BEASACBDAI AVAQBHAQBJAI@\RL{mas'alaT .durUriyaT}, principe}
\index{ATBCBC@\RL{^skk}!\RL{^skkuN}, doute}
\index{AQAUAO@\RL{r.sd}!AEAQAUAGAO AUAMBJAMAI@\RL{al-'ar.sAd al-.sa.hI.haT}, les observations}
\index{BGBJAC@\RL{haya'a}!BGBJACAI@\RL{hI'aT}, configuration, astronomie}
\addcontentsline{toc}{chapter}{Préface}
\includepdf[pages=1,pagecommand={\thispagestyle{plain}}]{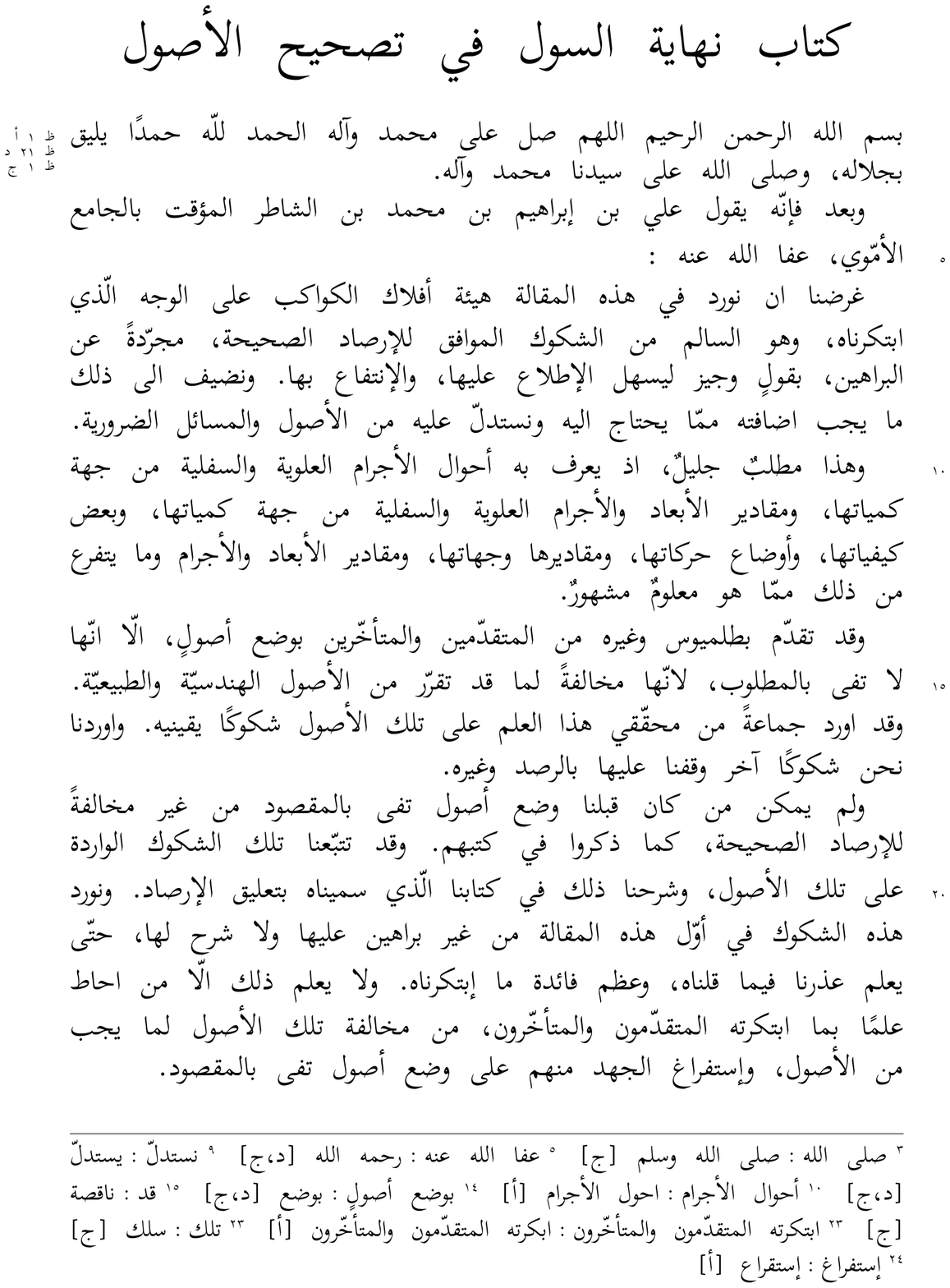}\phantomsection

\begin{center}
\LARGE L'achèvement de l'enquête et la correction des fondements
\end{center}
Au nom de Dieu, le Clément, le Miséricordieux~; louer Dieu
convient à sa Majesté~; que Dieu bénisse notre seigneur Mahomet et sa
descendance.

`Al{\=\i} b. Ibr\=ah{\=\i}m b. Mu\d{h}ammad b. al-\v{S}\=a\d{t}ir,
\emph{muwaqqit}\footnote{Qui observe et indique l'heure de la prière.}
à la Grande Mosquée des Omeyyades, dit~:

Dans ce traité, nous avons voulu exposer la configuration des orbes
des astres selon la méthode que nous avons inventée -- méthode sauve
des doutes et en accord avec les observations -- mais sans les
démonstrations, à la manière d'un abrégé, pour en simplifier
l'apprentissage et l'usage. Nous y ajoutons ce qu'il faut et nous le
déduisons des fondements et des principes.

C'est une noble intention, car elle fait connaître les états des corps
supérieurs et inférieurs (du point de vue de leurs quantités), leurs
grandeurs et leurs distances (du point de vue de leurs quantités et de
certaines de leurs qualités), les lieux de leurs mouvements, la
grandeur et le sens de ces mouvements, et ce qui, parmi les choses
connues, découle des grandeurs et des distances des corps.

Ptolémée, et d'autres parmi les anciens et les modernes, avaient
entrepris d'établir des fondements, mais ils ne suffisent pas car ils
contredisent les fondements déjà admis en géométrie et en
physique. Plusieurs spécialistes ont déjà exposé des doutes réfléchis
sur ces fondements. Nous en avons exposé d'autres que l'observation
(entre autres) nous a engagé à considérer.

Aucun de nos prédécesseurs n'a su établir de fondements suffisants
sans contredire les observations, et ils l'ont avoué dans leurs
livres.  Nous avons étudié ces doutes avec soin et nous avons expliqué
cela dans notre livre intitulé \emph{Commentaire des
  observations}\footnote{\emph{Ta`l{\=\i}q al-ar\d{s}\=ad}, ouvrage
  perdu d'Ibn al-\v{S}\=a\d{t}ir.}. Nous exposerons aussi ces doutes
au début du présent traité, sans démonstration ni explication~: nous
ferons ainsi savoir notre excuse pour avoir dit ce que nous y disons,
et nous montrerons la grande utilité de ce que nous avons
inventé. Seul le saura qui embrasse la connaissance de ce qu'ont
inventé anciens et modernes, la contradiction des fondements, et
l'épuisement de leur effort à établir des fondements suffisants.

\newpage\phantomsection  
\includepdf[pages=2,pagecommand={\thispagestyle{plain}}]{edit.pdf}\phantomsection

Nos prédécesseurs n'ont pu faire réussir cette entreprise, mais Dieu
(loué soit-il) l'a fait réussir~: nous avons indiqué toutes les
prémisses de ces fondements, leurs conséquences, et leurs
démonstrations, dans notre \emph{Commentaire des observations}. Dans
le présent traité, nous avons voulu détacher les fondements des
démonstrations, et exposer ces doutes comme nous avons dit -- nous
demandons à Dieu réussite et protection. Ce traité se compose de deux
parties. La première compte trente chapitres et une
conclusion. L'introduction comprend deux sections.

\newpage\phantomsection  
\index{AYBDBE@\RL{`lm}!AYAGBDBE@\RL{`Alam}, monde}
\index{AYBDBE@\RL{`lm}!AYAGBDBE AYBDBHBJ@\RL{`Alam `alwiy}, monde supérieur, \emph{i. e.} supralunaire)}
\index{AYBDBE@\RL{`lm}!AYAGBDBE ASBABDBJ@\RL{`Alam safliy}, monde inférieur, \emph{i. e.} sublunaire}
\index{AYBDBE@\RL{`lm}!AYAGBDBE BCBHBF BAASAGAO@\RL{`Alam al-kwn wa-al-fsAd}|see{\RL{`Alam safliy}}}
\index{AMAOAO@\RL{.hdd}!BEAMAOAOAO ALBGAGAJ@\RL{m.hddid al-jihAt}, le neuvième orbe}
\index{AHASAW@\RL{bs.t}!ALASBE AHASBJAW@\RL{jsm basI.t}, corps simple}
\index{AQBCAH@\RL{rkb}!ALASBE BEAQBCAH@\RL{jsm mrkb}, corps composé}
\index{AUBHAQ@\RL{.swr}!AUBHAQ@\RL{.sUraT j .suwar}, figure, constellation}
\index{AWAHAY@\RL{.tb`}!AWAHBJAYAI@\RL{.tabI`aT}, nature}
\index{ACAKAQ@\RL{'a_tr}!ACAKBJAQBJ@\RL{'a_tIriy}, éthéré}
\index{AYBFAUAQ@\RL{`n.sr}!AYBFAUAQ@\RL{`un.sar}, élément simple}
\index{BHBDAO@\RL{wld}!BEBHAGBDBJAO AKBDAGAKAI@\RL{al-mwAlId al-_talA_taT}, les trois règnes}
\index{ACAKAQ@\RL{'a_tr}!ACAKAQ AYBDBHBJ@\RL{'a_tr `alwiy}, météore}
\index{ASAMAH@\RL{s.hb}!ASAMAH@\RL{s.hAb j s.hb}, nuage, nébuleuse ??}
\index{ATBGAH@\RL{^shb}!ATBGAH APBHAGAJ APBFAGAH@\RL{^shb _dwAt al-a_dnAb}, étoile filante ??}
\index{AMAQBC@\RL{.hrk}!AMAQBCAI AYAQAVBJAI@\RL{.harakaT `r.diyaT}, mouvement par accident}
\index{AMAQBC@\RL{.hrk}!AMAQBCAI BBASAQBJAI@\RL{.harakaT qsriyaT}, mouvement violent}
\index{AMAQBC@\RL{.hrk}!AMAQBCAI AEAQAGAOBJAI@\RL{.harakaT 'irAdiyaT}, mouvement volitif}
\index{AMAQBC@\RL{.hrk}!AMAQBCAI AWAHBJAYBJAI@\RL{.harakaT .tabI`iyaT}, mouvement naturel}
\index{AMAQBC@\RL{.hrk}!AMAQBCAI AHASBJAWAI BEAQBCBCAHAI@\RL{.harakaT basI.taT mrkkbaT}, mouvement simple-composé}
\addcontentsline{toc}{chapter}{Introduction}
\includepdf[pages=3,pagecommand={\thispagestyle{plain}}]{edit.pdf}\phantomsection

\begin{center}
  \Large Introduction

  \large Première section
\end{center}
Sache que le \emph{monde} est le nom de ce qu'embrasse la
surface visible de l'orbe supérieur. Cet orbe est ce qui définit les
directions, car lui et son centre définissent deux directions~: le
haut et le bas. Il est [l'orbe] absolument simple.

Les corps sont classés en corps \emph{simples} et corps
\emph{composés}. Le corps simple est ce dont les parties et la nature
sont semblables, c'est-à-dire ce qui ne se divise pas en corps de
figures et de natures différentes~; au contraire, le corps simple a
une nature unique et ce qui en émane procède d'une manière
unique\footnote{\emph{d'une manière unique} ou \emph{d'une même
    manière} ou \emph{d'une seule manière} ou \emph{selon le même}~?
  Nous avons choisi la première traduction. C'est un attribut du
  mouvement dans la classification qui va suivre.}. Le corps composé,
c'est le contraire.

Parmi les corps simples, on distingue l'\emph{éther} et les
\emph{éléments simples}. L'éther, c'est les orbes et les astres situés
en eux. On appelle cela monde supérieur ou monde des orbes et des
cieux. Les autres corps simples sont les quatre éléments simples bien
connus. On appelle cela les éléments (y compris ce qui se compose
[d'éléments]). On appelle cela monde inférieur, monde de la génération
et de la corruption, lieu des éléments simples et des corps qui en
sont composés sous l'orbe de la Lune.

Dans les corps composés, on distingue une première division -- ce dont
la composition est parfaite et dont la forme se conserve un temps --
et une seconde division -- sans perfection ni forme qui se
conserve. La première division, c'est les trois règnes~: minéral,
végétal, animal. Le minéral est ce qui ne possède pas la faculté de
croissance, la plante est ce qui possède la faculté de croissance,
sans la perception, et l'animal est ce qui possède la faculté de
croissance avec la perception. La seconde division des corps composés,
c'est les météores\footnote{Les météores appartiennent-ils au monde
  sublunaire ou bien aux cieux~? Al-\d{T}\=us\={\i} laisse entendre
  qu'il s'agit d'éther (\emph{cf.} \cite{altusi1993} \S~II.1), bien
  qu'ils appartiennent au monde sublunaire (\cite{altusi1993}
  \S~II.2). Ils seraient parfois emportés par le mouvement des cieux,
  à cause du fait qu'ils sont assez loin de la Terre (contrairement à
  la partie de l'air <<~adjacente à la Terre~>> qui se conforme au
  mouvement rectiligne, \emph{cf.} \cite{altusi1993} \S~II.1).}~: les
nuages, les étoiles filantes, \emph{etc}. Les lieux des corps composés
sont les lieux des parties qui sont majoritaires en eux pourvu qu'il
n'y ait pas d'attraction en sens contraire.

Le mouvement est la condition de la chose mue entre l'origine et le
terme, en tant que son état à chaque instant diffère de ce qu'il sera
après et de ce qu'il était avant. 

\newpage\phantomsection  
\index{AMAQBC@\RL{.hrk}!AMAQBCAI AHASBJAWAI@\RL{.harakaT basI.taT}, mouvement simple}
\index{AMAQBC@\RL{.hrk}!AMAQBCAI BEANAJBDBAAI@\RL{.harakaT mu_htalifaT}, mouvement irrégulier (\textit{i. e.} non uniforme)}
\index{AMAQBC@\RL{.hrk}!AMAQBCAI BEAQBCBCAHAI@\RL{.harakaT mrkkbaT}, mouvement composé}
\index{AMAQBC@\RL{.hrk}!BEAJAMAQAQBC BFBAASBG AZBJAQBG@\RL{mt.hrrk binafsihi / bi.gayrihi}, mû par soi / par un autre}
\index{ASACBD@\RL{sa'ala}!BEASACBDAI AVAQBHAQBJAI@\RL{mas'alaT .durUriyaT}, principe}
\index{ANBDBH@\RL{_hlw}!ANBDAGAB@\RL{_hlA'}, vide}
\index{ASBCBF@\RL{skn}!ASBCBHBF@\RL{sukUn}, repos}
\label{var27}
\includepdf[pages=4,pagecommand={\thispagestyle{plain}}]{edit.pdf}\phantomsection

\noindent Le mouvement est ou par accident, ou violent, ou volitif, ou
naturel (ou bien il se compose de mouvements naturels, comme les
mouvements des astres, et c'est alors un composé simple que certains
disent aussi être un composé volitif). Le \emph{mouvement par
  accident} est comme le mouvement du passager du bateau. Le\label{bateau}
\emph{mouvement violent} est comme le mouvement du bateau et il lui
vient d'une cause extérieure. Le \emph{mouvement volitif} est comme le
mouvement des animaux qui ont une conscience. Le mouvement est
\emph{naturel} s'il est celui d'êtres qui n'ont pas conscience de ce
qui procède d'eux, comme les éléments, les orbes (à ce que certains
disent), et les plantes.\footnote{Ce paragraphe laisse déjà entendre
  qu'il y a deux doctrines, l'une affirmant que le mouvement des
  astres est un mouvement simple volitif, l'autre que le mouvement des
  orbes est un mouvement simple naturel. D'autre part, c'est la
  première mention d'un mouvement <<~composé simple~>>.}

On distingue le mouvement naturel \emph{qui n'agit pas d'une manière
  unique} (le mouvement des plantes) du mouvement naturel \emph{qui
  agit d'une manière unique}. Selon moi, il y a dans celui-ci deux
divisions~: le mouvement droit et le mouvement circulaire. Le
\emph{mouvement droit} tend vers le centre et il est propre aux
éléments, ainsi qu'à leurs composés. Le \emph{mouvement circulaire}
est selon moi comme le mouvement des orbes.

Il y a une doctrine distinguant le mouvement volitif \emph{qui ne
  procède pas d'une manière unique} (le mouvement des animaux) du
mouvement volitif \emph{qui procède d'une manière unique} comme le
mouvement des orbes autour de leurs centres. Selon cette doctrine, les
orbes et les astres ont une conscience, mais selon une autre doctrine
ce sont des corps simples qui n'ont pas de conscience.
Mon point de vue est qu'ils sont composés (sauf le neuvième orbe),
mais pas composés d'éléments.\footnote{\textit{cf.} l'analyse de ces
deux doctrines dans notre commentaire p.~\pageref{mouvement}.}

Le mouvement céleste individuel est ce qui émane d'un moteur simple
unique, et le mouvement composé ce qui émane d'une multiplicité de
moteurs simples (plus qu'un). Chaque mouvement individuel est simple,
mais chaque mouvement irrégulier est composé~; car il y a certains
mouvements célestes simples composés comme on le montrera si Dieu le
veut.\footnote{\textit{cf.} notre commentaire sur le mouvement
  simple-composé p.~\pageref{simple_compo}.}

On a besoin de sept principes\footnote{Le mot \emph{mas\=a'il} renvoie
  aux \emph{mas\=a'il \d{da}r\=uriyya} de la préface~: d'où ma
  traduction.} physiques~:

\emph{Premier principe.} Le vide est impossible, et dans les orbes, le
repos est impossible.\label{pas_de_repos}

\newpage\phantomsection  
\index{AHAWBDBEBJBHAS@\RL{b.talimayUs}, Ptolémée}
\index{BGBJAC@\RL{haya'a}!BGBJACAI@\RL{hI'aT}, configuration, astronomie}
\index{AHAOAB@\RL{bda'}!BEAHAOAB@\RL{mabda'}, principe, origine} 
\index{ALASBE@\RL{jsm}!ALASBE@\RL{jsm}, corps, solide}
\label{var24}
\includepdf[pages=5,pagecommand={\thispagestyle{plain}}]{edit.pdf}\phantomsection

\emph{Deuxième principe.} Tout mouvement a un principe. Si le principe
n'est pas séparé de la chose mue quant à sa position, de sorte que
l'indication sensible qui le montre est une, on dit qu'elle est
\emph{mue par soi}. Si au contraire il en est séparé, de sorte que le
mouvement se rapporte à la chose mue et que l'action de mouvoir se
rapporte à ce en quoi est son principe, on dit alors qu'elle est
\emph{mue par un autre}.

\emph{Troisième principe.} Un corps qui n'est pas mû par soi mais par
le mouvement d'un  autre corps mû par soi paraît se mouvoir par soi.

\emph{Quatrième principe.} Aucune chose ayant en elle un principe de
mouvement circulaire n'accepte le mouvement rectiligne, et
inversement, sauf par violence.

\emph{Cinquième principe.} Il ne peut y avoir de principe de deux
mouvements différents dans une chose mue simple car la différence des
mouvements exigerait une différence des choses mues\footnote{<<~choses
  mues~>>, \emph{sic}, à moins que \textit{mt\d{h}arrk} puisse aussi
  signifier moteur.}. \`A chaque fois qu'il y a plusieurs mouvements
différents dans un orbe, c'est qu'il a un mouvement par soi et un
mouvement par un autre.

\emph{Sixième principe.} Les orbes n'acceptent pas le mouvement
rectiligne, et il ne peut en être du mouvement de l'astre dans les
cieux comme du mouvement du poisson dans l'eau.

\emph{Septième principe.} Les mouvements des orbes agissent d'une
manière unique, donc il n'accélèrent ni ne ralentissent, ni ne
reviennent avant complétion d'un tour, ni ne s'arrêtent, ni ne sortent
de leur lieu, ni ne changent d'état. Au contraire, ils sont
constamment mus d'un mouvement simple, circulaire, dans la direction
vers laquelle ils tendent. Les mouvements que l'observation trouve
différents de cela sont des mouvements par accident dûs à la
composition de mouvements simples.
  
On ne peut sortir des limites imposées par ces fondements, or ils sont
incompatibles avec la configuration admise par Ptolémée et les autres,
même les meilleurs parmi les modernes, comme en témoignent leurs
livres.

\newpage\phantomsection  
\index{BGBJAC@\RL{haya'a}!BGBJACAI@\RL{hI'aT}, configuration, astronomie}
\index{AHAWBDBEBJBHAS@\RL{b.talimayUs}, Ptolémée}
\index{BABDBC@\RL{flk}!BABDBC ANAGAQAL BEAQBCAR@\RL{falak _hArij al-markaz}, orbe excentrique}
\index{AMAQBC@\RL{.hrk}!AMAQBCAI BEASAJBHBJAI@\RL{.harakaT mustawiyaT}, mouvement uniforme}
\index{AMAQBC@\RL{.hrk}!AMAQBCAI BEANAJBDBAAI@\RL{.harakaT mu_htalifaT}, mouvement irrégulier (\textit{i. e.} non uniforme)}
\index{AMAPBH@\RL{.h_dw}!BFBBAWAI BEAMAGAPACAI@\RL{nuq.taT al-m.hA_dAT}, point de prosneuse}
\includepdf[pages=6,pagecommand={\thispagestyle{plain}}]{edit.pdf}\phantomsection

\begin{center}\large Deuxième section \end{center}
Des doutes et des impossibilités qui nous ont arrêtés dans
les configurations connues.

\label{doutes_excentrique}
Parmi ces [doutes], il y a \emph{l'orbe excentrique}. Un orbe
excentrique rapporté au centre du monde, c'est impossible, car cela
entraîne nécessairement qu'il y ait dans les orbes rapportés au centre
du monde des formes sans la perfection de la circularité et qui ne
sont pas mues. S'il se meut autour du centre du monde cela entraîne
une irrégularité, et si l'excentrique se meut d'un mouvement uniforme
autour de son propre centre, alors il y a des mouvements qui ne sont
pas uniformes par rapport au centre du monde. Si l'on permettait cela,
alors il n'y aurait plus aucune nécessité de poser des orbes
multiples~: autant dire qu'à chaque astre appartient un orbe qui le
meut d'un mouvement irrégulier, qui stationne, rétrograde, accélère ou
ralentit. C'est impossible. Avec le mouvement de l'Apogée, c'est
encore plus impossible. Comme Ptolémée pensait que l'Apogée est fixe,
il a adopté l'excentrique. Quant à l'orbe épicycle\footnote{Ibn
  al-\v{S}\=a\d{t}ir fait ici allusion aux deux modèles équivalents
  proposés par Ptolémée pour le mouvement du Soleil~: l'un faisait
  intervenir un excentrique (de centre fixe distinct du centre du
  monde et dans la direction de l'Apogée, elle-même fixée, cet
  excentrique était en rotation uniforme autour de son propre centre),
  et l'autre faisait intervenir un épicycle.} pour le Soleil, il est
permis, mais il n'est pas en accord avec des observations précises,
comme tu en prendras connaissance à propos de la configuration du
Soleil. Car nous avons trouvé que l'irrégularité du Soleil,
c'est-à-dire l'équation du Soleil, n'est pas en accord avec des
observations précises dans les octants.

\label{doutes_lune}Parmi les doutes concernant les orbes de la Lune,
il y a aussi l'excentrique. L'excentrique dans les orbes de la Lune
est impossible, et l'uniformité du mouvement de l'excentrique autour
d'un point autre que son propre centre est aussi
impossible\footnote{Allusion aux deuxième et troisième modèles de la
  Lune dans l'\textit{Almageste}.}. Le fait que le rayon de l'épicycle
suive un point autre que le centre de l'orbe qui le porte est
impossible\footnote{Allusion au point de prosneuse dans le troisième
  modèle de la Lune dans l'\textit{Almageste}, mais cette critique
  concerne aussi son deuxième modèle.}~: l'orbe
épicycle paraîtrait alors avoir un mouvement multiple. Il y aurait en
lui deux mouvements, l'un uniforme et l'autre irrégulier, or le
mouvement irrégulier ne fait pas un tour complet et c'est impossible.
[Ces hypothèses] entraînent que le rayon
de la Lune dans les quadratures est double du rayon de la Lune dans
les syzygies, or c'est impossible~: on n'a jamais vu ça. 

\newpage\phantomsection  
\index{AWBHASBJ@\RL{n.sIr al-dIn al-.tUsI}, Na\d{s}{\=\i}r al-D{\=\i}n al-\d{T}\=us{\=\i}}
\index{AHAWBDBEBJBHAS@\RL{b.talimayUs}, Ptolémée}
\index{ALASBE@\RL{jsm}!ALASBE@\RL{jsm}, corps, solide}
\index{AMAPBH@\RL{.h_dw}!BFBBAWAI BEAMAGAPACAI@\RL{nuq.taT al-m.hA_dAT}, point de prosneuse}
\index{ANBDBH@\RL{_hlw}!ANBDAGAB@\RL{_hlA'}, vide}
\index{AMAQBC@\RL{.hrk}!AMAQBCAI BD-AYAQAV@\RL{.harakaT al-`r.d}, mouvement en latitude}
\index{AYAWAGAQAO@\RL{`u.tArid}, Mercure}
\index{BABDBC@\RL{flk}!BABDBC BEAYAOAOBD BEASBJAQ@\RL{falak mu`ddl al-msIr}, orbe équant}
\index{AQABAS@\RL{ra'asa}!AQABAS@\RL{ra's}, tête, n{\oe}ud ascendant}
\index{APBFAH@\RL{_dnb}!APBFAH@\RL{_dnb}, queue, n{\oe}ud descendant}
\index{AYAQAV@\RL{`r.d}!AYAQAV@\RL{`r.d}, latitude, \emph{i. e.} par rapport à l'écliptique}
\includepdf[pages=7,pagecommand={\thispagestyle{plain}}]{edit.pdf}\phantomsection

\noindent D'après des
observations précises, quand la distance de la Lune au Soleil est d'un
signe et demi, la Lune vraie déduite de ces principes contredit ce
qu'exige l'observation. D'ailleurs ni Ptolémée ni aucun autre n'a
mentionné d'observation de la Lune vraie dans ces
positions. La localisation de l'astre moyen\label{astre_moyen} et
du centre de l'astre dans un cercle semblable à l'orbe de l'écliptique
est impossible car ils ne sont pas pris dans le même cercle, ni par
rapport au même point.

L'orbe excentrique est aussi impossible pour les autres astres
errants, puisque nous avons dit que l'uniformité du mouvement des
excentriques autour de points autres que leurs centres est
impossible. Quand ils disent que l'excentrique est dans le plan de
l'orbe incliné et que le plan de l'excentrique coupe le parécliptique
en deux points opposés (la tête et la queue), c'est impossible, car
rien sauf un grand cercle ne coupe un grand cercle en deux moitiés, or
l'excentrique n'est pas un grand cercle~; il est donc impossible que
le parécliptique soit coupé en deux points opposés, et il y a là une
subtilité, prête attention. Que le rayon de l'épicycle suive un point
autre que le centre de l'orbe qui le porte, c'est impossible, or ils
l'ont posé ainsi dans les orbes des cinq planètes. L'orbe équant est
impossible. C'est une image fausse. Imaginer une droite dont l'un des
deux bouts est fixe au centre de l'orbe équant, et l'autre s'allonge
ou se raccourcit, passe par le centre des épicycles des astres et les
fait tourner autour de ce point, c'est impossible. La droite n'est pas
corporelle, il faudrait du vide, et d'autres choses impossibles qui en
dépendent. Le fait que Mercure soit plus proche ailleurs qu'au point
opposé à l'Apogée est impossible à cause des observations, bien que
cela ne soit pas impossible à concevoir. Concevoir des orbes dans les
orbes des épicycles des astres errants, de la manière indiquée par
Ptolémée dans les \emph{Hypothèses planétaires}, orbes qui inclinent
leurs plans par rapport au plan du zodiaque, pour le mouvement en
latitude des astres, c'est impossible.

La position mentionnée par Na\d{s}{\=\i}r
al-\d{T}\=us{\=\i} \label{doute_tusi_lune} dans la
\emph{Ta\b{d}kira} sur l'élimination des doutes des orbes de
la Lune est impossible parce qu'on y trouve l'excentrique et que le
diamètre de l'épicycle suit un point autre que le centre de l'orbe
portant l'épicycle.

\newpage\phantomsection  
\index{AYAQAVBJ@\RL{al-mwyd al-`r.dI}, al-Mu'ayyad al-`Ur\d{d}{\=\i}}
\index{ATBJAQAGARBJ@\RL{q.tb al-dIn al-^sIrAzI}, Qu\d{t}b al-D{\=\i}n al-\v{S}{\=\i}r\=az{\=\i}}
\index{ARAQBBAGBDAI@\RL{al-zarqAlaT}, al-Zarq\=ulla}
\index{BBAHBD@\RL{qbl}!AEBBAHAGBD@\RL{'iqbAl}, trépidation, accession}
\index{AOAHAQ@\RL{dbr}!AEAOAHAGAQ@\RL{'idbAr}, trépidation, récession}
\index{AMAPBH@\RL{.h_dw}!BFBBAWAI BEAMAGAPACAI@\RL{nuq.taT al-m.hA_dAT}, point de prosneuse}
\index{ANBDBA@\RL{_hlf}!ACANAJBDAGBA@\RL{i_htilAf}, irrégularité, anomalie, variation}
\index{AWBHASBJ@\RL{n.sIr al-dIn al-.tUsI}, Na\d{s}{\=\i}r al-D{\=\i}n al-\d{T}\=us{\=\i}}
\index{AQBJAN@\RL{ry_h}!BEAQAQBJAN@\RL{marrI_h}, Mars}
\index{ARAMBD@\RL{z.hl}!ARAMBD@\RL{zu.hal}, Saturne}
\index{ATAQBJ@\RL{^sry}!BEATAJAQBJ@\RL{al-mu^starI}, Jupiter}
\index{ARBGAQ@\RL{zhr}!ARBGAQAI@\RL{zuharaT}, Vénus}
\index{APAQBH@\RL{_drw}!APAQBHAI@\RL{_dirwaT}, sommet, apogée}
\index{BABDBC@\RL{flk}!BABDBC ANAGAQAL BEAQBCAR@\RL{falak _hArij al-markaz}, orbe excentrique}
\index{BABDBC@\RL{flk}!BABDBC BEAYAOAOBD BEASBJAQ@\RL{falak mu`ddl al-msIr}, orbe équant}
\index{BBBEAQ@\RL{qmr}!BBBEAQ@\RL{qamar}, Lune}
\index{AMAQBC@\RL{.hrk}!AMAQBCAI BD-AYAQAV@\RL{.harakaT al-`r.d}, mouvement en latitude}
\index{AYAQAV@\RL{`r.d}!AYAQAV@\RL{`r.d}, latitude, \emph{i. e.} par rapport à l'écliptique}
\includepdf[pages=8,pagecommand={\thispagestyle{plain}}]{edit.pdf}\phantomsection

La position mentionnée par Na\d{s}{\=\i}r al-\d{T}\=us{\=\i} dans la
correction de la configuration des quatre planètes Saturne, Jupiter,
Mars et Vénus est impossible, car y demeurent les excentriques, les
orbes équants, l'irrégularité du mouvement propre à cause du mouvement
de l'apogée\footnote{Première occurrence du mot \textit{dhirwa} que
  l'on traduit par <<~apogée~>> (sans majuscule)~: c'est une extrémité
  du diamètre de l'épicycle parallèle à la direction
  Terre--astre~moyen. On traduit \textit{'awj} par <<~Apogée~>> (avec
  une majuscule)~: c'est le point où le centre de l'épicycle est à
  distance maximale de la Terre. Voir la figure
  \ref{fig040}~p.~\pageref{fig040} \textit{infra}.} des épicycles, et
d'autres choses impossibles. Il l'a lui-même reconnu --- Dieu le
pardonne~: il lui manquait un principe suffisant.

La position conjecturée par al-Mu'ayyad al-`Ur\d{d}{\=\i} dans la
configuration des orbes de la Lune, concernant son inversion du sens
du mouvement de l'excentrique et la variation de l'apogée de
l'épicycle qui suit un point autre que le centre de l'orbe portant
l'épicycle, c'est impossible.\label{urdi_lune}

La position mentionnée par al-Mu'ayyad al-`Ur\d{d}{\=\i} dans la
réforme des orbes des astres errants est impossible car y demeurent
l'excentrique, les orbes équants, etc.

La position mentionnée par Qu\d{t}b al-D{\=\i}n
al-\v{S}{\=\i}r\=az{\=\i} dans la réforme de la configuration des
orbes de la Lune est impossible car y demeurent l'excentrique et le
point de prosneuse. L'argument qu'il a avancé (Dieu ait pitié de lui)
pour permettre le point de prosneuse de la Lune est la plus impossible
des impossibilités et c'est une image fausse. Il en est revenu dans la
\emph{Tu\d{h}fa} ; il y indique une autre manière, mais elle est aussi
impossible.\label{shirazi_lune}

Le principe qu'il a nommé l'invention concernant le désordre des
latitudes des astres est impossible.\label{ibdaai}

La position qu'il a mentionnée dans la réforme des orbes des astres
errants est impossible à cause qu'on y trouve les excentriques, les
équants, et l'irrégularité des deux apogées de l'épicycle.

Le mouvement de trépidation\footnote{Littéralement, <<~accession et
  récession~>>.} n'est pas vrai car cela diffère de ce qu'on vérifie
dans les observations anciennes et modernes. C'est d'ailleurs une idée
fausse bien qu'on puisse concevoir de poser des orbes qui causent ces
mouvements. S'ils étaient réels, j'ai déterminé le Soleil du point de
vue d'al-Zarq\=ulla et je l'ai trouvé inférieur à la vérité de plus de
trois degrés et demi --- ceci pendant l'accession. Pendant la
récession, la différence serait de plus de douze degrés s'il est en
fin de récession.

\newpage\phantomsection
\index{ACASAO@\RL{'asad}!BBBDAH ACASAO@\RL{qlb al-'asad}, le c{\oe}ur du Lion}
\index{AHAQAL@\RL{brj}!AHAQBHAL BEBCBHBCAHAI@\RL{brUj mkwkbaT}, constellations du zodiaque}
\index{BEALAQBJAWBJ@\RL{al-majrI.tI}, al-Majr{\=\i}\d{t}{\=\i}}
\index{AMAQBC@\RL{.hrk}!AMAQBCAI AKAGBFBJAI@\RL{.hrkaT _tAnyaT}, deuxième mouvement}
\index{AQAUAO@\RL{r.sd}!AEAQAUAGAO AUAMBJAMAI@\RL{al-'ar.sAd al-.sa.hI.haT}, les observations}
\includepdf[pages=9,pagecommand={\thispagestyle{plain}}]{edit.pdf}\phantomsection

\noindent J'ai déjà déterminé le c{\oe}ur du Lion avec ce mouvement,
or cela n'est pas en accord avec sa position [observée], et cela
contredit toutes les autres observations de cet astre, anciennes et
présentes. La différence est grande et ajoute souvent quatre degrés à
deux degrés. Dans les fondements, je n'ai rien trouvé en quoi il
puisse manquer [quelque chose] de la même grandeur. On a vérifié par
l'observation que le mouvement des fixes est un mouvement simple et
uniforme autour du centre du Monde et sur les pôles du zodiaque, vers
l'Est. Que cela soit différent est impossible.

Al-Majr{\=\i}\d{t}{\=\i} a affirmé que les mouvements de tous les
orbes des astres vont d'Est en Ouest, et que ce qu'on voit aller vers
l'Est vient de la <<~carence~>> des mouvements venant de l'Est à rejoindre
le mouvement diurne. C'est impossible.

Son affirmation que le pôle de l'orbe du zodiaque tourne en un cercle
(autour du pôle du Monde) dont l'azimut est de la grandeur de
l'inclinaison maximale, vers l'Est, par carence du mouvement du
huitième orbe, c'est impossible. Il affirme cela parce qu'il confond
entre les parts de l'équateur et les parts du zodiaque. En fait c'est
le contraire. Par exemple, si les pôles de l'orbe du zodiaque se
meuvent, par carence, vers l'Est, d'un quart de cercle, et si nous
traçons un grand cercle qui passe par les pôles de l'orbe du zodiaque,
par les pôles de l'équateur, et par les solstices, alors il passera
par la part de l'équateur qui suit le commencement du Bélier. On avait
vérifié par l'observation l'immobilité des quatre points de
l'équateur~; seulement il confond avec ces points les parties du
zodiaque étoilé qui n'est pas le zodiaque véritable. Ce qu'il a dit
est impossible à imaginer et cela n'existe pas.

De même qu'il a rendu possible un mouvement des pôles du zodiaque vers
l'Est autour des pôles de l'équateur, il a permis d'autres mouvements
vers l'Est d'une manière semblable~; mais ceci est différent de ce que
j'ai contesté [ci-dessus] sur l'existence des mouvements vers l'Est.
Son énoncé selon lesquel les mouvements vers l'Est sont dus à la
carence du huitième, \textit{etc.} est futile. La démonstration a déjà
établi l'existence de ces mouvements et a montré qu'une carence de
l'orbe pour un mouvement de cette grandeur est impossible.

\newpage\phantomsection
\index{ANBDBA@\RL{_hlf}!ACANAJBDAGBA@\RL{i_htilAf}, irrégularité, anomalie, variation}
\index{BEALAQBJAWBJ@\RL{al-majrI.tI}, al-Majr{\=\i}\d{t}{\=\i}}
\index{ATBEAS@\RL{^sms}!ATBEAS@\RL{^sams}, Soleil}
\index{AMAQBC@\RL{.hrk}!AMAQBCAI BDBHBDAHBJAI@\RL{.harakaT lawlbiyaT}, mouvement hélicoïdal}
\index{ACAQASAWBHAC@RL{'aris.tU'a}, Aristote}
\index{AMAQBC@\RL{.hrk}!AMAQBCAI BEASAJBHBJAI@\RL{.harakaT mustawiyaT}, mouvement uniforme}
\index{AMAQBC@\RL{.hrk}!AMAQBCAI BEAQBCBCAHAI@\RL{.harakaT mrkkbaT}, mouvement composé}
\index{AOBHAQ@\RL{dwr}!AJAOBHBJAQ@\RL{tadwIr}, épicycle}
\index{BABDBC@\RL{flk}!BABDBC AKAGBEBF@\RL{flk _tAmin}, huitième orbe (étoiles fixes)}
\index{BABDBC@\RL{flk}!BABDBC BEAYAOAOBD BFBGAGAQ@\RL{flk m`ddl al-nhAr}, orbe de l'équateur}
\index{AQAUAO@\RL{r.sd}!AEAQAUAGAO AUAMBJAMAI@\RL{al-'ar.sAd al-.sa.hI.haT}, les observations}
\index{AHASAW@\RL{bs.t}!ALASBE AHASBJAW@\RL{jsm basI.t}, corps simple}
\index{AQBCAH@\RL{rkb}!ALASBE BEAQBCAH@\RL{jsm mrkb}, corps composé}
\label{var1}
\includepdf[pages=10,pagecommand={\thispagestyle{plain}}]{edit.pdf}\phantomsection

L'énoncé qui dit qu'à chaque astre il y a un orbe unique avec lequel
il se meut autour des pôles du zodiaque vers l'Est, encore par
carence, et que ces deux mouvements produisent, dans les astres,
arrêt, retour, mouvement droit ou incliné de la trajectoire de part et
d'autre, et d'autres choses qu'on rencontre dans les astres, est un
énoncé futile. Je l'ai déjà étudié avec soin et j'ai établi par la
démonstration sa futilité dans mon livre intitulé \textit{Commentaire
  des observations}.

Tout ce qu'a suivi al-Majr{\=\i}\d{t}{\=\i} dans sa lettre sur la
modification des mouvements et l'établissement de l'orbe du Soleil
autour du centre du monde, son centre sortant des pôles du monde, et
que les pôles de l'orbe du Soleil tournent sur ce cercle excentrique
vers l'Est par carence de cet orbe, c'est une impossibilité. Et cela
que, quand l'axe de l'orbe du Soleil tournait dans l'excentrique sur
les pôles de l'équateur, le Soleil, voire tout point dans le plan de
l'orbe du Soleil, décrivait l'orbe excentrique~; c'est différent de ce
que j'indique. D'autant plus que ce qu'ils ont indiqué n'est pas fondé
sur l'observation~; or nous en avons déjà montré la futilité.

Le mouvement hélicoïdal qu'a mentionné Aristote ne fait apparaître ni
vitesse, ni lenteur, ni arrêt, ni retour dans les astres. Ce n'est
qu'une indication de leur trajectoire décrite par deux mouvements~: le
mouvement diurne et le mouvement qui leur est propre. Quand le Soleil
était à la tête du Capricorne puis se mouvait pendant une demi-année
jusqu'à la tête du Cancer, alors il décrivait une trajectoire
hélicoïdale autour du centre du monde. Aristote n'a indiqué que cela,
par le mouvement hélicoïdal~; mais qui entend cela pense que l'intention
d'Aristote était que, si le Soleil se meut d'un mouvement hélicoïdal,
alors l'irrégularité du mouvement du Soleil se produit, ainsi que
celle des autres astres, or c'est impossible.

Tout mouvement uniforme par rapport à un point est irrégulier par
rapport à un autre, et inversement. De plus, tout orbe dont le
mouvement est uniforme autour d'un point autre que son centre a un
mouvement composé.

L'existence de petits orbes comme les orbes des épicycles, sans
rapport au centre du monde, n'est pas interdit sauf dans la neuvième
sphère. En effet, c'est comme le fait qu'il existe dans chaque orbe un
astre, et dans le huitième orbe de nombreux astres sphériques, chacun
plus grand que les épicycles de certains des astres tandis que cet
astre lui-même est distinct du corps de l'orbe. Donc il n'est pas
interdit qu'il existe des orbes d'épicycles et de choses semblables,
et à partir de là on comprend qu'il y a dans les orbes une certaine
composition ; celui qui est absolument simple, c'est le neuvième orbe,
et l'on ne peut imaginer d'astre sur cet orbe ni rien d'autre
semblable.

\newpage\phantomsection
\index{AHAWBDBEBJBHAS@\RL{b.talimayUs}, Ptolémée}
\index{BABDAM@\RL{ibn afla.h}, Ibn Afla\d{h}}
\index{AYAQAVBJ@\RL{al-mwyd al-`r.dI}, al-Mu'ayyad al-`Ur\d{d}{\=\i}}
\index{ARBGAQ@\RL{zhr}!ARBGAQAI@\RL{zuharaT}, Vénus}
\index{ACBHBD@\RL{'awl}!ACBDAI@\RL{'AlaT j 'AlAt}, instrument}
\index{ACAQAV@\RL{'ar.d}!ACAQAV@\RL{'ar.d}, la Terre}
\index{BABDBC@\RL{flk}!BABDBC AHAQBHAL@\RL{flk al-burUj}, orbe de l'écliptique}
\index{AHAQAL@\RL{brj}!BEBFAWBBAI AHAQBHAL@\RL{min.taqaT al-burUj}, écliptique, ceinture de l'écliptique}
\index{BBBJAS@\RL{qys}!BEBBBJAGAS@\RL{miqyAs j maqAyIs}, gnomon}
\index{ATAYAY@\RL{^s``}!ATAYAGAY@\RL{^su`A`}, rayon (de lumière)}
\index{AQAUAO@\RL{r.sd}!AEAQAUAGAO AUAMBJAMAI@\RL{al-'ar.sAd al-.sa.hI.haT}, les observations}
\index{BGBFAO@\RL{hnd}!BGBFAO@\RL{al-hnd}, les Indiens}
\includepdf[pages=11,pagecommand={\thispagestyle{plain}}]{edit.pdf}\phantomsection
Les \emph{man\u{s}\=ura} qu'indique Ptolémée dans son livre intitulé
\emph{Hypothèses planétaires} sont futiles.

Ce qu'ont indiqué Ibn Afla\d{h} et al-Mu'ayyad al-`Ur\d{d}{\=\i} quant
au fait que Vénus est au-dessus du Soleil, et ce qu'ils en ont conclu,
c'est impossible.

Le rapprochement de la circonférence de l'écliptique et de la
circonférence de l'équateur, et leur éloignement, c'est impossible.

La variation de l'inclinaison maximale n'est pas due à une différence
de temps ni de lieu. Ce qu'on y trouve de variable est dû aux
instruments, à ce qui diffère dans leur érection, c'est dû à
l'irrégularité de l'ombre des extrémités des gnomons, c'est dû à une
variation de la lumière du Soleil dans le corps de la cible, c'est dû
au fait que le centre du volume de la Terre n'est pas fixe par rapport
au centre du monde, comme nous l'avons exposé dans mon livre intitulé
\textit{Commentaire des observations}.

D'après les Indiens et d'autres, on indique que l'inclinaison totale
du Soleil est vingt-quatre degrés. C'est futile et il n'y aucune
vérité là-dedans. C'est comme si cette doctrine indiquait le corps qui
contient l'inclinaison du Soleil des deux côtés~; son extension en
latitude est de quarante-huit degrés dont la moitié est vingt-quatre
degrés~; ceci est à peu près la latitude maximale du soleil mais
\emph{son centre} ne dépasse pas $23;31$. Si l'on prend cette
doctrine de ce point de vue, alors c'est vrai et sans défaut.

\newpage\phantomsection  
\index{BGBFAOASAI@\RL{hndsaT}, géométrie}
\index{ANAWAW@\RL{_h.t.t}!ANAWAW@\RL{_ha.t.t}, ligne}
\index{ANAWAW@\RL{_h.t.t}!ANAWAW BEASAJBBBJBE@\RL{_ha.t.t mustaqIm}, droite}
\index{ASAWAM@\RL{s.t.h}!ASAWAM@\RL{s.t.h}, surface}
\index{ASAWAM@\RL{s.t.h}!ASAWAM BEASAJBH@\RL{\vocalize s.t.h mstwiN}, plan}
\index{ALASBE@\RL{jsm}!ALASBE@\RL{jsm}, corps, solide}
\phantomsection
\addcontentsline{toc}{chapter}{I.1 Fondements posés comme axiomes en géométrie}
\label{var10}
\includepdf[pages=12,pagecommand={\thispagestyle{plain}}]{edit.pdf}\phantomsection

\begin{center}
  \Large Chapitre un

  \large Fondements posés comme axiomes en géométrie
\end{center}
Le point, la ligne, l'angle, la surface, le cercle et la
sphère sont des êtres de figures connues. Chaque point, droite, plan
ou angle s'ajuste sur son semblable. L'intersection de deux lignes
quelconques est un point, celle de deux surfaces est une ligne, et
celle de deux solides est une surface.

Nous pouvons supposer une droite donnée sur n'importe quel plan ou
bien passant par un point arbitraire. Nous pouvons fixer un point sur
n'importe quelle droite ou plan. Nous pouvons reporter une droite
finie sur n'importe quel plan ou solide. Nous pouvons tracer un cercle
en un point quelconque avec un rayon quelconque.

Les angles droits sont égaux. Lorsqu'une droite tombe sur une droite,
cela produit deux angles droits, ou bien un angle aigu et un angle
obtus~; la somme des deux angles est égale à deux droits. La somme des
angles de tout triangle est égale à deux droits.

\newpage\phantomsection  
\index{ASBEBH@\RL{smw}!ASBEAGAB@\RL{al-samA'}, les cieux}
\index{AKBBBD@\RL{_tql}!BEAQBCAR AKBBBD@\RL{markaz al-_tql}, centre de gravité}
\phantomsection
\addcontentsline{toc}{chapter}{I.2 Fondements déjà établis}
\includepdf[pages=13,pagecommand={\thispagestyle{plain}}]{edit.pdf}\phantomsection

\begin{center}
  \Large Chapitre deux

  \large Fondements déjà établis
\end{center}
La surface de la Terre et de l'eau est approximativement ronde.

La rotondité des cieux est comme celle de la sphère.

La Terre est au cieux comme le centre de la sphère est à son bord.

La grandeur de la Terre est sensible comparée à ce qui est au-dessous
de l'orbe de Mars\footnote{Donc la parallaxe sera négligeable dans
  l'observation des planètes supérieures (Mars, Jupiter, Saturne).}.

Le centre de gravité de la Terre coïncide avec le centre du monde (et
non le centre de son volume) où elle réside sans mouvement d'aucun
côté.

Le haut est ce qui est vers les cieux et s'éloigne du centre~; le bas
est ce qui est vers le centre. Tous les graves tendent vers le centre,
et ce qui est léger tend vers le bord.

\newpage\phantomsection  
\index{BABDBC@\RL{flk}!BABDBC BEAYAOAOBD BFBGAGAQ@\RL{flk m`ddl al-nhAr}, orbe de l'équateur}
\index{AYAOBD@\RL{`dl}!BEAYAOAOBD BFBGAGAQ@\RL{m`ddl al-nhAr}, équateur}
\index{BABDBC@\RL{flk}!BABDBC AHAQBHAL@\RL{flk al-burUj}, orbe de l'écliptique}
\index{BABDBC@\RL{flk}!BABDBC BEBCBHBCAH@\RL{flk mkwkb}, orbe des étoiles fixes (\emph{syn.} orbe de l'écliptique)}
\index{AHAQAL@\RL{brj}!AHAQBHAL BEBCBHBCAHAI@\RL{brUj mkwkbaT}, constellations du zodiaque}
\index{ARBEBF@\RL{zmn}!ARBEAGBF@\RL{'azmAn mu`addl al-nahAr}!temps équatoriaux}
\index{AWBDAS@\RL{.tls}!AWBDAS@\RL{a.tlas}, Atlas |see{\RL{flk m`ddl al-nhAr}}}
\index{BCBCAH@\RL{kkb}!BCBHBCAH AJBGAGAHAJ@\RL{kwkb thAbit}, étoile fixe}
\index{ALAR@\RL{jz}!ALARAB@\RL{juz' j 'ajzA'}, partie, segment, portion, part}
\index{AQAJAH@\RL{rtb}!AJAQAJBJAH@\RL{tartIb}, ordre, disposition}
\phantomsection
\addcontentsline{toc}{chapter}{I.3 Disposition des neuf orbes}
\label{var11}
\includepdf[pages=14,pagecommand={\thispagestyle{plain}}]{edit.pdf}\phantomsection

\begin{center}
  \Large Chapitre trois

  \large Disposition des neuf orbes
\end{center}
Le premier est celui qui entoure les huit. Sa forme est
sphérique et son centre est le centre du monde. Il a deux pôles sur
lesquels il tourne autour de son centre, centre du monde et centre du
tout, d'un mouvement uniforme qui fait, en un jour et une nuit, une
révolution plus la coascension du Soleil en un jour et une
nuit\footnote{``Un jour et une nuit'' désigne une jour solaire. Un
  jour solaire est en effet un peu plus long qu'un jour tropique.}. De
ce mouvement, l'ensemble des orbes se meut d'Est en Ouest autour de
son centre et sur ses pôles. La ceinture de ce mouvement, lieu
équidistant des pôles de cet orbe, est appelé \emph{équateur}~; le
pôle situé à gauche quand on est tourné vers l'Est est appelé pôle
Nord de l'équateur, et son homologue, pôle Sud. Dans le plan de
l'équateur on suppose qu'il y a un second cercle dont le bord est
divisé selon les parts et leurs fractions, au nombre de trois cent
soixante parts. Ces divisions s'appellent temps équatoriaux ou encore
parts de l'équateur~; à cause de l'immobilité de ce cercle, c'est par
lui qu'on mesure les mouvements. Cet orbe, le premier, est celui qui
est absolument simple~; il définit les directions et on l'appelle
l'\emph{Atlas}.

Le deuxième orbe est l'orbe des astres appelés étoiles fixes. Sa forme
est sphérique, sa partie convexe touche la partie concave du neuvième
et sa partie concave touche la partie convexe du septième qui entoure
les orbes de Saturne. Il a deux pôles fixes sur lesquels il tourne
autour du centre de monde d'un mouvement simple et uniforme d'Ouest en
Est, comme on l'a vérifié par la démonstration et par
l'observation. Dans le plan de la ceinture de ce mouvement (lieu
équidistant des pôles de cet orbe), sont dessinés les douze divisions
du zodiaque, leurs degrés et leurs fractions, en commençant par la
partie qui longe la figure du Bélier dans l'orbe des étoiles fixes
dans la première des désignations\footnote{Peut-être {\shatir} fait-il
  allusion ici à une ``première nomenclature''. Il expliquera plus
  loin qu'une première nomenclature consiste à nommer les divisions du
  zodiaque selon les constellations (les ``figures'' formées par les
  étoiles, donc). La seconde nomenclature consiste à les nommer en
  partant du point vernal : la première division s'appellerait alors
  ``Bélier'' par pure convention, et elle ne longerait pas toujours la
  constellation du Bélier, à cause du mouvement de précession.}, comme
on sait.

Le plan de ce lieu est incliné de $23;31$ degrés par rapport au plan
de l'orbe de l'équateur (c'est là l'inclinaison maximale). 

\newpage\phantomsection  
\index{BBBDAH@\RL{qlb}!BEBFBBBDAH AUBJBABJBJ@\RL{mnqlb .sayfiyy}, solstice d'été}
\index{BBBDAH@\RL{qlb}!BEBFBBBDAH ATAJBHBJ@\RL{mnqlb ^satawI}, solstice d'hiver}
\index{AYAOBD@\RL{`dl}!AYAJAOAGBD AQAHBJAYBJBJ@\RL{i`tidAl rabI`iyy}, équinoxe de printemps}
\index{AYAOBD@\RL{`dl}!AYAJAOAGBD ANAQBJBABJBJ@\RL{i`tidAl _hrIfiyy}, équinoxe d'automne}
\index{BCBCAH@\RL{kkb}!BCBHBCAH AJBGAGAHAJ@\RL{kwkb thAbit}, étoile fixe}
\index{AHAQAL@\RL{brj}!AHAQAL@\RL{brj}, signe (du zodiaque) |see{\RL{falak al-burUj}}}
\index{BABDBC@\RL{flk}!BABDBC AHAQBHAL@\RL{flk al-burUj}, orbe de l'écliptique}
\index{AHAQAL@\RL{brj}!AJBHAGBDBJ AHAQBHAL@\RL{`al_A tawAliy al-burUj}, dans le sens des signes}
\index{AHAQAL@\RL{brj}!ANBDAGBA AJBHAGBDBJ AHAQBHAL@\RL{il_A _hilAf tawAliy al-burUj}, en sens contraire de celui des signes}
\index{AOAQAL@\RL{drj}!AOAQALAI@\RL{drjaT j darjAt, drj}, degré}
\label{var28}
\includepdf[pages=15,pagecommand={\thispagestyle{plain}}]{edit.pdf}\phantomsection

\noindent Considérons le grand cercle qui touche les deux pôles de
l'équateur et les deux pôles de l'écliptique~; il coupe la ceinture de
l'écliptique en deux points. L'intersection proche du pôle Nord
s'appelle \emph{solstice d'été}, et la seconde, \emph{solstice
  d'hiver}.

Le plan de l'écliptique coupe l'équateur en deux points~: le point qui
précède le solstice d'été s'appelle l'\emph{équinoxe de printemps}, et celui
qui précède le solstice d'hiver, l'\emph{équinoxe d'automne}. Ces quatre
points sont immobiles ; ils ne se meuvent pas par rapport au lieu
qu'ils longent dans l'équateur. On découpe en trois portions égales
chacun des arcs
délimités par ces points~: nous imaginons six grands cercles qui
passent par les pôles de l'écliptique, et la ceinture est coupée par eux
en douze portions égales, de même que toute la figure du huitième
orbe. On appelle chacune de ces portions un \emph{signe du zodiaque}, et
chaque signe est divisé en trente portions égales qu'on appelle degrés
de l'écliptique.

Il y a très longtemps, la constellation du Bélier (dans le huitième
orbe, l'orbe des étoiles fixes) longeait la portion dont le
commencement est l'équinoxe de printemps. Cette portion s'appelle le
signe du Bélier. La deuxième portion s'appelle selon la constellation
qui la longeait, c'est-à-dire le Taureau, puis les Gémeaux, puis le
Cancer, \textit{etc.}, comme vous savez. Mais à la fin de l'an huit
cent vingt de l'hégire la première portion longera la constellation
des Poissons de l'orbe des étoiles fixes, donc les constellations se
déplaceront d'un signe dans le sens des signes\footnote{Le mot
  <<~signe~>> désigne désormais un arc de trente degrés de
  l'écliptique, et la locution <<~dans le sens des signes~>> indique
  l'ordre (conventionnel) des signes du zodiaque, c'est-à-dire d'Ouest
  en Est.}. Les portions immobiles longeront d'autres constellations
que celles qu'elles longeaient autrefois. On pourra alors appeler
Poissons la première portion, Bélier la deuxième, Taureau la
troisième, \textit{etc.} Ou bien on gardera leurs noms, mais on ne
considèrera plus les constellations qu'ils longent dans le zodiaque
étoilé. On a observé que ce mouvement est d'un degré en soixante-dix
années persanes. Chez Ptolémée il est d'un degré par siècle, chez les
modernes il est d'un degré et demi par siècle, et chez certains
d'entre eux il est d'un degré tous les soixante-dix ans, en années
persanes.

\newpage\phantomsection  
\index{BCBCAH@\RL{kkb}!BCBHBCAH AJBGAGAHAJ@\RL{kwkb thAbit}, étoile fixe}
\index{AMAQBC@\RL{.hrk}!AMAQBCAI BJBHBEBJBJAI@\RL{.harakaT yawumiyyaT}, mouvement diurne}
\index{BABDBC@\RL{flk}!BABDBC BFAGAQ@\RL{falak al-nAr}, orbe du feu}
\index{BABDBC@\RL{flk}!BABDBC BGBHAG@\RL{flk al-hwA}, orbe de l'Air}
\index{BABDBC@\RL{flk}!BABDBC BEAGAB@\RL{flk al-mA'}, orbe de l'Eau}
\index{BABDBC@\RL{flk}!BABDBC AJAQAGAH@\RL{flk al-trAb}, orbe de la Terre}
\index{ACAQAV@\RL{'ar.d}!BCAQAI ACAQAV@\RL{kuraT al-'ar.d}, globe terrestre}
\includepdf[pages=16,pagecommand={\thispagestyle{plain}}]{edit.pdf}\phantomsection

Ceci étant dit, sachez que l'orbe inférieur (le premier) appartient à
la Lune, le deuxième à Mercure, le troisième à Vénus, le quatrième au
Soleil, le cinquième à Mars, le sixième à Jupiter, le septième à
Saturne, et dans le huitième sont enchassées les étoiles fixes, appelées
ainsi car leurs positions les unes par rapport aux autres sont
fixes.

Le neuvième, mentionné ci-dessus, se meut sur ses pôles et autour de
son centre d'Est en Ouest. C'est le mouvement diurne, et les orbes des
sept astres sont mûs par ce mouvement parallèlement à l'équateur. Sous
l'orbe de la Lune, il y a l'orbe du Feu, puis l'Air, puis l'Eau, puis
la Terre, c'est-à-dire le globe terrestre. Nous avons déjà dénoncé la
faiblesse du raisonnement selon lequel Vénus est au-dessus du
Soleil. Rare sont les savants qui ont vu cela. Nous avons expliqué
cela ailleurs dans nos livres. Dieu est le plus savant.

\newpage\phantomsection 
\index{ASBHBI@\RL{sw_A}!ANAWAW ASAJBHAGAB@\RL{_ha.t.t al-istiwA'}, équateur terrestre}
\index{AYAOBD@\RL{`dl}!BEAYAOAOBD BFBGAGAQ@\RL{m`ddl al-nhAr}, équateur}
\index{AYAQAV@\RL{`r.d}!AYAQAV AHBDAO@\RL{`r.d al-bld}, latitude du pays, \emph{i. e.} par rapport à l'équateur}
\index{BFAUBA@\RL{n.sf}!AOAGAEAQAI BFAUBA BFBGAGAQ@\RL{dA'iraT n.sf al-nhAr}, méridien}
\index{AOBHAQ@\RL{dwr}!AOAGAEAQAI BFAUBA BFBGAGAQ@\RL{dA'iraT n.sf al-nhAr}|see{\RL{n.sf}}}
\index{AWBHBD@\RL{.twl}!AWBHBD AHBDAO@\RL{.tUl al-bld}, longitude du pays (par rapport à l'équateur)}
\addcontentsline{toc}{chapter}{I.4 Cercles remarquables et latitudes des pays}
\includepdf[pages=17,pagecommand={\thispagestyle{plain}}]{edit.pdf}\phantomsection

\begin{center}
  \Large Chapitre quatre

  \large Cercles remarquables et latitudes des pays
\end{center}
Le plus remarquable des grands cercles est le cercle de
l'\emph{équateur}. C'est celui qui est équidistant des pôles de
l'équateur, qui est ceinture du mouvement diurne et cercle des
équinoxes.  Qu'on conçoive son plan coupant le globe terrestre, la
ligne de coupe s'appelle équateur en raison de l'égalité du jour et de
la nuit pour tous les habitants qui s'y trouvent, parce que les cieux
y tournent en roue avec les pôles de l'équinoxe à l'horizon. Puisque
la Terre est sphérique et stable au centre du monde, pour qui avance
sur Terre de l'équateur en allant vers le Nord, le pôle Nord s'élève
au-dessus de l'horizon~; et pour qui avance vers le Sud, le pôle Sud
s'élève au-dessus de l'horizon. Plus l'un des deux pôles s'élève, plus
l'autre s'abaisse.

L'élévation du pôle s'appelle \emph{latitude du pays}.  C'est un arc
du grand cercle\footnote{Le cercle méridien décrit dans le paragraphe
  suivant.} passant par le zénith et par les pôles de l'équateur, cet
arc est égal à l'abaissement de l'autre pôle [sous l'horizon], et il
est aussi égal à la distance entre l'équateur et le zénith dans cette
région.

Le \emph{cercle méridien} ou longitude du pays coupe toutes les
trajectoires visibles en deux moitiés. Il passe par le pôle visible et
par le zénith. Il change quand on se déplace vers l'Est ou vers
l'Ouest mais ne change pas quand on suit le Nord ou le Sud.

Pour chaque pays, soit son méridien ; nous imaginons qu'il coupe la
sphère terrestre ou qu'il y est dessiné. Ce qui est compris entre le
méridien passant par l'extrémité Ouest du monde habité et le méridien
d'un lieu supposé, c'est la \emph{longitude du pays}. Son maximum est
cent quatre-vingt degrés. La latitude maximale est quatre-vingt-dix
degrés, et elle est atteinte là où l'un des pôles est au zénith. Les
cieux tournent alors comme une meule; aux autres latitudes, ils
tournent de manière oblique. 

\newpage\phantomsection 
\index{AHAQAL@\RL{brj}!BEBFAWBBAI AHAQBHAL@\RL{min.taqaT al-burUj}, écliptique, ceinture de l'écliptique}
\index{BABDBC@\RL{flk}!BABDBCBEBEAKAKBD@\RL{falak muma_t_tal}, parécliptique}
\index{AHAQAL@\RL{brj}!AHAQAL@\RL{brj}, signe (du zodiaque) |see{\RL{falak al-burUj}}}
\index{BABDBC@\RL{flk}!BABDBC AHAQBHAL@\RL{flk al-burUj}, orbe de l'écliptique}
\index{BEBJBD@\RL{myl}!AOAGAEAQAI BEBJBD@\RL{dA'iraT al-mIl}, cercle de déclinaison}
\index{AOBHAQ@\RL{dwr}!AOAGAEAQAI BEBJBD@\RL{dA'iraT al-mIl}|see{\RL{mIl}}}
\index{BEBJBD@\RL{myl}!BEBJBD ACBHBHBD@\RL{mIl 'awwal}, première déclinaison, \emph{i. e.} par rapport à l'équateur}
\index{AWBDAY@\RL{.tl`}!BEAWAGBDAY@\RL{m.tAl`}, coascension}
\index{BEAQAQ@\RL{mrr}!BEBEAQAQ@\RL{mamarr j At}, passage, transit (en général au méridien)}
\index{BEAQAQ@\RL{mrr}!AOAQALAI BEBEAQAQ@\RL{darjaT mamarr}, degré de transit}
\index{AYAQAV@\RL{`r.d}!AYAQAV@\RL{`r.d}, latitude, \emph{i. e.} par rapport à l'écliptique}
\index{AOBHAQ@\RL{dwr}!AOAGAEAQAI AYAQAV@\RL{dA'iraT al-`r.d}, cercle de latitude}
\index{BEBJBD@\RL{myl}!BEBJBD AKAGBFBJ@\RL{mIl _tAniy}, seconde déclinaison, \emph{i. e.} par rapport à l'écliptique}
\index{AYAQAV@\RL{`r.d}!AYAQAV BEAYAOAOBD@\RL{`r.d m`ddl}, latitude rapportée à l'équateur (mais mesurée sur un cercle de latitude)}
\includepdf[pages=18,pagecommand={\thispagestyle{plain}}]{edit.pdf}\phantomsection

\noindent Les régions habitées sur Terre vont : [en
  latitude], de l'équateur jusque là où le pôle Nord atteint une
hauteur égale au complément de l'inclinaison [de l'écliptique]
$66;29$, et en longitude, d'une longitude d'un degré jusqu'aux cent
quatre-vingt degrés.

Ceci étant admis, parmi les cercles remarquables je dis qu'il y a
aussi la \emph{ceinture de l'écliptique}. C'est la ceinture du
huitième orbe et on l'appelle aussi <<~voie du Soleil~>> car le Soleil
la suit toujours. Nous imaginons que ce cercle découpe tous les orbes
et dessine dans chacun un grand cercle appelé \emph{orbe
  parécliptique}. La moitié de la ceinture de l'écliptique qui commence
par l'équinoxe de printemps est au Nord de l'équateur, et l'autre
moitié est au Sud. Chacune des deux moitiés est divisée en six
portions appelées \emph{signes}, comme précédemment.

Parmi les grands cercles, il y a aussi le \emph{cercle de déclinaison}.
C'est le grand cercle qui passe par les pôles de
l'équateur céleste et par n'importe quel point ou par un astre donné~;
il coupe la ceinture de l'équateur céleste. Le plus petit des deux
arcs [du cercle de déclinaison] situés entre ce point et l'équateur
est la déclinaison de ce point appelée \emph{première déclinaison}. Si
ce cercle est celui qui passe par la tête du Cancer et du Capricorne, il
passe alors aussi par les pôles de l'orbe de l'écliptique et la
déclinaison [d'un astre sur l'écliptique] est ici maximale. Soit un
cercle de déclinaison passant par un astre, alors l'arc compris entre
l'astre et l'équateur céleste est la distance de l'astre à l'équateur
céleste. Le point où il rencontre l'équateur est la \emph{coascension}
du degré de transit de cet astre au méridien~; et le second point, où
il rencontre la ceinture de l'écliptique, est le \emph{degré de transit} de cet
astre.

Parmi les grands cercles, il y a aussi le \emph{cercle de
  latitude}. C'est le grand cercle imaginaire passant par les pôles de
l'écliptique et par un point donné de l'écliptique ou par un astre. On
le connaît aussi sous le nom de \emph{cercle de seconde
  déclinaison}. [S'il passe par un point de l'écliptique], l'angle
aigu situé sur le cercle [de latitude] entre la ceinture de
l'écliptique et l'équateur est la seconde déclinaison de ce point,
dite aussi \emph{latitude} de ce point.

\newpage\phantomsection 
\index{AYAQAV@\RL{`r.d}!AMAUAUAI AYAQAV@\RL{.hi.s.saT `r.d}, argument de latitude}
\index{AWBHBD@\RL{.twl}!AWBHBD@\RL{.tUl}, longitude (par rapport à l'écliptique)}
\index{AWBHBD@\RL{.twl}!AWBHBD BEAYAOAOBD@\RL{.tUl m`ddl}, longitude rapportée à l'équateur (mais mesurée sur un cercle de latitude)}
\index{BABHBB@\RL{fwq}!ADBABB@\RL{'ufuq}, horizon}
\index{ASBEAJ@\RL{smt}!ASBEAJ AQABAS@\RL{samt al-ra's}, zénith}
\index{ASBEAJ@\RL{smt}!ASBEAJ BBAOBE@\RL{samt al-qadam, samt al-rjl}, nadir}
\index{BFAUBA@\RL{n.sf}!AOAGAEAQAI BFAUBA BFBGAGAQ@\RL{dA'iraT n.sf al-nhAr}, méridien}
\index{AYAQAV@\RL{`r.d}!AYAQAV AHBDAO@\RL{`r.d al-bld}, latitude du pays, \emph{i. e.} par rapport à l'équateur}
\index{BFAUBA@\RL{n.sf}!ANAWAW BFAUBA BFBGAGAQ@\RL{_ha.t.t n.sf al-nahAr}, ligne méridienne}
\index{ANAWAW@\RL{_h.t.t}!ANAWAW BFAUBA BFBGAGAQ@\RL{_ha.t.t n.sf al-nahAr}|see{\RL{n.sf}}}
\includepdf[pages=19,pagecommand={\thispagestyle{plain}}]{edit.pdf}\phantomsection

\noindent S'il passe par un astre, l'arc
de ce cercle compris entre l'astre et la ceinture de l'écliptique du
côté le plus proche est la \emph{latitude} de cet astre, et le petit
arc de ce cercle compris entre l'astre et l'équateur est la latitude
de l'astre rapportée à l'équateur et on l'appelle aussi \emph{argument
  de latitude} de l'astre.\footnote{Dans les modèles planétaires,
  l'\emph{argument de latitude} désigne plutôt l'arc de l'écliptique
  entre le n{\oe}ud et la longitude moyenne de l'astre.  Dans tout ce
  paragraphe l'auteur semble viser une synthèse rationnelle des usages
  du mot \textit{`ar\d{d}} -- latitude. Pour un astre de latitude
  nulle, au sens où cet astre serait sur l'écliptique, l'usage devait
  prescrire à ce terme un autre sens que celui auquel nous sommes
  habitués aujourd'hui.} Si le cercle de latitude passe par un astre donné et
qu'il coupe la ceinture de l'écliptique, il la coupe en le degré de
\emph{longitude}
de l'astre. L'origine des longitudes de l'astre est l'équinoxe de
printemps, à la tête du Bélier, par convention.
Si le cercle de latitude [passant par un point donné] coupe
l'équateur, l'intersection est la longitude de l'astre rapportée à
l'équateur. Si le point donné est un des solstices, le cercle de
latitude [passant par ce point] est confondu avec le cercle de
première déclinaison~: c'est le cercle passant par les pôles de
l'écliptique et par les pôles de l'équateur, les première et seconde
déclinaisons sont réunies, et chacune est à son maximum $23;31$.

Parmi les cercles remarquables, il y a le \emph{cercle de
  l'horizon}. C'est celui qui sépare la partie visible de la partie
invisible des cieux.  Un de ses pôles est le \emph{zénith}, l'autre le
\emph{nadir}. Sur ce cercle ont lieu les levers et les couchers. [Ce
  cercle] est décrit par la droite issue de l'observateur, tangente à
la sphère terrestre et allant jusqu'au neuvième orbe~; or la partie
visible est plus grande que la partie cachée.

Il a tort celui qui dit que la partie visible est plus grande que la
partie cachée d'une grandeur de quatre minutes et vingt-six
secondes. J'ai trouvé par l'observation que les parts visibles du
Soleil et des astres sont plus grandes que ce qu'on a calculé, d'une
grandeur variable qui dépasse deux tiers de degré et qui est composée
de combien est descendu l'observateur et de combien a baissé le centre
du volume du globe terrestre par rapport au centre du monde. La raison
en est que c'est le centre de gravité de la Terre qui est confondu
avec le centre du monde, et non le centre de son volume. La différence
entre les centres est nécessaire, car [la Terre] est constituée d'eau
et de terre, et le centre de gravité est du côté de la terre~; ceci a
été établi par examens successifs. J'ai déjà consacré un écrit à ce
sujet.

Parmi les grands cercles, il y a le \emph{cercle méridien} dont on a
déjà décrit certaines propriétés. C'est le grand cercle imaginaire
passant par les pôles de l'horizon et par les pôles de l'équateur. Si
le pôle de l'horizon est pôle de l'équateur, [le méridien] est [un]
cercle perpendiculaire à l'horizon qui coupe les trajectoires diurnes
visibles en deux moitiés. 

\newpage\phantomsection 
\index{AQBAAY@\RL{rf`}!AQAJBAAGAY@\RL{irtifA`}, hauteur (coordonnées azimutales)}
\index{AMAWAW@\RL{.h.t.t}!BFAMAWAGAW@\RL{in.hi.tA.t}, abaissement|see{\RL{irtifA`}}}
\index{ASBEAJ@\RL{smt}!AOAGAEAQAI ACBHBHBD ASBEBHAJ@\RL{dA'iraT 'awwal al-sumUt}, cercle origine des azimuts}
\index{AOBHAQ@\RL{dwr}!AOAGAEAQAI ACBHBHBD ASBEBHAJ@\RL{dA'iraT 'awwal al-sumUt}|see{\RL{smt}}}
\index{ASBEAJ@\RL{smt}!ASBEAJ@\RL{smt}, direction, azimut}
\index{AYAQAV@\RL{`r.d}!AYAQAV BBBDBJBE AQBHBJAI@\RL{`r.d aqlIm al-rwyaT}, latitude de ce qui est visible sous un climat, \textit{i. e.} latitude du pôle de l'horizon par rapport à l'écliptique}
\index{AQBAAY@\RL{rf`}!AQAJBAAGAY@\RL{irtifA`}, hauteur (coordonnées azimutales)}
\index{AOBHAQ@\RL{dwr}!AOAGAEAQAI AQAJBAAGAY@\RL{dA'irat al-irtifA'}, cercle de hauteur}
\index{AOBHAQ@\RL{dwr}!AOAGAEAQAI ADBABB AMAGAOAK@\RL{dA'iraT al-'ufuq al-.hAdi_t}, cercle de l'horizon occurrent}
\index{ATAQBB@\RL{^srq}!ANAWAW BEATAQBB BEAZAQAH@\RL{_ha.t.t al-ma^sraq wa-al-ma.grab}, ligne de l'Est et de l'Ouest}
\label{var4}
\includepdf[pages=20,pagecommand={\thispagestyle{plain}}]{edit.pdf}\phantomsection

\noindent Il coupe la ceinture de l'écliptique
  au-dessus de la Terre en un point qu'on appelle le dixième, et
  au-dessous de la Terre en un point qu'on appelle le quatrième ou
  bien le piquet. Quand l'astre arrive sur ce cercle, il est au terme
de son ascension si c'est un astre pas toujours visible~; si c'est un
astre toujours visible, alors il a sur ce cercle deux hauteurs, et la
moitié de leur somme est la latitude du lieu. L'arc situé sur le
méridien entre les pôles de l'équateur et l'horizon vrai est la
latitude du pays. Il coupe la Terre en une ligne qui s'appelle
\emph{ligne méridienne}.

Parmi les grands cercles, il y a le cercle de l'Est et de
l'Ouest. C'est le grand cercle imaginaire passant par le zénith et le
nadir et par les pôles du méridien, donc il est perpendiculaire à
l'horizon, et on l'appelle \emph{cercle origine des azimuts}.
Si un astre ou un point donné y passe, il n'aura pas d'azimut lors du
passage. Si on l'imagine coupant le globe terrestre, [il y dessine] la
ligne de l'Est et de l'Ouest qui est perpendiculaire à la ligne
méridienne. La ligne de l'Est et de l'Ouest ne change pas quand on va
vers des lieux plus à l'Est ou plus à l'Ouest\footnote{Est-ce pour autant
  qu'Ibn al-\v{S}\=a\d{t}ir ne connaît pas la loxodromie~? Il était cependant
  connu que, pour un pays à même latitude que la Mecque et à l'Ouest, la
  direction de la \textit{qibla} n'est pas le point Est, et qu'elle est ``à gauche''
  du point Est (\textit{cf. infra} p.~\pageref{loxoqibla}).},
mais elle change quand on va vers le Nord ou vers le Sud.

Parmi les grands cercles, il y a le \emph{cercle de latitude de ce qui
  est visible sous un climat}.
C'est le grand cercle imaginaire passant par les pôles de l'écliptique et
de l'horizon, donc il est perpendiculaire à l'horizon, et il coupe en
deux moitiés la partie visible et la partie cachée de l'écliptique. L'arc
de ce cercle compris entre un des pôles de l'écliptique et l'horizon
s'appelle \emph{latitude de ce qui est visible sous le climat}.

Parmi les grands cercles, il y a le \emph{cercle de hauteur}.
C'est un grand cercle imaginaire qui passe par les pôles de l'horizon
et par le point donné, donc il est perpendiculaire à l'horizon et le
coupe en deux points dont la distance à la ligne de l'Est et de
l'Ouest est l'azimut de [ce cercle] de hauteur. L'arc de ce cercle
compris entre l'horizon et le point donné est la \emph{hauteur},
et la distance du point au zénith sur ce cercle est le complément de
la hauteur.

Parmi les grands cercles, il a le \uwave{\emph{cercle de l'horizon
  occurrent}}.
C'est un grand cercle qui passe par les deux points Nord et Sud et par
le point donné ou par un astre donné. La \emph{hauteur de l'horizon
  occurrent} est un arc du cercle origine des azimuts, compris entre
l'horizon du lieu et l'horizon occurrent. 

\newpage\phantomsection 
\index{AOBHAQ@\RL{dwr}!BEAOAGAQ@\RL{mdAr}, trajectoire, trajectoire diurne (\textit{i. e.} cercle parallèle à l'équateur)}
\index{BFBGAQ@\RL{nhr}!BBBHAS BFBGAGAQ@\RL{qws al-nhAr}, arc diurne, durée du jour}
\index{AOBHAQ@\RL{dwr}!BEAOAGAQAGAJ AYAQBHAV@\RL{mdArAt al-`rU.d}, trajectoires selon les latitudes (\textit{i. e.} cercles parallèles à l'écliptique)}
\index{BBBFAWAQ@\RL{qn.tr}!BEBBBFAWAQAGAJ@\RL{muqan.tarAt}, muqantars, cercles parallèles à l'horizon}
\includepdf[pages=21,pagecommand={\thispagestyle{plain}}]{edit.pdf}\phantomsection

\noindent La \emph{latitude de l'horizon occurrent} est un arc d'un
grand cercle passant par les pôles de l'équateur et perpendiculaire à
l'horizon occurrent.

Parmi les cercles remarquables, il y a les \emph{trajectoires
  [diurnes]}. Ce sont des cercles parallèles à l'équateur et passant
par le point donné. Leur route est réglée [par rapport à] la
déclinaison et au pôle. Les plus grands d'entre eux sont proches de la
ceinture de la sphère. Dans un pays où la latitude est telle que la
partie visible des trajectoires du Nord est plus grande que la moitié,
et que celle des trajectoires du Sud est plus petite que la moitié,
l'\emph{arc diurne} de l'astre est proportionné à la partie visible de
sa trajectoire~; si la distance de la trajectoire au pôle est
inférieure ou égale à la latitude du lieu, alors l'astre qui la suit
est toujours visible s'il est du côté du pôle visible, et toujours
caché sinon.

Parmi les cercles remarquables, il y a les \emph{trajectoires selon
  les latitudes}. Ce sont des cercles parallèles à la ceinture de
l'écliptique et passant par n'importe quel astre. Chacune des étoiles
fixes tourne autour des pôles de l'écliptique vers l'Est, et elle est
attachée à l'un de ces cercles sans jamais en dévier.

Parmi les cercles remarquables, il y a les \emph{muqantars}. Ce sont
des cercles parallèles à l'horizon~: ceux qui sont au-dessus de
l'horizon s'appellent muqantars de hauteur, et ceux qui sont
en-dessous, muqantars d'abaissement.  Les cercles de hauteur leur sont
perpendiculaires. Ils coupent tous les cercles de hauteur au même
niveau. Leurs pôles sont le zénith et le nadir. Leur forme varie si on
les projette, et elle varie en fonction de la
surface sur laquelle on les projette et du pôle.

\newpage\phantomsection 
\index{ALBHAR@\RL{jwz}!ALBHARAGAB@\RL{al-jawzA'}, les Gémeaux}
\index{BABDBC@\RL{flk}!BABDBC AKAGBEBF@\RL{flk _tAmin}, huitième orbe (étoiles fixes)}
\index{AKAHAJ@\RL{_tbt}!BCBHAGBCAH AKAGAHAJAI@\RL{al-kwAkb al-_tAbitaT}, les étoiles fixes}
\index{AEAHAQANAS@\RL{'ibr_hs}, Hipparque}
\index{ASAMAH@\RL{s.hb}!ASAMAH@\RL{s.hAb j s.hb}, nuage, nébuleuse ??}
\index{AXBDBE@\RL{.zlm}!BCBHBCAH BEAXBDBE@\RL{kwkb m.zlm}, étoile obsure}
\index{AUBHAQ@\RL{.swr}!AUBHAQ@\RL{.sUraT j .suwar}, figure, constellation}
\index{AMBEBD@\RL{.hml}!AMBEBD@\RL{.hml}, Bélier (zodiaque)}
\index{AMBHAJ@\RL{.hwt}!AMBHAJ@\RL{.hwt}, Poissons (zodiaque)}
\index{AKBHAQ@\RL{_twr}!AKBHAQ@\RL{_twr}, Taureau (zodiaque)}
\index{ASAQAW@\RL{sr.t}!ASAQAWAGBF@\RL{sr.tAn}, Cancer (zodiaque)}
\addcontentsline{toc}{chapter}{I.5 Le mouvement des étoiles fixes, et
  une notice sur leurs changements de position causés par les
  déplacements vers l'Est ou vers l'Ouest}
\includepdf[pages=22,pagecommand={\thispagestyle{plain}}]{edit.pdf}\phantomsection
\begin{center}
  \Large Chapitre cinq

  \large Le mouvement des étoiles fixes, et une notice sur leurs
  changements de position causés par les déplacements vers l'Est ou
  vers l'Ouest
\end{center}
Les astres appelés \emph{étoiles fixes} constellent la surface du
\emph{huitième orbe}. Elles ont diverses magnitudes. Il est impossible de les
dénombrer, sauf à Dieu, bien que d'une époque reculée jusqu'à Hipparque
ils en ont observé mille vingt-deux. Ils ont déterminé leurs
longitudes et leurs latitudes à une date donnée. Ils ont classé leurs
magnitudes selon six degrés~: le premier degré est la magnitude la
plus grande, et le sixième, la plus petite. Chaque degré se
différencie en trois genres~: majeur, moyen, mineur. Parmi les mille
vingt-deux étoiles, ils en ont trouvé quinze de première magnitude,
quarante-cinq de deuxième magnitude, deux cent huit de troisième
magnitude, quatre cent soixante-quatorze de quatrième magnitude, deux
cent dix-sept de cinquième magnitude, quarante-neuf de sixième
magnitude, et quatorze en dehors de ces degrés de magnitude. Parmi ces
quatorze, il y a neuf étoiles cachées que l'on appelle \emph{obscures}
et cinq \emph{étoiles nuageuses} comme si elles étaient des taches ou des
morceaux de nuage.

Comme il est difficile de marquer par des noms
cette multitude d'étoiles mais que leurs formes sont fixes, ils ont
imaginé, pour les reconnaître, des figures sur lesquelles ou à côté
desquelles sont situées les étoiles. On dit ainsi <<~le c{\oe}ur du
scorpion~>>, <<~la corne du taureau~>>, <<~la queue de l'ourse~>>
\textit{etc.} pour reconnaître telle étoile par cette nomination et
pour ne pas la confondre avec une autre. Ils ont donc imaginé quarante-huit
figures~: vingt-et-une au Nord de la ceinture du zodiaque, douze sur
la ceinture du zodiaque et quinze au Sud.

Les figures du Nord sont~: la
Petite Ourse, la Grande Ourse, le Dragon, Céphée, le Bouvier, la
Couronne Boréale, Hercule, la Lyre, le Cygne, Cassiopée, Persée, le
Cocher, Ophiuchus, le Serpent, la Flèche, l'Aigle, le Dauphin, le
Petit Cheval, Pégase, Andromède, le Triangle. Les étoiles observées
sur ces figures sont trois cent soixante.

Les figures du zodiaque sont~:
le Bélier, le Taureau, les Gémeaux, le Cancer, le Lion, la Vierge, la
Balance, le Scorpion, le Sagittaire, le Capricorne, le Verseau, les
Poissons. Les étoiles observées sur ces figures ou autour d'elles sont
trois cent quarante-six.

\newpage\phantomsection 
\index{ALAQAQ@\RL{jrr}!BEALAQAQAI@\RL{majarraT}, Voie Lactée}
\index{AHAWBDBEBJBHAS@\RL{b.talimayUs}, Ptolémée}
\index{BFARBD@\RL{nzl}!BEBFAGARBD BBBEAQ@\RL{manAzil al-qamar}, les maisons lunaires}
\index{AMAQBC@\RL{.hrk}!AMAQBCAI AKAGBFBJAI@\RL{.hrkaT _tAnyaT}, deuxième mouvement}
\index{AOBHAQ@\RL{dwr}!BEAOAGAQ@\RL{mdAr}, trajectoire, trajectoire diurne (\textit{i. e.} cercle parallèle à l'équateur)}
\includepdf[pages=23,pagecommand={\thispagestyle{plain}}]{edit.pdf}\phantomsection

Les figures du Sud sont~: la Baleine, Orion,
\'Eridan, le Lièvre, le Grand Chien, le Petit Chien, le Navire Argo,
l'Hydre, la Coupe, le Corbeau, le Centaure, le Loup, l'Autel, la
Couronne Australe, le Poisson Austral. Les étoiles observées sur ces
figures et autour d'elles sont trois cent seize.

Quant aux étoiles
nuageuses, il y en a une sur le poignet de Persée, une deuxième est la
tête d'Orion, une troisième \uwave{...}, une quatrième est celle qui
suit \uwave{le dard (le front ?)} du Scorpion, et une cinquième est un {\oe}il
du Sagittaire.

La \emph{Voie Lactée} est faite d'un très grand nombre de
petites étoiles très rapprochées. De par leur petitesse, leur nombre
et leur proximité, elles apparaissent comme des taches nuageuses dont
la couleur rappelle le lait, et c'est pour cela que Ptolémée l'a nommée Voie
Lactée.

Enfin, il y a vingt-huit \emph{maisons lunaires} dont les étoiles
sont bien identifiées et nommées.

\begin{center}
  \large Section
\end{center}
Le mouvement de ces étoiles est exclusivement en longitude
sans mouvement en latitude. C'est un mouvement d'Ouest en Est autour
des pôles du zodiaque, qu'on appelle \emph{deuxième mouvement}. Sa
grandeur\footnote{c'est-à-dire sa vitesse} est~: selon Ptolémée, un
degré en cent ans~; selon des modernes, un degré et demi en cent ans~;
et selon certains d'entre eux, un degré en soixante-dix ans. Ayant
vérifié ceci par l'observation, j'ai trouvé qu'elles se meuvent d'un
seul degré en soixante-dix années persanes. J'ai embrassé toutes les
observations des anciens et des modernes parvenues jusqu'à nous au
sujet de ces étoiles~; dans mon livre \textit{Ta`l{\=\i}q
  al-ar\d{s}\=ad}, j'ai consigné cela ainsi que tout ce que j'ai
observé et établi.

Cela est établi, et sache que les étoiles dites <<~fixes~>> ne
changent aucunement de trajectoires selon les latitudes, et qu'elles
ne changent pas non plus de positions les unes par rapport aux autres,
ni par rapport aux pôles du zodiaque. Cependant, à cause du deuxième
mouvement, elles changent de positions par rapport à l'équateur. Leurs
trajectoires diurnes changent en tant que leurs distances à l'équateur
augmente ou diminue. Relativement aux habitants du
climat\footnote{\textit{i. e.} relativement à un observateur situé sur
  Terre en un lieu dont la latitude par rapport à l'équateur serait
  donnée.}, elles changent de positions~: ce qui était élevé en
hauteur s'abaisse -- c'est ainsi lorsque sa trajectoire
diurne s'éloigne du zénith -- et inversement. Une étoile peut passer
par le zénith alors qu'elle ne le faisait pas avant, et c'est ainsi
lorsque sa distance à l'équateur céleste atteint la latitude
  du pays (du même côté que le pays) alors qu'elle était auparavant
supérieure ou bien inférieure.

\newpage\phantomsection 
\label{var21}
\includepdf[pages=24,pagecommand={\thispagestyle{plain}}]{edit.pdf}\phantomsection
\noindent Inversement, lorsque sa distance à
l'équateur céleste devient inférieure à la latitude du pays
(resp. supérieure) tandis qu'elle lui était égale et du même côté, alors sa
trajectoire diurne passe du côté du pôle caché (resp. visible) en
s'éloignant du zénith. Une étoile peut devenir toujours visible
(resp. cachée) alors qu'elle ne l'était pas auparavant, et c'est ainsi
lorsque le complément de sa distance à l'équateur céleste devient
inférieur ou égal à la latitude du pays et qu'elle est du côté du pôle
visible (resp. caché) alors qu'il lui était supérieur
auparavant. Auparavant elle avait un lever et un coucher. Quand
l'égalité est atteinte~: elle touche l'horizon à chaque révolution
quand elle passe par le méridien~; elle ne se couche plus si elle est
du côté du pôle visible~; et elle ne se lève plus si elle est du coté
du pôle caché. Dans ce cas sa distance maximale à l'horizon est égale
au double de la latitude du pays. Quand le complément de sa distance
devient inférieur, elle ne touche plus [l'horizon]~; sa distance à
l'horizon est, au minimum, la différence entre la latitude du pays et
le complément de sa distance à l'équateur, et au maximum, la somme de
la latitude du pays et de sa distance à l'équateur. Enfin, il peut
arriver qu'une étoile ait un lever et un coucher alors qu'elle était
toujours visible ou cachée auparavant, et c'est ainsi lorsque le
complément de sa distance à l'équateur dépasse la latitude du pays
alors qu'elle lui était inférieure ou égale auparavant. Ainsi, les
expressions <<~toujours visible~>> ou <<~toujours cachée~>> sont un
abus de langage. Dieu est le plus savant.

\newpage\phantomsection
\index{AMAQBC@\RL{.hrk}!AMAQBCAI AYAQAVBJAI@\RL{.harakaT `r.diyaT}, mouvement par accident}
\index{AMBHBI@\RL{.hw_A}!AMAGBH@\RL{\vocalize .hAwiN (al-.hAwiy)}, contenant}
\index{AMBHBI@\RL{.hw_A}!BEAMBHBI@\RL{m.hw_A}, contenu}
\index{AMAQBC@\RL{.hrk}!AMAQBCAI ANAGAUAI@\RL{.harakaT _hA.saT}, mouvement propre}
\index{ASBABF@\RL{sfn}!ASBABJBFAI@\RL{safInaT}, navire}
\index{ATAHAK@\RL{^sb_t}!AJATAHAHAK@\RL{ta^sabba_ta}, adhérer}
\index{BDARBB@\RL{lzq}!ACBDARBB@\RL{'alzaqa}, agglutiner, coller}
\index{AMAQBC@\RL{.hrk}!BEAJAMAQAQBC BFBAASBG AZBJAQBG@\RL{mt.hrrk binafsihi / bi.gayrihi}, mû par soi / par un autre}
\addcontentsline{toc}{chapter}{I.6 Les mouvements par accident et les autres mouvements}
\label{var25}
\includepdf[pages=25,pagecommand={\thispagestyle{plain}}]{edit.pdf}\phantomsection
\begin{center}
  \Large Chapitre six

  \large Les mouvements par accident et les autres mouvements
\end{center}
Quand le neuvième orbe se meut, les autres orbes se meuvent avec lui
d'Est en Ouest sur ses pôles et autour de son centre. Cela s'appelle
<<~mouvement du contenant au contenu~>>. Les anciens ont convenu qu'il
puisse y avoir mouvement du contenu par le mouvement du contenant, que
ce mouvement soit dans le même sens que celui du contenant ou bien en
sens inverse. Là où ils n'étaient pas d'accord, c'est quand le contenu
se meut \emph{de son mouvement propre autour de l'axe du contenant et
  de son centre}~: peut-il alors se mouvoir par le mouvement du
contenant, ou non~?\label{impuissance}

L'un penche en faveur de cette possibilité en raisonnant comme
suit. Le mouvement du contenant au contenu se fait avec attachement du
mû en son lieu dans celui qui meut, donc il est mû par accident à
cause du mouvement de son lieu, comme le passager du navire est mû à
cause du mouvement du navire. En plus de cela, il se meut de son
mouvement propre, comme le passager qui va et vient sur le navire,
tantôt dans le sens de son mouvement, tantôt en sens inverse. C'est
certes impossible à percevoir, car il est impossible de percevoir deux
mouvements différents dans une même sphère sur une ceinture et ses
pôles~: on perçoit plutôt un seul mouvement, composé de la somme des
deux s'ils vont dans le même sens, ou bien résultant de la différence
entre le plus rapide et le plus lent s'ils ne vont pas dans le même
sens (et on jugerait de même s'ils sont plus
nombreux)\footnote{Al-\d{T}\=us{\=\i} utilise précisément cet argument
  dans la \emph{Tadhkira} (II.2)~: selon lui, il est alors impossible
  que le mouvement de précession des étoiles fixes (huitième orbe) ait
  mêmes pôles que le mouvement diurne (neuvième orbe).}. \`A quoi l'on
répond~: ce n'est pas impossible à percevoir quand il y a
  aussi un astre sur l'orbe supérieur.

L'autre pense que c'est impossible en raisonnant comme suit. Si le
contenu adhère au contenant, alors ce n'est pas par les deux pôles car
ils sont fixes sur l'axe commun. Et ce n'est pas non plus par deux
autres points. S'il y avait toujours adhérence, alors il faudrait qu'il
y ait adhérence des surfaces entières, et ceci interdit au contenu le
mouvement qui lui est propre. Et il est impossible qu'il y ait des
points du contenu qui tantôt s'attachent à des points du contenant,
tantôt s'en séparent. Le mouvement du contenu est moindre que le
mouvement du contenant, non par soi mais à cause du fait qu'ils sont
\emph{collés}. Alors il n'est pas mû par le mouvement du
contenant tout entier. 

\newpage\phantomsection 
\index{ATBJAQAGARBJ@\RL{q.tb al-dIn al-^sIrAzI}, Qu\d{t}b al-D{\=\i}n al-\v{S}{\=\i}r\=az{\=\i}}
\index{ATBEAS@\RL{^sms}!ATBEAS@\RL{^sams}, Soleil}
\index{AMAQBC@\RL{.hrk}!AMAQBCAI BJBHBEBJBJAI@\RL{.harakaT yawumiyyaT}, mouvement diurne}
\index{ACBHAL@\RL{'awj}!AMAQBCAIACBHAL@\RL{.harakaT al-'awj}, mouvement des apogées}
\includepdf[pages=26,pagecommand={\thispagestyle{plain}}]{edit.pdf}\phantomsection

\noindent C'est comme dans la roue située sous une
  charge~; alors le mouvement devrait être sans ordre parce que ce
mouvement n'est dans la nature d'aucune des deux sphères, mais c'est
[en fait] une coïncidence, et il n'y a rien là qui contrevienne à
l'ordre. Pourtant le contenu n'est pas, dès lors qu'il adhère au
contenant, mû par celui-ci~; et il n'adhère pas quand il est mû par
lui~; donc il faut qu'il soit mû par lui sans être mû par lui, et
qu'il adhère sans adhérer~; c'est absurde.\footnote{Celui qui
  <<~pense que c'est impossible~>> fait donc une sorte de
  raisonnement par l'absurde~: la possibilité du mouvement en question
  ne peut s'expliquer que par l'hypothèse de l'<<~adhérence~>>, mais
  l'adhérence empêche le mouvement du contenu par le
  contenant. L'hypothèse de l'adhérence semble se décliner en trois
  hypothèses rejetées tour à tour~:
\begin{itemize}
\item il y a toujours adhérence
\item des points tantôt s'attachent, tantôt se séparent
\item le contenu est collé au contenant et cela freine le mouvement
\end{itemize}} 
[Le même dira~:] en ce qui concerne le mouvement du contenu sur l'axe
du contenant, d'un mouvement qui coïncide avec le mouvement du
contenant ou qui en diffère, je ne vois pas de raison qui empêche
cela. Si les deux mouvements s'équilibrent et que leur sens et leur
durée s'unissent, alors les séparer est inconcevable et inutile.

Voilà quelles sont les opinions avec les indices et les
doutes. Qu\d{t}b al-D{\=\i}n al-\v{S}{\=\i}r\=az{\=\i} a dit que la
preuve de la possibilité est plus flagrante et plus accessible mais à
condition de distinguer des cas. [Il a dit que] la preuve de
l'impossibilité n'échappe pas à l'objection parce qu'elle est fondée
sur le fait que quand les centres sont confondus il faut qu'il y ait
adhérence, or cela n'est pas nécessaire.

\begin{center}
  \large Section
\end{center}
Je ne vois pas comment il serait possible que le contenu soit mû par
un contenant, ainsi que par soi-même, ainsi que par un autre contenant
sur d'autres pôles et dans l'autre sens. Montrons cela. L'orbe
parécliptique du Soleil est mû par le mouvement diurne d'Est en Ouest,
et il se meut lui-même sur les deux pôles de l'écliptique d'Ouest en
Est, d'un déplacement égal au mouvement du centre du Soleil ou du
Soleil moyen. Or son mouvement par le mouvement du huitième orbe, en
vue du mouvement de l'Apogée, n'est pas permis~: il y a besoin d'un
autre orbe qui le meuve vers l'Est sur les pôles de l'écliptique, d'un
déplacement égal au mouvement de l'Apogée.

\newpage\phantomsection
\index{ATBEAS@\RL{^sms}!ATBEAS@\RL{^sams}, Soleil}
\includepdf[pages=27,pagecommand={\thispagestyle{plain}}]{edit.pdf}\phantomsection
\noindent \`A supposer qu'il soit permis, il faudrait que le reste des
orbes des autres astres se meuve du même mouvement, et le
parécliptique du Soleil ne pourrait se mouvoir comme cela sans que se
meuvent ainsi les autres orbes, ceux de Saturne et de Jupiter. C'est
subtil, prête attention~!\footnote{{\shatir } s'intéresse ici à la
  composition de trois mouvements concentriques~: un mouvement propre,
  le mouvement d'un contenant, et le mouvement d'un autre
  contenant. Le Soleil est entraîné par un mouvement annuel ayant pour
  origine l'Apogée, dû au mouvement de l'<<~orbe parécliptique du
  Soleil~>>. Mais il est aussi entraîné par le mouvement diurne du
  neuvième orbe le contenant, et par le mouvement de l'Apogée. Le
  mouvement de l'Apogée doit être causé par le mouvement d'un autre
  orbe le contenant. Si cet orbe est le huitième (animé du mouvement
  de précession des fixes), alors les orbes intermédiaires (dont ceux
  de Jupiter et Saturne) doivent aussi être entraînés par ce mouvement
  comme par le mouvement diurne. Et sinon, alors il faut concevoir un
  autre orbe.}

\emph{Remarque}. Si nous supposons que le Soleil a un épicycle et un
parécliptique comme dans la configuration connue\footnote{Le modèle du
  Soleil avec un épicycle, chez Ptolémée.}, alors le parécliptique du
Soleil se meut de lui-même sur les pôles de l'écliptique comme le
centre du Soleil, l'épicyle se meut en haut de cet orbe en sens
inverse, et le neuvième orbe meut ces deux orbes d'Est en Ouest.  Or
le huitième est impuissant à mouvoir le parécliptique d'Ouest en Est
d'un déplacement égal au mouvement de l'Apogée, comme nous l'avons
dit. Il y a donc besoin d'un autre orbe qui meuve le parécliptique,
d'un déplacement égal au mouvement de l'Apogée. Il n'y a personne qui
dise que l'Apogée n'est mû ni selon l'hypothèse que son mouvement est
semblable à celui du huitième orbe, ni selon l'hypothèse que le
huitième orbe le meut par accident dans le sens des signes, que le
neuvième le meut dans le sens inverse sur des pôles distincts de ceux
du premier mouvement, et qu'en plus de cela l'Apogée se meut par
lui-même~; alors il y a besoin d'ajouter un autre orbe qui meuve
l'Apogée.\footnote{{\shatir} semble simplement répéter le raisonnement
  du paragraphe précédent. La structure de tout ce paragraphe et du
  précédent est confuse.}

\emph{Autre remarque.} Toutes les observations vraies témoignent que le
mouvement de l'Apogée est plus rapide que le mouvement du huitième
orbe. Il faut dès lors ajouter un autre orbe qui meuve l'Apogée d'un
déplacement égal au mouvement observé.

\emph{Autre remarque.} Si on suppose que le Soleil a un parécliptique,
un épicyle et un autre orbe qui meut l'Apogée, avec ce principe, la
longitude calculée ne coïncide pas avec l'observation vraie, car nous
avons fait des observations précises dans les octants et nous avons
trouvé que cela diffère de la longitude calculée sur le principe
ci-dessus. Et nous avons déjà démontré cela dans le livre
\textit{Ta`l{\=\i}q al-ar\d{s}\=ad}. Donc il nous faut un principe
simple qui coïncide avec les observations~; Dieu a fait trouver un
principe qui suffit à l'effet visé, et [c'est un principe] que tu vas
connaître, si Dieu le veut.

\newpage\phantomsection
\index{AQBCAH@\RL{rkb}!ALASBE BEAQBCAH@\RL{jsm mrkb}, corps composé}
\index{AMAQBC@\RL{.hrk}!AMAQBCAI AHASBJAWAI BEAQBCBCAHAI@\RL{.harakaT basI.taT mrkkbaT}, mouvement simple-composé}
\includepdf[pages=28,pagecommand={\thispagestyle{plain}}]{edit.pdf}\phantomsection

\emph{Autre remarque}. Les astronomes n'ont pas été d'accord au sujet
du mouvement des petits orbes qui ne tournent pas autour du centre du
monde, tel l'orbe de l'épicycle et d'autres semblables. Ils ont tous
accepté qu'ils puissent tourner dans une direction ou une autre, car ils
ont signalé que l'orbe de l'épicycle avait une moitié supérieure et une
moitié inférieure~; si bien que, s'il tourne dans le sens des signes
dans la moitié supérieure, alors il tournera dans le sens contraire dans
la partie inférieure, et inversement. Son mouvement ne serait ni un
mouvement violent, ni un mouvement par accident, mais ce serait un
mouvement naturel. Ils ont tous accepté qu'il puisse y avoir des
épicycles dans les orbes autres que le neuvième. En effet, dans les
orbes, nous voyons des astres, or l'astre, dans son orbe, montre une
certaine \emph{composition}. A quiconque dirait que les orbes sont tous
simples, qu'il ne peut donc exister d'épicycle en eux, et qu'aucun
mouvement excentrique ne peut donc être un mouvement simple, je dirais
que l'existence et le mouvement des épicycles a déjà été prouvé. S'il
y avait une preuve décisive qui empêche cela, les orbes n'en
demeureraient pas moins des objets composés. Il n'y a pas de
simplicité dans les orbes. Je pense qu'il s'agit d'une composition de
choses simples et non d'éléments, sauf en ce qui concerne le neuvième
orbe, mais Dieu est le plus savant.\footnote{Justification des
  épicycles~: bien qu'ils causent une certaine <<~composition~>>, de
  toute façon les orbes ne sont jamais absolument <<~simples~>> (sauf
  le neuvième qui ne contient aucun astre...).}

\newpage\phantomsection  
\index{AOBHAQ@\RL{dwr}!BEAQAOAGAQAGAJBEAQAGBCARBCAQ@\RL{madArAt marAkaz al-akr}, trajectoires des centres des orbes}
\index{BABDBC@\RL{flk}!BABDBCBEBEAKAKBD@\RL{falak muma_t_tal}, parécliptique}
\index{BABDBC@\RL{flk}!BABDBCAMAGBEBD@\RL{falak .hAmil}, orbe déférent}
\index{BABDBC@\RL{flk}!BABDBCBEAOBJAQ@\RL{falak mdIr}, orbe rotateur}
\index{BABDBC@\RL{flk}!BABDBCATAGBEBD@\RL{falak ^sAmil}, orbe total}
\index{AQBCAR@\RL{rkz}!BEAQBCAR@\RL{markaz}, centre}
\index{ALASBE@\RL{jsm}!BCAQAGAJ BEALASBEAI AJAGBEBEAI@\RL{kraT mjsm, falak mjsm}, sphère solide, orbe solide}
\index{AMAQBC@\RL{.hrk}!AMAQBCAI BEAQBCAR ATBEAS@\RL{.harakaT markaz al-^sams}, mouvement du centre du Soleil}
\index{AMAQBC@\RL{.hrk}!AMAQBCAI BHASAW ATBEAS@\RL{.harakaT ws.t al-^sams}, mouvement du Soleil moyen}
\addcontentsline{toc}{chapter}{I.7 L'orbe du Soleil et ses mouvements}
\includepdf[pages=29,pagecommand={\thispagestyle{plain}}]{edit.pdf}\phantomsection

\begin{center}
  \Large Chapitre sept

  \large L'orbe du Soleil et ses mouvements
\end{center}
Nous évoquerons d'abord les trajectoires des centres des sphères
solides parce qu'elles sont plus faciles à représenter sur le
plan. Nous imaginons un orbe dans le plan de l'écliptique et en son
centre, et nous l'appelons le \emph{parécliptique}. Nous imaginons
un autre orbe dont le centre soit sur la périphérie de l'orbe
parécliptique et dont le rayon soit quatre degrés et trente-sept
minutes (si le rayon du parécliptique mesure soixante parts), et nous
l'appelons \emph{orbe déférent}. Sur la périphérie du déférent, il y
a le centre d'un autre orbe dont le rayon est de deux degrés et demi
des mêmes parts. On l'appelle \emph{orbe rotateur}. Et il y a un
autre orbe au centre du monde dont la périphérie est au-dessus de ces
orbes. Son rayon intérieur est la somme des rayons des trois orbes
mentionnés auxquels on ajoute le rayon du globe solaire. Le rayon du
globe solaire est un sixième d'une des mêmes parts. Donc le rayon [de
  l'orbe total] est soixante-sept degrés et dix-sept minutes, et son
épaisseur est d'une grandeur telle qu'il peut se mouvoir vers l'Est
comme le mouvement de l'Apogée. Nous supposons que son rayon est
soixante-huit degrés, alors son épaisseur est quarante-trois
minutes. On l'appelle \emph{orbe total}.

A présent, supposons que les trois orbes sont sur la ligne issue du
centre du zodiaque et passant par l'Apogée (voir figure
suivante). Quant aux mouvements, l'orbe parécliptique se meut d'un
mouvement simple autour du centre du zodiaque et sur ses pôles dans le
sens des signes. Ce mouvement est comme le mouvement du centre du
Soleil, c'est-à-dire de $0;59,8,9,51,46,57,32,3$ degrés par
jour\footnote{Si l'on utilise cette valeur pour calculer le Soleil
  moyen (Soleil moyen $=$ centre $+$ Apogée), les trois derniers rangs
  sexagésimaux diffèrent de la valeur donnée dans la suite du
  texte. On retiendra donc, comme valeur du mouvement du centre,
  $0;59,8,9,51,47$ degrés par jour. Le manuscrit C corrige d'ailleurs
  cette erreur (voir apparat critique dans l'arabe).}. Le centre du
déférent se déplace donc d'autant dans le sens des signes. L'orbe
déférent se meut d'un mouvement uniforme en son centre, dans le sens
inverse des signes quand on est dans sa partie supérieure. Ce
mouvement est aussi comme celui du centre du Soleil, c'est-à-dire
$0;59,8,9,51,46,57,32,3$. L'orbe rotateur se meut d'un mouvement
uniforme en son centre, dans le sens des signes quand on est dans sa
partie supérieure, comme deux fois le centre du Soleil, c'est-à-dire
$1;58,16,19,43,34$ degrés par jour. Le centre du corps du Soleil se
déplace donc d'autant, dans le sens des signes. 

\newpage\phantomsection  
\index{BHASAW@\RL{ws.t}!BHASAW@\RL{ws.t}, astre moyen}
\index{AQBCAR@\RL{rkz}!BEAQBCAR@\RL{markaz}, centre}
\index{ACBHAL@\RL{'awj}!ACBHAL@\RL{'awj}, Apogée}
\index{AMAVAV@\RL{.h.d.d}!AMAVBJAV@\RL{.ha.dI.d}, périgée}
\index{AYAOBD@\RL{`dl}!AJAYAOBJBD ASBEAS@\RL{ta`dIl al-^sams}, équation du Soleil}
\index{ACBHAL@\RL{'awj}!AMAQBCAIACBHAL@\RL{.harakaT al-'awj}, mouvement des apogées}
\index{AMAQBC@\RL{.hrk}!AMAQBCAI AHASBJAWAI BEAQBCBCAHAI@\RL{.harakaT basI.taT mrkkbaT}, mouvement simple-composé}
\index{ANBDBA@\RL{_hlf}!ACANAJBDAGBA@\RL{i_htilAf}, irrégularité, anomalie, variation}
\includepdf[pages=30,pagecommand={\thispagestyle{plain}}]{edit.pdf}\phantomsection

\noindent L'orbe total se meut d'un mouvement simple autour de son
centre, d'Ouest en Est, comme le mouvement de l'Apogée, c'est-à-dire
de $0;0,0,9,51,46,51$ degrés par jour, ce qui fait une minute par an,
ou un degré tous les soixante ans\footnote{Soixante années persanes.},
par l'observation véritable. Observe et corrige, alors la distance du
centre de l'orbe déférent
\footnote{Le texte porte ``épicycle'', mais il faut bien lire
  ``déférent''.}  au point fixe du zodiaque devient la somme des deux
mouvements. Leur somme est le [mouvement du Soleil] moyen, et c'est un
mouvement simple-composé.

Alors la distance du Soleil au centre du monde est maximale quand il
est à l'Apogée, minimale quand il est au périgée, et moyenne quand le
centre est $3$ signes et $7;7,1,6,8,13$ degrés\footnote{Ou bien
  faut-il lire <<~$7;7$ ou $8;13$~>>~? Mais ce résultat est erroné. La
  distance moyenne est atteinte lorsque le Soleil est sur la ceinture
  du parécliptique, c'est-à-dire lorsque la distance au centre du
  monde est de 60 parts. La résolution d'une équation de degré 2
  montre que cette distance est atteinte lorsque le centre est de 3
  signes et $0;18,2,25,27,47,11$ degrés. \`A moins qu'{\shatir } donne
  ici une valeur en signes et en \emph{heures} (avec 1 h =
  $0°2'30''$)~: si le cosinus du centre est arrondi à $0;0,18,40$, on
  trouve en effet 3 signes et $7;7,48$ heures.}. Trace une droite
joignant le centre du zodiaque au centre du déférent, et une autre
droite joignant le centre du zodiaque au corps du Soleil. Elles
cernent l'angle d'anomalie du Soleil, ou équation [du Soleil]. Elles
sont confondues quand le centre de l'orbe déférent est à l'Apogée ou
au périgée. Et le maximum de cette anomalie est $2;2,6$ degrés. Il est
atteint quand le centre vaut trois signes et sept degrés, ou huit
[signes] et vingt-trois degrés\footnote{Il suffit pour le voir de
  construire une table de l'équation du Soleil en fonction du centre,
  comme avait pu en construire {\shatir}. Si on construit une telle table
  pour des valeurs du centre variant de demi-degré en demi-degré,
  on constate que la correction maximale est atteinte
  quand le centre est entre $96;30$ et $97$ degrés, et qu'elle est
  alors supérieure à $2;2,5,10$. La comparaison des premières
  différences montre que la correction ne peut, sur cet intervalle,
  varier de plus de $0;0,0,38$ degrés quand le centre varie de $0;30$
  degrés. La correction maximale est donc inférieure à $2;2,5,48$. Une
  table plus précise montrerait que la correction maximale est de
  $2;2,5,13$. \textit{Cf.} aussi le graphe de l'équation du Soleil dans notre
  commentaire \textit{supra} figure \ref{fig005} p.~\pageref{fig005}.}.
La droite passant par le centre du déférent coupe
l'écliptique en le \emph{Soleil moyen}. C'est-à-dire que la distance
entre cette droite et la droite issue du zodiaque et passant par
l'équinoxe de printemps est le Soleil moyen. Sa distance à la
droite passant par l'Apogée est son \emph{centre}. La distance entre
la droite passant par l'Apogée et la droite passant par le
commencement du Bélier\footnote{Ibn al-Sh\=a\d{t}ir veut dire
  l'équinoxe de printemps. A cause du <<~deuxième mouvement~>> (le
  mouvement de précession), le Bélier se déplace par rapport à
  l'équinoxe de printemps...} est l'\emph{Apogée}. La somme de
l'Apogée et du centre est le Soleil moyen.

\newpage\phantomsection  
\index{BFAWBB@\RL{n.tq}!BEBFAWBBAI@\RL{min.taqaT}, ceinture}
\index{BBAWAH@\RL{q.tb}!BBAWAH@\RL{q.tb}, pôle}
\index{ACBHAL@\RL{'awj}!AMAQBCAIACBHAL@\RL{.harakaT al-'awj}, mouvement des apogées}
\index{ALASBE@\RL{jsm}!BCAQAGAJ BEALASBEAI AJAGBEBEAI@\RL{kraT mjsm, falak mjsm}, sphère solide, orbe solide}
\includepdf[pages=31,pagecommand={\thispagestyle{plain}}]{edit.pdf}\phantomsection

\noindent  La distance maximale du
centre du Soleil au centre du monde est soixante-sept degrés et sept
minutes, en parts telles que le rayon du parécliptique en compte
soixante. La distance minimale du centre du Soleil au centre du monde
est cinquante-deux degrés et cinquante-trois minutes, en les mêmes
parts.

Le Soleil moyen à midi du premier jour de l'année sept cent un de
l'ère de Yazdgard 9 signes et $10;9,0$ degrés à Damas et l'Apogée du
Soleil à la date indiquée est 2 signes et $29;52,3,1$ degrés. Le
mouvement de l'Apogée dans le sens des signes, en soixante années
persanes, est d'un seul degré.
Le mouvement du Soleil [moyen] en une
année persane est $11$ [signes] et $29;45,40$ degrés, en un mois
c'est-à-dire trente jours, $0$ [signe] et $29;34,9,51,46,50,57,32$
degrés, en un jour $0;59,8,19,43,33,41,55,4,6$, et en une heure égale
$0;2,27,50,49,18,54,14,47$. \^Otons du Soleil moyen l'Apogée, il reste
le centre.

Nous imaginons l'orbe parécliptique comme un orbe sphérique solide dont la
partie concave touche la partie convexe de la totalité des orbes de
Vénus et dont la partie convexe touche la partie concave de l'orbe
total du Soleil~; son centre est le centre du monde, il a deux pôles
alignés avec les pôles de l'écliptique, et sa ceinture est dans le
plan de la ceinture de l'écliptique.

\label{complement_sol}
Son épaisseur est de quatorze degrés et trente-quatre minutes, plus
un complément. Pous que [les orbes] soient mutuellement contigus
quand on y plonge l'orbe déférent, nous prenons un
complément d'un tiers de degré à l'extérieur et à l'intérieur, en
parts telles que son rayon extérieur est soixante-sept parts et
dix-sept minutes. Nous prenons une sphère dont le diamètre est quatorze
degrés et trente-quatre minutes, contenue en lui et le touchant à
l'intérieur et à l'extérieur. Nous prenons une autre sphère, de rayon
deux degrés, un demi-degré et un sixième de degré des mêmes parts,
contenue dans la
sphère du déférent de sorte qu'elle la touche en un point. Nous
prenons le corps du Soleil contenu dans cette sphère. Nous prenons les
ceintures de ces sphères dans le plan de l'orbe écliptique. C'est ce que nous
avons représenté dans la figure plane.

\newpage\phantomsection  
\index{AHAWBDBEBJBHAS@\RL{b.talimayUs}, Ptolémée}
\index{AOBHAQ@\RL{dwr}!AJAOBHBJAQ@\RL{tadwIr}, épicycle}
\index{ALAOBD@\RL{jdl}!ALAOBHBD@\RL{jdwl}, table}
\index{AYAOBD@\RL{`dl}!AJAYAOBJBD@\RL{ta`dIl}, équation}
\index{BBBHBE@\RL{qwm}!AJBBBHBJBE ATBEAS@\RL{taqwIm al-^sams}, Soleil vrai}
\label{var5}\label{var15}\label{var16}
\includepdf[pages=32,pagecommand={\thispagestyle{plain}}]{edit.pdf}\phantomsection

Les bordant, il y a une sphère dont le centre est le centre du monde,
dont la partie concave touche la partie convexe du parécliptique du
Soleil, dont l'épaisseur est de quarante-trois minutes, dont la partie
convexe touche la partie concave du parécliptique de Mars, dont la
circonférence est dans le plan de l'écliptique et dont les deux pôles
sont les pôles de l'écliptique. Les mouvements sont comme on l'a dit
précédemment. Ainsi le cas du Soleil est mené à bien par quatre
orbes. Il y en a deux qui entourent le centre du monde et sont dans le
plan de l'écliptique~; et il y a le déférent et le rotateur. Si la
correction s'accordait au principe de l'épicycle selon la
configuration de Ptolémée dans les octants, il faudrait ajouter un
autre orbe mouvant l'Apogée d'un degré tous les soixante ans. Mais
cela contredit l'observation et ne s'y accorde qu'en les équinoxes et
les solstices et non ailleurs. De plus ces orbes ont un même centre,
et le centre [du Soleil], le Soleil moyen et l'Apogée sont des
quantités prises dans un même cercle~; or dans la configuration
connue\footnote{un modèle du Soleil avec un excentrique},
le Soleil moyen est une
quantité prise dans un certain cercle, et le centre et l'Apogée sont
des quantités prises dans un autre cercle~; et ceci est erroné, ils en
ont fait un arc semblable, et c'est impossible\footnote{Dans les modèles avec
  un excentrique, il faut bien faire attention à quel centre se rapporte chaque
  paramètre angulaire (par rapport au centre de l'excentrique, ou non).
  Chez {\shatir} la question ne se pose pas, car à chaque orbe ne correspond
  qu'un unique centre. Mieux encore, les paramètres \textit{centre},
  \textit{Apogée} et \textit{Soleil moyen} désigneront tous des angles mesurés
  le long de cercles concentriques ; c'est plus compliqué chez Ptolémée.}.

\begin{center}
  \large Section
\end{center}
\emph{Comment trouver le Soleil vrai par les tables et par le calcul.}
Si tu prends le Soleil moyen à la date que tu veux, que tu
prends aussi l'Apogée, et que tu ôtes l'Apogée du Soleil moyen, il
reste le centre. Donc tu tires le Soleil moyen et l'Apogée de la
table, à la date que tu veux. Pour cela, tu prends ce qui appartient
aux années entières, aux années intercalaires, aux mois et aux jours
-- ne compte pas le jour dans lequel tu es -- puis ajoute rang à rang
et augmente ce qui doit être augmenté\footnote{addition sexagésimale
  avec retenues}~; rejette les circonférences entières~; reste alors
le Soleil moyen. Après cela, fait l'addition pour l'Apogée, puis ôtes
l'Apogée du Soleil moyen, il reste le centre. Entre-le alors dans la
table de l'équation du Soleil. Ce que tu y trouves en degrés et
fractions de degré, corrigé en interpolant ce qu'il y a entre les
lignes, c'est l'équation du Soleil. Ensuite, si le centre tombe
dans la partie supérieure de la table, tu soustrais l'équation au
Soleil moyen, et s'il tombe dans la partie inférieure de la table, tu
ajoutes l'équation au Soleil moyen.

\newpage\phantomsection  
\includepdf[pages=33,pagecommand={\thispagestyle{plain}}]{edit.pdf}\phantomsection

\noindent On obtient alors le Soleil
vrai, par rapport à l'orbe écliptique, à midi du jour pour lequel tu
as calculé, à Damas.

Et si tu veux calculer l'équation, par le calcul, sans table~:
si c'est cela que tu veux, alors multiplie le cosinus du centre du
Soleil par sept degrés et sept minutes.
Ce qu'on obtient, ajoute-le à soixante si le
centre est du côté de la distance maximale, c'est-à-dire quand le
centre est inférieur à trois signes, sept degrés et sept
minutes\footnote{il faut lire seulement~: trois signes, c'est-à-dire 90°.},
ou bien supérieur à huit signes, vingt-deux degrés et cinquante-trois
minutes\footnote{il faut lire plutôt~: neuf signes, c'est-à-dire 270°.}~;
sinon, il est du côté de la distance minimale, alors retranche ce qu'on
obtient de soixante. Le somme de l'addition, ou bien le reste de la
soustraction, on l'appelle le total.
Ensuite, multiplie le sinus du centre par deux
degrés et sept minutes, on appelle le résultat dividende. Mets au
carré le dividende et mets au carré le total, fais la somme des deux
carrés, prends-en la racine, c'est la distance du Soleil au centre du
monde (en parts telles que le rayon du parécliptique en compte
soixante). Divise le dividende par la distance du Soleil au centre du
monde, il en sort le sinus de l'équation du Soleil. Cherche son arc
dans la table des sinus, il en sort l'équation [du Soleil]. Nous
avons fait ce calcul pour chaque degré, et nous l'avons établi dans la
table. Ensuite, si le centre est inférieur à six signes, tu soustrais
l'équation au Soleil moyen, et s'il est supérieur à six signes, tu
ajoutes l'équation au Soleil moyen. Tu obtiens la position du
Soleil par rapport à l'écliptique.

\begin{center}
  \large Section
\end{center}
\label{diam_sol} Le \emph{diamètre du Soleil} à distance
moyenne (c'est-à-dire quand sa distance au centre du monde est comme
le rayon du parécliptique qui est par supposition de soixante parts)
est $32;32$. Son diamètre à l'Apogée est $29;5$, et au périgée
$36;55$.

\newpage\phantomsection  
\index{AHBGAJ@\RL{bht}, vitesse apparente}
\index{AHAWBDBEBJBHAS@\RL{b.talimayUs}, Ptolémée}
\index{AEAHAQANAS@\RL{'ibr_hs}, Hipparque}
\index{BGBFAO@\RL{hnd}!BGBFAO@\RL{al-hnd}, les Indiens}
\label{var6}
\includepdf[pages=34,pagecommand={\thispagestyle{plain}}]{edit.pdf}\phantomsection

La méthode pour calculer son diamètre ailleurs qu'en ces lieux est que
tu divises son diamètre à distance moyenne par sa distance au centre
du monde qui est le centre de l'écliptique\footnote{Notons $\alpha$ le
  diamètre apparent du Soleil (c'est un angle), et OS sa distance au
  centre du monde. \`A distance moyenne, OS = 60 et $\alpha=0°32'32''$
  (probablement une donnée de l'observation), or
  $\alpha\simeq\sin\alpha$ est inversement proportionnel à OS. On a
  donc, aux autres distances~: $\alpha=\dfrac{0°32'32''\times
    60}{\text{OS}}$.}. Son diamètre à distance moyenne est selon
Hipparque $32;45$, chez les Indiens $32;33$, chez les modernes
$32;35$, chez Ptolémée dans son chapitre\footnote{\textit{L'Almageste}, livre 5 chapitre 14.} concernant les grandeurs $31;20$,
et chez moi on a trouvé par l'observation
$32;32$. Pour cela il y a une méthode approchée qui est que tu
multiplies sa vitesse apparente \footnote{Le mot \textit{buht} que nous
  traduisons par ``vitesse apparente'' est utilisé pour chaque astre dans le
  même contexte et semble désigner la vitesse angulaire mesurée par
  rapport aux étoiles fixes, c'est-à-dire l'arc parcouru par l'astre le
  long d'un grand cercle du huitième orbe, par unité de temps.
  Il ne faut sûrement pas y voir un concept défini de manière rigoureuse.
  Cette grandeur est, grossièrement,
  inversement proportionnelle à OS ; elle est donc, grossièrement,
  proportionnelle au diamètre apparent de l'astre. Pour le Soleil,
  on s'en convainc aisément en pensant à un modèle (approximatif)
  avec un excentrique.} par trente-trois, il en
sort son diamètre approché~; mais la première méthode est exacte. Dieu
est le plus savant.

\newpage\phantomsection  
\includepdf[pages=35,pagecommand={\thispagestyle{plain}}]{edit.pdf}\phantomsection

\begin{center}\label{soleil_trajectoires}
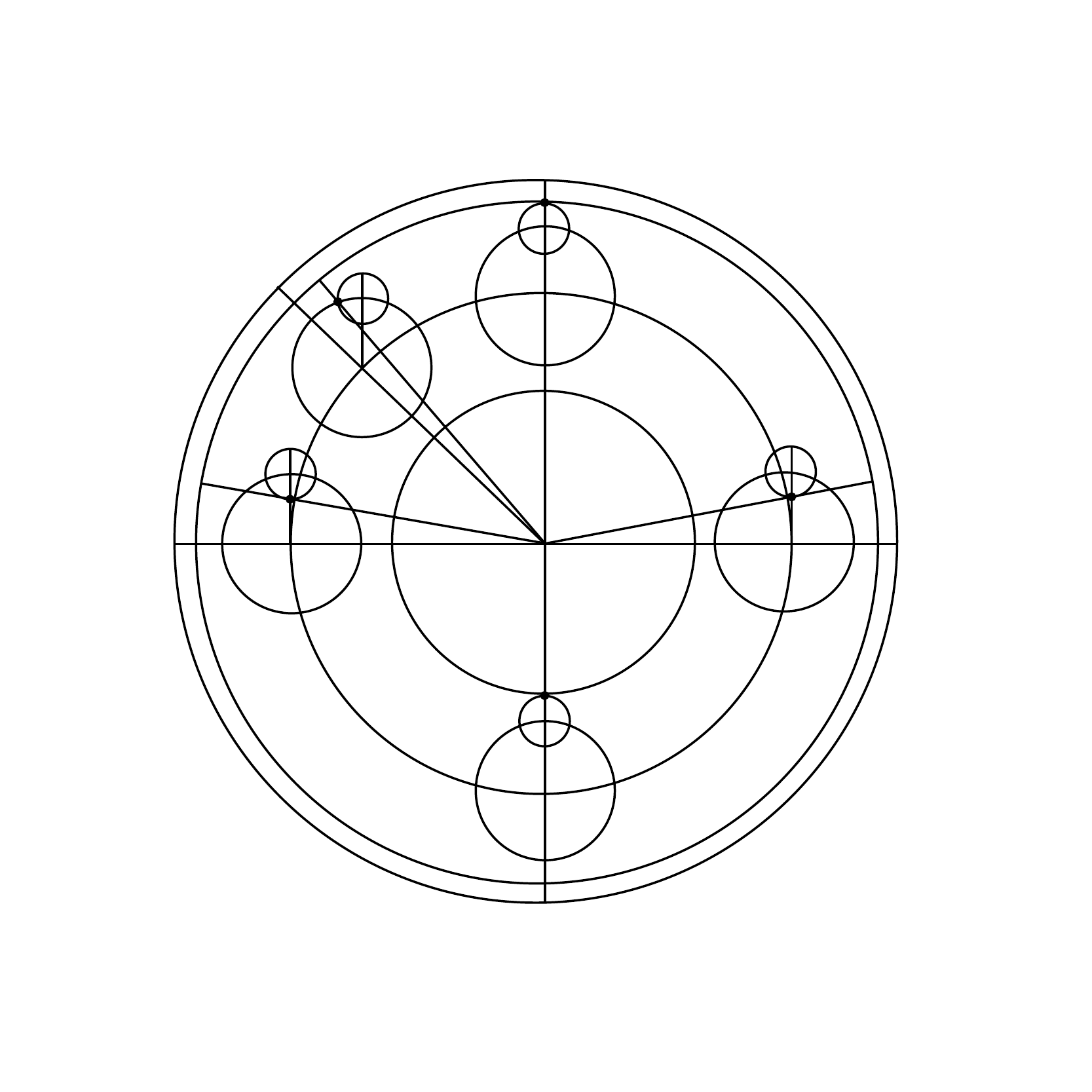
Les orbes du Soleil figurés par les trajectoires des centres des sphères
\end{center}

\newpage\phantomsection  
\includepdf[pages=36,pagecommand={\thispagestyle{plain}}]{edit.pdf}\phantomsection

\begin{center}\label{soleil_orbes_solides}
\begingroup%
  \makeatletter%
  \providecommand\color[2][]{%
    \errmessage{(Inkscape) Color is used for the text in Inkscape, but the package 'color.sty' is not loaded}%
    \renewcommand\color[2][]{}%
  }%
  \providecommand\transparent[1]{%
    \errmessage{(Inkscape) Transparency is used (non-zero) for the text in Inkscape, but the package 'transparent.sty' is not loaded}%
    \renewcommand\transparent[1]{}%
  }%
  \providecommand\rotatebox[2]{#2}%
  \ifx\svgwidth\undefined%
    \setlength{\unitlength}{480.1bp}%
    \ifx\svgscale\undefined%
      \relax%
    \else%
      \setlength{\unitlength}{\unitlength * \real{\svgscale}}%
    \fi%
  \else%
    \setlength{\unitlength}{\svgwidth}%
  \fi%
  \global\let\svgwidth\undefined%
  \global\let\svgscale\undefined%
  \makeatother%
  \begin{picture}(1,1)%
    \put(0,0){\includegraphics[width=\unitlength]{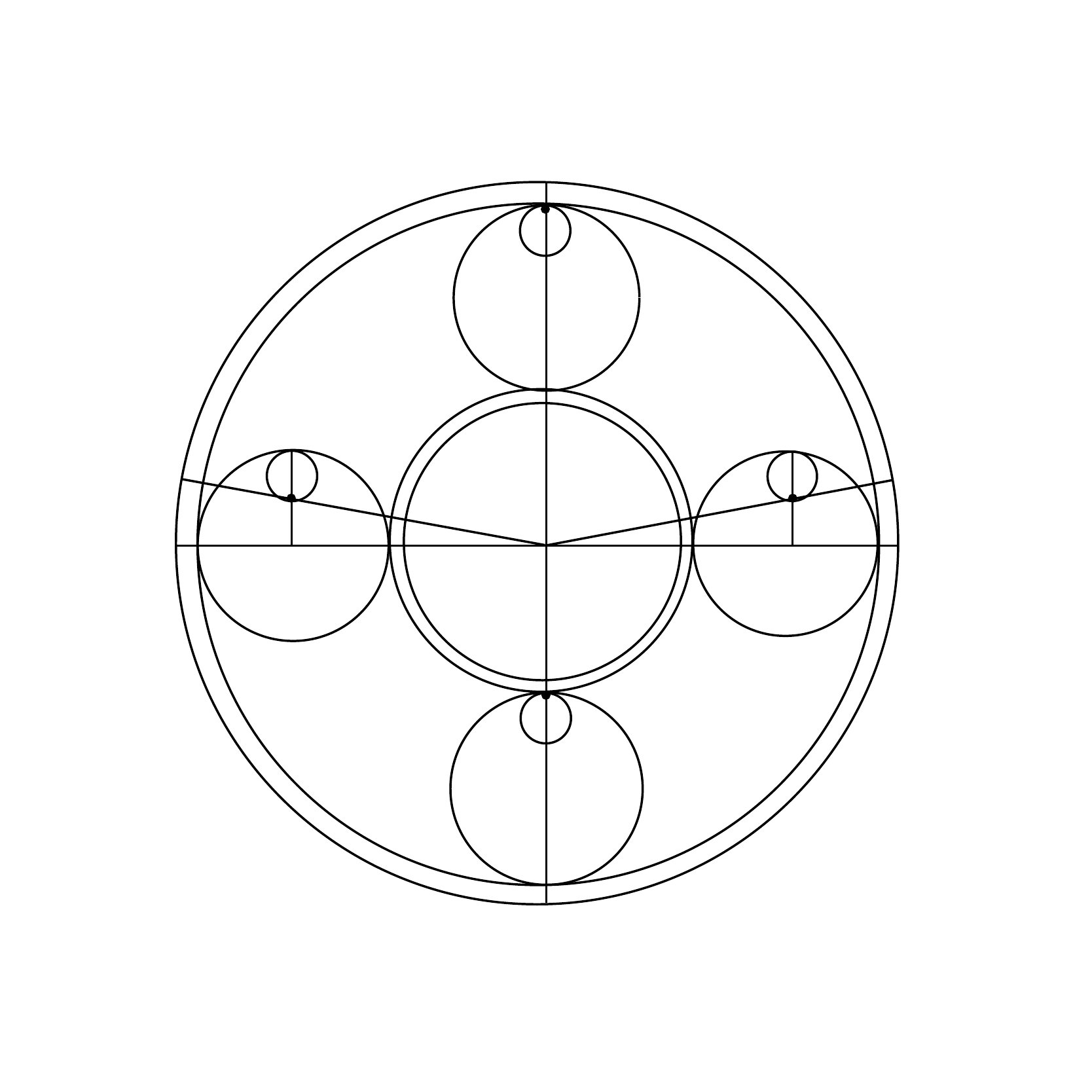}}%
    \put(0.35215843,0.54725943){\color[rgb]{0,0,0}\rotatebox{59.99999989}{\makebox(0,0)[lb]{\smash{partie concave}}}}%
    \put(0.82977474,0.49684524){\color[rgb]{0,0,0}\makebox(0,0)[lb]{\smash{Soleil moyen}}}%
    \put(0.03872766,0.49481337){\color[rgb]{0,0,0}\makebox(0,0)[lb]{\smash{Soleil moyen}}}%
    \put(0.73131559,0.52833371){\color[rgb]{0,0,0}\makebox(0,0)[lb]{\smash{Soleil}}}%
    \put(0.43253079,0.35124523){\color[rgb]{0,0,0}\makebox(0,0)[lb]{\smash{corps du Soleil}}}%
    \put(0.53946873,0.79214937){\color[rgb]{0,0,0}\rotatebox{-123.87356031}{\makebox(0,0)[lb]{\smash{orbe rotateur}}}}%
    \put(0.82304183,0.56164986){\color[rgb]{0,0,0}\makebox(0,0)[lb]{\smash{longitude}}}%
    \put(0.60638109,0.7300403){\color[rgb]{0,0,0}\rotatebox{-117.64843093}{\makebox(0,0)[lb]{\smash{orbe déférent}}}}%
    \put(0.20768515,0.59173845){\color[rgb]{0,0,0}\rotatebox{59.99999989}{\makebox(0,0)[lb]{\smash{partie convexe de l'orbe total}}}}%
    \put(0.17418188,0.61619639){\color[rgb]{0,0,0}\rotatebox{75.00000215}{\makebox(0,0)[lb]{\smash{orbe parécliptique}}}}%
    \put(0.47187343,0.38565528){\color[rgb]{0,0,0}\makebox(0,0)[lb]{\smash{périgée}}}%
    \put(0.48615614,0.8403225){\color[rgb]{0,0,0}\makebox(0,0)[lb]{\smash{Apogée}}}%
    \put(0.45402,0.48206377){\color[rgb]{0,0,0}\makebox(0,0)[lb]{\smash{centre de l'écliptique}}}%
    \put(0.3813091,0.54168651){\color[rgb]{0,0,0}\rotatebox{59.99999989}{\makebox(0,0)[lb]{\smash{de l'orbe total}}}}%
    \put(0.82418102,0.5415752){\color[rgb]{0,0,0}\makebox(0,0)[lb]{\smash{du Soleil}}}%
    \put(0.06156641,0.56981677){\color[rgb]{0,0,0}\makebox(0,0)[lb]{\smash{longitude}}}%
    \put(0.06270559,0.54974211){\color[rgb]{0,0,0}\makebox(0,0)[lb]{\smash{du Soleil}}}%
  \end{picture}%
\endgroup%

Les orbes solides du Soleil
\end{center}

\begin{center}
\includegraphics[width=\textwidth]{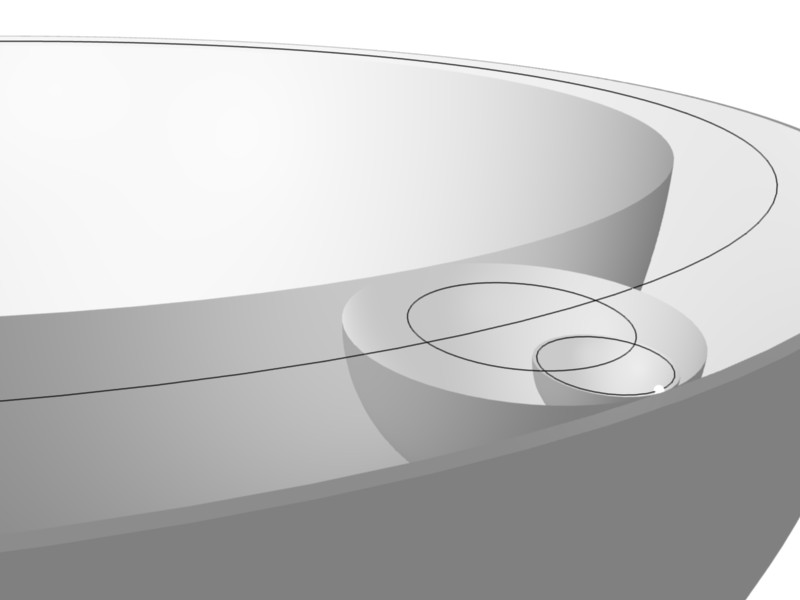}
Une représentation moderne des orbes du Soleil, 
librement inspirée d'Ibn al-\v{S}\=a\d{t}ir\label{soleil_blender}
\end{center}

\begin{center}
\includegraphics[width=\textwidth]{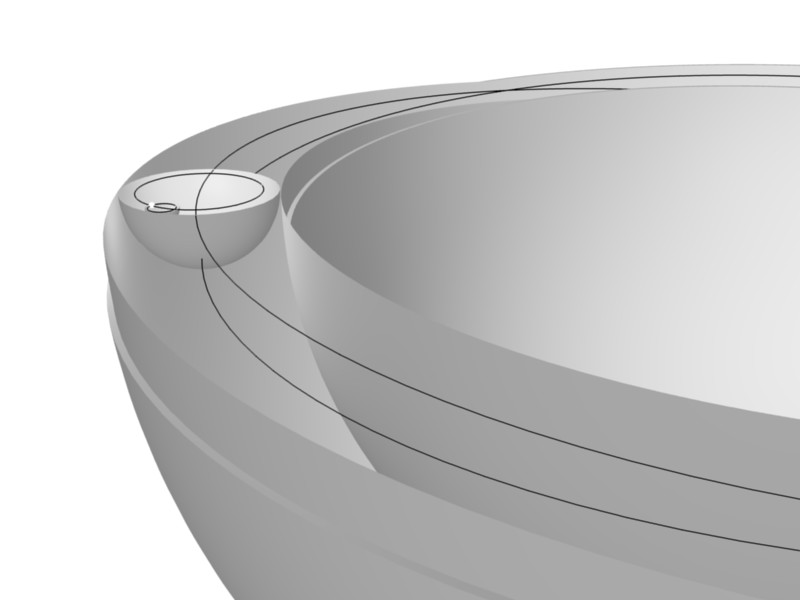}
Une représentation moderne des orbes de la Lune,
librement inspirée d'Ibn al-\v{S}\=a\v{t}ir
\end{center}

\newpage\phantomsection
\index{BJBHBE@\RL{ywm}!BJBHBE AHBDBJBDAJBG@\RL{al-yawm bilaylatihi j al-'ayyAm bilayAlIhA}, nychtémère (\textit{litt.} le jour avec sa nuit)}
\index{AWBDAY@\RL{.tl`}!BEAWAGBDAY ASAJBHAGAEBJBJAI@\RL{m.tAl` istiwA'iyyaT}, coascension calculée à l'équateur terrestre}
\index{AYAOBD@\RL{`dl}!AJAYAOBJBD ACBJBJAGBE@\RL{ta`dIl al-'ayyAm}, équation des nychtémères}
\index{AWBDAY@\RL{.tl`}!BEAWAGBDAY@\RL{m.tAl`}, coascension}
\index{AHBGAJ@\RL{bht}, vitesse apparente}
\addcontentsline{toc}{chapter}{I.8 Calcul de l'équation des nychtémères}
\includepdf[pages=37,pagecommand={\thispagestyle{plain}}]{edit.pdf}\phantomsection

\begin{center}
  \Large Chapitre huit

  \normalsize Calcul de l'équation des nychtémères
\end{center}
[Cette équation] est l'écart entre les nychtémères, rapporté
au jour [où] la quantité observée est [à la fois] l'astre moyen est
l'astre corrigé. Nous avons déjà posé que le jour où nous avons
observé à la fois l'astre moyen et l'astre corrigé est l'équinoxe de
printemps, et nous lui rapportons tous les autres nychtémères. La
méthode pour calculer cette équation est la suivante. Tu calcules le
Soleil moyen et le Soleil corrigé au jour souhaité. Ajoute au Soleil
moyen deux degrés, une minute et sept secondes. Prends la différence
entre cela et la coascension du Soleil corrigé (coascension calculée
par rapport à l'équateur terrestre en prenant comme origine le
commencement du Bélier). Multiplie-la par quatre. On obtient les
minutes de l'équation des nychtémères.  Si le Soleil moyen avec la
quantité ajoutée est supérieur à la coascension du Soleil corrigé,
alors retranche l'équation de la date, et s'il lui est inférieur,
alors ajoute-la à la date. On obtient la date corrigée par l'équation
des nychtémères. On en tire les [astres] moyens dans la table.

\emph{Remarque.} Multiplie la vitesse apparente de l'astre (par heure)
par les minutes de l'équation des nychtémères. Regarde ce qu'on
obtient~: si la position moyenne avec la quantité ajoutée est
supérieure à la coascension du Soleil corrigé, alors retranche le
résultat de cet astre corrigé, mais si elle leur est inférieure, alors
ajoute le résultat à cet astre corrigé. On obtient l'astre corrigé à
l'instant demandé. Nous avons déjà calculé une table pour l'y trouver
en fonction du Soleil moyen.

\newpage\phantomsection 
\index{AOBHAQ@\RL{dwr}!BEAOAGAQ@\RL{mdAr}, trajectoire, trajectoire diurne (\textit{i. e.} cercle parallèle à l'équateur)}
\index{ANASBA@\RL{_hsf}!ANASBHBAAI@\RL{_hsUf}, éclipse de Lune}
\index{BCASBA@\RL{ksf}!BCASBHBAAI@\RL{ksUf}, éclipse de Soleil}
\index{BABDBC@\RL{flk}!BABDBC BEAGAEBD@\RL{falak mA'il}, orbe incliné}
\index{ALBEAY@\RL{jm`}!AEALAJBEAGAYAI@\RL{'ijtimA`}, conjonction}
\index{BBAHBD@\RL{qbl}!AEASAJBBAHAGBDAI@\RL{'istiqbAl}, opposition}
\index{AHAWBDBEBJBHAS@\RL{b.talimayUs}, Ptolémée}
\index{AOBHAQ@\RL{dwr}!AJAOBHBJAQ@\RL{tadwIr}, épicycle}
\index{AEAHAQANAS@\RL{'ibr_hs}, Hipparque}
\index{AMAVAV@\RL{.h.d.d}!AMAVBJAV@\RL{.ha.dI.d}, périgée}
\index{APAQBH@\RL{_drw}!APAQBHAI@\RL{_dirwaT}, sommet, apogée}
\index{ANBDBA@\RL{_hlf}!ACANAJBDAGBA@\RL{i_htilAf}, irrégularité, anomalie, variation}
\addcontentsline{toc}{chapter}{I.9 Les orbes de la Lune et leurs mouvements}
\includepdf[pages=38,pagecommand={\thispagestyle{plain}}]{edit.pdf}\phantomsection

\begin{center}
  \Large Chapitre neuf

  \large Les orbes de la Lune et leurs mouvements selon la
méthode juste et sauve des doutes
\end{center}
Quand on a observé la Lune, on a découvert qu'elle se meut sur une
trajectoire différente de celle du Soleil et la coupant en deux lieux
opposés qui se déplacent le long de la trajectoire du Soleil en sens
inverse des signes. Cette trajectoire est inclinée par rapport à celle
du Soleil, et l'inclinaison maximale est toujours d'une même grandeur
de chaque côté (cinq parts), et telle que, dans une moitié de sa
trajectoire, la Lune est au Nord de la ceinture du zodiaque, et dans
l'autre moitié, au Sud. On a donc su qu'elle a un orbe
parécliptique, et un orbe incliné par rapport à celui-ci, et
que ces deux orbes se coupent en deux points opposés, donc que l'un
bissecte l'autre.

Quand on a découvert que le mouvement de la Lune sur l'orbe incliné
varie, en lenteur et en vitesse, autour du centre du monde qui est le
centre du zodiaque, on a su qu'elle a un orbe dont le centre est sur
la ceinture de l'orbe incliné. De sorte que, quand la Lune est dans la
partie haute de cet orbe, elle est plus éloignée du centre du monde,
et quand elle est dans sa partie basse, elle en est plus proche~; en
témoigne la variation de la grandeur du corps de la Lune, en
particulier lors des éclipses de Lune et des éclipses de Soleil. On
appelle cet orbe \emph{épicycle} depuis longtemps.

Ensuite, on a trouvé que la grandeur du rayon de l'épicycle, lors des
conjonctions et des oppositions, ne dépasse jamais cinq parts et un
sixième chez nous (cinq parts et un quart chez Ptolémée et
Hipparque)~; mais on a trouvé que sa grandeur, si la distance de la
Lune au Soleil est d'un quart de cercle des deux
côtés\footnote{C'est-à-dire lors des quadratures.}, est de huit
parts. Nous avons su ainsi qu'il y a un autre petit orbe sur la
ceinture de l'épicycle, de sorte que, quand la Lune est dans les
conjonctions ou les oppositions, elle est au périgée de cet orbe, et
que, quand elle est dans les quadratures, elle est à l'apogée -- la
distance au centre de l'épicycle est alors huit parts. Si la Lune est
au périgée de cette orbe, sa distance au centre de l'épicycle est cinq
parts et un sixième. Son mouvement est le double du mouvement de
l'orbe portant le centre de l'épicycle, de sorte que, quand l'orbe
portant l'épicycle tourne d'un quart de cercle, l'orbe portant le
corps de la Lune se meut d'un demi-cercle.

Ainsi savons-nous que la Lune a quatre orbes et quatre mouvements
simples.

\newpage\phantomsection  
\index{BABDBC@\RL{flk}!BABDBC BEAGAEBD@\RL{falak mA'il}, orbe incliné}
\index{BABDBC@\RL{flk}!BABDBC BFAGAQ@\RL{falak al-nAr}, orbe du feu}
\index{BABDBC@\RL{flk}!BABDBC AJAOBHBJAQ@\RL{falak al-tadwIr}, orbe de l'épicycle}
\index{BABDBC@\RL{flk}!BABDBCBEAOBJAQ@\RL{falak mdIr}, orbe rotateur}
\index{AYBBAO@\RL{`qd}!AYBBAOAI@\RL{`qdaT}, n{\oe}ud}
\index{ALARBGAQ@\RL{jzhr}!ALBHARBGAQ@\RL{jUzhr}|see{\RL{`qdaT}}}
\index{ALARBGAQ@\RL{jzhr}!ALBHARBGAQ BEALAGAR ATBEAGBDBJBJ@\RL{jwzhr mujAz ^simAliyy}|see{\RL{ra's}}}
\index{ALARBGAQ@\RL{jzhr}!ALBHARBGAQ BEALAGAR ALBFBHAHBJBJ@\RL{jwzhr mujAz junUbiyy}|see{\RL{_dnb}}}
\index{AQABAS@\RL{ra'asa}!AQABAS@\RL{ra's}, tête, n{\oe}ud ascendant}
\index{APBFAH@\RL{_dnb}!APBFAH@\RL{_dnb}, queue, n{\oe}ud descendant}
\index{AOBHAQ@\RL{dwr}!AJAOBHBJAQ@\RL{tadwIr}, épicycle}
\index{BABDBC@\RL{flk}!BABDBCBEBEAKAKBD@\RL{falak muma_t_tal}, parécliptique}
\includepdf[pages=39,pagecommand={\thispagestyle{plain}}]{edit.pdf}\phantomsection

Le premier est l'orbe \emph{parécliptique}. Ses deux pôles sont alignés
avec les pôles de l'écliptique. Sa partie convexe touche la partie
concave des orbes de Mercure, et sa partie concave touche la partie
convexe du deuxième de ses orbes. Son rayon est de soixante-neuf
parts.

Le deuxième orbe est un \emph{orbe incliné}. L'inclinaison de son
plan par rapport au plan du parécliptique est constante et son maximum
est de cinq parts. Son centre est le centre du monde, c'est le centre
du parécliptique~; sa partie convexe touche la partie concave du
parécliptique, et sa partie concave touche l'orbe du feu, comme on
sait. Son rayon est, par supposition, de soixante parts. La ceinture
qui est sur sa partie convexe coupe la ceinture qui est sur la partie
concave du parécliptique, en deux points opposés. On les appelle les
\emph{n{\oe}uds}. Le premier s'appelle \emph{tête}, ou n{\oe}ud
  ascendant, car la Lune le passe en s'inclinant vers le Nord par
rapport au parécliptique. L'autre, qui lui est opposé, s'appelle
\emph{queue}, ou n{\oe}ud descendant. Les pôles de l'orbe
incliné sont attachés à deux points de la partie concave du
parécliptique. Leur distance aux pôles du parécliptique mesure
l'inclinaison maximale de l'orbe incliné (c'est aussi la latitude
maximale de la Lune)~: cinq parts. Le centre de la partie concave de
l'orbe incliné est le centre du monde et son rayon est de
cinquante-et-une parts.

Quant au troisième orbe, nous imaginons une sphère de rayon huit
parts, seize minutes et vingt-sept secondes (des parts telles que le
rayon de l'orbe incliné en compte soixante). Nous la supposons plongée
dans l'orbe incliné qu'elle touche en un point de sa ceinture. Nous
l'appelons sphère de l'\emph{épicycle}.

Quant au quatrième orbe, nous imaginons une sphère de rayon une part,
quarante-et-une minutes et vingt-sept secondes (des mêmes parts). Nous
la supposons plongée dans le corps de l'épicycle~; elle le touche en
un point de sa ceinture qui est dans le plan de l'orbe incliné, et c'est
aussi en ce point que l'épicycle touche la ceinture de l'orbe incliné. On
appelle cet orbe le \emph{rotateur}. Nous supposons que le corps de
la Lune est plongé dans le corps du rotateur. Un point de la partie
convexe de la Lune touche la ceinture du rotateur qui est dans le plan
de l'orbe incliné. On a trouvé par nos observations que le diamètre du
globe lunaire est de trente-deux minutes cinquante-quatre secondes
(des mêmes parts).

\newpage\phantomsection  
\index{ALARBGAQ@\RL{jzhr}!AMAQBCAI ALBHARBGAQ@\RL{.harakaT al-jwzhr}, mouvement des n{\oe}uds}
\index{BHASAW@\RL{ws.t}!BHASAW@\RL{ws.t}, astre moyen}
\index{AMAQBC@\RL{.hrk}!AMAQBCAI BD-AYAQAV@\RL{.harakaT al-`r.d}, mouvement en latitude}
\index{AMAQBC@\RL{.hrk}!AMAQBCAI ANAGAUAI@\RL{.harakaT _hA.saT}, mouvement propre}
\index{APAQBH@\RL{_drw}!APAQBHAI BD-AJAOBHBJAQ BD-BEAQAEBJBJ@\RL{_dirwaT al-tadwIr al-mar'iyy}, apogée apparent}
\index{AHAYAO@\RL{b`d}!AHAYAO BEAG AHBJBF BD-BFBJAQBJBF AH-BD-BHASAW@\RL{bu`d mA bayna al-nIrayn bi-al-ws.t}, élongation moyenne}
\index{BBAWAQ@\RL{q.tr}!BFAUBA BBAWAQ BD-AJAOBHBJAQ BD-BEAQAEBJBJ@\RL{n.sf q.tr al-tadwIr al-mar'iyy}, rayon de l'épicycle apparent}
\includepdf[pages=40,pagecommand={\thispagestyle{plain}}]{edit.pdf}\phantomsection

Quant aux mouvements, le premier mouvement, mouvement du
parécliptique, est un mouvement simple autour de son centre (qui est
le centre du monde) et sur ses pôles, en sens inverse des signes, de
la grandeur du mouvement des n{\oe}uds, c'est-à-dire
$0;3,10,38,27$. Se déplacent de ce mouvement~: la tête, la queue, le
lieu où la latitude de la Lune est maximale, et tous les orbes de la
Lune. En vérité, à proprement parler, ce mouvement est l'excédent du
mouvement des n{\oe}uds sur le mouvement des étoiles fixes.\label{lab001}

Le deuxième mouvement est le mouvement de l'orbe incliné autour de son
centre (qui est le centre du parécliptique et le centre du monde) et
autour de ses pôles fixes que l'on a mentionnés précédemment. C'est un
mouvement simple et uniforme qui fait en des temps égaux des
[déplacements] égaux dans le sens des signes du zodiaque, de
$13;13,45,39,40$. Il vaut la somme de la Lune moyenne et du mouvements
des n{\oe}uds, et on l'appelle depuis longtemps \emph{mouvement en
  latitude}. Le centre de l'épicycle se déplace donc dans le sens des
signes, à partir d'un point qu'on imagine être un point fixe par rapport
à l'écliptique, d'un mouvement de grandeur égale à la Lune moyenne
c'est-à-dire $13;10,35,1,13,53$.

Le troisième mouvement est le mouvement de l'épicycle autour de son
centre et sur ses pôles situés à la perpendiculaire du plan de l'orbe
incliné. Ce mouvement est de $13;3,53,46,18$ en un jour et une nuit,
et il est en sens inverse des signes dans la partie supérieure [de
  l'épicycle]. On le connait depuis longtemps~: c'est le
\emph{mouvement propre} de la Lune. Ce mouvement commence à l'apogée
de l'épicycle apparent.

Le quatrième mouvement est le mouvement du rotateur. C'est un
mouvement simple autour de son centre, qui fait se mouvoir la Lune sur
la ceinture du rotateur dans le plan de la ceinture de l'orbe incliné,
dans le sens des signes quand elle est dans sa partie supérieure. Ce
mouvement est de $24;22,53,23$ en un jour et une nuit. C'est deux fois
l'excédent de la Lune moyenne sur le Soleil moyen. Ce mouvement
commence au périgée lors des conjonctions et des oppositions moyennes
et il passe par l'apogée lors des quadratures.

Pour montrer cela, supposons Soleil moyen et Lune moyenne en
conjonction en un point fixe de l'écliptique~; supposons que le centre
de l'épicycle et le centre du rotateur sont sur une droite passant par
le centre du monde et par ce point que l'on imagine être fixe par
rapport à l'écliptique~; supposons la Lune au périgée du rotateur, en
cet instant, sur cette droite, à distance minimale du centre de
l'épicycle (cinq parts et un sixième). 

\newpage\phantomsection  
\includepdf[pages=41,pagecommand={\thispagestyle{plain}}]{edit.pdf}\phantomsection

\noindent Puis chaque orbe se meut du
mouvement simple observé que nous avons assigné à chacun~:

-- Le Soleil se meut, en un jour et une nuit, de $0;59,8,20$ dans
  le sens des signes.

-- Le parécliptique se meut, en un jour et une nuit, de
  $0;3,10,38,27$ dans le sens inverses des signes~; il entraîne l'orbe
  incliné, l'épicycle et le rotateur dans ce mouvement en sens inverse
  des signes.

-- L'orbe incliné se meut dans le sens des signes, dans le même
  temps, de $13;13,45,39,40$. Laissant le point fixe de l'écliptique
  où les deux luminaires étaient en conjonction, l'élongation du
  centre de l'épicycle devient, en un jour et une nuit,
  $12;11,26,41,30$~; c'est l'excédent du mouvement de la Lune moyenne
  sur le Soleil moyen, et on l'appelle \emph{élongation moyenne entre
    les deux luminaires}.

-- Le rotateur se meut autour de son centre (donc la Lune se
  déplace du périgée du rotateur), en sens inverse des signes, du
  double de l'élongation, c'est-à-dire $24;22,53,23$ en un jour et une
  nuit.

Alors le corps de la Lune s'éloigne du centre de l'épicycle~; sa
distance s'appelle \emph{rayon de l'épicycle apparent}.

Si l'orbe incliné s'est mû (avec son inclinaison) d'un quart de cercle
en sept jours, dix-huit minutes, quatorze secondes, cinquante-deux
tierces et quinze quartes d'un jour et une nuit, alors l'orbe rotateur
s'est mû, dans le même temps, d'un demi-cercle (car son mouvement est
le double de l'élongation). La Lune se déplace donc du périgée de
l'épicycle vers son apogée, et la distance du centre du corps de la
Lune au centre de l'épicycle devient huit parts~: c'est le rayon
de l'épicycle apparent dans les première et seconde quadratures.

Si l'orbe incliné s'est mû, pendant une durée double de celle-ci, d'un
demi-cercle, alors le centre de l'épicycle est passé par l'opposé du
point fixe (sur la droite joignant ce point fixe au centre du monde),
du côté opposé au lieu de la conjonction. Il y a donc à présent
opposition entre le Soleil moyen et la Lune moyenne. Le rotateur se
meut d'un mouvement double de celui de l'orbe incliné, c'est-à-dire
d'un cercle entier, et il déplace la Lune vers le périgée du
rotateur.

\newpage\phantomsection  
\index{ANBDBA@\RL{_hlf}!ACANAJBDAGBA@\RL{i_htilAf}, irrégularité, anomalie, variation}
\index{AMAQBC@\RL{.hrk}!AMAQBCAI AHASBJAWAI@\RL{.harakaT basI.taT}, mouvement simple}
\index{AMAQBC@\RL{.hrk}!AMAQBCAI BEANAJBDBAAI@\RL{.harakaT mu_htalifaT}, mouvement irrégulier (\textit{i. e.} non uniforme)}
\index{AYAOBD@\RL{`dl}!AJAYAOBJBD@\RL{ta`dIl}, équation}
\label{var22}
\includepdf[pages=42,pagecommand={\thispagestyle{plain}}]{edit.pdf}\phantomsection

\noindent La distance de la Lune au centre de l'épicycle devient
minimale, et le rayon de l'épicycle apparent devient cinq parts et un
sixième ($5;10$), comme pendant la conjonction. Donc le rayon de
l'épicycle apparent est toujours $5;10$ dans les conjonctions et les
oppositions, et huit parts dans les quadratures. On l'a toujours
trouvé ainsi dans l'observation.

De plus, l'épicycle se meut de son mouvement propre. Le centre de la
Lune tourne alors autour du centre de l'épicycle en une trajectoire
dont le rayon est $5;10$, comme nous l'avons dit, dans les
conjonctions et les oppositions. L'arcsinus de ceci est
$4;56,24$. C'est l'équation maximale dans les conjonctions et les
oppositions, atteinte quand la Lune est sur la droite tangente à
l'épicycle apparent, d'un côté ou de l'autre. Mais quand le centre de
l'épicycle est dans les quadratures, la distance du centre du corps de
la Lune au centre de l'épicycle est de huit parts~: c'est le rayon
de l'épicycle apparent dans les quadratures. Son arcsinus est sept
parts et deux tiers. C'est ce que nous avons trouvé par l'observation
dans les quadratures quand la Lune est sur la droite tangente à
l'épicycle apparent d'un côté ou de l'autre.

Comme le rapport entre le mouvement propre de la Lune et le mouvement
de la Lune moyenne est inférieur au rapport entre la droite joignant
le centre du monde au périgée de l'épicycle apparent et le rayon de
l'épicycle apparent, alors la Lune n'a ni station ni rétrogradation
(mais son mouvement devient lent dans la moitié de l'apogée et rapide
dans la moitié du périgée)~: on ne la voit pas revenir en sens inverse
des signes à cause de la petitesse de l'épicycle et de la faiblesse de
son mouvement par rapport au mouvement de l'orbe incliné.

Avec le rayon de l'épicycle apparent variant entre cinq parts un
sixième et huit parts, c'est comme si les degrés de lenteur et de
vitesse n'étaient pas constants, mais variables, de sorte que la
lenteur, la vitesse, \emph{etc.} reviennent, tantôt moins, tantôt
plus, à cause de cette variation.

De ces orbes aux mouvements simples procèdent les mouvements
irréguliers -- mouvements en accord avec l'observation en tout temps.

Quant aux variations simples en longitude qui s'imposent en raison de
ces mouvements, il y en a quatre~; on les appelle \emph{équations}.

\newpage\phantomsection  
\label{var7}
\includepdf[pages=43,pagecommand={\thispagestyle{plain}}]{edit.pdf}\phantomsection

La première variation\footnote{notée $c_2$ dans notre commentaire, car
  elle sera désignée comme ``deuxième variation'' page \pageref{deuxieme}
  ci-dessous.} est ce qu'on obtient à cause du rayon de
l'épicycle dans les conjonctions et les oppositions. C'est l'angle
formé au centre du monde par deux droites qui en sont issues, l'une
passant par le centre de l'épicycle et l'autre par le centre du corps
de la Lune. Cet angle est maximal quand la droite passant par la Lune
touche l'épicycle apparent d'un côté ou de l'autre (sans le
couper). Ce maximum, dans les conjonctions et les oppositions, est de
$4;56,24$ (à condition que le rayon de l'orbe incliné soit de soixante
parts). Cet angle s'évanouit à l'apogée et au périgée~: c'est quand
les deux droites citées sont confondues. Cette variation est
soustraite à la Lune moyenne tant que la Lune descend de l'apogée au
périgée en sens inverse des signes, dans la moitié droite, jusqu'à ce
que la droite passant par la Lune coïncide avec la droite passant par
le centre de l'épicycle. Au contraire, on l'ajoute à la Lune moyenne
tant que la Lune monte du périgée vers l'apogée dans le sens des
signes, jusqu'à ce que l'une des droites coïncide avec l'autre. Cette
variation est appelée équation simple ou \emph{première équation}.

La deuxième variation est causée par l'augmentation de cette première
variation quand le centre de l'épicycle n'est ni en conjonction ni en
opposition. Elle est mêlée à la première et n'existe pas sans
elle. Elle est maximale quand le centre de l'épicycle est dans les
quadratures, c'est-à-dire à un quart de cercle du lieu où se situe la
conjonction moyenne ou l'opposition moyenne. Son maximum est $2;44$,
soit $7;40$ avec la première équation. Elle est toujours ajoutée à la
première équation, que celle-ci soit ajoutée ou soustraite [à la Lune
  moyenne]. On l'appelle \emph{deuxième équation} (et on l'appelait
autrefois variation de la distance minimale).

\emph{Remarque.} Nous avons déjà avancé que le rayon de l'épicycle apparent
est toujours $5;10$ dans les conjonctions et les oppositions, $8;0$
dans les quadratures. Entre ces positions, il augmente graduellement
du plus petit au plus grand. La première variation devient, à cause de
cette variation du rayon, la somme des deux variations. C'est subtil~:
prête attention.

La troisième variation\footnote{notée $c_1$ dans notre commentaire.}
est appelée \emph{équation de la Lune
  propre}.\label{c1} C'est l'angle formé au centre de l'épicycle par deux
droites qui en sont issues, l'une passant par le centre du rotateur et
l'autre par le centre de la Lune.

\newpage\phantomsection  
\index{BFBBBD@\RL{nql}!AJAYAOBJBD BFBBBD@\RL{ta`dIl al-naql}, équation du déplacement}
\index{BFASAH@\RL{nsb}!AOBBAGAEBB BFASAH@\RL{daqA'iq al-nisb}, coefficient d'interpolation}
\index{AYAQAV@\RL{`r.d}!AYAQAV@\RL{`r.d}, latitude, \emph{i. e.} par rapport à l'écliptique}
\index{ANBDBA@\RL{_hlf}!ACANAJBDAGBA@\RL{i_htilAf}, irrégularité, anomalie, variation}
\includepdf[pages=44,pagecommand={\thispagestyle{plain}}]{edit.pdf}\phantomsection

\noindent Si celle qui passe par le centre de
la Lune touche le rotateur, alors l'équation est maximale, et c'est,
selon nos observations, $12;26$. Cette équation s'évanouit lorsque la
Lune est à l'apogée ou au périgée par rapport au centre de
l'épicycle. On l'ajoute à la Lune propre quand l'élongation du centre
de l'épicycle au Soleil moyen est inférieure à un quart de cercle ou
supérieure à trois quarts de cercle, c'est-à-dire supérieure à neuf
signes ou inférieure à trois signes. On la soustrait quand
l'élongation est supérieure à un quart de cercle et inférieure à trois
quarts de cercle.

La quatrième variation est appelée \emph{équation du
  déplacement}. Sans cette équation, on obtient la Lune vraie par rapport
de l'orbe incliné~; mais le but est d'avoir sa position par rapport au
parécliptique pour connaître sa position par rapport à la ceinture du
zodiaque. Cette équation est l'écart entre la distance de son lieu au
n{\oe}ud sur la ceinture du parécliptique, et cette même distance sur
la ceinture de l'orbe incliné. J'en tiens compte quand je veux
convertir l'une de ces deux distances en l'autre. On appelle cela
<<~déplacer la Lune de l'orbe incliné à l'écliptique~>>. Cette
équation est maximale dans les octants de l'orbe incliné, et son
maximum est environ six minutes. Elle s'évanouit aux n{\oe}uds et
quand la Lune atteint sa latitude maximale de part et d'autre. Les
modernes en ont peu fait usage.

Le \emph{coefficient d'interpolation} en minutes est un nombre dont le
rapport à soixante est la proportion due à la variation de
l'élongation par rapport au Soleil. C'est une donnée qui facilite le
calcul de la Lune vraie par les tables.\footnote{Ce coefficient est noté $\chi$ dans notre commentaire. Ceci est mieux expliqué au chapitre 11.}

\emph{Remarque.} Si nous faisions le calcul selon la méthode de calcul
de la première équation en partant du rayon de l'épicycle apparent,
alors on obtiendrait la somme due des première et deuxième équations
correspondant à l'élongation pour laquelle on aurait fait le
calcul. C'est subtil~: prête attention.

Parlons de la latitude de la Lune. J'ai dit précédemment que ses
dernières latitudes\footnote{<<~ses dernières latitudes~>>,
  c'est-à-dire les latitudes (Nord et Sud) maximales.} sont égales au
Nord et au Sud. D'après ce qu'on a trouvé en observant, la latitude
maximale est de cinq parts. La Lune est au Nord quand elle va de la
tête vers la queue, et au Sud quand elle va de la queue vers la
tête. La Lune monte de sa dernière latitude Sud à sa dernière latitude
Nord, elle descend de sa dernière latitude Nord à sa dernière latitude
Sud.

\newpage\phantomsection  
\index{BFAXAQ@\RL{n.zr}!ANAJBDAGBA BEBFAXAQ@\RL{i_htilAf al-man.zr}, parallaxe}
\index{ATBCBD@\RL{^skl}!AJATBCBCBDAGAJ@\RL{t^skkulAt}, phases (de la Lune)}
\index{BEAMBH@\RL{m.hw}!BEAMBH@\RL{m.hw}, obscurcissement}
\index{BFBHAQ@\RL{nwr}!BFBHAQ@\RL{nwr}, lumière}
\index{ANBDBA@\RL{_hlf}!ACANAJBDAGBA@\RL{i_htilAf}, irrégularité, anomalie, variation}
\index{AHAWBDBEBJBHAS@\RL{b.talimayUs}, Ptolémée}
\includepdf[pages=45,pagecommand={\thispagestyle{plain}}]{edit.pdf}\phantomsection

\noindent Le dessein de l'ascension est de rapprocher la Lune du pôle
visible, et celui de la descente est de l'en éloigner. Sa parallaxe et
ses phases selon sa position par rapport au Soleil seront traitées
ultérieurement dans un chapitre particulier, si Dieu le veut.

Quant à la variation des parts de sa surface qui reçoivent la lumière et
qu'on appelle <<~obscurcissement~>>\footnote{les tâches dues au relief à la
  surface de la Lune, comme on le sait depuis Galilée.}, le plus vraisemblable
est que certaines sont placées face à la lumière et d'autres pas~; ou
bien c'est à cause de la variation de sa position par rapport au
Soleil. On dit aussi qu'il y a, dans l'épicycle de la Lune, avec elle,
des corps distincts qui ne sont pas disposés pareillement face à
l'éclairage. On dit aussi que la raison de cette variation est la
réflexion des rayons vers la Lune par la mer environnante, à cause de
son poli, et l'absence de leur réflexion par la surface du quart
habité et du reste du globe, à cause de sa rugosité. En effet les
lieux de sa face éclairés par les rayons venant à elle du Soleil après
réflexion par la surface de la mer seraient plus lumineux que les
lieux qui sont éclairés par les rayons venant seulement du Soleil.

\begin{center}
  \large Section
\end{center}
Quelle est la distance de la Lune au centre de la Terre~? Dans cette
configuration, lors des conjonctions et des oppositions, cette
distance vaut soixante parts en moyenne, $65;10$ au maximum et $54;50$
au minimum. Dans les quadratures, elle vaut $68;0$ au maximum et
$52;0$ au minimum. La distance maximale atteinte à l'intérieur de ses
orbes est $69;0$, et la distance minimale est $51;0$. Le tout, en
parts telles que le rayon de l'orbe incliné en compte soixante.

Sache que la distance du centre de l'épicycle au centre de la Terre
est soixante parts et qu'elle n'augmente ni ne diminue~; mais si tu
veux la distance de la Lune en parts telles que le rayon de la Terre
en compte une seule, ou bien telle autre quantité en usage pour les
distances et les volumes, alors voici la route à suivre selon mes
observations. Tu ôtes de chaque degré deux minutes et il restera ce
qu'on demandait~; en suivant la voie de Ptolémée, tu jetterais une
minute de chaque degré. En effet la distance du centre de l'épicycle
au centre du monde est selon nos observations cinquante-huit fois le
rayon de la Terre, et selon Ptolémée cinquante-neuf fois le rayon de
la Terre.

\newpage\phantomsection
\index{ALARBGAQ@\RL{jzhr}!BHASAW ALBHARBGAQ@\RL{ws.t al-jwzhr}, n{\oe}ud moyen}
\index{ALARBGAQ@\RL{jzhr}!AJBBBHBJBE ALBHARBGAQ@\RL{tqwIm al-jwzhr}, n{\oe}ud vrai}
\index{BHASAW@\RL{ws.t}!BHASAW BBBEAQ@\RL{ws.t al-qamar}, Lune moyenne}
\index{AHAYAO@\RL{b`d}!AHAYAO BEAVAGAYBA@\RL{b`d m.dA`f}, élongation double}
\index{ANAUAU@\RL{_h.s.s}!ANAGAUAUAI BBBEAQ@\RL{_hA.s.saT al-qamar}, Lune propre}
\index{APAQBH@\RL{_drw}!APAQBHAI AJAOBHBJAQ BHASAWBI@\RL{_dirwaT al-tadwIr al-ws.t_A}, apogée moyen de l'épicycle}
\index{APAQBH@\RL{_drw}!APAQBHAI BEAQAEBJBJAI@\RL{_dirwaT mar'iyyaT}, apogée apparent}
\index{BBBHBE@\RL{qwm}!AJBBBHBJBE BBBEAQ@\RL{tqwIm al-qamar}, Lune vraie}
\index{AYAQAV@\RL{`r.d}!AMAUAUAI AYAQAV@\RL{.hi.s.saT `r.d}, argument de latitude}
\includepdf[pages=46,pagecommand={\thispagestyle{plain}}]{edit.pdf}\phantomsection

Comment interpréter les choses qui se rapportent à la Lune~?

-- Le \emph{n{\oe}ud moyen} est [l'arc] entre le commencement du
  Bélier et le n{\oe}ud ascendant, dans l'orbe parécliptique, en sens
  inverse des signes.\label{lab002}

-- Le \emph{n{\oe}ud vrai} est [l'arc] entre ces deux points dans
  l'orbe parécliptique, dans le sens des signes.

-- La \emph{Lune moyenne} est [l'arc] entre le point face au
  commencement du Bélier (à condition qu'il soit fixe) et l'extrémité
  de la droite joignant le centre du monde au centre de l'épicycle. Et
  si j'avais dit, l'élongation du centre de l'épicycle au point
  mentionné dans le sens des signes, cela aurait été juste aussi.

-- Si l'on ôte le Soleil moyen de la Lune moyenne il reste
  l'\emph{élongation} du centre de l'épicycle au Soleil moyen, et le
  double de cette élongation s'appelle l'\emph{élongation
    double}. C'est égal au mouvement du rotateur de la Lune.

-- La \emph{Lune propre} est [l'arc] entre l'apogée moyen de l'épicycle
  et le centre du corps de la Lune, sur la ceinture de
  l'épicycle, dans le sens qu'on a supposé [pour l'épicycle]
  (c'est-à-dire, si on est dans la moitié haute, en sens inverse des
  signes).

La grandeur de ces arcs ne varie pas. Parmi les choses dont la
grandeur varie en croissant ou en décroissant, il y a~:

-- La \emph{Lune propre corrigée} est [l'arc] entre le centre du
  corps de la Lune et son apogée apparent sur la ceinture de son
  épicycle. C'est la véritable Lune propre, après qu'on l'a corrigée
  par la première équation.

-- La \emph{Lune vraie} est [l'arc] entre le commencement du
  Bélier et le point du parécliptique où tombe son cercle de latitude,
  dans le sens des signes.

-- Son \emph{argument de latitude} est l'excédent de la Lune vraie
  sur le n{\oe}ud vrai.

\newpage\phantomsection  
\index{ASBFBH@\RL{snw}!ASBFAI BBBEAQBJAI@\RL{sanaT qamariyaT}, année lunaire}
\index{ASBFBH@\RL{snw}!ASBFAI ATBEAS@\RL{sanaT al-^sams}, année solaire}
\index{ATBGAQ@\RL{^shr}!ATBGAQ BBBEAQBJBJ@\RL{^shr qamariyy}, mois lunaire}
\index{ALAOBD@\RL{jdl}!ALAOBHBD@\RL{jdwl}, table}
\addcontentsline{toc}{chapter}{I.10 Détermination des mouvements de la Lune à une date donnée}
\includepdf[pages=47,pagecommand={\thispagestyle{plain}}]{edit.pdf}\phantomsection

\begin{center}
  \Large Chapitre dix

  \large Détermination des mouvements de la Lune à une date donnée
\end{center}
Dans les tables, nous avons écrit que la Lune moyenne à midi du
premier jour de l'année sept cent un de l'ère de Yazdgard est sept
signes et $3;35,50$ degrés. Le mouvement de la Lune moyenne en vingt
années persanes est cinq signes et $13;7,40,49,48,21$\footnote{Le
  <<~cinq signes et $13;7,40,49,48,21$ degrés~>>, c'est notre correction
  (modulo 360 comme les valeurs suivantes)~; le texte arabe (y compris
  l'édition page ci-contre) donne <<~deux signes et $7;40,49,48,21$~>>.},
en une seule année quatre
signes et $9;23,2,29,25,1$ degrés, en trente jours un signe et
$5;17,30,36,56,18,9$ degrés, en un jour et une nuit
$13;10,35,1,13,52,36,18,16,7$ degrés, et en une heure égale
$0;32,56,27,33,41,31$ degrés.

Nous avons aussi écrit que la Lune propre à la date mentionnée est
quatre signes et $18;32,27$ degrés. Son mouvement en vingt années
persanes est onze signes et $4;22,30$ degrés, en une seule année deux
signes et $28;43,7$ degrés, en un mois un signe et $1;56,58$ degrés,
en un jour $13;3,53,56$ degrés, et en une heure $0;32,39,45$.

Nous avons établi que le n{\oe}ud moyen à la date mentionnée est
neuf signes et $5;7,35$ degrés. Son mouvement en sens inverse des
signes est toutes les vingt années persanes
de zéro signe et $26;34,37,20$ degrés, en une seule année
$19;19,43,52$ degrés, en un mois $1;35,19,13,18,54,15$, en un jour
$0;3,10,38,26,37,48,29,36$, et en une heure $0;0,7,56,36,6,35$.

Ceci étant admis, sache que l'année lunaire\footnote{Une année lunaire
  compte exactement douze mois lunaires.} est exactement
$354;22,1,35,8,50,24$, le mois lunaire moyen\footnote{Avec les notations
  adoptées dans notre commentaire mathématique, le mois lunaire moyen est $360°/(\omega_m-\omega_m^{\astrosun})$. On trouve bien la valeur donnée ici, jusqu'au cinquième rang sexagésimal, à partir des paramètres $\omega_m$ et $\omega_m^{\astrosun}$ adoptés par {\shatir}.} est
$29;31,50,7,55,44,12$, l'année solaire exacte est trois cent
soixante-cinq jours, cinq heures, quarante-neuf minutes et trois
cinquièmes de secondes, c'est-à-dire $365;5,49,0,36,32$.\footnote{Le
  texte arabe (y compris l'édition ci-contre) dit <<~trois cinquièmes de
  soixantièmes de secondes~>>, ce
  qui semble incohérent avec la valeur écrite en chiffres. Pourtant,
  en calculant l'année solaire avec un mouvement moyen par année
  persane de $359;45,40$, on trouve 365 jours et $5;49,0,33,36$
  heures. Le 36 exprime donc bien <<~trois cinquièmes de soixantièmes
  de secondes~>>, et l'erreur commise par le copiste serait alors
  d'avoir oublié trentre-trois secondes. Mais la durée de l'année en
  jours donnée ensuite confirme le choix que nous avons fait (sinon,
  elle serait en effet $365;14,32,31,24$).}

Supposons que le jour et la nuit comptent soixante parts, alors
l'année solaire est
$365;14,32,31,31,20$. L'excédent\footnote{Nous n'avons pas compris cette
  phrase. Est-ce une variation ?} par rapport au cercle est de $6;15,9,8$.

\newpage\phantomsection  
\index{AHBGAJ@\RL{bht}, vitesse apparente}
\index{AHAWBDBEBJBHAS@\RL{b.talimayUs}, Ptolémée}
\index{AEAHAQANAS@\RL{'ibr_hs}, Hipparque}
\index{BGBFAO@\RL{hnd}!BGBFAO@\RL{al-hnd}, les Indiens}
\index{ALBEBGAQ@\RL{jmhr}!ALBEBGBHAQ@\RL{al-jumhUr}, les Grecs}
\index{ACAKBJAQ BD-AOBJBF BD-AHBGAQBJ@\RL{'a_tIr al-dIn al-abharI}, Ath{\=\i}r al-D{\=\i}n al-Abhar{\=\i}}
\index{BCBHATBJAGAQ AHBF BDAHAHAGBF@\RL{kU^syAr bn labbAn}, K\=u\v{s}y\=ar b. Labb\=an}
\index{AWBHASBJ@\RL{n.sIr al-dIn al-.tUsI}, Na\d{s}{\=\i}r al-D{\=\i}n al-\d{T}\=us{\=\i}}
\index{AXBDBD@\RL{.zll}!BBAWAQ AXBDBD@\RL{qi.tr al-.zill}, diamètre de l'ombre}
\index{ANASBA@\RL{_hsf}!ANASBHBA@\RL{_husUf}, éclipse de Lune}
\label{var18}
\includepdf[pages=48,pagecommand={\thispagestyle{plain}}]{edit.pdf}\phantomsection

Ceci étant admis, sache que le diamètre de la Lune quand elle est à
distance moyenne de la Terre, c'est-à-dire quand la distance
Terre-Lune est de soixante parts (en parts du rayon de l'orbe portant
l'épicycle), est $0;32;54;33$. J'avais déjà découvert cela en suivant
une autre voie et j'avais trouvé $32;22$ mais je préfère la première
valeur. Le diamètre de la Lune est selon Hipparque $33;15$, selon
Ptolémée $33;5$, selon les Indiens $32;0$, selon les Grecs $32;20$, et
selon Ath{\=\i}r\footnote{Ath{\=\i}r al-D{\=\i}n al-Abhar{\=\i}},
K\=u\v{s}y\=ar\footnote{K\=u\v{s}y\=ar b. Labb\=an} et
al-\d{T}\=us{\=\i} $32;3$.\label{obs_diam_lune}

Si nous voulons le diamètre de la Lune à d'autres distances que la
distance moyenne, nous divisons son diamètre à distance moyenne par sa
distance au centre du monde (en parts du diamètre du parécliptique)~:
on obtient le diamètre de la Lune quand elle est à cette
distance.\footnote{C'est-à-dire que $0;32;54;33\times
  60=\text{diamètre apparent}\times\text{distance Terre-Lune}$, où le
  diamètre apparent est exprimé en degrés, et la distance Terre-Lune
  en soixantièmes du rayon de la trajectoire du centre de l'épicycle
  dans l'orbe incliné.}

La vitesse apparente de la Lune par heure, c'est aussi son diamètre, mais sache
que la première méthode est plus exacte.

Quant au diamètre de l'ombre, selon Ptolémée, il est deux fois et
trois cinquièmes de fois comme le diamètre de la Lune. Je l'ai trouvé,
par l'observation de nombreuses éclipses anciennes et récentes, deux
fois et demi et un cinquième de fois comme le diamètre de la Lune si
nous supposons que la latitude maximale de la Lune est cinq degrés~;
et si nous supposons qu'elle est de quatre degrés, deux tiers et un
quart de degrés, alors le diamètre de l'ombre est deux fois et deux
tiers de fois comme le diamètre de la Lune. Par l'observation, j'ai
trouvé que la durée des éclipses est supérieure à ce que prédit le
calcul sous l'hypothèse que le diamètre de l'ombre serait deux fois et
trois cinquièmes de fois comme le diamètre de la
Lune.\footnote{Al-\v{S}\=a\d{t}ir donne donc trois valeurs du diamètre
  de l'ombre : $2;36$, $2;42$ et $2;40$. La première, attribuée à Ptolémé, est rejetée sur la base d'observations. }

\newpage\phantomsection  
\includepdf[pages=49,pagecommand={\thispagestyle{plain}}]{edit.pdf}\phantomsection

\begin{center}
\begingroup%
  \makeatletter%
  \providecommand\color[2][]{%
    \errmessage{(Inkscape) Color is used for the text in Inkscape, but the package 'color.sty' is not loaded}%
    \renewcommand\color[2][]{}%
  }%
  \providecommand\transparent[1]{%
    \errmessage{(Inkscape) Transparency is used (non-zero) for the text in Inkscape, but the package 'transparent.sty' is not loaded}%
    \renewcommand\transparent[1]{}%
  }%
  \providecommand\rotatebox[2]{#2}%
  \ifx\svgwidth\undefined%
    \setlength{\unitlength}{480bp}%
    \ifx\svgscale\undefined%
      \relax%
    \else%
      \setlength{\unitlength}{\unitlength * \real{\svgscale}}%
    \fi%
  \else%
    \setlength{\unitlength}{\svgwidth}%
  \fi%
  \global\let\svgwidth\undefined%
  \global\let\svgscale\undefined%
  \makeatother%
  \begin{picture}(1,1)%
    \put(0,0){\includegraphics[width=\unitlength]{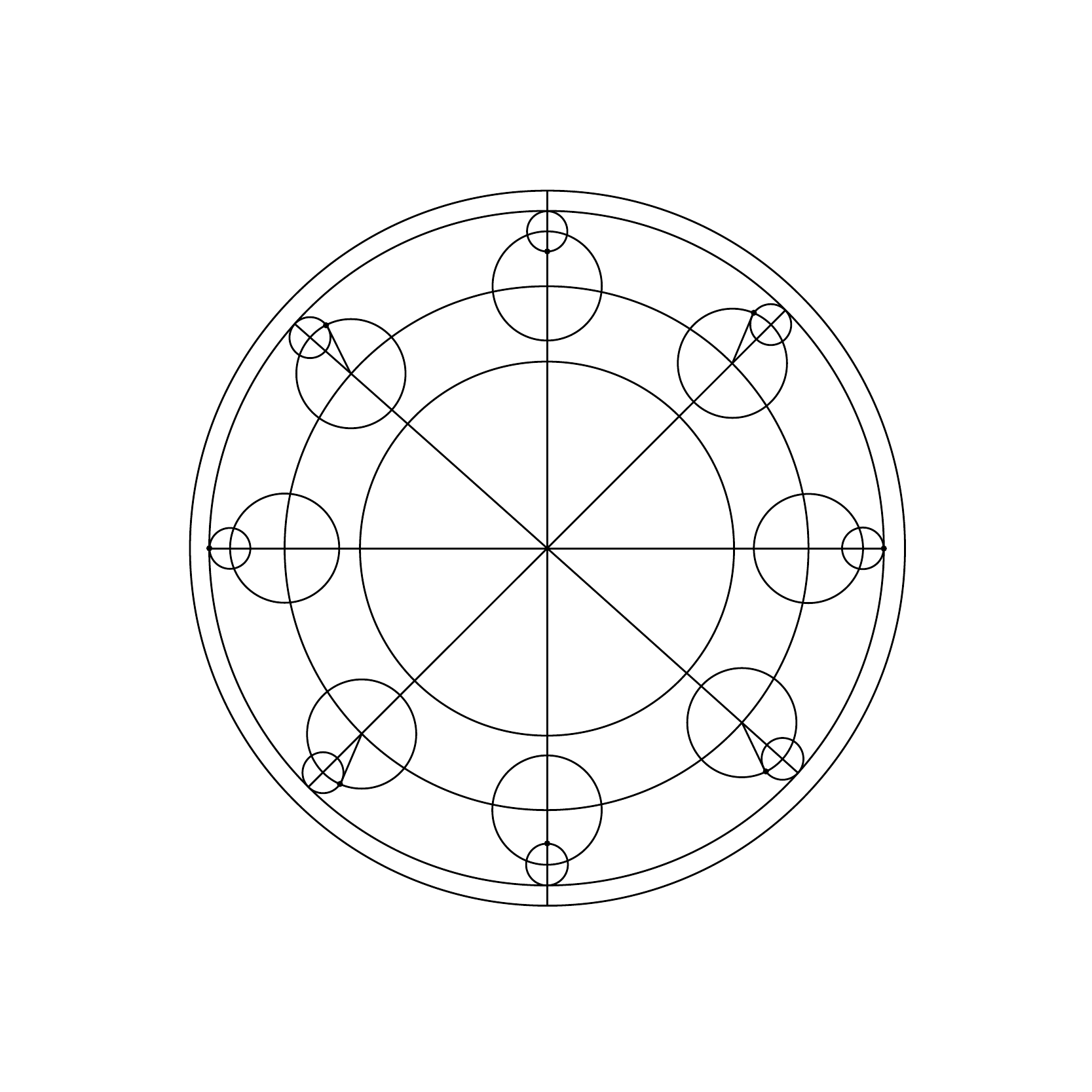}}%
    \put(0.41454449,0.50281993){\color[rgb]{0,0,0}\makebox(0,0)[lb]{\smash{centre}}}%
    \put(0.3963209,0.64856558){\color[rgb]{0,0,0}\rotatebox{15.00000022}{\makebox(0,0)[lb]{\smash{partie concave}}}}%
    \put(0.44880952,0.83392858){\color[rgb]{0,0,0}\makebox(0,0)[lb]{\smash{conjonction}}}%
    \put(0.5101033,0.82600407){\color[rgb]{0,0,0}\rotatebox{-59.99999989}{\makebox(0,0)[lb]{\smash{orbe rotateur}}}}%
    \put(0.37800503,0.69509599){\color[rgb]{0,0,0}\rotatebox{45}{\makebox(0,0)[lb]{\smash{orbe de l'épicycle}}}}%
    \put(0.58234099,0.73199126){\color[rgb]{0,0,0}\rotatebox{-30.00000011}{\makebox(0,0)[lb]{\smash{ceinture de l'orbe incliné}}}}%
    \put(0.51547619,0.50416667){\color[rgb]{0,0,0}\makebox(0,0)[lb]{\smash{du Monde}}}%
    \put(0.43259463,0.95602404){\color[rgb]{0,0,0}\rotatebox{-30.00000011}{\makebox(0,0)[lb]{\smash{parécliptique (dans un plan incliné}}}}%
    \put(0.74164621,0.72483371){\color[rgb]{0,0,0}\rotatebox{45}{\makebox(0,0)[lb]{\smash{un signe et demi}}}}%
    \put(0.83333335,0.49761907){\color[rgb]{0,0,0}\makebox(0,0)[lb]{\smash{quadrature}}}%
    \put(0.03452381,0.4940476){\color[rgb]{0,0,0}\makebox(0,0)[lb]{\smash{quadrature}}}%
    \put(0.49761903,0.15238093){\color[rgb]{0,0,0}\makebox(0,0)[lb]{\smash{opposition}}}%
    \put(0.12778012,0.83658445){\color[rgb]{0,0,0}\rotatebox{-45}{\makebox(0,0)[lb]{\smash{un signe et demi}}}}%
    \put(0.7491682,0.27417333){\color[rgb]{0,0,0}\rotatebox{-45}{\makebox(0,0)[lb]{\smash{sept signes et demi}}}}%
    \put(0.11634805,0.10086512){\color[rgb]{0,0,0}\rotatebox{45}{\makebox(0,0)[lb]{\smash{quatre signes et demi}}}}%
    \put(0.31511055,0.76029672){\color[rgb]{0,0,0}\rotatebox{30.00000011}{\makebox(0,0)[lb]{\smash{partie convexe}}}}%
    \put(0.30489828,0.73549801){\color[rgb]{0,0,0}\rotatebox{-74.99999978}{\makebox(0,0)[lb]{\smash{équation du mouvement propre}}}}%
    \put(0.32388492,0.74271219){\color[rgb]{0,0,0}\rotatebox{30.00000011}{\makebox(0,0)[lb]{\smash{de l'orbe incliné}}}}%
    \put(0.40394305,0.62890082){\color[rgb]{0,0,0}\rotatebox{15.00000022}{\makebox(0,0)[lb]{\smash{de l'orbe incliné}}}}%
    \put(0.47642107,0.90608954){\color[rgb]{0,0,0}\rotatebox{-30.00000011}{\makebox(0,0)[lb]{\smash{par rapport à l'orbe incliné)}}}}%
  \end{picture}%
\endgroup%

Les orbes de la Lune figurés par les trajectoires des centres des
sphères. Nous y avons représenté l'orbe de l'épicycle et le rotateur
en huit positions~: la conjonction, l'opposition, les quadratures et
les octants.
\end{center}

\newpage\phantomsection   
\includepdf[pages=50,pagecommand={\thispagestyle{plain}}]{edit.pdf}\phantomsection

\begin{center}
\begingroup%
  \makeatletter%
  \providecommand\color[2][]{%
    \errmessage{(Inkscape) Color is used for the text in Inkscape, but the package 'color.sty' is not loaded}%
    \renewcommand\color[2][]{}%
  }%
  \providecommand\transparent[1]{%
    \errmessage{(Inkscape) Transparency is used (non-zero) for the text in Inkscape, but the package 'transparent.sty' is not loaded}%
    \renewcommand\transparent[1]{}%
  }%
  \providecommand\rotatebox[2]{#2}%
  \ifx\svgwidth\undefined%
    \setlength{\unitlength}{480bp}%
    \ifx\svgscale\undefined%
      \relax%
    \else%
      \setlength{\unitlength}{\unitlength * \real{\svgscale}}%
    \fi%
  \else%
    \setlength{\unitlength}{\svgwidth}%
  \fi%
  \global\let\svgwidth\undefined%
  \global\let\svgscale\undefined%
  \makeatother%
  \begin{picture}(1,1)%
    \put(0,0){\includegraphics[width=\unitlength]{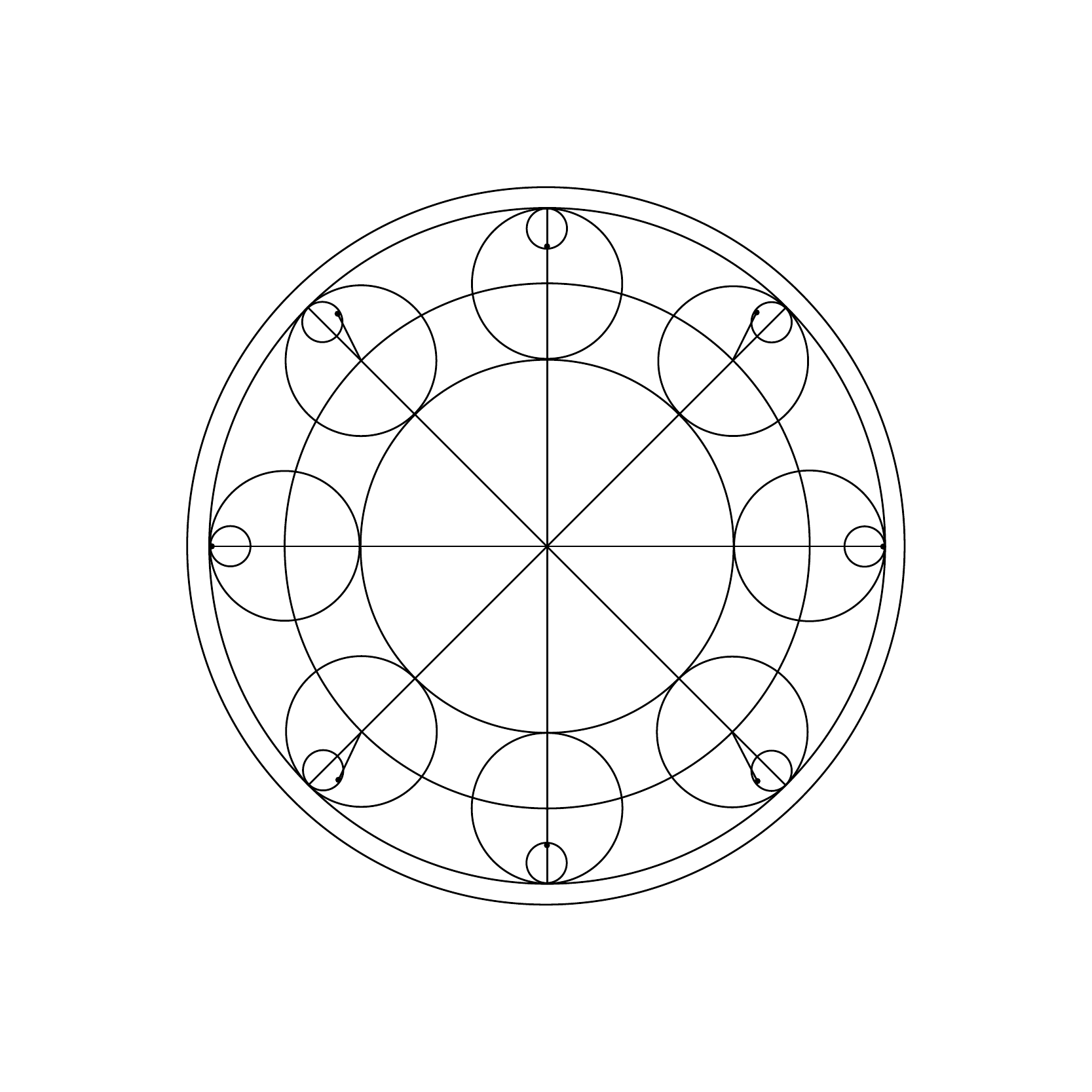}}%
    \put(0.44761906,0.83452382){\color[rgb]{0,0,0}\makebox(0,0)[lb]{\smash{conjonctions}}}%
    \put(0.83095229,0.49642855){\color[rgb]{0,0,0}\makebox(0,0)[lb]{\smash{quadratures}}}%
    \put(0.47633488,0.15040855){\color[rgb]{0,0,0}\makebox(0,0)[lb]{\smash{oppositions}}}%
    \put(0.04178632,0.49727033){\color[rgb]{0,0,0}\makebox(0,0)[lb]{\smash{quadratures}}}%
    \put(0.51309523,0.48690478){\color[rgb]{0,0,0}\makebox(0,0)[lb]{\smash{du Monde}}}%
    \put(0.43690474,0.48690478){\color[rgb]{0,0,0}\makebox(0,0)[lb]{\smash{centre}}}%
    \put(0.43558871,0.7209831){\color[rgb]{0,0,0}\rotatebox{-30.00000011}{\makebox(0,0)[lb]{\smash{orbe de l'épicycle}}}}%
    \put(0.458436,0.78080642){\color[rgb]{0,0,0}\rotatebox{-30.00000011}{\makebox(0,0)[lb]{\smash{rotateur}}}}%
    \put(0.54990743,0.81233856){\color[rgb]{0,0,0}\rotatebox{-15.00000022}{\makebox(0,0)[lb]{\smash{partie convexe de l'orbe incliné}}}}%
    \put(0.1871077,0.80565496){\color[rgb]{0,0,0}\rotatebox{-45}{\makebox(0,0)[lb]{\smash{octants}}}}%
    \put(0.36474536,0.58088626){\color[rgb]{0,0,0}\rotatebox{30.00000011}{\makebox(0,0)[lb]{\smash{partie concave}}}}%
    \put(0.27207276,0.74675741){\color[rgb]{0,0,0}\rotatebox{45}{\makebox(0,0)[lb]{\smash{dans un plan incliné par rapport}}}}%
    \put(0.12340479,0.59700192){\color[rgb]{0,0,0}\rotatebox{45}{\makebox(0,0)[lb]{\smash{l'orbe parécliptique}}}}%
    \put(0.38173447,0.56027133){\color[rgb]{0,0,0}\rotatebox{30.00000011}{\makebox(0,0)[lb]{\smash{de l'orbe incliné}}}}%
    \put(0.40014308,0.84057228){\color[rgb]{0,0,0}\rotatebox{45}{\makebox(0,0)[lb]{\smash{à l'orbe incliné}}}}%
  \end{picture}%
\endgroup%

Les orbes solides de la Lune mus de mouvements simples uniformes en
leurs centres. Le plan de la figure est le plan de l'orbe incliné~: on
a vu que le plan de l'orbe incliné est incliné par rapport au plan de
l'écliptique ou du parécliptique d'une inclinaison égale à la latitude
de la Lune.
\end{center}

\newpage\phantomsection  
\index{AYAOBD@\RL{`dl}!AJAYAOBJBD BEAMBCBE@\RL{ta`dIl mu.hakam}, équation principale}
\addcontentsline{toc}{chapter}{I.11 La Lune vraie par les tables ou par le calcul}
\label{var8}\label{var12}
\includepdf[pages=51,pagecommand={\thispagestyle{plain}}]{edit.pdf}\phantomsection

\begin{center}
  \Large Chapitre onze

  \large La Lune vraie par les tables ou par le calcul
\end{center}
Si tu veux [la Lune vraie], alors calcule à partir des tables la Lune
moyenne, la Lune propre et le Soleil moyen de la manière que je t'ai
indiquée pour le Soleil. Puis soustrais le Soleil moyen de la Lune
moyenne~; le double du résultat est le centre de la Lune, c'est-à-dire
l'élongation du centre de l'épicycle au point où la distance est
maximale\footnote{Le <<~centre de la Lune~>> est donc ce qu'il a
  appelé précédemment (chapitre 9) l'<<~élongation double~>>,
  c'est-à-dire la longueur de l'arc du rotateur entre le centre de la
  Lune et le point du rotateur où la distance (au centre de
  l'épicycle) est maximale.}. Entre, avec cela, dans la table des
équations de la Lune, et prends-y l'équation de la Lune propre. Si le
centre est plus petit que six signes\footnote{Six signes
  $=180^\circ$.}, alors ajoute l'équation à la Lune propre, et s'il
est plus grand que six signes, alors soustrais-la de la Lune
propre. On obtient la Lune propre corrigée\footnote{C'est donc la
  <<~troisième équation~>> du chapitre 9, notée $c_1$ dans
  notre commentaire.}. Puis entre avec la Lune
propre corrigée dans la table des équations de
la Lune, et prends-y la troisième\footnote{C'est $c_2(\alpha,0)$ dans notre
  commentaire, \textit{i. e.} la <<~première équation~>> du chapitre 9.}
et la quatrième équation\footnote{C'est $(c_2(\alpha,180°)-c_2(\alpha,0))$
  dans notre commentaire, \textit{i. e.}  la <<~deuxième équation~>>
  du chapitre 9.}. Puis entre, avec le
centre, dans la table des équations de la Lune, prends-y la deuxième
équation\footnote{C'est $\chi$ dans notre commentaire, appelé
  <<~coefficient d'interpolation~>> au chapitre 9.},
multiplie-la par la quatrième, et ajoute le
tout à la troisième~; on obtient la troisième équation 
\emph{principale}\footnote{Nous traduisons \textit{mu\d{h}km} par
  \emph{principale}. Pour la Lune comme pour les autres astres,
  l'équation principale est celle calculée au moyen des
  transformations planes, après rabattement du plan de l'orbe incliné
  dans le plan de l'écliptique. Elle néglige donc l'effet des
  inclinaisons des orbes. Mais peut-être ce nom renvoie-t-il aussi au
  fait que cette équation est toujours calculée par interpolation --
  alors aurait-on pu aussi traduire par \emph{équation
    interpolée}.}. Puis regarde si la Lune propre corrigée est
inférieure à six signes. Dans ce cas, retranche la principale de la
Lune moyenne~; mais si elle est plus grande que six signes, alors
ajoutes-la à la Lune moyenne. On obtient la Lune vraie par rapport à
la ceinture de l'orbe incliné. Puis retranche le n{\oe}ud vrai de la
Lune vraie. Il reste l'argument de latitude de la Lune. Entre-le dans
la table de l'équation du déplacement de la Lune\footnote{C'est la
  <<~quatrième équation~>> du chapitre 9.}. Retranche ce que tu y
trouves de la Lune vraie si l'argument de latitude de la Lune est
compris entre un et quatre-vingt-dix degrés ou bien entre cent
quatre-vingt et deux cent soixante-dix degrés, et ajoute-le sinon.

\newpage\phantomsection  
\label{var19}
\includepdf[pages=52,pagecommand={\thispagestyle{plain}}]{edit.pdf}\phantomsection

\noindent On obtient la Lune vraie par rapport à l'écliptique à la
date pour laquelle tu as fait le calcul, à midi, à Damas~; et si la
date n'a pas encore été corrigée par la correction des jours avec
leurs nuits, alors il faut la corriger comme je l'ai enseigné pour le
Soleil\footnote{Voir chapitre 8.}.

Si tu veux la latitude de la Lune, entre l'argument de latitude dans
la table des latitudes de la Lune, tu y trouves la latitude de la Lune
et son sens\footnote{son <<~sens~>>, c'est-à-dire si elle est orientée
  vers le Nord ou bien vers le Sud.}.

La manière d'y arriver par le calcul est que tu multiplies le sinus de
l'argument de latitude par la tangente de la latitude
entière\footnote{L'\emph{argument de latitude} est ici l'argument de
  latitude \emph{par rapport à l'écliptique} (de même que {\shatir} a
  distingué précédemment la Lune vraie par rapport à la ceinture de
  l'orbe incliné et la Lune vraie par rapport à l'écliptique). La
  \emph{latitude entière} est l'inclinaison de l'orbe incliné par
  rapport à l'écliptique, $5^\circ$. Voir commentaire mathématique.}
de la Lune, on obtient la tangente de la latitude partielle de la
Lune. Tu trouves son arc dans la table des tangentes inverses~: c'est
la latitude. On a déjà parlé de son sens.

\begin{center}
  \large Section

  \normalsize Détermination de l'équation de la Lune propre par le
  calcul
\end{center}
Détermine le centre de la Lune et multiplie son sinus par le rayon de
l'orbe portant le corps de la Lune, c'est-à-dire par $1;25$. Appelons
ce qu'on obtient \emph{le dividende}. Multiplie ensuite le cosinus du
centre par $1;25$. \^Ote le résultat du rayon de l'épicycle ($6;35$)
si le centre est compris entre un et quatre-vingt-dix degrés ou bien
entre deux cent soixante-dix et trois cent soixante degrés. Sinon,
ajoute-le au rayon de l'épicycle. Ajoute le carré de ce qu'on obtient
au carré du dividende. Prends-en la racine carrée, c'est le rayon
de l'épicycle apparent. Garde-le, tu en auras besoin ici et
ailleurs. Ensuite, divise le dividende par le rayon de l'épicycle
apparent, il en sort le sinus de l'équation de la Lune propre~;
écris-la en face de ce degré dans la table.

\newpage\phantomsection  
\includepdf[pages=53,pagecommand={\thispagestyle{plain}}]{edit.pdf}\phantomsection

\begin{center}
  \large Section

  \normalsize Calcul de la deuxième équation\footnote{\label{deuxieme}
    C'est $c_2$ dans notre commentaire. \`A nouveau, Ibn
    al-\v{S}\=a\d{t}ir a changé de numérotation. Il désigne ici par
    <<~deuxième équation~>> celle qu'il appelait <<~troisième
    équation~>> au début de ce chapitre, et <<~première équation~>> au
    chapitre neuf.} de la Lune.
\end{center}
Multiplie le sinus de la Lune propre corrigée et son cosinus par le
rayon de l'épicycle apparent lors d'une conjonction ou d'une
opposition, c'est-à-dire par l'excédent du rayon de l'épicycle sur le
rayon du rotateur, ce qui fait $5;10$. Ce qu'on obtient du cosinus,
ajoute-le toujours à soixante si la Lune propre corrigée est dans la
moitié supérieure de l'épicycle, et ôte-le si elle est dans sa moitié
inférieure. Ce qu'on obtient, ajoute son carré au carré du produit de
sinus du mouvement propre par $5;10$ (ce dernier produit est ce qu'on
appelle à présent le \emph{dividende}). Prends la racine carré de
cette somme. Divise le dividende par cela. On obtient le sinus de la
deuxième équation correspondant au degré pour lequel tu as fait le
calcul. Prends son arc dans la table des sinus~: tu trouves l'équation
correspondant à cette position.

La manière de calculer la troisième équation est la suivante. Tu
supposes que le rayon de l'épicycle apparent est de huit degrés,
c'est-à-dire le sinus de sept degrés et deux tiers, et tu suis ce que
je t'ai enseigné pour la deuxième équation. Ce qu'on obtient pour
chaque position, soustrais-en la deuxième équation pour cette même
position. \'Ecris le résultat en face de cette position.

\emph{Remarque.} Si tu utilises le rayon de l'épicycle apparent à
une position donnée au lieu de $5;10$, et que tu suis le procédé
indiqué, il en sort l'équation principale, c'est-à-dire la
deuxième équation corrigée.

\newpage\phantomsection  
\includepdf[pages=54,pagecommand={\thispagestyle{plain}}]{edit.pdf}\phantomsection

\begin{center}
  \large Section

  \normalsize Calcul du coefficient d'interpolation pour la Lune
\end{center}
Si tu veux cela, prends l'arc du rayon de l'épicycle apparent dans la
table des sinus~; ce qu'on obtient est la variation maximale pour ce
degré comme équation principale. Retranches-en $4;56,24$
qui est le maximum de la deuxième équation. Ce qui reste, divise-le
par $2;43,21$ (c'est la différence entre l'équation maximale dans les
quadratures $7;39,45,11$ et le maximum de la deuxième équation). Ce
qui sort de la division, c'est le coefficient d'interpolation~;
écris-le en face de ce degré.

\emph{Remarque.} Sache qu'il est possible de construire d'autres
tables pour la Lune vraie, autrement disposées que celles
connues. Pour cela, tu supposes que l'épicycle de la Lune est
l'épicycle véritable dont le rayon est $6;35$. Par rapport à ce
nombre, tu calcules une différence qui lui est retranchée à distance
maximale, et une différence qui lui est ajoutée à distance
minimale\footnote{<<~à distance minimale~>>, il veut dire ici quand le
  rayon de l'épicycle apparent est maximal, c'est-à-dire que la Lune
  peut être au plus proche de la Terre.} quand il vaut huit parts.
Par rapport à ce nombre-ci, tu calcules une autre différence
qui en est retranchée avec un coefficient d'interpolation (à l'opposé
de ce qu'on a vu ci-dessus\footnote{Veut-il dire simplement que
  l'interpolation se fait ici par soustraction plutôt que par
  addition~?}). Enfin, on peut aussi faire le calcul de la Lune vraie
au moyen d'une table unique\footnote{Dans ce qui suit, Ibn
  al-\v{S}\=a\d{t}ir semble expliquer comment construire, pour chaque
  valeur du centre, une <<~table unique~>> contenant l'équation de la
  Lune propre et toutes les valeurs de l'équation principale
  correspondant aux différentes valeurs de la Lune propre. Il conçoit
  donc une sorte de tableau à double entrée pour représenter
  l'équation principale comme fonction de deux variables~: le centre
  et la Lune propre. Chaque colonne de ce tableau est elle-même une
  table, et il explique dans le paragraphe suivant comment construire
  la colonne correspondant à un centre de $90^\circ$.} qui montrerait
à la fois la variation maximale à une position donnée dans l'orbe
portant [l'épicycle]\footnote{L'orbe portant l'épicycle, c'est-à-dire
  l'orbe incliné.} et la correction du mouvement propre dépendant de
cette position.

\newpage\phantomsection 
\label{var20}
\includepdf[pages=55,pagecommand={\thispagestyle{plain}}]{edit.pdf}\phantomsection

Par exemple, si nous supposons que le centre de la Lune (c'est-à-dire
l'élongation double) est de trois signes, alors la correction du
mouvement propre est $12;26$. D'autre part, si le centre est de trois
signes, la variation maximale est de six parts et un sixième de part~;
alors nous avons calculé la deuxième équation sous l'hypothèse que son
maximum est de six parts et un sixième\footnote{Veut-il dire que, pour
  un centre de $90^\circ$, il ne calcule de manière exacte qu'une
  seule valeur de la deuxième équation (son maximum, dénoté $\max\vert
  e(\cdot,90^\circ)\vert$ dans notre commentaire mathématique), et
  calcule toutes les autres par interpolation en utilisant cette
  unique valeur~? En effet, il serait bien fastidieux de calculer de
  manière exacte chaque cellule du tableau à double entrée~; mais
  l'explication reste sujette à interprétation.}. Nous avons écrit
ceci en face de $12;26$. Quand nous entrons dans cette table le
mouvement propre absolu (je veux dire, sans sa correction), nous
trouvons ainsi l'équation principale sans équation additionnelle~;
mais nous aurons alors besoin de nombreuses tables\footnote{Peut-être
  veut-il dire qu'il y aura de nombreuses colonnes dans ce tableau à
  double entrée : en effet, une colonne pour chaque valeur de
  l'élongation double.}.

\begin{center}
  \large Section
\end{center}
\emph{Si tu veux la distance de la Lune au centre du monde,} prends
connaissance du sinus de la Lune propre corrigée et de son cosinus~;
multiplie chacun par le rayon de l'épicycle apparent ($5;10$ lors
d'une conjonction ou d'une opposition, sinon, c'est ce que je t'ai
demandé de conserver dans un paragraphe antérieur). Ce qu'on obtient à
partir du sinus de la Lune propre, c'est le \emph{dividende} et
garde-le. Ce qu'on obtient à partir du cosinus, ajoute-le à soixante
si la Lune propre corrigée est inférieure à trois signes ou supérieure
à neuf signes, sinon ôte-le de soixante. Le total ou le reste, prends
son carré et ajoute-le au carré du dividende. Prends la racine carrée
de cela, c'est la distance de la Lune au centre du monde à condition
que le rayon de l'orbe portant [l'épicycle]\footnote{C'est-à-dire le rayon
  de l'orbe incliné.} soit soixante.

\emph{Remarque.} Si nous divisons par cette distance le diamètre de la
Lune à distance moyenne (c'est $0;32,54,33$), alors on obtient son
diamètre à cette distance. Si nous multiplions la vitesse apparente de
la Lune par jour par deux minutes et demi\footnote{Deux minutes et
  demi font un vingt-quatrième de degré.}, alors on obtient aussi le
diamètre de la Lune, mais la première méthode est plus précise. Dans
les conjonctions et les oppositions, le diamètre minimum de la Lune
est $0;30,18$, et son diamètre maximum $0;36,0$. Dans les quadratures
son diamètre minimum est $0;29,2,15$ et son diamètre maximum
$0;37,58,20$. Sache cela.

\emph{Seconde remarque.} Si nous divisions par cette distance le
dividende qu'on a gardé, il en sortirait le sinus de l'équation
principale de la Lune. Comprends cela.

\newpage\phantomsection 
\index{BABDBC@\RL{flk}!BABDBCBEBEAKAKBD@\RL{falak muma_t_tal}, parécliptique}
\index{BABDBC@\RL{flk}!BABDBC BEAGAEBD@\RL{falak mA'il}, orbe incliné}
\index{BABDBC@\RL{flk}!BABDBCAMAGBEBD@\RL{falak .hAmil}, orbe déférent}
\index{BABDBC@\RL{flk}!BABDBCBEAOBJAQ@\RL{falak mdIr}, orbe rotateur}
\index{BABDBC@\RL{flk}!BABDBC AJAOBHBJAQ@\RL{falak al-tadwIr}, orbe de l'épicycle}
\index{ACBHAL@\RL{'awj}!AMAQBCAIACBHAL@\RL{.harakaT al-'awj}, mouvement des apogées}
\index{AMAQBC@\RL{.hrk}!AMAQBCAI BEAQBCAR@\RL{.harakaT al-markaz}, mouvement du centre (en général, le centre d'un épicycle)}
\index{AMAQBC@\RL{.hrk}!AMAQBCAI ANAGAUAI@\RL{.harakaT _hA.saT}, mouvement propre}
\index{ANAUAU@\RL{_h.s.s}!ANAGAUAUAI BCBHBCAH@\RL{_hA.s.saT al-kwkab}, astre propre (traduit par ``anomalie'' par de nombreux auteurs)}
\index{ARAMBD@\RL{z.hl}!ARAMBD@\RL{zu.hal}, Saturne}
\addcontentsline{toc}{chapter}{I.12 Configuration des orbes de Saturne selon la vraie méthode}
\includepdf[pages=56,pagecommand={\thispagestyle{plain}}]{edit.pdf}\phantomsection

\begin{center}
  \Large Chapitre douze

  \large Configuration des orbes de Saturne selon la vraie méthode
\end{center}
Parmi les orbes de Saturne nous imaginons~:

-- un \emph{orbe parécliptique} représentant l'orbe de l'écliptique, dans
son plan, autour de son centre et sur ses pôles

-- un deuxième orbe, \emph{incliné} par rapport au parécliptique d'une
inclinaison constante de deux parts et demi et le coupant en deux
points opposés dont l'un s'appelle la tête et l'autre la queue

-- un troisième orbe dont le centre est sur le bord de l'orbe
  incliné et dont le rayon fait cinq parts et un huitième (en parts
  telles que le rayon de l'orbe incliné en compte soixante)~; il
  s'appelle \emph{déférent}
  
-- un quatrième orbe dont le centre est sur le bord du déférent et
  dont le rayon fait un degré, quarante-deux minutes et trente
  secondes~; il s'appelle \emph{rotateur}

-- un cinquième orbe dont le centre est sur le bord du rotateur et
  dont le rayon fait six parts et demi (en les mêmes parts)~; il
  s'appelle \emph{orbe de l'épicycle}

Le centre du corps de Saturne est attaché en un point de la ceinture
de l'épicycle.

Quant aux mouvements, le parécliptique se meut d'un mouvement simple
autour de son centre, dans le sens des signes. C'est le
\emph{mouvement des Apogées} qui fait chaque jour $0;0,9,52$. Les
intersections -- la tête et la queue -- se déplacent donc, et les deux
inclinaisons extrêmes aussi.

L'orbe incliné se meut d'un mouvement simple autour de son centre qui
est centre de l'univers, dans le sens des signes comme le
précédent. Sa mesure est le \emph{mouvement du centre de Saturne} et
c'est l'excédent du mouvement moyen sur le mouvement de l'Apogée~: il
fait, en un jour et une nuit, $0;2,0,26,17$.

Le déférent se meut en sens contraire dans sa partie supérieure, de la
même quantité que le mouvement du centre de Saturne, c'est-à-dire, en
un jour et une nuit, $0;2,0,26,17$.

Le rotateur se meut dans le sens des signes dans sa partie supérieure,
d'une quantité double de celle du mouvement du centre de Saturne,
c'est-à-dire, en un jour et une nuit, $0;4,0,52,34$. Une révolution
du rotateur s'achève donc avec une demi-révolution du déférent.

L'orbe de l'épicycle se meut d'un mouvement simple autour de son
centre dans le sens des signes dans sa partie supérieure. Sa mesure
est l'excédent du \emph{mouvement propre} de Saturne sur le mouvement de son
centre.

\newpage\phantomsection 
\index{AMAQBC@\RL{.hrk}!AMAQBCAI AHASBJAWAI BEAQBCBCAHAI@\RL{.harakaT basI.taT mrkkbaT}, mouvement simple-composé}
\index{APAQBH@\RL{_drw}!APAQBHAI@\RL{_dirwaT}, sommet, apogée}
\includepdf[pages=57,pagecommand={\thispagestyle{plain}}]{edit.pdf}\phantomsection

\noindent Ce mouvement est simple-composé. En effet, supposons que
l'épicycle est immobile\footnote{Immobile, au sein du référentiel
  constitué par l'orbe rotateur.}. L'orbe incliné se meut d'un quart
de cercle, l'orbe déférent d'un quart de cercle, et l'orbe rotateur
d'un demi cercle~: cela signifie que le diamètre passant par l'apogée
et le périgée de l'épicycle s'est déplacé d'un quart de cercle [en
  trop]. Nous devons donc supposer qu'il est animé d'un mouvement en
sens contraire égal au mouvement du déférent, c'est-à-dire au
mouvement du centre de Saturne. Ainsi le diamètre passant par l'apogée
et le périgée se déplacera vers son lieu.

Nous avons trouvé par l'observation que le mouvement propre dont se
meut l'épicycle fait, en un jour et une nuit, $0;57,7,43,34,22$, mais
on vient de voir que cet orbe se meut en sens contraire aux signes,
d'un mouvement simple autour de son centre. Or il semble se mouvoir
dans le sens des signes, d'un mouvement simple autour de son centre,
donc le plus petit des deux mouvements est soustrait du plus grand,
c'est-à-dire que nous soustrayons du mouvement propre le mouvement du
centre. Il reste le mouvement de l'orbe de l'épicycle dans le sens des
signes, égal à l'excédent du mouvement propre de Saturne sur le
mouvement de son centre~: par jour, $0;55,7,17,17,22$. C'est son vrai
mouvement simple, sauf que son apogée, dans le sens des signes, se
meut toujours d'autant que le mouvement propre de Saturne, somme de ce
mouvement simple et du mouvement de son centre. Il s'agit d'une notion
subtile, mais il faut la comprendre. Il en est de même pour les autres
planètes, donc on n'aura pas à répéter cette explication.\footnote{Le
  mot <<~apogée~>>, \textit{dhirwa}, désigne l'extrémité d'un diamètre
  de l'épicycle parallèle à la direction
  Terre--centre~de~l'épicycle. La grandeur appelée <<~mouvement de
  l'astre propre~>> désigne le mouvement de rotation de l'épicycle sur
  lui-même \emph{relativement à ce diamètre}~; mais le mouvement de
  l'épicycle \emph{relativement à l'orbe rotateur}, excédent du mouvement
  propre sur le mouvement du centre, intéresse {\shatir} davantage.}.
Quand nous voudrons l'élongation entre l'astre
et l'apogée, nous devrons sommer le mouvement de cet épicycle et le
mouvement du centre de Saturne~: le résultat sera comme le mouvement
propre de Saturne, ce sera l'élongation entre l'astre et l'apogée. Je
t'ai expliqué cette notion, elle contient une subtilité, examine-la
donc avec attention.

Supposons que les centres du déférent, du rotateur et de l'épicycle
soient alignés sur une droite passant par le centre du Monde et par
l'Apogée. Alors le centre de l'épicycle est à distance maximale de la
Terre, c'est-à-dire $63;25$. Maintenant, si l'orbe incliné se meut
d'un quart de cercle, que le déférent se meut, dans le même temps,
d'un quart de cercle aussi, et que le rotateur se meut, dans le même
temps, d'un demi-cercle, alors la distance du centre de l'épicycle au
centre du déférent devient six parts, une demi-part et un tiers de
part.

\newpage\phantomsection 
\index{AYAOBD@\RL{`dl}!AJAYAOBJBD ANAGAUAUAI@\RL{ta`dIl al-_hA.s.saT}, équation de l'astre propre}
\index{AYAOBD@\RL{`dl}!AJAYAOBJBD BEAQBCAR@\RL{ta`dIl al-markaz}, équation du centre}
\index{AMAQBC@\RL{.hrk}!AMAQBCAI ANAGAUAI@\RL{.harakaT _hA.saT}, mouvement propre}
\index{AYAOBD@\RL{`dl}!AJAYAOBJBD@\RL{ta`dIl}, équation}
\index{ANBDBA@\RL{_hlf}!ACANAJBDAGBA@\RL{i_htilAf}, irrégularité, anomalie, variation}
\index{APAQBH@\RL{_drw}!APAQBHAI BEAQAEBJBJAI@\RL{_dirwaT mar'iyyaT}, apogée apparent}
\index{APAQBH@\RL{_drw}!APAQBHAI AMBBBJBBBJBJAI@\RL{_dirwaT .haqIqiyyaT}, apogée vrai}
\includepdf[pages=58,pagecommand={\thispagestyle{plain}}]{edit.pdf}\phantomsection

Traçons une droite du centre du Monde au centre de l'épicycle, et une
autre droite du centre du Monde au centre du déférent. L'angle entre
ces deux droites est l'angle de la \emph{première équation de
  Saturne}. L'\emph{équation de l'astre propre} est de cette même
grandeur.

On explique ceci sur un exemple. Si l'orbe incliné se meut d'un
demi-cercle, que le déférent se meut d'un demi-cercle aussi, et que le
rotateur se meut, dans le même temps, d'un cercle entier, alors le
centre de l'épicycle vient à distance minimale du centre du Monde,
c'est-à-dire $56;35$, et la droite allant du centre du Monde au centre
de l'épicycle se confond avec celle allant du centre du Monde au
centre du déférent~; ainsi la première équation (qui est aussi bien 
équation du centre qu'équation de l'astre propre) s'évanouit.

Comme on a trouvé que le Soleil rejoint Saturne au milieu de son
mouvement direct et s'oppose à Saturne au milieu de sa rétrogradation,
on a su que le mouvement de l'orbe de l'épicycle vaut l'excédent du
mouvement de Soleil moyen sur Saturne moyen. Dès lors, nous avons su
qu'en soustrayant Saturne moyen du Soleil moyen, reste le mouvement
propre de Saturne. De même pour Jupiter et Mars.

Ceci étant admis, sache qu'à ces orbes et ces mouvements sont liées
des anomalies qu'on appelle \emph{équations}. La première équation
s'appelle \emph{équation du centre et de l'astre propre}. C'est
l'angle au centre du Monde (le centre de l'orbe incliné) entre deux
droites qui en sont issues, l'une passant par le centre du déférent,
et l'autre par le centre de l'épicycle. Cet angle est égal à l'angle
au centre de l'épicycle entre deux droites qui en sont issues, l'une
passant par le centre du Monde, et l'autre, parallèle à la droite
joignant le centre du déférent au centre du Monde. L'écart entre les
extrémités de ces deux droites le long de l'épicycle, c'est l'écart
entre l'apogée vrai et l'apogée apparent\footnote{\textit{Cf.} figure
  \ref{fig040} p.~\pageref{fig040} pour mieux saisir la terminologie
  adoptée par {\shatir} pour les planètes supérieures et en
  particulier pour Saturne. L'\emph{apogée apparent} et le
  \emph{périgée apparent} sont les points de l'épicycle alignés avec
  le centre du Monde et le centre de l'épicycle.}. Il est démontré que
l'équation de l'astre propre est égale à l'équation du centre à
cause de l'égalité des deux angles. On le voit clairement sur
l'exemple. Cette équation s'évanouit à l'Apogée et au périgée. Elle
s'ajoute au centre et se retranche de l'astre propre quand le centre
est supérieur à six signes~; elle se retranche du centre et s'ajoute
à l'astre propre quand le centre est inférieur à six signes. Son
maximum, pour Saturne, est de six parts, un tiers de part et un
cinquième de part, atteint quand l'élongation entre le centre du
déférent et l'Apogée est, dans un sens ou dans l'autre, trois signes
et quatre degrés.

La \emph{deuxième équation} est causée par l'épicycle. 

\newpage\phantomsection 
\index{APAQBH@\RL{_drw}!APAQBHAI BEAQAEBJBJAI@\RL{_dirwaT mar'iyyaT}, apogée apparent}
\index{APAQBH@\RL{_drw}!APAQBHAI AMBBBJBBBJBJAI@\RL{_dirwaT .haqIqiyyaT}, apogée vrai}
\index{BFASAH@\RL{nsb}!AOBBAGAEBB BFASAH@\RL{daqA'iq al-nisb}, coefficient d'interpolation}
\index{ACBHAL@\RL{'awj}!ACBHAL@\RL{'awj}, Apogée}
\index{AQBCAR@\RL{rkz}!BEAQBCAR@\RL{markaz}, centre}
\index{BHASAW@\RL{ws.t}!BHASAW@\RL{ws.t}, astre moyen}
\index{ANAUAU@\RL{_h.s.s}!ANAGAUAUAI BCBHBCAH@\RL{_hA.s.saT al-kwkab}, astre propre (traduit par ``anomalie'' par de nombreux auteurs)}
\includepdf[pages=59,pagecommand={\thispagestyle{plain}}]{edit.pdf}\phantomsection

\noindent C'est l'angle au centre du Monde entre deux droites qui en
sont issues~: l'une passant par le centre de l'épicycle, et l'autre
par l'astre. Son maximum est atteint quand celle-ci est tangente à
l'épicycle. Elle s'évanouit à l'apogée et au périgée apparents de
l'épicycle. Son maximum, pour Saturne, est six parts et treize
minutes. On l'ajoute à Saturne moyen quand l'astre propre corrigé
est inférieur à six signes, et on la retranche de Saturne moyen quand
l'astre propre corrigé est supérieur à cela.

La \emph{troisième équation} est la variation de la
deuxième\footnote{\textit{Cf.} chapitre 23, p. \pageref{troisieme_equation}
  \textit{infra}.}. En effet, la deuxième équation varie selon la plus
ou moins grande distance entre le centre de l'épicycle et le centre du
Monde. On ajoute ou on retranche cette troisième équation, selon qu'on
ajoute ou qu'on retranche la deuxième.

Le \emph{coefficient d'interpolation}, c'est un nombre dont le rapport à
soixante est comme le rapport entre ce qu'il faut de la troisième
anomalie et la totalité de cette anomalie. On clarifiera cela lors
du calcul des équations.

\begin{center}
  \large Section
\end{center}
Ceci étant admis, sache que l'\emph{Apogée de l'astre} est un arc de
l'orbe incliné entre deux droites issues de son centre~: l'une passant
par l'équinoxe de printemps\footnote{{\shatir} ne prend plus la peine
  de préciser <<~le point en face de l'équinoxe de printemps~>>,
  c'est-à-dire le point de l'orbe incliné situé à la même distance du
  n{\oe}ud que l'équinoxe de printemps sur l'écliptique (\textit{cf.}
  modèle de la Lune).}, et l'autre passant par le centre du déférent,
le centre du rotateur et le centre de l'épicycle quand le centre de
l'épicycle est à distance maximale du centre de l'orbe incliné. En
effet ces trois centres sont alignés à l'Apogée et au périgée~; et
j'appelle Apogée la distance maximale.

Le \emph{centre de l'astre} est un arc de l'orbe incliné entre deux
droites issues de son centre~: l'une passant par l'Apogée et l'autre
par le centre du déférent.

L'\emph{astre moyen} est la somme de l'Apogée et du centre. C'est un
arc de l'orbe incliné compris entre deux droites issues de son centre~:
l'une passant par le commencement du Bélier et l'autre par le centre
du déférent.

L'\emph{[astre] propre} est un arc de l'orbe de l'épicycle entre
l'astre et l'apogée vrai, dans le sens des signes.

L'\emph{[astre propre] corrigé} est un arc de cet orbe entre l'astre
et l'apogée apparent, dans le sens des signes.

L'\emph{apogée apparent} de cet orbe, c'est le point où le coupe la
droite passant par le centre de l'épicycle et le centre de l'orbe
incliné, dans la partie supérieure de l'épicycle (l'autre point est le
périgée apparent).

Menons une droite issue du centre du Monde (centre de l'écliptique et
centre de l'orbe incliné) et passant par l'astre, jusqu'à l'orbe de
l'écliptique.

\newpage\phantomsection 
\index{BBBHBE@\RL{qwm}!BEBBBHBE BCBHBCAH@\RL{mqwm al-kawkab}, astre vrai (= longitude vraie de l'astre)}
\index{BBBHBE@\RL{qwm}!AJBBBHBJBE@\RL{taqwIm}, 1) calcul de la longitude vraie, 2) longitude vraie d'un astre}
\includepdf[pages=60,pagecommand={\thispagestyle{plain}}]{edit.pdf}\phantomsection

\noindent Prenons un des cercles de latitude qui passe par
l'extrémité de cette droite. Il coupe l'orbe de l'écliptique en un
point qui est l'astre vrai\footnote{Nous traduisons \textit{mqwm
    al-kawkab} ainsi que \textit{taqw{\=\i}m al-kawkab} par <<~astre
  vrai~>>.}~; c'est-à-dire que l'élongation entre ce point et le
commencement du Bélier est l'astre vrai.

L'élongation entre ce point et l'astre, sur ce cercle, est la latitude
de l'astre. Nous avons mis la question des latitudes des astres dans
un chapitre séparé pour en faciliter la compréhension.

Ceci étant admis, sache qu'une première manière, c'est d'appliquer les
démonstrations géométriques aux trajectoires des centres des sphères
solides appartenant à Saturne.

Quant aux sphères, on les représente comme suit. Nous imaginons un
orbe circonscrit à l'ensemble des orbes de Saturne, dans le plan du
zodiaque et sur les pôles du zodiaque. C'est l'orbe parécliptique. Sa
surface extérieure touche la surface intérieure du huitième orbe, et
sa surface intérieure touche la surface extérieure de l'orbe incliné
de Saturne.

La surface intérieure de l'orbe incliné de Saturne touche la surface
extérieure de l'orbe parécliptique de Jupiter, et l'épaisseur de
l'orbe incliné est vingt-six degrés et trente minutes\footnote{Cette
  valeur est certainement corrompue~: l'épaisseur de l'orbe incliné ne
  peut être inférieure à $26;40$, diamètre de l'orbe déférent.}.  Son
plan est incliné par rapport au plan de l'orbe du zodiaque, de deux
parts et demi.

Ensuite, nous imaginons un globe entièrement sphérique, de rayon
treize degrés et un tiers de degré, appelé le déférent.

Nous imaginons une deuxième sphère, de rayon huit parts, douze minutes
et demi, plongée dans l'orbe déférent et le touchant en un point. Elle
s'appelle orbe rotateur.

Nous imaginons une troisième sphère, de rayon six parts et demi,
plongée dans l'orbe rotateur et le touchant en un point opposé au
point où le rotateur touche le déférent. Elle s'appelle orbe de
l'épicycle.

L'astre lui-même est plongé dans l'épicycle. Il le touche en un point
opposé au point où il touche le rotateur -- quand il est à l'Apogée.

Nous supposons que le déférent est plongé à l'intérieur de l'orbe
incliné et le touche au point de l'Apogée, puis que chaque orbe se
meut et que la grandeur et le sens de son mouvement sont comme nous
les avons supposés. Alors le parécliptique se meut autour des pôles du
zodiaque dans le sens des signes, et son mouvement est comme le
mouvement des Apogées. L'orbe incliné se meut dans le sens des signes,
d'un mouvement simple autour de son centre, de la grandeur du
mouvement du centre de Saturne, c'est-à-dire, par jour,
$0;2,0,26,17$. 

\newpage\phantomsection 
\index{AMAQBC@\RL{.hrk}!AMAQBCAI AHASBJAWAI BEAQBCBCAHAI@\RL{.harakaT basI.taT mrkkbaT}, mouvement simple-composé}
\includepdf[pages=61,pagecommand={\thispagestyle{plain}}]{edit.pdf}\phantomsection

\noindent L'orbe déférent se meut dans le sens contraire aux signes
dans sa partie supérieure~; son mouvement est de la grandeur du
mouvement du centre de Saturne. Le rotateur, dans sa partie
supérieure, se meut dans le sens des signes, d'une grandeur double de
celle du mouvement du centre de Saturne. Enfin, l'orbe de l'épicycle
se meut d'un mouvement simple autour de son centre, dans le sens des
signes dans sa partie supérieure~; son mouvement est l'excédent du
mouvement propre de Saturne sur le mouvement de son centre -- c'est le
mouvement simple-composé qu'on a déjà expliqué. Ainsi varient les
positions des orbes, et semblent aussi varier leurs mouvements qui
sont composés de mouvements simples.

A l'Apogée, l'astre est à distance maximale s'il est à l'apogée de
l'épicycle. De même, il est à distance minimale s'il est au périgée de
l'épicycle quand l'épicycle est au périgée du déférent. Entre Apogée
et périgée, il est à distance intermédiaire.\footnote{Quand il s'agit
  de l'orbe de l'épicycle, le mot \textit{\b{d}irwa} désigne
  en général un point situé \emph{dans la partie supérieure} de
  l'orbe, et non pas le point situé à l'apogée de l'orbe de l'épicycle
  \emph{relativement au centre de l'orbe portant} cet orbe. On peut
  aussi penser qu'{\shatir} fait alors abstraction des orbes
  intermédiaires entre l'orbe déférent et l'orbe de l'épicycle, comme
  il semblera le faire quand il parlera parfois d'un ``épicycle
  apparent''. L'orbe qui porte le centre de l'épicycle apparent est
  bien l'orbe déférent. L'apogée de l'épicycle est alors le point
  situé dans la ``partie supérieure'' de l'orbe de l'épicycle~: mais
  c'est le périgée de l'épicycle relativement au centre de l'orbe
  rotateur.}

Sache que l'astre Saturne n'atteint ni le point de ses orbes le plus
éloigné du centre du Monde, ni le point le plus proche. La distance
maximale de Saturne est, en parts, soixante-neuf, deux tiers et un
quart~; sa distance minimale est, en parts, cinquante et un
demi-sixième~; mais pour compléter les orbes et en faire 
des sphères circonscrites les unes aux autres, la distance maximale au
sein de l'orbe incliné de Saturne doit être (rayon de l'orbe incliné
plus rayon de l'orbe déférent), en parts, soixante-treize et un
tiers. En sus, il y a l'épaisseur du globe du parécliptique~; qu'on la
prenne égale à deux tiers de part, le total est soixante-quatorze
parts. La distance minimale est à la surface intérieure de l'orbe
incliné de Saturne, en parts, quarante-six et deux tiers~; en plus de
la réunion des orbes, il faut aussi compter le rayon de l'astre donc
cela fait quarante-six parts.\label{contiguite_sat}

Sache cela, et fais de manière analogue pour les autres astres. Ce sont
des modèles parfaits qui échappent aux doutes~; ce sont des sphères
parfaites, et chacun de leurs mouvements est un mouvement simple,
uniforme autour du centre de l'orbe qui se meut, un mouvement qui
décrit des portions égales en des temps égaux. Gloire à Dieu qui nous
a donné ceci et nous l'a appris.

\newpage\phantomsection
\includepdf[pages=62,pagecommand={\thispagestyle{plain}}]{edit.pdf}\phantomsection

\newpage\phantomsection
\begin{figure}[h!]
\centering
\noindent\hspace{-85mm}\hspace{.5\textwidth}\begin{minipage}{17cm}
\footnotesize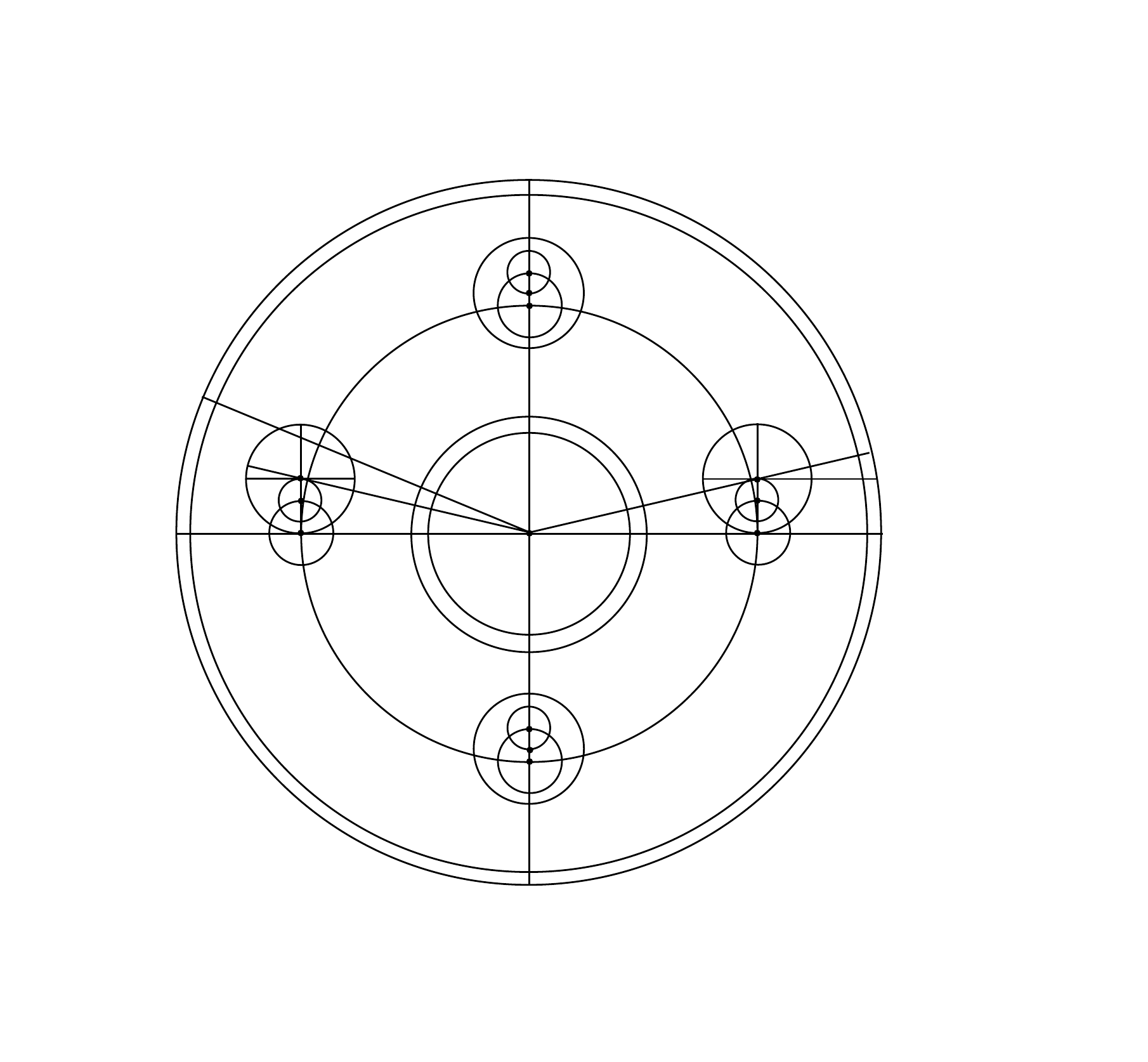
\end{minipage}
Les orbes de Saturne, représentées dans le plan par les ceintures
trajectoires des centres des orbes solides ou par les surfaces de ces
ceintures
\end{figure}

\newpage\phantomsection
\includepdf[pages=63,pagecommand={\thispagestyle{plain}}]{edit.pdf}\phantomsection

\newpage\phantomsection
\begin{figure}[h!]
\centering
\noindent\hspace{-85mm}\hspace{.5\textwidth}\begin{minipage}{17cm}
\footnotesize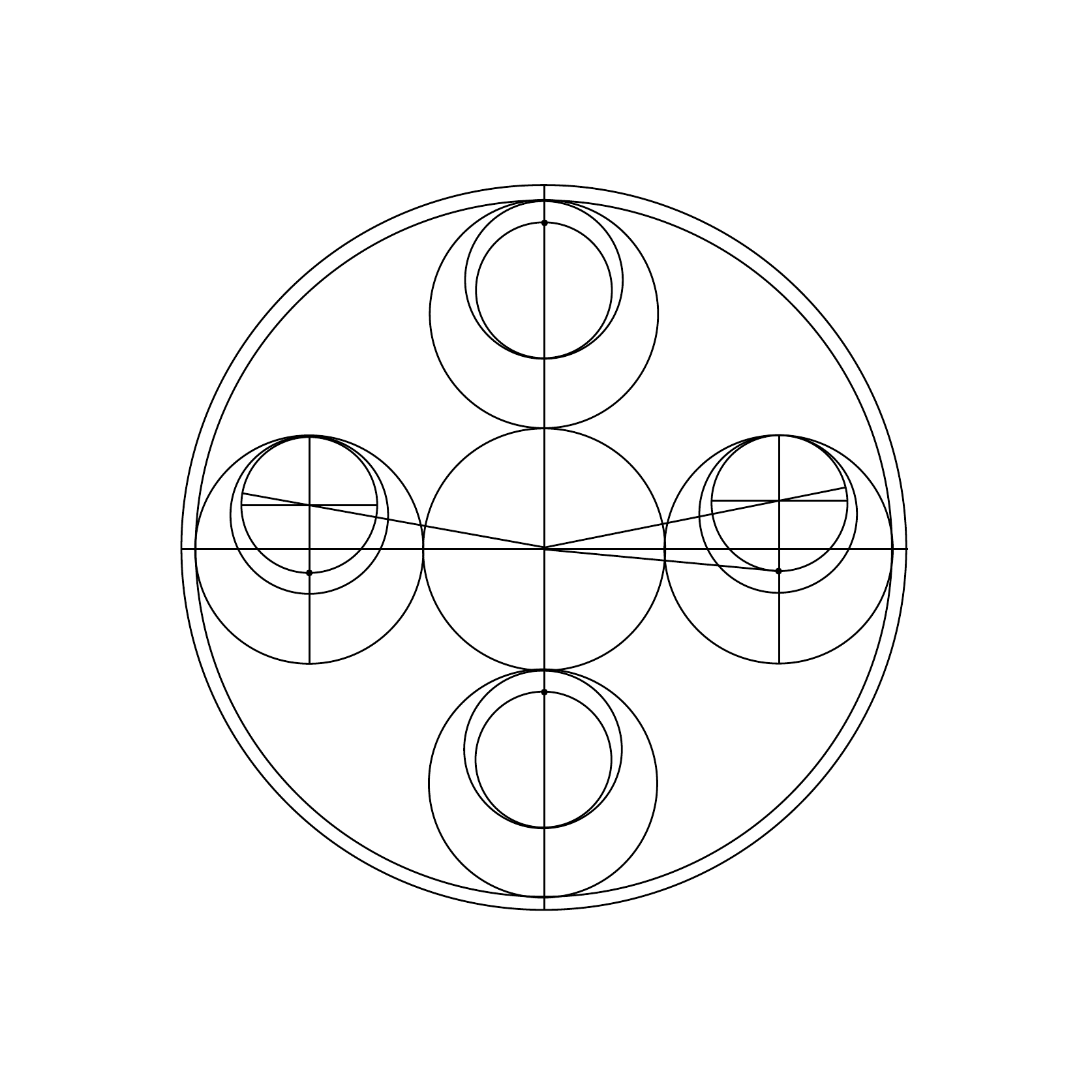
\end{minipage}
Les orbes de Saturne~: sphères entières représentées dans le plan à
l'Apogée, au périgée et aux élongations intermédiaires
\end{figure}

\newpage\phantomsection 
\index{ACAQAN@\RL{'ara_ha}!AJABAQBJAN BEAJBBAOAOBE@\RL{ta'rI_h mtqddm}, \'Epoque}
\index{ARAMBD@\RL{z.hl}!ARAMBD@\RL{zu.hal}, Saturne}
\addcontentsline{toc}{chapter}{I.13 La régulation des mouvements de Saturne et leurs conditions initiales à l'\'Epoque}
\includepdf[pages=64,pagecommand={\thispagestyle{plain}}]{edit.pdf}\phantomsection

\begin{center}
  \Large Chapitre treize

  \large La régulation des mouvements de Saturne et leurs conditions initiales à l'\'Epoque
\end{center}
Nous avons observé que Saturne moyen, à midi du premier jour de
l'année sept cent un de l'ère de Yazdgard, est 5 signes et $7;58,20$
degrés, et son Apogée 8 signes et $14;52$ degrés\footnote{Pour Saturne
  et pour les autres planètes, Kennedy et Roberts ont retrouvé les
  valeurs des Apogées à cette date par un calcul de précession au
  moyen des valeurs indiquées dans le \emph{Z{\=\i}j} de {\shatir }
  pour une autre date~; ceci confirme cette lecture (en particulier
  pour le chiffre 52) et explique aussi la coïncidence des rangs
  fractionnaires des Apogées des cinq planètes (voir
  \cite{roberts1959} p.~232).}. Le mouvement de Saturne moyen en vingt
années persanes est 8 signes et $4;33,20$ degrés~; son mouvement en
une année persane est $12;13,40$ degrés, en trente jour
$1;0,18,4,55,53,26$ degrés,
en un jour $0;2,0,36,9,51,46,50,57,32$, et en une
heure $0;0,5,1,30,25$ (en degré).

Le mouvement propre de chaque astre
est l'excédent du mouvement du Soleil moyen sur le mouvement de
l'astre moyen, et le mouvement de son centre est l'excédent du
mouvement de l'astre moyen sur le mouvement des Apogées, dont la
grandeur est d'un degré en soixante ans~; donc le mouvement propre [de
  Saturne] est, par jour, $0;57,7,43,34,21,15,5$,
et le mouvement de son centre est, par jour, $0;2,0,26,17$ (en
degré). Dieu est le plus savant.

\newpage\phantomsection 
\index{BBBHBE@\RL{qwm}!AJBBBHBJBE@\RL{taqwIm}, 1) calcul de la longitude vraie, 2) longitude vraie d'un astre}
\index{ALAOBD@\RL{jdl}!ALAOBHBD@\RL{jdwl}, table, tableau, catalogue}
\index{AYAOBD@\RL{`dl}!AJAYAOBJBD BEAMBCBE@\RL{ta`dIl mu.hakam}, équation principale}
\index{ARAMBD@\RL{z.hl}!ARAMBD@\RL{zu.hal}, Saturne}
\addcontentsline{toc}{chapter}{I.14 Détermination de Saturne vrai}
\includepdf[pages=65,pagecommand={\thispagestyle{plain}}]{edit.pdf}\phantomsection
\begin{center}
  \Large Chapitre quatorze

  \large Détermination de Saturne vrai
\end{center}
Calcule Saturne moyen, son Apogée, et le Soleil moyen à la date
souhaitée, puis retranche Saturne moyen du Soleil moyen, il reste
Saturne propre, puis retranche l'Apogée de Saturne moyen, il reste son
centre, puis consulte l'entrée de la table des équations de Saturne
correspondant au centre de Saturne, et prends-y la première
équation. Si le centre de Saturne est inférieur à six signes,
retranche l'équation de Saturne moyen et du centre, et ajoute-la au
mouvement propre~; mais si le centre est supérieur à six signes, alors
ajoute l'équation à l'astre moyen et au centre et retranche-la de l'astre
propre. On obtient l'astre moyen corrigé, le centre corrigé, et l'astre
propre corrigé. Puis consulte l'entrée correspondant au centre
corrigé dans la table des équations de Saturne, prends-y le coefficient
d'interpolation, regarde s'il est additif ou soustractif, et
mets-le de côté. Puis consulte l'entrée correspondant au mouvement
propre corrigé dans la table des équations de Saturne, et prends-y la
troisième et la quatrième équation si le coefficient d'interpolation
est additif -- sinon prends-y la deuxième et la troisième. Puis
multiplie la deuxième ou la quatrième par le coefficient d'interpolation, et
ajoute le résultat à la troisième équation si le coefficient d'interpolation
est additif -- sinon retranche-le de la troisième équation. On
obtient ainsi la troisième [équation] \emph{principale}~; ajoute-la à
l'astre moyen corrigé si le mouvement propre corrigé est supérieur à
six signes, mais s'il est inférieur à six signes alors retranche-la de
l'astre moyen. On obtient l'astre vrai par rapport à l'écliptique à la
date pour laquelle tu as fait ce calcul. Sache cela.\footnote{{\shatir
  } explique ainsi comment calculer la <<~troisième équation
  principale~>>, c'est-à-dire l'équation due au mouvement propre de
  l'épicycle, en s'aidant des tables. Il doit y avoir trois colonnes numérotées
  $C_2$, $C_3$, $C_4$, ainsi qu'un multiplicateur $\lambda$. Il dit de
  calculer $C_3+\lambda C_4$ ou bien $C_3-\lambda C_2$ suivant que
  $\lambda$ est <<~additif~>> ou <<~soustractif~>>.}

\newpage\phantomsection
\index{ATAQBJ@\RL{^sry}!BEATAJAQBJ@\RL{al-mu^starI}, Jupiter}
\index{BABDBC@\RL{flk}!BABDBCBEBEAKAKBD@\RL{falak muma_t_tal}, parécliptique}
\index{BABDBC@\RL{flk}!BABDBCAMAGBEBD@\RL{falak .hAmil}, orbe déférent}
\index{BABDBC@\RL{flk}!BABDBCBEAOBJAQ@\RL{falak mdIr}, orbe rotateur}
\index{BABDBC@\RL{flk}!BABDBC AJAOBHBJAQ@\RL{falak al-tadwIr}, orbe de l'épicycle}
\index{AQABAS@\RL{ra'asa}!AQABAS@\RL{ra's}, tête, n{\oe}ud ascendant}
\index{APBFAH@\RL{_dnb}!APBFAH@\RL{_dnb}, queue, n{\oe}ud descendant}
\addcontentsline{toc}{chapter}{I.15 Configuration nouvelle des orbes de Jupiter}
\includepdf[pages=66,pagecommand={\thispagestyle{plain}}]{edit.pdf}\phantomsection
\begin{center}
  \Large Chapitre quinze

  \large Configuration nouvelle des orbes de Jupiter
\end{center}
La configuration des orbes de Jupiter est comme celle des orbes de
Saturne (vue précédemment) sauf en ce qui concerne les grandeurs des
orbes et des mouvements.

Parmi les orbes de Jupiter, nous imaginons~:

-- un orbe dans le plan de l'écliptique, sur ses pôles et
en son centre~: on l'appelle \emph{parécliptique}.

-- un deuxième orbe, incliné par rapport au parécliptique,
d'inclinaison constante égale à une part et demi. Il coupe le
parécliptique en deux points opposés dont l'un s'appelle la
\emph{tête}, l'autre la \emph{queue}.

-- un troisième orbe dont le centre soit sur le bord de l'orbe
incliné et dont le rayon soit quatre parts et un huitième (en parts
telles que le rayon de l'orbe incliné en compte soixante)~; on
l'appelle \emph{déférent}.

-- un quatrième orbe dont le centre soit sur le bord du
déférent et dont le rayon soit une part, vingt-deux minutes et
demi. On l'appelle \emph{rotateur}.

-- un cinquième orbe dont le centre soit sur le bord du
rotateur et dont le rayon soit onze parts et demi (en les mêmes
parts). On l'appelle \emph{orbe de l'épicycle}.

-- le centre du corps de Jupiter est situé sur
la ceinture de l'épicycle.

Supposons à présent que les centres des orbes déférent, rotateur, de
l'épicycle et de l'astre soient tous alignés sur la droite issue du
centre du parécliptique dans la direction de l'Apogée. Supposons que
le parécliptique est mû d'un mouvement simple sur ses pôles, dans le
sens des signes, comme le mouvement des Apogées~: en un jour et une
nuit, $0;0,0,9,52$. La tête et la queue sont entraînées par ce
mouvement, ainsi que les deux parties d'inclinaison maximale. L'orbe
incliné se meut d'un mouvement simple autour de son centre (qui est
aussi le centre du parécliptique), dans le sens des signes, comme le
mouvement du centre de Jupiter, c'est-à-dire l'excédent de Jupiter
moyen sur le mouvement de l'Apogée~: en un jour et une nuit,
$0;4,59,6,14,35$. Le déférent se meut en sens inverses des signes dans
sa partie supérieure, d'autant que le mouvement du centre de Jupiter
aussi~: en un jour et une nuit, $0;4,59,6,14,35$. Le rotateur se meut,
dans sa partie supérieure, dans le sens des signes, d'autant que le
double du mouvement du centre de Jupiter~: en un jour et une nuit,
$0;9,58,12,29,10$. 

\newpage\phantomsection
\index{AMAQBC@\RL{.hrk}!AMAQBCAI AHASBJAWAI BEAQBCBCAHAI@\RL{.harakaT basI.taT mrkkbaT}, mouvement simple-composé}
\index{AMAQBC@\RL{.hrk}!AMAQBCAI BEANAJBDBAAI@\RL{.harakaT mu_htalifaT}, mouvement irrégulier (\textit{i. e.} non uniforme)}
\index{BABDBC@\RL{flk}!BABDBC@\RL{flk}, 1) orbe, 2) cieux}
\index{AYAOBD@\RL{`dl}!BEAYAOAOBD@\RL{mu`addal}, corrigé}
\index{AYAOBD@\RL{`dl}!BEAQBCAR BEAYAOAOBD@\RL{markaz mu`addal}, centre corrigé}
\index{AYAOBD@\RL{`dl}!ANAGAUAUAI BEAYAOAOBD@\RL{_hA.s.saT mu`addalaT}, astre propre corrigé}
\index{BBBHBE@\RL{qwm}!BEBBBHBE BCBHBCAH@\RL{mqwm al-kawkab}, astre vrai (= longitude vraie de l'astre)}
\index{BHASAW@\RL{ws.t}!BHASAW@\RL{ws.t}, astre moyen}
\index{AYAQAV@\RL{`r.d}!AYAQAV@\RL{`r.d}, latitude, \emph{i. e.} par rapport à l'écliptique}
\index{ALASBE@\RL{jsm}!BCAQAGAJ BEALASBEAI AJAGBEBEAI@\RL{kraT mjsm, falak mjsm}, sphère solide, orbe solide}
\index{APAQBH@\RL{_drw}!APAQBHAI AJAOBHBJAQ BHASAWBI@\RL{_dirwaT al-tadwIr al-ws.t_A}, apogée moyen de l'épicycle}
\includepdf[pages=67,pagecommand={\thispagestyle{plain}}]{edit.pdf}\phantomsection

\noindent L'orbe de l'épicycle se meut autour de son centre,
dans sa partie supérieure, dans le sens des signes, de la grandeur de
l'excédent du mouvement propre de Jupiter sur le mouvement de son
centre~: en un jour et une nuit, $0;49,9,57,22,24,10$. Ceci est le
mouvement simple-composé, donc l'astre s'éloigne de l'apogée moyen,
chaque jour, d'autant que le mouvement propre de Jupiter qui est
$0;54,9,3,36,59,10$ (comme nous l'avons expliqué à propos des orbes de
Saturne).

Si les orbes continuent à se mouvoir des mêmes mouvements, il
s'avère que l'astre aura un mouvement irrégulier composé de
mouvements simples.

Les tenants et les aboutissants sont comme nous l'avons expliqué dans
les orbes de Saturne.

Nous représentons la configuration des orbes de Jupiter réduites aux
cercles [des trajectoires des centres des orbes] selon ce qu'on
projette sur la figure comme nous l'avons représenté avec les orbes de
Saturne~; puis nous les représentons en tant que sphères solides comme
on l'imagine dans les cieux.\label{cieux}

Mentionnons à présent leurs grandeurs et distances. Sache que la
droite issue du centre du Monde dirigée vers le centre du déférent
passe par l'astre moyen. Celle issue du centre du Monde et passant par
le centre de l'épicycle passe, dans l'écliptique, par le centre
corrigé de l'astre, et elle est délimite, dans l'épicycle,
l'astre propre corrigé\footnote{Plus précisément, cette droite coupe
  l'épicycle en un point qui est son <<~apogée apparent~>> et qui sert
  d'origine pour mesurer l'astre propre corrigé.
  \textit{Cf.} fig. \ref{fig_term} p.~\pageref{fig_term}.}. La
droite issue du centre du Monde et passant par l'astre passe, dans
l'écliptique, par l'astre vrai. [Tout ceci], si l'astre est sans
latitude~; s'il a une latitude, alors les cercles passant par les
pôles de l'écliptique et par les extrémités des droites mentionnées
coupent l'écliptique en l'astre moyen, le centre, et l'astre vrai,
comme on l'a établi.

Ceci dit, sache que la distance maximale de l'astre de Jupiter au
centre du Monde est soixante-quatorze parts et un quart, et que sa
distance minimale est quarante-cinq parts, une demi-part et un quart~;
la distance maximale du centre de l'épicycle au centre du Monde est
$62;45$, et sa distance minimale est $57;15$~; sauf que Jupiter
n'atteint pas la distance maximale de ses orbes solides, ni leur
distance minimale.

\newpage\phantomsection
\index{BHAUBD@\RL{w.sl}!BEAJAJAUBD AHAYAVBG AHAHAYAV@\RL{mutta.sil ba`.dh biba`.d}, contigus}
\includepdf[pages=68,pagecommand={\thispagestyle{plain}}]{edit.pdf}\phantomsection

Quant aux grandeurs de ces orbes, je les décris comme suit. Le rayon
de la sphère du déférent est dix-sept parts. Le rayon de la sphère du
rotateur est douze parts et cinquante-deux minutes et demi. Le rayon
de la sphère de l'épicycle est onze parts et demi. Le rayon de l'orbe
incliné est soixante-dix-sept. Sa face externe touche la face interne
de son orbe parécliptique. L'épaisseur du parécliptique peut être
choisie arbitrairement, et on la pose égale à un degré. La distance
minimale dans l'orbe incliné est quarante-trois, et même inférieure à
cela de la grandeur du rayon de la sphère de l'astre avec [ce qu'il faut]
pour que les orbes soient contigus, comme nous l'avons dit~; on pose
ceci égal à un degré.

\newpage\phantomsection
\includepdf[pages=69,pagecommand={\thispagestyle{plain}}]{edit.pdf}\phantomsection

\newpage\phantomsection
\begin{figure}[h!]
\centering
\noindent\hspace{-85mm}\hspace{.5\textwidth}\begin{minipage}{17cm}
\footnotesize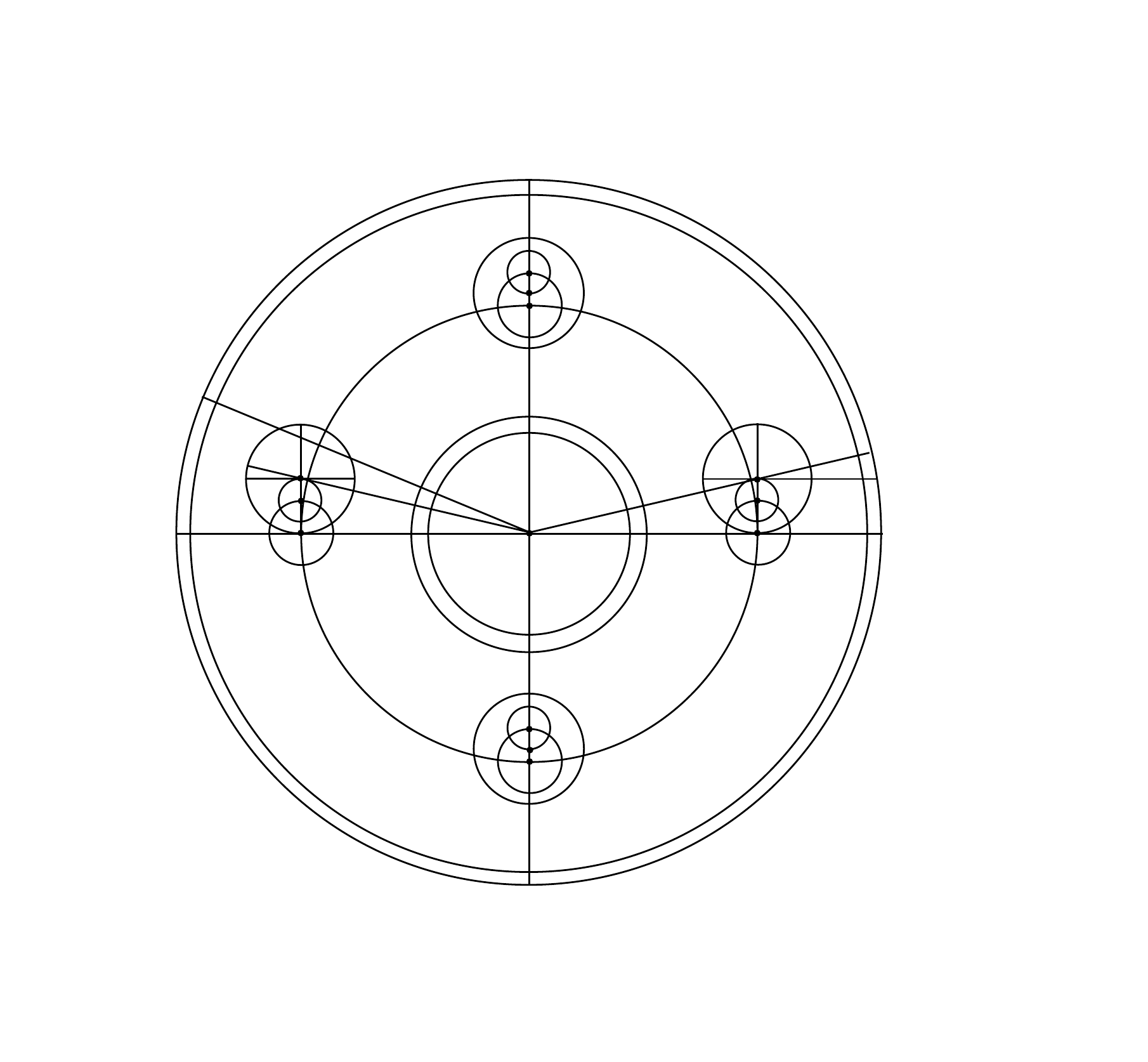
\end{minipage}
Les orbes de Jupiter représentées par les trajectoires des centres des
sphères solides à l'Apogée, au périgée, et aux élongations
intermédiaires
\end{figure}

\newpage\phantomsection
\includepdf[pages=70,pagecommand={\thispagestyle{plain}}]{edit.pdf}\phantomsection

\newpage\phantomsection
\begin{figure}[h!]
\centering
\noindent\hspace{-85mm}\hspace{.5\textwidth}\begin{minipage}{17cm}
\footnotesize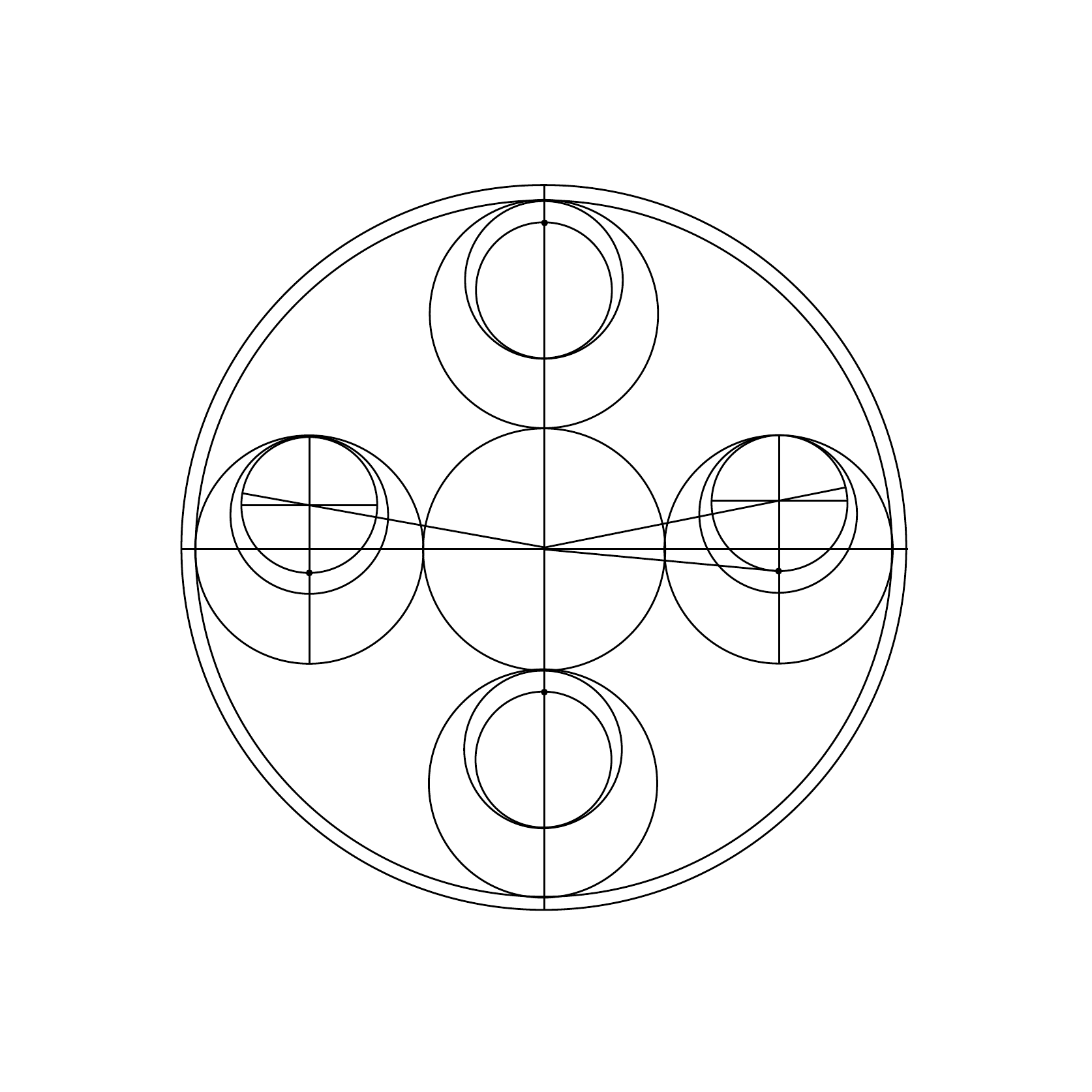
\end{minipage}
Configurations des orbes solides de Jupiter à l'Apogée, au périgée et
aux élongations intermédiaires
\end{figure}

\newpage\phantomsection
\index{BBBHBE@\RL{qwm}!AJBBBHBJBE@\RL{taqwIm}, 1) calcul de la longitude vraie, 2) longitude vraie d'un astre}
\index{AMAQBC@\RL{.hrk}!AMAQBCAI ANAGAUAI@\RL{.harakaT _hA.saT}, mouvement propre}
\index{AMAQBC@\RL{.hrk}!AMAQBCAI BHASAW@\RL{.harakaT al-ws.t}, mouvement de l'astre moyen}
\index{AMAQBC@\RL{.hrk}!AMAQBCAI BEAQBCAR@\RL{.harakaT al-markaz}, mouvement du centre (en général, le centre d'un épicycle)}
\index{ATAQBJ@\RL{^sry}!BEATAJAQBJ@\RL{al-mu^starI}, Jupiter}
\addcontentsline{toc}{chapter}{I.16 Régulation des mouvements de Jupiter}
\includepdf[pages=71,pagecommand={\thispagestyle{plain}}]{edit.pdf}\phantomsection
\begin{center}
  \Large Chapitre seize

  \large Régulation des mouvements de Jupiter
\end{center}
Leurs conditions initiales à l'\'Epoque, c'est-à-dire à midi du
premier jour de l'année
sept cent un de l'ère de Yazdgard, sont de 9 signes et $2;6,10$ degrés
pour Jupiter moyen et de 6 signes et $0;52$ degrés pour son Apogée. Le
mouvement de Jupiter moyen en vingt années persanes est 8 signes et
$6;51,0$ degrés, en une année 1 signe et $0;20,33$ degrés, en un mois
persan $2;29,38,3,17,16$, en un jour $0;4,59,16,6,34,32$, et en une
heure $0;0,12,28,10,16,26,20$. Le mouvement des Apogées est d'un degré
en soixante ans.

Le mouvement du centre est l'excédent du mouvement de l'astre
moyen sur [le mouvement des] Apogées. Le mouvement propre est, pour
les planètes
supérieures, l'excédent du mouvement du Soleil moyen sur le mouvement
de l'astre moyen~; c'est pour Jupiter, en un jour et une nuit,
$0;54,9,3,36,59,10$.

Jupiter vrai [se calcule] de la même manière que Saturne Vrai.

\newpage\phantomsection
\index{AQBJAN@\RL{ry_h}!BEAQAQBJAN@\RL{marrI_h}, Mars}
\index{BABDBC@\RL{flk}!BABDBCBEBEAKAKBD@\RL{falak muma_t_tal}, parécliptique}
\index{BABDBC@\RL{flk}!BABDBCAMAGBEBD@\RL{falak .hAmil}, orbe déférent}
\index{BABDBC@\RL{flk}!BABDBCBEAOBJAQ@\RL{falak mdIr}, orbe rotateur}
\index{BABDBC@\RL{flk}!BABDBC AJAOBHBJAQ@\RL{falak al-tadwIr}, orbe de l'épicycle}
\index{AQABAS@\RL{ra'asa}!AQABAS@\RL{ra's}, tête, n{\oe}ud ascendant}
\index{APBFAH@\RL{_dnb}!APBFAH@\RL{_dnb}, queue, n{\oe}ud descendant}
\addcontentsline{toc}{chapter}{I.17 La configuration des orbes de Mars
  d'après notre description, avec une notice sur les cercles
  trajectoires des centres des sphères solides}
\includepdf[pages=72,pagecommand={\thispagestyle{plain}}]{edit.pdf}\phantomsection
\begin{center}
  \Large Chapitre dix-sept

  \large La configuration des orbes de Mars d'après notre
description, avec une notice sur les cercles trajectoires des centres
des sphères solides
\end{center}
Nous imaginons~:

-- un orbe parécliptique représentant l'écliptique, sur ses pôles
  et dans son plan.

-- un deuxième orbe dont le plan est incliné d'une seule part par
  rapport au plan du parécliptique et qui le coupe en deux points
  opposés appelés la tête et la queue.
  
-- un troisième orbe dont le centre est sur le bord de l'orbe
  incliné et dont le rayon est neuf parts (en parts telles que le
  rayon de l'orbe incliné en compte soixante)~; on l'appelle déférent.

-- un quatrième orbe dont le centre est sur le bord du déférent et
dont le rayon est trois parts~; on l'appelle rotateur.

-- un cinquième orbe dont le centre est sur le bord du rotateur et
dont le rayon est trente-neuf parts et demi (en les mêmes parts)~;
on l'appelle orbe de l'épicycle.

-- le centre de l'astre est sur le bord de l'orbe de l'épicycle.

Supposons que les centres du déférent, du rotateur et de l'épicycle
sont à l'Apogée, c'est-à-dire sur la droite issue du centre du Monde
dans la direction de l'Apogée. Supposons que le parécliptique est mû
d'un mouvement simple autour de son centre et sur les pôles de
l'écliptique, dans le sens des signes, comme le mouvement des
Apogées~: en un jour et une nuit, $0;0,0,9,52$. La tête, la queue et
les deux parties d'inclinaison maximale sont entraînées par ce
mouvement. L'orbe incliné se meut dans le sens des signes, d'un
mouvement simple autour de son centre, comme le mouvement du centre de
Mars~: en un jour et une nuit $0;31,26,29,45$. Le déférent se meut, en
sens inverse des signes, d'autant que le mouvement du centre de
Mars. Le rotateur se meut, dans le sens des signes dans sa partie
supérieure, d'autant que le double du mouvement du centre de Mars~: en
un jour et une nuit, $1;2,52,59,30$. L'orbe de l'épicycle se meut
autour de son centre d'un mouvement simple, en sens inverse des signes
dans sa partie supérieure, d'autant que l'excédent du mouvement du
centre de Mars sur son mouvement propre.

\newpage\phantomsection
\index{AMAQBC@\RL{.hrk}!AMAQBCAI BEANAJBDBAAI@\RL{.harakaT mu_htalifaT}, mouvement irrégulier (\textit{i. e.} non uniforme)}
\index{AMAQBC@\RL{.hrk}!AMAQBCAI AHASBJAWAI BEAQBCBCAHAI@\RL{.harakaT basI.taT mrkkbaT}, mouvement simple-composé}
\index{AMAQBC@\RL{.hrk}!AMAQBCAI ANAGAUAI@\RL{.harakaT _hA.saT}, mouvement propre}
\index{AMAQBC@\RL{.hrk}!AMAQBCAI BHASAW@\RL{.harakaT al-ws.t}, mouvement de l'astre moyen}
\index{AQALAY@\RL{rj`}!AQALBHAY@\RL{rujU`}, rétrogradation}
\index{AOBHAQ@\RL{dwr}!BEAQAOAGAQAGAJBEAQAGBCARBCAQ@\RL{madArAt marAkaz al-akr}, trajectoires des centres des orbes}
\includepdf[pages=73,pagecommand={\thispagestyle{plain}}]{edit.pdf}\phantomsection

\noindent Ainsi l'orbe de l'épicycle
semble se mouvoir dans le sens des signes dans sa partie supérieure
d'autant que le mouvement propre de Mars, en un jour, $0;27,41,50$~;
car il tourne de par le mouvement du rotateur et de par son propre
mouvement qui est, autour de son centre, en un jour et une nuit,
$0;3,44,40$~; retranchant cette grandeur du mouvement du centre en un
jour ($0;31,26,29,45$), il reste $0;27,41,50$ et c'est le mouvement
apparent de l'astre par l'épicycle, dans le sens des signes, égal au
mouvement propre de Mars, c'est-à-dire l'excédent du mouvement du
Soleil moyen sur le mouvement de l'astre moyen. C'est le mouvement
simple-composé déjà mentionné.\footnote{Il y a une petite incohérence
  dans la
  valeur du mouvement propre, au niveau du troisième rang après la
  virgule~: le 50 devrait être un 40 si on calcule le mouvement propre
  comme différence entre Soleil moyen et Mars moyen. Cela peut venir
  d'une confusion entre Mars moyen et son centre (le mouvement des
  Apogées est environ $0;0,0,10$).}

Voyons un exemple~; supposons donc que les centres du déférent, du
rotateur et de l'épicycle sont sur la droite issue du centre du Monde
dans la direction de l'Apogée, et que se meuve chaque orbe d'un
mouvement simple autour de son centre de la grandeur et dans le sens
qu'on a supposés. Est engendré en l'astre un mouvement irrégulier
et il [lui] arrive de rétrograder (l'observation en rend compte).
C'est le mouvement irrégulier composé de mouvements simples comme
on l'a expliqué ailleurs.

Voyez la figure des orbes de Mars où les trajectoires des centres des
sphères solides sont projetées dans le plan pour qu'on se les
représente facilement.

Ceci étant dit, sache que la distance maximale du centre de l'épicycle
de Mars au centre du Monde est soixante-six parts, et que sa distance
minimale est cinquante-quatre parts. La distance maximale de l'astre
Mars au centre du Monde est cent cinq parts et demi, et sa distance
minimale est quatorze parts et demi (tout cela en parts telles que le
rayon de l'orbe incliné en compte soixante).

\newpage\phantomsection
\index{ALASBE@\RL{jsm}!BCAQAGAJ BEALASBEAI AJAGBEBEAI@\RL{kraT mjsm, falak mjsm}, sphère solide, orbe solide}
\index{BHAUBD@\RL{w.sl}!BEAJAJAUBD AHAYAVBG AHAHAYAV@\RL{mutta.sil ba`.dh biba`.d}, contigus}
\includepdf[pages=74,pagecommand={\thispagestyle{plain}}]{edit.pdf}\phantomsection
\noindent Cependant l'astre Mars
n'atteint pas la distance maximale de ses orbes, ni leur distance
minimale~; car le rayon du déférent est cinquante-et-une parts et
demi, le rayon de la sphère du rotateur est quarante-deux parts et
demi, et le rayon de la sphère de l'épicycle est trente-neuf parts et
demi, or la distance maximale de l'orbe incliné (c'est son rayon) est
cent onze degrés et demi, au-dessus desquels il y a l'épaisseur du
parécliptique, soit trente minutes, cela fait cent douze parts de
distance maximale, et la distance minimale [de l'orbe incliné] est
huit parts et trente minutes, et avec [ce qu'il faut] pour que
  les orbes soient contigus, soit trente minutes, il reste huit
parts. Conclusion~: la distance minimale des orbes de Mars est huit
parts, et leur distance maximale est cent douze parts. C'est ce sur
quoi nous nous sommes appuyés.

Voyez la figure des orbes solides de Mars tels qu'on les
imagine dans les cieux~; ce sont des sphères entières dont on a
présenté les rayons et les positions comme on l'a fait en l'expliquant
pour Saturne et les autres, donc il n'est pas besoin de répéter cela.

\newpage\phantomsection
\includepdf[pages=75,pagecommand={\thispagestyle{plain}}]{edit.pdf}\phantomsection

\newpage\phantomsection
\begin{figure}[h!]
\centering
\noindent\hspace{-85mm}\hspace{.5\textwidth}\begin{minipage}{17cm}
\footnotesize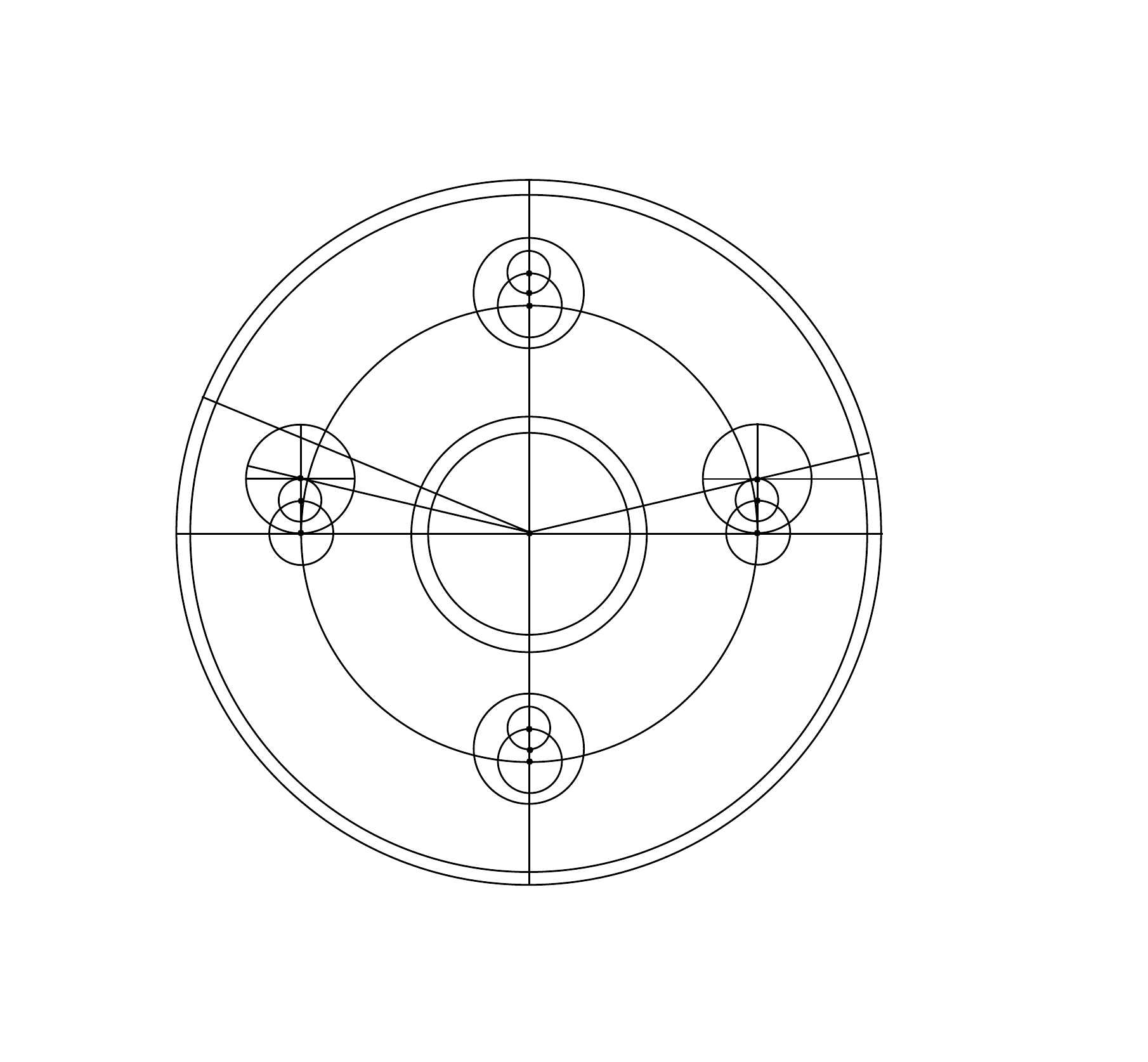
\end{minipage}
\label{mars_trajectoires}
Configuration des orbes de Mars représentées dans le plan par les
trajectoires des centres des sphères solides à l'Apogée, au périgée et
aux élongations intermédiares
\end{figure}

\newpage\phantomsection
\includepdf[pages=76,pagecommand={\thispagestyle{plain}}]{edit.pdf}\phantomsection

\newpage\phantomsection
\begin{figure}[h!]
\centering
\noindent\hspace{-85mm}\hspace{.5\textwidth}\begin{minipage}{17cm}
\footnotesize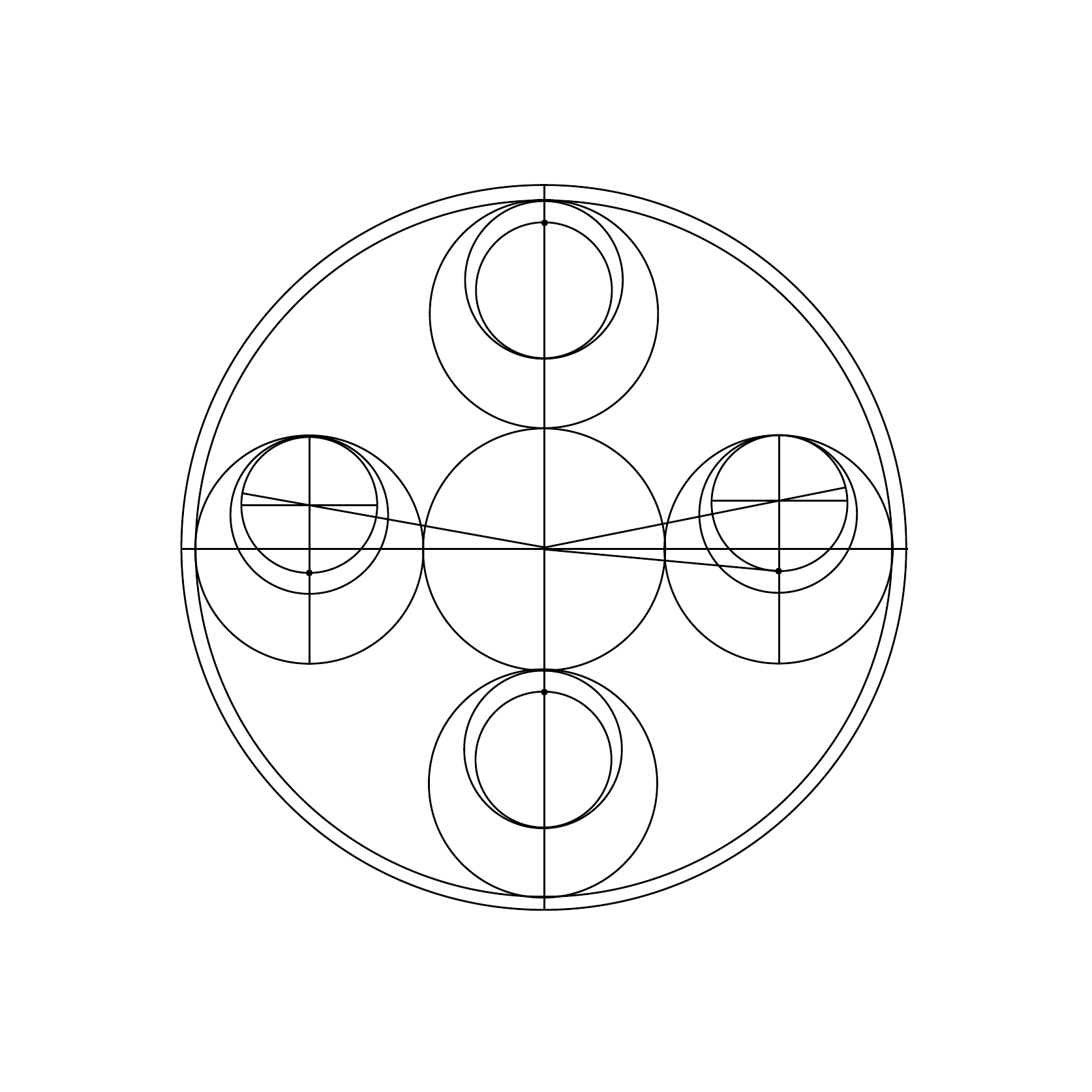
\end{minipage}
Configuration des orbes de Mars~: sphères entières représentées dans
le plan à l'Apogée, au périgée et aux élongations intermédiaires
\end{figure}

\newpage\phantomsection
\index{AQBJAN@\RL{ry_h}!BEAQAQBJAN@\RL{marrI_h}, Mars}
\index{BBBHBE@\RL{qwm}!AJBBBHBJBE@\RL{taqwIm}, 1) calcul de la longitude vraie, 2) longitude vraie d'un astre}
\addcontentsline{toc}{chapter}{I.18 Régulation des mouvements de Mars}
\includepdf[pages=77,pagecommand={\thispagestyle{plain}}]{edit.pdf}\phantomsection
\begin{center}
  \Large Chapitre dix-huit

  \large La régulation des mouvements de Mars
\end{center}
Mars moyen, mardi midi à Damas, le premier jour de l'année sept-cent
un de l'ère
de Yazdgard, est 9 [signes] et $22;0,0$ [degrés], et son Apogée 4 [signes]
et $17;52,0$ [degrés]. Le mouvement
de Mars moyen en vingt années persanes est 8 [signes] et $15;43,20$
[degrés]\footnote{modulo $360°$.}~; en une année
persane, 6 [signes] et $11;17,11,0$~; en trente jours,
$0;15,43,19,48,29,39$~; en un jour et une nuit,
$0;0,31,26,39,36,59,18$~; et en une heure $0;0,1,18,36,39$. Le
mouvement des apogées est d'un degré en soixante ans~; en un ans,
une minute~; en un jour, $0;0,0,9,52$. Le mouvement propre de Mars est
d'autant que l'excédent du mouvement du Soleil moyen sur Mars moyen,
c'est-à-dire par jour $0;0,27,41,40,7$. Mars vrai se calcule de la même
manière que Saturne vrai.

\newpage\phantomsection
\index{ARBGAQ@\RL{zhr}!ARBGAQAI@\RL{zuharaT}, Vénus}
\index{BABDBC@\RL{flk}!BABDBCBEBEAKAKBD@\RL{falak muma_t_tal}, parécliptique}
\index{BABDBC@\RL{flk}!BABDBCAMAGBEBD@\RL{falak .hAmil}, orbe déférent}
\index{BABDBC@\RL{flk}!BABDBCBEAOBJAQ@\RL{falak mdIr}, orbe rotateur}
\index{BABDBC@\RL{flk}!BABDBC AJAOBHBJAQ@\RL{falak al-tadwIr}, orbe de l'épicycle}
\index{AQABAS@\RL{ra'asa}!AQABAS@\RL{ra's}, tête, n{\oe}ud ascendant}
\index{APBFAH@\RL{_dnb}!APBFAH@\RL{_dnb}, queue, n{\oe}ud descendant}
\index{BEBJBD@\RL{myl}!AZAGBJAI BEBJBD@\RL{.gAyaT al-mIl}, partie d'inclinaison maximale (dans un orbe incliné)}
\index{AMAQBC@\RL{.hrk}!AMAQBCAI BEAQBCAR@\RL{.harakaT al-markaz}, mouvement du centre (en général, le centre d'un épicycle)}
\index{AMAQBC@\RL{.hrk}!AMAQBCAI BEAQBCAR ATBEAS@\RL{.harakaT markaz al-^sams}, mouvement du centre du Soleil}
\addcontentsline{toc}{chapter}{I.19 La configuration des orbes de Vénus
  selon notre méthode, avec une notice sur les cercles trajectoires
  des sphères solides} 
\includepdf[pages=78,pagecommand={\thispagestyle{plain}}]{edit.pdf}\phantomsection
\begin{center}
  \Large Chapitre dix-neuf

  \large La configuration des orbes de Vénus selon
notre méthode, avec une notice sur les cercles trajectoires des
sphères solides
\end{center}
Nous imaginons~:

-- un orbe parécliptique représentant l'écliptique, dans son plan
  et sur ses pôles

-- un deuxième orbe, au centre de l'écliptique, son plan étant
  incliné par rapport au plan du parécliptique, d'une inclinaison
  constante d'un sixième de part (selon l'opinion la plus juste), et
  le coupant en deux points opposés appelés la tête et la queue

-- un troisième orbe dont le centre est sur le bord de l'orbe
  incliné, et dont le rayon est une part et quarante-et-une minutes
  (en parts telles que le rayon de l'orbe incliné en compte
  soixante)~; on l'appelle déférent

-- un quatrième orbe dont le centre est sur le bord du déférent et
  dont le rayon est vingt-six minutes~; on l'appelle rotateur

-- un cinquième orbe dont le centre est sur le bord du rotateur et
  dont le rayon est, d'après nos observations, quarante-trois parts et
  trente-trois minutes (en les mêmes parts)~; on l'appelle orbe de
  l'épicycle

-- le centre de l'astre, sur le bord de ce dernier orbe, attaché en
  un de ses points

Supposons les centres du déférent, du rotateur et de l'épicycle
alignés sur la droite issue du centre de l'orbe incliné dans la
direction de l'Apogée. Supposons le parécliptique mû sur ses pôles,
dans le sens des signes, comme le mouvement des Apogées~: en un jour
et une nuit, $0;0,0,9,52$. Ce mouvement entraîne la tête et la queue,
et les deux parties d'inclinaison maximale. L'orbe incliné se meut dans
le sens des signes, comme le mouvement du centre de Vénus qui est comme
le mouvement du centre du Soleil~: en un jour et une nuit,
$0;59,8,9$. Le déférent se meut en sens inverse des signes, d'une
  grandeur égale au mouvement du centre de Vénus~; le rotateur se meut
  dans les sens des signes, d'une grandeur double du mouvement du
centre de Vénus~: en un jour et une nuit, $1;58,16,18$. 

\newpage\phantomsection
\index{AMAQBC@\RL{.hrk}!AMAQBCAI BEANAJBDBAAI@\RL{.harakaT mu_htalifaT}, mouvement irrégulier (\textit{i. e.} non uniforme)}
\index{AKAHAJ@\RL{_tbt}!AKAHAGAJ@\RL{_tibAt}, immobilité}
\index{AMAQBC@\RL{.hrk}!AMAQBCAI AHASBJAWAI BEAQBCBCAHAI@\RL{.harakaT basI.taT mrkkbaT}, mouvement simple-composé}
\index{AMAQBC@\RL{.hrk}!AMAQBCAI ANAGAUAI@\RL{.harakaT _hA.saT}, mouvement propre}
\index{APAQBH@\RL{_drw}!APAQBHAI@\RL{_dirwaT}, sommet, apogée}
\index{BHAUBD@\RL{w.sl}!BEAJAJAUBD AHAYAVBG AHAHAYAV@\RL{mutta.sil ba`.dh biba`.d}, contigus}
\index{ALASBE@\RL{jsm}!BCAQAGAJ BEALASBEAI AJAGBEBEAI@\RL{kraT mjsm, falak mjsm}, sphère solide, orbe solide}
\label{var26}
\includepdf[pages=79,pagecommand={\thispagestyle{plain}}]{edit.pdf}\phantomsection
\noindent L'orbe de l'épicycle se meut d'un mouvement simple autour de
son centre, en sens inverse des signes dans sa partie supérieure,
d'autant que l'excédent du mouvement du centre de Vénus sur son
mouvement propre~: en un jour et une nuit, $0;22,8,41,25,33,53$. Il
semble donc progresser dans le sens des signes d'autant que le
mouvement propre de Vénus qui est, en un jour et une nuit,
$0;36,59,28,26$. En effet, supposons l'épicycle immobile, et que se
meuvent l'orbe incliné, le déférent et le rotateur~: l'apogée se
déplace dans le sens des signes, d'autant que le mouvement du
centre. Supposons lui un autre mouvement, en sens inverse des signes,
d'autant que le mouvement du centre, alors l'apogée revient en son
lieu~; mais d'autre part, on a trouvé par l'observation que
[l'épicycle] se meut, par rapport à l'apogée, dans le sens des signes,
d'autant que le mouvement propre. Or on a démontré que, si on suppose
l'épicycle mû en sens inverse des signes, d'autant que l'excédent du
mouvement du centre sur le mouvement propre, alors l'apogée reste en
arrière d'autant que le mouvement propre~: c'est le mouvement simple
composé déjà démontré.

Ayant supposé les centres du déférent, du rotateur et de l'épicycle
alignés sur la droite issue du centre de l'écliptique dans la
direction de l'Apogée, si les orbes se meuvent d'autant qu'on l'a
supposé avec l'orientation ci-dessus, alors l'astre semblera avoir un
mouvement irrégulier composé de mouvements simples, et c'est ce qu'on
trouve par l'observation.

Ceci étant dit, sache que la distance maximale du centre de l'épicycle
de Vénus au centre du Monde est soixante-et-une parts et un quart de
part, et que sa distance minimale est cinquante-huit parts, une
demi-part et un quart de part. La distance maximale de Vénus au centre
du Monde est cent quatre parts et quarante-huit minutes, et sa
distance minimale est quinze parts et un cinquième de
part. En fait, l'astre de Vénus n'atteint ni la distance minimale de
ses orbes, ni leur distance maximale, car le rayon de la sphère du
déférent est quarante-cinq parts et deux tiers de parts (et le rayon
du rotateur est quarante-trois parts et cinquante-neuf minutes), donc
la distance maximale de l'orbe incliné est, pour Vénus, cent cinq
parts et deux tiers de part à quoi il faut ajouter l'épaisseur du
parécliptique supposée d'un tiers de part, et la distance minimale du
parécliptique est quatorze parts et un tiers de part dont il faut
retrancher un intervalle supposé d'un tiers de part~: la distance
minimale des orbes de Vénus semble donc être quatorze parts, et leur
distance maximale, cent six parts (peut-être plus, mais pas moins que
cela).

\newpage\phantomsection
\index{AOBHAQ@\RL{dwr}!BEAQAOAGAQAGAJBEAQAGBCARBCAQ@\RL{madArAt marAkaz al-akr}, trajectoires des centres des orbes}
\index{ALASBE@\RL{jsm}!BCAQAGAJ BEALASBEAI AJAGBEBEAI@\RL{kraT mjsm, falak mjsm}, sphère solide, orbe solide}
\index{BFAWBB@\RL{n.tq}!BEBFAWBBAI@\RL{min.taqaT}, ceinture}
\includepdf[pages=80,pagecommand={\thispagestyle{plain}}]{edit.pdf}\phantomsection
\begin{center}
  \large Remarque
\end{center}
Etant donné les rayons des trajectoires des centres des sphères (ce
sont les cercles représentés sur la figure), soit à calculer les
rayons des sphères solides pour tel astre. Voici la méthode~: la somme
du rayon de l'orbe de l'épicycle et du rayon du rotateur est le rayon
de la sphère du rotateur, ajoutes-y le rayon du déférent, on obtient
le rayon de la sphère du déférent. Si nous ajoutons le rayon de la
sphère du déférent à soixante, on obtient le rayon extérieur de l'orbe
incliné auquel il faut ajouter l'épaisseur du parécliptique~; et si
nous le retranchons de soixante, et que nous retranchons de ce qui reste
[ce qu'il faut] pour que les orbes soient contigus, alors on obtient
le rayon [intérieur] de l'orbe incliné.

Voyez la description des orbes de Vénus représentés par les trajectoires
des centres des sphères. [On obtient] la figure des orbes solides de Vénus
en considérant [les cercles] projetés dans le plan comme étant les ceintures
des sphères solides.

\newpage\phantomsection
\includepdf[pages=81,pagecommand={\thispagestyle{plain}}]{edit.pdf}\phantomsection

\newpage\phantomsection
\begin{figure}[h!]
\centering
\noindent\hspace{-85mm}\hspace{.5\textwidth}\begin{minipage}{17cm}
\footnotesize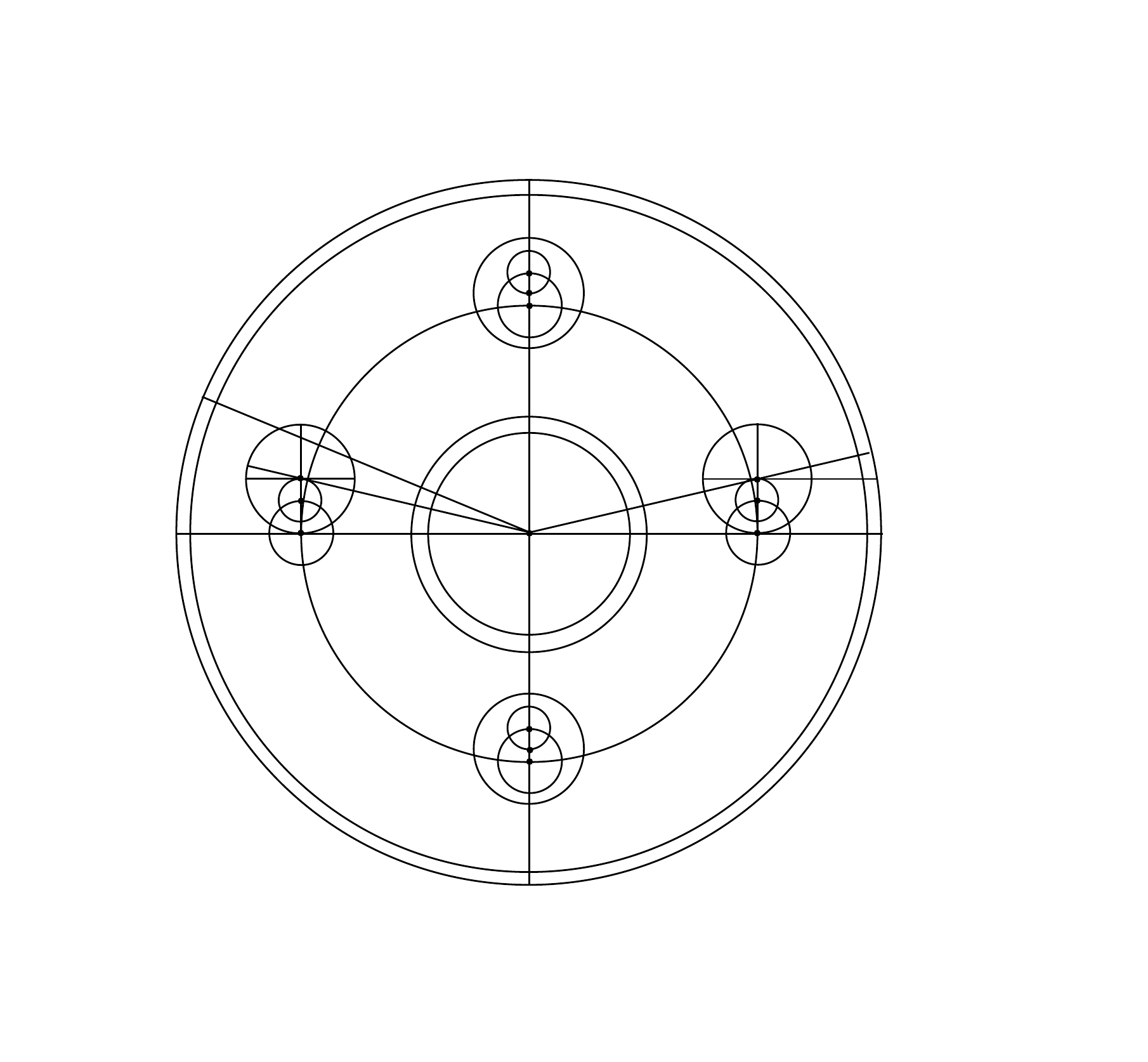
\end{minipage}
Configuration des orbes de Vénus représentées par les trajectoires des
centres des orbes solides à l'Apogée, au périgée et aux élongations
intermédiaires
\end{figure}

\newpage\phantomsection
\includepdf[pages=82,pagecommand={\thispagestyle{plain}}]{edit.pdf}\phantomsection

\newpage\phantomsection
\begin{figure}[h!]
\centering
\noindent\hspace{-85mm}\hspace{.5\textwidth}\begin{minipage}{17cm}
\footnotesize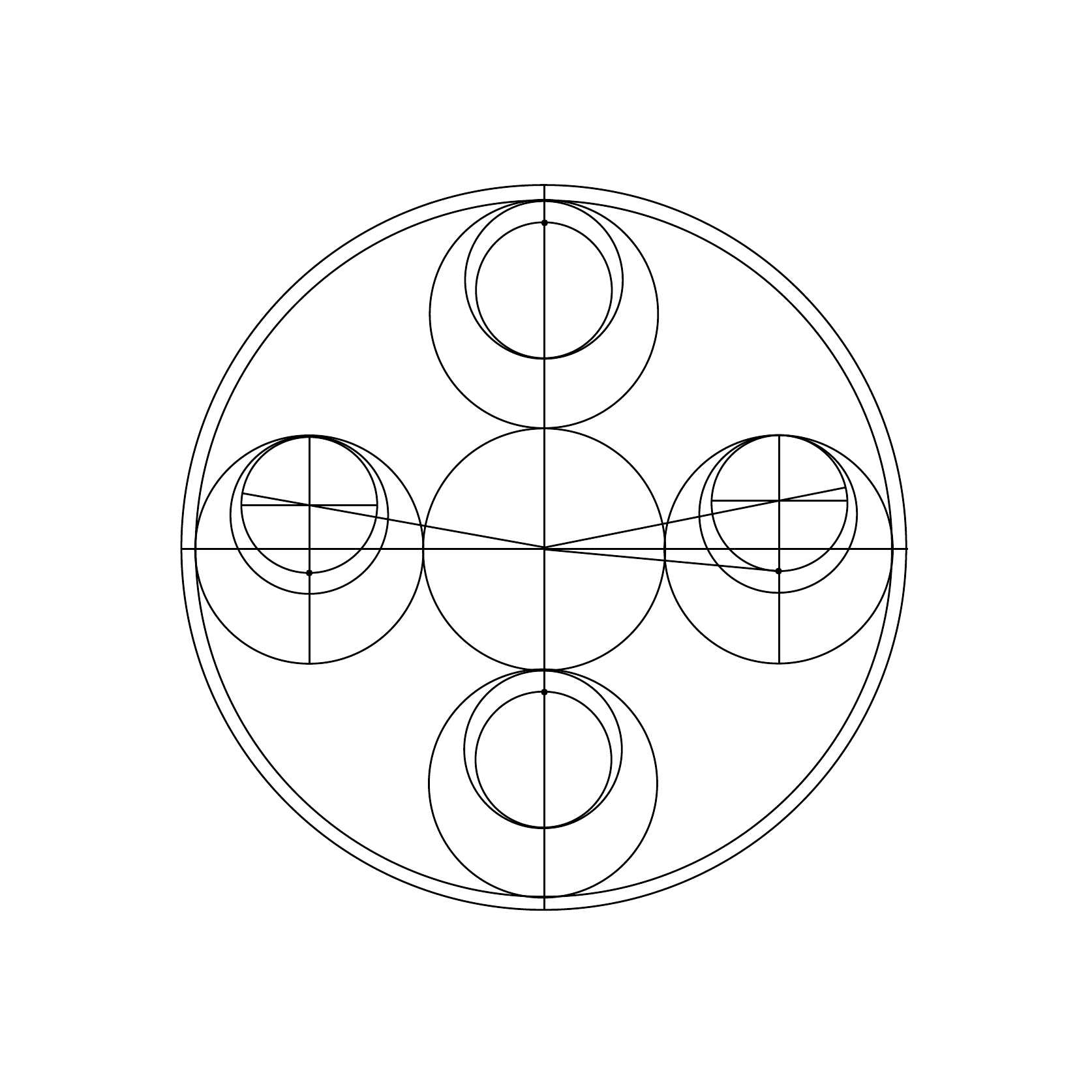
\end{minipage}
Configuration des orbes solides de Vénus~: sphères entières
représentées dans le plan à l'Apogée, au périgée et aux élongations
intermédiaires
\end{figure}

\newpage\phantomsection
\index{BBBHBE@\RL{qwm}!AJBBBHBJBE@\RL{taqwIm}, 1) calcul de la longitude vraie, 2) longitude vraie d'un astre}
\index{ARBGAQ@\RL{zhr}!ARBGAQAI@\RL{zuharaT}, Vénus}
\addcontentsline{toc}{chapter}{I.20 Régulation des mouvements de Vénus}
\includepdf[pages=83,pagecommand={\thispagestyle{plain}}]{edit.pdf}\phantomsection
\begin{center}
  \Large Chapitre vingt

  \large La régulation des mouvements de Vénus
\end{center}
Leur compte à midi du premier jour de l'année sept cent un de l'ère de
Yazdgard est de 9 signes et $10;9,0$ degrés pour Vénus moyen (comme le
Soleil moyen), et 2 signes $17;52,0$ degrés pour l'Apogée. Le
mouvement de Vénus moyen est comme le mouvement du Soleil moyen, soit
en vingt ans 11 signes et $25;13,20$ degrés, en une année persane 11
signes et $29;45,40,0$ degrés, en un mois persan 0 signe et
$29;34,9,51,46,50,57,32$ degrés, en un jour
$0;59,8,19,43,33,41,55,4,6$, et en une heure $0;2,27,50,50$. \`A la
date mentionnée, Vénus propre a atteint 10 signes et $20;50,19$
degrés~; son mouvement en vingt années persanes est 6 signes et
$0;36,0$ degrés, en une année persane 7 signes et $15;1,48,41$, en un
mois 0 signe et $18;29,44,13,9$ degrés, en un jour
$0;36,59,28,26,13,4,16$, et en une heure 0 signe et $0;1,32,28,41,6$
degrés. La correction [de Vénus] est faite selon la méthode de Saturne
et des autres astres errants.

\newpage\phantomsection
\index{AYAWAGAQAO@\RL{`u.tArid}, Mercure}
\index{BABDBC@\RL{flk}!BABDBCBEBEAKAKBD@\RL{falak muma_t_tal}, parécliptique}
\index{BABDBC@\RL{flk}!BABDBCAMAGBEBD@\RL{falak .hAmil}, orbe déférent}
\index{BABDBC@\RL{flk}!BABDBCBEAOBJAQ@\RL{falak mdIr}, orbe rotateur}
\index{AQABAS@\RL{ra'asa}!AQABAS@\RL{ra's}, tête, n{\oe}ud ascendant}
\index{APBFAH@\RL{_dnb}!APBFAH@\RL{_dnb}, queue, n{\oe}ud descendant}
\index{AMAQBC@\RL{.hrk}!AMAQBCAI BEAQBCAR@\RL{.harakaT al-markaz}, mouvement du centre (en général, le centre d'un épicycle)}
\index{AMAQBC@\RL{.hrk}!AMAQBCAI BEAQBCAR ATBEAS@\RL{.harakaT markaz al-^sams}, mouvement du centre du Soleil}
\index{AQAUAO@\RL{r.sd}!AEAQAUAGAO AUAMBJAMAI@\RL{al-'ar.sAd al-.sa.hI.haT}, les observations}
\index{AYAQAV@\RL{`r.d}!AYAQAV@\RL{`r.d}, latitude, \emph{i. e.} par rapport à l'écliptique}
\index{BABDBC@\RL{flk}!BABDBC BEAMBJAW@\RL{falak mu.hI.t}, orbe englobant}
\index{BABDBC@\RL{flk}!BABDBC AMAGBAAX@\RL{falak .hAfi.z}, orbe protecteur}

\addcontentsline{toc}{chapter}{I.21 Configuration des orbes de
  Mercure selon notre méthode confirmée par l'observation}
\includepdf[pages=84,pagecommand={\thispagestyle{plain}}]{edit.pdf}\phantomsection
\begin{center}
  \Large Chapitre vingt-et-un

  \large La configuration des orbes de Mercure selon notre méthode
  confirmée par l'observation.
\end{center}
Nous imaginons~:

-- un orbe dans le plan de l'écliptique, sur ses pôles et
en son centre, appelé le parécliptique.

-- un deuxième orbe dans un plan incliné par rapport au
plan du parécliptique d'un demi et un quart de degré à l'Apogée vers
le Sud, mais cette inclinaison n'est pas constante. Selon une autre
doctrine l'orbe est incliné d'un sixième de part et son inclinaison
est constante~: cette doctrine est plus juste. Le plan de l'orbe
incliné coupe le plan du parécliptique en deux points opposés appelés
la tête et la queue.

-- un troisième orbe dont le centre est sur la ceinture de
l'orbe incliné, de rayon quatre parts et cinq minutes (en parts telles
que le rayon de l'orbe incliné en compte soixante), appelé le
\emph{déférent}.

-- un quatrième orbe dont le centre est sur la ceinture du déférent,
de rayon un demi et un tiers de degré\footnote{Le rayon du rotateur
  serait donc $0;50$. Mais les calculs des rayons des orbes solides
  donnés plus loin laissent penser qu'il doit être $0;55$,
  c'est-à-dire une demi, un tiers et un demi-sixième de degré.},
appelé le \emph{rotateur}.

-- un cinquième orbe dont le centre est sur la ceinture du
rotateur, de rayon vingt-deux parts et quarante-six minutes (en les
mêmes parts), appelé orbe de l'épicycle.

-- un sixième orbe dont le centre est sur la ceinture de
l'épicycle, de rayon trente-trois minutes, appelé l'\emph{englobant}.

-- un septième orbe dont le centre est sur l'englobant, de rayon
semblable au rayon de l'englobant (trente-trois minutes), appelé le
\emph{protecteur}.

-- Mercure est centré sur la ceinture de ce [dernier] orbe.

Passons aux mouvements. Le parécliptique se meut sur
les pôles de l'écliptique, dans le sens des signes, d'un degré tous
les soixante ans, comme le mouvement des Apogées.

L'orbe incliné se meut dans le sens des signes, comme le mouvement du
centre de Mercure, c'est-à-dire comme le centre du Soleil~: en un jour
et une nuit, $0;59,8,10$.

Le déférent se meut en sens inverse des signes dans sa partie
supérieure, aussi d'autant que le mouvement du centre de Mercure.

L'orbe de l'épicycle se meut dans le sens des signes dans sa partie
supérieure, d'autant que l'excédent du mouvement propre de Mercure sur
le mouvement de son centre~: en un jour et une nuit,
$2;7,16,0,9,51,39,56,39$.

\newpage\phantomsection
\index{BBAWAQ@\RL{q.tr}!BFAUBA BBAWAQ BD-AJAOBHBJAQ BD-BEAQAEBJBJ@\RL{n.sf q.tr al-tadwIr al-mar'iyy}, rayon de l'épicycle apparent}
\index{BABDBC@\RL{flk}!BABDBC BEAMBJAW@\RL{falak mu.hI.t}, orbe englobant}
\index{BABDBC@\RL{flk}!BABDBCATAGBEBD@\RL{falak ^sAmil}, orbe total}
\index{AMAQBC@\RL{.hrk}!AMAQBCAI AHASBJAWAI BEAQBCBCAHAI@\RL{.harakaT basI.taT mrkkbaT}, mouvement simple-composé}
\index{AMAQBC@\RL{.hrk}!AMAQBCAI ANAGAUAI@\RL{.harakaT _hA.saT}, mouvement propre}
\index{APAQBH@\RL{_drw}!APAQBHAI@\RL{_dirwaT}, sommet, apogée}
\index{AHAOAB@\RL{bda'}!BEAHAJAOAB@\RL{mubtada'}|see{\RL{mabda'}}}
\index{AHAOAB@\RL{bda'}!BEAHAOAB@\RL{mabda'}, principe, origine}
\includepdf[pages=85,pagecommand={\thispagestyle{plain}}]{edit.pdf}\phantomsection
\noindent C'est un mouvement simple. Quant au
mouvement propre de Mercure, c'est un mouvement simple-composé qui
mesure, ensemble, le mouvement de cet épicycle
($2;7,16,0,9,51,39,56,39$) et le mouvement du centre de Mercure
($0;59,8,10$)~; car ces deux mouvements vont dans le même
sens. L'astre s'éloigne de l'apogée, d'autant que la somme des deux
mouvements, 0 signe et $3;6,24,10,1,38,37,28,42$~: c'est le mouvement
propre de Mercure, un mouvement composé mais uniforme par rapport au
centre de l'épicycle.

Voici une explication complémentaire. Si l'orbe incliné se meut d'un
quart de cercle, le déférent d'un quart de cercle, et le rotateur d'un
demi-cercle, alors se déplace l'apogée (origine du mouvement
propre) d'un quart de cercle dans le sens des signes. Cependant, on
trouve par l'observation qu'elle se déplace, dans le sens des signes,
comme le mouvement propre de Mercure (0 signe $3;6,24,10$), donc le
mouvement de l'épicycle autour de son centre est, dans le sens des
signes, l'excédent de ce mouvement propre sur le mouvement du centre
(car ils vont dans le même sens). J'ai déjà expliqué cela.

L'englobant\footnote{Cet orbe qu'{\shatir} avait appelé
  \textit{mu\d{h}\={\i}\d{t}} est maintenant désigne du nom de
  \textit{\v{s}\=amil}. Nous avons traduit ces deux termes par
  <<~englobant~>>.} se meut dans le sens des signes dans sa partie
supérieure, comme le double du mouvement du centre de Mercure~: en un
jour $1;58,16,20$.

Le protecteur se meut en sens inverse des signes dans sa partie supérieure,
quatre fois comme le mouvement du centre de Mercure~: en un
jour $3;56,32,40$.

Toujours située sur la droite issue du centre de l'épicycle et passant
par le centre de l'englobant, Mercure tantôt se rapproche, tantôt
s'éloigne du centre de l'épicycle~; et elle ne quitte jamais cette
droite\footnote{Ces deux orbes (englobant et protecteur) constituent un
  \emph{couple de \d{T}\=us{\=\i}}.}. Quand le centre de l'épicycle
est à l'Apogée ou au périgée, Mercure est à distance minimale du
centre de son épicycle, distance appelée rayon de l'épicycle apparent
qui vaut alors vingt-et-une parts et deux tiers de part. Quand le
centre est à trois signes [de l'Apogée], Mercure est à distance
maximale du centre de l'épicycle, à vingt-trois degrés et
cinquante-deux minutes.

La distance maximale de Mercure au centre du Monde est
quatre-vingt-six et deux tiers. Sa distance
minimale est trente-trois et un tiers\footnote{Il faut comprendre
  ces distances extrémales comme étant les distances extrémales de Mercure
  au centre du Monde \emph{quand Mercure est à l'Apogée ou au périgée},
  \textit{i. e.} quand $\overline{\kappa}=0°$ ou $180°$. Le rayon de
  l'épicycle apparent est alors minimal.}. Sauf que
Mercure n'atteint pas la distance minimale de ses orbes solides, comme
nous l'avons expliqué précédemment.

\newpage\phantomsection
\index{ALASBE@\RL{jsm}!BCAQAGAJ BEALASBEAI AJAGBEBEAI@\RL{kraT mjsm, falak mjsm}, sphère solide, orbe solide}
\index{BHAUBD@\RL{w.sl}!BEAJAJAUBD AHAYAVBG AHAHAYAV@\RL{mutta.sil ba`.dh biba`.d}, contigus}
\includepdf[pages=86,pagecommand={\thispagestyle{plain}}]{edit.pdf}\phantomsection
Quant aux dimensions des orbes solides~: le rayon de la sphère du
déférent est $28;52$, le rayon de la sphère du rotateur est $24;47$,
le rayon de la sphère de l'épicycle est $23;52$, le rayon de la sphère
de l'englobant est $1;6$, le rayon de la sphère du protecteur
est $0;33$ (tout cela en parts telles que le rayon du parécliptique en
compte soixante).

\label{rapport_merc}La distance maximale dans le parécliptique est donc quatre-vingt-huit
parts et cinquante-deux minutes, avec en plus de cela l'épaisseur du
parécliptique~; admettons qu'elle atteint 89. La distance minimale de
ses orbes est $31;8$, et même moins, de [ce qu'il faut] pour que les
orbes soient contigus~; on la pose égale à $31;0$. Dieu est le plus savant.

\newpage\phantomsection
\includepdf[pages=87,pagecommand={\thispagestyle{plain}}]{edit.pdf}\phantomsection

\newpage\phantomsection
\begin{figure}[h!]
\centering
\noindent\hspace{-85mm}\hspace{.5\textwidth}\begin{minipage}{17cm}
\footnotesize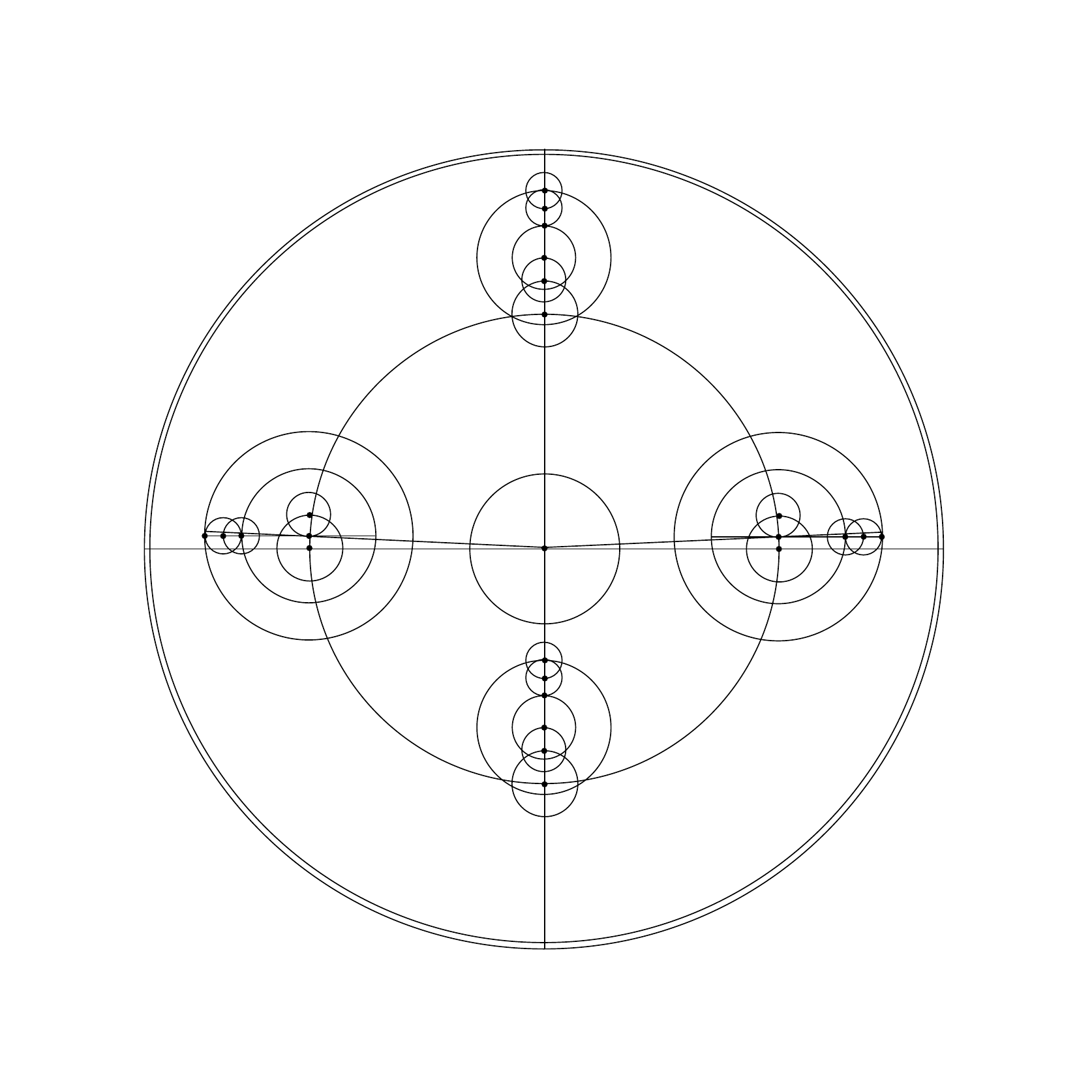
\end{minipage}
Configuration des orbes de Mercure représentées par les centres des
sphères solides dans la plan
\end{figure}

\newpage\phantomsection
\includepdf[pages=88,pagecommand={\thispagestyle{plain}}]{edit.pdf}\phantomsection

\newpage\phantomsection
\begin{figure}[h!]
\centering
\noindent\hspace{-85mm}\hspace{.5\textwidth}\begin{minipage}{17cm}
\footnotesize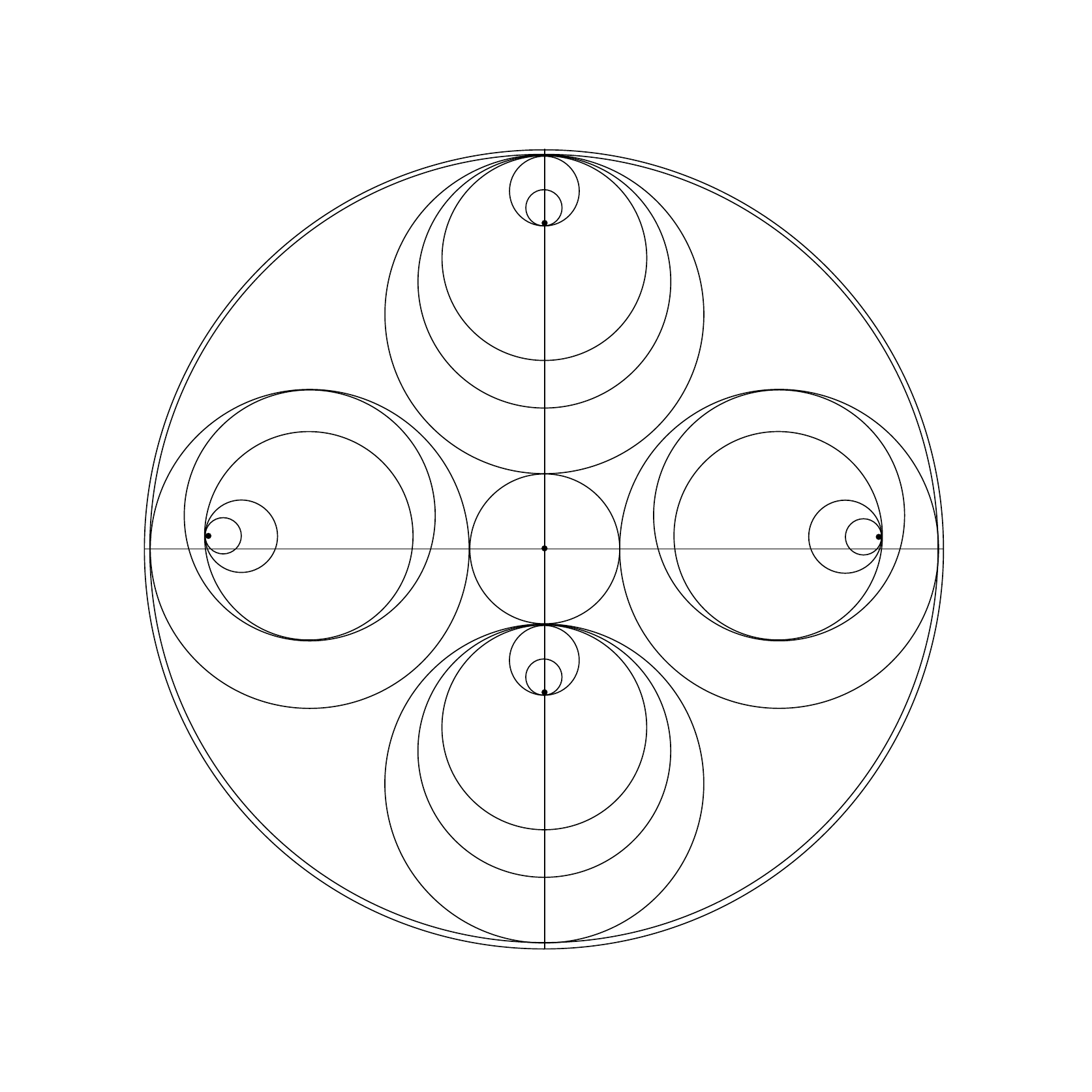
\end{minipage}
Configuration des orbes solides de Mercure, sphères entières
représentées dans le plan à l'Apogée, au périgée et aux élongations
intermédiaires
\end{figure}

\newpage\phantomsection
\index{AYAWAGAQAO@\RL{`u.tArid}, Mercure}
\index{ACAUBD@\RL{'a.sl}!ACAUBD@\RL{'a.sl}, 1) fondement, 2) origine}
\index{AMAQBC@\RL{.hrk}!AMAQBCAI BHASAW@\RL{.harakaT al-ws.t}, mouvement de l'astre moyen}
\index{AMAQBC@\RL{.hrk}!AMAQBCAI BHASAW ATBEAS@\RL{.harakaT ws.t al-^sams}, mouvement du Soleil moyen}
\index{BBBHBE@\RL{qwm}!AJBBBHBJBE@\RL{taqwIm}, 1) calcul de la longitude vraie, 2) longitude vraie d'un astre}
\addcontentsline{toc}{chapter}{I.22 Régulation des mouvements de Mercure}
\includepdf[pages=89,pagecommand={\thispagestyle{plain}}]{edit.pdf}\phantomsection
\begin{center}
  \Large Chapitre vingt-deux

  \large La régulation des mouvements de Mercure
\end{center}
Leurs conditions initiales à midi du premier jour de l'année sept
cent un de l'ère de Yazdgard sont de 9 signes et $10;9,0$ degrés pour
l'origine de Mercure
moyen (c'est comme l'origine du Soleil moyen et de Vénus moyen), et 7
signes et $2;52,0$ degrés pour l'Apogée. L'origine du mouvement propre
à cette date est de 5 signes $4;2,0$ degrés.

Le mouvement de Mercure moyen est comme le mouvement du Soleil moyen
(par an, 11 signes $29;45,40$ degrés). Le mouvement propre de Mercure
est\footnote{Les valeurs suivantes ont probablement été calculées à
  partir des valeurs observées par an ou en vingt années persanes. On
  obtient alors, par division, presque exactement la valeur donnée par
  jour $3;6,24,10,1,38,27,48,29$~: seul le huitième rang est
  erroné.}~: en vingt années persanes 11 signes $29;0,20,0$ degrés,
par an 1 signe $23;57,1$, par mois 3 signes $3;12,5,1$,
par jour 0 signe $3;6,24,10,1,38,37,29$, et par heure 0 signe $0;7,46,0,25$.

Si l'on retranche le mouvement de l'Apogée du mouvement moyen, on
obtient le mouvement du centre.

La [longitude] vraie de cet astre [se calcule] de la même manière que
Saturne vrai.

\newpage\phantomsection
\index{AYAOBD@\RL{`dl}!AJAYAOBJBD@\RL{ta`dIl}, équation}
\index{ALBJAH@\RL{jIb}!ALBJAH@\RL{jIb}, sinus}
\index{ALBJAH@\RL{jIb}!jIb tamAm@\RL{jIb al-tamAm}, cosinus}
\index{ALAOBD@\RL{jdl}!ALAOBHBD@\RL{jdwl}, table, tableau, catalogue}
\index{AQBCAR@\RL{rkz}!BEAQBCAR@\RL{markaz}, centre}
\index{AYAOBD@\RL{`dl}!ANAGAUAUAI BEAYAOAOBD@\RL{_hA.s.saT mu`addalaT}, astre propre corrigé}
\addcontentsline{toc}{chapter}{I.23 Calcul de l'équation des astres}
\includepdf[pages=90,pagecommand={\thispagestyle{plain}}]{edit.pdf}\phantomsection
\begin{center}
  \Large Chapitre vingt-trois

  \large Calcul de l'équation des astres
\end{center}
Pour calculer la \emph{première équation}\footnote{notée
  $c_1(\overline{\kappa})$ dans notre commentaire.}, pour Saturne,
Jupiter, Mars et Vénus, tu prends le sinus du centre et son cosinus,
et tu multiplies chacun par le rayon du déférent. Ce que tu obtiens à
partir du sinus, appelons-le \emph{base}~; ce que tu obtiens avec le
cosinus, appelons-le \emph{hauteur}. Puis multiplie le sinus, puis le
cosinus, chacun par le rayon du rotateur. Ce que tu obtiens à partir
du sinus, ajoute-le à la base si le centre est inférieur à trois
signes, et s'il est supérieur à trois signes alors
retranche-le\footnote{Ceci est erroné, il faut toujours ajouter
  quelque soit la valeur du centre (voir commentaire
  mathématique).}. Appelons \emph{dividende} leur somme ou
différence. Ce que tu obtiens à partir du cosinus, retranche-le de la
hauteur~; ajoute soixante à ce qui reste si le centre est inférieur à
trois signes, et s'il est supérieur à trois signes alors retranche ce
qui reste de soixante. Mets cela au carré, et ajoutes-y le carré du
dividende. Prends la racine carrée du total~: tu obtiens la distance
du centre de l'épicycle au centre du Monde. Divise le dividende que tu
as retenu par cela, il en sort le sinus de la première
équation. Trouve son arc dans la table des sinus, c'est la première
équation. \'Ecris-la en face du degré pour lequel tu l'as calculée,
c'est l'équation du centre et du mouvement propre.

Pour calculer la \emph{deuxième équation}\footnote{notée
  $c_2(0°,\alpha)$ dans notre commentaire.},
tu multiplies le sinus de l'astre propre corrigé
(puis son cosinus) par le rayon de l'orbe de l'épicycle de cet
astre. Ce que tu obtiens à partir du cosinus, tu l'ajoutes à la
distance du centre de l'épicycle au centre du Monde pour cet astre
($63;25$ pour Saturne, $62;45$ pour Jupiter, 66 pour Mars, $61;15$
pour Vénus, 65 pour Mercure) si l'astre propre corrigé est dans la
moitié supérieure de l'orbe de l'épicycle, et tu le retranches si
l'astre propre corrigé est dans sa moitié inférieure. Ce qu'on
obtient, ajoute son carré au carré du produit du sinus de l'astre
propre par le rayon de l'épicycle (ce produit est le
\emph{dividende})~; prends la racine carrée des deux carrés, et divise
le dividende par cela. On obtient le sinus de la deuxième équation
pour cet astre. Trouve son arc dans la table des sinus~; écris-le en
face de [la valeur] de l'astre propre corrigé pour laquelle tu l'as
calculé.

\newpage\phantomsection
\index{ANBDBA@\RL{_hlf}!ACANAJBDAGBA@\RL{i_htilAf}, irrégularité, anomalie, variation}
\index{BFASAH@\RL{nsb}!AOBBAGAEBB BFASAH@\RL{daqA'iq al-nisb}, coefficient d'interpolation}
\index{AYAOBD@\RL{`dl}!AJAYAOBJBD BEAMBCBE@\RL{ta`dIl mu.hakam}, équation principale}
\index{AHAYAO@\RL{b`d}!AHAYAO@\RL{b`d}, 1)~distance, 2)~élongation (\textit{i.e.}, écart en longitude)}
\includepdf[pages=91,pagecommand={\thispagestyle{plain}}]{edit.pdf}\phantomsection
\emph{Remarque}. La racine carrée des deux carrés est la distance de
l'astre au centre
du Monde à cette date, en parts telles que le rayon de l'orbe incliné
en compte soixante, et le rayon de l'épicycle, $6;30$ pour Saturne,
$11;30$ pour Jupiter, $39;30$ pour Mars, $43;33$ pour Vénus, et
$21;40$ pour Mercure (pour ce dernier, c'est le rayon de l'épicycle
apparent quand la distance est maximale ou bien
minimale\footnote{c'est-à-dire à l'Apogée ou au périgée, en lesquels
  le couple de \d{T}\=us{\=\i } rend l'épicycle apparent moindre que
  le vrai.}).

\label{troisieme_equation}La \emph{troisième équation}\footnote{notée
  $c_2(180°,\alpha)-c_2(0°,\alpha)$ dans notre commentaire.} est
la variation à distance minimale. Pour la calculer, tu fais
comme pour la deuxième équation, mais tu suppose le centre de
l'épicycle à distance minimale, c'est-à-dire $56;35$ pour Saturne,
$57;15$ pour Jupiter, 54 pour Mars, $58;45$ pour Vénus, et 55 pour
Mercure. Tu calcules ensuite selon la méthode pour la deuxième
équation. Tu retranches du résultat la deuxième équation pour le même
degré, et tu écris la différence en face du degré pour lequel tu as
fait ce calcul~: c'est la variation à distance minimale, et il faut
toujours l'ajouter à la deuxième équation (comme pour la Lune, j'ai
déjà indiqué cela).

La méthode pour calculer le \emph{coefficient d'interpolation}
\footnote{noté $\chi(\overline{\kappa})$ dans notre commentaire.} est
comme je t'ai enseigné pour le calcul du coefficient
d'interpolation pour la Lune. J'ai conçu les équations de sorte à ce que
la variation s'ajoute [toujours additivement] à la deuxième équation
(c'est plus facile ainsi). La méthode
pour calculer le coefficient d'interpolation est la suivante. Tu
calcules la distance du centre de l'épicycle au centre du Monde comme
je t'ai enseigné, puis tu calcules le rayon de l'épicycle
apparent\footnote{Pour toutes les planètes sauf Mercure, il s'agit
  simplement du rayon constant de l'épicycle.}. Tu divises celui-ci
par la distance du centre de l'épicycle que tu as calculée. Tu trouves
l'arcsinus du résultat dans la table des sinus~: c'est l'équation
principale maximale pour ce degré. Tu en soustrais l'équation
maximale à distance maximale pour cet astre. Tu divises la différence
par la différence entre l'équation à distance maximale et l'équation à
distance minimale. Le résultat est le coefficient d'interpolation.

\newpage\phantomsection
\index{ALAOBD@\RL{jdl}!ALAOBHBD@\RL{jdwl}, table, tableau, catalogue}
\includepdf[pages=92,pagecommand={\thispagestyle{plain}}]{edit.pdf}\phantomsection

Il faut bien conduire ce calcul, c'est subtil. Si tu préfères calculer
les équations de la manière la plus courante\footnote{Cette manière
  ``plus courante'' conduirait aux tables auxquelles le chapitre 14 fait
  implicitement référence. C'est ainsi que faisait Ptolémée.},
alors tu calcules quand
l'épicycle est à distance moyenne, et tu ajoutes à cela la variation à
distance minimale, ou bien tu en retranches la variation à distance
maximale. Je t'ai aussi expliqué dans le chapitre sur le calcul des
équations de la Lune que l'on peut calculer une table unique pour
corriger la Lune. De même pour chaque astre~: on peut faire une
équation unique pour corriger chacun sans équation du centre ni
équation de [l'astre] propre.

\newpage\phantomsection
\index{AYBBAO@\RL{`qd}!AYBBAOAI@\RL{`qdaT}, n{\oe}ud}
\index{ACBHAL@\RL{'awj}!ACBHAL@\RL{'awj}, Apogée}
\index{AYAQAV@\RL{`r.d}!AYAQAV@\RL{`r.d}, latitude, \emph{i. e.} par rapport à l'écliptique}
\index{BEBJBD@\RL{myl}!BEBJBD@\RL{myl}, inclinaison, inclinaison maximale (\textit{i. e.} mesure d'un angle dièdre)}
\index{ASAWAM@\RL{s.t.h}!ASAWAM BEASAJBH@\RL{\vocalize s.t.h mstwiN}, plan}
\index{APAQBH@\RL{_drw}!APAQBHAI@\RL{_dirwaT}, sommet, apogée}
\index{AMAVAV@\RL{.h.d.d}!AMAVBJAV@\RL{.ha.dI.d}, périgée}
\index{AQBJAN@\RL{ry_h}!BEAQAQBJAN@\RL{marrI_h}, Mars}
\index{ARAMBD@\RL{z.hl}!ARAMBD@\RL{zu.hal}, Saturne}
\index{ATAQBJ@\RL{^sry}!BEATAJAQBJ@\RL{al-mu^starI}, Jupiter}
\index{BABDBC@\RL{flk}!BABDBC BEAGAEBD@\RL{falak mA'il}, orbe incliné}
\index{BABDBC@\RL{flk}!BABDBC AJAOBHBJAQ@\RL{falak al-tadwIr}, orbe de l'épicycle}
\index{BBAWAQ@\RL{q.tr}!BBAWAQ@\RL{qi.tr}, diamètre}
\addcontentsline{toc}{chapter}{I.24 Latitudes des trois planètes supérieures~: Saturne, Jupiter et Mars}
\label{lat_debut}
\includepdf[pages=93,pagecommand={\thispagestyle{plain}}]{edit.pdf}\phantomsection
\begin{center}
  \Large Chapitre vingt-quatre

  \large Les latitudes des trois planètes supérieures~: Saturne,
  Jupiter et Mars
\end{center}
Dans la configuration de chacune d'elles, on a vu que l'inclinaison
maximale\footnote{Pour décrire la mesure de l'angle entre deux plans,
  on utilisait habituellement l'expression ``inclinaison maximale''
  qui désigne en effet l'arc maximal d'un cercle orthogonal à l'un des
  deux plans, de centre situé sur la droite commune (le cercle est alors
  orthogonal aux deux plans). Mais {\shatir} oubliera souvent le mot ``maximale'' pour parler simplement d'``inclinaison''.} de
l'orbe incliné par rapport au parécliptique est~: deux parts et demi
pour Saturne, une part et demi pour Jupiter, et une part pour Mars.

\`A présent, j'affirme qu'on a aussi trouvé ceci par l'observation~:
quand le centre de l'épicycle est à la moitié de l'arc d'orbe incliné
situé entre les n{\oe}uds (\textit{i. e.} là où l'inclinaison de l'orbe
incliné par rapport au parécliptique est maximale, vers le Nord),
l'inclinaison maximale du plan de l'épicycle par rapport à l'orbe incliné
vaut quatre parts et demi pour Saturne, deux parts et demi
pour Jupiter, et deux parts un quart pour Mars. Par l'observation, on a
trouvé qu'alors
le rayon passant par l'Apogée et le périgée de l'épicycle avait son
extrémité à l'Apogée comprise entre l'orbe incliné et le
parécliptique, et son extrémité au périgée incliné du même côté que
l'orbe incliné. Le rayon perpendiculaire au rayon mentionné est alors
parallèle au plan de l'écliptique.

Puis on a observé quand le centre de l'épicycle est situé aux
n{\oe}uds~: on a trouvé que les plans des trois cercles\footnote{Trois
  cercles~: le déférent, du rotateur et de l'épicycle.} étaient dans
le plan du parécliptique, c'est-à-dire le plan de l'écliptique. En effet,
on a trouvé que la latitude des astres était nulle quand le centre de
l'épicycle est à la tête ou bien à la queue, et que l'astre est en un
point quelconque de l'orbe de l'épicycle.

J'en déduis que les diamètres passant par les apogées ne sont
pas constamment situés dans les plans des orbes inclinés, ni dans les
plans des paracléptiques, sauf quand les centres des épicycles sont en
l'un des n{\oe}uds. Après [un passage à la tête] les apogées des
planètes supérieures s'inclinent toujours du côté du plan de
l'écliptique, et leurs périgées du côté opposé. Leur inclinaison
augmente jusqu'à un maximum atteint au milieu de l'arc compris entre
les deux n{\oe}uds, puis elle décroît jusqu'à disparaître au deuxième
n{\oe}ud. Puis l'apogée commence à s'incliner vers le plan de
l'écliptique, et le périgée du côté opposé~; cette inclinaison
augmente jusqu'à atteindre son maximum au milieu de l'arc compris
entre les deux n{\oe}uds, puis elle diminue jusqu'à disparaître à la
tête~; et la chose recommence depuis le début.

\newpage\phantomsection
\index{AHAWBDBEBJBHAS@\RL{b.talimayUs}, Ptolémée}
\index{AEAHAQANAS@\RL{'ibr_hs}, Hipparque}
\index{AQBJAN@\RL{ry_h}!BEAQAQBJAN@\RL{marrI_h}, Mars}
\index{ARAMBD@\RL{z.hl}!ARAMBD@\RL{zu.hal}, Saturne}
\index{ATAQBJ@\RL{^sry}!BEATAJAQBJ@\RL{al-mu^starI}, Jupiter}
\index{BFAWBB@\RL{n.tq}!BEBFAWBBAI@\RL{min.taqaT}, ceinture}
\index{AQABAS@\RL{ra'asa}!AQABAS@\RL{ra's}, tête, n{\oe}ud ascendant}
\index{AYAQAV@\RL{`r.d}!AYAQAV@\RL{`r.d}, latitude, \emph{i. e.} par rapport à l'écliptique}
\index{AWBHBD@\RL{.twl}!AWBHBD@\RL{.tUl}, longitude (par rapport à l'écliptique)}
\index{ACAUBD@\RL{'a.sl}!ACAUBD@\RL{'a.sl}, 1) fondement, 2) origine}
\includepdf[pages=94,pagecommand={\thispagestyle{plain}}]{edit.pdf}\phantomsection

L'apogée, après s'être éloigné des n{\oe}uds, reste toujours entre
les deux ceintures~; pour le périgée c'est le contraire. Cela a pour
effet de diminuer les latitudes des apogées et d'augmenter les
latitudes des périgées. Si la distance du centre de l'épicycle au
centre de l'écliptique varie, alors l'inclinaison de l'apogée et du
périgée par rapport au centre de l'écliptique varie aussi~: si la
distance du centre de l'épicycle augmente, alors l'inclinaison de
l'apogée décroît, et inversement.

\label{noeud_sup1}
Le n{\oe}ud ascendant de Saturne au premier jour de l'année sept cent
neuf du calendrier persan est à cinq parts du Lion, le n{\oe}ud
ascendant de Jupiter à vingt-et-une parts du Cancer, et le n{\oe}ud
ascendant de Mars à dix-huit parts du Taureau. Pour chaque astre, le
n{\oe}ud descendant est à l'opposé du n{\oe}ud ascendant. Un quart de
cercle après le n{\oe}ud ascendant se situe la latitude maximale vers
le Nord~; c'est pour Saturne $2;4$ quand il est à l'apogée de
l'épicycle et $3;3$ quand il est au périgée\footnote{Ces valeurs et
  les suivantes sont tirées des tables de latitudes de
  l'\textit{Almageste} (chapitre XIII).  Les données de l'observation
  rapportées par Ptolémée sont respectivement~: 2°, 3°, 1°, 2°, --,
  $4;20$°, 2°, 3°, 1°, 2°, --, 7° (\textit{cf.} \cite{swerdlow2005}
  p.~45-46).} pour Jupiter $1;6$ quand il est à l'apogée et $2;4$
quand il est au périgée, et pour Mars $0;7$ quand il est à l'apogée et
$4;21$ quand il est au périgée. Un quart de cercle après le n{\oe}ud
descendant se situe la latitude maximale vers le Sud~; c'est pour
Saturne $2;2$ quand il est à l'apogée de l'épicycle et $3;5$ quand il
est en son périgée, pour Jupiter $1;5$ quand il est à l'apogée de
l'épicycle et $2;8$ quand il est au périgée, et pour Mars $0;3$ quand
il est à l'apogée de l'épicycle et $7;7$ quand il est en son périgée.

\begin{center}
  \large Section
\end{center}

Sache que, concernant le problème des latitudes, ces fondements sont
des fondements nécessaires trouvés par l'observation, certains ayant
aussi été confirmés par le calcul~; on ne peut aucunement s'en
passer. Ni Ptolémée ni Hipparque, ni aucun autre des Anciens, ni aucun
des Modernes jusqu'à nos jours, n'a pu faire cela~: donner à
l'astronomie des fondements simples et vrais qui suffisent à la fois
pour les mouvements en longitude et pour les mouvements en
latitude. Ptolémée n'a pas pu imaginer de fondements suffisants pour
les mouvements en latitude~; que l'on ait ou pas d'indulgence pour la
théorie classique des mouvements en longitude, on ne peut avoir aucune
indulgence pour les mouvements en latitude. Le problème est le suivant.

\newpage\phantomsection
\index{AHAWBDBEBJBHAS@\RL{b.talimayUs}, Ptolémée}
\index{AEAHAQANAS@\RL{'ibr_hs}, Hipparque}
\index{AWBHASBJ@\RL{n.sIr al-dIn al-.tUsI}, Na\d{s}{\=\i}r al-D{\=\i}n al-\d{T}\=us{\=\i}}
\index{ATBJAQAGARBJ@\RL{q.tb al-dIn al-^sIrAzI}, Qu\d{t}b al-D{\=\i}n al-\v{S}{\=\i}r\=az{\=\i}}
\index{AYAQAVBJ@\RL{al-mwyd al-`r.dI}, al-Mu'ayyad al-`Ur\d{d}{\=\i}}
\index{BGBJAKBE@\RL{ibn al-hay_tam}, Ibn al-Haytham}
\index{AHAOAB@\RL{bda'}!BEAHAOAB@\RL{mabda'}, principe, origine}
\includepdf[pages=95,pagecommand={\thispagestyle{plain}}]{edit.pdf}\phantomsection

Quand le centre de l'épicycle est au milieu de l'arc compris
entre les n{\oe}uds, que le plan de l'épicycle est incliné par rapport
au plan de l'orbe incliné (comme je l'ai dit et comme on l'a trouvé
par l'observation) et que l'orbe incliné tourne en direction du
n{\oe}ud, alors il faut que l'épicycle ne cesse d'être incliné par
rapport au plan de l'orbe incliné, y compris au n{\oe}ud. Nous avons
certes vérifié par l'observation que le plan de l'épicycle n'est pas
incliné par rapport au plan de l'écliptique et qu'eux deux ont une
latitude nulle. Inversement, après que l'épicycle a été dans le plan
de l'écliptique et quand il parcourt l'arc compris entre les
n{\oe}uds, il présente une forte inclinaison.

\label{critique_lat1}
Le rapprochement de la ceinture de l'épicycle et de la ceinture de
l'orbe incliné, sans un mobile qui ne perturbe les mouvements en
longitude, dans l'astronomie classique de l'\textit{Almageste} et des
\emph{Hypothèses} concernant les corps mûs en latitude, est
impossible. Ce dont se sont efforcés Ptolémée et ses successeurs comme
Ibn al-Haytham dans son \textit{\'Epître}, Na\d{s}{\=\i}r
al-\d{T}\=us{\=\i} dans sa \textit{Tadhkira}, Mu'ayyad
al-`Ur\d{d}{\=\i} dans son \textit{Astronomie}, Qu\d{t}b
al-Sh{\=\i}r\=az{\=\i} dans sa \uwave{\textit{Muntaha 'adw\=ar}},
et ce que quiconque a repris de leurs propos dans son livre (sans
compter l'extension qu'on leur a donné, et la gloire qu'on s'est faite
du principe \textit{ibd\=a`{\=\i}}), rien de tout cela ne suffit au but
visé pour la longitude et la latitude, et quand cela convient pour
l'une, cela pêche pour l'autre. Que Dieu pardonne à ceux qui ont ici
admis leur faiblesse. \`A la fin de ce qu'a dit Qu\d{t}b al-D{\=\i}n,
il a dit (que Dieu aide celui qui examine son livre) qu'il avait trouvé
une manière parfaite de tout définir et d'éliminer les défauts qui restent
dans ce que nous avons mentionné (c'est Dieu qui inspire les choses
correctes)~; mais Ptolémée avait déjà dit des choses semblables~; la
connaissance appartient à ceux qui la mentionne [en premier].

Sache que la raison empêchant de se débarrasser des doutes dans
l'astronomie connue est que le principe de cette astronomie est faux~;
il est donc nécessaire de changer complètement ce principe pour l'amener
à une position sauve du doute et qui satisfasse au problème posé,
si c'est possible, comme nous l'avons fait.

Dès lors que nous avons conçu la configuration des orbes de la manière
sus-dite, que nous en sommes venus aux latitudes de ces astres, et que
nous avons rapporté ce qu'on a trouvé par l'observation concernant
leurs latitudes et leurs inclinaisons, nous avons trouvé que ceci
était en accord avec notre conception des orbes et de leurs
mouvements ; ceci se montre par ce que je vais expliquer ici.

\newpage\phantomsection
\index{AWBHASBJ@\RL{n.sIr al-dIn al-.tUsI}, Na\d{s}{\=\i}r al-D{\=\i}n al-\d{T}\=us{\=\i}}
\index{ATBJAQAGARBJ@\RL{q.tb al-dIn al-^sIrAzI}, Qu\d{t}b al-D{\=\i}n al-\v{S}{\=\i}r\=az{\=\i}}
\index{AYAQAVBJ@\RL{al-mwyd al-`r.dI}, al-Mu'ayyad al-`Ur\d{d}{\=\i}}
\index{AQBJAN@\RL{ry_h}!BEAQAQBJAN@\RL{marrI_h}, Mars}
\index{ARAMBD@\RL{z.hl}!ARAMBD@\RL{zu.hal}, Saturne}
\index{ATAQBJ@\RL{^sry}!BEATAJAQBJ@\RL{al-mu^starI}, Jupiter}
\index{BABDBC@\RL{flk}!BABDBC AJAOBHBJAQ@\RL{falak al-tadwIr}, orbe de l'épicycle}
\index{BABDBC@\RL{flk}!BABDBC BEAGAEBD@\RL{falak mA'il}, orbe incliné}
\index{BABDBC@\RL{flk}!BABDBCBEAOBJAQ@\RL{falak mdIr}, orbe rotateur}
\index{BABDBC@\RL{flk}!BABDBCAMAGBEBD@\RL{falak .hAmil}, orbe déférent}
\index{BEBJBD@\RL{myl}!BEBJBD@\RL{myl}, inclinaison, inclinaison maximale (\textit{i. e.} mesure d'un angle dièdre)}
\index{AQAUAO@\RL{r.sd}!AEAQAUAGAO AUAMBJAMAI@\RL{al-'ar.sAd al-.sa.hI.haT}, les observations}
\includepdf[pages=96,pagecommand={\thispagestyle{plain}}]{edit.pdf}\phantomsection

Supposons que, quand le centre du déférent et du rotateur sont au
milieu de l'arc compris entre les n{\oe}uds, le diamètre du déférent
est incliné par rapport à l'orbe incliné du côté où se trouve
l'inclinaison du périgée de l'épicycle et que son inclinaison est
chez Saturne trois parts et demi ; le rotateur de Saturne est incliné
par rapport au plan de son déférent d'une seule part. L'orbe déférent
de Mars est incliné par rapport au plan de l'orbe incliné d'un degré
et demi et un huitième, et l'épicycle de Mars est incliné par rapport
au plan de l'orbe déférent d'une demi-part et un huitième. Enfin, le
diamètre de l'orbe déférent de Jupiter est incliné par rapport au plan
de l'orbe incliné de deux parts, et le diamètre de l'épicycle de
Jupiter est incliné par rapport au plan du déférent de Jupiter d'une
demi-part. Faisons se mouvoir les orbes autour de leurs centres en
donnant à chacune son mouvement simple de la grandeur supposée dans le
sens supposé. Quand l'orbe incliné se meut d'un quart de cercle, le
déférent se meut d'un quart de cercle aussi, le rotateur se meut d'un
demi-cercle, et l'inclinaison change de côté: l'inclinaison de
l'épicycle passe du côté opposé à l'inclinaison du déférent et le plan
de l'épicycle vient dans le plan de l'écliptique\footnote{En fait,
  bien qu'ils deviennent parallèles, les deux plans ne sont jamais
  confondus. Ils le seraient si le rayon de la ceinture du déférent était
  nul, \textit{i. e.} si le centre du rotateur était confondu avec le
  centre du déférent. Comme le rayon du déférent est petit devant le
  rayon de l'orbe incliné, cette hypothèse peut se justifier~;
  {\shatir} l'a aussi adoptée dans le paragraphe suivant et sur la
  figure.}.
En effet on a posé
l'inclinaison du déférent égale à la somme de l'inclinaison du l'orbe
incliné et de l'inclinaison de l'épicycle. Or quand s'inverse
l'inclinaison de l'épicycle il faut soustraire l'inclinaison de
l'épicycle de l'inclinaison du déférent, il reste l'excédent de
l'inclinaison de l'orbe déférent sur l'inclinaison de l'épicycle, et
cet excédent est de la grandeur de l'inclinaison de l'orbe incliné; le
centre du déférent étant arrivé au n{\oe}ud, il a atteint le plan de
l'écliptique, et le plan de l'épicycle est donc dans le plan de
l'épicycle. Ceci est évident.

Prête attention et sache que la somme de l'inclinaison du déférent et
de l'inclinaison de l'épicycle a été posée égale à l'inclinaison,
observée, de l'épicycle par rapport au plan de l'orbe incliné, et que
la différence entre l'inclinaison du déférent et l'inclinaison de
l'épicycle a été posée de la grandeur de l'inclinaison de l'orbe
incliné. De là, tu sais aussi que, quand les deux inclinaisons
s'ajoutent, elles atteignent un maximum et c'est au milieu de l'arc
compris entre les n{\oe}uds; mais quand elles sont opposées, leur
différence est l'inclinaison de l'orbe incliné, on est aux n{\oe}uds
et l'épicycle est dans le plan de l'écliptique. Comme on a supposé
l'épicycle dans le plan de l'écliptique à cause de l'observation,
cette notion est donc juste.  

\newpage\phantomsection
\index{AYAOBD@\RL{`dl}!BEAQBCAR BEAYAOAOBD@\RL{markaz mu`addal}, centre corrigé}
\index{AQAUAO@\RL{r.sd}!AEAQAUAGAO AUAMBJAMAI@\RL{al-'ar.sAd al-.sa.hI.haT}, les observations}
\index{AWBHASBJ@\RL{n.sIr al-dIn al-.tUsI}, Na\d{s}{\=\i}r al-D{\=\i}n al-\d{T}\=us{\=\i}}
\index{ATBJAQAGARBJ@\RL{q.tb al-dIn al-^sIrAzI}, Qu\d{t}b al-D{\=\i}n al-\v{S}{\=\i}r\=az{\=\i}}
\index{AYAQAVBJ@\RL{al-mwyd al-`r.dI}, al-Mu'ayyad al-`Ur\d{d}{\=\i}}
\includepdf[pages=97,pagecommand={\thispagestyle{plain}}]{edit.pdf}\phantomsection

\noindent Pour démontrer cela\footnote{Pour
  transcrire les lettres \textit{\=a l\=a b j d h w z \d{h} \d{t}
    {\=\i} k l m n s ` f \d{s} q r \v{s}} utilisées dans la figure, on
  utilisera respectivement les lettres A A' B C D E F G H I J K L M N
  S O P \d{S} Q R \v{S}.}, menons la droite AB dans le plan de
l'écliptique et soit B le centre de l'écliptique. En soit issue la
droite BE inclinée par rapport à AB comme l'inclinaison de l'orbe
incliné par rapport au plan de l'écliptique. Traçons la droite CED
supposée diamètre de l'orbe déférent; du point E, menons aussi la
droite GI supposée diamètre de l'épicycle et formant un angle SEG
somme des deux inclinaisons quand elles sont du même côté (au milieu
de l'arc entre les n{\oe}uds). Quand le centre du déférent se déplace
vers le point B, le diamètre du déférent CD garde son inclinaison et
devient la droite Q\d{S}. Si le diamètre de l'épicycle gardait alors
son inclinaison, la droite GI deviendrait la droite OK; mais si le
rotateur tourne d'un demi-cercle, alors l'inclinaison du diamètre OK
par rapport au diamètre \d{S}Q passe de l'autre côté, l'angle QBK
devient l'angle A'BQ et le diamètre de l'épicycle (la droite OK)
coïncide avec le plan de l'écliptique. Il en est ainsi car la somme
des deux inclinaisons est comme l'inclinaison de l'épicycle, et
l'excédent de l'inclinaison du déférent sur l'inclinaison de
l'épicycle est comme l'inclinaison de l'orbe incliné. On a déjà répété
cela. Attention, quand on dit que le centre arrive aux n{\oe}uds, il
faut comprendre: le centre corrigé.

\begin{center}
  \large Remarque
\end{center}\label{critique_lat2}
Si nous avions posé l'inclinaison du déférent de la même grandeur que
l'inclinaison de l'épicycle, c'est-à-dire toutes deux de la grandeur
de l'inclinaison de l'épicycle trouvée par l'observation, alors après
rotations du déférent d'un quart de cercle et du rotateur d'un
demi-cercle le plan de l'épicycle serait dans le plan de l'orbe
incliné; c'est évident. Or l'observation confirme que c'est le plan de
l'écliptique et non le plan de l'orbe incliné. La plupart des Modernes
-- Na\d{s}{\=\i}r al-\d{T}\=us{\=\i}, Mu'ayyad al-`Ur\d{d}{\=\i},
Qu\d{t}b al-Sh{\=\i}r\=az{\=\i} -- ont pensé que la latitude
s'annulait quand c'est le plan de l'orbe incliné: c'est impossible
puisque le plan de l'épicycle reste incliné par rapport au plan de
l'écliptique. C'est très clair~; cependant les observateurs ont dit
qu'il est soit dans le plan de l'écliptique soit proche de celui-ci,
et ils n'ont pas vérifié avec certitude qu'il est exactement dans le
plan de l'écliptique.

\newpage\phantomsection
\index{AQBJAN@\RL{ry_h}!BEAQAQBJAN@\RL{marrI_h}, Mars}
\index{ARAMBD@\RL{z.hl}!ARAMBD@\RL{zu.hal}, Saturne}
\index{ATAQBJ@\RL{^sry}!BEATAJAQBJ@\RL{al-mu^starI}, Jupiter}
\index{AQAUAO@\RL{r.sd}!AEAQAUAGAO AUAMBJAMAI@\RL{al-'ar.sAd al-.sa.hI.haT}, les observations}
\index{AHAWBDBEBJBHAS@\RL{b.talimayUs}, Ptolémée}
\index{AYBBAO@\RL{`qd}!AYBBAOAI@\RL{`qdaT}, n{\oe}ud}
\index{ACBHAL@\RL{'awj}!ACBHAL@\RL{'awj}, Apogée}
\includepdf[pages=98,pagecommand={\thispagestyle{plain}}]{edit.pdf}\phantomsection

\begin{center}
  \large Deuxième remarque
\end{center}
Si l'on avait trouvé par l'observation que l'inclinaison de l'épicycle
\emph{par rapport à l'orbe incliné}\footnote{Nous avons mis en italique~:
  {\shatir} envisage ici un modèle où l'épicycle serait incliné par rapport
  à l'orbe incliné et non par rapport au déférent (les plans du déférent,
  du rotateur et de l'orbe incliné seraient confondus).} était de la
grandeur trouvée par
l'observation, alors après rotation du déférent d'un quart de cercle
et rotation, dans le même temps, du rotateur, d'un demi-cercle,
l'inclinaison de l'épicycle s'inverserait~; mais il passe alors dans le plan
de l'écliptique; or on a trouvé par l'observation que l'inclinaison
de l'épicycle est différente de l'inclinaison de l'orbe
incliné. \`A cause de cela, nous avons déduit les inclinaisons des
déférents et des épicycles de sorte que la somme des deux inclinaisons
soit comme l'inclinaison, observée, de l'épicycle, et que la
différence des deux soit comme l'inclinaison de l'orbe incliné.
Ce que nous avons accompli est correct.

\begin{center}
  \large Troisième remarque
\end{center}
Puisqu'on a découvert par l'observation que l'inclinaison maximale de
l'épicycle par rapport à l'orbe incliné est atteinte au milieu de
l'arc entre les n{\oe}uds, que là se trouve la dernière latitude Nord,
et que ce lieu est, chez Mars, le lieu de l'Apogée, alors les
diamètres du déférent et de l'épicycle passant par l'Apogée et le
périgée coïncident avec la droite issue du centre de l'écliptique et
passant par le lieu de la dernière latitude Nord (c'est-à-dire le
milieu de l'arc entre les n{\oe}uds): car ce lieu est le lieu de
l'Apogée chez Mars.

\label{noeud_sup2}
Puisque le lieu de la dernière latitude Nord qui est milieu de l'arc
entre les n{\oe}uds suit le lieu de l'Apogée chez Jupiter de vingt
parts d'après Ptolémée et vingt-huit parts d'après nos observations,
et qu'il précède le lieu de l'Apogée chez Saturne de cinquante parts,
alors le diamètre dans les dernières latitudes n'est pas celui passant
par l'Apogée du déférent et par l'Apogée du rotateur. Si nous
supposons que les orbes sont à l'Apogée, puis qu'ils se meuvent, le
déférent se mouvant de vingt parts, et le rotateur, du double, alors
l'Apogée du déférent se déplace de vingt parts dans le sens des signes
et devient le point D, et le diamètre EI devient celui où se situent
les dernières latitudes; le rotateur ayant tourné du double de cela,
son Apogée devient le point F tel que l'arc DE soit comme l'arc EF,
et le diamètre [du rotateur] coïncidant\footnote{Il faudrait plutôt
  dire ``parallèle à la droite EI''.} avec la droite EI se situe aux
dernières latitudes.

\newpage\phantomsection
\index{BBAWAQ@\RL{q.tr}!BBAWAQ AYAQAVBJ@\RL{qi.tr `r.diy}, diamètre des latitudes}
\includepdf[pages=99,pagecommand={\thispagestyle{plain}}]{edit.pdf}\phantomsection

On fait de même pour Saturne: posons que les
orbes se meuvent à l'envers et reviennent au lieu situé cinquante
parts avant l'Apogée de Saturne, lieu de ses dernières
latitudes. Parmi les diamètres du déférent et les diamètres du
rotateur, ceux qui coïncident avec la droite issue du centre de
l'écliptique dans la direction des dernières latitudes sont les deux
\emph{diamètres des latitudes} (j'entends par là, les diamètres en
lesquels sont les dernières latitudes). Voir la figure.


\begin{center}
  \large Remarque
\end{center}
Si l'on avait supposé que le déférent était dans le plan de l'orbe
incliné, que l'inclinaison du rotateur par rapport au plan du déférent
était comme l'inclinaison qu'on a supposée du déférent par rapport à
l'orbe incliné, et que l'inclinaison de l'épicycle par rapport au
rotateur était ce qu'elle est, alors le résultat de ces deux
hypothèses serait le même. Sache cela.

\begin{figure}[h!]
\centering\label{lat_inf_profil}
\footnotesize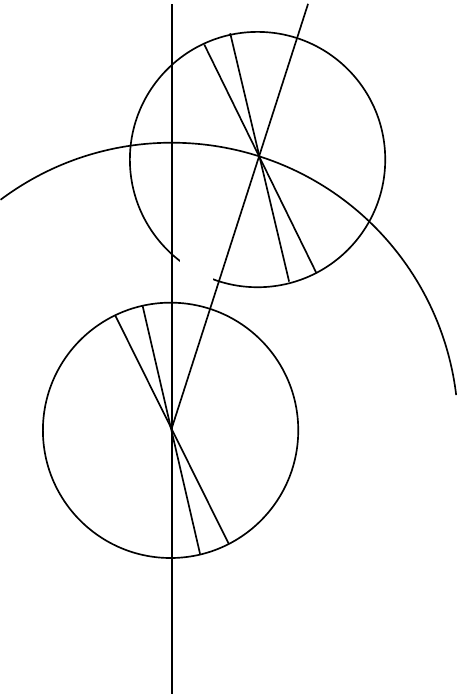
\end{figure}

\newpage\phantomsection
\includepdf[pages=100,pagecommand={\thispagestyle{plain}}]{edit.pdf}\phantomsection

\newpage\phantomsection
\begin{figure}[h!]
\centering
\noindent\hspace{-85mm}\hspace{.5\textwidth}\begin{minipage}{17cm}
\footnotesize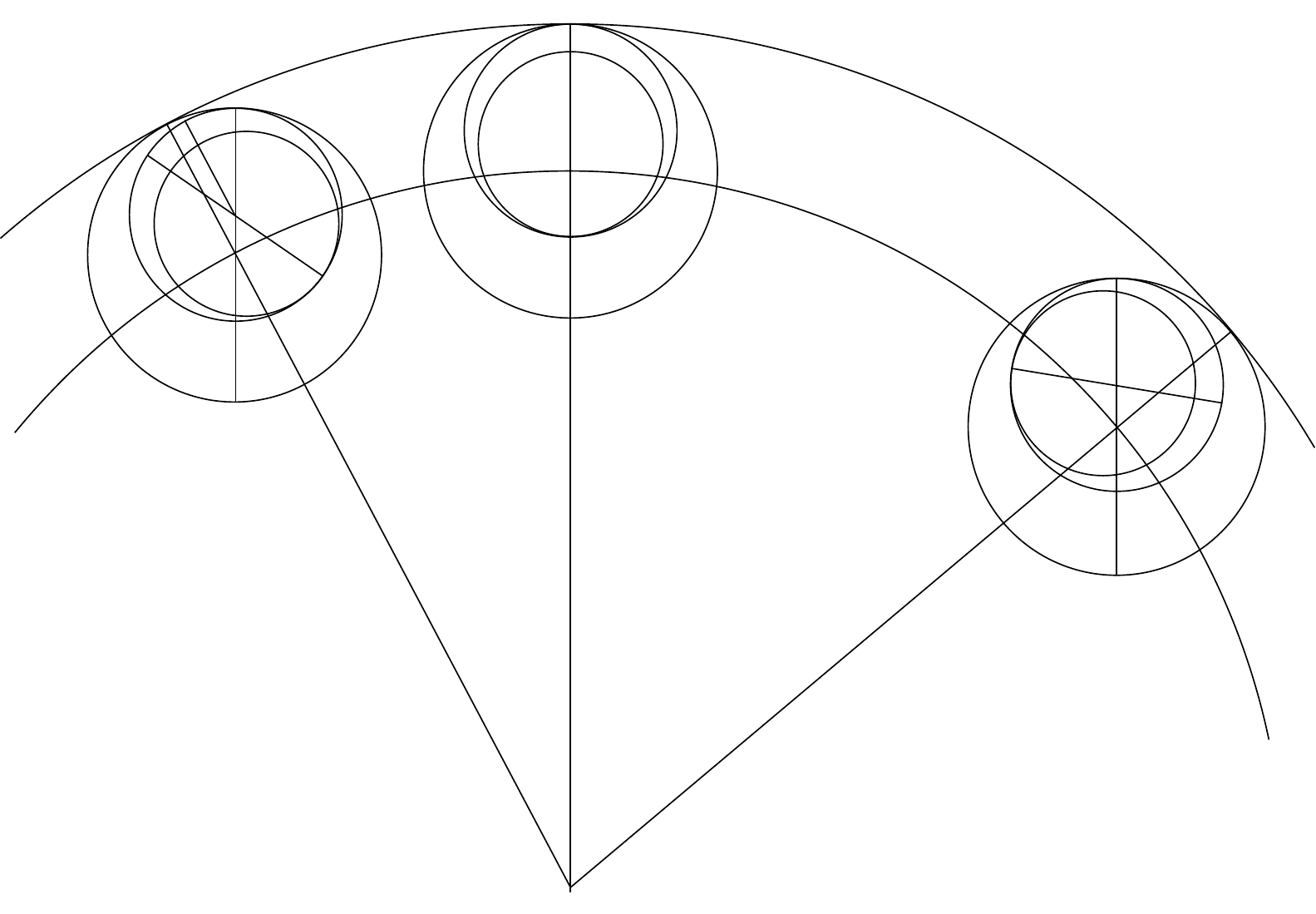
\end{minipage}
\end{figure}

\newpage\phantomsection
\index{AHAWBDBEBJBHAS@\RL{b.talimayUs}, Ptolémée}
\index{ARBGAQ@\RL{zhr}!ARBGAQAI@\RL{zuharaT}, Vénus}
\index{BABDBC@\RL{flk}!BABDBC BEAGAEBD@\RL{falak mA'il}, orbe incliné}
\index{AYBBAO@\RL{`qd}!AYBBAOAI@\RL{`qdaT}, n{\oe}ud}
\index{AYAQAV@\RL{`r.d}!AYAQAV@\RL{`r.d}, latitude, \emph{i. e.} par rapport à l'écliptique}
\index{ACBHAL@\RL{'awj}!ACBHAL@\RL{'awj}, Apogée}
\index{BBAWAQ@\RL{q.tr}!BBAWAQ@\RL{qi.tr}, diamètre}
\index{BEBJBD@\RL{myl}!BEBJBD BHAQAGAH@\RL{mIl al-wirAb}, inclinaison de biais}
\index{BEBJBD@\RL{myl}!BEBJBD AJAOBHBJAQ@\RL{mIl al-tadwIr}, inclinaison de l'épicycle}
\addcontentsline{toc}{chapter}{I.25 Mouvements de Vénus et Mercure en latitude, en deux parties}
\includepdf[pages=101,pagecommand={\thispagestyle{plain}}]{edit.pdf}\phantomsection
\begin{center}
  \Large Chapitre vingt-cinq

  \large Mouvements de Vénus et Mercure en latitude, en deux parties
\end{center}

\begin{center}
  \large Première partie

  \normalsize Les latitudes de Vénus
\end{center}
Chez Ptolémée, il y a deux théories à ce sujet.

\emph{Sa première théorie} est celle qu'il expose dans
l'\emph{Almageste} quand il dit que l'orbe incliné de Vénus tantôt
s'approche tantôt s'éloigne du plan de l'écliptique. Aux n{\oe}uds, le
centre de l'épicycle de Vénus est dans le plan de l'écliptique, mais
hors des n{\oe}uds, il est toujours au Nord du plan de l'écliptique
(au maximum, à un sixième de degré, inclinaison vers le Nord de
l'orbe incliné de Vénus quand elle est à l'Apogée). Quand le centre de
son épicycle est aux n{\oe}uds (c'est-à-dire aux quadratures de
l'Apogée), le diamètre perpendiculaire au diamètre passant par
l'Apogée de l'épicycle et par son périgée est dans le plan de
l'écliptique; quand il est à la queue, l'Apogée de l'épicycle est
incliné vers le Nord et le périgée vers le Sud; quand il est à la
tête, l'Apogée de l'épicycle est incliné vers le Sud et son périgée
vers le Nord; la grandeur de cette inclinaison est deux degrés et demi
dans l'\emph{Almageste}, et elle s'appelle \emph{inclinaison de
  l'épicycle}.  Quand le centre de l'épicycle est à l'Apogée, le
diamètre passant par l'Apogée de l'épicycle et par son périgée est
dans le plan de l'orbe incliné, et le diamètre perpendiculaire est
incliné par rapport au plan de l'écliptique, de sorte que la moitié
allant de l'Apogée de l'épicycle à son périgée est au Nord de l'orbe
incliné, et l'autre moitié, au Sud; la grandeur de cette inclinaison
par rapport au plan de l'écliptique est de trois degrés et demi, et
elle s'appelle \emph{inclinaison de biais}. C'est la théorie de
Ptolémée dans l'\emph{Almageste}. Certains parmi les Modernes l'ont
amendée par l'observation : ils ont trouvé que l'inclinaison de
l'épicycle était égale à l'inclinaison de biais. Les savants ont
calculé les tables de latitude de Vénus selon cette théorie.

Sa \emph{seconde théorie} est celle qu'il expose dans son livre
connu sous le nom des
\emph{Hypothèses} et qui vient après le premier: l'orbe incliné de
Vénus serait d'inclinaison constante de part et d'autre. Son
inclinaison maximale est un sixième de degré vers le Nord, et de même
vers le Sud. L'inclinaison de l'épicycle est, dans les deux cas, trois
degrés et demi.

\newpage\phantomsection
\index{AHAWBDBEBJBHAS@\RL{b.talimayUs}, Ptolémée}
\index{ANBDBA@\RL{_hlf}!ACANAJBDAGBA@\RL{i_htilAf}, irrégularité, anomalie, variation}
\index{BABDBC@\RL{flk}!BABDBCAMAGBEBD@\RL{falak .hAmil}, orbe déférent}
\index{BABDBC@\RL{flk}!BABDBCBEAOBJAQ@\RL{falak mdIr}, orbe rotateur}
\includepdf[pages=102,pagecommand={\thispagestyle{plain}}]{edit.pdf}\phantomsection

Ni Ptolémée ni aucun autre n'avait encore pu établir
les principes du mouvement de l'épicycle et de l'orbe incliné selon la
description ci-dessus à cause de la variation d'inclinaison de
l'épicycle dans les deux cas, et de la variation des deux diamètres
dont l'un est incliné et l'autre est tantôt dans le plan de
l'écliptique, tantôt dans le plan de l'orbe incliné. On y est parvenu
grâce à Dieu.

\`A cet effet, on supposera que l'inclinaison maximale de l'orbe
incliné vers le Nord est atteinte à l'Apogée de Vénus et que cette
inclinaison est un sixième de degré. C'est une inclinaison de grandeur
constante vers le Nord, et de même vers le Sud~; mais elle se meut,
elle se déplace avec l'Apogée, et les n{\oe}uds (la tête et la queue)
se déplacent aussi de par son mouvement. Supposons le centre de
l'épicycle et du rotateur sur la droite passant par
l'Apogée. Supposons le déférent dans le plan de l'orbe
incliné. Supposons l'apogée moyen du rotateur incliné vers le Sud par
rapport au plan de l'orbe incliné, d'un angle de cinq minutes d'arc,
et la moitié du rotateur commençant à l'Apogée (c'est la << première >>
moitié) inclinée vers le Nord par rapport au plan de l'orbe
incliné, formant un angle en son centre d'une grandeur de trois
degrés. Supposons l'apogée de l'épicycle incliné par rapport à
l'apogée du rotateur, vers le Nord, de cinq minutes d'arc. Le diamètre
de l'épicycle [passant par] l'Apogée reste donc dans le plan de l'orbe
incliné. Supposons la première moitié de l'épicycle (celle qui
commence à l'Apogée) inclinée d'un demi-degré vers le Nord par rapport
à la moitié du rotateur inclinée vers le Nord.

Ceci étant admis, que les orbes tournent. Que l'orbe incliné tourne
d'un quart de cercle, le déférent d'un quart de cercle aussi, et le
rotateur d'un demi-cercle dans le même temps; alors l'inclinaison de
l'épicycle s'inverse, le diamètre perpendiculaire au diamètre de
l'Apogée arrive dans le plan de l'écliptique, et l'autre diamètre est
incliné par rapport au plan de l'écliptique d'un angle d'une grandeur
de deux parts et demi; ceci, à condition que la variation
d'inclinaison de l'épicycle soit comme l'indique l'\emph{Almageste}.

Quant à l'amendement fait par les Modernes (pour qui l'épicycle est
incliné aux n{\oe}uds, à l'Apogée et au périgée de trois degrés et
demi), et la théorie qu'expose Ptolémée dans les \emph{Hypothèses},
nous supposerons l'orbe incliné comme on l'a indiqué dans la première
théorie. Le déférent est dans le plan de l'orbe incliné.

\newpage\phantomsection
\index{AQABAS@\RL{ra'asa}!AQABAS@\RL{ra's}, tête, n{\oe}ud ascendant}
\index{APBFAH@\RL{_dnb}!APBFAH@\RL{_dnb}, queue, n{\oe}ud descendant}
\index{ACBHAL@\RL{'awj}!ACBHAL@\RL{'awj}, Apogée}
\index{AYAWAGAQAO@\RL{`u.tArid}, Mercure}
\index{AHAWBDBEBJBHAS@\RL{b.talimayUs}, Ptolémée}
\index{BABDBC@\RL{flk}!BABDBC BEAGAEBD@\RL{falak mA'il}, orbe incliné}
\label{var2}
\includepdf[pages=103,pagecommand={\thispagestyle{plain}}]{edit.pdf}\phantomsection

\noindent \`A l'Apogée,
l'apogée du rotateur est incliné par rapport à l'orbe incliné, de cinq
minutes d'arc vers le Sud~; l'apogée de l'épicycle est incliné par
rapport à l'apogée du rotateur, de cinq minutes d'arc vers le Nord~;
la moitié du rotateur\footnote{Ici, le manuscrit porte le mot
  <<~épicycle~>>, \textit{cf.} texte arabe. C'est pourtant le plan
  du rotateur qui doit subir l'inclinaison de $3°30$ ; l'épicyle en
  hérite par composition.}
commençant à l'Apogée est inclinée vers le
Nord de trois degrés et demi. C'est la voie sur laquelle s'appuie
l'amendement fait par les Modernes.

L'inclinaison maximale de l'épicycle est, vers le Sud, $8;40$, et vers
le Nord, $1;3$, quand l'inclinaison maximale de biais\footnote{Il
  est difficile de donner un sens précis à ces deux valeurs. Peut-être
  $1°3'$ désigne-t-elle la latitude Nord maximale atteinte par
  Vénus quand le centre de l'épicycle est à la queue~: le diamètre de
  l'épicycle passant par son apogée est alors incliné de $2°30'$ par
  rapport au plan de l'écliptique, et on a
$$1°3'\simeq 2°30'\times\frac{43;33}{60+43;33}.$$
  La valeur $8°40'$ semble en revanche être erronée, bien que les
  latitudes de Vénus puissent dépasser cette valeur dans le modèle
  fondé sur l'<<~amendement fait par les Modernes~>>.} est $2;30$.
La dernière latitude
Nord est toujours, pour l'orbe incliné, un sixième de degré. La tête
de Vénus précède l'Apogée d'un quart de cercle, et sa queue le suit
d'un quart de cercle. Les dernières latitudes Nord et Sud, pour l'orbe
incliné, sont à l'Apogée et au périgée.

\emph{Remarque}. Concernant la latitude de Vénus, tout ce qu'ont fait
les Anciens et les Modernes en concevant des ceintures qui
s'approchent ou une inclinaison variable de l'épicycle est
impossible. Qui est bien versé dans cet art ne l'ignore pas.

\begin{center}
  \large Deuxième partie

  \normalsize Les latitudes de Mercure
\end{center}
Chez Ptolémée, il y a deux théories, la
première dans l'\emph{Almageste} et la seconde dans les
\emph{Hypothèses}.

\emph{Première théorie}.  Nous supposons un orbe incliné par rapport
au parécliptique, d'inclinaison maximale atteinte en l'Apogée, vers le
Sud, un demi-degré et un quart. Cette situation à l'Apogée est
réciproque de celle au périgée, où l'inclinaison est vers le Nord. De
plus, la ceinture de l'orbe incliné tend à se rapprocher de la
ceinture de l'écliptique jusqu'à ce qu'elles se confondent, puis elle
s'incline à nouveau d'autant dans l'autre sens (et de même pour le
lieu situé à l'opposé, [au périgée]). Or ceci est impossible: on ne
peut le concevoir dans l'astronomie de Ptolémée, et ni lui ni aucun
autre n'a abordé le problème du mobile de ce mouvement.

\newpage\phantomsection
\index{BEBJBD@\RL{myl}!BEBJBD BHAQAGAH@\RL{mIl al-wirAb}, inclinaison de biais}
\index{BEBJBD@\RL{myl}!BEBJBD AJAOBHBJAQ@\RL{mIl al-tadwIr}, inclinaison de l'épicycle}
\includepdf[pages=104,pagecommand={\thispagestyle{plain}}]{edit.pdf}\phantomsection

D'autre part, quand le centre de l'épicycle est au milieu de l'arc
entre les n{\oe}uds, à l'Apogée, à la dernière latitude Sud, il faut
supposer que le diamètre passant par l'Apogée et le périgée de
l'épicycle est dans le plan de l'orbe incliné, et que l'autre diamètre
est incliné de sorte que la première moitié de l'épicycle soit au Sud,
l'autre moitié au Nord, et que l'angle de cette inclinaison (dite
\emph{de biais}) soit sept degrés. Quand
le centre de l'épicycle est aux n{\oe}uds, le diamètre perpendiculaire
à l'Apogée est dans le plan de l'écliptique, et l'autre diamètre est
incliné par rapport au plan de l'écliptique, son périgée vers le Nord
et son Apogée vers le Sud; l'angle de cette inclinaison est six degrés
et un quart, et on l'appelle \emph{inclinaison de
  l'épicycle}. Ainsi, le centre de
l'épicycle est constamment au Sud du plan de l'écliptique, ou bien
dans ce plan (aux n{\oe}uds). C'est cette théorie qu'on trouve dans
l'\textit{Almageste}.

\emph{Seconde théorie}. C'est la théorie des \emph{Hypothèses}. Là,
l'inclinaison de Mercure est constante, vers le Sud à l'Apogée, d'une
grandeur de six degrés (et l'autre extrémité, au périgée, est donc
inclinée vers le Nord de six degrés). Quand le centre de l'épicycle
est entre les deux n{\oe}uds, le diamètre passant par l'Apogée
coïncide avec l'orbe incliné, et l'autre diamètre est incliné (ainsi
que la deuxième moitié de l'épicycle) vers le Nord d'un angle de six
degrés et demi; on l'appelle \emph{inclinaison de biais}. Quand le
centre de l'épicycle est [aux] n{\oe}uds, le diamètre perpendiculaire
au diamètre de l'Apogée est dans le plan de l'écliptique, et le
diamètre de l'Apogée est incliné, son périgée vers le Nord et son
Apogée vers le Sud, d'un angle égal à l'inclinaison de biais,
c'est-à-dire six parts et demi. C'est la seconde théorie; mais ni
Ptolémée ni aucun autre n'a pu établir des principes permettant de
concevoir ces inclinaisons sans perturber les mouvements en longitude.

Nous avons pu concevoir ces deux théories~; gloire à Dieu.

Voici \emph{la théorie sur laquelle on s'est appuyé}. Nous
supposons que l'inclinaison maximale de l'orbe incliné est vers le Sud
à l'Apogée et vers le Nord au périgée (cette inclinaison maximale est
de six degrés dans les \emph{Hypothèses}).

\newpage\phantomsection
\index{BABDBC@\RL{flk}!BABDBCAMAGBEBD@\RL{falak .hAmil}, orbe déférent}
\index{BABDBC@\RL{flk}!BABDBCBEAOBJAQ@\RL{falak mdIr}, orbe rotateur}
\index{AHAWBDBEBJBHAS@\RL{b.talimayUs}, Ptolémée}
\index{AQAUAO@\RL{r.sd}!AEAQAUAGAO AUAMBJAMAI@\RL{al-'ar.sAd al-.sa.hI.haT}, les observations}
\index{BEBJBD@\RL{myl}!BEBJBD BHAQAGAH@\RL{mIl al-wirAb}, inclinaison de biais}
\includepdf[pages=105,pagecommand={\thispagestyle{plain}}]{edit.pdf}\phantomsection

\noindent Nous supposons que l'apogée
du rotateur est incliné vers le Nord de cinq minutes d'arc, que la
seconde moitié du rotateur est inclinée vers le Sud d'un angle de six
degrés, un demi-degré et un huitième de degré par rapport au plan de
l'orbe incliné, que l'apogée de l'épicycle est incliné par rapport à
l'apogée du rotateur vers le Sud de cinq minutes d'arc, et que sa
seconde moitié est inclinée de sept degrés vers le Nord par rapport au
plan de l'orbe incliné (elle est donc inclinée d'un quart et un
huitième de degré par rapport au plan du rotateur). Quand les orbes se
meuvent, l'orbe incliné d'un quart de cercle, le déférent aussi, et
l'épicycle d'un demi-cercle, alors l'inclinaison du diamètre de
l'Apogée devient six degrés et un quart, et le diamètre qui lui est
perpendiculaire passe dans le plan de l'écliptique. Or il en est ainsi
à l'observation, c'est donc la voie sur laquelle on s'est appuyé; elle
ne prend pas en considération que, selon Ptolémée, le centre de
l'épicycle serait toujours au Sud, mais il est certes revenu de cette
opinion dans les \emph{Hypothèses}. Sache cela.

La seconde théorie est conçue comme suit. Nous supposons l'orbe
incliné comme on l'a décrit ci-dessus, de sorte qu'il coupe le
parécliptique aux n{\oe}uds et que l'inclinaison maximale soit vers le
Sud à l'Apogée et vers le Nord au périgée. Nous supposons que le
centre de l'épicycle est à l'Apogée au milieu de l'arc compris entre
les n{\oe}uds. Nous supposons que l'Apogée du rotateur est incliné de
cinq minutes d'arc vers le Nord, et que sa seconde moitié est inclinée
vers le Nord, et l'autre moitié vers le Sud, d'un angle de six degrés
et demi; on appelle cela \emph{inclinaison de biais}. Nous supposons
l'apogée de l'épicycle incliné par rapport à l'apogée du rotateur de
cinq minutes d'arc vers le Sud; alors le diamètre de l'Apogée est dans
le plan de l'orbe incliné. Nous supposons que la seconde moitié de
l'épicycle est inclinée de six degrés et demi par rapport au plan de
l'orbe incliné. Que l'on se représente cela, puis que les orbes se
meuvent jusqu'à ce que le centre de l'épicycle soit au n{\oe}ud. Le
diamètre perpendiculaire au diamètre de l'Apogée coïncide avec le plan
de l'écliptique, et le diamètre de l'Apogée reste incliné de six
degrés et demi par rapport au plan de l'écliptique. La moitié [de
  l'épicycle] contenant son Apogée est vers le Sud, et l'autre moitié
vers le Nord. C'est ainsi qu'il faut concevoir les mouvements des
orbes en latitude, de façon à ne pas perturber les mouvements en
longitude.

\newpage\phantomsection
\index{AWBHBD@\RL{.twl}!AWBHBD@\RL{.tUl}, longitude (par rapport à l'écliptique)}
\index{BFBBBD@\RL{nql}!AJAYAOBJBD BFBBBD@\RL{ta`dIl al-naql}, équation du déplacement}
\index{AOBHAQ@\RL{dwr}!BEAOAGAQAGAJ AYAQBHAV@\RL{mdArAt al-`rU.d}, trajectoires selon les latitudes (\textit{i. e.} cercles parallèles à l'écliptique)}
\index{BBBHBE@\RL{qwm}!BEBBBHBE BCBHBCAH@\RL{mqwm al-kawkab}, astre vrai (= longitude vraie de l'astre)}
\index{BBBHBE@\RL{qwm}!AJBBBHBJBE@\RL{taqwIm}, 1) calcul de la longitude vraie, 2) longitude vraie d'un astre}
\index{BBBHBE@\RL{qwm}!AJBBBHBJBE BBBEAQ@\RL{tqwIm al-qamar}, Lune vraie}
\label{var3}
\includepdf[pages=106,pagecommand={\thispagestyle{plain}}]{edit.pdf}\phantomsection

\begin{center}
  \large Remarque
\end{center}
Supposer, au besoin, que les orbes sont dans les plans des [orbes qui]
les portent ne perturbe pas les mouvements en longitude de manière
trop importante. Tu sais que ceci est analogue à l'équation du
déplacement\footnote{Le raisonnement peu
  rigoureux d'Ibn al-\v{S}\=a\d{t}ir est à peu près le suivant. Pour
  chaque astre, aussi bien que pour la Lune qu'il a traitée dans un
  chapitre antérieur, l'inclinaison de l'orbe incliné a peu
  d'influence sur le mouvement en longitude (l'équation du
  déplacement, qu'il sait calculer, est faible). Donc les inclinaisons
  des autres orbes auront \textit{a fortiori} une influence
  négligeable sur les mouvements en longitude.} de la Lune vraie de
l'orbe incliné au parécliptique. Puisque la dernière latitude de la
Lune est de cinq degrés, et que la variation maximale due au
déplacement est de six minutes d'arc et deux tiers, l'omission de
l'équation du déplacement des trajectoires selon les latitudes à la
ceinture de l'écliptique aura pour effet une perturbation d'une minute
d'arc et un tiers, par degré de latitude, au plus: c'est une grandeur
négligeable.

Calcule le déplacement si tu le souhaites. En voici la méthode. Prends
l'élongation entre l'astre vrai et la tête de cet astre. Prends
l'équation qui correspond à cette élongation dans la table de
l'équation du déplacement de la Lune. Multiplie cela par la latitude
de cet astre, et divise le produit par la latitude maximale de la Lune
(c'est cinq degrés). En sort l'équation qu'il faut ajouter à l'astre
vrai si le nombre entré dans la table était compris
entre quatre-vingt-dix et cent quatre-vingt ou entre deux cent
soixante-dix et trois cent soixante degrés, et qu'il faut sinon
soustraire de l'astre vrai; reste [l'astre] vrai véritable,
rapporté à l'écliptique. Cette équation est maximale dans les
octants, et nulle aux n{\oe}uds et aux dernières latitudes de part et
d'autre. Dieu est le plus grand et le plus savant.
\label{lat_fin}

\newpage\phantomsection
\index{ASAQAY@\RL{sr`}!ASAQAYAI@\RL{sur`aT}, vitesse}
\index{AHAWBH@\RL{b.tw}!AHAWAB@\RL{bu.t'}, lenteur}
\index{BHBBBA@\RL{wqf}!BHBBBHBA@\RL{wuqUf}, station|see{\RL{rujU`}}}
\index{AQALAY@\RL{rj`}!AQALBHAY@\RL{rujU`}, rétrogradation}
\index{AOBHAQ@\RL{dwr}!AJAOBHBJAQ@\RL{tadwIr}, épicycle}
\index{AMAQBC@\RL{.hrk}!ARAGBHBJAI AMAQBCAI@\RL{zAwiyaT al-.harakaT}, mouvement angulaire}
\addcontentsline{toc}{chapter}{I.26 Cause de la vitesse, de la lenteur, des stations et des rétrogradations des astres}
\includepdf[pages=107,pagecommand={\thispagestyle{plain}}]{edit.pdf}\phantomsection
\begin{center}
  \Large Chapitre vingt-six

  \large Cause de la vitesse, de la lenteur, des stations et
des rétrogradations des astres.
\end{center}
On a vu que le centre de l'épicycle était sur le bord de l'orbe
incliné, que le mouvement de l'orbe incliné était dans le sens des
signes, et que le mouvement de l'épicycle était aussi dans le sens des
signes dans sa partie supérieure. Si l'astre est dans la partie
supérieure de l'épicycle, les mouvements de l'épicycle et de l'orbe
incliné sont dans le sens des signes, et il en résulte que l'astre est
rapide à cause de l'addition des deux mouvements qui vont dans le même
sens. Si l'astre est dans la partie inférieure de l'épicycle, le
mouvement de l'épicycle est en sens inverse du mouvement de l'orbe
incliné. Dans ce cas, si le mouvement angulaire de l'orbe incliné par
jour est égal au mouvement angulaire de l'épicycle, l'astre paraît
stationnaire~; si le [mouvement] angulaire de l'épicycle est supérieur
au [mouvement] angulaire de l'orbe incliné, l'astre paraît rétrograder
d'autant que l'excédent~; et si le [mouvement] angulaire de l'épicycle
par jour est inférieur au [mouvement] angulaire de l'orbe incliné,
l'astre paraît avancer d'autant que la différence entre les deux
[mouvements] angulaires.

Si nous représentons un orbe d'épicycle, et que nous menons en son
centre une droite issue du centre de l'écliptique, alors le rapport du
mouvement de l'épicycle au mouvement de l'orbe incliné est soit
inférieur, soit égal, soit supérieur au rapport de la droite joignant
le centre de l'orbe incliné et le périgée de l'épicycle au rayon de
l'épicycle.

\emph{Si le rapport est inférieur}, alors les seuls effets causés en
l'astre par les deux mouvements sont la vitesse dans la partie
lointaine et la lenteur dans la partie proche, parce que le mouvement
dans la partie lointaine est la somme des deux mouvements, et qu'il
est, dans la partie proche, la différence entre le mouvement de l'orbe
incliné et le mouvement de l'épicycle. L'astre n'a alors ni station ni
rétrogradation, car il y a station quand le rapport des droites
mentionnées est égal au rapport des deux mouvements, et la
rétrogradation est établie quand le rapport est inférieur (ce rapport
est, [au périgée de l'épicycle], le plus petit de ces rapports, donc
un [autre] rapport, égal ou inférieur, ne se peut trouver).

\emph{Si le rapport est égal}, alors l'astre aura une station au
milieu de la période de lenteur, c'est-à-dire quand il sera à distance
minimale et sur la droite mentionnée.

\newpage\phantomsection
\index{APAQBH@\RL{_drw}!APAQBHAI BEAQAEBJBJAI@\RL{_dirwaT mar'iyyaT}, apogée apparent}
\index{BBBHBE@\RL{qwm}!BEBBBJBE AQALBHAY@\RL{muqIm lil-rujU` / lil-istiqAmaT}, en station rétrograde / en station directe}
\index{BBBHAS@\RL{qws}!BBBHAS AQALBHAY@\RL{qws al-rujU` / al-istiqAmaT}, arc rétrograde / arc direct}
\includepdf[pages=108,pagecommand={\thispagestyle{plain}}]{edit.pdf}\phantomsection
\noindent Il n'y aura pas rétrogradation, car la rétrogradation est
établie quand le rapport des droites est inférieur au rapport des deux
mouvements, mais ce rapport, étant égal au plus petit de tous les
rapports, ne peut en avoir un qui le minore~: en effet le plus petit
des rapports est atteint dans l'épicycle quand la droite passe par son
centre. L'astre n'aura donc pas de rétrogradation.

\emph{Si le rapport est supérieur}, alors l'astre aura une
rétrogradation dans la partie proche, entre deux stations. Les
stations se situent des deux côtés de l'épicycle, sur deux droites
issues du centre de l'écliptique telles que, pour chacune, le rapport
du mouvement de l'épicycle au mouvement de l'orbe incliné soit égal au
rapport du segment de cette droite entre le centre de l'écliptique et
le bord de l'épicycle près du périgée à la demi-corde de l'épicycle
située sur cette droite. \`A l'arrivée de l'astre dans la partie
proche, à la première de ces deux droites (appellée première station),
on dit que l'astre est en \textit{station rétrograde}~: il s'arrête
après un ralentissement progressif, et de là, il rétrograde jusqu'à
son arrivée à la seconde droite. Il rétrograde d'abord de plus en plus
vite, la vitesse maximale étant en vigueur quand l'astre est à
distance minimale, puis il ralentit à nouveau jusqu'à atteindre la
seconde droite, lieu de la seconde station, où l'on dit que l'astre
est en \textit{station directe} et où il subit un second arrêt. L'arc
compris entre les deux droites du côté du périgée s'appelle
\textit{arc rétrograde}, et l'arc compris entre les deux droites du
côté de l'apogée de l'épicycle apparent s'appelle \textit{arc
  direct}. Après le second arrêt, l'astre avance en allant
progressivement de l'arrêt à une course lente, puis à une course
moyenne, puis à une course rapide. Les courses moyennes sont en
vigueur entre la lenteur et la vitesse, quand l'astre est située à
distance moyenne le long de l'épicycle. L'avance la plus rapide a lieu
à l'apogée apparent au milieu de l'arc direct. La rétrogradation la
plus rapide a lieu au périgée apparent, au milieu de l'arc rétrograde.

\newpage\phantomsection
\begin{center}
\begingroup%
  \makeatletter%
  \providecommand\color[2][]{%
    \errmessage{(Inkscape) Color is used for the text in Inkscape, but the package 'color.sty' is not loaded}%
    \renewcommand\color[2][]{}%
  }%
  \providecommand\transparent[1]{%
    \errmessage{(Inkscape) Transparency is used (non-zero) for the text in Inkscape, but the package 'transparent.sty' is not loaded}%
    \renewcommand\transparent[1]{}%
  }%
  \providecommand\rotatebox[2]{#2}%
  \ifx\svgwidth\undefined%
    \setlength{\unitlength}{219.17660712bp}%
    \ifx\svgscale\undefined%
      \relax%
    \else%
      \setlength{\unitlength}{\unitlength * \real{\svgscale}}%
    \fi%
  \else%
    \setlength{\unitlength}{\svgwidth}%
  \fi%
  \global\let\svgwidth\undefined%
  \global\let\svgscale\undefined%
  \makeatother%
  \begin{picture}(1,1.45520107)%
    \put(0,0){\includegraphics[width=\unitlength]{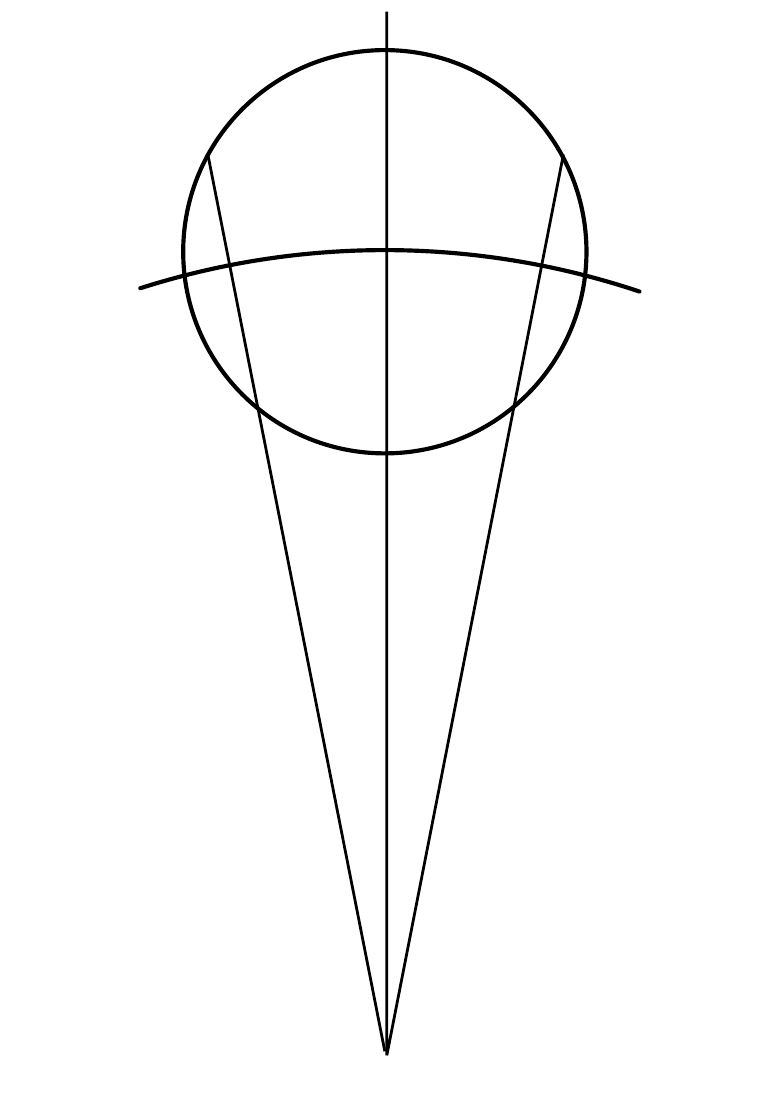}}%
    \put(0.48946017,0.00787749){\color[rgb]{0,0,0}\makebox(0,0)[lb]{\smash{\RL{mrkz al-mA'il}}}}%
    \put(0.79331395,1.1166264){\color[rgb]{0,0,0}\rotatebox{-15.70877313}{\makebox(0,0)[lb]{\smash{\RL{al-mA'il}}}}}%
    \put(0.21188973,1.11793567){\color[rgb]{0,0,0}\rotatebox{13.79720126}{\makebox(0,0)[rb]{\smash{\RL{al-mA'il}}}}}%
    \put(0.44579148,1.43541935){\color[rgb]{0,0,0}\rotatebox{29.81732902}{\makebox(0,0)[rb]{\smash{\RL{flk al-tadwIr}}}}}%
    \put(0.48286453,1.08401568){\color[rgb]{0,0,0}\rotatebox{14.66539518}{\makebox(0,0)[rb]{\smash{\RL{mrkz al-tadwIr}}}}}%
  \end{picture}%
\endgroup%

\end{center}

\newpage\phantomsection
\begin{center}
\begingroup%
  \makeatletter%
  \providecommand\color[2][]{%
    \errmessage{(Inkscape) Color is used for the text in Inkscape, but the package 'color.sty' is not loaded}%
    \renewcommand\color[2][]{}%
  }%
  \providecommand\transparent[1]{%
    \errmessage{(Inkscape) Transparency is used (non-zero) for the text in Inkscape, but the package 'transparent.sty' is not loaded}%
    \renewcommand\transparent[1]{}%
  }%
  \providecommand\rotatebox[2]{#2}%
  \ifx\svgwidth\undefined%
    \setlength{\unitlength}{202.73088015bp}%
    \ifx\svgscale\undefined%
      \relax%
    \else%
      \setlength{\unitlength}{\unitlength * \real{\svgscale}}%
    \fi%
  \else%
    \setlength{\unitlength}{\svgwidth}%
  \fi%
  \global\let\svgwidth\undefined%
  \global\let\svgscale\undefined%
  \makeatother%
  \begin{picture}(1,1.55912002)%
    \put(0,0){\includegraphics[width=\unitlength]{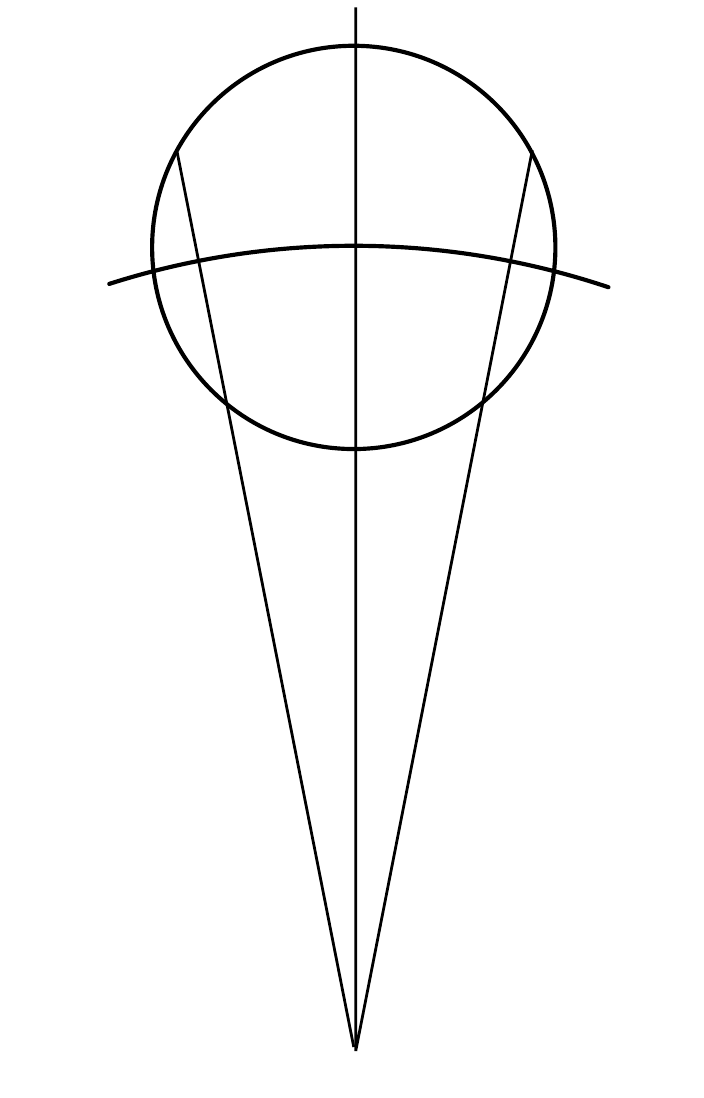}}%
    \put(0.48501184,0.00053951){\color[rgb]{0,0,0}\makebox(0,0)[lb]{\smash{centre de l'orbe incliné}}}%
    \put(0.81351452,1.19923121){\color[rgb]{0,0,0}\rotatebox{-15.70877313}{\makebox(0,0)[lb]{\smash{orbe incliné}}}}%
    \put(0.18492461,1.20064669){\color[rgb]{0,0,0}\rotatebox{13.79720126}{\makebox(0,0)[rb]{\smash{orbe incliné}}}}%
    \put(0.4378007,1.54388495){\color[rgb]{0,0,0}\rotatebox{29.81732902}{\makebox(0,0)[rb]{\smash{orbe de l'épicycle}}}}%
    \put(0.47788115,1.16397508){\color[rgb]{0,0,0}\rotatebox{14.66539518}{\makebox(0,0)[rb]{\smash{centre de l'épicycle}}}}%
  \end{picture}%
\endgroup%

\end{center}

\newpage\phantomsection
\index{ALBEAY@\RL{jm`}!AEALAJBEAGAYAI@\RL{'ijtimA`}, conjonction}
\index{BBAHBD@\RL{qbl}!AEASAJBBAHAGBDAI@\RL{'istiqbAl}, opposition}
\includepdf[pages=109,pagecommand={\thispagestyle{plain}}]{edit.pdf}\phantomsection

Le rapport de la droite entre le centre de l'écliptique et le périgée
de l'épicycle de la Lune ou du Soleil au rayon de l'épicycle est
supérieur au rapport du mouvement de l'épicycle au mouvement de
l'orbe incliné, donc aucun des deux n'a de rétrogradation ni de
station, mais seulement de la lenteur et de la vitesse, d'après ce que
je t'ai expliqué. Le rapport précédent est donc un critère. Sache
qu'avance et rétrogradation sont composés à partir des angles
parcourus lors de mouvements uniformes. On soustrais le plus petit du
plus grand~: si domine ce qui est dans le sens des signes, alors
l'astre avance, sinon il rétrograde. Si additif et soustractif sont
égaux, alors l'astre semble s'arrêter.

\emph{Remarque}\footnote{L'objet de cette remarque est une coïncidence
  importante présente dans tous les modèles médiévaux~; ce sera un
  argument essentiel en faveur de l'héliocentrisme copernicien aux
  yeux de Kepler.}.  Les trois planètes supérieures sont en
conjonction avec le Soleil aux apogées de leurs épicycles (on appelle
cela \uwave{<<~foyer~>>}), et en opposition avec le Soleil aux
périgées de leurs épicycles.  Les planètes inférieures sont en
conjonction avec le Soleil aux apogées et au périgées de leurs
épicycles (on appelle ça aussi un \uwave{<<~foyer~>>}). Dieu est le
plus savant.

\newpage\phantomsection
\index{BBBHAS@\RL{qws}!BBBHAS AQABBJAI@\RL{qws al-ru'yaT}, arc de visibilité}
\index{AXBGAQ@\RL{.zhr}!AXBGBHAQ@\RL{.zuhUr}, visibilité ($\neq$ \RL{i_htifA'})}
\index{ANBABJ@\RL{_hafiya}!ANAJBAAGAB@\RL{i_htifA'}, invisibilité}
\index{AWBDAY@\RL{.tl`}!AWBDBHAY@\RL{.tulU`}, lever (d'un astre)}
\index{AQBJAN@\RL{ry_h}!BEAQAQBJAN@\RL{marrI_h}, Mars}
\index{ARAMBD@\RL{z.hl}!ARAMBD@\RL{zu.hal}, Saturne}
\index{ATAQBJ@\RL{^sry}!BEATAJAQBJ@\RL{al-mu^starI}, Jupiter}
\index{ARBGAQ@\RL{zhr}!ARBGAQAI@\RL{zuharaT}, Vénus}
\index{AYAWAGAQAO@\RL{`u.tArid}, Mercure}
\index{BABHBB@\RL{fwq}!ADBABB@\RL{'ufuq}, horizon}
\index{AQAUAO@\RL{r.sd}!AEAQAUAGAO AUAMBJAMAI@\RL{al-'ar.sAd al-.sa.hI.haT}, les observations}
\index{AMAWAW@\RL{.h.t.t}!BFAMAWAGAW@\RL{in.hi.tA.t}, abaissement|see{\RL{irtifA`}}}
\index{AQBAAY@\RL{rf`}!AQAJBAAGAY@\RL{irtifA`}, hauteur (coordonnées azimutales)}
\addcontentsline{toc}{chapter}{I.27 Visibilité et invisibilité des cinq astres errants}
\includepdf[pages=110,pagecommand={\thispagestyle{plain}}]{edit.pdf}\phantomsection
\begin{center}
  \Large Chapitre vingt-sept

  \large Visibilité et invisibilité des cinq astres errants
\end{center}
Quand l'abaissement du Soleil sous l'horizon au lever\footnote{Mesuré
  le long d'un cercle de hauteur, cet abaissement maximal du
  Soleil sous l'horizon s'appelle <<~arc de visibilité~>> dans les
  paragraphes suivants.} est, pour Saturne \uwave{$11;0$}, pour
Jupiter $10;0$, pour Mars $11;30$, \emph{les trois planètes
  supérieures} sont à la limite de leur visibilité.

Quand elles sont en opposition avec le Soleil, l'arc de visibilité
mesure la moitié des arcs ci-dessus.

L'arc de visibilité de \emph{Vénus et Mercure}, quand ils commencent à
être visibles le soir (et ne le sont plus le matin), est sept parts
pour Vénus et douze parts pour Mercure. Quand ils sont visibles le
matin (et plus le soir), leur arc de visibilité est cinq parts pour
Vénus et sept parts pour Mercure. La cause nécessaire de la diminution
des arcs est la grandeur [apparente] des deux corps parce que les deux
astres sont alors proches du périgée de l'épicycle.

Ces limites sont telles lorsque le centre de l'épicycle est à distance
moyenne~; elles varient sinon en fonction de cela, de la latitude des
planètes, et d'autres raisons, comme nous l'avons expliqué avec
l'élaboration de ce livre dans notre livre \emph{Commentaire des
  observations}.

\newpage\phantomsection
\index{BFBHAQ@\RL{nwr}!BFBJBJAQ@\RL{nayyir}, lumineux}
\index{BFBHAQ@\RL{nwr}!BFBJAQAGBF@\RL{al-nayrAn}, les deux luminaires (Lune et Soleil)}
\index{ALBEAY@\RL{jm`}!AEALAJBEAGAYAI@\RL{'ijtimA`}, conjonction}
\index{BBAHBD@\RL{qbl}!AEASAJBBAHAGBDAI@\RL{'istiqbAl}, opposition}
\index{AXBDBD@\RL{.zll}!BEANAQBHAW AXBDBD@\RL{ma_hrU.t al-.zill}, cône d'ombre}
\index{BBBEAQ@\RL{qmr}!BBBEAQ@\RL{qamar}, Lune}
\index{ATBEAS@\RL{^sms}!ATBEAS@\RL{^sams}, Soleil}
\index{ANASBA@\RL{_hsf}!ANASBHBA@\RL{_husUf}, éclipse de Lune}
\index{BGBDBD@\RL{hll}!BGBDAGBD@\RL{hilAl}, croissant (de Lune)}
\index{AHAUAQ@\RL{b.sr}!AHAUAQ@\RL{ba.sar}, 1)~vision, 2)~observateur}
\index{AXBDBE@\RL{.zlm}!BEAXBDBE@\RL{mu.zlim}, obscur ($\neq$ \RL{nayyir})}
\index{ATAYAY@\RL{^s``}!ATAYAGAY@\RL{^su`A`}, rayon (de lumière)}
\index{AVBHAB@\RL{.daw'}!AVBHAB@\RL{.daw'}, clarté}
\index{BFAXAQ@\RL{n.zr}!ANAJBDAGBA BEBFAXAQ@\RL{i_htilAf al-man.zr}, parallaxe}
\addcontentsline{toc}{chapter}{I.28 Cause des éclipses de Lune et de Soleil, en deux sections}
\includepdf[pages=111,pagecommand={\thispagestyle{plain}}]{edit.pdf}\phantomsection
\begin{center}
  \Large Chapitre vingt-huit

  \large Cause des éclipses de Lune et de Soleil, en deux sections
\end{center}

\begin{center}
  \large Première section

  \normalsize La cause des éclipses de Lune
\end{center}
On a établi que la Lune est un corps sphérique opaque obscur, sous la
lumière du Soleil. Ses rayons ne sont pas transmis par elle~; malgré
son obscurité elle est lisse et les rayons du Soleil sont réfléchis
par elle. La clarté lui est attribuée, mais c'est en réalité les
rayons du Soleil refléchis par le corps de la Lune. La face devant le
Soleil est toujours lumineuse et l'autre face est toujours obscure.

Quand le croissant apparaît, sa partie visible est de la grandeur dont
est déviée la moitié éclairée qui est devant le Soleil. La déviation
augmente chaque nuit jusqu'à l'opposition~; la face devant le Soleil
devient alors visible entièrement car l'observateur s'interpose entre
les deux luminaires. Puis le processus s'inverse~: la partie visible
se soustrait de la lumière jusqu'à la conjonction des deux
luminaires. Puis le processus recommence depuis le début.

La cause des éclipses de Lune est son entrée dans le cône d'ombre de
la Terre. Cela ne peut se produire que pendant l'opposition des deux
luminaires.

Une autre condition est que la latitude [de la Lune] doit être
inférieure à soixante-trois minutes. Si sa latitude est égale à cela,
elle touche le cône d'ombre de la Terre sans y tomber. [Si sa latitude
  est] inférieure à cela, elle entre dans le cône d'ombre de la
Terre~: une partie de la Lune s'éclipse, de la grandeur de ce qui
entre dans le cône d'ombre de la Terre. [Si sa latitude est supérieure
  à un degré et trois minutes, il n'y a pas d'éclipse, la clarté du
  Soleil se propage jusqu'à sa surface, et elle apparaît comme lumière
  en étant réfléchie.

Le fait décisif concerne la latitude de la Lune quand elle est en
opposition avec le Soleil. Si elle est supérieure à la somme du rayon
de la Lune et du rayon de l'ombre de la Terre alors il n'y a pas
éclipse de Lune. Si la latitude est égale aux deux rayons alors la
Lune touche le bord du disque d'ombre, de l'extérieur, du côté de la
latitude de la Lune, et il n'y a pas d'éclipse non plus. Si la
latitude est inférieure aux deux rayons alors il y a bien
éclipse.

\newpage\phantomsection
\index{BCASBA@\RL{ksf}!BCASBHBAAI@\RL{ksUf}, éclipse de Soleil}
\index{ALBEAY@\RL{jm`}!AEALAJBEAGAYAI@\RL{'ijtimA`}, conjonction}
\index{AHAUAQ@\RL{b.sr}!AHAUAQ@\RL{ba.sar}, 1)~vision, 2)~observateur}
\index{BCAKBA@\RL{k_tf}!BCAKAGBAAI@\RL{ki_tAfaT}, opacité}
\label{var23}
\includepdf[pages=112,pagecommand={\thispagestyle{plain}}]{edit.pdf}\phantomsection

\noindent Dans ce dernier cas, si la latitude est pourtant supérieure
au rayon de l'ombre, alors est éclipsée une partie moindre que la
moitié~; si elle lui est égale, alors le bord du disque d'ombre passe
par le centre de cette face de la Lune et la moitié est éclipsée~; si
la latitude est inférieure au rayon de l'ombre mais supérieure à la
différence entre le rayon de la Lune et le rayon de l'ombre, alors est
éclipsée une partie plus grande que la moitié~; si la latitude est
égale à la différence entre le rayon de la Lune et le rayon de
l'ombre, alors la Lune touche le bord du disque d'ombre, de
l'intérieur, du côté de sa latitude, et elle est totalement éclipsée
mais elle ne peut le demeurer~; si la latitude est inférieure à cette
différence, alors l'éclipse est totale et dure en proportion de
combien la Lune s'est engagée dans le disque d'ombre. La durée est
maximale quand le centre de la Lune passe, à mi-temps de l'éclipse, au
centre de disque de l'ombre -- une éclipse dure quand la Lune reste
éclipsée et sombre pendant un certain temps, et cette durée est
proportionnelle au rapport entre la portion de l'écliptique contenue
dans le disque de l'ombre et l'excédent de sa vitesse sur la vitesse
du Soleil.

Le commencement de l'assombrissement suivi de l'éclaircissement est
vers le Sud-Est si la latitude de la Lune est au Nord, il est vers le
Nord-Est si la latitude est au Sud, et si la latitude est nulle
alors [cela dépend de son sens de variation].

La partie sombre de la Lune comporte toujours deux parties
convexes, l'une (à l'opposé de la latitude) appartenant à la
Lune, l'autre appartenant au disque de l'ombre. Sa partie
lumineuse a toujours la forme d'un croissant dont la partie convexe
appartient à la Lune et dont la partie concave appartient au disque de
l'ombre.

\begin{center}
  \large Deuxième section

  \normalsize La cause des éclipses de Soleil
\end{center}
[Le Soleil est éclipsé] quand il cesse d'éclairer
les éléments proches de nous, à un instant où d'habitude il les éclaire.

[Ces éclipses] sont causées par l'interposition de la Lune entre
l'observateur et le Soleil~: elle dérobe alors la lumière du Soleil
aux regards de l'observateur à cause de son opacité. Dès qu'elle
traverse l'intervalle entre l'observateur et le Soleil, l'éclipse
survient. Cela se produit lors des conjonctions qui ont lieu en plein jour,
pourvu que la Lune soit alignée avec le Soleil et l'observateur. Si elle
n'entre pas dans l'alignement avec l'observateur par rapport
au Soleil, celui-ci ne s'éclipsera pas pendant cette conjonction.

\newpage\phantomsection
\index{BFBHAQ@\RL{nwr}!AMBDBBAI BFBHAQAGBFBJAI@\RL{.halaqaT nUrAniyaT}, anneau lumineux (lors d'une éclipse annulaire)}
\index{BGBDBD@\RL{hll}!BGBDAGBD@\RL{hilAl}, croissant (de Lune)}
\index{BFAXAQ@\RL{n.zr}!ANAJBDAGBA BEBFAXAQ@\RL{i_htilAf al-man.zr}, parallaxe}
\includepdf[pages=113,pagecommand={\thispagestyle{plain}}]{edit.pdf}\phantomsection

Ce qui importe est la conjonction des astres tels qu'ils apparaissent
après correction avec la parallaxe~: grâce à cela on sait qu'il
peut y avoir éclipse chez certains peuples et non chez les autres bien
que le Soleil soit au-dessus de l'horizon chez l'un comme chez les
autres, à cause de la variation de la position relative de la Lune par
rapport aux observateurs en fonction du lieu qu'ils habitent. Elle
dépend du pays~: il peut y avoir éclipse dans un pays et non dans
l'autre, ou bien éclipse dans l'un et éclipse dans l'autre mais de
grandeurs différentes, dans des directions différentes, ou de durées
différentes.

Sache que les diamètres des deux luminaires sont égaux en apparence,
ou différents. Qu'ils soient égaux, et que la latitude apparente soit
inférieure à la moitié de la somme des rayons des deux luminaires~;
alors, si elle est égale au rayon du Soleil, il est éclipsé à moitié~;
si elle est supérieure, il est éclipsé moins qu'à moitié~; si elle est
inférieure, il est éclipsé plus qu'à moitié~; si la latitude apparente
est nulle, il est entièrement éclipsé mais cela ne dure pas. Et si les
rayons des deux luminaires étaient différents en apparence~? Que le
rayon du Soleil soit plus grand en apparence que rayon de la Lune~; si
la latitude apparente déjà mentionnée est supérieure au rayon du
Soleil, il est éclipsé moins qu'à moitié~; mais s'il y a égalité, il
est aussi éclipsé moins qu'à moitié, d'autant moins que l'excédent de
son rayon sur le rayon de la Lune est grand~; si la latitude est
inférieure et qu'elle est égale à la différence entre le rayon du
Soleil et le rayon de la Lune, la Lune touche de l'intérieur le bord
du disque du Soleil dont il ne reste plus qu'un anneau lumineux de la
forme d'un croissant~; si la latitude de la Lune est exactement nulle,
à mi-temps de l'éclipse l'anneau deviendra circulaire, entourant le
corps la Lune, de rondeur uniforme~; si la latitude est entre ces deux
états, la Lune est sculptée différemment d'un croissant ou d'un anneau
circulaire parfait. Dans le dernier cas, ainsi que dans le cas d'un
croissant, la partie lumineuse la plus épaisse est du côté opposé à la
latitude apparente. Enfin, que le rayon de la Lune soit plus grand en
apparence que le rayon du Soleil~; si la latitude apparente déjà
mentionnée est égale au rayon de la Lune, il est éclipsé à moitié, le
bord de la Lune passant par le centre du Soleil~; si la latitude est
supérieure à cela, il est éclipsé moins qu'à moitié~; si la latitude
est inférieure à cela et qu'elle est égale à la différence entre le
rayon du Soleil et le rayon de la Lune, il est éclipsé entièrement
mais cela ne dure pas~; si elle est encore inférieure à la différence,
l'éclipse totale dure en proportion de cela~; la durée maximale est
atteinte quand la latitude de la Lune est exactement nulle à mi-temps
de l'éclipse, et cette durée est proportionnelle au rapport entre la
différence des deux diamètres et la vitesse de la Lune.

\newpage\phantomsection
\index{ANASBA@\RL{_hsf}!ANASBHBA@\RL{_husUf}, éclipse de Lune}
\index{AXBDBD@\RL{.zll}!BBAWAQ AXBDBD@\RL{qi.tr al-.zill}, diamètre de l'ombre}
\index{ATBGAQ@\RL{^shr}!ATBGAQ BBBEAQBJBJ BHASAW@\RL{^shr qamariyy}, mois lunaire}
\addcontentsline{toc}{chapter}{I.29 Les intervalles de temps entre les éclipses de Lune et de Soleil}
\includepdf[pages=114,pagecommand={\thispagestyle{plain}}]{edit.pdf}\phantomsection
\begin{center}
  \Large Chapitre vingt-neuf

  \large Les intervalles de temps entre les éclipses de Lune et de Soleil, en deux sections
\end{center}

\begin{center}
  \large Première section
  
  \normalsize Les intervalles de temps entre les éclipses de Lune
\end{center}
Sache que les limites des éclipses sont quand, à l'opposition des deux
luminaires, la latitude de la Lune est inférieure à un degré et quatre
minutes. Cette latitude requiert une distance au n{\oe}ud (quelque
soit le n{\oe}ud, et que la distance soit dans le sens des signes ou
en sens contraire) de douze degrés.  Si la latitude dépasse cela, elle
dépasse la somme des rayons de la Lune et de l'ombre de la Terre.

Sachant cela, sache qu'il ne peut y avoir deux éclipses de Lune à un
mois d'intervalle, car entre les limites des éclipses de Lune, de part
et d'autre, il n'y a pas plus que vingt-cinq parts, or le Soleil
décrit plus que cette grandeur en un mois lunaire et il sort donc des
limites de l'éclipse.

Il ne peut pas non plus y avoir deux éclipses de Lune à sept mois
d'intervalle si l'opposition lors de l'éclipse a lieu avant l'arrivée
au premier n{\oe}ud, à l'extrême limite [du segment entre les
  limites], et que l'autre opposition est après le passage du second
n{\oe}ud, sept mois plus tard, car la seconde opposition ne tombe pas
entre les limites d'une éclipse et elle dépasse la grandeur limite
pour une éclipse, en sens contraire des signes, au delà du second
n{\oe}ud. En effet, le Soleil se meut en sept mois lunaires d'environ
deux cent cinq degrés. Si il est lors de la première opposition à la
limite de l'éclipse (comme on l'a supposé), il dépassera la limite de
l'éclipse d'un degré lors de la seconde opposition~: après une durée
de douze degrés il arrive au premier n{\oe}ud, puis après cent
quatre-vingt degrés il arrive au second n{\oe}ud, et après treize
jours il passe la limite de l'éclipse d'un degré. Cela suppose que le
n{\oe}ud est immobile, mais en fait il se meut de onze degrés pendant
cette durée en sens contraire des signes, donc la distance entre le
Soleil et la limite de l'éclipse est douze degré~; ainsi il ne peut y
avoir deux éclipses de Lune à sept mois d'intervalle.

Le plus souvent, il y a six mois entre deux éclipses car pendant cette
durée le Soleil se déplace du voisinage d'un des n{\oe}uds au
voisinage de l'autre n{\oe}ud.

\newpage\phantomsection
\index{BCASBA@\RL{ksf}!BCASBHBAAI@\RL{ksUf}, éclipse de Soleil}
\index{BFAXAQ@\RL{n.zr}!ANAJBDAGBA BEBFAXAQ@\RL{i_htilAf al-man.zr}, parallaxe}
\index{ASBCBF@\RL{skn}!BEASBCBHBF@\RL{al-maskUn}, la partie habitée (du globe)}
\includepdf[pages=115,pagecommand={\thispagestyle{plain}}]{edit.pdf}\phantomsection

Il peut aussi y avoir cinq mois entre deux éclipses, mais plus
rarement. Il faut qu'une opposition lors d'une éclipse ait lieu après
le passage au n{\oe}ud à l'extrême limite [du segment entre les
  limites], puis qu'une autre opposition ait lieu cinq mois plus tard
avant d'atteindre l'autre n{\oe}ud. Ceci peut arriver à la limite de
l'éclipse, à cause du mouvement du n{\oe}ud qui est pendant cette
durée de huit parts en sens contraire des signes, mais aucune des deux
éclipses ne peut alors être totale (tandis que des éclipses à six mois
d'intervalle peuvent être toutes deux totales, toutes deux partielles,
ou bien l'une totale et l'autre partielle).

\begin{center}
  \large Deuxième section

  \normalsize Les intervalles de temps entre les éclipses de Soleil
\end{center}
On détermine ces intervalles de temps au moyen des limites des
éclipses de Soleil mais elles ne sont pas égales des deux côtés comme
pour les éclipses de Lune, car c'est la latitude apparente qui importe
ici (et non la latitude réelle comme pour les éclipses de Lune), or la
parallaxe doit tantôt être ajoutée à la latitude de la Lune, tantôt
retranchée, pour obtenir la latitude apparente, et la limite par
rapport au n{\oe}ud de part et d'autre différe selon la parallaxe du
lieu.
 
Au centre du monde habité, le sinus de la latitude est de trente-six
parts, et le lieu des éclipses s'étend jusqu'à dix-huit degrés après
le n{\oe}ud ascendant ou avant le n{\oe}ud descendant, mais seulement
jusqu'à neuf degrés avant le n{\oe}ud ascendant ou après le n{\oe}ud
descendant. En effet, la parallaxe maximale à cette latitude est de
soixante-quatre minutes~; or les rayons des deux luminaires font, au
plus, un total de trente-quatre minutes~; si la latitude est
quatre-vingt-dix-huit minutes Nord, retranches-en la parallaxe, il
reste trente-quatre degrés de latitude apparente, or c'est égal aux
deux rayons des luminaires, et la limite de l'éclipse au Nord est donc
bien le premier nombre ci-dessus (la latitude de la Lune est
quatre-vingt-dix-huit minutes quand la Lune est environ dix-huit
degrés après le n{\oe}ud ascendant ou avant le n{\oe}ud descendant).
Quand la Lune a une latitude Sud, son maximum est trente-quatre
minutes, et elle a cette latitude quand elle est six degrés et demi
après le n{\oe}ud descendant ou avant le n{\oe}ud ascendant.

\newpage\phantomsection
\index{BCASBA@\RL{ksf}!BCASBHBAAI@\RL{ksUf}, éclipse de Soleil}
\index{ANASBA@\RL{_hsf}!ANASBHBA@\RL{_husUf}, éclipse de Lune}
\includepdf[pages=116,pagecommand={\thispagestyle{plain}}]{edit.pdf}\phantomsection

Ceci étant admis, sache que~:

-- il ne peut y avoir deux éclipses de Soleil à un mois
d'intervalle au même endroit,

-- mais cela peut arriver dans deux endroits différents, l'un dans
les latitudes Nord, l'autre dans les latitudes Sud.

Le premier point résulte du fait que les limites dans lesquelles peut
avoir lieu l'éclipse, avant le n{\oe}ud et après le n{\oe}ud, sont
distantes de vingt-cinq degrés, or le Soleil se meut d'environ un
signe par mois lunaire, c'est-à-dire trente degrés, donc il sort des
limites de l'éclipse et celle-ci ne peut avoir lieu. Pour le second
point, si l'on fait avec les latitudes Sud dans les endroits situés
dans les latitudes Sud comme on l'a fait avec les latitudes Nord dans
les endroits situés dans les latitudes Nord, on trouve que les deux
limites des éclipses sont distantes de trente-six degrés, or le Soleil
décrit trente degrés par mois, donc il peut retomber dans les limites
de l'éclipses.

Il arrive qu'il y ait deux éclipses de Soleil au même endroit à cinq
mois d'intervalle, l'une après le n{\oe}ud ascendant, puis l'autre
avant le n{\oe}ud descendant. Cela arrive aussi à sept mois
d'intervalle, l'une avant le n{\oe}ud descendant, puis l'autre avant
le n{\oe}ud ascendant, mais le premier cas est plus fréquent. Et
l'occurrence de deux éclipses de Soleil à six mois d'intervalle n'est
pas sujette à caution~: c'est le cas le plus fréquent. Enfin, il peut
y avoir une éclipse de Soleil la moitié d'un mois après une éclipse de
Lune.

Les limites de l'éclipse de Lune sont comptées le long d'un cercle
passant par le centre de l'ombre, perpendiculaire à l'orbe incliné~;
et la limite de l'éclipse de Soleil est comptée le long d'un arc
passant par le centre de la Lune, perpendiculaire à l'écliptique. Dieu
est le plus savant.

\newpage\phantomsection
\index{AKBBBD@\RL{_tql}!BEAQBCAR AKBBBD@\RL{markaz al-_tql}, centre de gravité}
\index{AMALBE@\RL{.hjm}!BEAQBCAR AMALBE@\RL{markaz al-.hjm}, centre du volume}
\index{AOBHAQ@\RL{dwr}!ASAJAOAGAQAI@\RL{istidAraT}, sphéricité (de la Terre, des cieux)}
\index{BAAQASAN@\RL{frs_h}!BAAQASAN@\RL{farsa_h j frAs_h}, parasange (= 3 milles)}
\index{BEBJBD@\RL{myl}!BEBJBD BEBJAGBD@\RL{mIl j amyAl}, mille (= 4000 coudées modernes)}
\index{APAQAY@\RL{_dr`}!APAQAGAY@\RL{_dirA` j a_drA`}, coudée (coudée moderne = 24 doigts)}
\index{AUAHAY@\RL{.sb`}!AUAHAY@\RL{a.sb`}, doigt (= 6 grains d'orge)}
\index{ATAYAQ@\RL{^s`r}!ATAYBJAQAI BEAYAJAOBDAI@\RL{^sa`IraT mu`tadalaT}, grain d'orge (unité de longueur)}
\index{BFAUBA@\RL{n.sf}!ANAWAW BFAUBA BFBGAGAQ@\RL{_ha.t.t n.sf al-nahAr}, ligne méridienne}
\index{BABDBC@\RL{flk}!BABDBC@\RL{flk}, 1) orbe, 2) cieux}
\index{ALAHBD@\RL{jbl}!ALAHBD@\RL{jabal j jibAl}, montagne}
\index{ACAQAV@\RL{'ar.d}!ACAQAV@\RL{'ar.d}, la Terre}
\index{AZBHAQ@\RL{.gwr}!AZBHAQ@\RL{.gawr j 'a.gwAr}, terrain encaissé}
\addcontentsline{toc}{chapter}{I.30 Science de la forme de la Terre et de sa surface}
\includepdf[pages=117,pagecommand={\thispagestyle{plain}}]{edit.pdf}\phantomsection
\begin{center}
  \Large Trentième et dernier chapitre

  \large Science de la forme de la Terre et de sa surface
\end{center}
En astronomie, on a déjà établi des preuves que l'ensemble de la Terre
et de l'eau est sphérique, qu'elle est au centre des cieux, que son
centre de gravité coïncide avec le centre du Monde, et qu'elle est
immobile au centre [des cieux]. Elle ne s'incline d'aucun côté
autrement que ce qu'exige la différence entre son centre de gravité et
le centre de son volume. Ainsi la sphéricité de la surface de la Terre
et de l'eau est parallèle à la sphéricité des cieux~; les reliefs qui
proviennent des montagnes et des terrains encaissés ne la font pas
sortir des limites de la sphéricité, car ils sont insensibles par
comparaison avec la grandeur de la Terre. En effet, soit une montagne
haute d'une demi-parasange\footnote{La parasange est une unité de
  longueur. On va voir que~:

-- 1 parasange = 3 milles

-- 1 mille = 4000 coudées modernes

-- 1 coudée moderne = 24 doigts (1 coudée ancienne = 32 doigts)

-- 1 doigt = 6 grains d'orge

Dans ce paragraphe, {\shatir} compte en coudées modernes.},
c'est-à-dire six mille coudées. Par rapport à l'ensemble de la Terre,
elle est moindre qu'un cinquième d'un septième de la largeur d'un
grain d'orge\footnote{Si la circonférence terrestre est de 8000 parasanges, et que la coudée (moderne) mesure 144 grains d'orge, alors on a bien~:
  $$\frac{\dfrac{1}{2}\times 144}{8000\times \dfrac{7}{22}}=\frac{198}{1000}\times\frac{1}{7}<\frac{1}{5}\times\frac{1}{7}.$$
} par rapport à une sphère d'une coudée de diamètre~; c'est
bien insensible par comparaison avec une telle sphère.

Sache que les grands cercles à la surface de la Terre sont parallèles
aux grands cercles célestes. Ceux qui sont à sa surface se subdivisent
comme eux en trois cent soixante parts, et chaque part du cercle
terrestre est nommée comme son homologue du cercle céleste.

Si un voyageur se déplace le long du méridien, tout droit, sur une
terre uniforme, jusqu'à ce que le pôle s'élève ou bien s'abaisse d'un
degré exactement, alors la grandeur dont il s'est déplacé le long de
ce cercle est d'un degré, et le périmètre entier [du cercle] fait
trois cent soixante fois la grandeur dont il s'est déplacé.

\newpage\phantomsection
\index{AEAHAQANAS@\RL{'ibr_hs}, Hipparque}
\index{AHAWBDBEBJBHAS@\RL{b.talimayUs}, Ptolémée}
\index{AOAQAL@\RL{drj}!AOAQALAI@\RL{drjaT j darjAt, drj}, degré}
\index{ALBEBGAQ@\RL{jmhr}!ALBEBGBHAQ@\RL{al-jumhUr}, les Grecs}
\index{ACAQAV@\RL{'ar.d}!AOBHAQ ACAQAV@\RL{dwr al-'ar.d}, circonférence terrestre}
\includepdf[pages=118,pagecommand={\thispagestyle{plain}}]{edit.pdf}\phantomsection
\noindent Hipparque et Ptolémée (auteur de l'\emph{Almageste}) ont déjà
montré cela et ils ont trouvé qu'une portion d'un degré du grand
cercle imaginaire sur la Terre mesure soixante-six milles et deux
tiers de mille, si le mille mesure trois mille
coudées\footnote{Coudées anciennes}, la coudée\footnote{Coudée
  ancienne} mesure trente-deux doigts, et le doigt six grains d'orge
alignés \uwave{dans le sens de la longueur}. Une école de savants
s'est mise à observer cela au temps d'al-Ma'm\=un et sous ses
ordres\uwave{...} Ils ont trouvé qu'une portion d'un degré faisait
cinquante-six milles et deux tiers de mille, si le mille mesure quatre
mille coudées, la coudée vingt-quatre doigts, et le doigt six grains
d'orge.
J'ai certes enseigné qu'il y a
eu parmi nos prédécesseurs des conventions différentes quant à la
coudée, au mille et à la parasange (pas de différence quant au
doigt)~: l'ancienne coudée mesure une coudée moderne et un tiers, et
le mille vaut quatre-vingt-seize mille doigts pour les deux écoles,
c'est-à-dire trois mille coudées anciennes, mais quatre mille coudées
modernes. La différence pour les coudées est bien réelle, mais pour
les milles elle n'est que verbale~; quant à la parasange, il n'y a pas
de différence puisqu'elle vaut trois milles chez les anciens comme
chez les modernes, et une parasange fait douze mille coudées chez les
modernes~; ainsi la différence dont j'ai parlé, dans la détermination
du degré, entre les anciens et les modernes, est bien une différence
réelle et résulte d'un écart des procédés. La grandeur de la
différence est dix milles. Ceci étant admis, sache qu'un degré compte,
pour les anciens, vingt-deux parasanges et deux neuvièmes de
parasange, et pour les modernes, dix-huit parasanges et huit neuvièmes
de parasange. Et la mesure selon les anciens est en accord avec les Grecs.

Multiplions une portion d'un degré par trois cent soixante
(c'est-à-dire le nombre de parts que compte la circonférence
terrestre), on obtient la grandeur de la circonférence terrestre~:
huit mille parasanges selon les anciens, six mille huit cents
parasanges selon les modernes (la différence est de mille deux cents
parasanges).

\newpage\phantomsection
\index{AMBHAW@\RL{.hw.t}!BEAMBJAW@\RL{m.hI.t}, contour, périmètre}
\index{ASBHAM@\RL{sw.h}!BEASAGAMAI@\RL{masA.haT al-basI.t / al-jarim}, mesure d'une surface / d'un volume}
\index{ALAQBE@\RL{jrm}!ALAQBE@\RL{jarm j 'ajrAm}, corps, volume|see{\RL{masA.haT}}}
\index{ASBHAM@\RL{sw.h}!BEASAGAMAI@\RL{masA.haT al-basI.t / al-jarim}, mesure d'une surface / d'un volume}
\index{AHASAW@\RL{bs.t}!AHASBJAW@\RL{basI.t}, plan, surface|see{\RL{masA.haT}}}
\index{ASBCBF@\RL{skn}!BEASBCBHBF@\RL{al-maskUn}, la partie habitée (du globe)}
\includepdf[pages=119,pagecommand={\thispagestyle{plain}}]{edit.pdf}\phantomsection
\noindent Si l'on divise le périmètre par trois et un
septième\footnote{Ici, {\shatir } prend donc $\pi\simeq
  3+\dfrac{1}{7}$.}, on obtient le diamètre du globe terrestre (avec
l'eau)~: deux mille cinq cent quarante-cinq parasanges et cinq
onzièmes de parasange selon les anciens, deux mille cent
soixante-trois parasanges et sept onzièmes de parasange selon les
modernes.

Si l'on multiplie le diamètre par la circonférence, on obtient la
surface du globe terrestre (avec l'eau)~: 20363636 parasanges et
quatre onzièmes de parasange selon les anciens, 14712727 parasanges et
trois onzièmes de parasange selon les modernes.
 
On dit que la partie habitée est un quart du globe terrestre~: une
demi-circonférence suivant les longitudes et un quart de circonférence
suivant les latitudes.

La volume du globe terrestre entier avec l'eau est, selon les anciens,
8630116345 parasanges cubes et vingt-et-un cent-vingt-et-unièmes et
deux tiers de cent-vingt-et-unième de parasange cube (ces fractions
font environ un sixième)~; selon les modernes, c'est 5033856498
parasanges cubes et quarante-deux cent-vingt-et-unièmes de parasange
cube\footnote{Ces deux mesures de volume ne sont pas parfaitement
  exactes si l'on adopte les valeurs données ci-dessus pour la surface
  et le diamètre du globe. Si, comme pour les anciens, la longueur du
  méridien est 8000 parasanges, on devrait avoir~:
$$\text{volume de la Terre} =
  \frac{\text{surface}\times\text{diamètre}}{6} =
  8639118457+\frac{36+\dfrac{1}{3}}{121} \text{ parasanges cubes}$$ 
Si la longueur du méridien est 6800 parasanges, on devrait avoir un
  volume de ${5305498622+\dfrac{71+\dfrac{1}{3}}{121}}$ parasanges
  cubes.}.

Si tu veux connaître cela en coudées, alors multiplies le nombre de
coudées cubes d'une parasange par le nombre de parasanges du
volume. Si tu veux connaître la surface de la Terre en coudées, alors
multiplies le nombre de coudées carrées d'une parasange par le nombre
de parasanges de la surface de la Terre~; que ce soit selon les
anciens ou selon les modernes, on obtient ainsi la surface de la Terre
en coudées carrées (la coudée carrée des anciens ou bien celle des
modernes). C'est ce que nous voulions montrer.

\newpage\phantomsection
\index{AHAWBDBEBJBHAS@\RL{b.talimayUs}, Ptolémée}
\index{BFAXAQ@\RL{n.zr}!ANAJBDAGBA BEBFAXAQ@\RL{i_htilAf al-man.zr}, parallaxe}
\index{AQAUAO@\RL{r.sd}!AEAQAUAGAO AUAMBJAMAI@\RL{al-'ar.sAd al-.sa.hI.haT}, les observations}
\index{BBAWAQ@\RL{q.tr}!BFAUBA BBAWAQ BD-AJAOBHBJAQ BD-BEAQAEBJBJ@\RL{n.sf q.tr al-tadwIr al-mar'iyy}, rayon de l'épicycle apparent}
\index{BBBEAQ@\RL{qmr}!BBBEAQ@\RL{qamar}, Lune}
\index{ATBEAS@\RL{^sms}!ATBEAS@\RL{^sams}, Soleil}
\index{BBAWAQ@\RL{q.tr}!BFAUBA BBAWAQ ATBEAS BBBEAQ@\RL{ni.sf qi.tr al-^sams / al-qamar}, rayon apparent du Soleil / de la Lune}
\addcontentsline{toc}{chapter}{I.31 Conclusion~: distances et [grandeurs] des corps}
\label{conclusion}
\includepdf[pages=120,pagecommand={\thispagestyle{plain}}]{edit.pdf}\phantomsection
\begin{center}
  \Large Conclusion

  \large Distances et [grandeurs] des corps, en quatre sections
\end{center}

\begin{center}
  \large Première section

  \normalsize Correction des distances des deux luminaires au centre de la Terre
\end{center}
Ptolémée avait indiqué dans l'\emph{Almageste} que le rayon de
l'orbe incliné de la Lune était cinquante-neuf fois le rayon de la
Terre -- et c'est la distance maximale du centre de l'épicycle. Des
chercheurs en cet art ont corrigé ceci pendant les longs intervalles
de temps les séparant de Ptolémée, et en suivant des voies plus
précises que les siennes~: ils ont trouvé qu'il y avait un défaut dans
la valeur de la parallaxe utilisée par Ptolémée, or ils ont évité
fortement ce défaut en supposant que le rayon du parécliptique faisait
cinquante-huit fois le rayon de la Terre.
\`A partir de là, ils ont trouvé une parallaxe en accord avec
l'observation. Ils ont répété cela~: ils ont eu raison, et nous nous
sommes [ensuite] appuyés [sur ce résultat]. Nous avons déja montré que
le rayon de l'épicycle apparent pendant les oppositions et les
conjonctions est cinq parts et un sixième, en parts telles que le
rayon du parécliptique en compte soixante. Cela fait cinq parts, en
parts telles que son rayon en compte cinquante-huit. Si l'on ajoute
cela au rayon du parécliptique, on obtient soixante-trois parts,
distance maximale de la Lune au centre de la Terre pendant les
conjonctions et les oppositions. Il est déjà démontré que le diamètre
[apparent] de la Lune à cette distance de soixante-trois diamètres
terrestres (c'est la distance maximale pendant les conjonctions et les
oppositions) est $0;30,18$.

\label{diam_sol2}De nombreuses observations ont révélé que le diamètre
du Soleil à distance maximale est $0;29,5$, à distance moyenne
$0;32,32$, et à distance minimale $0;36,54,41$. Nous avons déjà
indiqué dans la configuration des orbes du Soleil que sa distance
minimale est $52;53$, sa distance moyenne $60$, et sa distance
maximale $67;7$. Si nous divisons son diamètre à distance moyenne
par la distance du Soleil, on obtient le diamètre du Soleil à cette
distance~; et si nous divisons le diamètre du Soleil à distance
moyenne par son diamètre à un instant donné, on obtient la distance du
Soleil\footnote{Autrement dit, la distance du Soleil et son diamètre
  apparent sont inversement proportionnels~: $\dfrac{\text{diamètre
      apparent du Soleil}}{60}=\dfrac{0;32,32}{\text{distance du
      Soleil}}$.}.

\newpage\phantomsection
\index{AXBDBD@\RL{.zll}!BBAWAQ AXBDBD@\RL{qi.tr al-.zill}, diamètre de l'ombre}
\index{AXBDBD@\RL{.zll}!BEANAQBHAW AXBDBD@\RL{ma_hrU.t al-.zill}, cône d'ombre}
\index{ANASBA@\RL{_hsf}!ANASBHBA@\RL{_husUf}, éclipse de Lune}
\includepdf[pages=121,pagecommand={\thispagestyle{plain}}]{edit.pdf}\phantomsection

D'après ce que nous avons dit, on peut vérifier que si nous
divisons le diamètre du Soleil à distance moyenne, $0;32,32$, par le
diamètre de la Lune à distance maximale, $0;30,18$, il sort $0;62,30$
comme distance du Soleil\footnote{Il y a erreur~:
  $\dfrac{0;32,32}{0;30,18}=0;64,25,20$. Hélas la valeur $0;62,30$
  sera utilisée dans la suite du texte.} au centre de l'écliptique à
l'instant où le rayon de la Lune est égal au rayon du Soleil quand la
Lune est à distance maximale lors des conjonctions et des oppositions.

Ceci étant admis, voyons le tracé de la figure en forme de sapin. Soit
AB le cercle représentant le Soleil, et C son centre. Soit DE le
cercle représentant la Terre, et R son centre. Soit ANB le cône qui la
contient, et RC sa hauteur. Soit JIH le cercle représentant la Lune,
et I son centre. Soit ARB le cône contenant la Lune et le
Soleil. Reportons la longueur RI en RL, et traçons les
perpendiculaires CB, JI, RE et LM à la hauteur de ce cône c'est-à-dire
la droite NC. Si on les prolonge de l'autre côté [de NC], on obtient
les diamètres de ces sphères. Sache que les droites passant par les
points de tangence avec chacun de ces cercles ne sont pas leurs
diamètres, mais ils en diffèrent peu car la distance est
longue. L'angle ARB contient le Soleil et la Lune à l'instant où la
distance [de la Lune] à la Terre est soixante-trois. Nous avons dit
que le diamètre de la Lune dans cette position est $0;30,18$, et pour
sa moitié IJ c'est l'angle IRJ~; or le rapport de RI (soixante-trois)
à IJ est comme le rapport du sinus de l'angle JIR au sinus de l'angle
IRJ~; donc IJ est connu.

Voici le calcul. Prenons le sinus de l'angle IRJ, c'est-à-dire du
rayon de la Lune, c'est $0;0,15,51,54$. Multiplions-le par RI,
c'est-à-dire soixante-trois. On obtient $0;16,39,29,42$, longueur de
la droite IJ, et c'est le rayon de la Lune si le rayon terrestre est
l'unité. Multiplions-le par le coefficient de l'ombre qui vaut à cette
distance, comme je l'ai vérifié par de nombreuses observations, deux
et quarante-trois minutes. On obtient le rayon de l'ombre\footnote{En
  effectuant ce calcul, on trouve que le rayon de l'ombre vaut
  $0;45,15,17,41,6$ rayon terrestre~; hélas, la suite révèle que {\shatir}
  a mené le calcul en supposant que le rayon de l'ombre est
  $0;45,30,19,41$.}.

La droite LR est comme la droite RI, et les droites LM, RE, IK sont
perpendiculaires à la droite NC donc elles sont parallèles entre
elles~; ainsi la somme des deux droites LM et IK fait le double de
RE~; mais RE est le rayon de la Terre et son double fait cent vingt
minutes~; si l'on en retranche la somme du rayon de l'ombre et du
rayon de la Lune, $62;9,49,23$, on obtient la droite JK qui fait
$0;57,50,10,37$.\footnote{Deuxième erreur~: en refaisant les calculs à
  partir de la valeur de rayon de l'ombre mentionnée dans la note
  précédente, on a trouvé $0;58,5,12,36,54$.}

\newpage\phantomsection
\includepdf[pages=122,pagecommand={\thispagestyle{plain}}]{edit.pdf}\phantomsection

La droite RE est parallèle à la droite JK donc le rapport de RE à JK
est comme le rapport de BR à BJ et comme le rapport de RC à CI~; alors
si RC valait soixante minutes, CI (distance entre les centres des deux
luminaires) vaudrait $57;50,10,37$~; reste RI (distance du centre de
la Terre à la Lune) égal à $0;2,9,49,23$~; et le rapport de ceci aux
soixante minutes de RC est comme le rapport de RI (distance du centre
de la Terre à la Lune), soixante-trois, à RC qui est la distance du
Soleil~; donc la distance du Soleil est connue.

Voici le calcul. Divisons la distance de la Lune, soixante-trois, par
$0;2,9,49,23$, il sort $29;7$. Multiplions-le par soixante, on obtient
1747, distance du Soleil au centre de la Terre à l'instant où son
diamètre est comme le diamètre de la Lune~; or c'est soixante-deux et
demi~; si nous multiplions la distance du Soleil ci-dessus par
soixante et que nous divisons le résultat par soixante-deux et demi,
on obtient 1677, sept minutes et un cinquième de minute, et c'est la
distance moyenne du Soleil si le rayon terrestre est
l'unité.\footnote{Ce calcul de la distance Terre-Soleil comprend donc
  deux erreurs~: l'estimation du rayon de l'ombre et la
  valeur $0;62,30$. Si l'on corrige ces deux erreurs, la distance moyenne
  du Soleil est 1840 rayons terrestres. Mais il est peu sensé de ``corriger''
  ainsi les erreurs de calcul de l'auteur, car ce procédé pour calculer la distance Terre-Soleil est très sensible à la précision des données de l'observation. Notons $u$ le sinus du rayon apparent de la Lune à distance maximale, et $C$ le ``coefficient de l'ombre'', $C=2;43$. La distance Terre-Soleil lors de syzygies où les deux rayons apparents sont égaux est, selon {\shatir}~:
  $$\frac{63}{63(1+C)u-1}.$$
Comme le dénominateur (c'est RI dans le texte) est de l'ordre de $0;2$, on voit aisément qu'une erreur relative de $1 \%$ sur $u$ ou sur $C$ conduit à une erreur relative de $30 \%$ sur la distance Terre-Soleil~!\label{sensibilite}}

Divisons-le par soixante. Le quotient de la division est $27;57,7,12$,
c'est une portion d'un degré du rayon du parécliptique, c'est-à-dire
que chaque degré du rayon du parécliptique fait cette grandeur. On a
montré que la distance maximale du Soleil faisait $7;7$ de plus que sa
distance moyenne, donc si nous multiplions la portion d'un degré par
$7;7$ (on obtient 198 $55;30,14,24$) et qu'on l'ajoute à la distance
moyenne du Soleil, on obtient la distance maximale du Soleil égale à
1876 $2;42$~; et si l'on retranche cela de la distance moyenne, il
reste la distance minimale du Soleil qui est 1478 et un cinquième.

Le rapport de JI à IR est comme le rapport de CB à RC donc CB est
connu. Pour le calculer, multiplions $29;7$ par le rayon de la Lune
$0;16,39,29,42$ puis divisons le résultat par la distance de la Lune
c'est-à-dire soixante-trois, on obtient $7;42,53$ et c'est le rayon du
Soleil si le rayon de la Terre est l'unité\footnote{Il y a
  une erreur de calcul ou de copie : on trouve $7;41,56$ et non
  $7;42,53$.}.\label{obs_diam_sol} 

\newpage\phantomsection
\index{ASBHAM@\RL{sw.h}!BEASAGAMAI@\RL{masA.haT al-basI.t / al-jarim}, mesure d'une surface / d'un volume}
\includepdf[pages=123,pagecommand={\thispagestyle{plain}}]{edit.pdf}\phantomsection
\noindent On vient de montrer que le rayon du Soleil est sept
fois et deux tiers de fois et un vingtième de fois le rayon de la
Terre~; le cube de ceci est 459 $8;6,47,46,25,17$~; or on sait que le
rapport de la boule à la boule est comme le rapport du cube du
diamètre au cube du diamètre~; donc le volume du Soleil mesure le
volume de la Terre quatre cent cinquante-neuf fois et deux tiers de
cinquième de fois environ.

Pour calculer la distance du sommet du cône d'ombre de la Terre,
traçons la droite LS parallèle à NE~; alors le triangle RNE est
semblable au triangle LRS donc le rapport de NR à LR est comme le
rapport de ER à RS, donc RN est connu. Voici le calcul. La droite LM
est comme la droite ES, et la droite ER qui est le rayon de la Terre
vaut soixante minutes~; donc si l'on en retranche $0;45,30,19,41$, il
reste $0;14,29,40,19$ et c'est la droite SR. Divisons par cela la
distance de la Lune qui vaut soixante-trois, on trouve $260;47,18$ et
c'est la distance du sommet du cône d'ombre de la Terre au centre de
la Terre quand le Soleil est à la distance mentionnée ci-dessus. Cela
augmente quand le Soleil s'éloigne du centre de la Terre, et cela
diminue quand il se rapproche.

\newpage\phantomsection
\index{ALASBE@\RL{jsm}!ALASBE@\RL{jsm}, corps, solide} 
\includepdf[pages=124,pagecommand={\thispagestyle{plain}}]{edit.pdf}\phantomsection

\begin{center}
  \normalsize Section
\end{center}
Si nous doublons le rayon de la Lune ($0;16,39,29,42$) on obtient
$0;33,18,59,24$ qu'on arrondit à $0;33,19$. Si nous divisons le
diamètre de la Terre (cent vingt minutes) par cela, on obtient
$3;36,6$, et le diamètre du globe terrestre est donc trois fois et
trois cinquièmes de fois le diamètre de la Lune (plus des fractions
très petites). Si nous divisons le diamètre du Soleil par le diamètre
de la Lune, on trouve vingt-sept et quatre cinquièmes~; donc le
diamètre du Soleil est vingt-sept fois et quatre cinquièmes de fois le
diamètre de la Lune. Si l'on pose que le diamètre de la Lune est une
unité, alors son cube est une unité, le cube formé sur le diamètre de
la Terre (trois et trois cinquièmes) est quarante-six et deux-tiers,
et le cube formé sur vingt-sept et quatre cinquièmes est vingt-et-un
mille quatre cent quatre-vingt-un et un quart environ.

Cela implique que, si le corps de la Lune est l'unité, alors le corps
de la Terre est quarante-six et deux tiers, et le corps du Soleil est
vingt-et-un mille quatre cent quatre-vingt-un et un quart. Le corps du
Soleil est quatre cent cinquante-neuf fois et un quinzième de fois le
corps de la Terre environ. Ces grandeurs ni ces distances ne peuvent
être moins que ce que nous avons dit~; mais elles pourraient être plus
grandes.

\newpage\phantomsection
\begin{figure}[h!]
  \begin{center}
    \tiny
    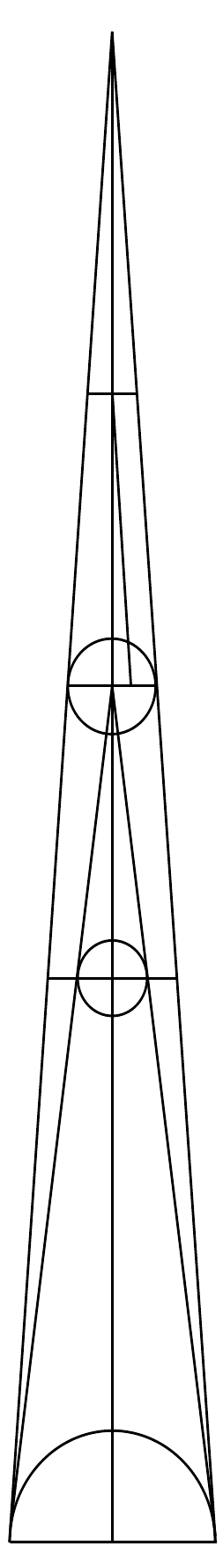
  \end{center}
\end{figure}

\newpage\phantomsection
\includepdf[pages=125,pagecommand={\thispagestyle{plain}}]{edit.pdf}\phantomsection

\newpage\phantomsection
\index{AHAWBDBEBJBHAS@\RL{b.talimayUs}, Ptolémée}
\index{ALAQBE@\RL{jrm}!ALAQBE@\RL{jarm j 'ajrAm}, corps, volume|see{\RL{masA.haT}}}
\index{ASBHAM@\RL{sw.h}!BEASAGAMAI@\RL{masA.haT al-basI.t / al-jarim}, mesure d'une surface / d'un volume}
\index{AXBDBD@\RL{.zll}!BBAWAQ AXBDBD@\RL{qi.tr al-.zill}, diamètre de l'ombre}
\index{BHAUBD@\RL{w.sl}!BEAJAJAUBD AHAYAVBG AHAHAYAV@\RL{mutta.sil ba`.dh biba`.d}, contigus}
\index{BBAWAQ@\RL{q.tr}!BFAUBA BBAWAQ ATBEAS BBBEAQ@\RL{ni.sf qi.tr al-^sams / al-qamar}, rayon apparent du Soleil / de la Lune}
\includepdf[pages=126,pagecommand={\thispagestyle{plain}}]{edit.pdf}\phantomsection

\begin{center}
  \normalsize Section
\end{center}
Pour Ptolémée, si le rayon de la Lune est l'unité, alors le rayon de
la Terre est trois parts et deux cinquièmes de part, et le rayon du
Soleil est dix-huit fois et quatre cinquièmes de fois comme le rayon
de la Lune~; le rayon du Soleil est alors cinq fois et demi le rayon
de la Terre, et son volume est cent soixante-dix fois le volume de la
Terre. La distance maximale de la Lune est soixante-quatre et un
sixième, et la distance moyenne du Soleil est mille deux cent vingt.

\begin{center}
  \normalsize Remarque
\end{center}
La cause de ces différences de grandeurs est que nous
différons de lui légèrement quant au diamètre de la Lune, au diamètre
de l'ombre et au diamètre du Soleil. Chez lui, le diamètre de la Lune
à distance maximale est $31;20$, chez nous c'est $30;18$. Nous
différons de son diamètre d'une minutes et deux secondes. Chez lui, le
diamètre du Soleil est comme le diamètre de la Lune dans son chapitre
concernant les grandeurs\footnote{\emph{L'Almageste}, livre 5,
  chapitre 14.}, et le rapport du diamètre de l'ombre de la Terre au
diamètre de la Lune est chez lui deux fois et trois cinquièmes de
fois, alors qu'il est chez nous deux fois et quarante-trois
minutes. Des chercheurs en cet art ont mis à l'épreuve de manière
répétée ce en quoi nous différons de lui, et nous aussi~; et si nous
utilisions sa valeur du rapport de l'ombre, la différence serait
encore plus grande, donc la différence provient des diamètres du
Soleil et de la Lune. Sache cela.

\begin{center}
  \normalsize Remarque
\end{center}
Si nous retranchons de la distance minimale du Soleil le rayon du
Soleil $7;43$ (avec pour unité le rayon terrestre), il reste $1470;29$
et c'est la distance minimale atteinte au sein des orbes du Soleil.
Si nous ajoutons cela à sa distance maximale, $1883;46$, on obtient la
distance maximale au sein de ses orbes~; en sus, il y a l'épaisseur du
parécliptique, et en deça, [ce qu'il faut] pour que les orbes soient
contigus.

\begin{center}
  \normalsize Remarque
\end{center}
Nous avons indiqué dans la configuration des orbes de la Lune que sa
distance maximale pendant les conjonctions et les oppositions est
$65;10$, que sa distance minimale est alors $54;50$, que sa distance
maximale pendant les quadratures est $68;0$, et que sa distance
minimale est alors $52;0$, toujours sous l'hypothèse que la distance
moyenne est soixante.

\newpage\phantomsection
\index{AYAQAVBJ@\RL{al-mwyd al-`r.dI}, al-Mu'ayyad al-`Ur\d{d}{\=\i}}
\index{ARBGAQ@\RL{zhr}!ARBGAQAI@\RL{zuharaT}, Vénus}
\index{BABDAM@\RL{ibn afla.h}, Ibn Afla\d{h}}
\index{BHAUBD@\RL{w.sl}!BEAJAJAUBD AHAYAVBG AHAHAYAV@\RL{mutta.sil ba`.dh biba`.d}, contigus}
\index{AYAWAGAQAO@\RL{`u.tArid}, Mercure}
\includepdf[pages=127,pagecommand={\thispagestyle{plain}}]{edit.pdf}\phantomsection
\noindent Si nous convertissons cela en parts telles que
le rayon de la Terre en soit une seule, on trouve que sa distance
minimale pendant les conjonctions et les oppositions est $53;0$, sa
distance maximale est alors $63;0$, sa distance minimale pendant les
quadratures est $50;16$, et sa distance maximale est alors $65;44$. Il
faut ajouter à cette dernière distance le rayon de la Lune,
$0;16,39,30$, on trouve $66;0,39,30$~; en sus, il y a l'épaisseur du
parécliptique. Et si nous retranchons [de la distance minimale pendant
  les quadratures] le rayon de la Lune, il reste $49;59,21$~; en deça,
il y a le [ce qu'il faut] pour que les orbes soient contigus.

\begin{center}
  \large Deuxième section

  \normalsize Les distances de Mercure et Vénus au centre de la Terre (en rayons terrestres)
\end{center}
Nous avons montré dans la configuration des orbes de Mercure que le
rapport de la distance minimale de ces orbes à leur distance maximale
est comme le rapport de trente-et-un à quatre-vingt-neuf et
demi\footnote{\textit{Cf.} cependant p.~\pageref{rapport_merc}, où le
  rapport est $\dfrac{31}{89}$.}, et que le rapport de la distance
minimale des orbes de Vénus à leur distance maximale est comme le
rapport de quatorze à cent six. Nous avons montré que l'orbe le plus
distant de la Lune est son orbe incliné, et qu'il est à soixante-six
rayons terrestres. En sus, il y a l'épaisseur du parécliptique que
l'on pose d'un degré~; donc la distance maximale des orbes de la Lune
est soixante-sept, et c'est aussi la distance minimale des orbes de
Mercure. Si nous la multiplions par le rapport propre à Mercure,
c'est-à-dire que nous la multiplions par quatre-vingt-neuf et demi,
puis que nous divisons par trente-et-un, il en sort cent
quatre-vingt-treize et un tiers et un dixième environ~: c'est la
distance maximale des orbes de Mercure. Si nous la multiplions par le
rapport propre à Vénus, c'est-à-dire que nous la multiplions par cent
six, puis que nous divisons le résultat par quatorze, il en sort mille
quatre cent soixante-quatre et quatre septièmes~: c'est la distance
maximale des orbes de Vénus. Elle est inférieure à la distance
minimale des orbes du Soleil, de cinq parts, deux tiers de part et un
quart de part~: c'est [ce qu'il faut] pour que les orbes du Soleil
soient contigus. Il est ainsi démontré que les lieux de Vénus et de
Mercure sont sous le Soleil, contrairement à l'opinion de Mu'ayyad
al-`Ur\d{d}{\=\i} et de Ibn Afla\d{h}.

\newpage\phantomsection
\index{ATAQBJ@\RL{^sry}!BEATAJAQBJ@\RL{al-mu^starI}, Jupiter}
\index{AQBJAN@\RL{ry_h}!BEAQAQBJAN@\RL{marrI_h}, Mars}
\index{ARAMBD@\RL{z.hl}!ARAMBD@\RL{zu.hal}, Saturne}
\includepdf[pages=128,pagecommand={\thispagestyle{plain}}]{edit.pdf}\phantomsection
\begin{center}
  \large Troisième section

  \normalsize Les distances des trois planètes supérieures au centre de la Terre (en rayons terrestres)
\end{center}
Nous avons montré dans la configuration des orbes des trois
planètes [supérieures] que le rapport de la distance minimale des
orbes de Mars à leur distance maximale est comme le rapport de huit à
cent douze, que le rapport de la distance minimale des orbes de
Jupiter à leur distance maximale est comme le rapport de quarante-deux
à soixante-dix-huit, et que le rapport de la distance minimale des
orbes de Saturne à leur distance maximale est comme le rapport de
quarante-six à soixante-quatorze. On vient de voir que la distance
maximale des orbes portant le Soleil est mille huit cent quatre-vingt-trois,
deux tiers et un dixième, avec en sus l'épaisseur du parécliptique que
nous supposons égale à $6;14$, et on atteint ainsi mille huit cent
quatre-vingt-dix~: c'est aussi la distance minimale des orbes de
Mars. Si nous multiplions cela par le numérateur du rapport propre à
Mars, puis que nous divisons le produit par son dénominateur, on
obtient vingt-six mille quatre cent soixante~: c'est la distance
maximale des orbes de Mars. Si nous multiplions cela par le numérateur
du rapport propre à Jupiter, puis que nous divisons le produit par son
dénominateur, on obtient quarante-neuf mille cent quarante environ~:
c'est la distance maximale des orbes de Jupiter. Si nous multiplions
cela par le numérateur du rapport propre à Saturne, puis que nous
divisons le produit par son dénominateur, on obtient soixante-dix-neuf
mille cinquante-et-un, un sixième et un huitième\footnote{$49140\times
  \dfrac{74}{46}=79051+\dfrac{7}{23}\simeq 79051+\dfrac{7}{24}
  =79051+\dfrac{1}{6}+\dfrac{1}{8}$.}~: c'est la distance
maximale des orbes de Saturne, et c'est aussi la distance minimale des
étoiles fixes. Moindre que cela, c'est impossible~; mais cela pourrait
être plus. Dieu est le plus savant.

\newpage\phantomsection
\index{AEAHAQANAS@\RL{'ibr_hs}, Hipparque}
\index{AHAWBDBEBJBHAS@\RL{b.talimayUs}, Ptolémée}
\index{BBAWAQ@\RL{q.tr}!ACBBAWAGAQ BCBHAGBCAH@\RL{'aq.tAr al-kawAkib}, diamètres apparents des planètes}
\index{AYAWAGAQAO@\RL{`u.tArid}, Mercure}
\index{ARBGAQ@\RL{zhr}!ARBGAQAI@\RL{zuharaT}, Vénus}
\index{ATAQBJ@\RL{^sry}!BEATAJAQBJ@\RL{al-mu^starI}, Jupiter}
\index{AQBJAN@\RL{ry_h}!BEAQAQBJAN@\RL{marrI_h}, Mars}
\index{ARAMBD@\RL{z.hl}!ARAMBD@\RL{zu.hal}, Saturne}
\includepdf[pages=129,pagecommand={\thispagestyle{plain}}]{edit.pdf}\phantomsection
\begin{center}
  \large Quatrième section

  \normalsize Mesures des corps des cinq astres errants
\end{center}
Hipparque et Ptolémée ont observé les diamètres de ces planètes, et
ils ont trouvé que le diamètre de Mercure est un quinzième du diamètre
du Soleil, celui de Vénus un dixième, celui de Mars un vingtième,
celui de Saturne un dix-huitième, et celui de Jupiter un douzième~;
mais ils n'ont pas dit à quelles distances se rapportent leurs
observations. Ptolémée les a adoptées dans les distances moyennes,
mais parmi les chercheurs en cet art, les Modernes les ont adoptées
dans les distances minimales, sauf pour Mercure où ils l'ont adoptée
dans les distances moyennes.

Ceci étant admis, si nous divisons la distance moyenne du Soleil par
son diamètre, on obtient cent huit parts et demi et un cinquième de
part~; gardons ce nombre de côté~; puis divisons la distance moyenne
de Mercure c'est-à-dire 130 et treize minutes par le rapport du Soleil
à son rayon c'est-à-dire quinze, on obtient huit parts et quarante
minutes que l'on divise par le nombre qu'on a mis de côté, et on
obtient quatre minutes et quatre cinquièmes de minute~: c'est le
diamètre de Mercure si le diamètre de la Terre est l'unité.

Si nous prenons de même le dixième de la distance moyenne de Vénus
(huit cent vingt-neuf) et que nous le divisons par le nombre mis de
côté, on obtient deux tiers de part et un dixième de part~: c'est le
diamètre de Vénus si le diamètre de la Terre est l'unité.

Quant à Mars, si nous supposons que ce rapport est le sien dans les
distances moyennes, alors il devrait être à distance minimale plus
grand que Vénus~; or ce n'est pas le cas, donc nous avons pris ce
rapport comme étant le sien quand il est entre sa distance moyenne et
sa distance minimale, c'est-à-dire quand sa distance à la Terre vaut
huit mille trente-deux fois le rayon de la Terre\footnote{8032 est la
  moyenne entre la distance minimale du système d'orbes de Mars (1890)
  et sa distance moyenne à la Terre $\dfrac{1890+26460}{2}=14175$
  rayons terrestres.}. Prenons un vingtième de ce nombre, divisons-le
par cent huit, un demi et un cinquième, on obtient trois parts et deux
tiers de part environ~: c'est le diamètre de Mars si le diamètre de la
Terre est l'unité.

Quant à Jupiter, si nous supposons que ce rapport est le sien dans les
distances moyennes, alors divisons sa distance moyenne (trente-sept
mille huit cents) par douze, on obtient trois mille cent cinquante~;
divisons cela par cent huit, un demi et un cinquième, on obtient
vingt-neuf environ~: le diamètre de Jupiter est donc vingt-neuf fois
le diamètre de la Terre.

\newpage\phantomsection
\index{BBAOAQ@\RL{qdr}!BBAOAQ@\RL{qdr}!magnitude (d'une étoile)}
\includepdf[pages=130,pagecommand={\thispagestyle{plain}}]{edit.pdf}\phantomsection
\noindent Si nous supposons au contraire que ce rapport
est le sien à distance minimale, c'est-à-dire à vingt-huit mille huit
cent onze environ (car la planète Jupiter n'atteint pas la distance
minimale de ses orbes qui est égale à la distance maximale de Mars),
alors le diamètre de Jupiter est vingt-trois fois le diamètre de la
Terre. Au cube, cela fait douze mille cent soixante-sept, et c'est le
rapport de son volume au volume de la Terre~: son volume mesure douze
mille cent soixante-sept fois le volume de la Terre environ.

Quant à Saturne, supposons que le rapport mentionné ci-dessus est le
sien à distance minimale, c'est-à-dire à la distance maximale de
Jupiter à laquelle on ajoute l'orbe déférent de Saturne~: cette
distance est quarante-neuf mille neuf cent
soixante-huit\footnote{\emph{Erreur}~: à la fin du chapitre 12,
  {\shatir} a montré que la distance minimale \emph{atteinte} par
  Saturne, rapportée au rayon de la ceinture de son orbe incliné, est
  $50;5$. En rayons terrestres, cela fait donc~:
  $\dfrac{49140+79051+\dfrac{1}{6}+\dfrac{1}{8}}{2}\times 0;50,5\simeq
  53502$, et non $49968$.}. Divisons-la par le rapport du diamètre du
Soleil au diamètre de Saturne (dix-huit), on obtient deux mille sept
cent soixante-seize. Divisons cela par cent huit, un demi et un
cinquième (le nombre que nous avions mis de côté), on trouve
vingt-cinq, un tiers et un cinquième environ~: c'est le diamètre de
Saturne si le diamètre de la Terre est l'unité.

Quant aux étoiles fixes, celles qui sont de première magnitude ont
pour diamètre un vingtième du diamètre du Soleil, celles qui sont de
deuxième magnitude ont pour diamètre un vingt-deuxième, celles qui
sont de troisième magnitude ont pour diamètre un vingt-quatrième,
celles qui sont de quatrième magnitude ont pour diamètre un
vingt-sixième, celles qui sont de cinquième magnitude ont pour
diamètre un vingt-huitième, et celles qui sont de sixième magnitude
ont pour diamètre un trentième du diamètre du Soleil s'il est à
distance moyenne.  Sache donc ces nombres, puis divise la distance
maximale de Saturne (soixante-dix-neuf mille cinquante-et-un) par cent
huit, un demi et un cinquième (le nombre qu'on avait mis de côté), on
obtient sept cent vingt-sept et un quart~: on appelle cela la base. Si
nous divisons ce nombre par vingt, on obtient le diamètre des étoiles
de première magnitude ; si nous le divisons par vingt-deux, on obtient
le diamètre de celles de deuxième grandeur ; si nous le divisons par
vingt-quatre, on obtient le diamètre de celles de troisième magnitude
; si nous le divisons par vingt-six, on obtient le diamètre de celles
de quatrième magnitude ; si nous le divisons par vingt-huit, on
obtient le diamètre de celles de cinquième magnitude ; et si nous le
divisons par trente, on obtient le diamètre de celles de sixième
magnitude.

\newpage\phantomsection
\includepdf[pages=131,pagecommand={\thispagestyle{plain}}]{edit.pdf}\phantomsection
\noindent Divisons donc par les nombres qu'on a mentionnés ; on
trouve que le diamètre de celles de première magnitude est trente-six
fois, un cinquième de fois et un sixième de fois [le diamètre de la
  Terre], que le diamètre de celles de deuxième magnitude est
trente-trois et un quatorzième environ, que le diamètre de celles de
troisième magnitude est trente, un cinquième et un sixième, que le
diamètre de celles de quatrième magnitude est vingt-huit environ, que
le diamètre de celles de cinquième magnitude est vingt-six environ, et
que le diamètre de celles de sixième magnitude est vingt-quatre, un
cinquième et un sixième.

Si nous mettons au cube les diamètres des astres errants et des fixes,
ceci sera rapporté à l'unité qu'est le cube du diamètre de la
Terre. On a déjà traité le cas de la Lune et du Soleil. Quant à
Mercure, son diamètre est quatre minutes et quatre cinquièmes de
minutes, et le cube de ceci est une seconde, quarante-six tierces,
quarante-cinq quartes, et sept quintes ; donc le rapport de son volume
au volume de la Terre est comme le rapport de l'unité à deux mille
vingt-cinq environ.

Le diamètre de Vénus est quarante-six minutes ; le cube de ceci est
vingt-sept minutes, [Vénus] est donc deux cinquièmes et un vingtième du
globe terrestre.

Nous avons montré que le diamètre de Mars est trois parts et deux tiers
de part si le diamètre de la Terre est l'unité ; son cube est
quarante-neuf, un sixième et un huitième ; le volume de Mars mesure
donc quarante-neuf fois, un sixième et un huitième de fois la Terre
environ.

Le diamètre de Jupiter est vingt-neuf fois le diamètre de la Terre ;
le cube de son diamètre est vingt-quatre mille trois cent
quatre-vingt-neuf, et son volume mesure donc vingt-quatre mille trois
cent quantre-vingt neuf fois le volume de la Terre, à condition que le
rapport du diamètre du Soleil à son diamètre vaille quand [Jupiter]
est à distance moyenne (à distance minimale, voir ce qui
précède). Sache cela.

Quant au diamètre de Saturne, il est vingt-cinq fois, un tiers et un
quart de fois comme le diamètre de la Terre ; son cube est seize mille
deux cent quatre-vingt-quatorze, un demi et un huitième environ ; donc
le volume de Saturne mesure seize mille deux cent
quatre-vingt-quatorze fois le volume de la Terre.

\newpage\phantomsection
\index{AHAWBDBEBJBHAS@\RL{b.talimayUs}, Ptolémée}
\index{BAAQASAN@\RL{frs_h}!BAAQASAN@\RL{farsa_h j frAs_h}, parasange (= 3 milles)}
\includepdf[pages=132,pagecommand={\thispagestyle{plain}}]{edit.pdf}\phantomsection

Quant à la grandeur
des étoiles fixes, le diamètre de celles qui sont de première
magnitude mesure trente-six fois, un cinquième de fois et un sixième
de fois le diamètre de la Terre environ ; le cube de ceci est
vingt-sept mille six cent cinquante-quatre arrondi ; et son volume est
d'autant par rapport au volume de la Terre, c'est à dire
quarante-sept-mille six cent cinquante-quatre fois la Terre environ.

Quant aux étoiles fixes les plus petites, de sixième magnitude, le
diamètre de chacune d'entre elles est vingt-quatre fois, un cinquième
et un sixième de fois le diamètre de la Terre ; le cube de ceci est
quatorze mille quatre cent soixante-trois environ ; donc son volume
est quatorze mille quatre cent soixante-trois fois le volume de la
Terre.

En vertu des fondements véritables, ces grandeurs ne peuvent être
moindres, mais peut-être sont-elles encore plus grandes, et cela ne
nous étonnerait point. La grandeur est à Dieu, il est le plus grand,
sa puissance est immense, il est puissant sur toute chose. Tous les
Anciens étaient d'accord quant à la portée des observations qui ont
conduit à ces choses. De manière répétée, Ptolémée disait déjà que les
distances ne peuvent être moindres que ce qu'il avait indiqué~; quant
aux volumes, ils sont fonctions des distances, et si les distances
sont supérieures alors les volumes le sont aussi~; cela n'est pas
exhaustif dans cette espèce~; et il se contredit à propos de la
configuration et des volumes des astres dans son livre intitulé
\emph{Les Hypothèses planétaires}. Nous demandons le succès à Dieu le
Très Haut.

\begin{center}
  \normalsize Remarque
\end{center}
Pour toute sphère dont le diamètre est connu, on connaît ses grands
cercles, ainsi que les portions de ses grands cercles, sa surface et
son volume. Toute grandeur connue par rapport à une certaine grandeur,
l'est aussi par rapport à tout ce par rapport à quoi cette grandeur-ci
l'est. Par exemple, les distances des orbes et les volumes des
astres sont connus par rapport à la grandeur du rayon de la Terre, et
le rayon de la Terre est connu par rapport aux parasanges, aux milles,
aux coudées, aux doigts et au grain d'orge ; ainsi les distances et
les volumes sont connus aussi par rapport à ces grandeurs-ci.

\newpage\phantomsection
\index{BABDBC@\RL{flk}!ACBBAQAH BBAQAH AJAGASAY@\RL{'aqrab qurb al-tAsi`}, rayon de la cavité du neuvième orbe}
\index{BABDBC@\RL{flk}!BABDBC AJAGASAY@\RL{flk tAsi`}, neuvième orbe|see{\RL{flk m`ddl al-nhAr}}}
\index{ASAQAY@\RL{sr`}!ASAQAYAI@\RL{sur`aT}, vitesse}
\index{ARBEBF@\RL{zmn}!ARBEAGBF@\RL{'azmAn mu`addl al-nahAr}!temps équatoriaux}
\includepdf[pages=133,pagecommand={\thispagestyle{plain}}]{edit.pdf}\phantomsection

Sache que si nous ajoutons à la distance maximale des orbes de Saturne
les diamètres des étoiles fixes, cela fait la distance maximale
atteinte au sein des étoiles fixes ; il est possible que cela soit
davantage, si Dieu le veut, mais il est impossible que cela soit
moins, comme je l'ai enseigné. La distance maximale du huitième orbe
ou la distance minimale du neuvième n'est pas moins que
soixante-dix-neuf mille quatre-vingt-huit fois la grandeur du rayon de
la Terre ; c'est son rayon, et si je le double puis que je multiplie
cela par trois plus un septième, on obtient la circonférence de la
cavité du neuvième orbe ; si nous divisons ensuite la circonférence
par trois-cent-soixante, on obtient une portion d'un degré du bord de
la cavité du neuvième orbe ; si nous divisons le degré par soixante,
on en obtient une portion d'une minute d'arc. Nous avons multiplié le
rayon de la partie convexe de la sphère des fixes (soixante-dix-neuf
mille quatre-vingt-huit) par trois plus un septième, puis nous avons
doublé le résultat, on obtient la circonférence d'un grand cercle
situé à la surface du huitième orbe (qui coïncide avec la cavité du
neuvième orbe)~: 497124 et quatre septièmes. Nous avons divisé cela
par trois cent soixante, on obtient 1380, deux tiers et un quart :
c'est une portion d'un degré d'un cercle de la cavité du neuvième
orbe. Nous divisons cela par soixante, on obtient vingt-trois et des
fractions (en rayons terrestes) : c'est une minute d'arc.
Cette durée\label{quatre_secondes}
suffit à un homme pour compter rapidement jusqu'à six~; le temps qu'il
dise <<~un~>>, l'orbe se meut deux fois comme le diamètre de la
Terre environ. Nous avons montré que quand le neuvième orbe se meut d'un
degré, c'est à dire en un quinzième d'heure, alors il se meut d'une
grandeur mesurant mille trois cent quatre-vingt fois, deux tiers de
fois et un quart de fois le rayon de la Terre. S'il se meut d'une
minute, alors alors il se meut d'une grandeur mesurant vingt-trois
fois le rayon de la Terre. Cela concerne sa face concave~; quant à son
épaisseur, personne ne la connaît sauf Dieu le Très Haut~; une portion
d'un degré de sa face convexe est plus grande qu'une telle portion de
sa face concave, mais personne ne sait de combien, sauf Dieu le Très
Haut.

\newpage\phantomsection
\index{AHAWBDBEBJBHAS@\RL{b.talimayUs}, Ptolémée}
\index{ARBGAQ@\RL{zhr}!ARBGAQAI@\RL{zuharaT}, Vénus}
\index{AQBJAN@\RL{ry_h}!BEAQAQBJAN@\RL{marrI_h}, Mars}
\index{BBAWAQ@\RL{q.tr}!ACBBAWAGAQ BCBHAGBCAH@\RL{'aq.tAr al-kawAkib}, diamètres apparents des planètes}
\index{BFAXAQ@\RL{n.zr}!ANAJBDAGBA BEBFAXAQ@\RL{i_htilAf al-man.zr}, parallaxe}
\index{AHAUAQ@\RL{b.sr}!AHAUAQ@\RL{ba.sar}, 1)~vision, 2)~observateur}
\index{BFBHAQ@\RL{nwr}!BFBHAQ@\RL{nwr}, lumière}
\index{BFBHAQ@\RL{nwr}!BFBJBJAQ@\RL{nayyir}, lumineux} 
\includepdf[pages=134,pagecommand={\thispagestyle{plain}}]{edit.pdf}\phantomsection

\begin{center}
  \normalsize Remarque
\end{center}
Sache que l'apparence des diamètres de Vénus et de Mars n'est pas
fonction de leurs distances. La distance maximale de chacun est le
double de sa distance minimale~; or les diamètres devraient suivre
les distances, c'est-à-dire que chacune de ces deux planètes devrait
paraître avoir, à distance minimale, un diamètre double de son
diamètre à distance maximale~; mais ce n'est pas ce qu'on voit.

Ptolémée a dit~: <<~La raison pour laquelle la vision voit et croit
que les grandeurs des corps ne sont pas dans le même rapport que leurs
distances, c'est l'erreur qui entre dans l'observateur à cause de la
parallaxe. On s'en aperçoit en toute chose vue par l'{\oe}il (il n'y a
pas de différence entre choses de grandeurs différentes, avec la
connaissance de leur diminution continuelle). On voit toute planète
beaucoup plus proche de nous que son état en vérité, à cause du fait
que la vision se réduit au distances auxquelles elle s'est accoutumée
et habituée. Nous avons montré que l'augmentation et la diminution
subie par la grandeur est en raison de l'augmentation et de la
diminution des distances, mais ce rapport est inférieur au rapport
réel, à cause de la faiblesse de la vision, d'après ce que nous avons
dit à propos de la distinction de la quantité de la différence pour
chaque espèce de ce que nous avons mentionné.~>>
Ceci est la réponse de Ptolémée. Nous pensons pouvoir analyser cela ainsi~:

-- Quand une sphère se rapproche de l'observateur, ce qui en est
visible est moindre que ce qui en est visible à distance.

-- L'observateur ne peut percevoir les diamètres réels des corps qui
ont de l'éclat à cause de la réception de la vision par les rayons
lumineux.
  
-- Avec la distance, l'observateur ne voit pas les choses comme elles
sont.

Expliquons cela. Prenons deux corps sphériques, le diamètre de l'un étant
double du diamètre de l'autre, et regardons-les de loin. L'{\oe}il ne
pourra percevoir la grandeur réelle de leurs diamètres, ni ne pourra
voir que le diamètre de l'un est double du diamètre de l'autre, ni ne
pourra le vérifier par l'aide d'aucun instrument, même s'il sait que
le diamètre de l'un est double du diamètre de l'autre.

\newpage\phantomsection
\index{AVBHAB@\RL{.daw'}!AVBHAB@\RL{.daw'}, clarté}
\index{BCAKBA@\RL{k_tf}!BCAKAGBAAI@\RL{ki_tAfaT}, opacité}
\index{AYBCAS@\RL{`ks}!AEBFAYBCAGAS@\RL{in`ikAs}, reflet}
\index{ATAYAY@\RL{^s``}!ATAYAGAY@\RL{^su`A`}, rayon (de lumière)}
\includepdf[pages=135,pagecommand={\thispagestyle{plain}}]{edit.pdf}\phantomsection
\noindent Que l'observateur estime la grandeur de l'un des
diamètres par rapport à
l'autre~; supposons que l'un des deux reste à sa place et que l'autre
se déplace au double de sa distance~; l'observateur ne pourra plus
admettre la différence qu'il avait estimée. D'autre part, si le regard
se porte sur la lumière d'une lampe à distance donnée et que le regard
imagine qu'elle est d'une certaine grandeur, puis qu'on s'approche [de
  cette lampe] à mi-distance ou qu'on s'en éloigne d'autant, alors le
regard n'aura pas l'impression d'un tel changement de distance. \`A
partir de là et de ce qui précède, on a la réponse à la question
indépendamment de l'incapacité de la vision à saisir la grandeur des
choses dans leur état. Nous affirmons que si l'{\oe}il se porte sur Paul et
Pierre, il ne peut déterminer lequel des deux est plus grand (même par la
distance), ni ne peut dire de quelle grandeur l'un est plus grand que
l'autre~; ceci, même s'il n'y a aucune distance entre eux deux. ll est donc
clair que la distance maximale de Mars est vingt-six mille quatre cent
soixante comparée à la grandeur du rayon de la Terre. La science de la
vision échoue pour les grandeurs à des distances bien moindres que
cela ; de plus, elle répugne à la détermination des grandeurs même si
nous percevons que les diamètres des astres diffèrent aux autres
distances et que nous le sentons (sauf que cela ne dépend pas du
rapport des distances comme nous l'avons expliqué). De là, on comprend
que les astres sont des corps qui ont un éclat et que leur lumière
leur est propre, non comme la Lune qui est un corps opaque dont la
lumière vient du reflet des rayons du Soleil à sa surface. Pour cette
raison, l'estimation du diamètre [de la Lune] aux autres distances de
la Terre est possible, mais non celle des diamètres des autres astres.

\newpage\phantomsection
\index{ALAOBD@\RL{jdl}!ALAOBHBD@\RL{jdwl}, table, tableau, catalogue}
\includepdf[pages=136,pagecommand={\thispagestyle{plain}}]{edit.pdf}\phantomsection

Ceci achève ce que nous voulions produire dans cette première
partie. Va suivre une deuxième partie sur la configuration de la Terre
et ce qui en dépend. Nous continuerons cela par le calcul des tables
des équations des astres selon ces fondements dans une troisième
\label{troisieme1}
partie si Dieu le veut. On espère que le lecteur de ce livre ne
rejettera pas ce qui n'est pas parfaitement connu, qu'il empruntera le
chemin de la vérité et de la raison, et qu'il échappera à la voie de
l'entêtement. Le bien-fondé de ce que nous avons produit et la clarté
de la voie que nous avons empruntée paraîtra à qui suivra la voie de
la vérité et confrontera ces fondements avec les Anciens et les
Modernes.

Merci à Dieu, qu'Il soit notre refuge, nous nous en remettons à Lui.

\newpage\phantomsection
\index{ASBHBI@\RL{sw_A}!ANAWAW ASAJBHAGAB@\RL{_ha.t.t al-istiwA'}, équateur terrestre}
\index{ANBDBA@\RL{_hlf}!ACANAJBDAGBA@\RL{i_htilAf}, irrégularité, anomalie, variation}
\index{AYAOBD@\RL{`dl}!BEAYAOAOBD BFBGAGAQ@\RL{m`ddl al-nhAr}, équateur, plan de l'équateur}
\index{ASBCBF@\RL{skn}!BEASBCBHBF@\RL{al-maskUn}, la partie habitée (du globe)}
\index{BGBJAC@\RL{haya'a}!BGBJACAI@\RL{hI'aT}, configuration, astronomie}
\index{AOBHAQ@\RL{dwr}!ASAJAOAGAQAI@\RL{istidAraT}, sphéricité (de la Terre, des cieux)}
\index{AHAMAQ@\RL{b.hr}!AHAMAQ@\RL{b.hr}, mer}
\index{BFBGAQ@\RL{nhr}!BBBHAS BFBGAGAQ@\RL{qws al-nhAr}, arc diurne, durée du jour}
\includepdf[pages=1,pagecommand={\thispagestyle{plain}}]{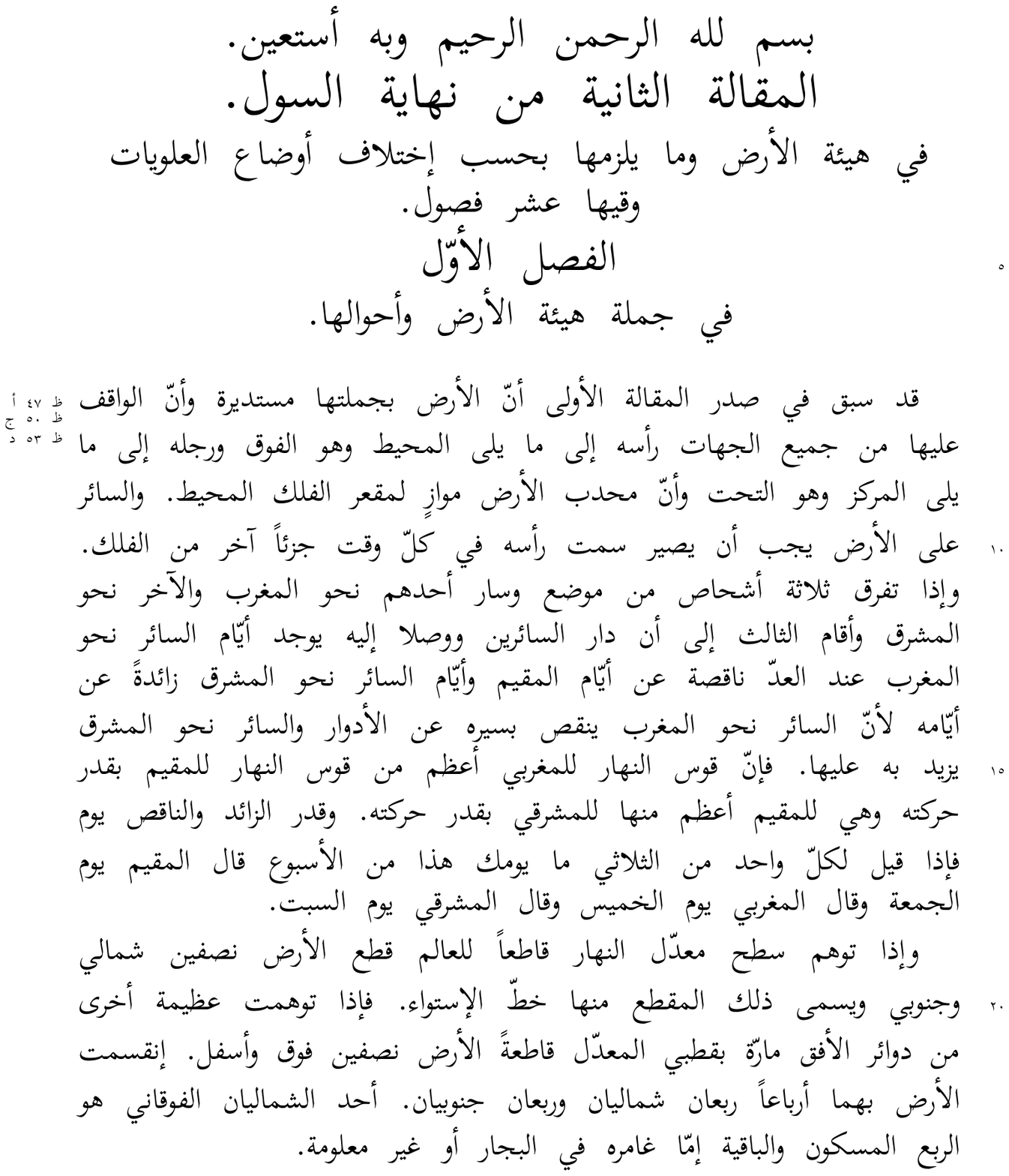}
\addcontentsline{toc}{chapter}{II.1 L'ensemble de la configuration de la Terre et ce qui l'entoure}
\begin{center}
\Large Au nom de Dieu clément et miséricordieux~; j'implore son secours

\LARGE Seconde partie de l'Achèvement de l'enquête

\normalsize La configuration de la Terre et de ce qui en dépend concernant la variation de position des corps célestes. En dix sections.

\Large Première section

\normalsize L'ensemble de la configuration de la Terre et ce qui
l'entoure.
\end{center}

\noindent Nous avons déjà dit au début de la première partie que la
Terre, globalement, est ronde, qu'un individu debout où que ce soit, a
la tête vers la périphérie (c'est-à-dire vers le haut) et les pieds
vers le centre (c'est-à-dire vers le bas), et que la surface convexe
de la Terre est parallèle à la surface concave de l'orbe qui
l'entoure. Le sommet de la tête d'un voyageur, sur Terre, indique sans
cesse une portion différente de cet orbe. Supposons que trois
personnes se séparent en un lieu, que l'un se dirige vers l'Ouest,
l'autre vers l'Est, et que le troisième demeure immobile jusqu'à ce
que les deux autres aient fait un tour et reviennent à lui~; on trouve
que les jours mis par le voyageur allant vers l'Ouest sont en nombre
inférieur aux jours de celui qui est immobile, et que les jours mis
par le voyageur allant vers l'Est sont en nombre supérieur aux
siens. En effet le voyageur allant vers l'Ouest retranche [quelque
  chose] des circonférences à cause de son voyage, et celui qui va
vers l'Est y ajoute [quelque chose] à cause de son voyage.  La durée
du jour pour celui qui va vers l'Ouest est plus grande que pour celui
qui est immobile, à raison de son mouvement~; et la durée du jour pour
celui qui est immobile est plus grande que pour celui qui va vers
l'Est, à raison du mouvement de celui-ci. La quantité ajoutée ou
retranchée est d'une journée, donc si l'on dit à chacun des trois
``quel jour de la semaine est-on ?'' et que celui qui est immobile
répond vendredi, alors celui qui va vers l'Ouest dira jeudi, et celui
qui va vers l'Est dira samedi.

Si l'on imagine que le plan de l'équateur coupe le Monde, il coupe la
Terre en deux moitiés Nord et Sud, et cette coupure s'appelle ligne de
l'équateur.  Imaginons encore un autre grand cercle passant par les
pôles de l'équateur~; il coupe la Terre en deux moitiés dessus comme 
dessous. Ils divisent donc la Terre en quatre~; deux quarts au Nord et
deux quarts au Sud. L'un des quarts Nord, le plus élevé, est le quart
habité. Le reste est soit en grande partie sous les mers, soit inconnu.

\newpage\phantomsection
\index{AMAQAQ@\RL{.harara}!AMAQAGAQAI@\RL{.hrAraT}, chaleur}
\index{AHAMAQ@\RL{b.hr}!AHAMAQ@\RL{b.hr}, mer}
\index{ARBFAL@\RL{znj}!ARBFAL@\RL{znj}, Noirs}
\index{AMAHAT@\RL{.hb^s}!AMAHATAI@\RL{.hb^saT}, \'Ethiopie}
\index{BFBJBD@\RL{nIl}, Nil}
\index{BEAUAQ@\RL{mi.sr}, \'Egypte}
\index{BBBJAS@\RL{qys}!BEBBBJAGAS@\RL{miqyAs j maqAyIs}, gnomon}
\index{AYBEAQ@\RL{`mr}!BEAYBEBHAQ@\RL{al-ma`mUr, al-`imAraT}, partie cultivée de la Terre~; par extension, le quart habité (\RL{maskUn})}
\index{AHAQAO@\RL{brd}!AHAQAO@\RL{brd}, froid}
\index{BDBHBF@\RL{lwn}, couleur}
\index{ASBCBF@\RL{skn}!BEASBCBHBF@\RL{al-maskUn}, la partie habitée (du globe)}
\index{ANASBA@\RL{_hsf}!ANASBHBA@\RL{_husUf}, éclipse de Lune}
\index{AMAQBB@\RL{.hrq}!AWAQBJBBAI BEAMAJAQBBAI@\RL{.tarIqaT mu.htariqaT}, routes brûlantes}
\index{BBBEAQ@\RL{qmr}!BBBEAQ@\RL{qamar}, Lune}
\index{ALAHBD@\RL{jbl}!ALAHBD@\RL{jabal j jibAl}, montagne}
\index{AWBHASBJ@\RL{n.sIr al-dIn al-.tUsI}, Na\d{s}{\=\i}r al-D{\=\i}n al-\d{T}\=us{\=\i}}
\includepdf[pages=2,pagecommand={\thispagestyle{plain}}]{edit2.pdf}

\uwave{La science ... est propre à ...} Certains spécialistes de cette
science disent que les contrées du Sud, parce qu'elles sont proches du
périgée solaire, sont très chaudes~; leur grande chaleur empêche toute
culture. Si c'en était la cause, les contrées du Sud dont la latitude
dépasse l'inclinaison maximale ne seraient pas incultes, et certaines
contrées du Nord seraient inappropriées à la culture, bien que
d'autres ne le sont pas, à latitude et longitude égale. Mais selon
certains de ces spécialistes, la couleur de ceux qui en sont proches
montre la durée pendant laquelle le périgée est au Sud de
l'écliptique, pour ces contrées très chaudes.

La chaleur attire l'humidité, comme le montre la flamme d'une lampe quand
elle attire l'huile ; ainsi les mers seraient-elles entraînées
vers la moitié Sud, et la moitié Nord de la Terre serait-elle laissée
à découvert ? L'existence de mers au Nord empêche ce jugement.

Selon certains, celles de ces contrées qui ne sont pas habitées sont
celles où vont se coucher les luminaires, et on les appelle ``routes
brûlantes'' ; mais c'est un conte invraisemblable.

L'auteur de la \emph{Ta\b{d}kira} a dit : si la raison n'en
était la Providence, pas un des deux quarts Nord ne serait approprié à
la culture. Mais pour atteindre la science de ce qui rend ce quart
approprié à la culture, on sait que la partie habitée est le quart le
plus élevé, dans la moitié Nord, parce qu'en longitude, il n'existe pas,
aux deux extrémités de la partie habitée, de lieux où la limite de
l'avancée des éclipses et de leur retard dépasse douze heures, et
parce qu'en latitude, il n'existe nulle part, dans toutes les contrées habitées,
où l'extrémité de l'ombre du gnomon soit côté Sud à midi des équinoxes, sauf
dans quelques lieux aux confins des pays des Noirs, de l'Ethiopie, etc.

Les premières latitudes habitées au Sud sont là où le pôle Sud est
haut de $16;25$ degrés, ce sont les lieux qu'on vient de
mentionner. Les dernières au Nord sont là où le pôle Nord est haut de
66 degrés ; il est impossible de vivre au delà à cause du froid dû à
l'éloignement du Soleil par rapport au zenith là-bas. Donc les
latitudes habitées mesurent $82;25$ degrés ; les longitudes, $177;15$
degrés. La mer entoure de presque tout côté la part habitable. Le côté
Ouest, le côté Nord, et presque tout le Sud, notamment le Sud-Est,
sont connus. Le Sud-Ouest n'est pas connu. Il est dit que des
voyageurs, en haut du cours du Nil en Egypte, sont arrivés en des
lieux dont la latitude Sud dépasse quelques dix degrés, et qu'ils ont
vu les montagnes blanches, comparables à la Lune, où sont les sources
du Nil ; même plus loin au Sud, ils n'ont pas atteint de mer.  Nous
n'avons pas non plus de témoignage au sujet d'une mer au
Nord-Est.

Dans la partie découverte et habitée, il y a de nombreuses
mers. 

\newpage\phantomsection
\index{AHAWBDBEBJBHAS@\RL{b.talimayUs}, Ptolémée}
\index{ASBHBI@\RL{sw_A}!ANAWAW ASAJBHAGAB@\RL{_ha.t.t al-istiwA'}, équateur terrestre}
\index{ALBEBGAQ@\RL{jmhr}!ALBEBGBHAQ@\RL{al-jumhUr}, les Grecs}
\index{ACBBBDBE@\RL{'aqlm}!AEBBBDBJBE@\RL{'iqlIm j 'aqAlIm}, climat}
\index{ANBDAL@\RL{_hlj}!ANBDBJAL AHAQAHAQBI@\RL{halIj brbr_A}, golfe de Berbera (golfe d'Aden)}
\index{ANBDAL@\RL{_hlj}!ANBDBJAL ACANAVAQ@\RL{_halIj a_h.dar}, <<~golfe vert~>> (golfe d'Oman)}
\index{ANBDAL@\RL{_hlj}!ANBDBJAL BAAGAQAS@\RL{_halIj fAris}, golfe Persique}
\index{ANBDAL@\RL{_hlj}!ANBDBJAL ACAMBEAQ@\RL{_halIj a.hmar}, mer Rouge}
\index{AMBHAW@\RL{.hw.t}!AHAMAQ BEAMBJAW@\RL{al-ba.hr al-mu.hI.t}, l'Océan}
\index{AHAMAQ@\RL{b.hr}!AHAMAQ BHAQBFBC@\RL{b.hr warank}, mer Baltique}
\index{AHAMAQ@\RL{b.hr}!AHAMAQ AWAHAQASAJAGBF@\RL{ba.hr .tabaristAn}, mer Caspienne}
\index{AHAMAQ@\RL{b.hr}!AHAMBJAQAI ANBHAGAQARBE@\RL{bu.hayraT _hawArizm}, mer d'Aral}
\index{AJBHBDBI@\RL{tUl_A}, île de Thulé}
\index{AUBBBDAH@\RL{.saqlab}, Slave}
\index{ANBDAO@\RL{_hld}!ALARAGAEAQ ANAGBDAOAGAJ@\RL{jazA'ir al-_hAlidAt}, îles Canaries}
\index{BBAHAH@\RL{qbb}!BBAHAHAI ACAQAV@\RL{qubbaT al-'ar.d}, dôme de la Terre}
\includepdf[pages=3,pagecommand={\thispagestyle{plain}}]{edit2.pdf}

\noindent Certaines mers touchent l'Océan, comme celle qui part du
Maghreb et d'Al-Andalous et qui est entre Al-Andalous et la Syrie, ou
comme la mer du Sud qui touche au côté Est d'où sortent quatre golfes
vers le milieu des terres habitées~: le golfe de
Berbera\footnote{\textit{i. e.} le golfe d'Aden où se trouve la ville
  de Berbera.}, celui qui est le plus proche du Maghreb, le golfe
Rouge\footnote{\textit{i. e.} la mer Rouge.}, le golfe Persique et le
golfe Vert\footnote{\textit{i. e.} le golfe d'Oman.}.  Chacun d'eux a
une longitude et une latitude avantageuses. [Parmi les mers qui
  touchent l'Océan], il y a aussi la mer des
Warank\footnote{\textit{i. e.} la mer Baltique.} du côté Nord.
Certaines mers ne touchent pas l'Océan, comme la mer du
Tabaristan\footnote{\textit{i. e.} mer Caspienne.}, le lac du
Khwarizm\footnote{\textit{i. e.} mer d'Aral.}, et d'autres qui
  sont des marais et des p\^echeries.

Outre les mers, parmi les choses empêchant l'habitation, il y a les
plaines, les montagnes, les sables, les collines, les maquis, et
beaucoup d'autres que savent les spécialistes des routes, le voyageur,
\emph{etc.} ; la connaissance détaillée de cela est seulement
  dans les livres variés de cette catégorie. La grandeur de la
partie habitée, du côté Nord, tombe entre un peu plus de dix degrés de
latitude jusqu'à la limite des cinquante degrés.

Les hommes de métier l'ont divisée en sept climats longitudinalement,
de sorte que chaque climat couvre l'étendue d'Est en Ouest sous un
certain parallèle ; les conditions des lieux qui s'y trouvent sont
semblables ; son [étendue] en latitude se calcule par des incréments
d'une demi-heure par rapport à la durée du jour le plus long dans le
climat moyen. Du septième climat jusqu'à une latitude de 66 parts, on
trouve peu de constructions, et leurs habitants ressemblent à des
sauvages. Au delà de 63 parts de latitude, il y a l'\^ile de
Thulé\footnote{les Hébrides, l'Islande ou la Norvège peut-être.} dont
les habitants vivent dans les bains à cause du grand froid ; le jour
le plus long y dure vingt heures. Ptolémée a dit que là où la latitude
est de 64 parts, il y a un peuple de Slaves qui ne savent
  pas. Les Grecs placent la longitude en partant de l'Ouest de sorte
que sa quantité augmente dans le sens de signes, ils placent le début
des latitudes sur l'équateur pour les délimiter de manière naturelle,
et ils placent le début des lieux habités aux îles Canaries qui sont
aujourd'hui entourées d'eau.  Un peuple le place au littoral ouest du
Maghreb. Entre eux deux, il y a dix degrés. La dernière limite de la
partie habitée est comme nous avons rapporté -- selon certains, c'est
la demeure des Démons -- et c'est le début de la partie habitée pour
ceux qui le placent du côté de l'Est. Ils appellent ce qui est [au
  milieu] entre les deux limites, sur l'équateur, le Dôme de la
Terre.

\newpage\phantomsection
\index{ACBBBDBE@\RL{'aqlm}!AEBBBDBJBE@\RL{'iqlIm j 'aqAlIm}, climat}
\includepdf[pages=4,pagecommand={\thispagestyle{plain}}]{edit2.pdf}

\noindent C'est situé à un quart de circonférence du commencement Ouest ;
il y a donc des variantes concernant ce nom, en raison des variantes
concernant le commencement Ouest.

\emph{Le premier climat} commence là où le jour le plus long dure douze
heures trois quarts et où la latitude est douze degrés et
\uwave{[deux ?]} tiers de degré.  En son milieu, le jour dure treize
heures et la latitude est seize degrés, un demi-degré et un huitième
de degré.\footnote{Début $12;20$ ou $12;40$. Milieu $16;37,30$.}

\emph{Le deuxième climat} commence là où le jour dure treize heures et
quart et où la latitude est vingt degrés, un quart et un cinquième de
degré.  En son milieu, le jour dure treize heures et demi, et la
latitude est vingt-quatre degrés, un demi et un sixième de degré.
\footnote{Début $20;27$. Milieu $24;40$.}

\emph{Le troisième climat} commence là où le jour dure treize heures
trois quarts et où la latitude est vingt-sept degrés et demi. En son
milieu, le jour dure quatorze heures, et la latitude est trente degrés
et deux tiers de degré.\footnote{Début $27;30$. Milieu $30;40$.}

\emph{Le quatrième climat} commence là où le jour dure quatorze heures
et quart, et où la latitude est trente-trois degrés et demi et un
huitième de degré.  En son milieu, le jour dure quatorze heures et
demi, et la latitude est trente-six degrés, un cinquième et un sixième
de degré.\footnote{Début $33;37,30$. Milieu $36;22$.}

\emph{Le cinquième climat} commence là où le jour dure quatorze heures
trois quarts et où la latitude est trente-neuf degré moins un dixième
de degré. En son milieu, le jour dure quinze heures, et la latitude
est quarante-et-un degrés et un quart de degré.\footnote{Début
  $38;54$.  Milieu $41;15$.}

\emph{Le sixième climat} commence là où le jour dure quinze heures et
quart et où la latitude est quarante-trois degrés, un quart et un
huitième de degré.  En son milieu, le jour dure quinze heures et demi,
et la latitude est quarante-cinq degrés, un quart et un dixième de
degré.\footnote{Début $43;22,30$. Milieu $45;21$.}

\emph{Le septième climat} commence là où le jour dure quinze heures
trois quarts et où la latitude est quarante-sept degrés et un
cinquième de degré.  En son milieu, le jour dure seize heures et
quart, et la latitude est cinquante degrés et un tiers de
degré.\footnote{Début $47;12$. Milieu $50;20$.}

La fin de chaque climat est le commencement de celui qui le
suit. D'autres gens posent le commencement du premier climat à
l'équateur, et la fin du septième là où s'achève la partie habitée,
non là où s'en achève la majeure partie.

\newpage\phantomsection
\index{BFBGAQ@\RL{nhr}!BFBGAQ@\RL{nahr j anhAr}, rivière}
\index{ALAHBD@\RL{jbl}!ALAHBD@\RL{jabal j jibAl}, montagne}
\index{BDBHBF@\RL{lwn}, couleur}
\index{BFBGAQ@\RL{nhr}!BBBHAS BFBGAGAQ@\RL{qws al-nhAr}, arc diurne, durée du jour}
\index{ANBDAL@\RL{_hlj}!ANBDBJAL BAAGAQAS@\RL{_halIj fAris}, golfe Persique}
\label{var13}
\includepdf[pages=5,pagecommand={\thispagestyle{plain}}]{edit2.pdf}

[figure des sept climats contenant les montagnes et les grandes
  rivières]\label{carte_fleuves_montagnes}

\emph{Dans le premier climat} il y a vingt montagnes et trente
rivières, et la plupart de ses peuples sont noirs. \emph{Dans le
  deuxième climat} il y a vingt-sept montagnes et vingt-sept rivières,
et les couleurs de ses peuples sont noire et brune. \emph{Dans le
  troisième climat}, il y a trente-trois montagnes et vingt-deux
rivières, et la plupart de ses peuples sont bruns.
\emph{Dans le quatrième climat}, il y a vingt-cinq montagnes
et vingt-deux rivières, et la couleur de la plupart de ses peuples est
entre brun et blanc ; ce sont les hommes les plus justes
\uwave{khlqan w khlqan}, et c'est pourquoi [ce climat] a été une mine
\uwave{de peuples}, de saints et de savants ; d'ailleurs les
peuples du troisième, du cinquième et des autres climats manquent
\uwave{khlqan w khlqan}. \emph{Dans le cinquième climat}, il y a
trente montagnes et quinze rivières, et la plupart de ses peuples sont
blancs.  \emph{Dans le sixième climat}, il y a onze montagnes et
quarante rivières, et souvent la couleur de ses peuples est
jaune. \emph{Dans le septième climat}, il y a autant de rivières et de
montagnes, et la couleur de ses peuples est entre jaune et blanc ; la
plus grande part en est inoccupée à cause de l'intensité du froid et
de la fréquence des neiges et des pluies, et les peuples de certaines
de ses contrées habitent dans des bains [chauds] pendant six mois de
l'année.  Entre la fin du septième climat et la fin de la partie
habitée, il y a beaucoup moins d'habitations qu'en amont ; là, le
peuple ressemble à des bêtes.

Le jour le plus long dure 17 heures quand la latitude est 54 degrés
environ.  Il dure 18 heures à 58 degrés de latitude. Il dure 19 heures
à 61 degrés de latitude. Il dure 20 heures à 63 degrés de latitude. Il
dure 21 heures à $64;30$ degrés de latitude. Il dure 22 heures à
environ 65 degrés de latitude. Il dure 23 heures à 66 degrés de
latitude. Il dure 24 heures quand la latitude atteint le complément de
l'inclinaison [de l'écliptique]. Il dure un mois à 67 degrés de
latitude. Il dure deux mois à $69;45$ degrés de latitude.  Il dure
trois mois à $73;30$ degrés de latitude. Il dure quatre mois à $78;30$
degrés de latitude. Il dure cinq mois à 84 degrés de latitude. Il dure
six mois quand la latitude atteint un quart de circonférence.

Passons aux pays connus des Anciens dans chaque climat, en allant de
l'extrémité Est vers l'extrémité Ouest. \emph{Le premier climat}
commence à l'Est aux confins de la Chine et il passe par la partie Sud
de la Chine où il y a la cité impériale, puis par le Sud de l'Inde et
du Sind, \uwave{l'île Krk Baha'ir}, le royaume du Yemen puis le golfe
Persique et la péninsule arabique qui fait à peu près cinq cents
stades et où se trouve tout l'état arabe, ses tribus et ses peuplades
près du Hedjaz, Taëf, le Yemen, Bahrein, Nejd, Tihama, la Mecque, et
Médine ou Yathrib, mais toutes ne sont pas dans le premier
climat. Parmi les villes célèbres, celles qui s'y trouvent sont la
ville de Zafar, Oman, Hadramaout, Aden, Ta`izz, Sanaa du Yemen,
M\=ar\=a, Zabid, Qalhât, Shihr, Huras, Mahra, et Saba. 

\newpage\phantomsection
\index{AMAHAT@\RL{.hb^s}!AMAHATAI@\RL{.hb^saT}, \'Ethiopie}
\index{BFBJBD@\RL{nIl}, Nil}
\index{AMBHAW@\RL{.hw.t}!AHAMAQ BEAMBJAW@\RL{al-ba.hr al-mu.hI.t}, l'Océan}
\includepdf[pages=6,pagecommand={\thispagestyle{plain}}]{edit2.pdf}

\noindent Puis [le premier climat] passe par l'Ethiopie et coupe le
Nil d'Egypte, le Soudan et la Nubie, \uwave{comme \d{H}arr\=a}, palais
du roi d'Ethiopie, et \uwave{Dahlamat} une ville de Nubie, et le Ghana
d'où vient l'or du Soudan, dans le Maghreb, puis il passe par le Sud
de la Berbérie, jusqu'à l'Océan à l'Ouest.

\emph{Le deuxième climat} commence à l'Est de la Chine, puis il
traverse la majeure partie de l'Inde, le Sind, comme la ville de
Mansura, Narwar, et \uwave{Danbek}, et
il arrive à Oman et Bassorah. Il coupe la péninsule arabique sur les
terres de Nejd et Tihama ; il compte, parmi les villes
de là-bas : Al-Yam\=ama, Bahrein, Hajar, Yathrib, Harr, La Mecque,
Taëf et Jeddah. Il coupe la mer de Qulzum\footnote{\textit{i. e.} la
  mer Rouge.} et passe par la Haute-Egypte, puis il coupe le Nil~;
il compte, parmi les villes de là-bas, Qûs, Akhmîm, \uwave{S\=abir},
\uwave{\d{S}\=a}, et Assouan, puis il prend sur les terres du
Maghreb et passe au c{\oe}ur de l'Afrique, puis par la Berbérie,
pour finir à l'Océan.

\emph{Le troisième climat} commence à l'Est et passe par la Chine, par
le Royaume Indien et le Gandh\=ara qui est dans les hautes-terres de
l'Inde, par Multan dans le Sind, par \uwave{Z\=a'il}, Bust,
S{\=\i}statm\=an\footnote{Sistan ?}, Kerman, Sirjan, Fars, Ispahan,
Ahvaz, Askar, Koufa, Bassorah, Wasit, Baghdad, Al-Anbâr, Hit ; puis il
passe par la Syrie, où sont les villes de \uwave{\d{H}iy\=ar},
Salamiya, \d{H}om\d{s}, Hama, Baalbek, Maarret en Nouman, Damas, Tyr,
Acre, Tibériade, Qays\=ariya\footnote{aujourd'hui Baniyas.},
Jérusalem, Ascalon, \uwave{Madiyan}, 'Ars\=uf\footnote{Apollonie.},
Ramla, Gaza, puis il passe par la Basse-Egypte où sont les villes de
Farama\footnote{Péluse.}, Tunis, Damiette, Le Caire, Alexandrie ; puis
il passe en Ifriqiya où sont les villes de Kerouan et Sousse, puis il
passe par les peuplades de la Berbérie, puis par Tanger, pour finir à
l'Océan.

\newpage\phantomsection
\index{AMBHAW@\RL{.hw.t}!AHAMAQ BEAMBJAW@\RL{al-ba.hr al-mu.hI.t}, l'Océan}
\includepdf[pages=7,pagecommand={\thispagestyle{plain}}]{edit2.pdf}

\emph{Le quatrième climat} commence au Nord de la Chine, il passe par
\uwave{Batout}, \uwave{Khoulkh{\=\i}z}, Cathay, Khotan, les montagnes
du Cachemire, Bul\=ur, Badakhsh\=an, Kaboul, Balkh, Hérât, Fergh\=ana,
Samarcande, Boukhara, \=Amul\footnote{aujourd'hui Türkmenabat au
  Turkmenistan.}, Merv, Sarakhs, Tus, Nishapur, Gorgan, Esfarayen,
Qûhistân, \uwave{Q\=umz}, Daylam, Tabaristan, Qom, Hamadan,
Azerbaïdjan, Qazvin, Dinavar, Hulwân, Shahraz\=ur, Mossoul,
S\=amarr\=a, Nusaybin, \=Amid\footnote{aujourd'hui Diyarbak{\i}r en
  Turquie.}, Mardin, Ra's al-`Ayn, Samosate, Malatya, Alep, Qinnasrîn,
Antioche, Tripoli\footnote{Tripoli du Liban.}, Tarse,
Ammuriye\footnote{la ville antique d'Amorium, aujourd'hui Hisarköy en
  Turquie.}, Lattaquié, puis il passe par l'île de Chypre, par Rhodes,
et par les terres occidentales de l'Europe, et par Tanger, pour finir
à l'Océan au détroit entre Al-Andalous et le Maghreb.

\emph{Le cinquième climat} commence à l'Est à la frontière du
Turkestan, près de Gog et Magog ; il passe par les célèbres
  nations du Turkestan avec ses peuplades, jusqu'à la limite de
Kachgar\footnote{dans le Xinjiang.}, Taraz, 'Isf{\=\i}j\=ab\footnote{
  aujourd'hui Sayram, au Kazakhstan.}, \uwave{\d{H}a\d{h}},
Usrushana\footnote{aujourd'hui Istaravchan au Tadjikistan.}, Ferghana,
Samarcande, Soghd, Boukhara, Khwarezm, \uwave{Shams}, la contrée
d'Arménie, Barda, Mayy\=af\=ariq{\=\i}n\footnote{Silvan, en Turquie.},
la Séleucie,
Erzurum, Ahlat\footnote{au bord du lac de Van.}, et il passe dans
l'Empire Romain par \uwave{\d{H}arshafat}, \uwave{Qrt}, et
\uwave{Rome}, puis par les côtes de la mer de
Syrie\footnote{\textit{i. e.} la mer Méditerranée.}, l'Italie et
Al-Andalous, jusqu'à l'Océan.

\emph{Le sixième climat} commence à l'Est par les terres des Turcs
orientaux et leurs peuplades, il traverse la mer de Gorgan\footnote{la
  mer Caspienne.}, il passe par les Khazars et les \uwave{M\=ut\=ani ?},
par les Slaves\footnote{D'après Kazimirski, nom générique donné aux
  peuples du Nord-Est de l'Europe.}, les \uwave{Aryens ?}, B\=ab
al-Abw\=ab\footnote{aujourd'hui Derbent en Russie.}, la Russie, puis
l'Empire Romain dont Constantinople, et le Nord d'Al-Andalous, pour
finir à l'Océan.

\newpage\phantomsection
\index{AMBHAW@\RL{.hw.t}!AHAMAQ BEAMBJAW@\RL{al-ba.hr al-mu.hI.t}, l'Océan}
\includepdf[pages=8,pagecommand={\thispagestyle{plain}}]{edit2.pdf}

\emph{Le septième climat} commence à l'Est en passant par l'extrémité
des Turcs orientaux et de leurs peuplades, par le Nord (Gog et Magog),
puis il passe par Ghiy\=a\d{d} où des montagnes recèlent des turques
vivant comme des bêtes. Puis il passe par Bolghar, la Russie et les
Slaves, et il traverse la Mer Méditerranée et les pays slaves pour
finir à l'Océan.

J'ai placé ces pays dans leurs climats, de manière détaillée, selon leur
longitude et leur latitude, dans les tables de la troisième partie,
avec les autres tables déjà mentionnées. Dieu approuve la raison.
\label{troisieme2}

\newpage\phantomsection
\index{ASBHBI@\RL{sw_A}!ANAWAW ASAJBHAGAB@\RL{_ha.t.t al-istiwA'}, équateur terrestre}
\index{AYAOBD@\RL{`dl}!BEAYAOAOBD BFBGAGAQ@\RL{m`ddl al-nhAr}, équateur, plan de l'équateur}
\index{AMAQBC@\RL{.hrk}!AMAQBCAI AKAGBFBJAI@\RL{.hrkaT _tAnyaT}, deuxième mouvement}
\index{AXBGAQ@\RL{.zhr}!AXBGBHAQ@\RL{.zuhUr}, visibilité ($\neq$ \RL{i_htifA'})}
\index{ASBFBH@\RL{snw}!BAAUBHBD ASBFAI@\RL{fu.sUl al-sanaT}, les saisons}
\index{BABHBB@\RL{fwq}!ADBABB@\RL{'ufuq}, horizon}
\index{AHAQAL@\RL{brj}!BEBFAWBBAI AHAQBHAL@\RL{min.taqaT al-burUj}, écliptique, ceinture de l'écliptique}
\includepdf[pages=9,pagecommand={\thispagestyle{plain}}]{edit2.pdf}
\addcontentsline{toc}{chapter}{II.2 Particularités de l'équateur terrestre}
\begin{center}
  \Large Deuxième section
  
  \normalsize Particularités de l'équateur terrestre
\end{center}
\noindent L'équateur passe par les zéniths des lieux qui se trouvent
sur l'équateur terrestre, et il coupe leurs horizons orthogonalement,
et leurs cercles horizontaux bissectent chaque trajectoire diurne
puisqu'ils passent par les pôles de l'équateur. C'est pourquoi en ces
lieux, le jour et la nuit sont toujours égaux. L'origine de leurs
azimuts est l'équateur, et ses pôles sont aux points Nord et Sud du
cercle horizontal. Ainsi tous les astres ont un lever et un coucher,
sauf [ceux qui sont] aux pôles, car une moitié [de leur trajectoire
  diurne] est visible et l'autre cachée. La durée de visibilité de
chaque point de l'orbe est égale à la durée pendant laquelle il est
caché ; il y a certes un écart à cause de l'irrégularité des
trajectoires due au second mouvement\footnote{Le second mouvement est
  le mouvement du huitième orbe, c'est-à-dire le mouvement de
  précession des fixes.}, mais cet écart est insensible.

Le Soleil passe deux fois l'an au zénith, et c'est quand il est aux
points des équinoxes ; il n'y alors, à midi, pas d'ombre projetée
horizontale. Il ne s'éloigne pas du zénith d'un angle supérieur à
l'inclinaison [de l'écliptique], et sa hauteur n'est jamais inférieure
au complément de l'inclinaison ; il passe la moitié de l'année de
chaque côté, et l'ombre à midi est alors située de l'autre côté. Les
ombres au début de l'été et au début de l'hiver sont égales. Les pôles
de l'écliptique sont à l'horizon chaque fois qu'un des deux points des
équinoxes est au zénith, la ceinture de l'écliptique est alors
orthogonale à l'horizon, et le méridien du lieu bissecte la partie
visible de l'écliptique. Pendant que la moitié Nord de la ceinture de
l'écliptique passe par le méridien du lieu, c'est son pôle Sud qui est
visible. Pendant que la moitié Sud de la ceinture y passe, c'est son
pôle Nord qui est visible. La hauteur de chacun ne dépasse pas
l'inclinaison [de l'écliptique].

Ils ont, dans chaque année, huit saisons : deux printemps, deux étés,
deux automnes et deux hivers. Le début du premier printemps, c'est
quand le Soleil est au milieu du Verseau, et le début du second, c'est
quand il est au milieu du Lion ; le début du premier été, quand il est
au commencement du Bélier, et le début du second quand il est au
commencement de la Balance ; le début du premier automne quand il est
au milieu du Taureau, et le début du second quand il est au milieu du
Scorpion ; le début du la premier hiver quand il est au commencement
du Cancer et le début du second quand il est au commencement du
Capricorne ; comme il est montré dans la figure suivante. Dieu est le
plus savant.

\newpage\phantomsection
\begin{center}
\begingroup%
  \makeatletter%
  \providecommand\color[2][]{%
    \errmessage{(Inkscape) Color is used for the text in Inkscape, but the package 'color.sty' is not loaded}%
    \renewcommand\color[2][]{}%
  }%
  \providecommand\transparent[1]{%
    \errmessage{(Inkscape) Transparency is used (non-zero) for the text in Inkscape, but the package 'transparent.sty' is not loaded}%
    \renewcommand\transparent[1]{}%
  }%
  \providecommand\rotatebox[2]{#2}%
  \ifx\svgwidth\undefined%
    \setlength{\unitlength}{405.85902868bp}%
    \ifx\svgscale\undefined%
      \relax%
    \else%
      \setlength{\unitlength}{\unitlength * \real{\svgscale}}%
    \fi%
  \else%
    \setlength{\unitlength}{\svgwidth}%
  \fi%
  \global\let\svgwidth\undefined%
  \global\let\svgscale\undefined%
  \makeatother%
  \begin{picture}(1,0.75749861)%
    \put(0,0){\includegraphics[width=\unitlength]{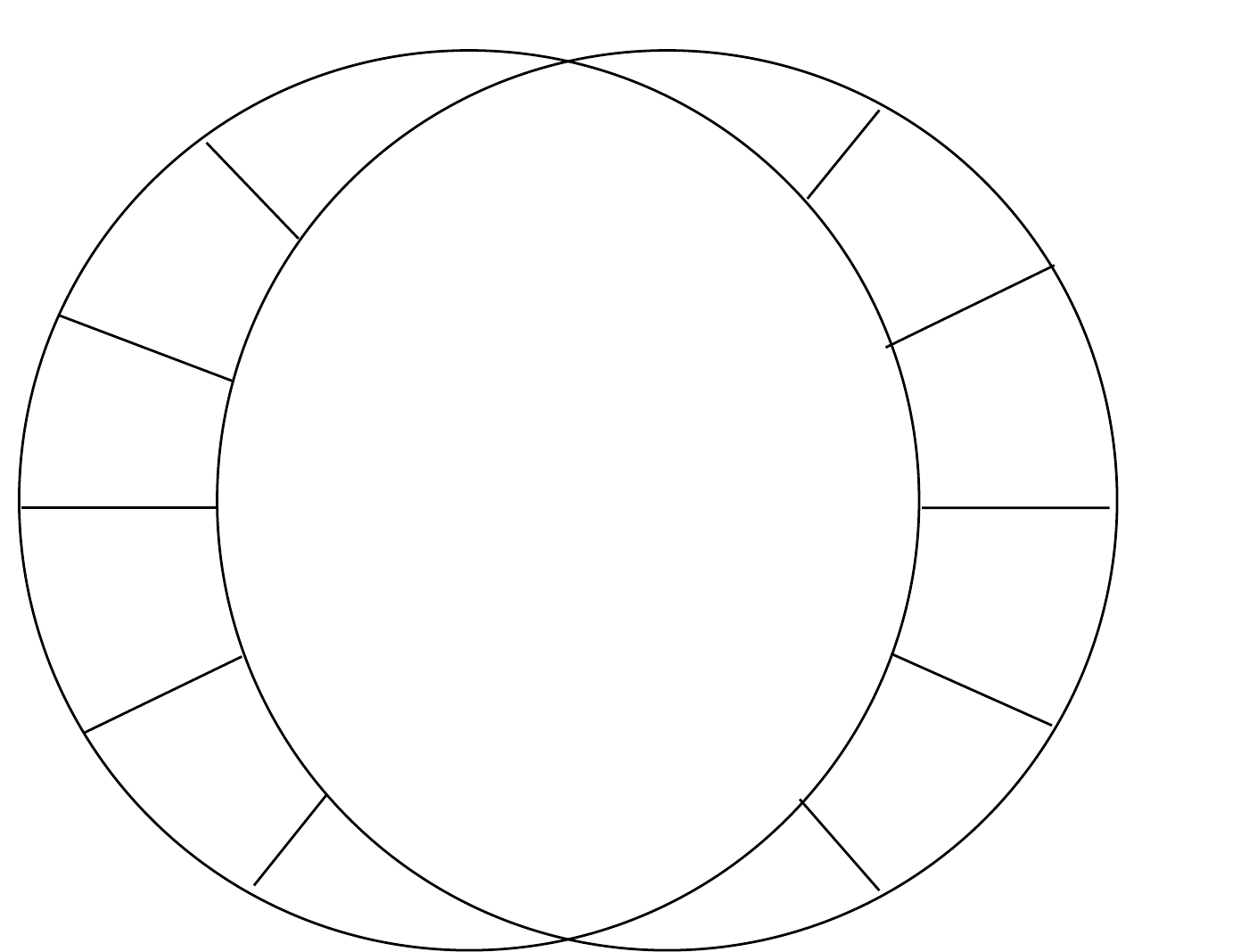}}%
    \put(0.66880868,0.66980252){\color[rgb]{0,0,0}\makebox(0,0)[rb]{\smash{\RL{snblaT}}}}%
    \put(0.74068501,0.57920566){\color[rgb]{0,0,0}\makebox(0,0)[rb]{\smash{\RL{asad}}}}%
    \put(0.81689766,0.42688346){\color[rgb]{0,0,0}\makebox(0,0)[rb]{\smash{\RL{sr.tAn}}}}%
    \put(0.83229938,0.28551255){\color[rgb]{0,0,0}\makebox(0,0)[rb]{\smash{ \RL{jUzA}}}}%
    \put(0.76055365,0.13617713){\color[rgb]{0,0,0}\makebox(0,0)[rb]{\smash{\RL{_tUr}}}}%
    \put(0.61050744,0.03064672){\color[rgb]{0,0,0}\makebox(0,0)[rb]{\smash{\RL{.hml}}}}%
    \put(0.24597591,0.04856697){\color[rgb]{0,0,0}\makebox(0,0)[lb]{\smash{\RL{.hwt}}}}%
    \put(0.13845434,0.14414164){\color[rgb]{0,0,0}\makebox(0,0)[lb]{\smash{\RL{dlw}}}}%
    \put(0.05283535,0.27157462){\color[rgb]{0,0,0}\makebox(0,0)[lb]{\smash{\RL{jdy}}}}%
    \put(0.03690623,0.40697211){\color[rgb]{0,0,0}\makebox(0,0)[lb]{\smash{\RL{qws}}}}%
    \put(0.11655181,0.54037845){\color[rgb]{0,0,0}\makebox(0,0)[lb]{\smash{\RL{`qrb}}}}%
    \put(0.23402902,0.65984685){\color[rgb]{0,0,0}\makebox(0,0)[lb]{\smash{\RL{mIzAn}}}}%
    \put(0.89188624,0.19936897){\color[rgb]{0,0,0}\rotatebox{68.45049768}{\makebox(0,0)[lb]{\smash{\RL{mn.tqaT al-burUj al-^simAliy}}}}}%
    \put(0.0118382,0.49193955){\color[rgb]{0,0,0}\rotatebox{59.99999989}{\makebox(0,0)[lb]{\smash{\RL{qws m`dl al-nhAr al-junUbiyaT}}}}}%
    \put(0.45814362,0.69553835){\color[rgb]{0,0,0}\rotatebox{90}{\makebox(0,0)[rb]{\smash{\RL{mabdA al-.sayf}}}}}%
    \put(0.67389382,0.54809704){\color[rgb]{0,0,0}\rotatebox{36.81728888}{\makebox(0,0)[rb]{\smash{\RL{mabdA al-rbI`}}}}}%
    \put(0.723188,0.34823337){\color[rgb]{0,0,0}\makebox(0,0)[rb]{\smash{\RL{mabdA al-^sitA}}}}%
    \put(0.6724752,0.18582849){\color[rgb]{0,0,0}\rotatebox{-36.09699436}{\makebox(0,0)[rb]{\smash{\RL{mbdA al-_harIf}}}}}%
    \put(0.45615246,0.02702749){\color[rgb]{0,0,0}\rotatebox{90}{\makebox(0,0)[lb]{\smash{\RL{mbdA al-.sayf}}}}}%
    \put(0.23792439,0.17717288){\color[rgb]{0,0,0}\rotatebox{38.94088684}{\makebox(0,0)[lb]{\smash{\RL{mbdA al-rabI`}}}}}%
    \put(0.18425052,0.34922904){\color[rgb]{0,0,0}\makebox(0,0)[lb]{\smash{\RL{mbdA al-^sitA}}}}%
    \put(0.21599137,0.51168428){\color[rgb]{0,0,0}\rotatebox{-32.14993087}{\makebox(0,0)[lb]{\smash{\RL{mbdA al-_harIf}}}}}%
  \end{picture}%
\endgroup%

\end{center}

\newpage\phantomsection
\begin{center}
\begingroup%
  \makeatletter%
  \providecommand\color[2][]{%
    \errmessage{(Inkscape) Color is used for the text in Inkscape, but the package 'color.sty' is not loaded}%
    \renewcommand\color[2][]{}%
  }%
  \providecommand\transparent[1]{%
    \errmessage{(Inkscape) Transparency is used (non-zero) for the text in Inkscape, but the package 'transparent.sty' is not loaded}%
    \renewcommand\transparent[1]{}%
  }%
  \providecommand\rotatebox[2]{#2}%
  \ifx\svgwidth\undefined%
    \setlength{\unitlength}{407.46595893bp}%
    \ifx\svgscale\undefined%
      \relax%
    \else%
      \setlength{\unitlength}{\unitlength * \real{\svgscale}}%
    \fi%
  \else%
    \setlength{\unitlength}{\svgwidth}%
  \fi%
  \global\let\svgwidth\undefined%
  \global\let\svgscale\undefined%
  \makeatother%
  \begin{picture}(1,0.71556368)%
    \put(0,0){\includegraphics[width=\unitlength]{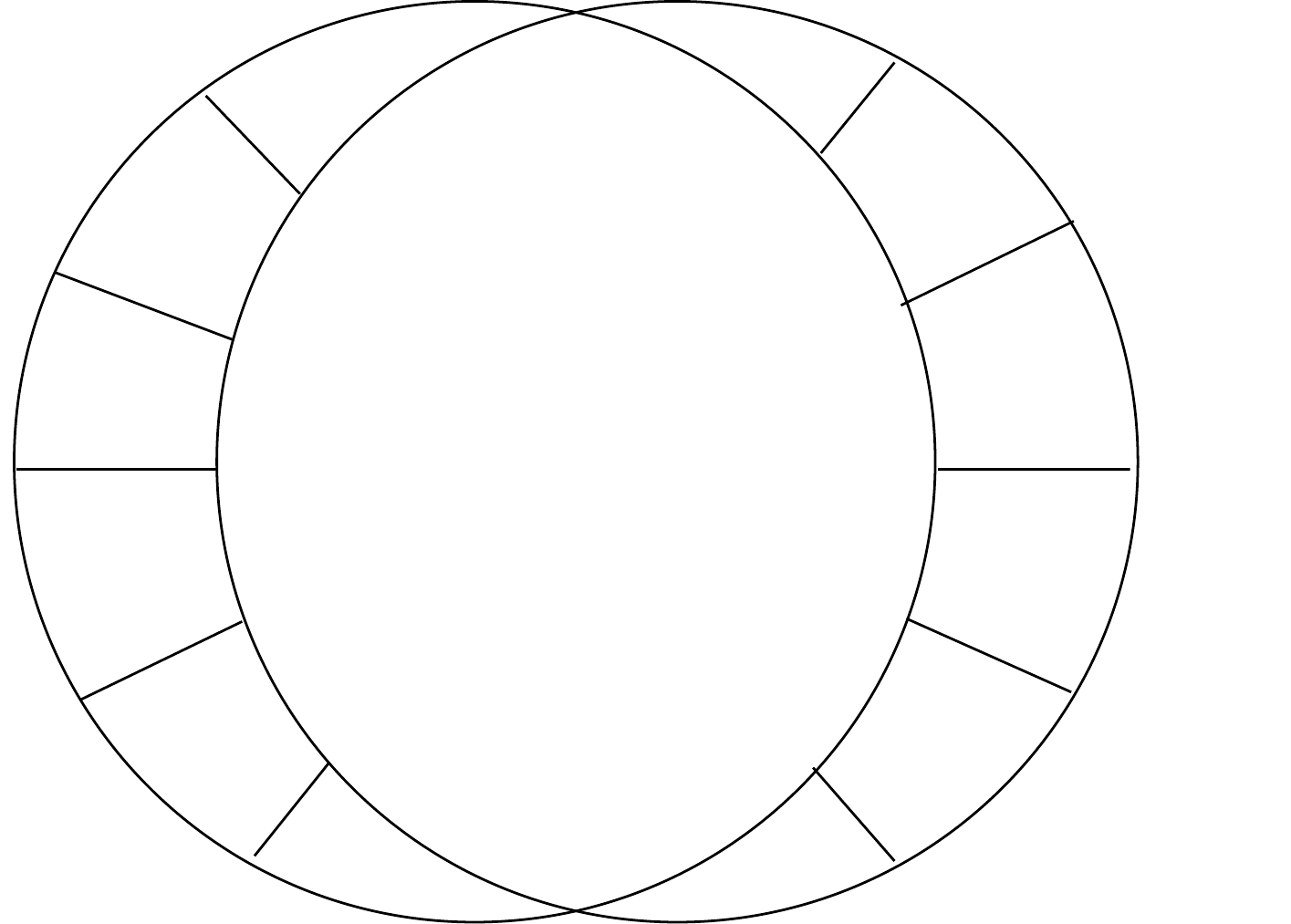}}%
    \put(0.55218574,0.67707742){\color[rgb]{0,0,0}\makebox(0,0)[lb]{\smash{Vierge}}}%
    \put(0.66919968,0.57592981){\color[rgb]{0,0,0}\makebox(0,0)[lb]{\smash{Lion}}}%
    \put(0.75448102,0.42519995){\color[rgb]{0,0,0}\makebox(0,0)[lb]{\smash{Cancer}}}%
    \put(0.75051441,0.28438657){\color[rgb]{0,0,0}\makebox(0,0)[lb]{\smash{Gémeaux}}}%
    \put(0.69101587,0.13564008){\color[rgb]{0,0,0}\makebox(0,0)[lb]{\smash{Taureau}}}%
    \put(0.55615229,0.03052586){\color[rgb]{0,0,0}\makebox(0,0)[lb]{\smash{Bélier}}}%
    \put(0.22296009,0.04440882){\color[rgb]{0,0,0}\makebox(0,0)[lb]{\smash{Poisson}}}%
    \put(0.11982912,0.14357319){\color[rgb]{0,0,0}\makebox(0,0)[lb]{\smash{Verseau}}}%
    \put(0.0325645,0.27050361){\color[rgb]{0,0,0}\makebox(0,0)[lb]{\smash{Capricorne}}}%
    \put(0.0325645,0.40536713){\color[rgb]{0,0,0}\makebox(0,0)[lb]{\smash{Sagittaire}}}%
    \put(0.09007982,0.5422139){\color[rgb]{0,0,0}\makebox(0,0)[lb]{\smash{Scorpion}}}%
    \put(0.22890989,0.6572446){\color[rgb]{0,0,0}\makebox(0,0)[lb]{\smash{Balance}}}%
    \put(0.88417271,0.19858272){\color[rgb]{0,0,0}\rotatebox{68.45049768}{\makebox(0,0)[lb]{\smash{partie Nord de la ceinture de l'écliptique}}}}%
    \put(0.00759534,0.48999948){\color[rgb]{0,0,0}\rotatebox{59.99999989}{\makebox(0,0)[lb]{\smash{un arc de l'équateur}}}}%
    \put(0.45214065,0.68679258){\color[rgb]{0,0,0}\rotatebox{90}{\makebox(0,0)[rb]{\smash{début de l'été}}}}%
    \put(0.6626638,0.53998528){\color[rgb]{0,0,0}\rotatebox{36.81728888}{\makebox(0,0)[rb]{\smash{début du printemps}}}}%
    \put(0.71057422,0.34586841){\color[rgb]{0,0,0}\makebox(0,0)[rb]{\smash{début de l'hiver}}}%
    \put(0.66562697,0.18509563){\color[rgb]{0,0,0}\rotatebox{-36.09699436}{\makebox(0,0)[rb]{\smash{début de l'automne}}}}%
    \put(0.45015735,0.0269209){\color[rgb]{0,0,0}\rotatebox{90}{\makebox(0,0)[lb]{\smash{début de l'été}}}}%
    \put(0.2327899,0.17647416){\color[rgb]{0,0,0}\rotatebox{38.94088684}{\makebox(0,0)[lb]{\smash{début du printemps}}}}%
    \put(0.17932771,0.34785178){\color[rgb]{0,0,0}\makebox(0,0)[lb]{\smash{début de l'hiver}}}%
    \put(0.21094338,0.50966634){\color[rgb]{0,0,0}\rotatebox{-32.14993087}{\makebox(0,0)[lb]{\smash{début de l'automne}}}}%
  \end{picture}%
\endgroup%

\end{center}

\newpage\phantomsection
\index{ACBHAL@\RL{'awj}!ACBHAL ATBEAS@\RL{'awj al-^sams}, Apogée du Soleil}
\index{AMAQAQ@\RL{.harara}!AMAQAQ@\RL{.hrr}, chaud}
\index{AHAQAO@\RL{brd}!AHAQAO@\RL{brd}, froid}
\index{AYAOBD@\RL{`dl}!AYAOBD@\RL{`dl}, tempéré}
\index{BHASAY@\RL{ws`}!ASAYAI@\RL{s`aT m^srq / m.grb}, amplitude Est / Ouest}
\index{AHBF ASBJBFAG@\RL{ibn sInA}, Avicenne}
\index{AQAGARBJ@\RL{al-rAzI}, Fa\b{h}r al-D{\=\i}n al-R\=az{\=\i}}
\index{AMAHAT@\RL{.hb^s}!AMAHATAI@\RL{.hb^saT}, \'Ethiopie}
\index{BABDBC@\RL{flk}!BABDBC BEASAJBBBJBE BCAQAI BEBFAJAUAHAI@\RL{falak mustaqIm, kuraT munta.sibaT}, sphère droite}
\index{ASBEAJ@\RL{smt}!ASBEAJ AQABAS@\RL{samt al-ra's}, zénith}
\index{BBBDAH@\RL{qlb}!BEAOAGAQ BEBFBBBDAHBJBF@\RL{madAr al-munqalbayn}, les deux tropiques}
\index{ATAYAY@\RL{^s``}!ATAYAGAY@\RL{^su`A`}, rayon (de lumière)}
\includepdf[pages=10,pagecommand={\thispagestyle{plain}}]{edit2.pdf}

Là-bas, l'orbe tourne comme une roue puisque l'horizon coupe toutes
les trajectoires diurnes orthogonalement. C'est pourquoi on appelle
ses horizons <<~horizons de l'orbe droit~>> ou <<~horizons de la
sphère droite~>>. Comme leur horizon est un des cercles de déclinaison
(il passe par les pôles), alors l'\emph{amplitude Est} de tout
point\footnote{L'\emph{amplitude Est d'un astre} mesure donc, à
  l'horizon, l'arc entre le point où se lève les équinoxes et le point
  où se lève l'astre. L'amplitude Ouest est relative aux
  couchers. Ragep : \textit{ortive amplitude}, \textit{occasive
    amplitude}.}, c'est-à-dire l'arc à l'horizon entre son lever et le
lever des équinoxes, est égal à sa déclinaison. De même pour son
amplitude Ouest.


Les gens s'accordent à dire que les lieux les plus chauds en été sont
ceux qui sont sous le tropique, si aucune raison terrestre ou
atmosphérique n'enlève rien à sa chaleur, car le Soleil arrive au
zénith et y prolonge son séjour pendant à peu près deux mois, à cause
de l'amoindrissement de l'écart en déclinaison.
Pour cette raison, il semble alors ne pas avoir de mouvement en
déclinaison pendant des jours. Le jour en été ne dure pas plus
longtemps, ni la nuit moins longtemps, mais la longueur du jour
n'influe pas sur l'intensité de la chaleur quand le Soleil est proche
du zénith, à cause de l'épaississement des rayons quand cela arrive,
épaississement dû au fait qu'ils forment un angle aigu.

Le Grand Maître\footnote{Avicenne} conteste cela et pense que les
contrées les plus tempérées sont à l'équateur terrestre, parce que le
Soleil n'y reste pas longtemps au zénith, parce que son mouvement en
déclinaison est y est le plus rapide. Bien que le passage au zénith
cause l'échauffement, c'est la durée du passage au zénith qui y
conduit. C'est ainsi que, [en général], l'été est plus chaud que le
printemps, et l'après-midi plus chaud qu'avant midi. De plus, à cause
de l'égalité entre la durée du jour et la durée de la nuit, la dureté
des conditions de chacun sera rapidement rompue par l'autre ; ainsi
les deux seront tempérés.

L'Imam\footnote{Al-R\=az{\=\i}} réfute cela, en disant que
l'échauffement du Soleil en hiver, à l'équateur terrestre, est comme
l'échauffement du Soleil en été dans un pays dont la latitude est le
double de l'inclinaison de l'écliptique, et que cet échauffement est
très intense.  Que penser alors de la chaleur estivale à l'équateur
terrestre où le Soleil peut être considéré au zénith toute l'année !

On réfute l'Imam en contredisant le fait que la chaleur hivernale à
l'équateur est comme la chaleur estivale dans ce pays~: la chaleur
estivale de ce pays serait plus intense, à cause de la longueur du
jour et de la courte durée de ses nuits, au contraire de ce qui se
passe à l'équateur terrestre. De plus les choses habituelles n'ont pas
autant d'influence que les choses inhabituelles.

\newpage\phantomsection
\index{BDBHBF@\RL{lwn}, couleur}
\includepdf[pages=11,pagecommand={\thispagestyle{plain}}]{edit2.pdf}
\noindent \`A cause de
l'accoutumance de la constitution des peuples de l'équateur à la
chaleur, et de leur manque d'accoutumance au froid, il est donc
possible qu'ils ne ressentent pas la chaleur de l'air quand le Soleil
est au zénith, ni sa froideur quand il est aux solstices,
contrairement au peuple de cet [autre] pays aux mêmes moments.

L'Imam juge que le climat le plus tempéré est le quatrième. Il le
déduit du fait qu'un lieu où prospèrent les habitations se trouve sans
aucun doute plus proche du milieu des sept climats que de leurs
extrémités. Ainsi, le fait de griller, ou celui de rester immature,
causés par les conditions de chaleur ou de froid, sont clairement
visibles aux extrémités [des sept climats].

La vérité est que, si l'on entend par \emph{tempéré} l'uniformité des
états, il n'y a aucun doute qu'elle est atteinte davantage à
l'équateur terrestre ; mais si l'on entend par là l'équilibre entre
les deux conditions, alors il n'y a aucun doute qu'il est atteint
davantage dans le quatrième climat. Le montrent les faits suivants. La
couleur noire des habitants de l'équateur terrestre dans les Pays des
Noirs et en \'Ethiopie est très intense, et leurs cheveux sont très
crépus~: c'est causé par la chaleur de l'air. Le contraire a lieu dans
le quatrième climat~; car l'air y est tempéré.

\newpage\phantomsection
\index{BABHBB@\RL{fwq}!ADBABB@\RL{'ufuq}, horizon}
\index{BABHBB@\RL{fwq}!ADBABB BEAGAEBD@\RL{'ufuq mA'il}, horizon incliné}
\index{AHAQAL@\RL{brj}!BEBFAWBBAI AHAQBHAL@\RL{min.taqaT al-burUj}, écliptique, ceinture de l'écliptique}
\includepdf[pages=12,pagecommand={\thispagestyle{plain}}]{edit2.pdf}
\addcontentsline{toc}{chapter}{II.3 Particularités des lieux qui ont une latitude [non nulle]}
\begin{center}
  \Large Troisième section

  \normalsize Particularités des lieux qui ont une latitude [non nulle]
\end{center}
On les appelle \emph{horizons inclinés}~; ce sont ceux qui sont entre
l'équateur terrestre et l'un des pôles de l'équateur. Là-bas, l'orbe
tourne [sous un angle égal à l'éloignement] de l'équateur par rapport
à ces horizons, du côté du pôle caché. Chacun de ces horizons est
incliné par rapport à [l'équateur] du côté du pôle visible. C'est
pourquoi on les appelle horizons inclinés. La hauteur du pôle est du
même côté, et elle est égale à la latitude du pays.

La distance des trajectoires diurnes toujours visibles ou toujours
cachées, par rapport à l'équateur, est supérieure ou égale au
complément de la latitude du pays. La plus grande d'entre elles est
celle dont la distance est égale au complément de la latitude, elle
touche l'horizon, et elle n'est jamais invisible, ou bien jamais
visible. Les trajectoires diurnes sont coupées par l'horizon en deux
arcs inégaux. Le plus grand [arc] est celui qui est visible, quand il
s'agit d'une trajectoire proche du pôle visible et située du même côté
que lui ; c'est celui qui est caché, quand il s'agit d'une trajectoire
proche du pôle caché et située du même côté que lui. Pour tout couple
de trajectoires diurnes équidistantes de l'équateur mais situées de
part et d'autre de l'équateur, les arcs alternés sont deux à deux
égaux\footnote{L'absence du duel en français rend cette phrase un peu
  obscure. L'arc visible de l'une des deux trajectoire est égal à
  l'arc caché de l'autre, et inversement.}.

Plus le Soleil s'éloigne de l'équateur du côté du pôle visible, plus
l'excès du jour sur la nuit croît, et inversement du côté du pôle
caché -- mais on augmente la latitude du pays, alors l'écart entre
l'arc visible et l'arc caché [de chaque trajectoire diurne] augmente
aussi. Le jour s'allonge jusqu'à ce que [le Soleil] atteigne le sommet
du solstice qui est proche du pôle visible, et il raccourcit jusqu'à
l'autre solstice. En chaque partie [du circuit du Soleil sur
  l'écliptique], le jour est égal à la nuit en la partie opposée, et
inversement. Les deux temps ne sont égaux que lorsque le Soleil est à
l'équinoxe à son lever ou à son coucher. S'il se lève quand il est à
l'équinoxe, alors la nuit de ce lever est égale à son jour ; s'il se
couche quand il est à l'équinoxe, alors le jour où il se couche est
égal à la nuit du lever [suivant]. D'où apparaît l'impossibilité que
[le jour et la nuit] soient égaux dans toutes les contrées à une même
date, puisqu'il est impossible que [le Soleil] soit situé en même
temps sur tous les horizons.

Soit deux cercles de déclinaison, l'un à l'Est, l'autre à l'Ouest,
tels que chacun passe par l'un des deux points d'intersection de
l'horizon et de la trajectoire diurne du Soleil, d'un astre, ou d'un
point quelconque. Il y a deux triangles~: l'un à l'Est, l'autre à
l'Ouest, sous l'horizon si l'on est du côté du pôle visible, et sur
l'horizon si l'on est du côté du pôle caché.

\newpage\phantomsection
\index{BHASAY@\RL{ws`}!ASAYAI@\RL{s`aT m^srq / m.grb}, amplitude Est / Ouest} 
\index{BFBGAQ@\RL{nhr}!AJAYAOBJBD BFBGAGAQ@\RL{ta`dIl al-nahAr}, équation du jour}
\index{BFBGAQ@\RL{nhr}!BBBHAS BFBGAGAQ@\RL{qws al-nhAr}, arc diurne, durée du jour}
\index{ASBEAJ@\RL{smt}!ASBEAJ AQABAS@\RL{samt al-ra's}, zénith}
\index{ASBEAJ@\RL{smt}!ASBEAJ BBAOBE@\RL{samt al-qadam, samt al-rjl}, nadir}
\includepdf[pages=13,pagecommand={\thispagestyle{plain}}]{edit2.pdf}

\noindent Le premier côté de chacun de ces deux triangles est un arc
du cercle de déclinaison, et sa longueur est la déclinaison du Soleil,
de l'astre, ou du point quelconque. Le deuxième côté est un arc du
cercle de l'horizon, et c'est l'amplitude Est quand on est du côté
[Est], ou bien l'amplitude Ouest quand on est du côté [Ouest] ; en
effet l'amplitude Est est un arc du cercle de l'horizon entre la
trajectoire de l'astre (ou du point quelconque) et le lever des
équinoxes\footnote{\textit{i. e.} le point Est, à l'horizon.}, et
l'amplitude Ouest est l'arc entre la trajectoire mentionnée et le
coucher des équinoxes\footnote{\textit{i. e.} le point Ouest, à
  l'horizon.}. Le troisième côté est l'arc de l'équateur entre ce
lever -- ou ce coucher -- et le cercle de déclinaison passant par
l'intersection de la trajectoire et de l'horizon~: c'est
l'\emph{équation du jour} du Soleil ou de l'astre. L'excédent de l'arc
diurne (de l'astre ou du point quelconque) sur son arc diurne en
l'équateur terrestre est le double de l'équation du jour, lorsque
l'astre (ou le point quelconque) est du côté du pôle visible ; et
c'est son défaut qui est le double, lorsque l'astre est du côté du
pôle caché. Il est clair qu'on obtient la moitié de l'arc diurne en
ajoutant l'équation, pendant un quart de circonférence, puis en la
retranchant ; il est donc bien permis de l'appeler ``équation du
jour'', car c'est l'équation pour la moitié [du jour].

Il y a des gens qui posent un seul cercle de déclinaison passant par
le lever de l'équinoxe et par son coucher. Avec l'horizon et avec
toute trajectoire diurne, [ce cercle de déclinaison] forme deux
triangles, l'un à l'Est et l'autre à l'Ouest ; mais ils sont sur
l'horizon si [la trajectoire diurne] est du côté du pôle visible, et
sous l'horizon si elle est du côté du pôle caché. L'équation du jour
est alors un arc de la trajectoire de l'astre ou du point quelconque,
entre le cercle de l'horizon et le cercle de déclinaison passant par
le lever de l'équinoxe et son coucher. Cet arc de la trajectoire
diurne est semblable à l'arc de l'équateur mentionné ci-dessus.

Toute trajectoire diurne dont la distance à l'équateur est égale à la
latitude du pays touche le cercle origine des azimuts, au zénith si
elle est du côté du pôle visible, et au nadir si elle est du côté du
pôle caché. Si la distance est supérieure à la latitude du pays, alors
[la trajectoire diurne] ne rencontrera pas le cercle origine des
azimuts ; si elle est du côté du pôle visible, alors elle passera près
du zénith mais loin du nadir du côté caché. Si en revanche la distance
est inférieure [à la latitude du pays], alors [la trajectoire diurne]
coupera le cercle origine des azimuts en deux points, l'un à l'Est et
l'autre à l'Ouest. Dès lors que l'astre est sur l'arc [de trajectoire
  diurne] situé entre le cercle origine des azimuts et l'équateur, il
s'éloigne du cercle origine des azimuts quand il est du côté caché si
cette trajectoire diurne est située du côté du pôle visible, et
inversement.

\newpage\phantomsection
\begin{center}
  \small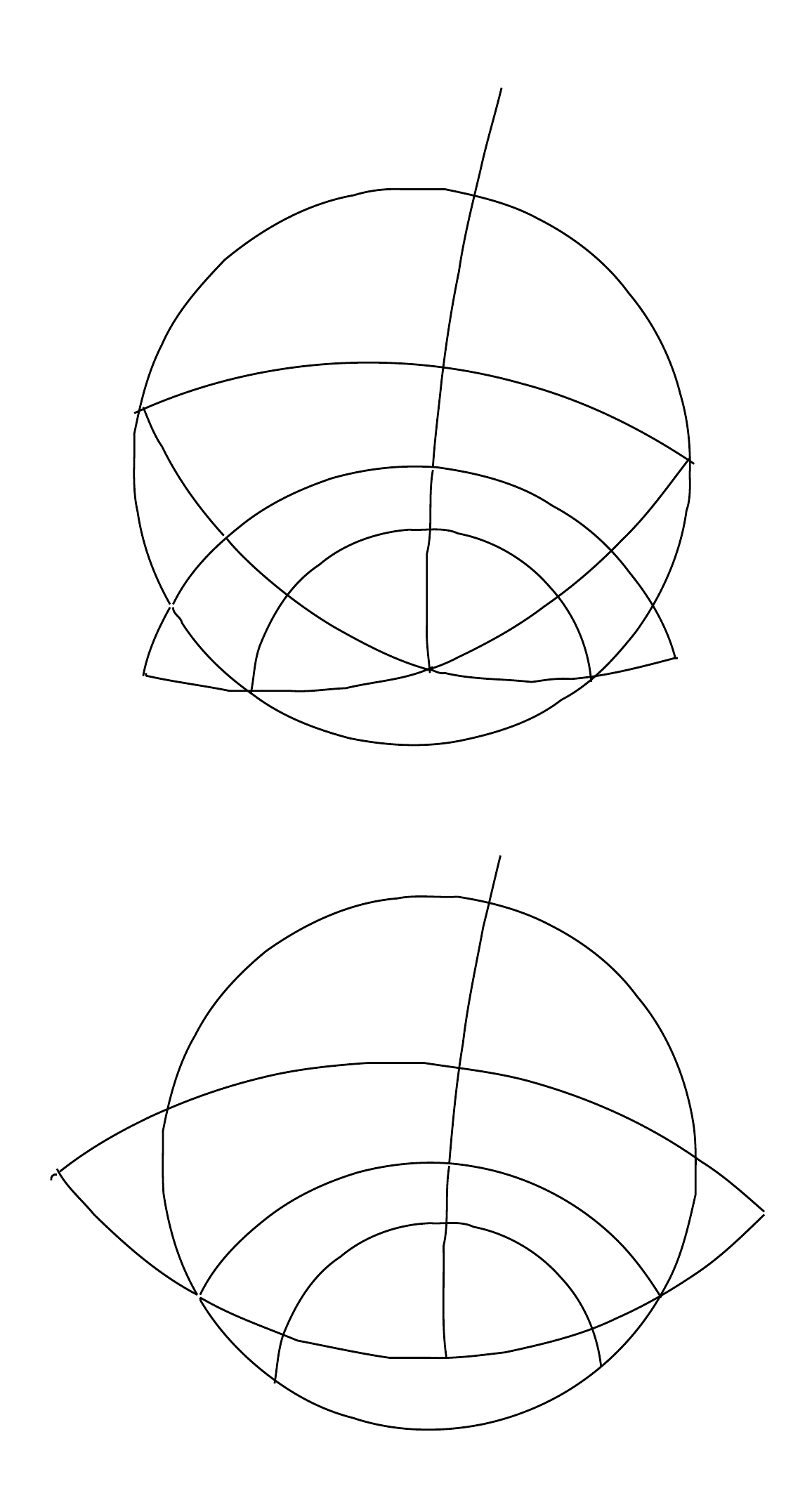\normalsize
\end{center}

\newpage\phantomsection
\begin{center}
  \small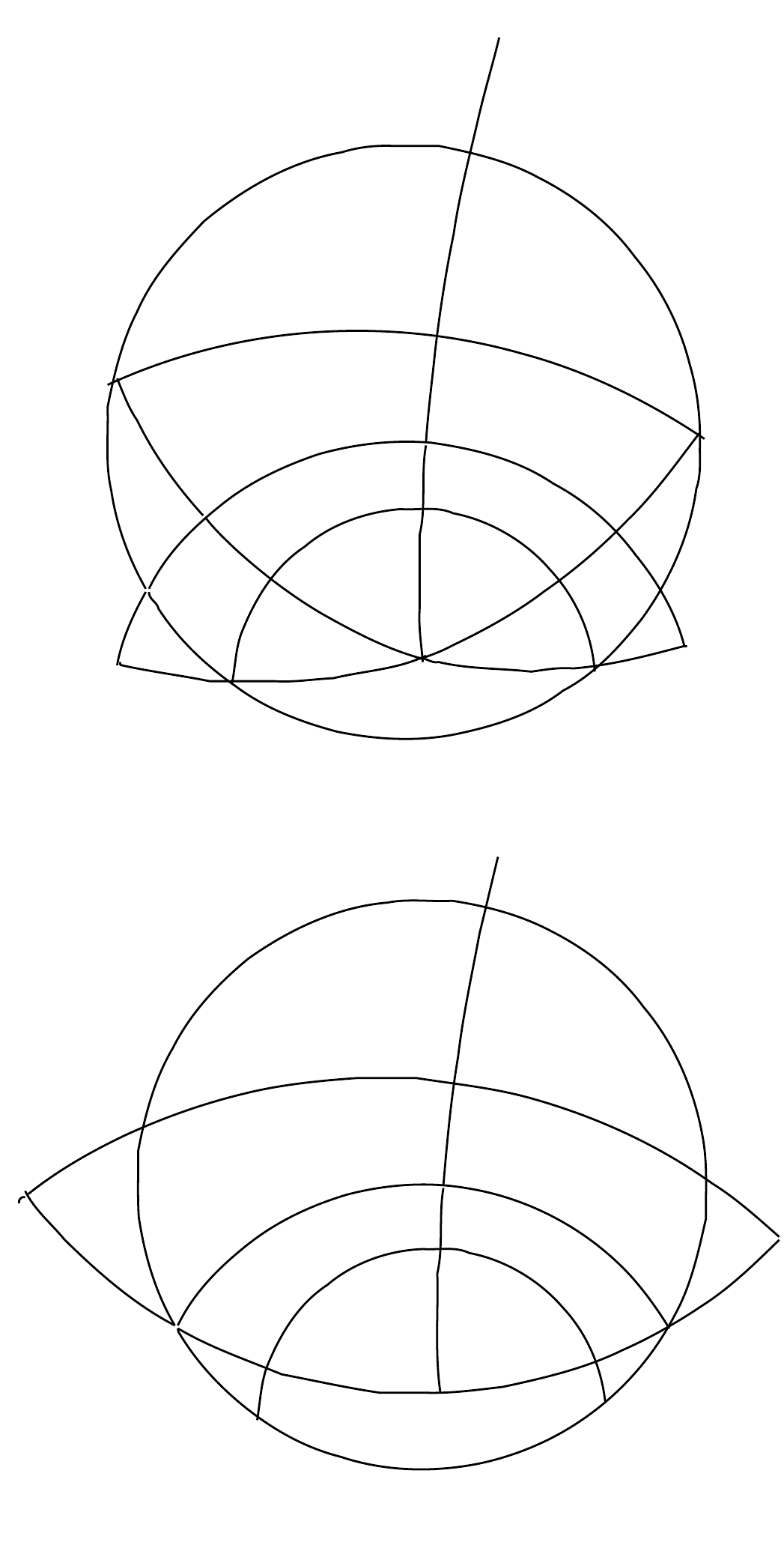\normalsize
\end{center}

\newpage\phantomsection
\index{AHAYAO@\RL{b`d}!AHAYAO AYBF ASBEAJ@\RL{b`d `an al-samt}, distance zénithale}
\index{ASBFBH@\RL{snw}!BAAUBHBD ASBFAI@\RL{fu.sUl al-sanaT}, les saisons}
\index{ASBEAJ@\RL{smt}!ASBEAJ AQABAS@\RL{samt al-ra's}, zénith}
\index{ACBHAL@\RL{'awj}!ACBHAL ATBEAS@\RL{'awj al-^sams}, Apogée du Soleil}
\index{AMAQAQ@\RL{.harara}!AMAQAQ@\RL{.hrr}, chaud}
\index{AHAQAL@\RL{brj}!BEBFAWBBAI AHAQBHAL@\RL{min.taqaT al-burUj}, écliptique, ceinture de l'écliptique}
\index{BABHBB@\RL{fwq}!ADBABB@\RL{'ufuq}, horizon}
\includepdf[pages=14,pagecommand={\thispagestyle{plain}}]{edit2.pdf}
\addcontentsline{toc}{chapter}{II.4 Particularités des lieux dont la latitude ne dépasse pas le complément de l'inclinaison [de l'écliptique]}
\begin{center}
  \Large Quatrième section

  \normalsize Particularités des lieux dont la latitude ne dépasse pas le complément de l'inclinaison [de l'écliptique]
\end{center}

\noindent Ces lieux se répartissent en quatre catégories.

\emph{Première catégorie :} là où la latitude est inférieure à
l'inclinaison de l'écliptique. En ces lieux, le Soleil passe au zénith
en deux points [de l'écliptique] dont les déclinaisons sont égales à
la latitude du lieu et qui sont situés du même côté que le pôle
visible ; la ceinture de l'écliptique est alors perpendiculaire à
l'horizon, ses pôles sont à l'horizon, et un homme debout n'a pas
d'ombre à midi. Tant que le Soleil est dans l'arc compris entre ces
deux points du côté du pôle visible, l'ombre sera portée du côté du
pôle caché, et le pôle visible de l'écliptique sera celui des deux qui
est proche du pôle caché de l'équateur, et son pôle caché sera celui
qui est proche du pôle visible [de l'équateur]. Tant que le Soleil est
dans l'arc opposé compris entre ces deux points du côté du pôle caché,
l'ombre sera portée du côté du pôle visible de l'équateur, et le pôle
visible de l'écliptique sera celui qui est proche du pôle visible de
l'équateur, et son pôle caché celui qui est proche du pôle caché. Plus
la latitude croît, plus les deux points se rapprochent l'un de
l'autre. L'arc entre les deux points, ainsi que les pôles de
l'écliptique, ont des levers et des couchers. Les saisons ne sont pas
égales. Ainsi leurs étés sont plus longs. Puisque le Soleil atteint le
zénith deux fois l'an, les saisons ne seront pas égales même si l'on
en comptait plus que quatre, car le Soleil atteint deux distances
zénithales maximales différentes, suivant qu'il est d'un côté ou de
l'autre de l'équateur.

\emph{Deuxième catégorie :} là où la latitude est égale à
l'inclinaison de l'écliptique. En ces lieux, le Soleil passe au zénith
une seule fois par an, l'un des pôles de l'écliptique est toujours
visible, l'autre toujours caché, et ils ne touchent pas l'horizon dans
leurs trajectoires sauf quand le point solsticial situé du côté du
pôle visible atteint le zénith -- et alors la ceinture de l'écliptique
est perpendiculaire à l'horizon. Les ombres sont toute l'année du côté
du pôle visible, sauf quand le Soleil est au zénith -- alors il n'y en
a pas. La hauteur du Soleil croît de l'un des deux points solsticiaux
à l'autre, puis la variation change de sens et elle décroît jusqu'à
revenir à sa valeur initiale. Les saisons sont quatre, ni plus ni
moins, et elles sont égales et analogues alternativement pour les
contrées du Nord et pour celles du Sud ; mais là où la latitude Sud
est égale à l'inclinaison [de l'écliptique] il fait plus chaud que là
où la latitude Nord est égale à l'inclinaison de l'écliptique, parce
que l'Apogée du Soleil est au Nord, et que son périgée est au Sud.

\newpage\phantomsection
\index{AMAWAW@\RL{.h.t.t}!BFAMAWAGAW@\RL{in.hi.tA.t}, abaissement|see{\RL{irtifA`}}}
\index{AQBAAY@\RL{rf`}!AQAJBAAGAY@\RL{irtifA`}, hauteur (coordonnées azimutales)}
\includepdf[pages=15,pagecommand={\thispagestyle{plain}}]{edit2.pdf}

\emph{Troisième catégorie :} là où la latitude est supérieure à
l'inclinaison de l'écliptique mais inférieure à son complément. Le
Soleil n'y atteint pas le zénith. Deux hauteurs lui sont propres~:
l'une est la somme de l'inclinaison de l'écliptique et du complément
de la latitude du lieu, et l'autre est plus petite, c'est l'excédent
du complément de la latitude sur l'inclinaison de l'écliptique. Il en
est de même du pôle toujours visible, car il n'atteint jamais
l'horizon, et sa hauteur maximale est atteinte lors du passage du
solstice caché au méridien du lieu, et sa hauteur minimale lors du
passage du solstice visible au méridien du lieu. Quant au pôle caché,
d'après le même raisonnement, il atteint deux hauteurs
négatives\footnote{hauteur négative, \textit{in\d{h}i\d{t}\=a\d{t}}
  que nous avons ailleurs traduit par <<~abaissement~>>, désigne une
  hauteur mesurée le long d'un grand cercle perpendiculaire à
  l'horizon mais \emph{sous} l'horizon.}. Les autres conditions, quant
au fait que les ombres tombent du côté du pôle visible et quant à la
longueur du jour et de la nuit, sont comme nous avons montré
ci-dessus. Quant aux planètes, celles dont la latitude dépasse
l'excédent de la latitude du lieu sur l'inclinaison de l'écliptique
passeront par le zénith deux fois ; celles dont la latitude est égale
à l'excédent y passeront une fois~; et celles dont la latitude est
moindre ne passeront pas par le zénith. Plus la latitude augmente,
plus l'équation du jour et les amplitudes Est et Ouest augmentent.

\emph{Quatrième catégorie :} là où la latitude est égale au complément
de l'inclinaison de l'écliptique. Là, la trajectoire du solstice situé
du côté du pôle visible est toujours visible, et la trajectoire de
l'autre solstice est toujours cachée. La trajectoire du pôle visible
de l'écliptique passe par le zénith, et la trajectoire de l'autre pôle
passe par le point opposé. Donc quand le point solsticial visible
arrive à l'horizon, il le touche un un point qui est le pôle du cercle
origine des azimuts du côté du pôle visible, le point solsticial
invisible touche l'autre pôle du cercle origine des azimuts, les deux
pôles [de l'écliptique] viennent au zénith et au [nadir], la ceinture
de l'écliptique coïncide avec le cercle de l'horizon, la tête du
Bélier avec le point Est, le commencement de la Balance avec le point
Ouest, la tête du Cancer avec le point Nord, et la tête du Capricorne
avec le point Sud ; et le point de l'équateur en face de la tête du
Capricorne sera situé sur le méridien du lieu au Sud au dessus de la
Terre\footnote{c'est-à-dire au dessus de l'horizon.}, et le point de
l'équateur en face de la tête du Cancer sera situé sur le méridien du
lieu au Nord au dessous de la Terre\footnote{c'est-à-dire sous
  l'horizon.}, si le pôle visible est au Nord. De cela, on connaît
aussi la position des deux ceintures relativement à l'horizon quand
[le pôle visible] est au Sud.

\newpage\phantomsection
\index{AJAGBHAOBHASBJBHAS@\RL{tAwdUsyUs}, Théodose}
\index{BHASAY@\RL{ws`}!ASAYAI@\RL{s`aT m^srq / m.grb}, amplitude Est / Ouest}
\index{BFBGAQ@\RL{nhr}!AJAYAOBJBD BFBGAGAQ@\RL{ta`dIl al-nahAr}, équation du jour}
\index{BJBHBE@\RL{ywm}!BJBHBE AHBDBJBDAJBG@\RL{al-yawm bilaylatihi j al-'ayyAm bilayAlIhA}, nychtémère (\textit{litt.} le jour avec sa nuit)}
\index{AVBHAB@\RL{.daw'}!AVBHAB@\RL{.daw'}, clarté}
\index{BBBJAS@\RL{qys}!BEBBBJAGAS@\RL{miqyAs j maqAyIs}, gnomon}
\index{ASBCBF@\RL{skn}!BEASBCBHBF@\RL{al-maskUn}, la partie habitée (du globe)}
\index{AYBEAQ@\RL{`mr}!BEAYBEBHAQ@\RL{al-ma`mUr, al-`imAraT}, partie cultivée de la Terre~; par extension, le quart habité (\RL{maskUn})}
\includepdf[pages=16,pagecommand={\thispagestyle{plain}}]{edit2.pdf}

Quand le pôle s'incline vers l'Ouest par rapport au zénith dans le
sens direct\footnote{Le sens direct est ici le sens du mouvement
  diurne, c'est-à-dire d'Est en Ouest ; ce n'est donc pas le sens
  ``des signes''.}, le point solsticial et la moitié Est de la
ceinture de l'écliptique (dont le centre est l'équinoxe de printemps)
s'élèvent [au dessus de l'horizon], si le pôle visible est au Nord ;
par là même, la moitié Ouest s'abaissera aussi. Les deux
cercles\footnote{Ceinture de l'écliptique et horizon.} se couperont en
deux points opposés proches des solstices et des points Nord et Sud :
il y a tangence en ces quatre points, donc nécessairement ils se
couperont en d'autres points. Alors la partie suivant le solstice
caché près du pôle du cercle origine des azimuts tendra à se coucher,
et la partie suivant le solstice visible près de l'autre pôle du
cercle origine des azimuts tendra à se lever~; puis les parties de la
moitié cachée se lèveront l'une après l'autre en chacune des parties
de la moitié Est de l'horizon, et les parties de la moitié visible se
coucheront l'une après l'autre en chacune des parties de la moitié
Ouest de l'horizon, pendant la durée d'un nychtémère, jusqu'à ce que
la position de l'orbe revienne à son état initial.

L'amplitude Est et l'équation du jour atteint [en ces lieux] un quart
de cercle. Le jour s'allonge jusqu'à ce que la durée d'un nychtémère
soit un jour entier lorsque le Soleil arrive au solstice visible,
parce que sa trajectoire est alors toujours visible comme nous l'avons
dit, puis apparaissent des nuits, et elles s'allongent jusqu'à ce que
la durée [d'un nychtémère] soit une nuit entière lorsqu'il arrive au
solstice caché ; c'est ainsi si l'on considère le début de la nuit ou
du jour comme étant l'arrivée du centre du Soleil à l'horizon. Si l'on
considère le début du jour comme étant l'apparition de la clarté et la
disparition des étoiles fixes, et le début de la nuit comme étant la
disparation [de la clarté] et l'apparition [des étoiles fixes], alors
la durée [maximale] de chacun des deux est \emph{un mois} d'après ce
qu'a montré Théodose dans les \textit{Habitations}. La hauteur du
Soleil croît jusqu'au double de l'inclinaison de l'écliptique, puis
elle décroît jusqu'à s'évanouir au point de tangence avec
l'horizon. Après son lever au pôle du cercle origine des azimuts, il
s'élève à l'Est, il franchit [le grand cercle] correspondant à la
droite Est-Ouest, il atteint sa hauteur maximale à son arrivée au
méridien du lieu, au Sud, et c'est le double de l'inclinaison [de
  l'écliptique] ; puis sa hauteur décroît jusqu'à toucher l'horizon au
pôle du cercle origine des azimuts. Puisque le Soleil tourne autour du
gnomon, et que l'ombre est toujours du côté opposé, l'ombre tourne
aussi autour du gnomon. Comme nous l'avons montré, en ces lieux, une
moitié du l'orbe de l'écliptique se lève avec une révolution de
l'équateur, tandis que l'autre moitié se lève instantanément. La
partie habitée s'arrête, du côté Nord, aux environs de cet horizon,
comme on sait.

\newpage\phantomsection
\index{AHAQAL@\RL{brj}!BEBFAWBBAI AHAQBHAL@\RL{min.taqaT al-burUj}, écliptique, ceinture de l'écliptique}
\index{BABHBB@\RL{fwq}!ADBABB@\RL{'ufuq}, horizon}
\index{AUAHAM@\RL{.sb.h}!AUAHAM@\RL{.sb.h}, aurore}
\index{ATBABB@\RL{^sfq}!ATBABB@\RL{^sfq}, crépuscule}
\index{ASBFBH@\RL{snw}!BAAUBHBD ASBFAI@\RL{fu.sUl al-sanaT}, les saisons}
\index{AMAQBC@\RL{.hrk}!AMAQBCAI ABBHBDBI@\RL{.harakaT 'Ul_A}, premier mouvement|see{\RL{.harakaT yawumiyyaT}}}
\index{AMAQBC@\RL{.hrk}!AMAQBCAI BJBHBEBJBJAI@\RL{.harakaT yawumiyyaT}, mouvement diurne}
\includepdf[pages=17,pagecommand={\thispagestyle{plain}}]{edit2.pdf}
\addcontentsline{toc}{chapter}{II.5 Particularités des lieux dont la latitude dépasse le complément de l'inclinaison [de l'écliptique] sans atteindre un quart de cercle}
\begin{center}
  \Large Cinquième section

  \normalsize Particularités des lieux dont la latitude dépasse le complément de l'inclinaison [de l'écliptique] sans atteindre un quart de cercle
\end{center}

\noindent En ces lieux, la trajectoire du pôle de l'écliptique dévie
un peu du zénith, du côté du pôle caché. La plus grande trajectoire
diurne toujours visible est plus grande que la trajectoire du
solstice, et elle coupe la ceinture de l'écliptique en deux points de
déclinaisons égales, du côté du pôle visible ; et la plus grande
trajectoire diurne toujours cachée la coupe en deux points qui sont
leurs opposés, du côté du pôle caché. La ceinture de l'écliptique se
divise en quatre portions. \emph{La première}, toujours visible, est
celle dont le centre est le solstice du côté du pôle visible, et le
Soleil y reste pendant un jour de leur été. \emph{La deuxième},
toujours cachée, est celle dont le centre est le solstice du côté du
pôle caché, et le Soleil y reste pendant une nuit de leur hiver. Les
extrémités du premier arc touchent l'horizon en l'un des pôles du
cercle origine des azimuts, du côté du pôle visible, et ils ne se
couchent pas avec le mouvement diurne ; les extrémités du deuxième arc
touche [l'horizon] en l'autre pôle du cercle origine des azimuts, et
ils ne se lèvent pas. \emph{La troisième} [portion] est celle dont le
centre est le commencement du Bélier ; elle se lève à rebours,
c'est-à-dire que sa fin se lève avant son début, et elle se couche à
l'endroit, c'est-à-dire que son début se couche avant sa fin, quand le
pôle visible est au Nord ; mais elle se lève à l'endroit et elle se
couche à rebours quand le pôle visible est Sud. \emph{La quatrième}
est celle dont le centre est le commencement de la Balance, et son
régime est à l'inverse de la troisième.

Le solstice visible a deux hauteurs. Sa hauteur maximale est la somme
de l'inclinaison de l'écliptique et du complémentaire de la latitude
du lieu, et [elle est atteinte] sur le méridien du lieu, du côté du
pôle caché. Sa hauteur minimale est l'excédent de la latitude du lieu
sur le complémentaire de l'inclinaison de l'écliptique, et [elle est
  atteinte] sur le méridien du lieu, du côté du pôle visible. Le pôle
visible de l'écliptique a aussi deux hauteurs. Sa hauteur maximale est
la somme du complémentaire de la latitude du lieu et du complémentaire
de l'inclinaison de l'écliptique. Sa hauteur minimale est l'excédent
de la latitude du lieu sur l'inclinaison de l'écliptique. Le pôle et
le solstice passent ensemble au méridien du lieu, de part et d'autre
du zénith, à des hauteurs complémentaires. On déduit de manière
analogue l'état du solstice et du pôle cachés. En ces horizons,
l'aurore et le crépuscule [peuvent] durer longtemps. L'ombre [peut]
tomber de tous les côtés, mais elle est plus longue du côté du pôle
caché.

\newpage\phantomsection
\index{BHASAY@\RL{ws`}!ASAYAI@\RL{s`aT m^srq / m.grb}, amplitude Est / Ouest}
\includepdf[pages=18,pagecommand={\thispagestyle{plain}}]{edit2.pdf}

Pour se représenter facilement les positions, supposons que la
latitude du lieu est $70$ degrés Nord, que l'arc toujours visible
[contient] les Gémeaux et le Cancer, que l'arc toujours caché
[contient] le Sagittaire et le Capricorne, que l'arc qui se lève à
rebours et se couche à l'endroit va du début du Verseau à la fin du
Taureau, et que celui qui se lève à l'endroit et se couche à rebours
va du début du Lion à la fin du Scorpion. Si la tête du Cancer est au
Sud, à sa hauteur maximale $43;35$, alors le pôle visible est au Nord
sur le méridien du lieu à sa hauteur minimale $46;25$, le début de la
Balance s'apprête à se lever au point du lever des équinoxes, le début
du Bélier s'apprête à se coucher au point de leur coucher, et la
moitié visible de l'orbe de l'écliptique est au Sud, d'Ouest en Est,
comme on voit sur la figure.

Que les orbes soient mûs par le premier mouvement, alors le début du
Cancer s'abaisse vers l'Ouest, le pôle de l'écliptique s'élève vers
l'Est, le début du printemps se couche, le début de l'automne se lève,
et de même pour les deux arcs qu'ils délimitent ; [le long de
  l'horizon], la distance entre le lever de chaque partie -- resp. le
coucher de la partie opposée -- et le lever de l'équinoxe -- resp. son
coucher -- croît jusqu'à ce que vienne le tour des deux parties dont
l'une touche l'horizon sans [jamais] se coucher et l'autre le touche
sans [jamais] se lever. Alors la Balance et le Scorpion se sont levés
à l'endroit, et leurs amplitudes Est ont occupé tout le quart
Sud-Est. Le Bélier et le Taureau se sont couchés à l'endroit, et leurs
amplitudes Ouest ont occupé tout le quart Nord-Ouest. Le début du
Sagittaire est au point Sud où il touche l'horizon ; et le début des
Gémeaux touche l'horizon au Nord. Le pôle visible de l'écliptique est
du côté Est, entre ses hauteurs maximale et minimale, sur le cercle
origine des azimuts ; son autre pôle est à l'opposé. La moitié visible
de la ceinture de l'écliptique est maintenant à l'Ouest, du Sud au
Nord ; et sa moitié cachée est à l'opposé. L'intersection de la
ceinture de l'écliptique et de l'horizon est aux points Nord et
Sud. Voir la figure.

Que les orbes [continuent] de se mouvoir, alors le début des Gémeaux
prend de la hauteur en allant vers l'Est. Adjacente, la fin du Taureau
se lève peu à peu -- de sorte que le lever de chaque partie soit plus
proche du lever de l'équinoxe que les levers des parties qui se lèvent
avant elle. Le Taureau se lève, puis la fin du Bélier, puis son
début. Le quart Nord-Est est entièrement occupé par les amplitudes Est
de ces signes. Finalement le début du Bélier atteint son lever. En
face de cela, le début du Sagittaire se met à s'abaisser sous
l'horizon, la fin du Scorpion -- qui lui est adjacente -- se couche,
il se couche peu à peu, puis la Balance se couche de sa fin jusqu'à
son début, et leurs amplitudes Ouest occupent tout le quart
Sud-Ouest. Le début de la Balance atteint son coucher.

\newpage\phantomsection
\begin{center}
  \small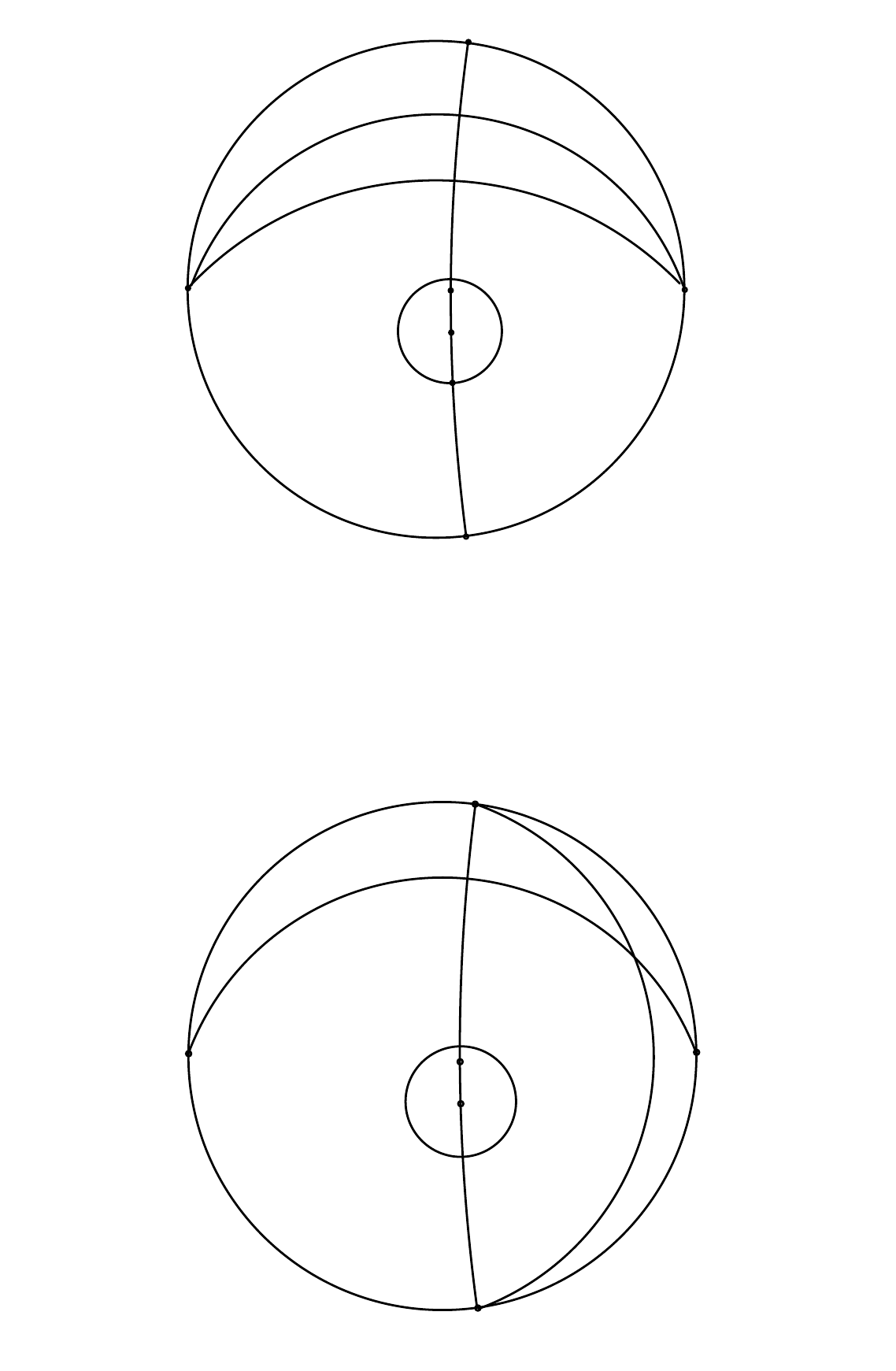\normalsize
\end{center}

\newpage\phantomsection
\begin{center}
  \small
\begingroup%
  \makeatletter%
  \providecommand\color[2][]{%
    \errmessage{(Inkscape) Color is used for the text in Inkscape, but the package 'color.sty' is not loaded}%
    \renewcommand\color[2][]{}%
  }%
  \providecommand\transparent[1]{%
    \errmessage{(Inkscape) Transparency is used (non-zero) for the text in Inkscape, but the package 'transparent.sty' is not loaded}%
    \renewcommand\transparent[1]{}%
  }%
  \providecommand\rotatebox[2]{#2}%
  \ifx\svgwidth\undefined%
    \setlength{\unitlength}{327.39820557bp}%
    \ifx\svgscale\undefined%
      \relax%
    \else%
      \setlength{\unitlength}{\unitlength * \real{\svgscale}}%
    \fi%
  \else%
    \setlength{\unitlength}{\svgwidth}%
  \fi%
  \global\let\svgwidth\undefined%
  \global\let\svgscale\undefined%
  \makeatother%
  \begin{picture}(1,1.46287372)%
    \put(0,0){\includegraphics[width=\unitlength]{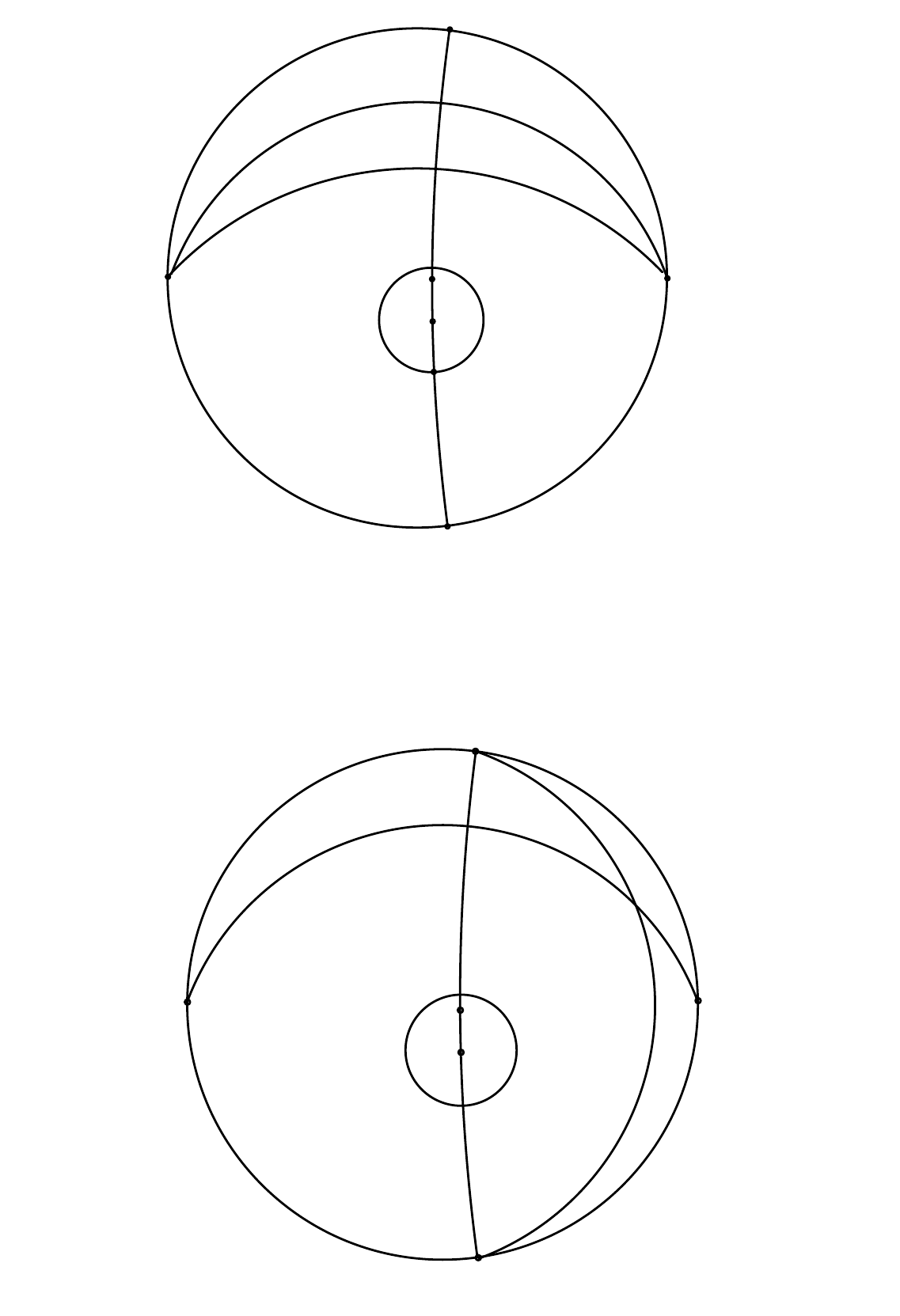}}%
    \put(0.48858761,1.44590757){\color[rgb]{0,0,0}\makebox(0,0)[lb]{\smash{point Sud}}}%
    \put(0.75044307,1.15222617){\color[rgb]{0,0,0}\makebox(0,0)[lb]{\smash{point Ouest}}}%
    \put(0.74869768,1.12255501){\color[rgb]{0,0,0}\makebox(0,0)[lb]{\smash{tête du Bélier}}}%
    \put(0.67371753,0.92678643){\color[rgb]{0,0,0}\rotatebox{41.26721321}{\makebox(0,0)[lb]{\smash{horizon}}}}%
    \put(0.49728858,0.85104354){\color[rgb]{0,0,0}\makebox(0,0)[lb]{\smash{point Nord}}}%
    \put(0.49066006,1.03123052){\color[rgb]{0,0,0}\makebox(0,0)[lb]{\smash{pôle de l'écliptique}}}%
    \put(0.48902757,1.10528){\color[rgb]{0,0,0}\makebox(0,0)[lb]{\smash{pôle Nord}}}%
    \put(0.4876056,1.14678894){\color[rgb]{0,0,0}\makebox(0,0)[lb]{\smash{zénith}}}%
    \put(0.17298867,1.15875842){\color[rgb]{0,0,0}\makebox(0,0)[rb]{\smash{point Est}}}%
    \put(0.17262045,1.13083265){\color[rgb]{0,0,0}\makebox(0,0)[rb]{\smash{tête de la Balance}}}%
    \put(0.29924781,1.24500902){\color[rgb]{0,0,0}\rotatebox{21.87109053}{\makebox(0,0)[lb]{\smash{ceinture de l'écliptique}}}}%
    \put(0.29758675,1.30515714){\color[rgb]{0,0,0}\rotatebox{33.67090665}{\makebox(0,0)[lb]{\smash{équateur}}}}%
    \put(0.51331485,0.6429357){\color[rgb]{0,0,0}\makebox(0,0)[lb]{\smash{point Sud}}}%
    \put(0.5126734,0.67409423){\color[rgb]{0,0,0}\makebox(0,0)[lb]{\smash{début du Sagittaire}}}%
    \put(0.78739814,0.34608697){\color[rgb]{0,0,0}\makebox(0,0)[lb]{\smash{point Ouest}}}%
    \put(0.78739814,0.3094345){\color[rgb]{0,0,0}\makebox(0,0)[lb]{\smash{coucher de l'équinoxe}}}%
    \put(0.69099574,0.10977405){\color[rgb]{0,0,0}\rotatebox{41.26721321}{\makebox(0,0)[lb]{\smash{horizon}}}}%
    \put(0.52632194,0.03815326){\color[rgb]{0,0,0}\makebox(0,0)[lb]{\smash{point Nord}}}%
    \put(0.52568034,0.00036987){\color[rgb]{0,0,0}\makebox(0,0)[lb]{\smash{début des Gémeaux}}}%
    \put(0.52111561,0.28579926){\color[rgb]{0,0,0}\makebox(0,0)[lb]{\smash{pôle Nord}}}%
    \put(0.51867144,0.33410358){\color[rgb]{0,0,0}\makebox(0,0)[lb]{\smash{zénith}}}%
    \put(0.19232223,0.35189626){\color[rgb]{0,0,0}\makebox(0,0)[rb]{\smash{point Est}}}%
    \put(0.19369932,0.31349841){\color[rgb]{0,0,0}\makebox(0,0)[rb]{\smash{lever de l'équinoxe}}}%
    \put(0.30252341,0.49061299){\color[rgb]{0,0,0}\rotatebox{33.67090665}{\makebox(0,0)[lb]{\smash{équateur}}}}%
    \put(0.62160168,0.128601){\color[rgb]{0,0,0}\rotatebox{63.25093676}{\makebox(0,0)[lb]{\smash{ceinture de l'écliptique}}}}%
  \end{picture}%
\endgroup%
\normalsize
\end{center}

\newpage\phantomsection
\index{AMAQBC@\RL{.hrk}!AMAQBCAI AKAGBFBJAI@\RL{.hrkaT _tAnyaT}, deuxième mouvement}
\index{AWBDAY@\RL{.tl`}!AWBDAY@\RL{.tl`}, se lever ($\neq$ \RL{.grb})}
\index{AZAQAH@\RL{.grb}!AZAQAH@\RL{.grb}, se coucher ($\neq$ \RL{.tl`})}
\includepdf[pages=19,pagecommand={\thispagestyle{plain}}]{edit2.pdf}

\noindent Le début du Cancer atteint le méridien du côté Nord où il
est à sa hauteur minimale $3;35$ et où le pôle de l'écliptique est à
sa hauteur maximale du côté Sud, $86;25$. La moitié visible de l'orbe
de l'écliptique est du côté Nord entre le lever des équinoxes et leur
coucher, mais dans un ordre différent de l'ordre habituel. Sa moitié
cachée est à l'opposé. La ceinture de d'écliptique et l'horizon se
coupent aux points Est et Ouest. Voir la figure.

Que les orbes [continuent] de se mouvoir, alors les Poissons se lèvent
de leur fin jusqu'à leur début, puis le Verseau, de sa fin jusqu'à son
début, de sorte que leurs amplitudes Est occupent tout le quart
Sud-Est. En face d'eux, la Vierge se couche de sa fin jusqu'à son
début, puis le Lion, de sa fin jusqu'à son début, leurs amplitudes
Ouest occupant tout le quart Nord-Ouest, et le début du Verseau vient
toucher l'horizon au point Sud. La moitié visible du cercle de
l'écliptique est entre les deux, du côté Est. Le début du Cancer a
pris de la hauteur à l'Est, et le pôle s'est mis à descendre à
l'Ouest. Voir la figure.

Que les orbes [continuent] de se mouvoir, alors le début du Lion
s'élève de l'horizon en se dirigeant vers l'Est, et les parts [du
  Lion] se lèvent dans l'ordre jusqu'à sa fin, puis [jusqu'à] la fin
de la Vierge, de même. Leurs amplitudes Est occupent tout le quart
Nord-Est. En face de cela, le début du Verseau s'abaisse sous terre à
partir de l'horizon, et le Verseau se couche puis les Poissons dans
l'ordre [des signes], occupant de leurs amplitudes Ouest tout le quart
Sud-Ouest, le lever atteignant le début de la Balance, et le coucher
atteignant le début du Bélier ; car la proximité des levers des parts
[de l'écliptique] -- resp. de leurs couchers -- et du lever --
resp. du coucher -- de l'équinoxe augmente. Pendant ce temps-là, le
début du Cancer a atteint sa hauteur maximale au méridien, et le pôle
visible de l'écliptique a atteint sa hauteur minimale au méridien. La
moitié visible de l'orbe de l'écliptique est du côté du Sud. Le cercle
est complet, et on retrouve les positions supposées au commencement.

En ces horizons, si la latitude du pays se rapproche de la latitude
maximale et que la hauteur de l'équateur par rapport à l'horizon
devient petite, il peut arriver qu'un astre ayant une trajectoire
[diurne] très proche de l'horizon se déplace vers une autre
trajectoire [diurne] à cause du deuxième mouvement. Il [peut] alors
disparaître après avoir été visible bien qu'il soit dans la moitié Est
; ou il peut apparaître après avoir été caché bien qu'il soit dans la
moitié Ouest. Il se serait alors couché à l'Est, ou levé à
l'Ouest. Ceci compte parmi les questions singulières.

\newpage\phantomsection
\begin{center}
  \small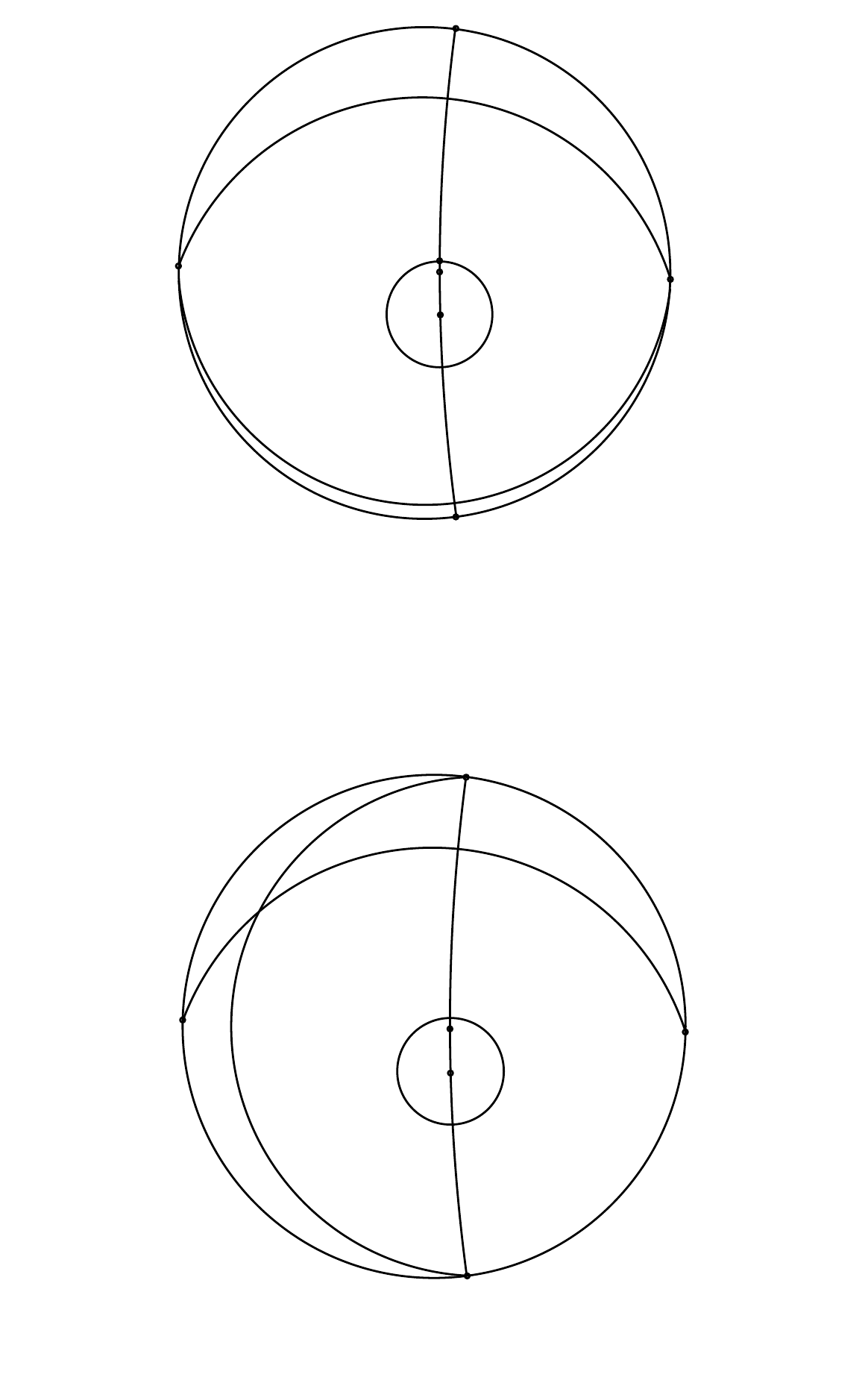\normalsize
\end{center}

\newpage\phantomsection
\begin{center}
  \small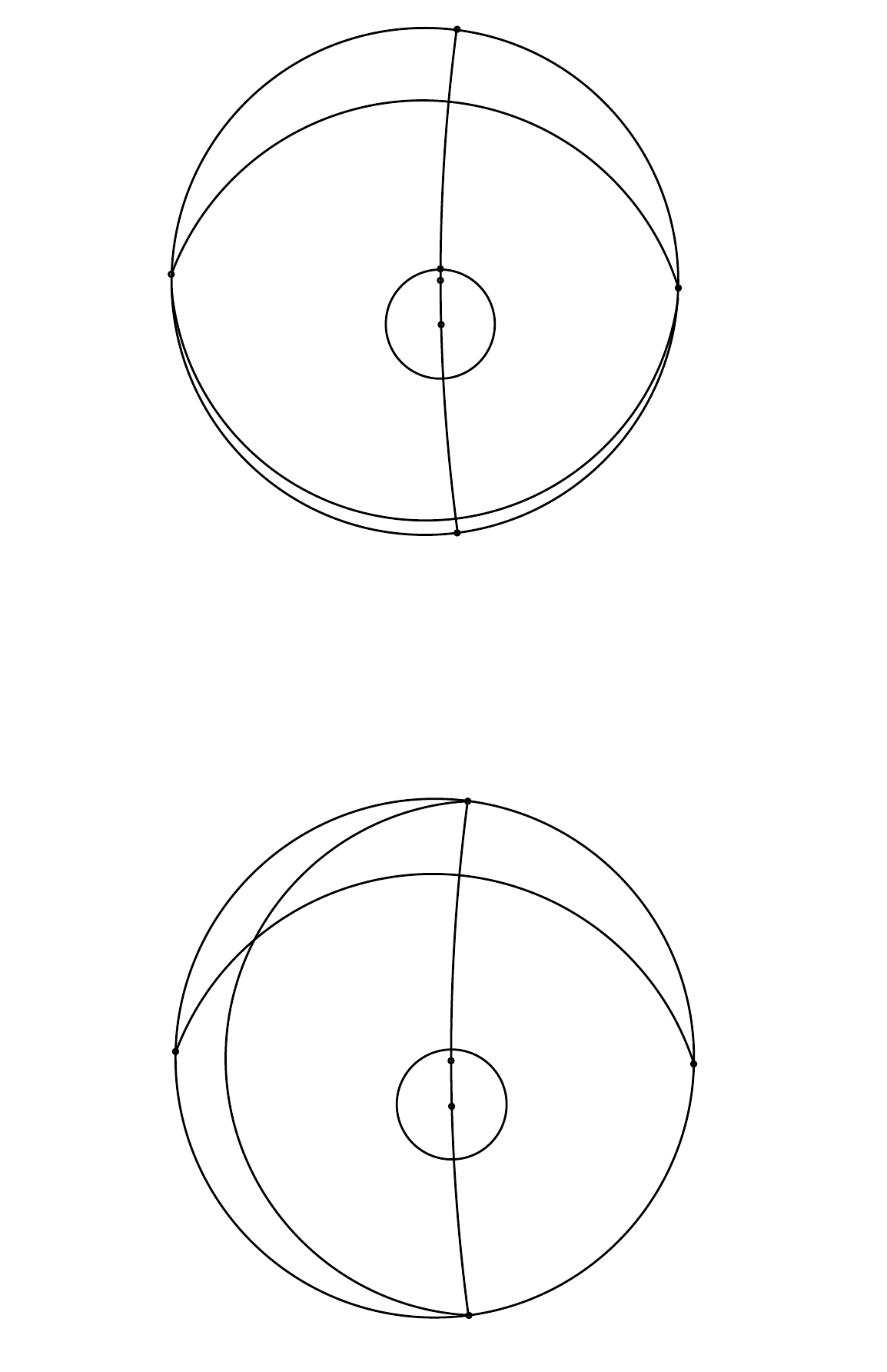\normalsize
\end{center}

\newpage\phantomsection
\index{AJAGBHAOBHASBJBHAS@\RL{tAwdUsyUs}, Théodose} 
\index{AMAQBC@\RL{.hrk}!AMAQBCAI AKAGBFBJAI@\RL{.hrkaT _tAnyaT}, deuxième mouvement}
\index{AWBDAY@\RL{.tl`}!AWBDAY@\RL{.tl`}, se lever ($\neq$ \RL{.grb})}
\index{AZAQAH@\RL{.grb}!AZAQAH@\RL{.grb}, se coucher ($\neq$ \RL{.tl`})}
\index{ATAQBB@\RL{^srq}!BEATAQBB@\RL{ma^sraq}, Est}
\index{AZAQAH@\RL{.grb}!BEAZAQAH@\RL{ma.grab}, Ouest}
\index{AOBHAQ@\RL{dwr}!BEAOAGAQ@\RL{mdAr}, trajectoire, trajectoire diurne (\textit{i. e.} cercle parallèle à l'équateur)}
\index{AHAQAL@\RL{brj}!BEBFAWBBAI AHAQBHAL@\RL{min.taqaT al-burUj}, écliptique, ceinture de l'écliptique}
\index{BABHBB@\RL{fwq}!ADBABB@\RL{'ufuq}, horizon}
\index{BFBGAQ@\RL{nhr}!BFBGAGAQ@\RL{nahAr}, jour}
\index{ACBHAL@\RL{'awj}!ACBHAL ATBEAS@\RL{'awj al-^sams}, Apogée du Soleil}
\index{AVBHAB@\RL{.daw'}!AVBHAB@\RL{.daw'}, clarté}
\index{AUAHAM@\RL{.sb.h}!AUAHAM@\RL{.sb.h}, aurore}
\index{ATBABB@\RL{^sfq}!ATBABB@\RL{^sfq}, crépuscule}
\addcontentsline{toc}{chapter}{II.6 Particularités des horizons dont la latitude est égale à un quart de cercle}
\includepdf[pages=20,pagecommand={\thispagestyle{plain}}]{edit2.pdf}

\newpage\phantomsection
\begin{center}
  \Large Sixième section

  \normalsize Particularités des horizons dont la latitude est égale à un quart de cercle
\end{center}

En ces lieux, un des pôles de l'équateur est au zénith, et l'équateur
et l'horizon sont confondus. Par le premier mouvement, les cieux
tournent comme une meule. L'Est ne se distingue pas de l'Ouest,
puisque lever et coucher sont possibles dans toutes les directions. La
hauteur maximale est égale à l'inclinaison de l'écliptique, et de même
pour l'abaissement maximal\footnote{\textit{i. e.} la hauteur minimale
  comptée négativement sous l'horizon.}. La moitié de l'orbe de
l'équateur située du côté du pôle visible est toujours visible, et
l'autre moitié est toujours cachée. De même pour l'orbe de
l'écliptique ; ainsi, il fait jour tant que le Soleil est dans la
moitié visible [de l'écliptique], et il fait nuit tant qu'il est dans
l'autre moitié.

L'année est un jour et une nuit, mais les deux diffèrent à cause de la
vitesse ou de la lenteur du mouvement [du Soleil]. \`A la date
présente, au pôle Nord, le jour dure sept nychtémères de plus que la
nuit, parce que l'Apogée du Soleil est vers la fin des Gémeaux et que
son périgée est vers la fin du Sagittaire. C'est ainsi si l'on
considère que le jour va du lever [du Soleil] à son coucher et que la
nuit va de son coucher à son lever ; mais si l'on considère que le
jour va de l'apparition de la clarté -- qui fait disparaître les fixes
-- jusqu'à sa disparition, et que c'est le contraire pour la nuit,
alors le jour dure plus de sept mois et la nuit dure presque cinq mois
selon Théodose dans les \textit{Habitations}. Si l'on considère que le
jour va du lever de l'aurore jusqu'au coucher du crépuscule, alors le
jour dure neuf mois et sept jours environ -- de nos jours. En effet le
lever de l'aurore et le coucher du crépuscule durent chacun cinquante
jours -- de nos jours -- comme le montrera leur description
ultérieurement.

En ces lieux, le lever ou le coucher des astres, causé par le deuxième
mouvement, ne se fait pas en une position déterminée de
l'horizon. Toute étoile fixe dont la latitude est nulle sera pendant
douze mille ans au dessus de l'horizon, et pendant la même durée en
dessous. [Toute étoile fixe] dont la latitude est inférieure à
l'inclinaison [de l'écliptique] aura un lever et un coucher, mais la
durée de sa visibilité et de son invisibilité dépendra de la distance
de sa trajectoire par rapport à la ceinture de
l'écliptique\footnote{Ici, ``trajectoire'' ne désigne pas une
  trajectoire diurne bien qu'il s'agisse du même mot. C'est plutôt la
  trajectoire de l'étoile par rapport au référentiel attaché au
  neuvième orbe~: il s'agit d'une trajectoire circulaire dans un plan
  parallèle à la ceinture de l'écliptique. Rappelons que pour
  {\shatir} le mouvement de précession, \textit{i. e.} le ``deuxième
  mouvement'', est le mouvement de rotation du huitième orbe~--~les
  étoiles fixes~--~autour d'un axe incliné par rapport à l'axe
  Nord-Sud au sein du référentiel attaché au neuvième orbe.}.

\newpage\phantomsection
\includepdf[pages=21,pagecommand={\thispagestyle{plain}}]{edit2.pdf}
\noindent [Toute étoile fixe] dont la latitude est égale
à [l'inclinaison] touche
l'horizon à chaque révolution due au deuxième mouvement ; une telle
étoile -- de même qu'une étoile ayant une latitude supérieure à
l'inclinaison -- n'aura ni lever ni coucher, et elle sera toujours
visible ou toujours cachée. Qu'on se souvienne de ce qui a été
mentionné au sujet de la position de l'orbe à cause des deux premiers
mouvements~; ceci dépend de cela.

Ici s'achève le compte-rendu des caractéristiques des contrées [quant
  aux] trajectoires diurnes \textit{etc.}

\newpage\phantomsection
\index{AOAQAL@\RL{drj}!AOAQALAI@\RL{drjaT j darjAt, drj}, degré}
\index{AOAQAL@\RL{drj}!AOAQAL ASBHAGAB@\RL{drj al-sawA'}, dégrés égaux ($\neq$ leur coascension)}
\index{ALBEBGAQ@\RL{jmhr}!ALBEBGBHAQ@\RL{al-jumhUr}, les Grecs}
\index{AWBDAY@\RL{.tl`}!BEAWAGBDAY@\RL{m.tAl`}, coascension}
\index{AYAOBD@\RL{`dl}!BEAYAOAOBD BFBGAGAQ@\RL{m`ddl al-nhAr}, équateur, plan de l'équateur}
\index{AHAQAL@\RL{brj}!BEBFAWBBAI AHAQBHAL@\RL{min.taqaT al-burUj}, écliptique, ceinture de l'écliptique}
\index{BABHBB@\RL{fwq}!ADBABB@\RL{'ufuq}, horizon}
\index{AKBDAK@\RL{_tl_t}!BEAKBDAK@\RL{m_tl_t}, triangle (éventuellement sphérique)}
\index{AHAQAL@\RL{brj}!AHAQAL@\RL{brj}, signe (du zodiaque) |see{\RL{falak al-burUj}}}
\index{BABDBC@\RL{flk}!BABDBC AHAQBHAL@\RL{flk al-burUj}, orbe de l'écliptique}
\addcontentsline{toc}{chapter}{II.7 Les coascensions de l'écliptique}
\includepdf[pages=22,pagecommand={\thispagestyle{plain}}]{edit2.pdf}

\newpage\phantomsection
\begin{center}
  \Large Septième section

  \normalsize Les coascensions de l'écliptique
\end{center}

La \emph{coascension}, ce sont les parts d'arc de l'équateur qui se
lèvent en même temps que se lèvent des parts d'arc de l'écliptique. On
dit que l'arc de l'écliptique [est mesuré] en \emph{degrés égaux}. Les
coascensions varient selon les horizons. Leur origine est l'équinoxe
de printemps chez les Grecs. Chez certains [autres savants] c'est le
solstice d'hiver~; ce qu'ils font montre qu'ils adoptent cette
hypothèse.

\`A l'équateur terrestre, chaque quart délimité par deux des quatres
points des équinoxes et des solstices se lève avec un quart [de
  l'équateur]. En effet, quand le point de l'équinoxe -- qui est l'une
des deux bornes [de chacun] des deux quarts des ceintures de
l'équateur et de l'écliptique -- atteint le zénith, alors les
solstices atteignent l'horizon. Et comme le [grand cercle] passant par
les quatre pôles est confondu avec lui, alors [l'horizon] coupe
perpendiculairement les ceintures de l'écliptique et de l'équateur. On
raisonne de manière analogue pour chacun des quatre quarts. Ce qui se
lève, de l'équateur, avec le signe [du zodiaque] suivant un des quatre
points, n'est pas [un arc] égal, c'est-à-dire trente parts. Un signe
adjacent à un points des équinoxes est plus grand que sa coascension,
parce que l'arc d'écliptique dont il est question est plus grand que
l'arc correspondant de l'équateur. En effet dans le triangle formé par
eux et par l'horizon, cet arc-là est le côté opposé à l'angle droit
formé par l'équateur et l'horizon, alors que la coascension est le
côté opposé à l'angle aigu formé par l'écliptique et l'horizon ; et le
côté opposé à l'angle droit est plus grand que le côté opposé à
l'angle aigu. Ce qu'on vient de dire se produit aussi pour deux
signes\footnote{\textit{i. e.} un arc de 60° adjacent à l'équinoxe le
  long de l'écliptique.} suivant le même point et pour leur
coascension. Quant à [l'arc] complément du quart [de la
  circonférence] de l'écliptique et adjacent au solstice, sa
coascension est plus grande que lui. On le démontre ainsi~: ce qui
reste de l'équateur pour compléter le quart [de l'équateur] est
supérieur à ce qui reste à partir du point de l'écliptique pour
compléter le quart [de l'écliptique], et l'excédent de ceci est égale
au manque de cela. De l'excédent et du manque en coascension, et de
l'égalité de cet excédent et de ce manque, il apparaît aussi que deux
arcs égaux et à distances égales de l'un des quatre points ont des
coascensions égales à l'équateur terrestre.

\newpage\phantomsection
\index{ASAWAM@\RL{s.t.h}!ASAWAM BEASAJBH@\RL{\vocalize s.t.h mstwiN}, plan}
\index{BEAQAQ@\RL{mrr}!AOAQALAI BEBEAQAQ@\RL{darjaT mamarr}, degré de transit}
\index{BEAQAQ@\RL{mrr}!BEBEAQAQ@\RL{mamarr j At}, passage, transit (en général au méridien)}
\index{BEAQAQ@\RL{mrr}!BEAQBHAQ@\RL{murUr}|see{\RL{mamarr}}}
\index{BEBJBD@\RL{myl}!AOAGAEAQAI BEBJBD@\RL{dA'iraT al-mIl}, cercle de déclinaison}
\index{AKBDAK@\RL{_tl_t}!BEAKBDAK@\RL{m_tl_t}, triangle (éventuellement sphérique)}
\index{AZAQAH@\RL{.grb}!BEAZAGAQAH@\RL{ma.gArib}, co-coucher (<<~codescension~>>)}
\includepdf[pages=23,pagecommand={\thispagestyle{plain}}]{edit2.pdf}

Si on connaît la coascension d'un quart [de l'écliptique], alors on
connaît aussi la coascension des autres quarts. De la même manière,
de la coascension d'un signe on peut connaître la coascension des
suivants, jusqu'à tout connaître. Divisons la ceinture de l'écliptique
en quatre segments dont les extrémités sont les milieux des quarts
[ci-dessus]. Chaque segment dont le milieu est une équinoxe est plus
grand que sa coascension, et chaque segment dont le milieu est un
solstice est plus petit que sa coascension.

Il en est des transits de la ceinture de l'écliptique et de l'équateur
aux méridiens et aux cercles de déclinaisons dans toutes les contrées
comme il en est de leurs coascensions à l'équateur terrestre, puisque
chaque [cercle de déclinaison] est un horizon parmi les horizons de
l'équateur terrestre. Et puisque le co-coucher est [l'arc] de
l'équateur qui se couche en même temps que se lève [l'arc mesurant] la
coascension, alors le co-coucher de chaque signe est égal à sa
coascension.

Quand l'horizon passe par les pôles de l'équateur et de l'écliptique,
c'est-à-dire quand le point de l'équinoxe est au zénith [ou au nadir],
le fait que l'horizon coupe perpendiculairement ces deux cercles se
démontre par la contraposée, comme suit. Chacun des deux arcs (celui
de l'équateur et celui de l'écliptique) est un côté opposé à un [des
  deux angles] angles droits. Il y a deux angles droits dans le
triangle, car s'il n'y avait pas deux angles droits, alors un quart de
la ceinture de l'écliptique serait nécessairement plus grand qu'un
quart de l'équateur ; or l'hypothèse est qu'ils sont égaux. Certes
l'impossibilité [de deux angles droits dans un triangle] n'a lieu que
dans les triangles plans.

Passons aux horizons inclinés. Là, une moitié [se lève] en même temps
qu'une moitié quand elles sont délimitées par les équinoxes -- et non
un quart en même temps qu'un quart. En effet l'équateur n'est pas
perpendiculaire à l'horizon. Si un quart adjacent à une équinoxe et
situé du côté du pôle visible se lève, alors il sera plus grand que sa
coascension, d'autant plus grand que l'angle obtus sera plus grand
qu'un angle droit ; car [ce quart] est le côté opposé à un angle obtus
dans le triangle mentionné ci-dessus, et sa coascension est le côté
opposé à un angle aigu. S'il est au contraire situé du côté du pôle
caché, alors sa coascension sera plus grande que lui, puisqu'on sera
dans le régime inverse. De là, on démontre que les coascensions de
deux portions égales situées à des distances égales d'une des deux
équinoxes sont égales. Si l'on divise l'orbe de l'écliptique en deux
moitiés au milieu de chacune desquelles se situe un des deux points
des équinoxes, alors celle dont le milieu est le passage du côté du
pôle visible sera plus grande que son lever, et inversement. 

\newpage\phantomsection
\index{BFBGAQ@\RL{nhr}!AJAYAOBJBD BFBGAGAQ@\RL{ta`dIl al-nahAr}, équation du jour}
\index{AWBDAY@\RL{.tl`}!AJAYAOBJBD BEAWAGBDAY@\RL{ta`dIl ma.tAli`}, équation de la coascension}
\includepdf[pages=24,pagecommand={\thispagestyle{plain}}]{edit2.pdf}

\noindent La coascension des arcs au Nord dans les horizons du Nord
est égale à la coascension des arcs opposés au Sud dans les horizons
du Sud ; de même pour les horizons du Sud. Quel que soit l'horizon, le
co-coucher d'un arc est toujours égal à la coascension de l'arc
opposé.

En un lieu dont la latitude est égale au complément de l'inclinaison
de l'écliptique, on l'a déjà montré plus haut, la moitié de
l'écliptique se lève avec un tour complet de l'équateur, et l'autre
moitié se lève instantanément ; quant aux couchers, il faut échanger
les deux moitiés.

Là où la latitude dépasse le complément de l'inclinaison sans
toutefois atteindre un quart de circonférence, il y a des arcs
toujours visibles de l'orbe de l'écliptique, et ce qui est à l'opposé
est toujours caché. On divise l'équateur en des parties qui se lèvent
en même temps que les parties de l'écliptique se levant à rebours, et
des parties qui se lèvent en même temps que les parties de
l'écliptique se levant à l'endroit. Ainsi à $70$ degrés de latitude,
les Gémeaux et le Cancer sont toujours visibles, et le Sagittaire et
le Capricorne toujours cachés~; quand l'équinoxe de printemps se lève,
les Poissons se lèvent après elle, à rebours, de leur fin jusqu'à leur
début, puis le Verseau [se lève] à rebours aussi, puis le Lion
commence à se lever, à partir de son début, à l'endroit, puis la
Vierge puis la Balance puis le Scorpion, de même ; lorsqu'est atteint
le début du Sagittaire, la fin du Taureau commence à se lever à
rebours, puis le Taureau et le Bélier se lèvent à rebours, jusqu'à ce
que le point de l'équinoxe de printemps revienne à l'horizon. On
traite de la même manière les autres horizons, ainsi que les
couchers. Les portions de l'écliptique ou de l'équateur qui ne se
lèvent ni ne se couchent n'ont ni coascension ni co-coucher. C'est
pourquoi il n'y a ni coascension ni co-coucher à la latitude de $90$
degrés.

\emph{L'équation de la coascension} d'un pays est ce qu'on ajoute ou
retranche à la coascension de la sphère
droite\footnote{\textit{i. e.} aux coascensions calculées à
  l'équateur terrestre.}. C'est l'équation du jour, et le co-coucher
de l'arc est comme sa coascension\footnote{En effet, quand le Soleil
  a un lever et un coucher par nychtémère, la durée du jour est la
  coascension d'une moitié de l'écliptique.}. La partie
\emph{ascendante} de l'orbe de l'écliptique est celle qui vient à
l'horizon vers l'Est, et sa partie \emph{descendante} est la partie
opposée, à l'horizon Ouest.

\newpage\phantomsection
\index{BJBHBE@\RL{ywm}!BJBHBE AHBDBJBDAJBG@\RL{al-yawm bilaylatihi j al-'ayyAm bilayAlIhA}, nychtémère (\textit{litt.} le jour avec sa nuit)}
\index{AWBDAY@\RL{.tl`}!BEAWAGBDAY@\RL{m.tAl`}, coascension}
\index{BJBHBE@\RL{ywm}!BJBHBE AMBBBJBBBJ@\RL{ywm .hqIqy}, jour vrai}
\index{BJBHBE@\RL{ywm}!BJBHBE BHASAWBJ@\RL{ywm ws.ty}, jour moyen}
\index{BHASAW@\RL{ws.t}!BHASAW ATBEAS@\RL{ws.t al-^sams}, Soleil moyen}
\index{ANBDBA@\RL{_hlf}!ACANAJBDAGBA@\RL{i_htilAf}, irrégularité, anomalie, variation}
\index{BEAQAQ@\RL{mrr}!BEBEAQAQ@\RL{mamarr j At}, passage, transit (en général au méridien)}
\index{ACBHAL@\RL{'awj}!ACBHAL ATBEAS@\RL{'awj al-^sams}, Apogée du Soleil}
\addcontentsline{toc}{chapter}{II.8 Comment déterminer la durée des jours avec leurs nuits, et leurs équations}
\includepdf[pages=25,pagecommand={\thispagestyle{plain}}]{edit2.pdf}

\newpage\phantomsection
\begin{center}
  \Large Huitième section

  \normalsize Comment déterminer la durée des jours avec leurs nuits,
  et leurs équations
\end{center}

Le \emph{jour} désignera le nychtémère\footnote{littéralement, le
  <<~jour avec sa nuit~>>.}. Il y a le jour \emph{vrai} et le jour
\emph{moyen}. Le jour vrai est la durée entre le lever du Soleil, son
coucher ou son transit au méridien, d'une part, et son retour à ce
qu'on a posé comme origine parmi les trois, d'autre part, en raison du
premier mouvement ; sa grandeur est une circonférence de l'équateur,
plus l'[arc d'équateur] qui se lève en même temps que l'arc parcouru
par le Soleil en raison de son mouvement propre entre le début et la
fin [du jour].

Cet arc [de l'écliptique] est plus petit [quand le Soleil est] dans la
moitié la plus distante, et il est plus grand [quand le Soleil est]
dans la moitié la plus proche ; et l'arc d'équateur qui se lève en
même temps qu'un arc [donné] de l'écliptique est tantôt plus petite,
tantôt plus grande. En raison de ces deux variations, la grandeur d'un
nychtémère varie. L'écart est petit, donc cette variation est
insensible sur une durée d'un petit nombre de jours (un ou deux jours
par exemple), mais elle est sensible sur une durée d'un grand nombre
de jours.

Les astronomes sont contraints d'utiliser des nychtémères de
grandeurs égales pour déterminer les mouvements moyens et d'autres
choses dont ils ont besoin~; aussi ont-ils posé que
l'excès\footnote{\textit{i. e.} la différence entre la durée du jour
  et la période du mouvement de l'équateur.} est égal au mouvement
moyen du Soleil pendant un nychtémère, et ils ont appelé ces jours
``jours moyens''. Chacun est égal à une rotation de l'équateur, plus
le Soleil moyen.

Ils l'ont appelé ``moyen'' car il implique le parcours moyen [du
  Soleil], de même qu'ils ont appelé le premier ``jour vrai'' car il
implique le parcours vrai [du Soleil], c'est-à-dire son mouvement en
longitude. Les deux [durées ainsi définies] sont parfois égales,
parfois différentes, parce que la coascension de son mouvement en
longitude est parfois supérieure, parfois inférieure [au mouvement en
  longitude]~; outre ces deux cas, le mouvement en longitude est
tantôt égal au mouvement moyen (quand [le Soleil] est à une position
où son mouvement est moyen), tantôt supérieur (quand [le Soleil] est
dans la moitié [de sa trajectoire passant par] le périgée), tantôt
inférieur (quand [le Soleil] est dans la moitié [de sa trajectoire
  passant par] l'Apogée). On aura donc six divisions.

\newpage\phantomsection
\index{AYAOBD@\RL{`dl}!AJAYAOBJBD ACBJBJAGBE@\RL{ta`dIl al-'ayyAm}, équation des nychtémères}
\index{BFAUBA@\RL{n.sf}!AOAGAEAQAI BFAUBA BFBGAGAQ@\RL{dA'iraT n.sf al-nhAr}, méridien}
\includepdf[pages=26,pagecommand={\thispagestyle{plain}}]{edit2.pdf}
Les deux [durées] sont égales quand coascension et [Soleil] moyen sont
égaux, c'est-à-dire quand sont égaux les deux excédents ou défauts par
rapport au mouvement en longitude. L'écart [entre les deux durées]
s'appelle \emph{équation des nychtémères}. Elle n'est sensible que sur
une durée d'un grand nombre de jours -- pas un jour ni deux. Pour la
déterminer quantitativement, il faut connaître les deux variations.

L'écart dû [à la variation du] parcours du Soleil est quatre fois
l'anomalie maximale qui est environ deux degrés, parce
que le mouvement en longitude dans la moitié [passant par] l'Apogée
est égal au mouvement moyen moins deux fois le maximum de ce qu'on lui
ajoute, et qu'il en est de même dans la moitié [passant] par le
périgée. Il n'y a pas de contradiction entre le fait que l'arc en
longitude soit supérieur à l'arc moyen -- comme il arrive dans les
orbes du Soleil -- et le fait que le mouvement en longitude soit
inférieur au mouvement moyen -- comme nous venons de dire.

L'écart dû à la coascension dépend des horizons, donc ça n'est pas la
même chose pour toutes les contrées. Si l'on pose que les jours
commencent quand le Soleil arrive à l'horizon, de sorte que l'origine
est l'arrivée du Soleil à l'horizon Est, alors on compte l'écart entre
les degrés égaux et leur coascension par rapport au lieu [donné]. Si
l'on pose que l'origine est l'arrivée du Soleil à l'horizon Ouest,
alors on compte l'écart entre les degrés égaux et la coascension [de
  l'arc] opposé, par rapport au lieu [donné]. Si l'on pose que
l'origine est l'arrivée [du Soleil] au méridien, alors l'écart est le
même pour tous les horizons, et les jours ont la même durée dans tous
[les horizons]~: on compte le déplacement et la coascension par
rapport à l'équateur terrestre. Et c'est pourquoi le méridien a été
choisi comme origine.

Ci-dessus\footnote{\textit{cf.} section précédente.}, on a divisé
l'orbe de l'écliptique en quatre. Parmi les quatre portions, il y en a
deux qui ont chacune un équinoxe en son milieu et qui sont chacune
supérieure à sa coascension. L'une va du milieu du Verseau au milieu
du Taureau, l'autre va du milieu du Lion au milieu du
Scorpion. L'excédent de chacune sur sa coascension par rapport à
l'équateur terrestre est cinq degrés. Les deux autres portions ont
chacune un solstice en son milieu et elles sont chacune inférieure à
sa coascension. L'une va du milieu du Taureau au milieu du Lion, et
l'autre va du milieu du Scorpion au milieu de Verseau. Le défaut de
chacune à sa coascension par rapport à l'équateur terrestre est aussi
cinq degrés.

\newpage\phantomsection
\index{ASBFBH@\RL{snw}!ASBFAI ATBEAS@\RL{sanaT al-^sams}, année solaire}
\includepdf[pages=27,pagecommand={\thispagestyle{plain}}]{edit2.pdf}

En combinant les deux écarts par excès et par défaut, cette variation
fera la grandeur totale de l'écart entre les jours moyens et les jours
vrais. On ne peut faire autrement que choisir un jour de l'année comme
origine. Les autres jours lui seront rapportés, de sorte que midi de
ce jour soit l'origine des jours moyens ainsi que des jours
vrais. Quel que soit le jour de l'année choisi comme origine du temps,
l'écart entre les jours moyens et les jours vrais écoulés depuis ce
jour sera tantôt additif, tantôt soustractif ; mais si l'on prend
l'origine vers la fin du Verseau, alors les jours vrais seront
toujours inférieurs au jours moyens~; et si l'on prend l'origine vers
le début du Scorpion, alors les jours vrais seront toujours supérieurs
aux jours moyens. Les hommes de métier se sont accordés à choisir
[l'origine] vers la fin du Verseau.

Voici la figure des six divisions, sous l'hypothèse que l'Apogée est à
la fin des Gémeaux.

L'écart dû à l'anomalie solaire change à cause du mouvement de
l'Apogée, mais seulement sur de longues durées.

Ainsi s'achève l'explication de l'écart dans la grandeur des
jours. L'explication de [leur] grandeur à un instant donné [figure]
dans les livres pratiques. Cet écart s'appelle équation des
nychtémères. Au cours d'une révolution complète, [la durée totale des]
jours vrais et [la durée totale des] jours moyens sont égales, et la
considération précédente tombe.

\newpage\phantomsection
\begin{center}
  \small
\begingroup%
  \makeatletter%
  \providecommand\color[2][]{%
    \errmessage{(Inkscape) Color is used for the text in Inkscape, but the package 'color.sty' is not loaded}%
    \renewcommand\color[2][]{}%
  }%
  \providecommand\transparent[1]{%
    \errmessage{(Inkscape) Transparency is used (non-zero) for the text in Inkscape, but the package 'transparent.sty' is not loaded}%
    \renewcommand\transparent[1]{}%
  }%
  \providecommand\rotatebox[2]{#2}%
  \ifx\svgwidth\undefined%
    \setlength{\unitlength}{364.88588507bp}%
    \ifx\svgscale\undefined%
      \relax%
    \else%
      \setlength{\unitlength}{\unitlength * \real{\svgscale}}%
    \fi%
  \else%
    \setlength{\unitlength}{\svgwidth}%
  \fi%
  \global\let\svgwidth\undefined%
  \global\let\svgscale\undefined%
  \makeatother%
  \begin{picture}(1,1.03286227)%
    \put(0,0){\includegraphics[width=\unitlength]{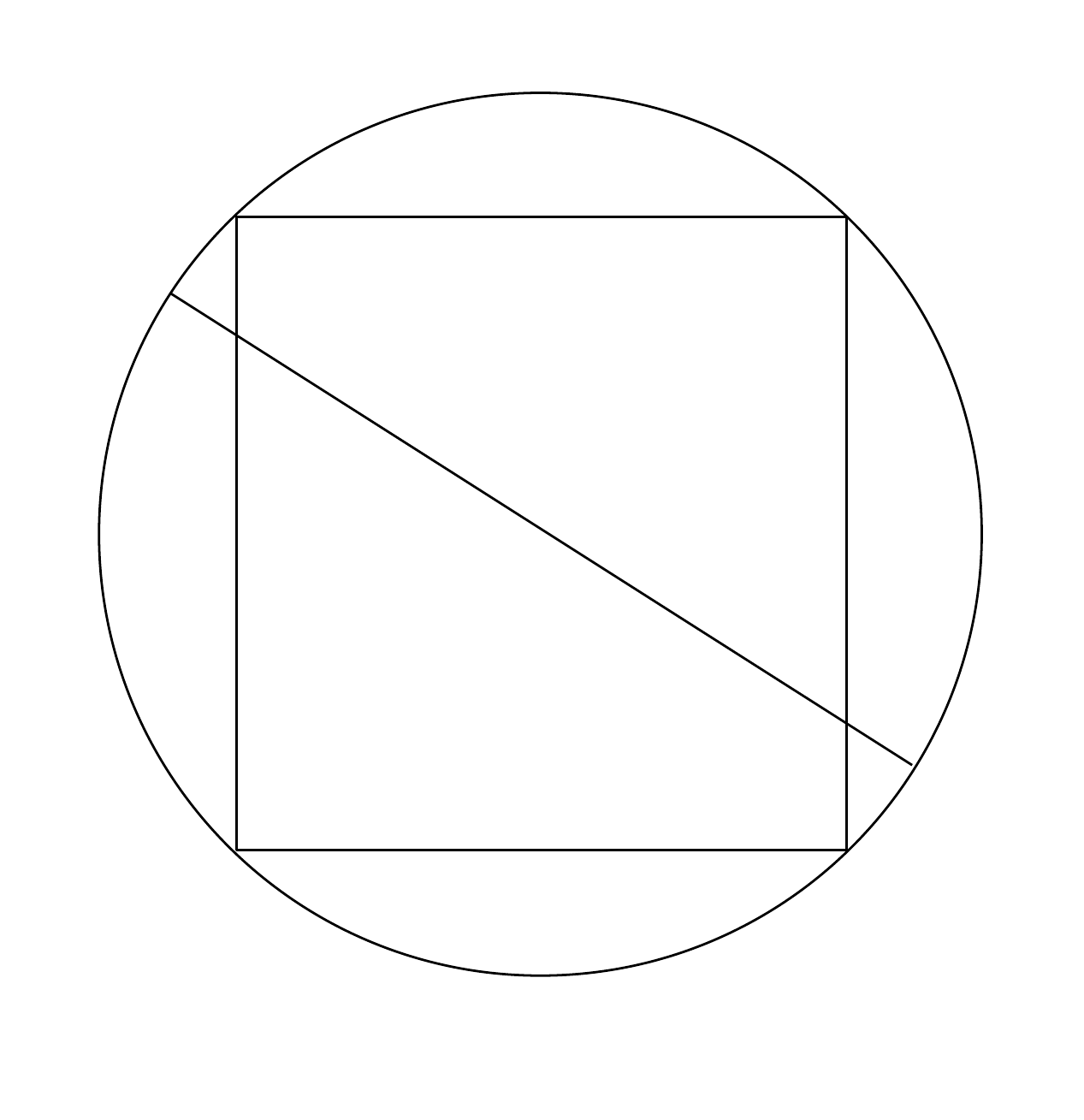}}%
    \put(0.2112979,0.23841569){\color[rgb]{0,0,0}\rotatebox{24.53529745}{\makebox(0,0)[rb]{\smash{\RL{ws.t al-`qrb}}}}}%
    \put(0.20365265,0.82880465){\color[rgb]{0,0,0}\rotatebox{-33.298821}{\makebox(0,0)[rb]{\smash{\RL{ws.t al-asad}}}}}%
    \put(0.14452474,0.75304672){\color[rgb]{0,0,0}\rotatebox{-28.46168182}{\makebox(0,0)[rb]{\smash{\RL{tarbI` al-awj}}}}}%
    \put(0.28131236,0.17452173){\color[rgb]{0,0,0}\rotatebox{31.92440824}{\makebox(0,0)[rb]{\smash{\RL{.hadId al-^sams}}}}}%
    \put(0.73498463,0.22057398){\color[rgb]{0,0,0}\makebox(0,0)[rb]{\smash{\RL{al-.hqIqyaT tzId `l_A al-ws.t_A bisabab al-m.tAl`}}}}%
    \put(0.6125275,0.48752617){\color[rgb]{0,0,0}\rotatebox{-31.50258743}{\makebox(0,0)[rb]{\smash{\RL{al-q.t`aT al-b`IdaT}}}}}%
    \put(0.59608109,0.4503291){\color[rgb]{0,0,0}\rotatebox{-32.21726869}{\makebox(0,0)[rb]{\smash{\RL{al-q.t`aT al-qrIbaT}}}}}%
    \put(0.74375851,0.84644537){\color[rgb]{0,0,0}\makebox(0,0)[rb]{\smash{\RL{al-.hqIqyaT tzId `l_A al-ws.t_A bisabab al-m.tAl`}}}}%
    \put(0.79008404,0.2394636){\color[rgb]{0,0,0}\rotatebox{-42.24247551}{\makebox(0,0)[lb]{\smash{\RL{ws.t al-dlw}}}}}%
    \put(0.78765049,0.8331687){\color[rgb]{0,0,0}\rotatebox{12.58834452}{\makebox(0,0)[lb]{\smash{\RL{ws.t al-_twr}}}}}%
    \put(0.51627765,0.11804792){\color[rgb]{0,0,0}\rotatebox{47.69252529}{\makebox(0,0)[rb]{\smash{\RL{awl al-jd_A}}}}}%
    \put(0.51168531,0.95472264){\color[rgb]{0,0,0}\rotatebox{22.74466614}{\makebox(0,0)[lb]{\smash{\RL{awl al-sr.tAn}}}}}%
    \put(0.70815817,0.89372063){\color[rgb]{0,0,0}\rotatebox{21.32351932}{\makebox(0,0)[lb]{\smash{\RL{awj al-^sams}}}}}%
    \put(0.80740052,0.79406091){\color[rgb]{0,0,0}\rotatebox{90}{\makebox(0,0)[rb]{\smash{\RL{al-.hqIqyaT tnq.s min al-ws.t_A bisabab al-m.tAl`}}}}}%
    \put(0.54352873,0.26105063){\color[rgb]{0,0,0}\rotatebox{-51.10087642}{\makebox(0,0)[rb]{\smash{\RL{al-.hqIqyaT tzId `l_A al-ws.t_A bisabab al-i_htilAf}}}}}%
    \put(0.19101864,0.29597814){\color[rgb]{0,0,0}\rotatebox{-90}{\makebox(0,0)[rb]{\smash{\RL{al-.hqIqyaT tnq.s min al-ws.t_A bisabab al-m.tAl`}}}}}%
    \put(0.85067846,0.31426666){\color[rgb]{0,0,0}\rotatebox{-28.46168182}{\makebox(0,0)[lb]{\smash{\RL{tarbI` al-awj}}}}}%
  \end{picture}%
\endgroup%
\normalsize
\end{center}

\begin{center}
  \small
\begingroup%
  \makeatletter%
  \providecommand\color[2][]{%
    \errmessage{(Inkscape) Color is used for the text in Inkscape, but the package 'color.sty' is not loaded}%
    \renewcommand\color[2][]{}%
  }%
  \providecommand\transparent[1]{%
    \errmessage{(Inkscape) Transparency is used (non-zero) for the text in Inkscape, but the package 'transparent.sty' is not loaded}%
    \renewcommand\transparent[1]{}%
  }%
  \providecommand\rotatebox[2]{#2}%
  \ifx\svgwidth\undefined%
    \setlength{\unitlength}{390.41563339bp}%
    \ifx\svgscale\undefined%
      \relax%
    \else%
      \setlength{\unitlength}{\unitlength * \real{\svgscale}}%
    \fi%
  \else%
    \setlength{\unitlength}{\svgwidth}%
  \fi%
  \global\let\svgwidth\undefined%
  \global\let\svgscale\undefined%
  \makeatother%
  \begin{picture}(1,0.96602366)%
    \put(0,0){\includegraphics[width=\unitlength]{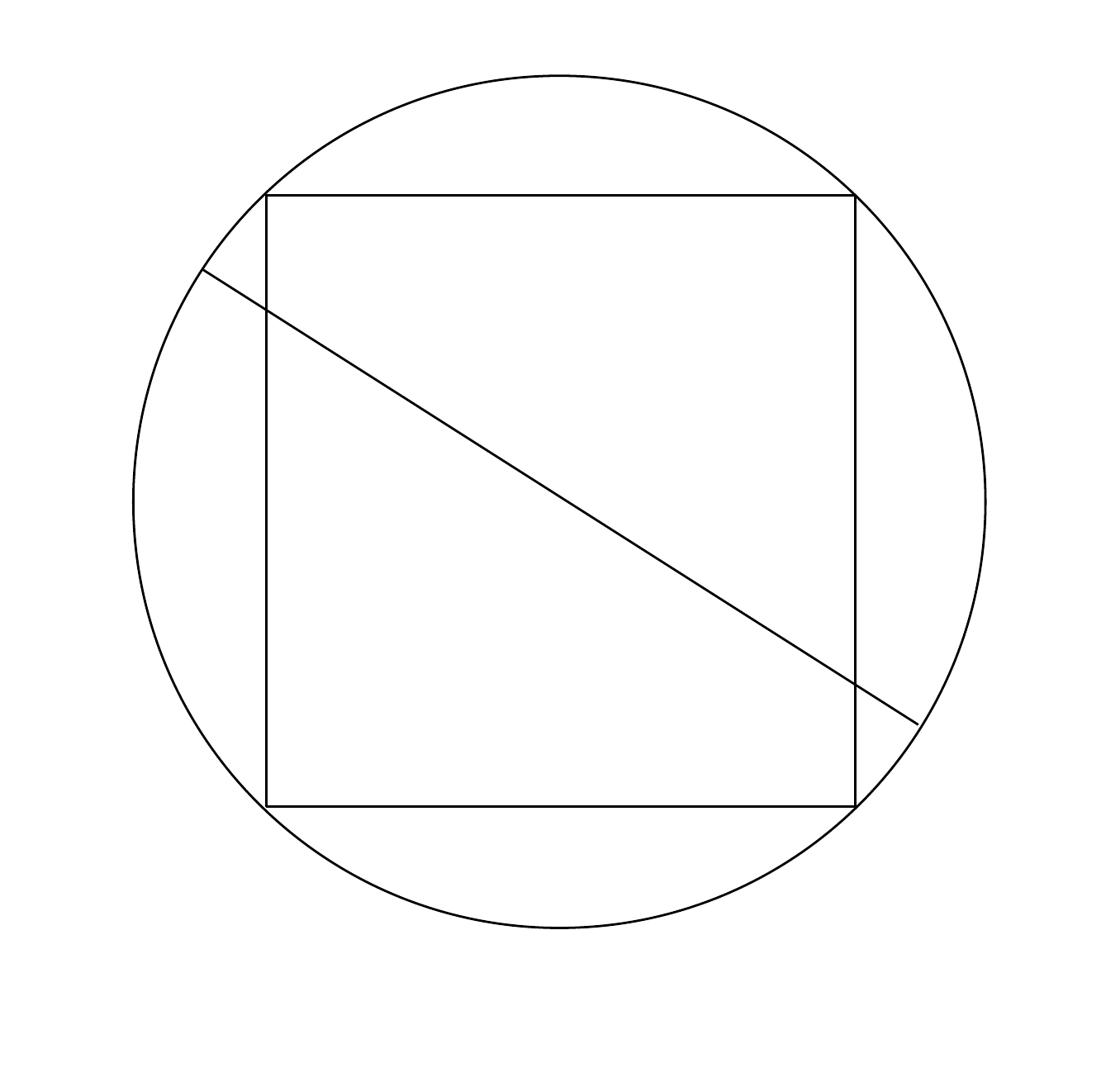}}%
    \put(0.23116334,0.23590004){\color[rgb]{0,0,0}\rotatebox{24.53529745}{\makebox(0,0)[rb]{\smash{milieu du Scorpion}}}}%
    \put(0.224018,0.78768279){\color[rgb]{0,0,0}\rotatebox{-33.298821}{\makebox(0,0)[rb]{\smash{milieu du Lion}}}}%
    \put(0.16875653,0.71687883){\color[rgb]{0,0,0}\rotatebox{-28.46168182}{\makebox(0,0)[rb]{\smash{quadrature de l'Apogée}}}}%
    \put(0.29659942,0.1761843){\color[rgb]{0,0,0}\rotatebox{31.92440824}{\makebox(0,0)[rb]{\smash{périgée du Soleil}}}}%
    \put(0.72060552,0.21922513){\color[rgb]{0,0,0}\makebox(0,0)[rb]{\smash{coascension $\Rightarrow$ jour vrai $>$ jour moyen}}}%
    \put(0.60615602,0.46872086){\color[rgb]{0,0,0}\rotatebox{-31.50258743}{\makebox(0,0)[rb]{\smash{la portion lointaine}}}}%
    \put(0.59078504,0.43395629){\color[rgb]{0,0,0}\rotatebox{-32.21726869}{\makebox(0,0)[rb]{\smash{la portion proche}}}}%
    \put(0.72880567,0.80417002){\color[rgb]{0,0,0}\makebox(0,0)[rb]{\smash{coascension $\Rightarrow$ jour vrai $>$ jour moyen}}}%
    \put(0.77210193,0.23687943){\color[rgb]{0,0,0}\rotatebox{-42.24247551}{\makebox(0,0)[lb]{\smash{milieu du Verseau}}}}%
    \put(0.7698275,0.79176142){\color[rgb]{0,0,0}\rotatebox{12.58834452}{\makebox(0,0)[lb]{\smash{milieu du Taureau}}}}%
    \put(0.51620006,0.12340336){\color[rgb]{0,0,0}\rotatebox{47.69252529}{\makebox(0,0)[rb]{\smash{début du Capricorne}}}}%
    \put(0.51190805,0.90536684){\color[rgb]{0,0,0}\rotatebox{22.74466614}{\makebox(0,0)[lb]{\smash{début du Cancer}}}}%
    \put(0.69553328,0.84835382){\color[rgb]{0,0,0}\rotatebox{21.32351932}{\makebox(0,0)[lb]{\smash{Apogée du Soleil}}}}%
    \put(0.78828605,0.75521103){\color[rgb]{0,0,0}\rotatebox{90}{\makebox(0,0)[rb]{\smash{coascension $\Rightarrow$ jour vrai $<$ jour moyen}}}}%
    \put(0.54166916,0.25705488){\color[rgb]{0,0,0}\rotatebox{-51.10087642}{\makebox(0,0)[rb]{\smash{anomalie $\Rightarrow$ jour vrai $>$ jour moyen}}}}%
    \put(0.21221016,0.28969853){\color[rgb]{0,0,0}\rotatebox{-90}{\makebox(0,0)[rb]{\smash{coascension $\Rightarrow$ jour vrai $<$ jour moyen}}}}%
    \put(0.82873398,0.30679121){\color[rgb]{0,0,0}\rotatebox{-28.46168182}{\makebox(0,0)[lb]{\smash{quadrature de l'Apogée}}}}%
  \end{picture}%
\endgroup%
\normalsize
\end{center}

\newpage\phantomsection
\index{ATAYAY@\RL{^s``}!ATAYAGAY@\RL{^su`A`}, rayon (de lumière)}
\index{BFBHAQ@\RL{nwr}!BFBHAQ@\RL{nwr}, lumière}
\index{AUAHAM@\RL{.sb.h}!AUAHAM@\RL{.sb.h}, aurore}
\index{ATBABB@\RL{^sfq}!ATBABB@\RL{^sfq}, crépuscule}
\index{AVBHAB@\RL{.daw'}!AVBHAB@\RL{.daw'}, clarté}
\index{AXBDBD@\RL{.zll}!BEANAQBHAW AXBDBD@\RL{ma_hrU.t al-.zill}, cône d'ombre}
\index{AUAHAM@\RL{.sb.h}!AUAHAM BCAGAPAH@\RL{.sb.h kA_dib}, aurore trompeuse}
\index{AUAHAM@\RL{.sb.h}!AUAHAM AUAGAOBB@\RL{.sb.h .sAdiq}, aurore franche}
\index{APBFAH@\RL{_dnb}!APBFAH ASAQAMAGBF@\RL{_danb al-sr.hAn}, queue du loup|see{\RL{.sb.h kA_dib}}}
\index{BFAXAQ@\RL{n.zr}!BFAXAQ@\RL{n.zr}, {\oe}il (de l'observateur)}
\index{ATBEAS@\RL{^sms}!ATBEAS@\RL{^sams}, Soleil}
\addcontentsline{toc}{chapter}{II.9 L'aurore et le crépuscule}
\includepdf[pages=28,pagecommand={\thispagestyle{plain}}]{edit2.pdf}

\begin{center}
  \Large Neuvième section

  \normalsize L'aurore et le crépuscule.
\end{center}

L'aurore est l'illumination de l'horizon Est à l'approche du Soleil,
et le crépuscule est l'illumination de l'horizon Ouest après que le
Soleil s'y soit retiré. La cause de l'aurore est la suivante~: pendant
la nuit, de dessous la Terre, quand le Soleil approche l'horizon Est,
le cône d'ombre de la Terre (en quoi consiste la nuit) penche vers
l'Ouest. Plus il penche, plus le rayon de lumière à son bord
s'approche du regard, du côté de l'Est, jusqu'à être vu, petit à
petit.

Pour se le représenter, imaginons un triangle à angles aigus dont la
base soit l'horizon et dont les côtés soient à la surface du cône dans
un plan passant par le centre du Soleil, le centre de la Terre et le
sommet du cône. Imaginons un autre triangle ayant un angle droit, dont
la base soit entre l'horizon et le \emph{premier point vu de l'aurore}
(c'est-à-dire un point du côté du premier triangle), et dont les deux
autres côtés soient issus de l'{\oe}il, l'un des deux passant par une
des extrémités de la base et l'autre par l'autre, l'une [des deux
  extrémités de la base] étant le point mentionné, et l'autre un point
de l'horizon. Comme [le côté] issu de l'{\oe}il et passant par le
point [mentionné] est perpendiculaire [en ce point], alors [l'autre]
côté [issu de l'{\oe}il] vers l'horizon est plus long car c'est
l'hypoténuse\footnote{Parmi les directions situées dans le même cercle
  de hauteur que le Soleil, celle l'observateur voit la clarté en
  premier est normale au cône d'ombre, car la distance à la surface du
  cône est minimale le long de la normale. C'est mieux expliqué chez
  \d{T}\=us{\=\i} \cite{altusi1993} p.~294-297.}. Ainsi ce qu'on voit
[en premier] de l'aurore est comme une droite située au dessus de
l'horizon~; ce qui est plus proche de l'horizon est [encore]
obscur. On appelle cela \emph{ce qui en est vu de l'aurore en premier}
parce que cela apparaît en premier à l'aurore, et \emph{aurore
  trompeuse} parce que cela [semble] être en contradiction avec
l'obscurité de ce qui est plus proche encore du Soleil. [On l'appelle
  aussi] \emph{queue du loup} à cause du déploiement de son début et
de la ponctualité de sa fin.

Ensuite, si le Soleil s'approche beaucoup de l'horizon, la lumière se
développe à l'horizon comme une demi-disque, et l'aurore devient alors
\emph{franche}. L'horizon Est s'emplit de clarté, la moitié du ciel est
atteinte, et cette clarté ne cesse plus de croître jusqu'à ce que
rougeoie l'horizon, que le Soleil se lève, et que le cône d'ombre se
couche. Plus le Soleil monte, plus le cône s'abaisse, jusqu'à midi ;
alors le Soleil commence à se diriger vers l'Ouest, et l'ombre et ses
limites vont vers l'Est, jusqu'à ce que le Soleil arrive à l'horizon
Ouest, qu'il le fasse rougir, et que le cône [d'ombre] arrive à
l'horizon Est. 

\newpage\phantomsection
\begin{center}
\begingroup%
  \makeatletter%
  \providecommand\color[2][]{%
    \errmessage{(Inkscape) Color is used for the text in Inkscape, but the package 'color.sty' is not loaded}%
    \renewcommand\color[2][]{}%
  }%
  \providecommand\transparent[1]{%
    \errmessage{(Inkscape) Transparency is used (non-zero) for the text in Inkscape, but the package 'transparent.sty' is not loaded}%
    \renewcommand\transparent[1]{}%
  }%
  \providecommand\rotatebox[2]{#2}%
  \ifx\svgwidth\undefined%
    \setlength{\unitlength}{414.53145082bp}%
    \ifx\svgscale\undefined%
      \relax%
    \else%
      \setlength{\unitlength}{\unitlength * \real{\svgscale}}%
    \fi%
  \else%
    \setlength{\unitlength}{\svgwidth}%
  \fi%
  \global\let\svgwidth\undefined%
  \global\let\svgscale\undefined%
  \makeatother%
  \begin{picture}(1,0.82623714)%
    \put(0,0){\includegraphics[width=\unitlength]{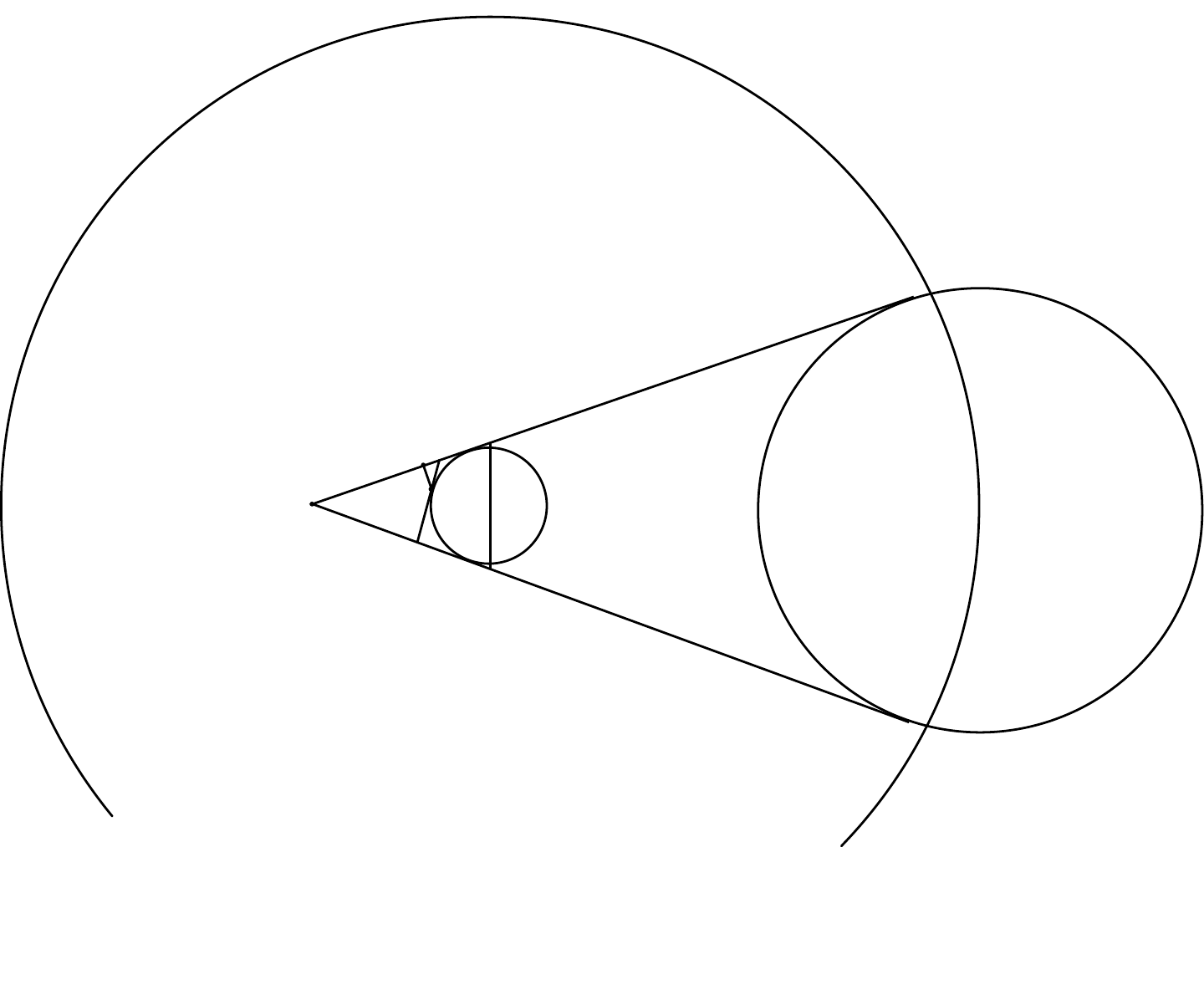}}%
    \put(0.25179414,0.40515976){\color[rgb]{0,0,0}\makebox(0,0)[rb]{\smash{\RL{ra's al-m_hrU.t}}}}%
    \put(0.29852407,0.38879388){\color[rgb]{0,0,0}\rotatebox{70.77317342}{\makebox(0,0)[rb]{\smash{\RL{m_tl_t al-m_hrU.t}}}}}%
    \put(0.34556334,0.36719478){\color[rgb]{0,0,0}\rotatebox{71.64733953}{\makebox(0,0)[rb]{\smash{\RL{s.t.h al-'ufuq al-mr_A}}}}}%
    \put(0.42882648,0.4413856){\color[rgb]{0,0,0}\rotatebox{90}{\makebox(0,0)[rb]{\smash{\RL{qA`daT al-m_hrU.t}}}}}%
    \put(0.66037691,0.3798044){\color[rgb]{0,0,0}\makebox(0,0)[lb]{\smash{\RL{al-^sams}}}}%
    \put(0.35145554,0.41132584){\color[rgb]{0,0,0}\rotatebox{-19.45913348}{\makebox(0,0)[rb]{\smash{\RL{mw.d` al-nA.zr}}}}}%
    \put(0.37470406,0.44685622){\color[rgb]{0,0,0}\rotatebox{74.85894298}{\makebox(0,0)[lb]{\smash{\RL{.trf al-'ufuq whw mw.d` istiqbAlihi bi-al-.dl`}}}}}%
    \put(0.1301665,0.12476204){\color[rgb]{0,0,0}\rotatebox{-58.51254587}{\makebox(0,0)[rb]{\smash{\RL{qws al-^sams}}}}}%
    \put(0.39870835,0.00416509){\color[rgb]{0,0,0}\makebox(0,0)[b]{\smash{\RL{w h_dh .sUraT al-'ufuq wa-al-m_tl_t wa-al-`mUd wa-al-^sams wa-al-'ar.d}}}}%
  \end{picture}%
\endgroup%

\end{center}

\newpage\phantomsection
\begin{center}
\begingroup%
  \makeatletter%
  \providecommand\color[2][]{%
    \errmessage{(Inkscape) Color is used for the text in Inkscape, but the package 'color.sty' is not loaded}%
    \renewcommand\color[2][]{}%
  }%
  \providecommand\transparent[1]{%
    \errmessage{(Inkscape) Transparency is used (non-zero) for the text in Inkscape, but the package 'transparent.sty' is not loaded}%
    \renewcommand\transparent[1]{}%
  }%
  \providecommand\rotatebox[2]{#2}%
  \ifx\svgwidth\undefined%
    \setlength{\unitlength}{414.53145082bp}%
    \ifx\svgscale\undefined%
      \relax%
    \else%
      \setlength{\unitlength}{\unitlength * \real{\svgscale}}%
    \fi%
  \else%
    \setlength{\unitlength}{\svgwidth}%
  \fi%
  \global\let\svgwidth\undefined%
  \global\let\svgscale\undefined%
  \makeatother%
  \begin{picture}(1,0.89970455)%
    \put(0,0){\includegraphics[width=\unitlength]{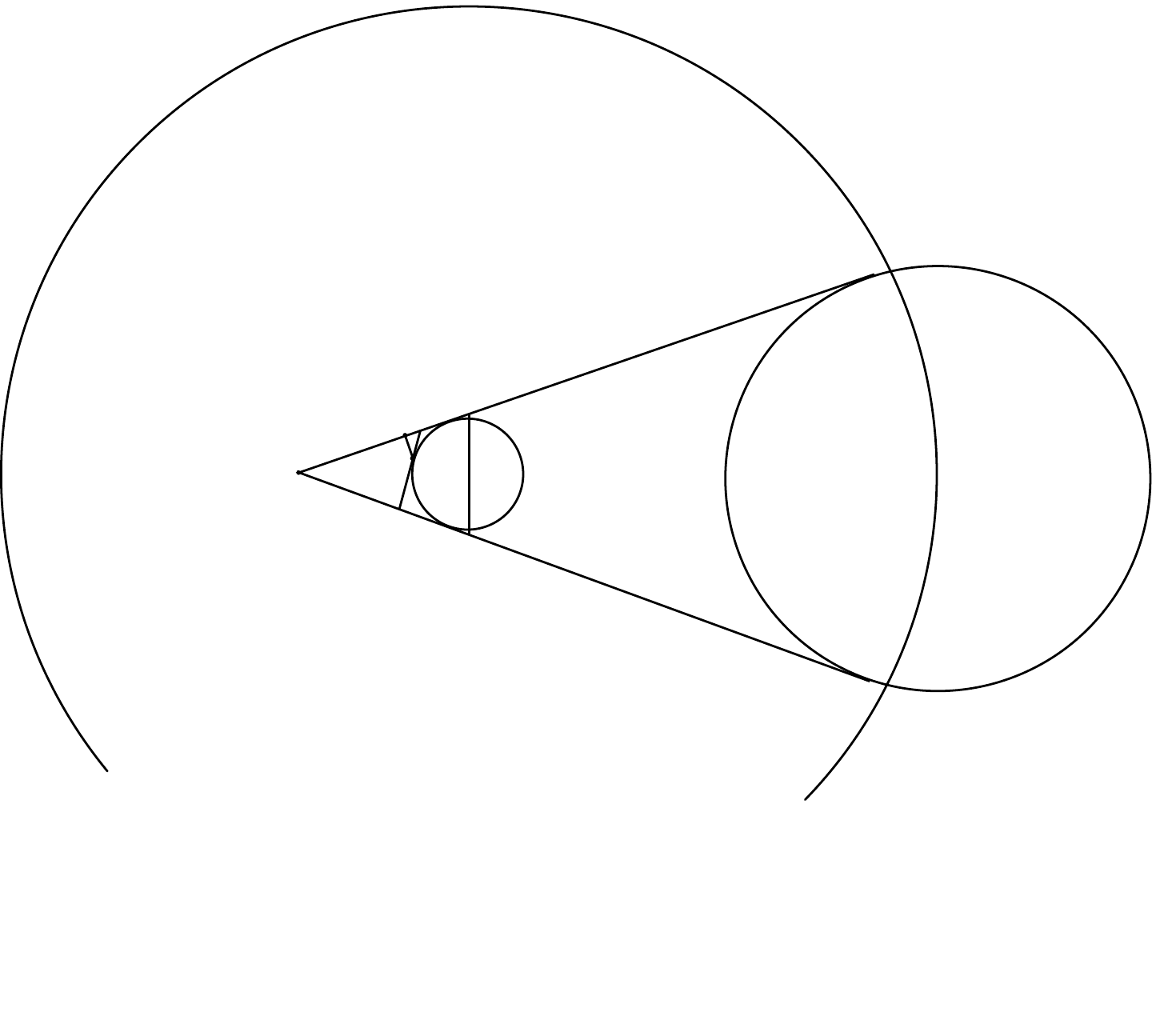}}%
    \put(0.25179414,0.48703802){\color[rgb]{0,0,0}\makebox(0,0)[rb]{\smash{sommet du cône}}}%
    \put(0.29852407,0.47067214){\color[rgb]{0,0,0}\rotatebox{70.77317342}{\makebox(0,0)[rb]{\smash{triangle [dans le plan] du cône}}}}%
    \put(0.34556334,0.44907304){\color[rgb]{0,0,0}\rotatebox{71.64733953}{\makebox(0,0)[rb]{\smash{plan de l'horizon apparent}}}}%
    \put(0.42882648,0.52326386){\color[rgb]{0,0,0}\rotatebox{90}{\makebox(0,0)[rb]{\smash{base du cône}}}}%
    \put(0.66037691,0.46168266){\color[rgb]{0,0,0}\makebox(0,0)[lb]{\smash{Soleil}}}%
    \put(0.35145554,0.4932041){\color[rgb]{0,0,0}\rotatebox{-19.45913348}{\makebox(0,0)[rb]{\smash{position de l'observateur}}}}%
    \put(0.37470406,0.52873448){\color[rgb]{0,0,0}\rotatebox{74.85894298}{\makebox(0,0)[lb]{\smash{extrémité de l'horizon à l'intersection avec le côté}}}}%
    \put(0.1301665,0.2066403){\color[rgb]{0,0,0}\rotatebox{-58.51254587}{\makebox(0,0)[rb]{\smash{arc d'écliptique}}}}%
    \put(0.39870835,0.00416509){\color[rgb]{0,0,0}\makebox(0,0)[b]{\smash{La figure avec l'horizon, le triangle, la perpendiculaire, le Soleil et la Terre}}}%
  \end{picture}%
\endgroup%

\end{center}

\newpage\phantomsection
\includepdf[pages=29,pagecommand={\thispagestyle{plain}}]{edit2.pdf}
\noindent Si le Soleil se couche et que se lève le cône d'ombre,
alors le rougeoiment appelé crépuscule commence à décroître. \`A sa
suite, la blancheur croît, jusqu'à ce que la blancheur soit complète
et que la rougeur cesse, puis la blancheur
commence à décroître, jusqu'à ce qu'il n'en reste plus qu'une clarté
prolongée, comme l'aurore trompeuse. Puis commence la nuit.

Les gens sont rendus favorables au dîner et au sommeil à l'injonction
du crépuscule, et au déploiement dans l'accomplissement des besoins à
l'aurore. Injonction du crépuscule~: qu'ils se distraient des
[travaux] pointilleux. Conditions de l'aurore~: qu'ils s'occupent de
[travaux] pointilleux.

On a appris empiriquement que l'abaissement du Soleil sous l'horizon
au commencement du lever de l'aurore et à la fin du coucher du
crépuscule est de dix-huit degrés [mesurés] sur le cercle de hauteur
passant par le centre du Soleil~; mais comme la coascension de [cet
  arc] est variable, la [durée] en heures entre le lever de l'aurore
et le lever du Soleil est variable, de même que [la durée en heures]
entre le coucher du Soleil et le coucher du crépuscule. \`A l'équateur
terrestre, quand le Soleil est à l'un des deux équinoxes, cette
grandeur de dix-huit degrés [équivaut] à une durée d'une heure et un
cinquième d'heure. Quand le Soleil est à l'un des deux équinoxes, il
n'y aucun lieu sur Terre où la durée de l'aurore ou du crépuscule soit
inférieure à cela. \`A l'équateur terrestre, en deux positions [de
  l'écliptique] équidistantes d'un des deux équinoxes, les durées de
l'aurore ou du crépuscule sont égales.

La durée de l'aurore et du crépuscule dans le quatrième climat est
deux heures quand le Soleil est vers le début du Cancer, et c'est
environ une heure et un tiers d'heure quand il est vers le début du
Capricorne.

Aux lieux dont la latitude vaut $48;30$, quand le Soleil est au
solstice d'été, l'aurore est contigüe au crépuscule. Là où la latitude
est encore supérieure, leur contigüité [a lieu] selon que
l'abaissement [maximal du Soleil] sous l'horizon est inférieur à cette
grandeur, et alors le lever de l'aurore est avant le coucher du
crépuscule\footnote{Pendant la période de l'année où l'abaissement
  maximal du Soleil sous l'horizon est inférieur à $18°$, aurore et
  crépuscule seront contigus.}.

Lorsque la latitude [du lieu] est égale au complément de l'inclinaison
[de l'écliptique], au solstice situé du côté opposé à la latitude, le
Soleil culmine quand il est à l'horizon~; il ne se lève pas, et la durée [de
  l'aurore et du crépuscule] est cinq heures et un tiers d'heure. Ce
qui reste des vingt-quatre heures est treize heures et un tiers~:
c'est la durée de l'obscurité. Le Soleil touche [aussi] l'horizon [au
  solstice] situé du même côté que la latitude~; alors il ne disparaît
pas, et il n'y a donc ni aurore ni crépuscule à ces dates.

\newpage\phantomsection
\includepdf[pages=30,pagecommand={\thispagestyle{plain}}]{edit2.pdf}

Lorsque la latitude [du lieu] est $90$ (à l'horizon où c'est comme une
meule), la durée de chacun est de $50$ jours -- de nos jours à nous.

L'explication des grandeurs que nous avons mentionnées et de celles
que nous n'avons pas mentionnées figure dans les livres pratiques.

\newpage\phantomsection
\index{BFBGAQ@\RL{nhr}!AJAYAOBJBD BFBGAGAQ@\RL{ta`dIl al-nahAr}, équation du jour}
\index{BFBGAQ@\RL{nhr}!BBBHAS BFBGAGAQ@\RL{qws al-nhAr}, arc diurne, durée du jour}
\index{BDBJBD@\RL{lyl}!BBBHAS BDBJBD@\RL{qaws al-layl}, arc nocturne, durée de la nuit}
\index{ASBHAY@\RL{sw`}!ASAGAYAGAJ BEAYBHALAI@\RL{sA`At ma`UjaT}, heures courbes}
\index{ASBHAY@\RL{sw`}!ASAGAYAGAJ BEASAJBHBJAI@\RL{sA`At mustawiyaT}, heures égales}
\index{AUAHAM@\RL{.sb.h}!AUAHAM AUAGAOBB@\RL{.sb.h .sAdiq}, aurore franche}
\index{ATBGAQ@\RL{^shr}!ATBGAQ BBBEAQBJBJ BHASAW@\RL{^shr qamariyy}, mois lunaire}
\index{ATBCBD@\RL{^skl}!AJATBCBCBDAGAJ@\RL{t^skkulAt}, phases (de la Lune)}
\index{ATBEAS@\RL{^sms}!ATBEAS@\RL{^sams}, Soleil}
\index{BBBEAQ@\RL{qmr}!BBBEAQ@\RL{qamar}, Lune}
\index{AWBDAY@\RL{.tl`}!BEAWAGBDAY@\RL{m.tAl`}, coascension}
\index{BJBHBE@\RL{ywm}!BJBHBE AHBDBJBDAJBG@\RL{al-yawm bilaylatihi j al-'ayyAm bilayAlIhA}, nychtémère (\textit{litt.} le jour avec sa nuit)}
\addcontentsline{toc}{chapter}{II.10 Détermination des parties du jour, les heures, et des
  multiples des heures, les mois et les années}
\includepdf[pages=31,pagecommand={\thispagestyle{plain}}]{edit2.pdf}

\begin{center}
  \Large Dixième section

  \normalsize Détermination des parties du jour, les heures, et des
  multiples des heures, les mois et les années
\end{center}

Il est bien connu que l'\emph{arc diurne} est la moitié d'une
circonférence de l'équateur si l'équation est nulle, et qu'il faut y
ajouter ou en retrancher le double de l'équation si elle n'est pas
nulle. En vérité, [l'arc diurne devrait mesurer] combien tourne
l'équateur à partir de l'instant où s'est levé à l'horizon la moitié
du corps du Soleil jusqu'à l'instant où s'est couché la moitié~; et
selon cette interprétation, il est supérieur à la grandeur ci-dessus,
d'un excédent égal à la coascension de l'arc parcouru par le Soleil le
long de l'écliptique pendant ce jour par rapport à cette
contrée. L'arc nocturne dépend [aussi] de ce qu'on vient de dire
concernant l'arc diurne~; et l'équation du jour exigée est la moitié
de l'excédent de l'arc diurne sur une demi-circonférence de
l'équateur.

Tantôt, on divise chacun des deux arcs en quinze portions, et on
appelle ces portions \emph{heures égales} à cause du fait que leurs
parts sont toujours [en nombre] égal. Tantôt, [on les divise] en douze
portions, et on appelle ces portions \emph{heures courbes} à cause du
fait que leurs parts sont [en nombre] variable et dépendent de la
variation des deux temps\footnote{Les deux <<~temps~>> sont le jour et
  la nuit.} suivant qu'ils sont courts ou longs. Si les deux temps
sont égaux, chacun sera douze heures, égales et courbes. Si les deux
temps diffèrent, alors le nombre d'heures égales augmente et [le
  nombre] de parts de chaque heure courbe augmente, car l'une est plus
que l'autre.

Dans l'esprit des astronomes, le début du jour et de l'arc diurne est
au lever du Soleil, et c'est le point de vue naturel. Dans l'esprit
des légistes, c'est au lever de l'aurore franche ; alors [l'arc
  diurne] dépasse l'autre, de la différence entre les deux levers. Le
début de la nuit est au coucher du Soleil. Chez certains
astronomes le début du nychtémère est à midi, chez d'autres
c'est à minuit, chez d'autres encore qui sont légistes c'est au début
de la nuit, et chez d'autres c'est au début du jour ; [il dure]
jusqu'au prochain [événement] semblable.

Quand on conçoit le mois à partir des phases de la Lune selon sa
position relative au Soleil, et selon que l'excès du mouvement de la
Lune sur le mouvement du Soleil -- en mouvements vrais -- fait une
circonférence, alors son nombre [de jours] est variable, à cause de
l'irrégularité des deux mouvements.

\newpage\phantomsection
\index{BGBDBD@\RL{hll}!BGBDAGBD@\RL{hilAl}, croissant (de Lune)}
\index{BCAHAS@\RL{kbs}!BJBHBE ASBFAI BCAHBJASAI@\RL{yawm / sanaT kabIsaT}, jour / année intercalaire}
\index{BJBHBE@\RL{ywm}!BJBHBE BCAHBJASAI@\RL{yawm kabIsaT / mu^stariqaT / lawA.hiq}, jour intercalaire / volé / supplémentaire}
\index{ASBFBH@\RL{snw}!ASBFAI ATBEAS@\RL{sanaT al-^sams}, année solaire}
\index{ASBFBH@\RL{snw}!BAAUBHBD ASBFAI@\RL{fu.sUl al-sanaT}, les saisons}
\index{AYAOBD@\RL{`dl}!AYAJAOAGBD AQAHBJAYBJBJ@\RL{i`tidAl rabI`iyy}, équinoxe de printemps}
\includepdf[pages=32,pagecommand={\thispagestyle{plain}}]{edit2.pdf}

\noindent \uwave{Les observateurs} disent que la période va d'un jour
de conjonction au [suivant], ou bien d'une nuit de visibilité du
croissant à la [suivante], ou bien de la dernière phase [de la Lune] à
la [dernière phase] suivante, par convention. Les astronomes
parlent de la différence entre les deux mouvements \emph{moyens}, au
sens où ils soustraient le Soleil moyen de la Lune moyenne, et ils
trouvent ainsi que l'intervalle entre deux conjonctions est vingt-neuf
jours et demi et une [petite] fraction. Puisque deux mois successifs
réunis font cinquante-neuf jours, ils posent qu'un mois [sur deux]
vaut trente jours et l'autre vingt-neuf jours. Les fractions
accumulées en sus de [chaque période de vingt-neuf jours et] demi
ajoutent, tous les trente ans, onze jours. Onze mois parmi les mois
devant compter vingt-neuf jours, tous les trente ans, deviendront
ainsi onze mois de trente [jours]. Ces jours s'appellent des jours
\emph{intercalaires}~: un jour intercalaire est un jour ajouté à une
année et constitué des fractions en sus des jours entiers, et ils
ajoutent les jours intercalaires à [certains] mois.


Autrement dit, ils posent que le premier mois de l'année
(\textit{mu\d{h}arram}) fait trente jours, que le deuxième fait
vingt-neuf jours, et ainsi de suite jusqu'à la fin de l'année, mais
\textit{\b{d}\=u al-\d{h}ijja} compte alors vingt-neuf jours
et un cinquième et un sixième de jour, c'est-à-dire $0;22$ jour,
résultat de la multiplication par douze du $0;1;50$ en sus du
demi-jour. [Ce mois] devient [un mois de] trente [jours] lors de
l'\emph{année intercalaire} qui est l'année comptant le jour
intercalaire, trentième jour de \textit{\b{d}\=u
  al-\d{h}ijja}, constitué des fractions [accumulées]. C'est ainsi que
sont les mois {lunaires}~; on en tire les mois vrais et les mois
moyens.

L'année est le retour du Soleil en son lieu sur l'écliptique. Ce
retour implique le retour des conditions de l'année selon les
saisons. Ceci se produit en trois cent soixante-cinq jours un quart,
et une fraction. [L'année] consiste en douze mois parmi les mois
lunaires, et onze jours, sans compter les fractions. On considère
souvent que l'année va du jour où le Soleil arrive à l'équinoxe de
printemps jusqu'au [jour où il revient au] même point. Ses mois
commencent aux jours où il arrive aux points semblables de
l'écliptique~; ou bien les mois comptent trente jours, trente jours,
etc. mais ils ajoutent au dernier cinq ou six jours. Les cinq jours
sont dits \emph{volés} ou \emph{supplémentaires}, et le sixième est
dit \emph{intercalaire} (comme on dit de l'année). L'année [ainsi
  conçue] est l'année \emph{solaire vraie}, et ces mois sont des mois
solaires vrais, ou bien des mois conventionnels qu'on peut faire
commencer à n'importe quel jour ne coïncidant peut-être avec aucune
observation de la position du Soleil. On choisit des mois aux
alentours de trente jours parce que les mois lunaires sont proches de
cela.

\newpage\phantomsection
\index{ASBFBH@\RL{snw}!ASBFAI BBBEAQBJAI@\RL{sanaT qamariyaT}, année lunaire}
\includepdf[pages=33,pagecommand={\thispagestyle{plain}}]{edit2.pdf}

Il arrive souvent qu'on prenne la fraction en sus des trois cent
soixante-cinq jours égale à un quart exactement~: on intercale un jour
tous les quatre ans. Il arrive souvent qu'on l'élimine complètement,
et ces années sont les années solaires conventionnelles. Ceux qui
prennent en considération les mois lunaires posent que l'année est
solaire, et, tous les trois ans ou tous les deux ans, ils ajoutent un
mois à l'année pour rassembler les [intervalles] de onze jours -- sans
les fractions -- suivant leur convention.

Certains peuples posent que douze mois lunaires font une année, et ils
l'appellent \emph{année lunaire}. Chaque peuple a une origine à
laquelle il rapporte toutes les années de son histoire. La
connaissance détaillée de ceci relève des livres de pratique.

\newpage\phantomsection
\addcontentsline{toc}{chapter}{II.11 Les degrés de transit des astres au méridien, leurs levers, et leurs couchers}
\index{BEAQAQ@\RL{mrr}!AOAQALAI BEBEAQAQ@\RL{darjaT mamarr}, degré de transit}
\index{BEBJBD@\RL{myl}!AOAGAEAQAI BEBJBD@\RL{dA'iraT al-mIl}, cercle de déclinaison}
\index{AOBHAQ@\RL{dwr}!AOAGAEAQAI AYAQAV@\RL{dA'iraT al-`r.d}, cercle de latitude}
\index{BEAQAQ@\RL{mrr}!AEANAJBDAGBA BEBEAQAQ@\RL{i_htilAf al-mamarr}, anomalie du transit}
\index{AYAOBD@\RL{`dl}!AJAYAOBJBD AOAQALAI BEBEAQAQ@\RL{ta`dIl darjaT al-mamarr}, équation du degré de transit}
\index{BFAUBA@\RL{n.sf}!AOAGAEAQAI BFAUBA BFBGAGAQ@\RL{dA'iraT n.sf al-nhAr}, méridien}
\includepdf[pages=34,pagecommand={\thispagestyle{plain}}]{edit2.pdf}

\begin{center}
  \Large Onzième section

  \normalsize Les degrés de transit des astres au méridien, leurs levers, et leurs couchers
\end{center}

Le \emph{degré de transit} de l'astre est le [point] du cercle de
l'écliptique qui traverse le méridien en même temps que l'astre le
traverse, et le cercle de déclinaison [de l'astre] indique [ce
  point]. On a déjà montré que le degré de longitude est indiqué par
le cercle de latitude~; donc si ces deux cercles sont confondus, comme
il arrive lorsque l'astre est sur le cercle passant par les quatre
pôles, ou bien si [l'astre] a une latitude nulle, alors [le degré de
  transit et le degré de longitude] coïncident. Ailleurs qu'en ces
deux positions, ils diffèrent.

La différence s'appelle \emph{anomalie du transit}, et l'arc de
l'équateur situé entre son intersection avec le cercle de latitude de
l'astre et son intersection avec son cercle de déclinaison s'appelle
\emph{équation du degré de transit}. L'anomalie est maximale aux
     [positions] proches des commencements du Bélier et de la Balance,
     et elle est minimale aux [positions] proches des commencements du
     Cancer et du Capricorne.

Si le pôle de l'écliptique est sur le méridien, c'est-à-dire quand les
deux points des solstices [y] sont aussi et quand les deux points des
équinoxes sont à l'horizon, alors le passage de l'astre [au méridien]
se fait en même temps que [celui du point indiquant] sa longitude,
parce que le méridien est son cercle de latitude.

Si celui de deux pôles [de l'écliptique] qui est visible est à l'Est
du méridien -- ceci a lieu quand la moitié de l'écliptique centrée en
l'équinoxe d'automne passe [par le méridien] et quand sa moitié Sud se
lève, si le pôle visible est au Nord, ou bien quand la moitié de
l'écliptique centrée en l'équinoxe de printemps passe [par le
  méridien] et quand sa moitié Nord se lève, si le pôle visible est au
Sud -- alors tout astre dont la latitude est du côté du pôle visible
passe au méridien \emph{après} son degré de transit. En effet, son
cercle de latitude issu du pôle le plus proche rencontre au méridien
l'astre en premier, puis il rencontre son degré de transit seulement
après avoir traversé le méridien et être passé à l'Ouest~;
[inversement], si son degré de transit vient au méridien, l'astre est
encore du côté Est [du méridien]. [Dans les mêmes conditions], tout
astre dont la latitude est du côté opposé au pôle visible passe au
méridien \emph{avant} son degré de transit. En effet, le cercle de
latitude en question rencontre au méridien le degré [de transit] de
l'astre en premier, puis il rencontre l'astre seulement après avoir
traversé le méridien et être passé à l'Ouest.

\newpage\phantomsection
\index{AWBDAY@\RL{.tl`}!AWBDBHAY@\RL{.tulU`}, lever (d'un astre)}
\index{AZAQAH@\RL{.grb}!AZAQBHAH@\RL{.gurUb}, coucher (d'un astre)}
\index{AZAQAH@\RL{.grb}!AOAQALAI AZAQBHAH@\RL{darjaT al-.gurUb}, degré du coucher}
\index{AWBDAY@\RL{.tl`}!AOAQALAI AWBDBHAY@\RL{darjaT al-.tulU`}, degré du lever}
\index{BABHBB@\RL{fwq}!ADBABB@\RL{'ufuq}, horizon}
\includepdf[pages=35,pagecommand={\thispagestyle{plain}}]{edit2.pdf}

Si le pôle visible est à l'Ouest -- ceci a lieu quand la moitié de
l'écliptique centrée en l'équinoxe de printemps passe [par le
  méridien] et quand sa moitié Nord se lève, si le pôle visible est au
Nord, ou bien quand la moitié de l'écliptique centrée en l'équinoxe
d'automne passe [par le méridien] et quand sa moitié Sud se lève, si
le pôle visible est au Sud -- alors tout astre dont la latitude est du
côté du pôle visible passe au méridien \emph{avant} son degré de
transit, et [tout astre] dont la latitude est du côté opposé passe [au
  méridien] \emph{après} son degré de transit, pour les mêmes raisons
que dans le paragraphe précédent.

Le degré du lever ou du coucher d'un astre est [le point] du cercle de
l'écliptique qui se lève ou qui se couche en même temps que cet
astre. Si l'astre a une latitude nulle, alors son degré de longitude
est égal à son degré de lever et de coucher. Il en est de même d'un
astre qui vient à l'horizon en même temps que le pôle de
l'écliptique~: l'horizon est alors son cercle de latitude puisqu'il
passe par l'astre et par le pôle de l'écliptique.

Si l'horizon est à latitude nulle, comme il arrive à l'équateur
terrestre, alors le lever et le coucher d'un astre sont comme les
transits au méridien pour les autres horizons. [Tout astre] qui vient
à l'horizon en même temps que le pôle et le solstice se lève et se
couche en même temps que son degré [de longitude], comme on vient de
le dire. [Tout astre] situé du côté du pôle visible de l'écliptique se
lève avant son degré [de longitude] et se couche après. [Tout astre]
situé du côté du pôle caché se lève après son degré [de longitude] et
se couche avant. En effet, le cercle de latitude issu du pôle visible
atteint l'astre situé du même côté, [quand l'astre est] à l'horizon,
avant [d'atteindre] son degré de longitude~; et [il atteint] son degré
de longitude, [quand celui-ci est] à l'horizon, avant [d'atteindre]
l'astre situé du côté du pôle caché. Là, le pôle Nord [de
  l'écliptique] est visible pendant le lever de la moitié centrée en
l'équinoxe de printemps et pendant le passage de la moitié Sud au
méridien au dessus [de l'horizon]. Le pôle Sud est visible pendant le
lever de la moitié centrée en l'équinoxe d'automne et pendant le
passage de la moitié Nord au méridien.

Dans les horizons inclinés, la règle du lever et du coucher des astres
est comme nous l'avons décrite pour l'équateur terrestre. Quant au
passage [au méridien] des moitiés de l'écliptique, il y a certes une
différence. Il est possible qu'un des deux pôles [de l'écliptique]
soit visible, et que [la partie] qui transite ou qui se lève soit un
arc inférieur à une moitié, ou bien supérieur. La règle ci-dessus suit
son cours quand la latitude du lieu est supérieure à l'inclinaison [de
  l'écliptique], sans avoir à distinguer [de cas], parce qu'un des
deux pôles de l'écliptique est alors toujours visible.

\newpage\phantomsection
\includepdf[pages=36,pagecommand={\thispagestyle{plain}}]{edit2.pdf}

\noindent La règle reste
identique là où la latitude n'est pas supérieure~: si celui des deux
pôles de l'écliptique situé au Nord est visible, alors la règle passe,
mais s'il est caché, alors c'est la règle inverse. En effet, l'astre
se lève après son degré [de longitude] et il se couche avant son degré
[de longitude] si sa latitude est au Nord, et inversement si sa
latitude est au Sud.

Après réflexion, il ne vous échappera pas que le degré du lever de
l'astre se lève de jour s'il est entre le Soleil et le point opposé
[au Soleil], et qu'il se lève de nuit s'il est entre le point opposé
et le Soleil. Le degré de coucher se couche de nuit s'il est entre ces
deux-là, et il se couche de jour s'il est entre ces deux-ci. Parmi les
astres situés sur un même grand cercle\footnote{sous-entendu, un
  cercle de latitude.} coupé par la trajectoire diurne toujours
visible maximale, ceux qui sont plus proches du pôle se lèvent avant
ceux qui sont plus loins, et ils se couchent après ceux qui sont plus
loins. C'est pourquoi l'écart entre les degrés de longitude et de
lever est plus grand pour un astre proche du pôle que pour un astre
loin [du pôle].

\newpage\phantomsection
\index{BFAUBA@\RL{n.sf}!ANAWAW BFAUBA BFBGAGAQ@\RL{_ha.t.t n.sf al-nahAr}, ligne méridienne}
\index{BFAUBA@\RL{n.sf}!AOAGAEAQAI BFAUBA BFBGAGAQ@\RL{dA'iraT n.sf al-nhAr}, méridien}
\index{ATAQBB@\RL{^srq}!ANAWAW BEATAQBB BEAZAQAH@\RL{_ha.t.t al-ma^sraq wa-al-ma.grab}, ligne de l'Est et de l'Ouest}
\index{ASBEAJ@\RL{smt}!AOAGAEAQAI ACBHBHBD ASBEBHAJ@\RL{dA'iraT 'awwal al-sumUt}, cercle origine des azimuts}
\index{BBAHBD@\RL{qbl}!BBAHBDAI@\RL{qiblaT}, \textit{qibla}, direction de la \textit{Ka`ba}}
\index{BBBJAS@\RL{qys}!BEBBBJAGAS@\RL{miqyAs j maqAyIs}, gnomon}
\index{BGBFAO@\RL{hnd}!AOAGAEAQAI BGBFAOBJAI@\RL{dA'iraT hindiyaT}, cercle indien}
\index{ASBEAJ@\RL{smt}!ASBEAJ@\RL{smt}, direction, azimut}
\index{BEBCBCAI@\RL{makkaT}, La Mecque}
\addcontentsline{toc}{chapter}{II.12 Détermination du méridien du lieu et de l'azimut de la \textit{qibla}}
\includepdf[pages=37,pagecommand={\thispagestyle{plain}}]{edit2.pdf}

\begin{center}
  \Large Douzième section

  \normalsize Détermination du méridien du lieu et de l'azimut de la \textit{qibla}
\end{center}

\noindent On détermine la \emph{ligne méridienne} en prenant le Soleil
à deux hauteurs égales de part et d'autre de sa hauteur maximale~; on
trace les deux directions de son ombre portées par un même gnomon sur
un terrain plan~; puis on bissecte l'angle formé entre elles. Pour ce
faire, on prend le sommet de l'angle comme centre, on trace un arc de
cercle qui coupe les deux ombres, et on mène une droite entre le
milieu [des extrémités] de cet arc et le centre. Cette droite est dans
le plan du cercle méridien et on l'appelle ligne méridienne. La ligne
passant par le centre et perpendiculaire à cette ligne est dans le
plan du cercle origine des azimuts, et on l'appelle ligne de l'Est et
de l'Ouest.

Voici une autre méthode pour déterminer [la ligne méridienne]. On
installe un gnomon sur un terrain extrêmement plan~; on y trace un
cercle de rayon deux fois plus long que le gnomon~; on observe
l'entrée de l'ombre du gnomon dans le cercle avant midi, puis sa
sortie [du cercle] après [midi]~; une fois connues l'entrée et la
sortie, on mène une droite entre le milieu de l'arc joignant les
extrémités des deux ombres et le centre~; c'est la ligne
méridienne. La ligne qui lui est perpendiculaire et qui passe par le
centre du cercle est la ligne de l'Est et de l'Ouest. Ces deux lignes
font quatre quarts de cercle, et on découpe chaque quart en $90$
subdivisions égales afin d'y lire les azimuts là où l'ombre tombe sur
la circonférence~: [le nombre] de ces subdivisions entre le point de
l'Est ou de l'Ouest et l'ombre est un azimut. Ce cercle s'appelle
cercle \emph{indien}.

L'azimut de la \textit{qibla} est l'intersection de l'horizon du lieu
et du cercle azimutal passant par le zénith du lieu et par le zénith
de la Mecque. La ligne menée entre ce point [d'intersection] et le
centre de l'horizon est la ligne d'azimut de la \textit{qibla}~; c'est
la \uwave{tangente} à l'arc [d'un grand cercle] sur lequel est la base
du \textit{mi\d{h}r\=ab}. Si l'orant en prière met cette
\uwave{tangente} entre ses deux jambes en s'asseyant dessus, alors il
peut prier du lieu où il s'agenouille sur la circonférence d'un cercle
à la surface de la Terre. Ce qui est entre ses deux pieds et le centre
de la Maison\footnote{la maison sacrée, c'est-à-dire la
  \textit{Ka`ba}.} va à la rencontre de la droite menée entre la
Maison et le point qui culmine au dessus d'elle au zénith de la
Mecque. [La \uwave{tangente}] peut certes aller à la rencontre de la
Maison, mais puisque l'horizon de la Mecque est au dessous de
l'horizon de l'orant, la vision [de l'orant] ne peut pas viser la
Maison.

\newpage\phantomsection
\begin{center}
  \small
\begingroup%
  \makeatletter%
  \providecommand\color[2][]{%
    \errmessage{(Inkscape) Color is used for the text in Inkscape, but the package 'color.sty' is not loaded}%
    \renewcommand\color[2][]{}%
  }%
  \providecommand\transparent[1]{%
    \errmessage{(Inkscape) Transparency is used (non-zero) for the text in Inkscape, but the package 'transparent.sty' is not loaded}%
    \renewcommand\transparent[1]{}%
  }%
  \providecommand\rotatebox[2]{#2}%
  \ifx\svgwidth\undefined%
    \setlength{\unitlength}{226.275bp}%
    \ifx\svgscale\undefined%
      \relax%
    \else%
      \setlength{\unitlength}{\unitlength * \real{\svgscale}}%
    \fi%
  \else%
    \setlength{\unitlength}{\svgwidth}%
  \fi%
  \global\let\svgwidth\undefined%
  \global\let\svgscale\undefined%
  \makeatother%
  \begin{picture}(1,2.47471444)%
    \put(0,0){\includegraphics[width=\unitlength]{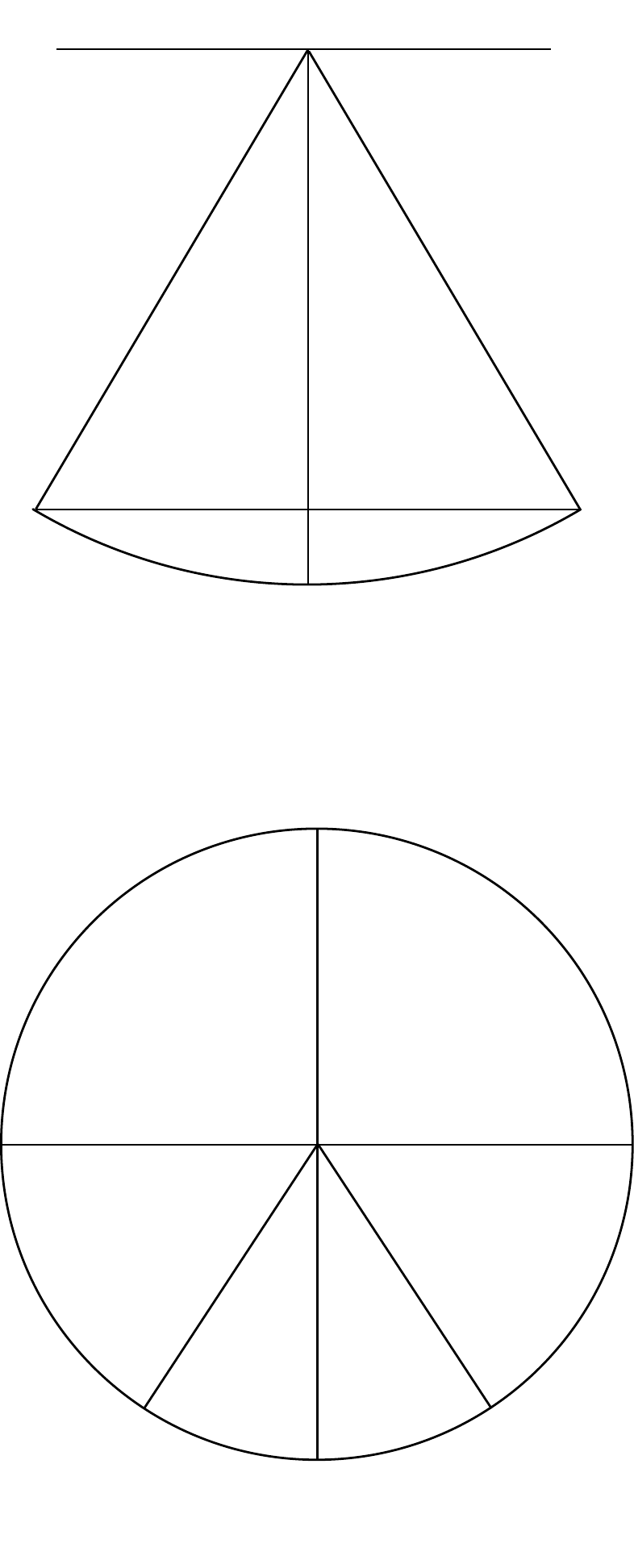}}%
    \put(0.84715185,2.45016605){\color[rgb]{0,0,0}\makebox(0,0)[rb]{\smash{\RL{al-`mUd alla_dI fI smt dA'iraT 'awl al-smUt}}}}%
    \put(0.52892292,2.32064453){\color[rgb]{0,0,0}\rotatebox{90}{\makebox(0,0)[rb]{\smash{\RL{al-zAwyaT}}}}}%
    \put(0.64943289,2.22650052){\color[rgb]{0,0,0}\rotatebox{122.04528609}{\makebox(0,0)[rb]{\smash{\RL{al-.zIl al-'awl ^srqy}}}}}%
    \put(0.33574667,2.22085954){\color[rgb]{0,0,0}\rotatebox{59.99999989}{\makebox(0,0)[rb]{\smash{\RL{al-.zIl al-_tAnI .grby}}}}}%
    \put(0.53130029,1.69592079){\color[rgb]{0,0,0}\rotatebox{90}{\makebox(0,0)[lb]{\smash{\RL{_h.t.t n.sf al-nhAr}}}}}%
    \put(0.6501728,0.11672223){\color[rgb]{0,0,0}\makebox(0,0)[rb]{\smash{\RL{nq.taT al-^smAl}}}}%
    \put(0.38753413,0.69503238){\color[rgb]{0,0,0}\makebox(0,0)[rb]{\smash{\RL{wa-al-m^srq}}}}%
    \put(0.62996979,1.19253093){\color[rgb]{0,0,0}\makebox(0,0)[rb]{\smash{\RL{nq.taT al-jnUb}}}}%
    \put(0.84576335,0.69414979){\color[rgb]{0,0,0}\makebox(0,0)[rb]{\smash{\RL{_h.t.t al-m.grb}}}}%
    \put(0.55004078,1.13115973){\color[rgb]{0,0,0}\rotatebox{90}{\makebox(0,0)[rb]{\smash{\RL{_h.t.t n.sf al-nhAr}}}}}%
    \put(0.21951173,0.23622022){\color[rgb]{0,0,0}\rotatebox{50.20901334}{\makebox(0,0)[rb]{\smash{\RL{m_hrj al-.zIl}}}}}%
    \put(0.78524017,0.22423414){\color[rgb]{0,0,0}\rotatebox{-55.3332966}{\makebox(0,0)[lb]{\smash{\RL{md_hl al-.zIl}}}}}%
    \put(0.49568873,0.63662594){\color[rgb]{0,0,0}\rotatebox{90}{\makebox(0,0)[rb]{\smash{\RL{mw.d` al-mqyAs}}}}}%
  \end{picture}%
\endgroup%
\normalsize
\end{center}

\newpage\phantomsection
\begin{center}
  \small
\begingroup%
  \makeatletter%
  \providecommand\color[2][]{%
    \errmessage{(Inkscape) Color is used for the text in Inkscape, but the package 'color.sty' is not loaded}%
    \renewcommand\color[2][]{}%
  }%
  \providecommand\transparent[1]{%
    \errmessage{(Inkscape) Transparency is used (non-zero) for the text in Inkscape, but the package 'transparent.sty' is not loaded}%
    \renewcommand\transparent[1]{}%
  }%
  \providecommand\rotatebox[2]{#2}%
  \ifx\svgwidth\undefined%
    \setlength{\unitlength}{226.275bp}%
    \ifx\svgscale\undefined%
      \relax%
    \else%
      \setlength{\unitlength}{\unitlength * \real{\svgscale}}%
    \fi%
  \else%
    \setlength{\unitlength}{\svgwidth}%
  \fi%
  \global\let\svgwidth\undefined%
  \global\let\svgscale\undefined%
  \makeatother%
  \begin{picture}(1,2.44549322)%
    \put(0,0){\includegraphics[width=\unitlength]{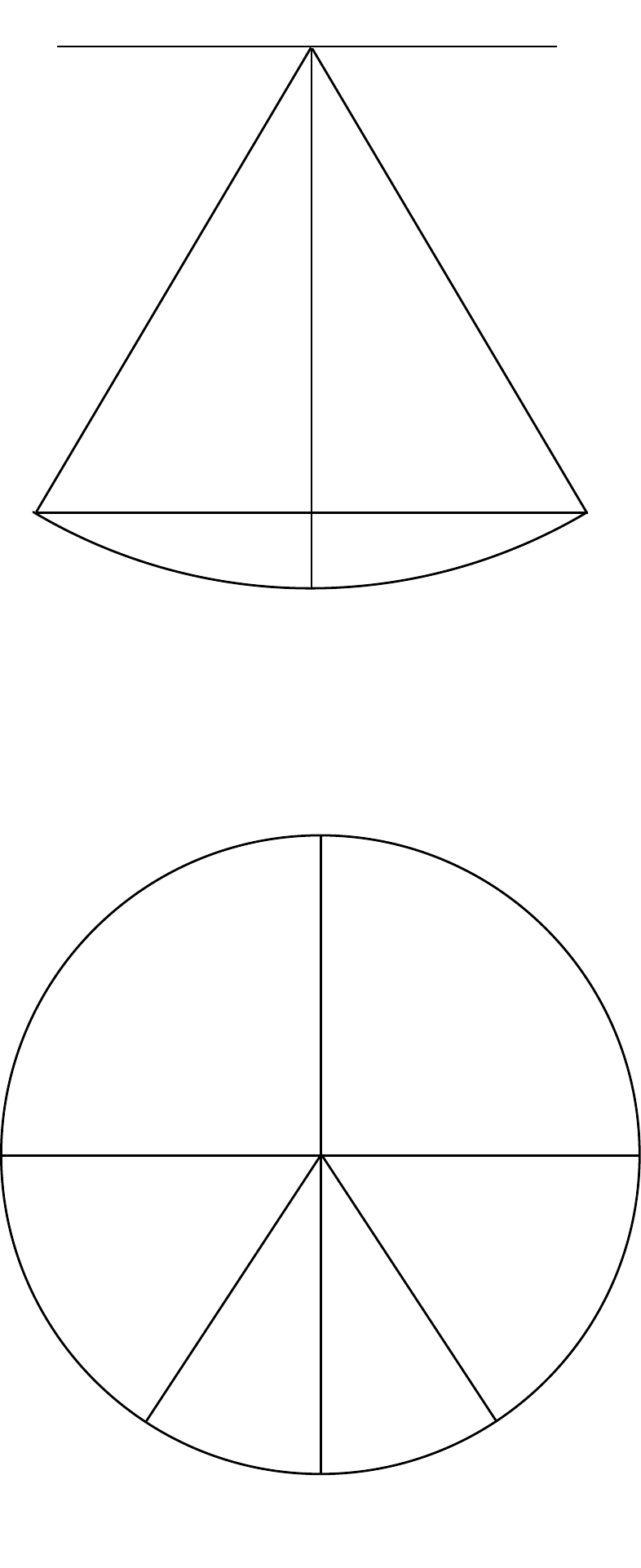}}%
    \put(0.99109799,2.42094483){\color[rgb]{0,0,0}\makebox(0,0)[rb]{\smash{la perpendiculaire, dans la direction du cercle origine des azimuts}}}%
    \put(0.52892292,2.29647406){\color[rgb]{0,0,0}\rotatebox{90}{\makebox(0,0)[rb]{\smash{l'angle}}}}%
    \put(0.65170505,2.2035667){\color[rgb]{0,0,0}\rotatebox{121.30954677}{\makebox(0,0)[rb]{\smash{la première ombre est à l'Est}}}}%
    \put(0.39096595,2.27765452){\color[rgb]{0,0,0}\rotatebox{57.78398425}{\makebox(0,0)[rb]{\smash{la deuxième ombre est à l'Ouest}}}}%
    \put(0.53130029,1.67175033){\color[rgb]{0,0,0}\rotatebox{90}{\makebox(0,0)[lb]{\smash{la ligne méridienne}}}}%
    \put(0.6501728,0.09255177){\color[rgb]{0,0,0}\makebox(0,0)[rb]{\smash{point Nord}}}%
    \put(0.38753413,0.67086191){\color[rgb]{0,0,0}\makebox(0,0)[rb]{\smash{la ligne de l'Est}}}%
    \put(0.62996979,1.16836046){\color[rgb]{0,0,0}\makebox(0,0)[rb]{\smash{point Sud}}}%
    \put(0.92657516,0.66745393){\color[rgb]{0,0,0}\makebox(0,0)[rb]{\smash{et de l'Ouest}}}%
    \put(0.55004078,1.10698927){\color[rgb]{0,0,0}\rotatebox{90}{\makebox(0,0)[rb]{\smash{la ligne méridienne}}}}%
    \put(0.21951173,0.21204975){\color[rgb]{0,0,0}\rotatebox{50.20901334}{\makebox(0,0)[rb]{\smash{sortie de l'ombre}}}}%
    \put(0.78524017,0.20006367){\color[rgb]{0,0,0}\rotatebox{-55.3332966}{\makebox(0,0)[lb]{\smash{entrée de l'ombre}}}}%
    \put(0.49568873,0.61245547){\color[rgb]{0,0,0}\rotatebox{90}{\makebox(0,0)[rb]{\smash{position du gnomon}}}}%
  \end{picture}%
\endgroup%
\normalsize
\end{center}

\newpage\phantomsection
\index{BBAHBD@\RL{qbl}!BBAHBDAI 2@\RL{qaws in.hirAf al-qiblaT}, arc de déviation de la \textit{qibla}}
\index{BEBCBCAI@\RL{makkaT}, La Mecque}
\index{AOBHAQ@\RL{dwr}!BEAOAGAQ@\RL{mdAr}, trajectoire, trajectoire diurne (\textit{i. e.} cercle parallèle à l'équateur)}
\index{ATAQBB@\RL{^srq}!BEATAQBB@\RL{ma^sraq}, Est}
\index{ATAQBB@\RL{^srq}!ANAWAW BEATAQBB BEAZAQAH@\RL{_ha.t.t al-ma^sraq wa-al-ma.grab}, ligne de l'Est et de l'Ouest}
\index{AZAQAH@\RL{.grb}!BEAZAQAH@\RL{ma.grab}, Ouest}
\label{var14}\label{var17}
\includepdf[pages=38,pagecommand={\thispagestyle{plain}}]{edit2.pdf}

\emph{L'arc de déviation} de la \textit{qibla} est l'arc de l'horizon
entre son point d'intersection avec le cercle azimutal et son point
d'intersection avec le méridien~; il mesure combien l'orant doit se
détourner de la direction d'un des quatre points [cardinaux] -- le
Nord, le Sud, l'Est ou l'Ouest -- pour faire face à la Maison.

Pour déterminer les azimuts [de la \textit{qibla}], il faut déterminer
la longitude de la Mecque, celle du lieu supposé, et leurs
latitudes. En prenant comme origine les \^Iles Canaries, la longitude
de la Mecque vaut soixante-dix-sept parts et un sixième~; en prenant
comme origine la côte de la mer, elle vaut soixante-sept parts et un
sixième~; sa latitude vaut vingt-et-une parts et deux tiers. La Mecque
est à l'Est de chaque pays dont la longitude est moindre que la
longitude de la Mecque, et elle est à l'Ouest de chaque pays dont la
longitude est supérieure à sa longitude. Si leur deux longitudes sont
égales, alors la Mecque est sur le méridien [du pays], vers le Sud si
sa latitude est inférieure à la latitude [du pays], et vers le Nord si
sa latitude est supérieure à la latitude [du pays]. Si leurs deux
latitudes sont égales, alors le pays et la Mecque sont ensemble sous
la même trajectoire diurne~; si la longitude [du pays] est moindre que
la longitude de la Mecque, alors la Mecque est à gauche du point où se
lève l'équinoxe pour ce pays~; et si la longitude [du pays] est
supérieure à sa longitude, alors la Mecque est à droite du point où se
couche l'équinoxe~; mais elle n'est pas au point Est dans le premier
cas, ni au point Ouest dans le second cas, ni sur la ligne de l'Est et
de l'Ouest comme on peut croire.

En effet, puisque leurs longitudes diffèrent, le cercle origine des
azimuts de l'un [des deux pays] est tangent à la trajectoire diurne
mentionnée en un point distinct du point où le cercle origine des
azimuts de l'autre [pays] lui est tangent~; donc [la trajectoire
  diurne et l'horizon] se coupent en des points distincts des points
Est et Ouest, et les lignes de l'Est et de l'Ouest [des deux pays] ne
sont les mêmes. Ce qu'on a cru est vrai seulement dans les contrées de
l'équateur terrestre, quelque soit la longitude~; car tous les zéniths
y sont sur l'équateur, et [l'équateur] est donc le cercle origine des
azimuts de toutes [ces contrées]. \label{loxoqibla}

\newpage\phantomsection
\includepdf[pages=39,pagecommand={\thispagestyle{plain}}]{edit2.pdf}
Parmi les moyens les plus simples pour déterminer l'azimut de la
\textit{qibla} en différentes longitudes ou même en différentes
longitudes et latitudes, en voici un. On observe le transit du Soleil
au zénith des gens de la Mecque -- c'est-à-dire quand il est là-bas à
midi dans le huitième ou le septième degré des Gémeaux ou bien dans le
vingt-troisième degré du Cancer. On prend l'écart entre les deux
longitudes. On pose une heure pour chaque [portion] de quinze parts
d'écart, et quatre minutes pour chaque part. Le total, ce sont des
heures à compter de midi. Puis on fait une observation, le même jour,
à cet instant avant midi si la Mecque est à l'Est, après midi si elle
est à l'Ouest~; alors l'azimut de l'ombre du gnomon à l'heure [où tu
  observes] est l'azimut de la \textit{qibla}.

On a besoin de retrancher l'arc de déviation aussi bien quand les deux
[pays] diffèrent seulement en longitude que quand ils diffèrent en
longitude et en latitude~; mais quand il ne diffèrent qu'en latitude,
alors leurs azimuts sont sur la ligne méridienne, et l'orant fait donc
face à point Sud si la latitude de la Mecque est inférieure, et il
fait face au point Nord si sa latitude est supérieure. Dieu est le
plus savant.

\begin{center}
  Avec cette section s'achève la deuxième partie du livre
  \textit{L'achèvement de l'enquête et la correction des fondements}.
  Dieu fait réussir ce qui est juste, il nous a donné suffisamment et
  il a combé de bienfaits celui qui s'en remet à lui.
\end{center}

\label{txt_fin}

\part*{Commentaire mathématique}
\addcontentsline{toc}{part}{Commentaire mathématique}
\label{comm_debut}

\paragraph{Théorème d'Apollonius, lemme de `Ur\d{d}{\=\i} et couple de
  \d{T}\=us{\=\i}} On démontre ici trois propositions concernant les rotations
affines du plan. Elles nous seront utiles pour comparer entre eux les modèles
planétaires.

\paragraph{Proposition 1} Soient $P_1$, $P_2$, $P_3$, $P_4$ quatre points
alignés et $\alpha\in\lbrack 0,2\pi\rbrack$. Si $\overrightarrow{P_1P_3}=\overrightarrow{P_2P_4}$, alors $R_{P_2,\alpha}=R_{P_1,\alpha}R_{P_3,-\alpha}R_{P_4,\alpha}$.
\label{apollonius}

\emph{Démonstration}. $R_{P_2,\alpha}R_{P_4,-\alpha}$ est une transformation
affine directe dont l'application linéaire associée est l'identité, et c'est
donc une translation. C'est la translation de vecteur~:
$$\overrightarrow{P_4R_{P_2,\alpha}(P_4)}=R_\alpha(\overrightarrow{P_2P_4})-\overrightarrow{P_2P_4},$$
où $R_\alpha$ désigne la rotation \emph{vectorielle} d'angle $\alpha$. De même $R_{P_1,\alpha}R_{P_3,-\alpha}$ est une translation de vecteur~:
$$\overrightarrow{P_3R_{P_1,\alpha}(P_3)}=R_\alpha(\overrightarrow{P_1P_3})-\overrightarrow{P_1P_3},$$
mais $\overrightarrow{P_1P_3}=\overrightarrow{P_2P_4}$, donc $R_{P_1,\alpha}R_{P_3,-\alpha}=R_{P_2,\alpha}R_{P_4,-\alpha}$, \textit{q. e. d.}

\paragraph{Remarque} Si $P_2$ est situé entre $P_1$ et $P_3$, une conséquence de cette proposition était bien connue de Ptolémée qui l'attribuait à Apollonius. Si $P_1$ est le centre du Monde, $P_2$ le centre d'un excentrique, et $P_4$ un astre à son apogée sur l'excentrique, on a~:
$$R_{P_2,\alpha}(P_4)=R_{P_1,\alpha}R_{P_3,-\alpha}(P_4).$$
Ceci montre que l'excentrique peut être remplacé par un modèle équivalent constitué d'un déférent de centre $P_1$ et d'un épicycle de centre $P_3$.

\paragraph{Proposition 2} Soient $P_1$, $P_2$, $P_3$, $P_4$, $P_5$ cinq points alignés, et $\alpha\in\lbrack 0,2\pi\rbrack$. Si $\overrightarrow{P_3P_4}=-\overrightarrow{P_1P_2}$ et $\overrightarrow{P_3P_5}=\overrightarrow{P_1P_2}$, alors
$$R_{P_1,\alpha}R_{P_3,\alpha}R_{P_4,-\alpha}=R_{P_2,\alpha}R_{P_5,-\alpha}R_{P_3,2\alpha}R_{P_4,-\alpha}.$$

\emph{Démonstration}. Comme $\overrightarrow{P_2P_5}=\overrightarrow{P_1P_3}$, d'après la proposition 1, on a~:
$$R_{P_1,\alpha}=R_{P_2,\alpha}R_{P_5,-\alpha}R_{P_3,\alpha},$$
d'où le résultat escompté si l'on compose à droite avec $R_{P_3,\alpha}R_{P_4,-\alpha}$.

\paragraph{Remarque} $R_{P_5,-\alpha}R_{P_3,2\alpha}R_{P_4,-\alpha}$ est une translation parallèlement à la droite $(P_1P_2)$. En effet~:\label{couple_transl}
$$(R_{P_5,-\alpha}R_{P_3,\alpha})(R_{P_3,\alpha}R_{P_4,-\alpha})=t_{\overrightarrow{P_3R_{P_5,-\alpha}(P_3)}}\circ t_{\overrightarrow{P_4R_{P_3,\alpha(P_4)}}},$$
mais
$$\overrightarrow{P_3R_{P_5,-\alpha(P_3)}}+\overrightarrow{P_4R_{P_3,\alpha(P_4)}}=R_{-\alpha}(\overrightarrow{P_3P_4})+R_\alpha(\overrightarrow{P_3P_4})-2\overrightarrow{P_3P_4},$$
or
$$(R_{-\alpha}(\overrightarrow{P_3P_4})+R_\alpha(\overrightarrow{P_3P_4})-2\overrightarrow{P_3P_4})\wedge\overrightarrow{P_3P_4}=\Vert\overrightarrow{P_3P_4}\Vert^2(\sin\alpha+\sin(-\alpha))=0.$$
Cette composée de rotations apparaît dans les modèles planétaires utilisant un ``couple de \d{T}\=us{\=\i}''.

\paragraph{Proposition 3} Soit $P_1$, $P_2$, $P_3$ trois points alignés avec $\overrightarrow{P_1P_3}=-\overrightarrow{P_1P_2}$, alors $R_{P_3,-\alpha}R_{P_1,2\alpha}R_{P_2,-\alpha}=R_{P_3,\alpha}R_{P_1,-2\alpha}R_{P_2,\alpha}$. En particulier, le sens de rotation des orbes d'un couple de \d{T}\=us{\=\i} n'importe pas.

\emph{Démonstration}.
$$R_{P_3,2\alpha}R_{P_1,-2\alpha}=t_{\overrightarrow{P_1R_{P_3,2\alpha}(P_1)}}=t_{\overrightarrow{P_2R_{P_1,2\alpha}(P_2)}}=R_{P_1,2\alpha}R_{P_2,-2\alpha},$$
  d'où le résultat.

\paragraph{Préliminaires} \label{preliminaires}
Nous n'analyserons pas en détail ce qu'est un orbe pour Ibn
al-\v{S}\=a\d{t}ir. Qu'il suffise de le concevoir comme un corps à
symétrie sphérique. Chaque orbe est doté d'un centre. Le centre de
certains orbes est aussi centre du Monde. Enfin, et surtout, Ibn
al-\v{S}\=a\d{t}ir utilise chaque orbe comme un référentiel solide en
mouvement. Pour qu'un corps solide à symétrie sphérique puisse faire
office de référentiel spatial, il faut distinguer un plan attaché à
l'orbe et passant par son centre, et une direction dans ce plan. Dans
tout ce qui suit, on représentera plus commodément chaque orbe par une
sphère, un cercle et un point~: à savoir, une sphère dont le centre
est le centre de l'orbe, un grand cercle section de cette sphère par
son plan, et un point de ce cercle situé dans la direction distinguée
par rapport au centre.

Chaque orbe est mobile par rapport à l'orbe qui le porte (c'est-à-dire
au sein du référentiel constitué par l'orbe qui le porte). Un seul
mouvement est admis pour chaque orbe~: un mouvement de rotation
uniforme autour d'un axe passant par son centre. Autrement dit, le
centre de chaque orbe est immobile au sein de l'orbe qui le porte. Il
n'y a aucune restriction quant à l'inclinaison de l'axe par rapport au
plan de l'orbe~: tel orbe a son axe perpendiculaire à son plan, tel
autre a son axe perpendiculaire au plan de l'orbe qui le porte. Mais
tous ces mouvements sont relatifs, et chaque orbe hérite aussi <<~par
accident~>> (au sens aristotélicien du terme) du mouvement de l'orbe
qui le porte, qui hérite lui-même du mouvement de l'orbe qui le porte,
et ainsi de suite. Le mouvement de chaque orbe est donc composé d'un
mouvement propre de rotation uniforme par rapport à l'orbe qui le
porte, et du mouvement absolu de l'orbe qui le porte. Seul le neuvième
orbe, porté par aucun autre, n'est animé que de son mouvement propre
de rotation uniforme.

Ce cadre conceptuel donne lieu à trois combinaisons possibles, toutes
présentes dans le texte d'Ibn al-\v{S}\=a\d{t}ir. Nous les décrirons
au moyen des trois figures (i), (ii) et (iii)
page~\pageref{fig2}. Dans chacune de ces trois figures il y a deux
orbes representés chacun par une sphère, un grand cercle de cette
sphère et un point sur ce cercle. Le premier orbe dont le point est
$P$ porte le second dont le point est $Q$. Dans les figures (i) et
(ii), le point $P$ de l'orbe portant est choisi par commodité au
centre de l'orbe porté. Au contraire, dans la figure (iii), les deux
orbes ayant même centre $O$, il a fallu choisir un autre point~: le
point $P$ est alors un point à l'intersection des plans des deux
orbes. Dans tout ce paragraphe, on désignera chaque orbe par son
point. L'orbe de $P$ porte donc l'orbe de $Q$. Ainsi, le lieu occupé
par l'orbe de $Q$ est constant au sein du référentiel solide constitué
par l'orbe de $P$. Mais l'orbe de $Q$ est animé d'un mouvement de
rotation uniforme sur lui-même, autour d'un axe représenté par un
vecteur sur chaque figure. Peu importe ici le sens de cette rotation~;
seule importe la direction de son axe, perpendiculaire au plan de
l'orbe de $P$ dans la figure (i), mais perpendiculaire au plan de
l'orbe de $Q$ dans les figures (ii) et (iii). La figure (i) présente
cette particularité que la trajectoire du point $Q$ au sein du
référentiel constitué par l'orbe de $P$ n'est pas dans le plan de
l'orbe de $Q$~: le point $Q$ décrit un petit cercle parallèle au plan
de l'orbe de $P$. Enfin, dans chacune des trois figures, l'orbe de $P$
est lui-même animé d'un mouvement de rotation uniforme autour d'un axe
représenté par un vecteur d'origine $O$. Peu importe que ce vecteur
soit perpendiculaire au plan de l'orbe de $P$ -- comme il l'est dans
ces figures -- ou non. Ce dernier mouvement entraîne le point $P$ et
aussi l'orbe de $Q$.

Chacune des planètes supérieures (Mars, Jupiter, Saturne) est portée
par un système d'orbes: un <<~orbe parécliptique~>>, un <<~orbe
incliné~>>, un <<~orbe déférent~>>, un <<~orbe rotateur~>> et un
<<~orbe de l'épicycle~>>. La figure (i) représente alors la relation
entre l'orbe incliné (orbe de $P$) et l'orbe déférent (orbe de $Q$),
ainsi que la relation entre l'orbe déférent (orbe de $P$) et l'orbe
rotateur (orbe de $Q$). Pour Vénus, les orbes portent les mêmes noms,
mais la figure (ii) représente alors la relation entre l'orbe déférent
(orbe de $P$) et l'orbe rotateur (orbe de $Q$), ainsi que la relation
entre l'orbe rotateur (orbe de $P$) et l'orbe de l'épicycle (orbe de
$Q$). Mercure compte davantage d'orbes sans toutefois introduire de
combinaison nouvelle.

Pour chaque planète, la figure (iii) représente la relation entre
l'orbe parécliptique (orbe de $P$) et l'orbe incliné (orbe de $Q$)~:
le point $P$ est alors l'un des n{\oe}uds. Mais cette figure
représente aussi la relation entre le neuvième orbe (orbe de $P$) et
le huitième, dit de l'écliptique (orbe de $Q$)~: le point $P$ est
alors le point vernal.

\begin{figure}
\begin{center}
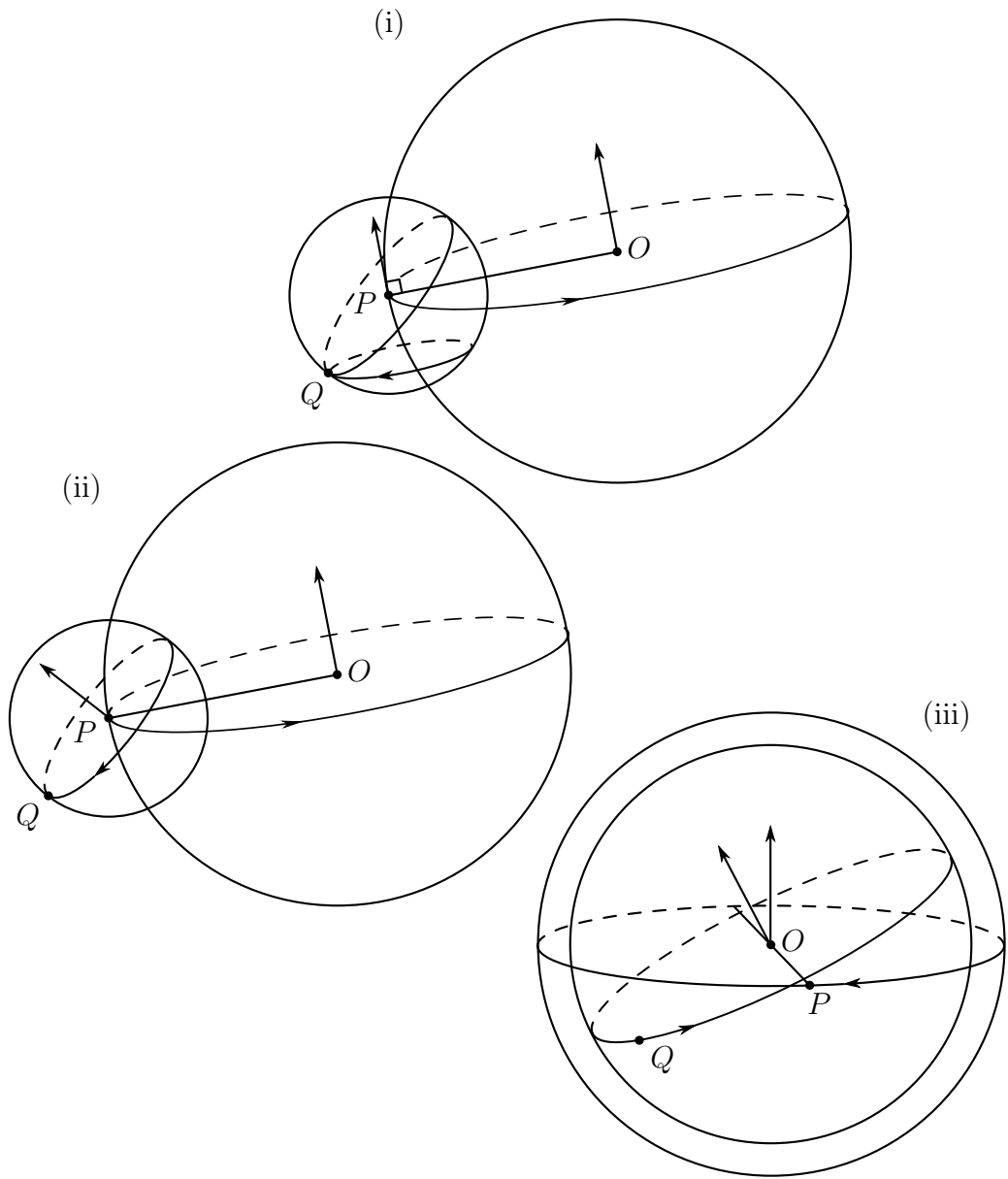
\caption{\label{fig2}Un orbe est un référentiel solide en mouvement}
\end{center}
\end{figure}

\paragraph{Précession des équinoxes}
Observons à nouveau la figure (iii). Chez Ptolémée, les étoiles fixes
sont attachées au dernier orbe et le mouvement relatif de cet orbe par
rapport aux orbes inférieurs rend compte du mouvement de précession
des équinoxes. Pour les astronomes arabes, depuis le \emph{Traité sur
  l'année solaire} certainement dû aux Banu Musa (IXème
siècle)\footnote{\textit{cf.} \cite{qurra1987}
  p.~\textit{xlvi}-\textit{lxxv} et 27-67.}, il existe un neuvième
orbe au delà des étoiles fixes situées sur le huitième orbe, et c'est
le mouvement relatif du huitième orbe par rapport au neuvième qui rend
alors compte de la précession. Tous les orbes inférieurs des planètes,
du Soleil et de la Lune héritent de ce mouvement comme les étoiles
fixes (tandis que chez Ptolémée, le Soleil en était
exempt\footnote{\textit{cf.} \cite{pedersen1974} p.~147.}).

Chez Ibn al-\v{S}\=a\d{t}ir, l'orbe parécliptique de chaque planète
est un orbe dont le centre coïncide avec le centre du Monde, comme
l'écliptique. L'orbe parécliptique reproduit, en plus petit, l'orbe de
l'écliptique et y est attaché en son centre. Le parécliptique est
animé d'un petit mouvement de rotation sur lui-même, par rapport au
référentiel constitué par l'orbe de l'écliptique~; mais comme l'axe de
rotation du parécliptique coïncide avec l'axe de l'écliptique, et que
la composée de deux rotations de même axe en est encore une, le
parécliptique est aussi en rotation uniforme autour de son axe au sein
du référentiel constitué par le neuvième orbe. Le mouvement de l'orbe
de l'écliptique par rapport au neuvième est d'un degré toutes les
soixante-dix années persanes. Le mouvement de l'orbe parécliptique par
rapport au neuvième orbe est un peu plus rapide, comme l'avait vérifié
al-Zarqalluh (XIème siècle)~: il est d'un degré toutes les soixante
années persanes, c'est-à-dire $0°1'$ par année
persane\footnote{Al-Zarqalluh (Azarquiel) avait observé que l'apogée
  solaire se meut d'un degré tous les 279 ans par rapport aux étoiles
  fixes, elles-mêmes étant animées du mouvement de précession
  (\textit{cf.} \cite{rashed1997} p.~285-287). Or
  $\dfrac{1}{70}+\dfrac{1}{279}\simeq\dfrac{1}{60}$.}.

\paragraph{Unité de temps}
Le temps absolu est mesuré par le mouvement du neuvième orbe, seul
orbe dont le mouvement propre est absolu (les mouvements propres des
autres orbes sont relatifs à l'orbe qui les porte).

En un lieu donné, la durée entre deux culminations du Soleil (de midi
à midi du lendemain) est approximativement constante.

Considérons deux unités possibles pour la mesure du temps~:
\begin{itemize}
\item le \emph{jour sidéral}~: durée d'une rotation complète du ciel
  visible autour du pôle nord (temps de retour d'une étoile fixe en
  son lieu par rapport au pôle nord et à l'horizon)
\item le \emph{jour solaire}~: durée entre deux culminations du Soleil
  (de midi à midi du lendemain)
\end{itemize}
Le neuvième orbe étant supposé mû d'un mouvement uniforme, le jour
sidéral est de durée constante (sauf à tenir compte du mouvement de
précession, le mouvement propre du huitième orbe qui est infime). En
revanche, la durée du jour solaire varie au fil des saisons, ce n'est
donc pas une bonne unité. En effet, soit midi au Soleil. Un jour
sidéral plus tard, le Soleil aura presque effectué une rotation
complète d'est en ouest à cause du mouvement du neuvième orbe~; pas
tout à fait cependant, à cause du mouvement du parécliptique du Soleil
qui l'aura légèrement entraîné vers l'est, d'un peu moins qu'un
degré. Il ne sera donc pas encore midi au Soleil~: le jour sidéral est
plus court que le jour solaire. Cet effet se fera sentir d'autant plus
que l'angle entre le parécliptique et la trajectoire diurne du Soleil
sera faible (ainsi la différence sera plus importante lors des
solstices que lors des équinoxes). De toute façon, le mouvement du
Soleil le long de son parécliptique n'est pas uniforme.

Tout calendrier civil adoptant pourtant le jour solaire comme unité,
il est commode de choisir comme unité astronomique un <<~jour solaire
moyen~>>~: durée entre deux culminations, non plus du Soleil, mais
d'un point imaginaire mû d'un mouvement uniforme le long de l'équateur
et faisant le tour du ciel exactement en même temps que le
Soleil\footnote{c'est-à-dire en une année tropique}. C'est l'unité
adoptée par Ibn al-\v{S}\=a\d{t}ir~: \textit{al-yawm
  bi-laylatihi}. Quand une date est donnée en temps civil
(c'est-à-dire un certain nombre d'heures à compter de midi au Soleil),
il faut la corriger en lui ajoutant une certaine <<~équation du
temps~>>, \textit{ta`d{\=\i}l al-'ay\=am}, qui s'annule lors de
l'équinoxe de printemps quand le Soleil et le point imaginaire
coïncident sur l'équateur\footnote{En fait, comme on va le voir, le
  point imaginaire suit le Soleil \emph{moyen} et il est donc situé un
  peu à l'Ouest du Soleil au moment de l'équinoxe.}. L'heure est un
vingt-quatrième de jour solaire moyen.

Soit $t_0$ la date d'un équinoxe de printemps, et $t$ un instant donné
quelconque, dates exprimées en jours solaires moyens. Les mêmes dates
exprimées en temps civil (jours solaires, et heures comptées à partir
de midi au soleil) seront notées $\tau_0$ et $\tau$, et en jours
sidéraux $T_0$ et $T$.

Le point imaginaire se déplace le long de l'équateur à la vitesse du
Soleil moyen~; son ascension droite varie donc comme la longitude
$\overline{\lambda}$ du Soleil moyen. \`A l'instant $t$, elle vaut
$\overline{\lambda}(t)-\overline{\lambda}(t_0)$.

De combien le Soleil s'est-il déplacé vers l'Est depuis le dernier
équinoxe de printemps~? C'est précisément son ascension droite
$\alpha(t)$.

Le jour sidéral est plus court que le jour solaire, d'autant qu'il
faut au mouvement diurne pour récupérer le déplacement du Soleil vers
l'Est. La vitesse du mouvement diurne est d'environ 15° par heure
(elle est de 360° par jour sidéral, mais le jour sidéral compte un peu
moins que 24 heures). On a donc~:
$$(T-T_0)-(\tau-\tau_0)\simeq\frac{\alpha(t)}{15°}\text{ heures}$$ 
Il ne s'agit que d'une approximation, puisque pendant cette durée, le
Soleil se déplace légèrement davantage vers l'Est...

De même, le jour sidéral est plus court que le jour solaire moyen,
d'autant qu'il faut au mouvement diurne pour récupérer le déplacement
du point imaginaire vers l'Est le long de l'équateur~:
$$(T-T_0)-(t-t_0)\simeq\frac{\overline{\lambda}(t)-\overline{\lambda}(t_0)}{15°}\text{ heures}$$

Ibn al-\v{S}\=a\d{t}ir donne $\overline{\lambda}(t_0)=-2°1'7''$. L'équation du
temps est donc\footnote{On adopte ici la convention française quant au signe de l'équation du temps $E$.}~:
$$E=(t-t_0)-(\tau-\tau_0)
\simeq\frac{-\overline{\lambda}(t)-2°1'7''+\alpha(t)}{15°}\text{ heures}$$

Chez Ibn al-\v{S}\=a\d{t}ir, l'unité de temps est le jour solaire
moyen, ou l'année persane de 365 jours solaires moyens. L'<<~époque~>>
(c'est-à-dire l'instant $t=0$) est à midi du 24 décembre 1331. Il
s'agit probablement d'une date en temps solaire moyen (les valeurs des
paramètres à cette date ayant été obtenus par le calcul à partir d'une
observation antérieure\footnote{Une remarque de Kennedy et Roberts
  (\cite{roberts1959} p.~232) au sujet du \textit{z{\=\i}j} d'Ibn
  al-\v{S}\=a\d{t}ir le confirme.}, un jour d'équinoxe de
printemps). Le temps GMT (Greenwich Mean Time) est un temps solaire
moyen mesurant comme précédemment le déplacement uniforme d'un point
imaginaire le long de l'équateur, mais ce point ne coïncide pas avec
le Soleil lors de l'équinoxe de printemps~: l'équation du temps n'a
donc pas la même origine que chez Ibn al-\v{S}\=a\d{t}ir, elle est
d'environ 8 min à l'équinoxe de printemps et elle s'annule vers le 15
avril. Il faut aussi tenir compte du décalage horaire à la longitude
de Damas, 36°18'23''. L'époque d'Ibn al-\v{S}\=a\d{t}ir est donc, en
temps GMT, le 24 décembre 1331 à 12 h 8 min
$-\dfrac{36°18'23''}{15°}\simeq$ 9 h 43 min.

Pourquoi {\shatir} choisit-il le 24 décembre 1331 comme époque ? Le premier jour de la première année de l'ère de Yazdgard est le 12 juin 632 (date julienne). Sept cents années persanes plus tard (soit $700\times 365=255500$ jours), on est le premier jour de l'année 701 de l'ère de Yazdgard, un mardi (\textit{yawm al-thal\=ath\=a'}), le 23 \textit{rab{\=\i}` al-'awwal} 732 de l'hégire (date julienne~: 24 décembre 1331). {\shatir } avait alors 24 ans. 

\paragraph{Le Soleil : figure initiale}
Soit $(\mathbf{i},\mathbf{j},\mathbf{k})$ la base canonique de $\mathbb{R}^3$ qu'on identifiera à l'espace physique. Selon {\shatir}, la trajectoire du Soleil est contenue dans le plan de l'écliptique $(\mathbf{i},\mathbf{j})$. Dans les modèles d'{\shatir}, on chercherait en vain une description mathématique du mouvement en termes de transformations ou de différences entre deux instants. L'objet de la description mathématique est ici l'opération transformant une figure initiale (qui n'a d'existence qu'imaginaire) en une autre figure représentant la configuration des astres à un instant donné. Mais pour décrire le mouvement des astres, {\shatir} n'admet que des mouvements circulaires uniformes ou bien des mouvements composés de mouvements circulaires uniformes~: les seules opérations autorisées seront donc des rotations spatiales dont les paramètres dépendent linéairement de la variable temps.

\begin{figure}
  \begin{center}
    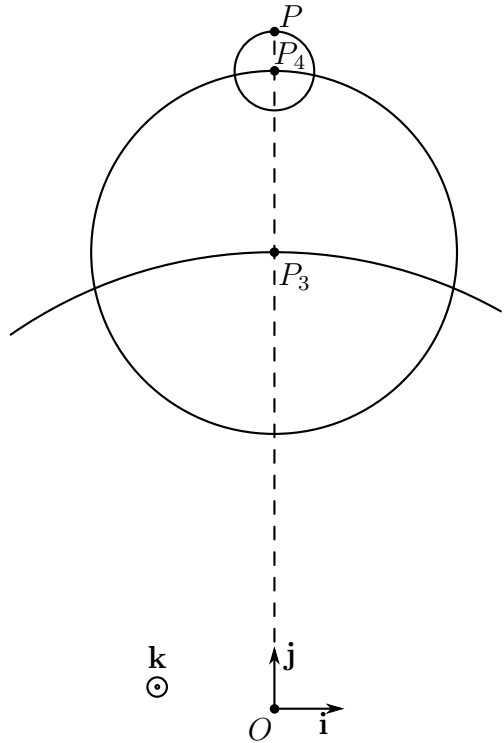
  \end{center}
  \caption{\label{fig003}Les orbes du Soleil : figure initiale}
\end{figure}

Dans la figure initiale, les centres des orbes du Soleil sont tous alignés dans la direction du vecteur $\mathbf{j}$ elle-même confondue avec la direction du point vernal. Le point $O$ est le centre du Monde, $P_1$ est le centre de l'orbe total (\textit{al-falak al-\v{s}\=amil}\footnote{La dénomination des orbes du Soleil peut prêter à confusion. Pour décrire le mouvement en longitude des planètes, chacune possède deux orbes concentriques, parécliptique et orbe incliné. Le Soleil n'a pas d'orbe incliné ; son ``orbe total'' joue à peu près le rôle du parécliptique des autres planètes, et son ``parécliptique'', celui de l'orbe incliné. Mais pour le Soleil, ces deux orbes ayant même centre et même axe de rotation, on pourrait aussi bien les confondre en un orbe unique ; {\shatir} le savait, et il laisse planer une certaine ambigu\"ité, par exemple dans la figure p. \pageref{soleil_orbes_solides} où le parécliptique semble contenir l'orbe total, tandis que son texte semble affirmer le contraire.}), et $P_2$ le centre du parécliptique du Soleil. Ces trois points sont confondus $O=P_1=P_2$. Le point $P_3$ est le centre de l'orbe déférent, $P_4$ le centre de l'orbe rotateur, et $P$ est le centre du Soleil, cet astre étant lui-même un corps sphérique. La figure \ref{fig003} précise la position de ces points
\footnote{Attention, notre figure \ref{fig003} n'est pas présente dans le traité d'{\shatir}. Il a certes représenté les centres des orbes du Soleil dans cinq positions distinctes, page \pageref{soleil_trajectoires} ; mais le Soleil y atteint son apogée au premier point du Cancer (c'était à peu près le cas à l'époque d'{\shatir}), et non au point vernal. Dans notre commentaire, nous préférons suivre l'usage adopté par {\shatir} pour d'autres planètes, car cet usage colle davantage avec l'architecture conceptuelle.}. Elle n'est pas à l'échelle. Nous donnons une représentation à l'échelle et en perspective p.~\pageref{soleil_blender}. Posons $\overrightarrow{OP_3}=60\ \mathbf{j}$, alors
$$\begin{array}{l}
\overrightarrow{P_3P_4}=4;37\ \mathbf{j}\\
\overrightarrow{P_4P}=2;30\ \mathbf{j}
\end{array}$$
où le point-virgule sépare partie entière et partie
fractionnaire\footnote{Toutes les valeurs numériques seront données en
  sexagésimal, et la virgule sert alors à séparer les rangs. Par
  exemple, $359;45,40=359+\dfrac{45}{60}+\dfrac{40}{60^2}$.}. 

Insistons encore sur ce point déjà soulevé~: \textit{a priori} la figure initiale ne représente la position des orbes à aucun instant de l'histoire passée ou future du Monde. La situation représentée peut se produire dans la réalité, s'il arrive que le Soleil soit à l'apogée de sa trajectoire à la traversée du point vernal, c'est-à-dire à l'équinoxe de printemps ; mais {\shatir} ne se pose même pas la question de savoir si le Soleil a jamais été à l'apogée de sa trajectoire à l'instant précis de l'équinoxe de printemps. Peu importe. Comme nous l'avons dit, il faut appliquer au point $P$ une suite de transformations géométriques pour obtenir la position du Soleil prédite ou observée à un instant donné. 

\paragraph{Le Soleil : transformations géométriques}
Les transformations appliquées au Soleil $P$ sont des rotations dans le plan $(\mathbf{i},\mathbf{j})$ paramétrées par deux angles $\lambda_A$, $\overline{\alpha}$ appelés l'\textit{apogée} et le \textit{centre}\footnote{Ici comme ailleurs, notre choix de notations est dicté par la commodité dans les comparaisons futures avec les modèles d'autres auteurs.}. Voici la liste de ces transformations~:
$$R(P_4,\mathbf{k},2\overline{\alpha}),\quad R(P_3,\mathbf{k},-\overline{\alpha}),
\quad R(P_2,\mathbf{k},\overline{\alpha}),\quad R(P_1,\mathbf{k},\lambda_A)$$
où $R(Q,\mathbf{k},\theta)$ désigne la rotation de centre $Q$ et d'angle $\theta$, le vecteur $\mathbf{k}$ désignant la direction de l'axe de rotation, et le sens de rotation. La description d'Ibn al-\v{S}\=a\d{t}ir ne laisse aucun
doute quant à l'ordre dans lequel appliquer ces quatres transformations.

Pour obtenir la configuration des orbes à un instant donné, il faut
donc appliquer toutes ces rotations aux points $P_3$, $P_4$, $P$.
L'image du point $P_3$ entraîné par les mouvements de l'orbe
total et de l'orbe parécliptique est~:
$$P_3'=R(P_1,\mathbf{k},\lambda_A)\circ R(P_2,\mathbf{k},\overline{\alpha})\ (P_3)$$
Quant au point $P_4$, il est aussi
entraîné par le mouvement de l'orbe déférent et devient~:
$$P_4'=R(P_1,\mathbf{k},\lambda_A)\circ R(P_2,\mathbf{k},\overline{\alpha})
\circ R(P_3,\mathbf{k},-\overline{\alpha})\ (P_4)$$
Enfin le point $P$, aussi entraîné par le mouvement de l'orbe rotateur,
devient~:
$$P'=R(P_1,\mathbf{k},\lambda_A)
\circ R(P_2,\mathbf{k},\overline{\alpha})
\circ R(P_3,\mathbf{k},-\overline{\alpha})
\circ R(P_4,\mathbf{k},2\overline{\alpha})\ (P)$$

\paragraph{Le Soleil : trajectoire}
Dans le cadre théorique que l'on vient de décrire, l'on ne peut guère
parler de trajectoire dans l'espace. Il y a bien toutefois une
trajectoire dans l'ensemble des valeurs des paramètres~:
$$\lambda_A=\dot{\lambda}_At+\lambda_A(0)$$
$$\overline{\alpha}=(\dot{\overline{\lambda}}-\dot{\lambda}_A)t+(\overline{\lambda}(0)-\lambda_A(0))$$
On a~:
$$\overline{\lambda}(0)=280;9,0,\quad\dot{\overline{\lambda}}=359;45,40
\text{ par année persane}$$
$$\lambda_{A}(0)=89;52,3,1,\quad \dot{\lambda}_A=0;1\text{ par année persane}$$

Si l'observation et la tradition délivrent une vitesse moyenne par année persane de $359;45,40$ degrés, une simple division permet d'en déduire la vitesse par jour solaire moyen, par mois de trente jours, ou bien par heure égale. C'est dans cet ordre qu'a dû procéder {\shatir}, comme le révèle la précision absurde des valeurs données~:
\begin{align*}
  \dot{\overline{\lambda}}&=0;59,8,19,43,33,41,55,4,6\text{ par jour solaire moyen,}\\
  &=29;34,9,51,46,50,57,32\text{ par mois de trente jours,}\\
  &=0;2,27,50,49,18,54,14,47\text{ par heure égale.}
\end{align*}
On en déduit aussi~:
$$\dot{\overline{\lambda}}-\dot{\lambda}_A=0;59,8,9,51,47\text{ par jour.}$$

\begin{figure}
  \begin{center}
    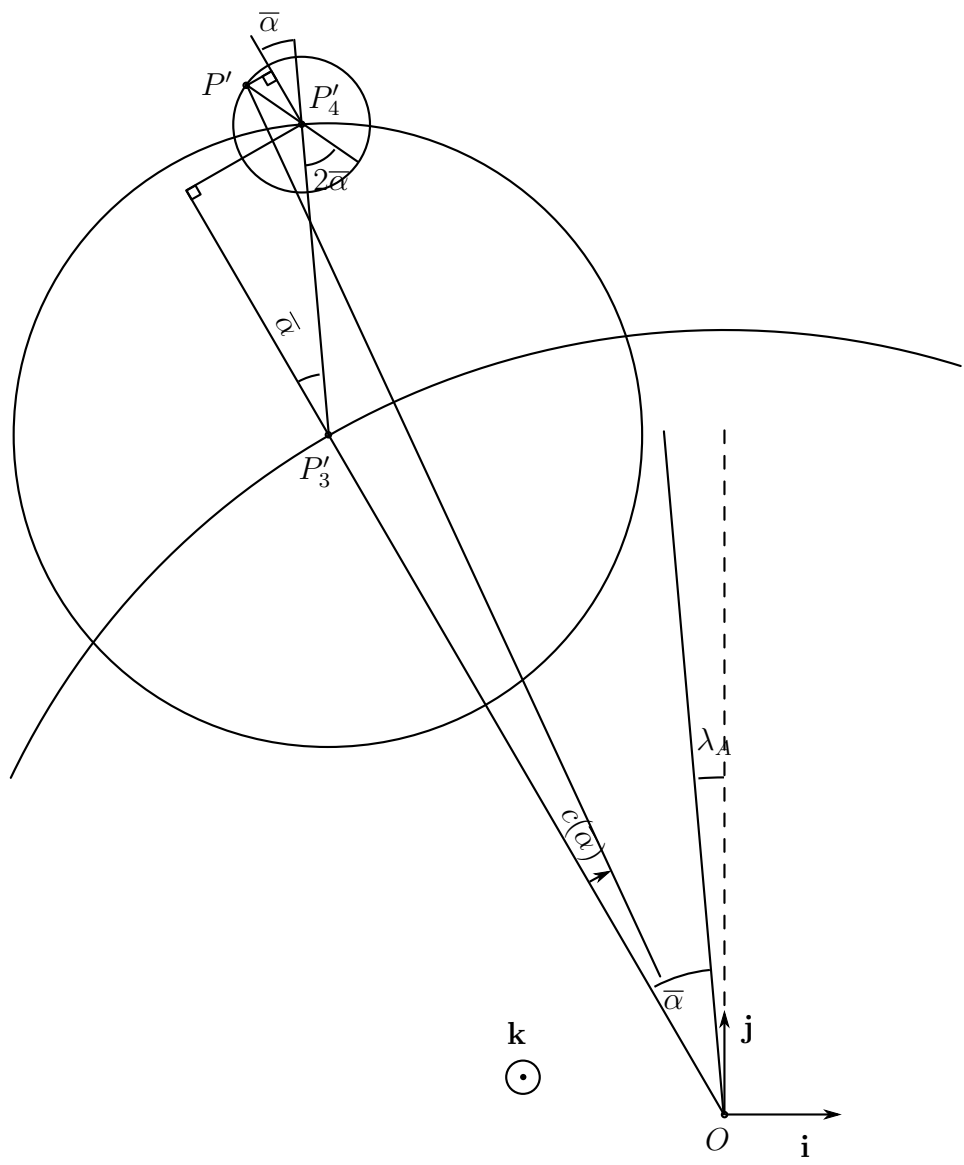
    \caption{\label{fig004}Le Soleil : transformations planes}
  \end{center}
\end{figure}

\begin{figure}
  \begin{center}
    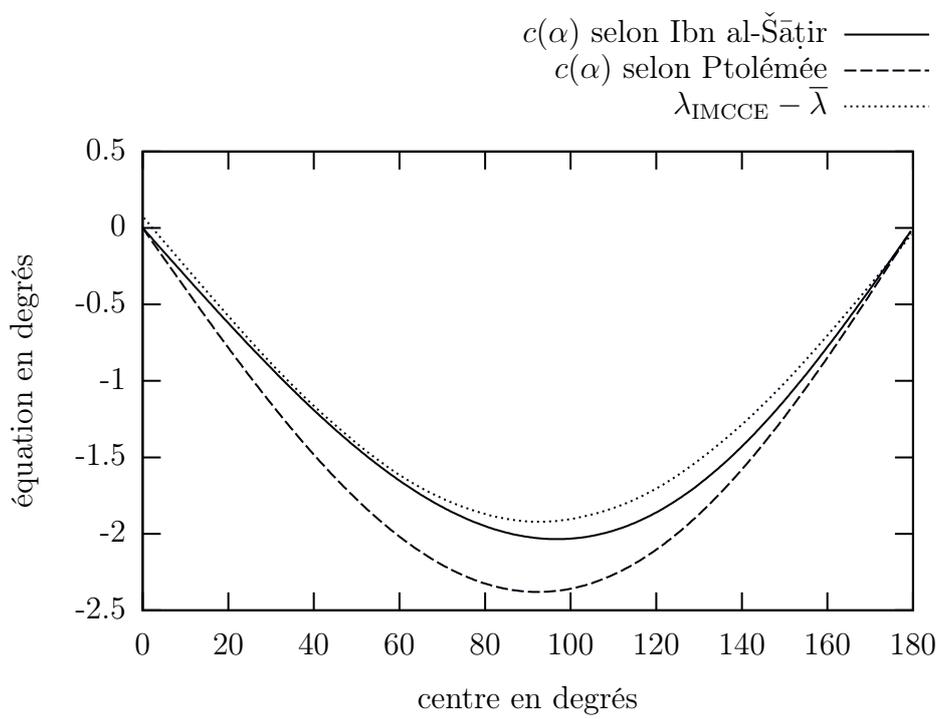
    \caption{\label{fig005}L'équation du Soleil}
  \end{center}
\end{figure}

\paragraph{L'équation du Soleil}
Le calcul de la position du point $P'$ revient à résoudre deux
triangles rectangles (\textit{cf}. fig.~\ref{fig004} p.~\pageref{fig004})
pour calculer l'angle $c(\overline{\alpha})$ traditionnellement
appelé <<~équation du Soleil~>>~:
$$c(\overline{\alpha})=(\overrightarrow{OP_3'},\overrightarrow{OP'})=
-\arcsin\left(\frac{P_3P_4\sin\overline{\alpha}-P_4P\sin\overline{\alpha}}{OP'}\right)$$
où
$$OP'=\sqrt{(P_3P_4\sin\overline{\alpha}-P_4P\sin\overline{\alpha})^2+
(OP_3+P_3P_4\cos\overline{\alpha}+P_4P\cos\overline{\alpha})^2}$$
La longitude de $P'$, depuis le point vernal, est alors~:
$$\left(\mathbf{j},\overrightarrow{OP'}\right)=\lambda_A+\overline{\alpha}+c(\overline{\alpha})$$
On remarque que~:
$$c(360°-\overline{\alpha})=-c(\overline{\alpha})$$
On a tracé le graphe de l'équation du Soleil en fonction de $\overline{\alpha}$, \textit{cf.} figure \ref{fig005}, où l'on a aussi représenté à des fins de comparaison l'équation calculée suivant la méthode de Ptolémée\footnote{\textit{cf.} \cite{pedersen1974} p.~150.}~:
$$c_{\text{Ptolémée}}(\overline{\alpha})=-\arcsin\left(\frac{2;29,30\times\sin\overline{\alpha}}{\sqrt{(2;29,30\times\sin\overline{\alpha})^2+(60+2;29,30\times\cos\overline{\alpha})^2}}\right)$$
La figure \ref{fig005} montre aussi l'équation telle qu'aurait pu l'observer {\shatir} en 1331, en retranchant à la longitude vraie observée le soleil moyen calculé avec les paramètres $\overline{\lambda}_0,\dot{\overline{\lambda}}$ mentionnés ci-dessus.

\paragraph{Application} On a vu que {\shatir } donnait la longitude du Soleil moyen à $t=0$ : c'est $\overline{\lambda}_0=280°9'0''$. Pour avoir la position exacte du Soleil le long de l'écliptique à cette date, il faut corriger cette valeur par l'équation du Soleil comme on l'explique dans la section précédente. On trouve $\lambda=280°33'31''$. \`A titre de comparaison, l'Institut de Mécanique Céleste et de Calcul des \'Ephémérides de l'Observatoire de Paris donne $\lambda=280°29'$ le 24 décembre 1331 à 9 h 43 min.

\begin{figure}
  \begin{center}
    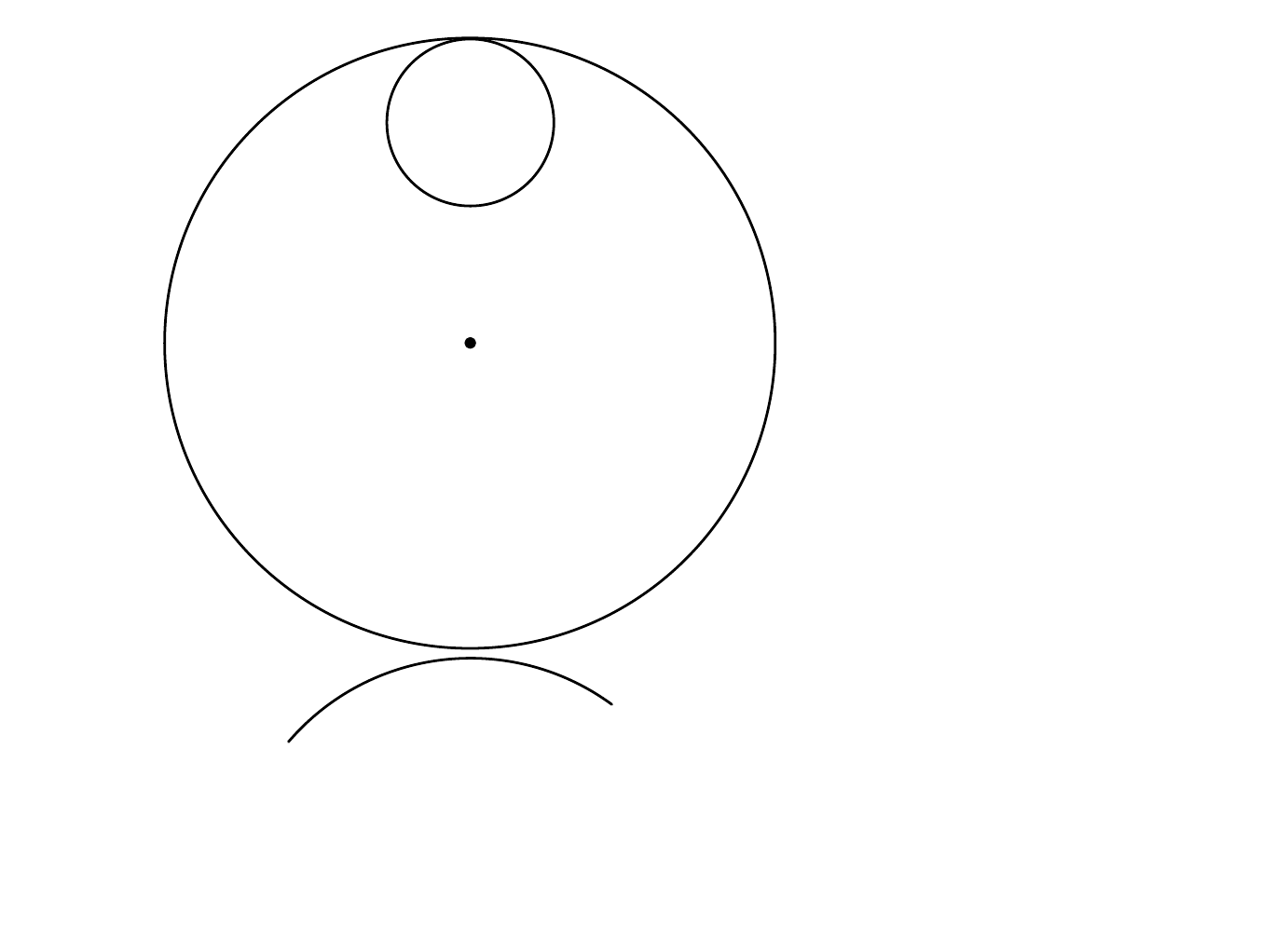
    \caption{\label{fig006}Les orbes solides du Soleil}
  \end{center}
\end{figure}

\paragraph{Distances extrémales et dimensions des orbes}
Tandis que la figure \ref{fig003} représentait, pour chaque orbe, la trajectoire circulaire de son centre au sein de l'orbe le portant, la figure \ref{fig006} est une section des corps à symétrie sphérique constituant les orbes : chaque cercle y représente le bord d'un orbe. Cette figure n'est pas à l'échelle, mais on y a écrit les dimensions des orbes, en unités telles que $OP_3=60$. Il est alors facile de calculer la distance maximale du Soleil à la Terre, $OP=67;7$, et sa distance minimale, $52;53$. Soit $r$ le rayon du globe solaire lui-même. On a~:
\begin{align*}
  \text{rayon du rotateur} &=2;30+r\\
  \text{rayon du déférent} &=4;37+2;30+r=7;7+r\\
  \text{rayon extérieur du parécliptique} &=60+7;7+r=67;7+r\\
  \text{rayon intérieur du parécliptique} &=60-(7;7+r)=52;53-r
\end{align*}
L'épaisseur du parécliptique serait donc $14;14+2r$. Or {\shatir} nous dit, p.~\pageref{complement_sol}~:
\begin{quote}
  ``Son épaisseur est quatorze degrés et trente-quatre minutes, plus un complément. Pour que [les orbes] soient mutuellement contigus quand on y plonge l'orbe déférent, nous prenons un complément d'un tiers de degré à l'extérieur et à l'extérieur, en parts telles que son rayon extérieur est soixante-sept parts et dix-sept minutes.''
\end{quote}
Une épaisseur de $14;34$ inclut donc déjà un ``complément'' de $0;20$ égal au diamètre supposé du globe solaire dans le chapitre 7~; mais cette valeur est provisoire, comme on le voit dans la conclusion de la première partie de la \textit{Nihaya}. D'une part, l'observation permet de déterminer le diamètre du globe solaire (\textit{cf.} p.~\pageref{obs_diam_sol}). D'autre part, il faudra prendre un ``complément'' suffisant pour que la partie concave de l'orbe du Soleil touche exactement la partie convexe du système d'orbes de Vénus (\textit{cf.} p.~\pageref{contiguite_venus_sol}). Enfin, l'épaisseur de $0;43$ choisie pour l'orbe total dans le chapitre 7 a aussi un caractère provisoire~; elle n'influe en rien sur la trajectoire décrite par le Soleil dans ce modèle~; et {\shatir} l'oubliera complètement dans sa conclusion.

Pour chaque astre, {\shatir} trace toujours deux figures, analogues des figures \ref{fig003} et \ref{fig006}. Celle représentant les trajectoires des centres des orbes est très utile à la compréhension du modèle mathématique : elle permet de situer la position des centres des orbes et donc aussi la position de l'astre. L'autre figure, représentant les ``orbes solides'', relève davantage de l'astronomie physique et de la nature des orbes ; elle nous offre une image du Monde par ses sections planes, et elle permet aussi d'en calculer les dimensions sous l'hypothèse que les systèmes d'orbes appartenant aux différents astres s'emboitent les uns dans les autres sans qu'il n'y ait de vide.

\paragraph{Diamètre apparent} Le diamètre apparent du Soleil quand il est à distance moyenne ($OP'=60$) est, selon {\shatir}, $0;32,32$ degrés. \`A partir de cette observation et des distances extrémales calculées précédemment, on trouve facilement le diamètre apparent à l'apogée~:
$$\frac{0;32,32\times 60}{67;7}\simeq 0;29,5$$
et au périgée~:
$$\frac{0;32,32\times 60}{52;53}\simeq 0;36,55$$
D'après Saliba\footnote{\textit{cf.} \cite{rashed1997} p.~103.}, {\shatir} est le premier auteur étudiant sérieusement la variation du diamètre apparent du Soleil, en lien avec la possibilité des éclipses annulaires que n'avait pas envisagées Ptolémée ; il y avait là motivation à offrir un nouveau modèle du Soleil. Ces observations suffisent en effet à déterminer $OP_3+P_3P_4+P_4P$ et $OP_3-P_3P_4-P_4P$, d'où l'on tire $P_3P_4+P_4P$, à un facteur d'échelle près permettant de choisir $OP_3=60$. Pour régler entièrement les paramètres du modèle, il suffirait enfin de déterminer $P_3P_4-P_4P$. Saliba suggère qu'{\shatir} a pu le faire en observant l'équation maximale. {\shatir} mentionne des observations concernant le diamètre apparent du Soleil, p.~\pageref{diam_sol} et \pageref{diam_sol2}.

Sans remettre en cause l'analyse de Saliba, nous aimerions suggérer que ce scénario n'est pas le seul possible. Un premier ensemble d'observations aurait pu permettre de déterminer l'équation maximale ou bien l'équation dans les quadratures. On montre facilement que l'équation dans les quadratures vaut :
$$c(90°)=-\arctan\left(\frac{P_3P_4-P_4P}{OP_3}\right)$$
Quant à l'équation maximale $c_{\max}$, on montre facilement qu'elle est atteinte lorsque
$$\cos\overline{\alpha}=-\frac{P_3P_4+P_4P}{OP_3},$$
et on a alors :
$$(OP_3^2-(P_3P_4+P_4P)^2)\tan^2c_{\max}=(P_3P_4-P_4P)^2$$
{\shatir} aurait ensuite pu utiliser des observations concernant l'équation du Soleil dans les octants : il mentionne lui-même ces observations comme étant un argument en faveur de ses modèles\footnote{\textit{cf. supra} l'introduction d'{\shatir} p.~\pageref{doutes_excentrique}.}. On montre facilement que :
$$\frac{P_3P_4-P_4P}{OP_3}=-\left(\sqrt{2}+\frac{P_3P_4+P_4P}{OP_3}\right)\tan c(45°)$$
Finalement, les observations concernant le diamètre apparent auraient très bien pu n'être utilisées que comme test d'un modèle déjà achevé. Dans tous les cas, il reste à étudier \emph{comment} de telles observations -- diamètre apparent, équation maximale, équation dans les quadratures, dans les octants -- ont pu être menées, et à quelle précision, avec les moyens de l'époque.


\paragraph{La Lune}
La trajectoire de la Lune est contenue dans un plan incliné par rapport au plan de l'écliptique. {\shatir} conçoit donc un \emph{orbe incliné}, dont le plan est immobile au sein du référentiel constitué par l'orbe parécliptique de la Lune. Ce plan incliné coupe le plan du parécliptique suivant une droite, la direction des <<~n{\oe}uds~>>. Les n{\oe}uds sont deux points du parécliptique, diamètralement opposés, dans cette direction. Comme l'orbe parécliptique de la Lune entraîne l'orbe incliné dans son mouvement, alors les n{\oe}uds se déplacent de manière uniforme le long de l'écliptique, par rapport au point vernal.

\begin{figure}
\begin{center}
  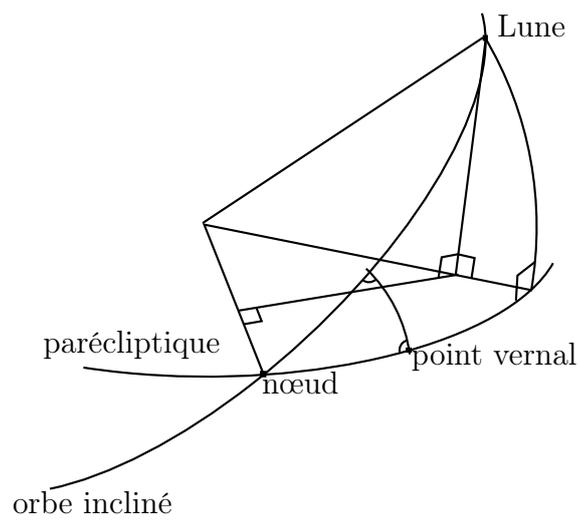
  \caption{\label{latitude}Mouvement en latitude de la Lune}
\end{center}
\end{figure}

Pour décrire le mouvement des points situés sur le plan de l'orbe incliné, {\shatir} utilise comme point de référence un point situé <<~en face du point vernal~>>~: l'arc de l'orbe incliné compris entre ce point et le n{\oe}ud ascendant est égal à l'arc de l'écliptique compris entre le point vernal et le n{\oe}ud ascendant (voir figure \ref{latitude}). {\shatir} devait en effet rabattre le plan de l'orbe incliné dans le plan de l'écliptique au moyen d'une rotation dont l'axe est la direction des n{\oe}uds. Le point <<~en face du point vernal~>> devient alors une origine naturelle pour les mesures angulaires.

\begin{figure}
\begin{center}
  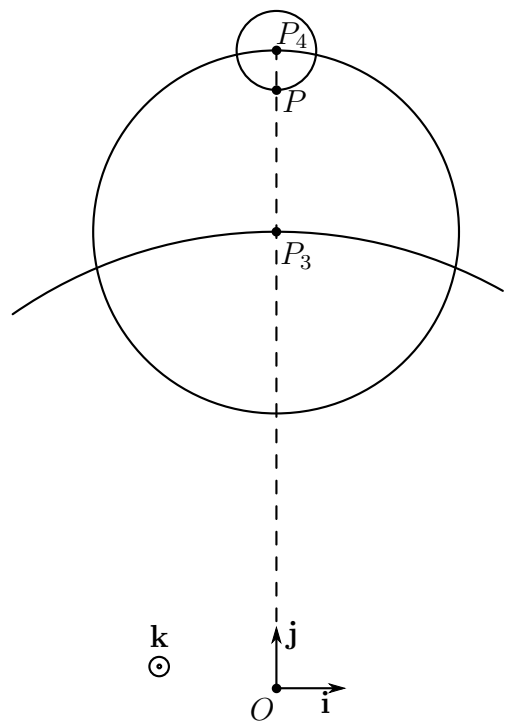
  \caption{\label{fig008}Les orbes de la Lune : figure initiale}
\end{center}
\end{figure}

\paragraph{La Lune : figure initiale}
Soit $(\mathbf{i},\mathbf{j},\mathbf{k})$ la base canonique de
$\mathbb{R}^3$. Dans la figure initiale, les plans des orbes sont tous
rabattus dans le plan de l'écliptique $(\mathbf{i},\mathbf{j})$ par
des rotations, les centres des orbes sont tous alignés dans la
direction du vecteur $\mathbf{j}$ elle-même confondue avec la
direction du point vernal~; enfin, la ligne des n{\oe}uds, à l'intersection
des plans de l'orbe incliné et du parécliptique, est orientée dans la
direction du vecteur $\mathbf{j}$, le n{\oe}ud ascendant étant du même
côté que $\mathbf{j}$. Le point $O$ est le centre du Monde, $P_1$
le centre du parécliptique de la Lune, et $P_2$ le centre de son orbe
incliné. Ces trois points sont confondus $O=P_1=P_2$. Le point $P_3$ est
le centre de l'orbe de l'épicycle, $P_4$ le centre de l'orbe rotateur,
et $P$ le centre de la Lune, cet astre étant lui-même un corps sphérique.
La figure \ref{fig008} précise la position de ces points, mais elle n'est
pas à l'échelle. Bien qu'ils soient alignés, on remarque que les rayons
vecteurs $\overrightarrow{P_nP_{n+1}}$ ne sont pas tous dans le même sens.
On pose~:
$\overrightarrow{OP_3}=60\ \mathbf{j}$, alors
$$\begin{array}{l}
\overrightarrow{P_3P_4}=6;35\ \mathbf{j}\\
\overrightarrow{P_4P}=-1;25\ \mathbf{j} 
\end{array}$$
Il faut appliquer au point $P$ une suite de transformations géométriques
pour obtenir la position de la Lune prédite ou observée à un instant
donné.

\paragraph{La Lune : transformations géométriques}
Les transformations appliquées à la Lune $P$ sont des rotations paramétrées par quatre angles~: le \emph{n{\oe}ud} $\lambda_{\ascnode}$, la \emph{Lune propre} $\overline{\alpha}$, le \emph{centre de la Lune} $2\overline{\eta}$, et l'argument de latitude\footnote{Cette grandeur mesure un arc le long de l'orbe incliné entre le n{\oe}ud ascendant et le centre de l'épicycle~; il ne faut donc pas la confondre avec une latitude.} moyen $(\overline{\lambda}+\lambda_{\ascnode})$. Enfin, $\overline{\lambda}$ désigne la \emph{Lune moyenne}. Voici la liste des transformations géométriques utilisées\footnote{Il règne une légère ambiguïté concernant $\dot{\lambda}_{\tiny\ascnode}$ dans le texte d'{\shatir}. Quant il décrit le mouvement du parécliptique de la Lune, qu'il dit être $0;3,10,38,27$ par jour, il ajoute <<~en vérité, à proprement parler, ce mouvement est l'excédent du mouvement des n{\oe}uds sur le mouvement des étoiles fixes~>> (p.~\pageref{lab001} \textit{supra}). Est-ce à dire que la valeur donnée par {\shatir} est la différence entre le mouvement du parécliptique de la Lune et le mouvement des fixes (précession), tous deux relatifs au référentiel solide du neuvième orbe~? Ce serait donc leur somme algébrique puisque ces deux mouvements vont en sens contraires~; c'est l'interprétation que nous avons retenue dans notre commentaire mathématique. Dans ce cas, la rotation $R(P_1,\mathbf{k},-\lambda_{\tiny\ascnode})$ rend en fait compte du mouvement des fixes \textit{et} du mouvement du parécliptique~; le <<~n{\oe}ud moyen~>> et le <<~mouvement du n{\oe}ud moyen~>> désignent alors des grandeurs mesurées par rapport au référentiel solide du neuvième orbe et décrivant le mouvement du parécliptique au sein de ce référentiel. {\shatir} dit plus loin (p.~\pageref{lab002}) que le n{\oe}ud moyen est l'arc <<~entre le commencement du Bélier et le n{\oe}ud ascendant, dans l'orbe parécliptique~>>. Si notre interprétation est juste, il faut donc penser qu'il commet un abus de langage en désignant le point vernal par la locution <<~commencement du Bélier~>>, puisque le commencement du Bélier désigne un point de l'écliptique qui subit le mouvement de précession, et non un point du neuvième orbe. Mais notre interprétation est confortée par le sort analogue dévolu au mouvement des apogées pour les autres astres. Par exemple, dans le modèle du Soleil vu ci-dessus, la grandeur $\dot{\lambda}_A$ du mouvement de l'apogée (un degré toutes les soixante années persanes) est bien la somme du mouvement de précession (un degré toutes les soixante-dix années persanes) et du mouvement de l'orbe total.}~:
$$R(P_4,\mathbf{k},2\overline{\eta}),\quad R(P_3,\mathbf{k},-\overline{\alpha}),\quad
R(P_2,\mathbf{k},\overline{\lambda}+\lambda_{\ascnode}),\quad R(P_2,\mathbf{j},5°), \quad
R(P_1,\mathbf{k},-\lambda_{\ascnode}).$$
La rotation dont l'axe est dans la direction du vecteur $\mathbf{j}$ a pour effet d'incliner le plan de l'orbe incliné par rapport au plan du parécliptique. La description d'Ibn al-\v{S}\=a\d{t}ir ne laisse aucun
doute quant à l'ordre dans lequel appliquer ces quatres transformations.

Pour obtenir la configuration des orbes à un instant donné, il faut
donc appliquer toutes ces rotations aux points $P_3$, $P_4$, $P$.
L'image du point $P_3$ entraîné par les mouvements du parécliptique et
de l'orbe incliné est~:\label{composees_lune}
$$R(P_1,\mathbf{k},-\lambda_{\ascnode})\circ R(P_2,\mathbf{j},5°)
\circ R(P_2,\mathbf{k},\overline{\lambda}+\lambda_{\ascnode})\ (P_3)$$
Quant au point $P_4$, il est aussi
entraîné par le mouvement de l'orbe de l'épicycle et devient~:
$$R(P_1,\mathbf{k},-\lambda_{\ascnode})\circ R(P_2,\mathbf{j},5°)
\circ R(P_2,\mathbf{k},\overline{\lambda}+\lambda_{\ascnode})\circ R(P_3,\mathbf{k},-\overline{\alpha})\ (P_4)$$
Enfin le point $P$, aussi entraîné par le mouvement de l'orbe rotateur,
devient~:
$$R(P_1,\mathbf{k},-\lambda_{\ascnode})\circ R(P_2,\mathbf{j},5°)
\circ R(P_2,\mathbf{k},\overline{\lambda}+\lambda_{\ascnode})\circ R(P_3,\mathbf{k},-\overline{\alpha})
\circ R(P_4,\mathbf{k},2\overline{\eta})\ (P)$$

\paragraph{La Lune : trajectoire paramétrée}
Le modèle de la Lune est couplé au mouvement du Soleil moyen dont on notera désormais $\dot{\overline{\lambda}}_{\astrosun}$ et $\overline{\lambda}_{\astrosun}(0)$ les paramètres notés précédemment $\dot{\overline{\lambda}}$ et $\overline{\lambda}(0)$, puisque ces derniers symboles concernent désormais la Lune moyenne. Pour la Lune, la trajectoire dans l'ensemble des valeurs des paramètres est~:
$$\lambda_{\ascnode}=\dot{\lambda}_{\ascnode}t+\lambda_{\ascnode}(0)$$
$$\overline{\lambda}+\lambda_{\ascnode}=(\dot{\overline{\lambda}}+\dot{\lambda}_{\ascnode}) t+\overline{\lambda}(0)+\lambda_{\ascnode}(0)$$
$$\overline{\alpha}=\dot{\overline{\alpha}}t+\overline{\alpha}(0)$$
$$2\overline{\eta}=2(\dot{\overline{\lambda}}-\dot{\overline{\lambda}}_{\astrosun})t+2(\overline{\lambda}(0)-\overline{\lambda}_{\astrosun}(0))$$
On a, d'après {\shatir}~:
$$\overline{\lambda}_{\astrosun}(0)=280;9,0,\quad\dot{\overline{\lambda}}_{\astrosun}=0;59,8,19,43,33,41,55,4,6\text{ par jour solaire moyen}$$
$$\overline{\lambda}(0)=213;35,50,\quad\dot{\overline{\lambda}}=13;10,35,1,13,53\text{ par jour solaire moyen}$$
$$\overline{\alpha}(0)=138;32,27,\quad\dot{\overline{\alpha}}=13;3,53,56\text{ par jour solaire moyen}$$
$$\lambda_{\ascnode}(0)=275;7,35,\quad\dot{\lambda}_{\ascnode}=0;3,10,38,27\text{ par jour solaire moyen}$$
d'où $\dot{\overline{\lambda}}+\dot{\lambda}_{\ascnode}=13;13,45,39,40$ par jour solaire moyen. On appelle <<~élongation moyenne~>> la grandeur $\overline{\eta}=\overline{\lambda}-\overline{\lambda}_{\astrosun}$. 

\paragraph{La Lune : transformations planes}
Si $R$ et $S$ sont deux rotations dans l'espace, il existe une
rotation $T$ telle que $R\circ S=T\circ R$, à savoir, la rotation de
même angle que $S$ et dont l'axe est l'image par $R$ de l'axe de
$S$. Appliquons cette relation de commutation à toutes les composées
de rotations décrites p.~\pageref{composees_lune} ci-dessus, de sorte à
réécrire à droite toutes les rotations dont l'axe est dans la
direction du vecteur $\mathbf{k}$. Il suffit de remarquer que~:
$$R(P_1,\mathbf{k},-\lambda_{\ascnode})\circ R(P_2,\mathbf{j},5°)=
R(P_2,\mathbf{u},5°)\circ R(P_1,\mathbf{k},-\lambda_{\ascnode})$$ 
où le vecteur $\mathbf{u}$ est l'image de $\mathbf{j}$ par
$R(P_1,\mathbf{k},-\lambda_{\ascnode})$. L'image de $P$ est donc~:
$$R(P_1,\mathbf{u},5°)(P')$$
où
$$P'=R(P_1,\mathbf{k},-\lambda_{\ascnode})\circ R(P_2,\mathbf{k},\overline{\lambda}+\lambda_{\ascnode})
\circ R(P_3,\mathbf{k},-\overline{\alpha})\circ R(P_4,\mathbf{k},2\overline{\eta})\ (P).$$
On introduit de même les points $P_3'$ et $P_4'$ suivants~:
$$P_3'=R(P_1,\mathbf{k},-\lambda_{\ascnode})\circ R(P_2,\mathbf{k},\overline{\lambda}+\lambda_{\ascnode})\ (P_3)$$
$$P_4'=R(P_1,\mathbf{k},-\lambda_{\ascnode})
\circ R(P_2,\mathbf{k},\overline{\lambda}+\lambda_{\ascnode})\circ R(P_3,\mathbf{k},-\overline{\alpha})\ (P_4)$$
Tous ces points sont dans le plan de la figure initiale~: calculer
leurs positions relève entièrement de la géométrie plane.

\begin{figure}
\begin{center}
  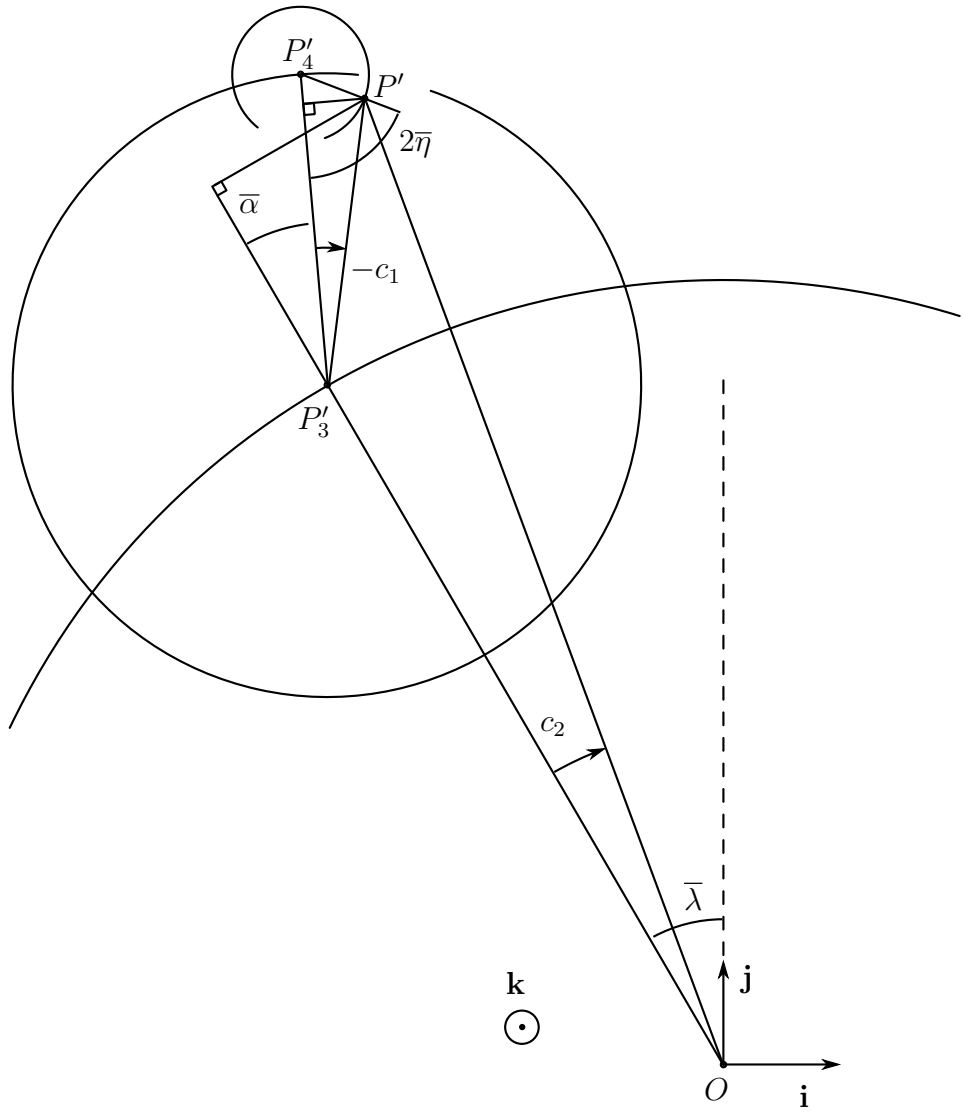
  \caption{\label{fig009}La Lune : transformations planes}
\end{center}
\end{figure}

\paragraph{Equations de la Lune}\label{equ_lune}
Le calcul de la position du point $P'$ revient à résoudre les triangles
rectangles représentés sur la figure \ref{fig009}. On calcule d'abord
l'angle $c_1=-(\overrightarrow{P_3'P_4'},\overrightarrow{P_3'P'})$ appelé <<~équation de la Lune propre~>> (cf.~p.~\pageref{c1}). C'est une fonction de $2\overline{\eta}$~:
$$c_1=\arcsin\left(\frac{P_4P\sin2\overline{\eta}}{P_3'P'}\right)$$
où $P_3'P'$ est le <<~rayon de l'épicycle apparent~>>~:
$$P_3'P'=\sqrt{(P_3P_4-P_4P\cos2\overline{\eta})^2+(P_4P\sin2\overline{\eta})^2}$$
On calcule ensuite une autre <<~équation~>>
$c_2=(\overrightarrow{OP_3'},\overrightarrow{OP'})$~:
$$c_2=\arcsin\left(\frac{P_3'P'\sin(-\overline{\alpha}-c_1)}{OP'}\right)$$
où
$$OP'=\sqrt{(OP_3+P_3'P'\cos(\overline{\alpha}+c_1))^2+(P_3'P'\sin(\overline{\alpha}+c_1))^2}$$
La position de $P'$ est alors donnée par $OP'$ et par l'angle suivant~:
$$\left(\mathbf{j},\overrightarrow{OP'}\right)=\overline{\lambda}+c_2$$
La fonction $c_1$ est une fonction d'une seule variable, et $c_2$ fonction des deux variables $\overline{\alpha}$ et $2\overline{\eta}$. Dans ce qui suit, $c_2$ sera plutôt conçue comme une fonction de $\alpha=\overline{\alpha}+c_1$ et de $2\overline{\eta}$, notée $c_2(\alpha,2\overline{\eta})$. Plus précisément,
$$c_2(\alpha,2\overline{\eta})=\arcsin\left(\frac{P_3'P'\sin(-\alpha)}{\sqrt{(OP_3+P_3'P'\cos\alpha)^2+(P_3'P'\sin\alpha)^2}}\right)$$
où $P_3'P'$ est la fonction de $2\overline{\eta}$ vue ci-dessus.
Ces fonctions sont bien sûr $2\pi$-périodiques par rapport à chacune des variables, mais on remarque aussi que~:
$$c_1(360°-2\overline{\eta})=-c_1(2\overline{\eta})$$
$$c_2(360°-\alpha,\ 360°-2\overline{\eta})=-c_2(\alpha,2\overline{\eta})$$

\paragraph{Fonction de deux variables et interpolation}
Les valeurs de $c_1$ et $c_2$ devront être reportées dans des tables ; mais $c_2$ est fonction de deux variables, et il faudrait donc construire une table pour chaque valeur de $2\overline{\eta}$. Les astronomes, à partir de Ptolémée\footnote{\textit{cf.} \cite{pedersen1974}, p.~84-89.}, utilisaient dans ce contexte une méthode d'interpolation visant à calculer une valeur approchée de $c_2$ au moyen d'un produit d'une fonction de $\alpha$ par une fonction de $2\overline{\eta}$. \`A $\alpha$ donné, {\shatir} interpole les valeurs de $c_2$ entre $c_2(\alpha,0)$ et $c_2(\alpha,180°)$ au moyen de la formule suivante (qui \emph{n'est pas} linéaire en $2\overline{\eta}$)~:
$$c_2(\alpha,2\overline{\eta})\simeq c_2(\alpha,0)
+\chi(2\overline{\eta})(c_2(\alpha,180°)-c_2(\alpha,0)).$$
Le coefficient d'interpolation $\chi$ est défini par~:
$$\chi(2\overline{\eta}) =\frac{\max\vert c_2(\cdot,2\overline{\eta})\vert-\max\vert
  c_2(\cdot,0)\vert} {\max\vert c_2(\cdot,180°)\vert-\max\vert
  c_2(\cdot,0)\vert}$$ 
Ce coefficient est fonction d'une seule variable,
on peut donc le calculer pour une série de valeurs de $2\overline{\eta}$ et en
faire une table. Pour ce faire, on remarque que le
maximum $\max\vert c_2(\cdot,2\overline{\eta})\vert$ est atteint quand $(\tan(c_2(\alpha,2\overline{\eta})))^2$ est maximum. En posant $z^{-1}=OP_3+P_3'P'\cos\alpha$, on montre que $(\tan(c_2(\alpha,2\overline{\eta})))^2$ est un polynôme de degré 2 en $z$. On calcule facilement son maximum\footnote{{\shatir} procède géomètriquement, en considérant une droite tangente à l'<<~épicycle apparent~>> qui est un cercle de rayon $P_3'P'$.}, il est atteint quand
$$\cos\alpha=-\frac{P_3'P'}{OP_3}$$
$$\vert\sin\alpha\vert=\sqrt{1-\left(\frac{P_3'P'}{OP_3}\right)^2}.$$
On reporte ces valeurs dans $\vert c_2(\alpha,2\overline{\eta})\vert$, et on trouve~:
$$\max\vert c_2(\cdot,2\overline{\eta})\vert=\arcsin\frac{P_3'P'}{OP_3}.$$

\begin{figure}
\begin{center}
  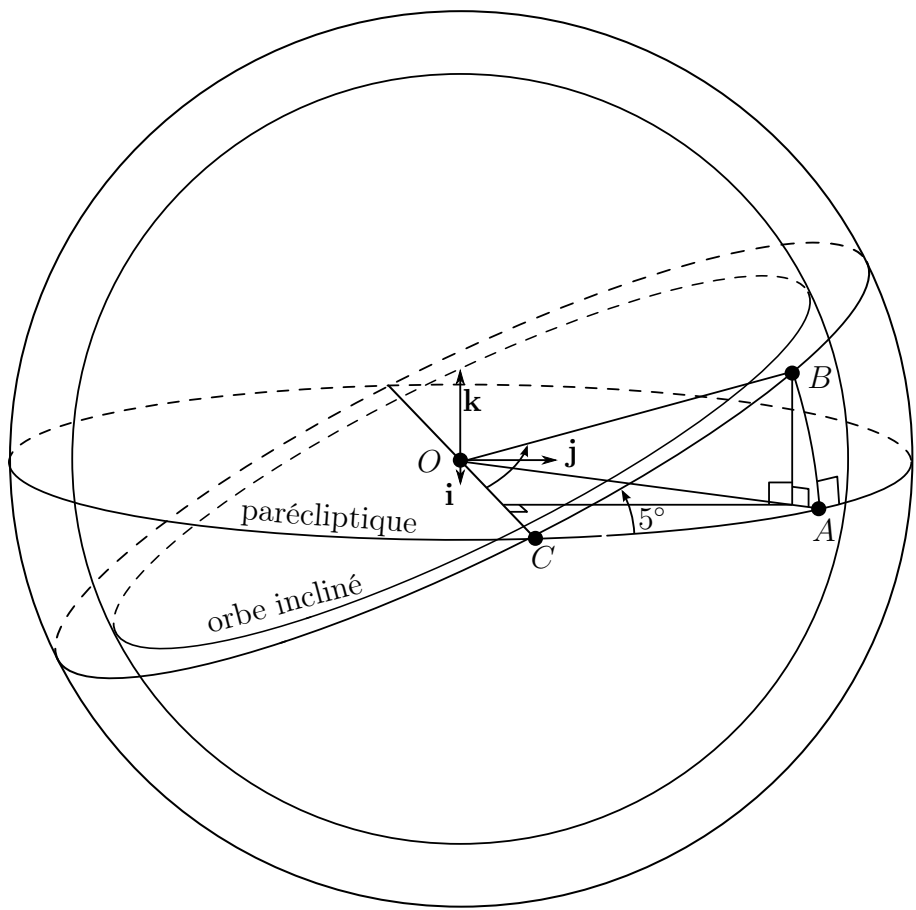
  \caption{\label{fig010}La Lune : coordonnées de $R(P_1,\mathbf{u},5°)(P')$}
\end{center}
\end{figure}

\paragraph{Trigonométrie sphérique} 
On va à présent calculer les coordonnées sphériques du point
$R(P_1,\mathbf{u},5°)(P')$ par rapport à l'écliptique. L'angle
formé entre le vecteur $\mathbf{u}$ et la direction du
point $P'$ vaut (modulo 360°)~:
$$(\mathbf{u},\overrightarrow{OP'})=(\mathbf{j},\overrightarrow{OP'})-(\mathbf{j},\mathbf{u})=\overline{\lambda}+c_2+\lambda_{\ascnode}$$
Les égalités suivantes seront aussi prises modulo 360°. Sur la figure \ref{fig010}, on a représenté sur la sphère de l'écliptique le point $B$ dans la direction du point $R(P_1,\mathbf{u},5°)(P')$, et le point $C$ dans la direction du vecteur $\mathbf{u}$. En particulier,
$$(\overrightarrow{OC},\overrightarrow{OB})=\overline{\lambda}+\lambda_{\ascnode}+c_2.$$
L'étude du triangle sphérique $ABC$ donne~:
$$\tan(\overrightarrow{OC},\overrightarrow{OA})
=\cos(5°)\times\tan(\overline{\lambda}+\lambda_{\ascnode}+c_2)$$
La longitude du point $R(P_1,\mathbf{u},5°)(P')$ par rapport
à l'écliptique, en prenant la direction du point vernal, c'est-à-dire 
$\mathbf{j}$, comme origine, est donc, quand
$\overline{\lambda}+\lambda_{\ascnode}+c_2\in\ ]-90°,\ 90°[$~:
$$(\mathbf{j},\overrightarrow{OA})=(\overrightarrow{OC},\overrightarrow{OA})-(\overrightarrow{OC},\mathbf{j})=\arctan(\cos(5°)\times\tan(\overline{\lambda}+\lambda_{\ascnode}+c_2))-\lambda_{\ascnode}$$
    Quand au contraire $\overline{\lambda}+\lambda_{\ascnode}+c_2\in\ ]90°,\ 270°[$, on a~:
    $$(\mathbf{j},\overrightarrow{OA})=180°+\arctan(\cos(5°)\times\tan(\overline{\lambda}+\lambda_{\ascnode}+c_2))-\lambda_{\ascnode}.$$
    On rassemble ces deux cas dans la formule suivante, valable pour toutes les valeurs des paramètres :
    $$(\mathbf{j},\overrightarrow{OA})=\overline{\lambda}+c_2+e_n(\overline{\lambda}+\lambda_{\ascnode}+c_2)$$
    où l'<<~équation du déplacement~>> $e_n(x)$ est définie comme suit sur l'intervalle $[-90°,270°]$, et ailleurs par périodicité :
    $$e_n(x)=\left\lbrace\begin{array}{l}
    \arctan(\cos(5°)\times\tan(x))-x,\text{ si }x\in\ ]-90°,\ 90°[\\
    180°+\arctan(\cos(5°)\times\tan(x))-x,\text{ si }x\in\ ]90°,\ 270°[\\
    0°,\text{ si }x=\pm 90°
    \end{array}\right.$$
Enfin, la latitude de la Lune est~:
$$(\overrightarrow{OA},\overrightarrow{OB})=
\arcsin(\sin(5°)\times\sin(\overline{\lambda}+\lambda_{\ascnode}+c_2).$$
On peut aussi l'obtenir ainsi (c'est la solution retenue par {\shatir} au chapitre 11)~:
$$(\overrightarrow{OA},\overrightarrow{OB})=
\arctan(\tan(5°)\times\sin(\overline{\lambda}+\lambda_{\ascnode}+c_2+e_n(\overline{\lambda}+\lambda_{\ascnode}+c_2))).$$

\paragraph{La Lune : critique des modèles antérieurs} Les critiques
adressées par {\shatir} à ses prédécesseurs p.~\pageref{doutes_lune}
concernent~:
\begin{itemize}
\item l'usage d'un orbe \emph{excentrique}
\item la question de l'uniformité du mouvement de rotation d'un orbe
  excentrique par rapport à un axe passant par son centre
\item ``le fait que le rayon de l'épicycle suive un point autre que le
  centre de l'orbe qui le porte'' (en particulier le point de
  \textit{prosneuse} dans le troisième modèle de Ptolémée)
\item la variation du \emph{diamètre apparent} de la Lune, vu de la
  Terre
\item l'équation de la Lune dans les \emph{octants}
\end{itemize}
Pour mieux évaluer les critiques adressées par {\shatir} à ses
prédécesseurs, et la précision de ses propres modèles et des
observations dont il a pu les déduire, on procèdera, pour chaque
astre, au choix d'une période de référence commençant à l'époque
choisie par {\shatir}. Pour la Lune, on a choisi trente jours à partir
du 24 décembre 1331. \`A des fins de comparaison, il faut remarquer
que les paramètres des modèles planétaires exprimés en ``mouvements
moyens'' ont toujours une valeur phénoménologique indépendante du
modèle géométrique considéré. Par exemple pour la Lune,
$\overline{\lambda}$ est la direction d'un point se mouvant
uniformément sur l'écliptique et dont la vitesse est la vitesse
\emph{moyenne} de la Lune en longitude\footnote{à ne pas confondre
  avec la vitesse instantanée, d'où notre notation pour la vitesse
  moyenne $\dot{\overline{\lambda}}$ qui est une constante, à ne pas
  confondre avec $\dot{\lambda}$.}. La longitude vraie $\lambda$ de la
Lune oscille autour de $\overline{\lambda}$ ; en théorie, on peut donc
déterminer $\overline{\lambda}$ statistiquement, indépendammant d'un
quelconque modèle géométrique. Il est encore plus facile de calculer
la vitesse moyenne $\dot{\overline{\lambda}}$~: c'est l'inverse d'une
période de retour aux mêmes longitudes. Les Babyloniens savaient déjà
la calculer en moyennant sur un grand nombre de périodes.

Pour comparer les modèles géométriques d'{\shatir} et ceux de ses
prédecesseurs, on pourra donc les utiliser pour calculer la
trajectoire d'un \emph{même} astre sur une \emph{même} période de
référence, en utilisant un \emph{même} jeu de paramètres concernant
les ``mouvements moyens''. Les résultats d'une telle comparaison
seront nécessairement partiels car aucun astre -- pas même le Soleil
si l'on tient compte du mouvement de l'Apogée -- n'a un mouvement
strictement périodique dans les modèles les plus élaborés.

\begin{figure}
\begin{center}
  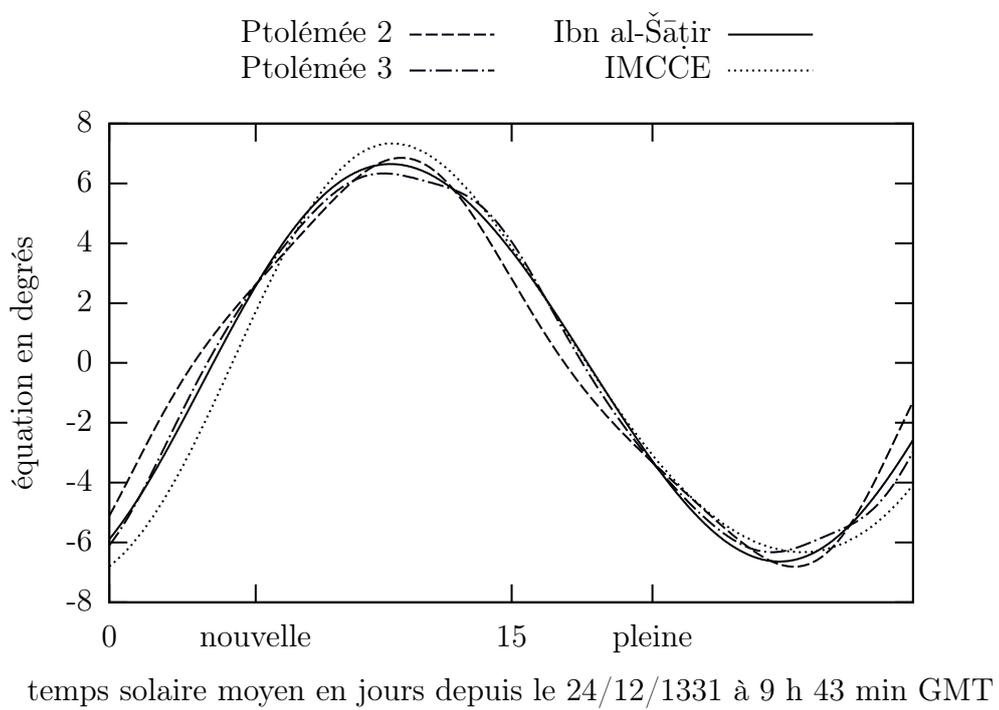
  \caption{\label{fig011}L'équation de la Lune en 1331, sur 30 jours}
\end{center}
\end{figure}

\begin{figure}
  \begin{center}
    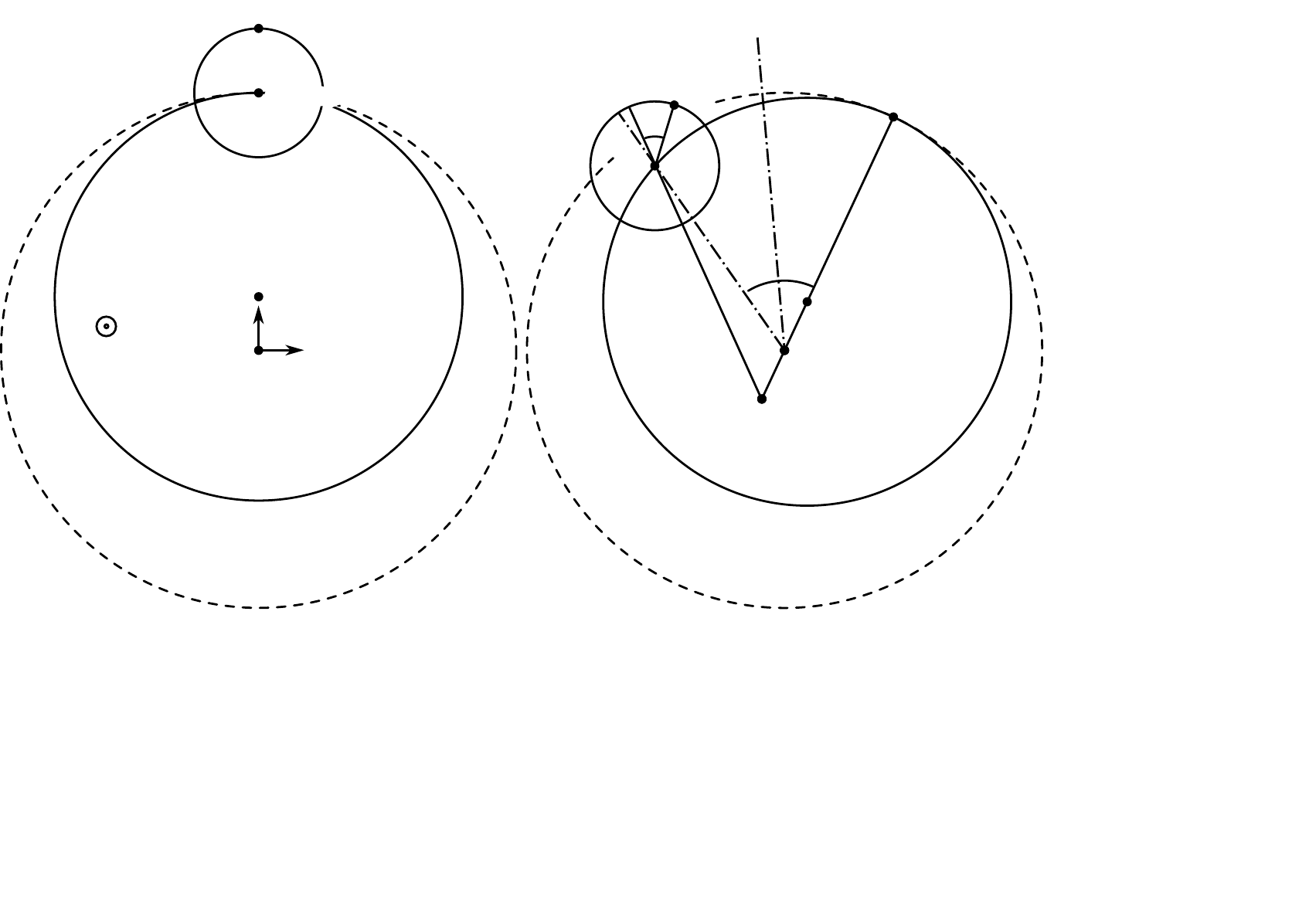
    \caption{\label{fig015}Le troisième modèle de la Lune dans l'\textit{Almageste}}
  \end{center}
\end{figure}

Ainsi pour la Lune, la figure \ref{fig011} compare le modèle
d'{\shatir} et les deuxième et troisième modèles proposés par
Ptolémée.  On a tracé plusieurs graphes. Le premier représente
l'équation de la Lune
$c_2+e_n(\overline{\lambda}+\lambda_{\ascnode}+c_2)$ calculée au moyen
du modèle d'{\shatir} ci-dessus, pendant trente jours à partir du 24
décembre 1331. Sur le même repère, on a aussi représenté l'équation
telle que pouvait l'\emph{observer} {\shatir}. Pour la connaître, à
une date donnée, on calcule la longitude précise de la Lune au moyen
d'outils modernes --~éphémérides de l'IMCCE~-- puis on en retranche la
Lune moyenne
$\overline{\lambda}=\overline{\lambda}(0)+\dot{\overline{\lambda}}t$
calculée avec les paramètres d'{\shatir} pour cette date. Enfin,
toujours sur le même repère, on a représenté l'équation de la Lune
calculée au moyen des deuxième et troisième modèles de Ptolémée dans
l'\textit{Almageste}.  Remarquons que le modèle décrit par Ptolémée
dans les \textit{Hypothèses planétaires} est identique au deuxième
modèle de l'\textit{Almageste}.

\`A des fins de référence, la figure \ref{fig015} rappelle les
caractéristiques du troisième modèle de la Lune dans
l'\textit{Almageste}. Rappelons que le premier modèle de la Lune dans
l'\textit{Almageste} est le modèle ``naïf'' avec un épicycle pour
rendre compte de la première anomalie lunaire. Le deuxième modèle de
la Lune dans l'\textit{Almageste} fait intervenir un excentrique dont
le centre n'est pas fixe. Cet excentrique porte un point, le centre de
l'épicycle, à vitesse angulaire constante par rapport au centre du
monde. Le troisième modèle de la Lune précise que la direction à
partir de laquelle est comptée le mouvement propre de rotation de
l'épicycle, $\overline{\alpha}$ dans la figure \ref{fig015}, passe par
le centre $P_4'$ de l'épicycle et par un point de ``prosneuse'' qui
est distinct du centre $P_3'$ de l'excentrique et du centre du Monde
$O$.

La figure \ref{fig011} montre que tous les modèles donnent des
résultats corrects aux syzygies, \textit{i.e.}, pleine lune et
nouvelle lune, indiquées en abscisse. D'ailleurs, c'est aussi le cas
du premier modèle de Ptolémée --- non représenté ici --- puisqu'il est
déduit des observations aux syzygies.

Le deuxième modèle de Ptolémée corrige le premier en tenant compte
d'observations dans les quadratures, \textit{i.e.}, à mi-chemin entre
les syzygies. Ce faisant, il ouvre le flanc à la première critique
d'{\shatir}, puisqu'il utilise un excentrique. Mais il encourt aussi
son troisième reproche, puisque la direction de l'``apogée'' de
l'épicycle, c'est-à-dire la direction à partir de laquelle est mesuré
le mouvement propre de rotation de l'épicycle, passe par le centre de
l'épicycle et le centre du monde, et non par le centre de l'orbe
excentrique portant le centre de l'épicycle.  De plus, sur la figure
\ref{fig011}, le deuxième modèle de Ptolémée présente un défaut
évident aux octants. 

Dans son troisième modèle, Ptolémée corrige le deuxième modèle en
supposant que la direction de l'``apogée'' de l'épicycle passe par le
centre de l'épicycle et par un autre point, qu'on désignera par la
locution ``point de prosneuse''. Ce point est distinct du centre de
l'excentrique et du centre du monde. Le modèle se prête donc toujours
aux critiques mentionnées ci-dessus, sauf concernant les
octants. Etrangement, p.~\pageref{doutes_lune}, {\shatir} semble
ignorer que Ptolémée a précisément utilisé des observations dans les
octants pour déduire son troisième modèle\footnote{\textit{Cf.}
  \textit{Almageste} V.5, \textit{in} \cite{ptolemy1952}
  p.~149.}. Comment est-ce possible ? Peut-être n'avait-il qu'une
connaissance indirecte de l'\textit{Almageste}, par le biais d'un des
nombreux commentaires qui circulaient à l'époque, peut-être même
seulement par l'intermédiaire de la \textit{Ta\b{d}kirat}
\cite{altusi1993}, qu'il connaissait, et où \d{T}\=us{\=\i} expose le
troisième modèle de Ptolémée en mentionnant des observations dans les
\textit{sextants} (et non les \textit{octants})\footnote{En revanche,
  `Ur\d{d}{\=\i}, dans son \textit{Kit\=ab al-hay'a}
  \cite{saliba1990}, commente les observations de Ptolémée dans les
  octants. Mais {\shatir} a pu négliger cette partie du traité de
  `Ur\d{d}{\=\i} après avoir compris le défaut du modèle lunaire de
  `Ur\d{d}{\=\i} -- \textit{cf.}  notre analyse ci-dessous.}.

\begin{figure}
\begin{center}
  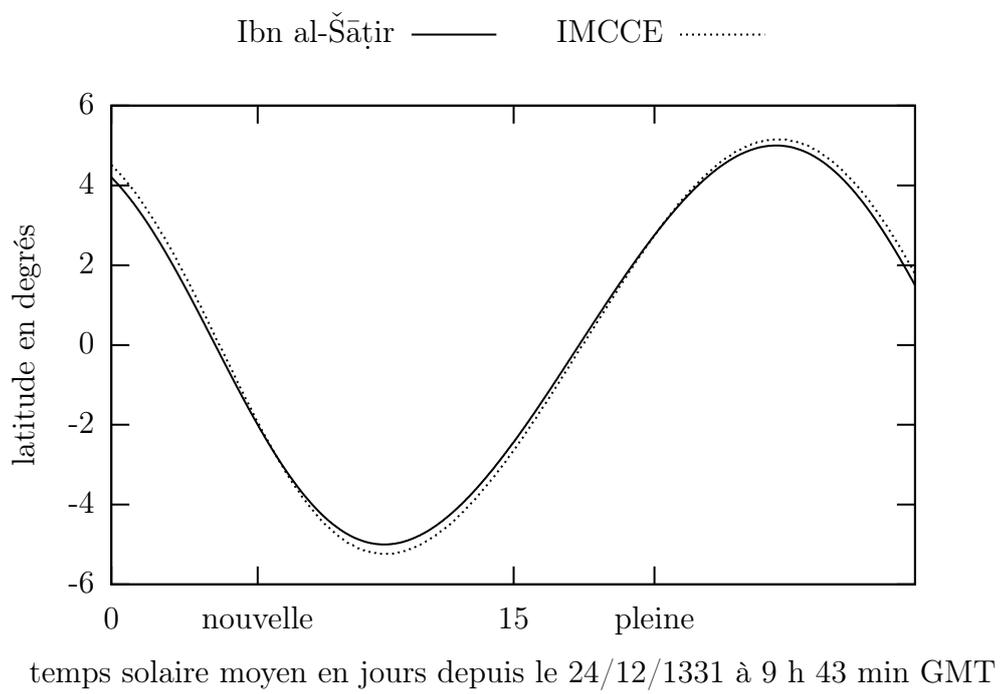
  \caption{\label{fig012}Latitudes de la Lune en 1331, sur 30 jours}
\end{center}
\end{figure}

\begin{figure}
  \begin{center}
    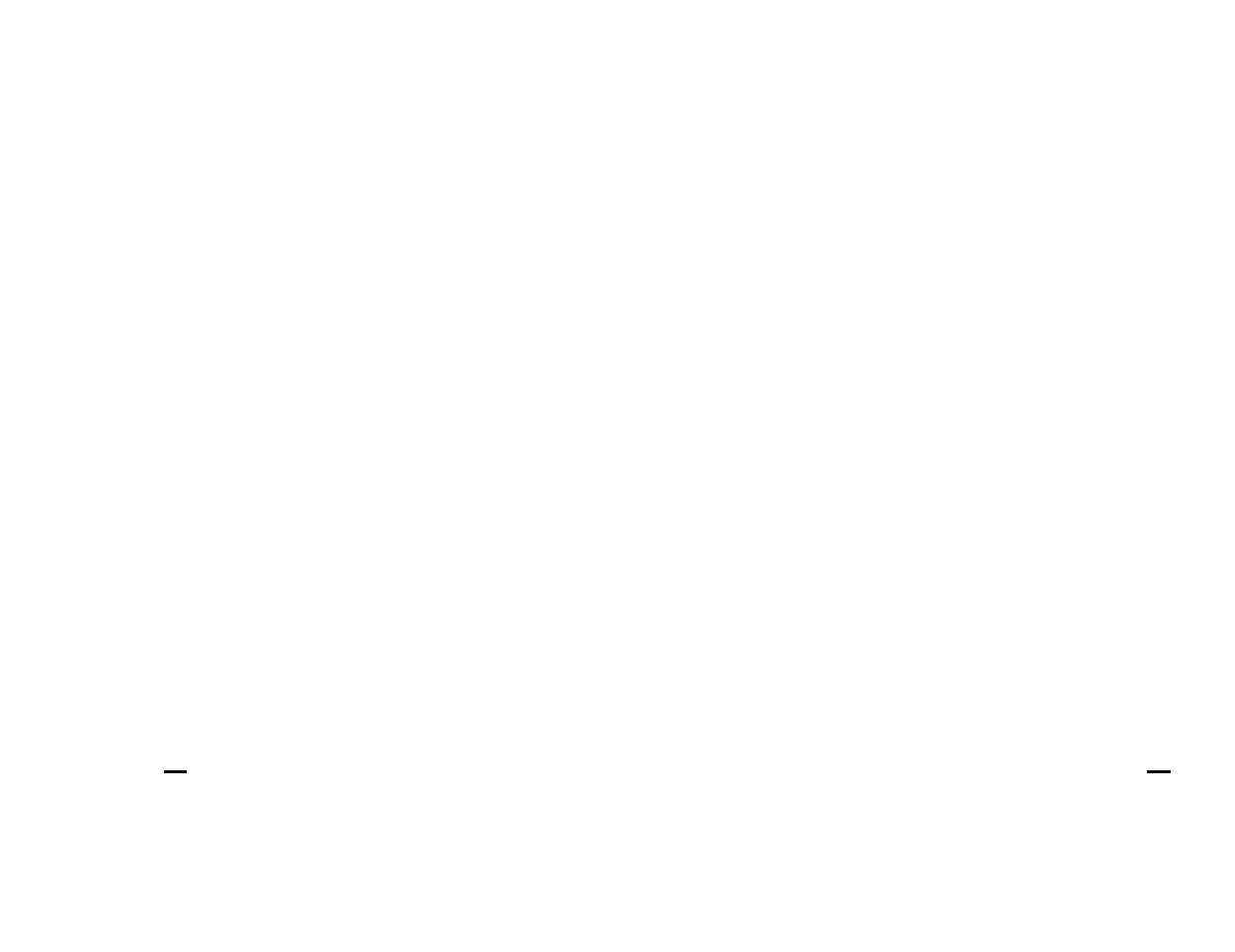
    \caption{\label{fig016}Distance Terre-Lune en 1331, sur 30 jours}
  \end{center}
\end{figure}

On a aussi représenté les latitudes de Lune prédites par le modèle
d'{\shatir} pour les mêmes dates, \textit{cf.}
fig. \ref{fig012}. {\shatir} n'adresse aucune objection à ses
prédécesseurs concernant les latitudes de la Lune, et son modèle
reproduit de très près les prédictions en latitude du troisième modèle
de Ptolémée~: à l'échelle de notre tracé, les deux courbes seraient
indiscernables. Enfin, on a représenté la distance Terre-Lune à la
figure \ref{fig016}, où l'on voit le progrès radical accompli par
{\shatir}. Le troisième modèle de Ptolémée\footnote{Son deuxième
  modèle présente le même défaut.} prédit en effet des valeurs
aberrantes pour la variation de la distance Terre-Lune, et donc aussi
pour la variation du diamètre apparent de la Lune -- il semble varier
du simple au double selon ce modèle. Le modèle d'{\shatir} offre des
prédictions correctes en ordre de grandeur.


\paragraph{La Lune chez al-`Ur\d{d}{\=\i}}
{\shatir} adresse une critique plus circonstanciée à son prédécesseur
Mu'ayyad al-Din al-`Ur\d{d}{\=\i} (?-1266). Il écrit en effet~:
\begin{quote}
  La position conjecturée par al-Mu'ayyad al-`Ur\d{d}{\=\i} dans la
  configuration des orbes de la Lune, concernant son inversion du sens
  du mouvement de l'excentrique et la variation de l'apogée de
  l'épicycle qui suit un point autre que le centre de l'orbe portant
  l'épicycle, est impossible.\footnote{\textit{Cf. supra}
    p.~\pageref{urdi_lune}.}
\end{quote}
Pour une fois, {\shatir} ne critique pas l'usage même d'un
excentrique~; le problème est ailleurs.  On peut penser que le modèle
en question est celui décrit par `Ur\d{d}{\=\i} dans son
\textit{Kit\=ab al-Hay'a}\footnote{\textit{Cf.} l'édition critique du
  texte arabe par Saliba \cite{saliba1990}, accompagnée d'une
  introduction en anglais où Saliba décrit le modèle de
  `Ur\d{d}{\=\i}, p.~50-55.}.

\begin{figure}
  \begin{center}
    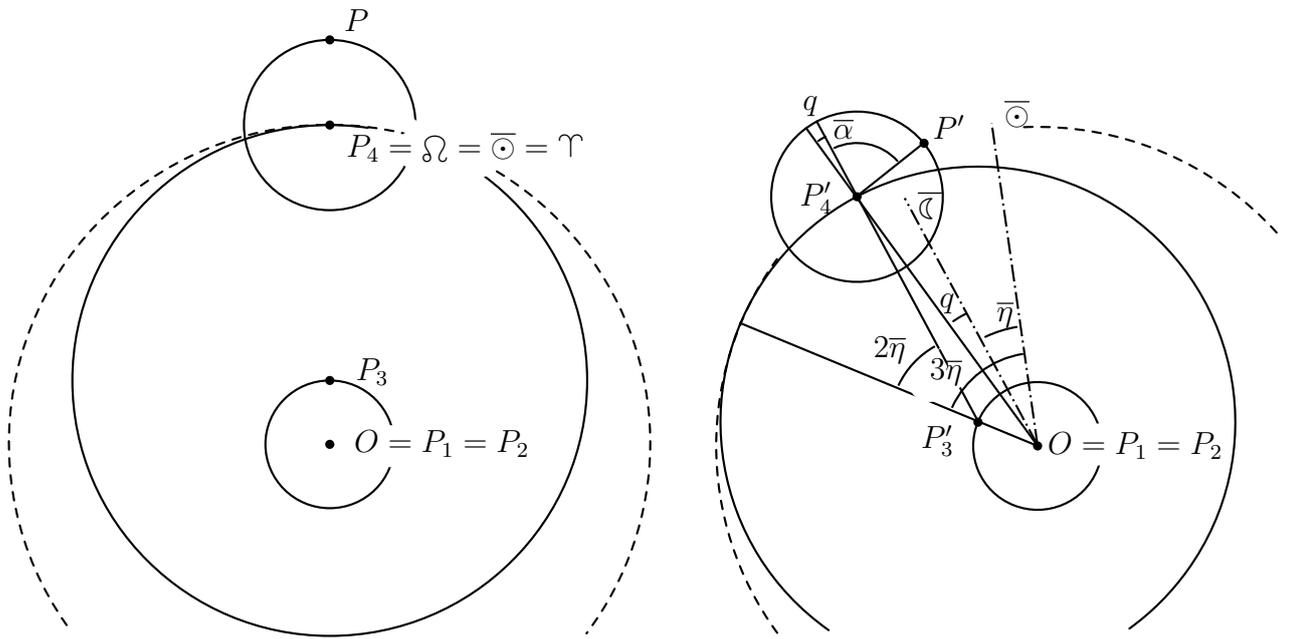
    \caption{\label{fig017}La Lune dans le \textit{Kitab al-Hay'at} de `Ur\d{d}{\=\i}}
  \end{center}
\end{figure}

Suivons la description donnée par al-`Ur\d{d}{\=\i} :
\begin{quote}
  Supposons que l'apogée de la Lune, le centre de son épicycle et le
  Soleil coïncident en un point donné de l'écliptique, à un instant
  donné. Qu'ils se meuvent, chacun selon son mouvement.\footnote{
    \textit{Cf.} \cite{saliba1990} p.~ 118.}
\end{quote}
On reconnaît la figure initiale que nous avons dessiné
fig.~\ref{fig017}, à gauche, où $P_1$ est le centre de l'``orbe des
n{\oe}uds'', $P_2$ le centre de l'orbe incliné, $P_3$ le centre du déférent,
$P_4$ le centre de l'épicycle, et $P$ la Lune. On a~:
$$P_2P_3=10;19\qquad P_3P_4=49;41\qquad P_4P=5;15$$
Al-`Ur\d{d}{\=\i} continue~:
\begin{quote}
  Imaginons l'orbe des n{\oe}uds se mouvoir autour du centre du Monde,
  sur ses pôles et sur son axe qui est l'axe de l'écliptique, en sens contraire
  aux signes, d'un mouvement uniforme égal à l'excès du mouvement en latitude
  sur le mouvement en longitude. Avec soi, il déplace aussi l'orbe incliné et
  tout ce qu'il contient [...]

  Imaginons l'apogée se mouvoir, uniformément, autour de pôles situés à la surface de l'orbe incliné. L'axe de cet orbe rencontre l'axe de l'écliptique au centre du Monde. Son mouvement par jour, dans le sens des signes, est $37;36,39,2$ [$=2\overline{\eta}+\overline{\lambda}+\lambda_{\ascnode}$] [...]

  L'orbe déférent se meut aussi et emporte l'épicycle en sens contraire aux signes, avec un mouvement par jour égal au double de l'élongation. Si l'on soustrait du mouvement de l'apogée (dans le sens des signes) les mouvements du n{\oe}ud et du déférent (en sens contraire), ce qui reste est le mouvement en longitude par jour.\footnote{\textit{Cf.} \cite{saliba1990} p.~118.}
\end{quote}

\begin{figure}
  \begin{center}
    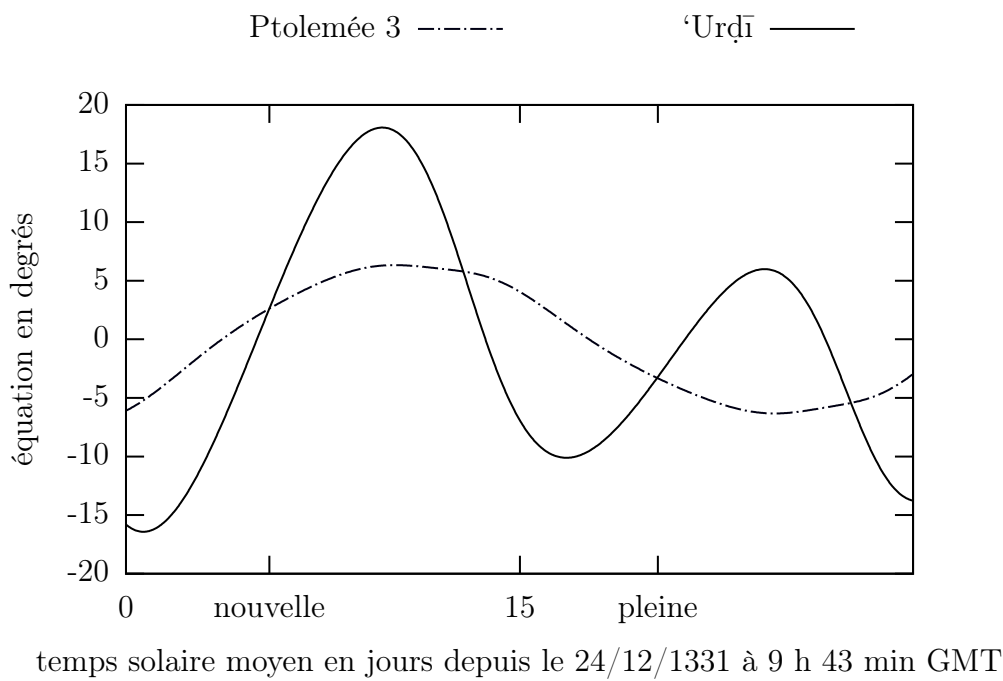
    \caption{\label{fig018}L'équation de la Lune selon al-`Ur\d{d}{\=\i}.}
  \end{center}
\end{figure}

Autrement dit, on passe de la figure initiale à la figure à droite en appliquant une composée de rotations. Si l'on néglige l'inclinaison de l'orbe incliné, on obtient ainsi~:
$$P'=R_{O,3\overline{\eta}+\overline{\lambda}_{\mbox{\tiny\astrosun}}}\,R_{P_3,-2\overline{\eta}}\,R_{P_4,-\overline{\alpha}}(P)$$
et en tenant compte de l'inclinaison de l'orbe incliné :
$$P'=R_{P_1,-\lambda_{\mbox{\tiny\ascnode}}}
\,R_{P_1,\mathbf{j},5°}
\,R_{P_2,\lambda_{\mbox{\tiny\ascnode}}+3\overline{\eta}+\overline{\lambda}_{\mbox{\tiny\astrosun}}}
\,R_{P_3,-2\overline{\eta}}
\,R_{P_4,-\overline{\alpha}}(P)$$
En résolvant des triangles rectangles, il est alors facile de calculer la position de point $P'$. En particulier~:
$$OP_4'=\sqrt{(P_3P_4\sin(2\overline{\eta}))^2+(P_3P_4\cos(2\overline{\eta})+OP_3)^2}$$
$$q=\arcsin\left(\frac{P_3'P_4'\sin(2\overline{\eta})}{OP_4'}\times\frac{OP_3}{P_3'P_4'}\right)=\arcsin\left(\frac{OP_3\sin(2\overline{\eta})}{OP_4'}\right)$$
$$OP'=\sqrt{(P_4P\sin(-(\overline{\alpha}+q)))^2+(OP_4'+P_4P\cos(-(\overline{\alpha}+q)))^2}$$
Si l'équation de la Lune était l'angle $\widehat{P_4'OP'}$, elle serait donc égale à~:
$$\arcsin\left(\dfrac{P_4P\sin(-(\overline{\alpha}+q))}{OP'}\right)$$
Al-`Ur\d{d}{\=\i} s'estimait heureux du fait que cet arc est presque égal à l'équation de la Lune dans le troisième modèle de Ptolémée, et l'avait vérifié, à une bonne précision, en particulier aux octants. Hélas, \emph{$P_4'$ n'est pas dans la direction de la Lune moyenne}. Comme on voit sur la figure, l'équation de la Lune vaut donc~:
$$\arcsin\left(\dfrac{P_4P\sin(-(\overline{\alpha}+q))}{OP'}\right)+q$$
\'Etrangement, comme l'a remarqué Saliba, al-`Ur\d{d}{\=\i} se défend par anticipation contre d'éventuels contradicteurs~:
\begin{quote}
  Qu'on ne dise pas que ces grandeurs ne sont pas prises dans le même cercle~; parce qu'on les entend comme des mouvements moyens, ce sont des arcs de l'écliptique, semblables aux arcs de l'excentrique coupés par le centre de l'épicycle. C'est comme on a coutûme de faire avec le mouvement du Soleil le long de son excentrique, où l'élongation de centre par rapport à l'Apogée est \emph{approximativement} égale à l'élongation par jour.\footnote{Nos italiques. \textit{Cf.} \cite{saliba1990} p. 119, et la remarque de Saliba qui évite cependant de juger ce modèle, p.~54.}
\end{quote}
On a donc le choix entre deux interprétations possibles~:
\begin{itemize}
\item Soit le modèle de `Ur\d{d}{\=\i} n'est pas un modèle cinématique, mais seulement un dispositif de calcul visant à déterminer l'équation du troisième modèle de Ptolémée en n'utilisant que des rotations, en évitant donc les difficultés propres au point de prosneuse. Ce serait alors une sorte d'équatoire.
\item Soit il s'agit d'une erreur de `Ur\d{d}{\=\i}, car l'équation de la Lune dans ce modèle est très loin de reproduire celle du modèle de Ptolémée (\textit{cf.} notre figure \ref{fig018}).
\end{itemize}
Comme `Ur\d{d}{\=\i} exprime par ailleurs clairement que son modèle doit remplacer celui de Ptolémée, la première interprétation nous semble exclue.

{\shatir} avait sûrement vu l'erreur de `Ur\d{d}{\=\i}. \`A la fin du paragraphe énumérant les critiques générales sur la Lune qu'{\shatir} adresse à ses prédécesseurs, celui-ci émet une mise en garde qui fait curieusement écho à la défense par anticipation de `Ur\d{d}{\=\i}~:
\begin{quote}
  La localisation de l'astre moyen et
  du centre de l'astre dans un cercle semblable à l'orbe de l'écliptique
  est impossible car ils ne sont pas pris dans le même cercle, ni par
  rapport au même point.\footnote{\textit{cf. supra} p.~\pageref{astre_moyen}.}
\end{quote}
Cette remarque était probablement adressée à `Ur\d{d}{\=\i} ou à un autre auteur, avant ou après `Ur\d{d}{\=\i}, ayant commis une erreur analogue.

\begin{figure} 
  \begin{center} 
    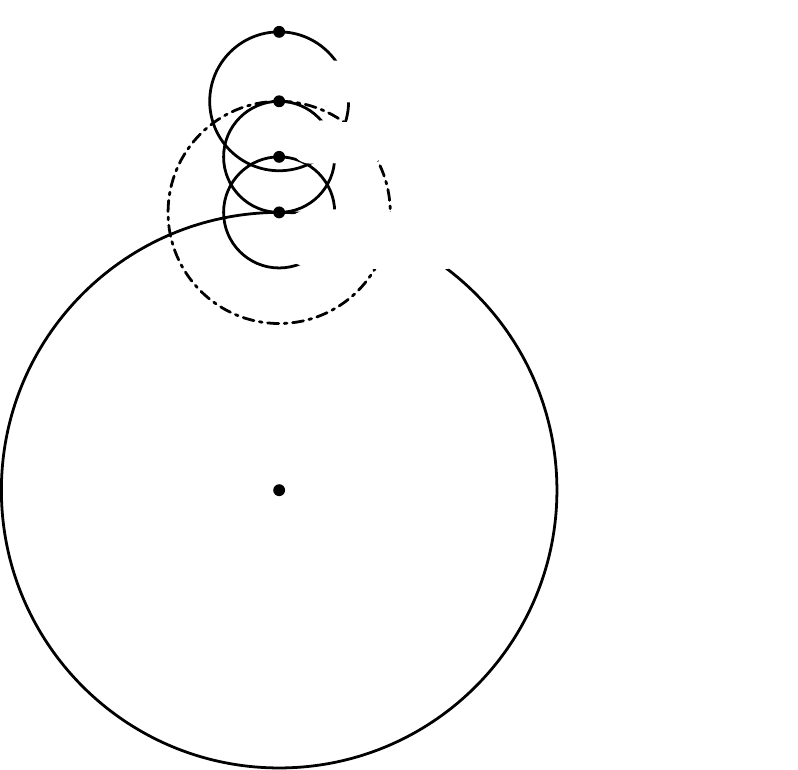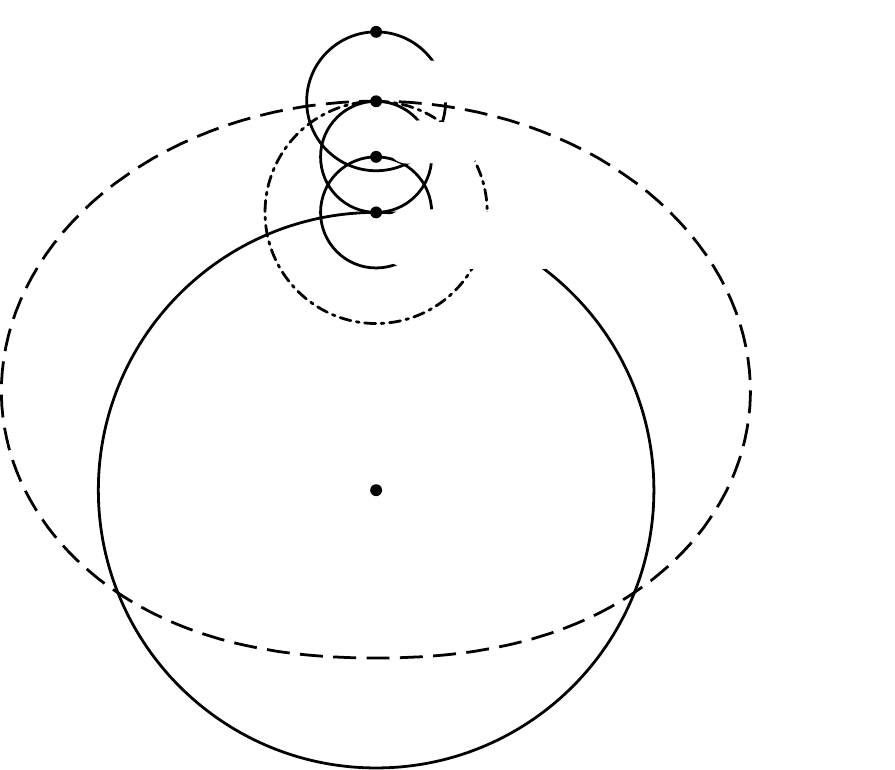
    \caption{\label{fig021}La Lune dans la \textit{Ta\b{d}kirat} de \d{T}\=us{\=\i}}
  \end{center}
\end{figure}

\paragraph{La Lune chez \d{T}\=us{\=\i} et Sh{\=\i}r\=az{\=\i}} Décrivons à présent le modèle exposé par \d{T}\=us{\=\i} dans sa \emph{Ta\b{d}kirat f{\=\i} `ilm al-hay'at}. On a seulement tracé la figure initiale, \textit{cf.} figure \ref{fig021}, à gauche, où $P_1$ est le centre du parécliptique, $P_2$ le centre de l'orbe incliné, $P_3$ le centre d'un orbe déférent, $P_4$ le centre d'une ``grande sphère'', $P_5$ le centre d'une ``petite sphère'', $P_6$ le centre d'une ``sphère englobante'', $P_7$ le centre de l'``épicycle'', et $P$ la Lune\footnote{Comme le dit \d{T}\=us{\=\i}, on pourrait remplacer le couple ``orbe incliné'' et ``orbe déférent'' par un orbe unique puisqu'ils ont même centre et même axe, mais il préfère garder deux orbes distincts. De même pour le couple ``grande sphère'' et ``épicycle''.}. On a~:
$$P_3P_4=49;41\qquad P_4P_5=P_5P_6=\frac{10;19}{2}\qquad P_6P=5;15$$
Le couple ``petite sphère'' et ``grande sphère'' forme un dispositif que les historiens ont appelé ``couple de \d{T}\=us{\=\i}''. Que les orbes tournent, et la position du point $P$ devient~:
$$P'=R_{P_1,-\lambda_{\mbox{\fontencoding {U}\fontfamily {wasy}\fontsize{6pt}{8pt}\selectfont \char 19}}}
\circ R_{P_1,\mathbf{j},5°}
\circ R_{P_2,-(\overline{\eta}-\overline{\lambda}_{\mbox{\tiny\astrosun}}-\lambda_{\mbox{\tiny\ascnode}})}
\circ R_{P_3,2\overline{\eta}}
\circ R_{P_4,2\overline{\eta}}
\circ R_{P_5,-4\overline{\eta}}
\circ R_{P_6,2\overline{\eta}}
\circ R_{P_7,-\overline{\alpha}}(P)$$

Dans le référentiel de l'orbe déférent centré en $P_3$, la trajectoire du point $P_6$ \emph{ressemble} à un cercle de rayon $49;41$. On a représenté cette trajectoire ovale sur la figure \ref{fig021}, en pointillés à droite, en adoptant pour $P_3P_4$, $P_4P_5$ et $P_5P_6$ les longueurs représentées sur la figure. Comme cette figure n'est pas à l'échelle -- le rapport $P_4P_5:P_3P_4$ est bien supérieur à $\dfrac{10;19}{2}:49;41$ pour la rendre plus lisible -- alors notre ovale est beaucoup plus applati que dans le modèle véritable. \d{T}\=us{\=\i} trace et étudie la ressemblance de cette trajectoire avec la ceinture de l'orbe excentrique du modèle de Ptolémée. Dans ses effets, ce modèle reproduit d'assez près le deuxième modèle de Ptolémée.

Afin de copier l'effet de la prosneuse, \d{T}\=us{\=\i} propose d'intercaler un autre couple d'orbes entre la ``sphère englobante'' et l'épicycle. Les historiens ont appelé ce dispositif un ``couple de \d{T}\=us{\=\i} curviligne''. Il s'agit de deux sphères homocentriques, centrées en $P_7$. Leurs rayons n'importent pas à la description du mouvement. Dans la composée de rotations qu'on fait subir au point $P$, \d{T}\=us{\=\i} introduit donc de nouvelles rotations agissant sur $R_{P_7,-\overline{\alpha}}(P)$. Les axes de ces rotations ne sont pas orthogonaux au plan de la figure initiale, et il faut ajouter des rotations pour les rabattre dans ce plan. La position du point $P'$ devient~:
\begin{align*}
  P'=&R_{P_1,-\lambda_{\mbox{\fontencoding {U}\fontfamily {wasy}\fontsize{6pt}{8pt}\selectfont \char 19}}}
  \circ R_{P_1,\mathbf{j},5°}
  \circ R_{P_2,-(\overline{\eta}-\overline{\lambda}_{\mbox{\tiny\astrosun}}-\lambda_{\mbox{\tiny\ascnode}})}
  \circ R_{P_3,2\overline{\eta}}
  \circ R_{P_4,2\overline{\eta}}
  \circ R_{P_5,-4\overline{\eta}}
  \circ R_{P_6,2\overline{\eta}}\\
  &\circ R_{P_7,\mathbf{i},\frac{c_{3\max}}{2}}
  \circ R_{P_7,\mathbf{j},-4\overline{\eta}}
  \circ R_{P_7,\mathbf{i},-\frac{c_{3\max}}{2}}
  \circ R_{P_7,\mathbf{j},2\overline{\eta}}
  \circ R_{P_7,-\overline{\alpha}}(P)
\end{align*}
On expliquera les rotations autour d'axes non orthogonaux au plan de la figure, de manière plus détaillée, quand on commentera l'usage qu'en fait {\shatir} pour les mouvements en latitude des planètes. Ici $c_{3\max}=13;9°$ est l'inclinaison maximale de l'apogée moyen de l'épicycle due à la prosneuse, dans le troisième modèle de Ptolémée.

\begin{figure}
  \begin{center}
    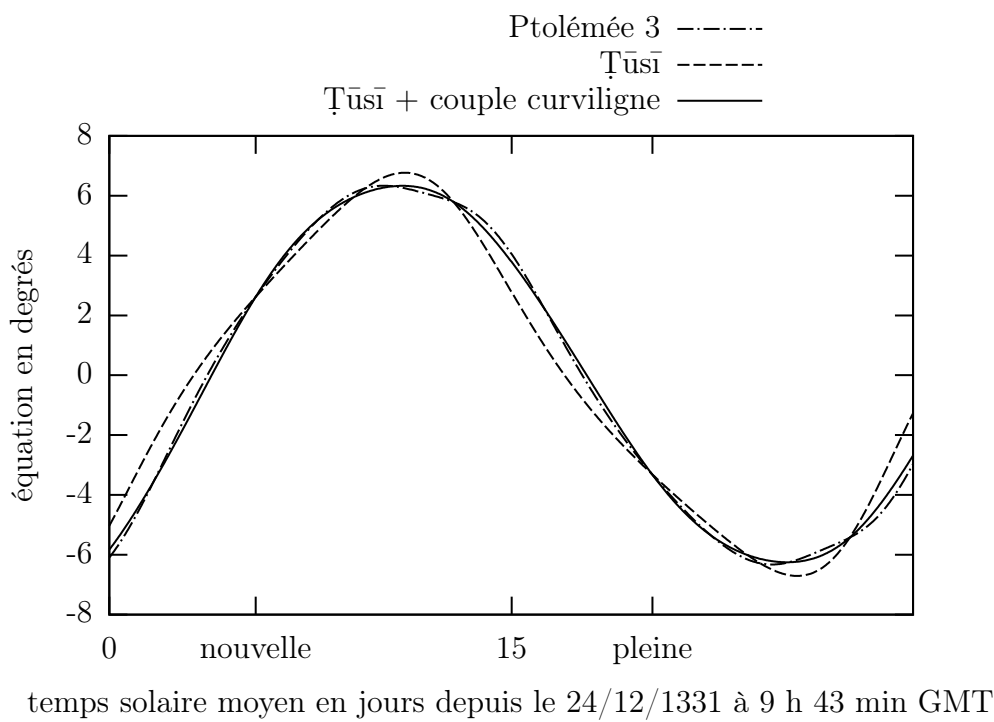
    \caption{\label{fig023}L'équation de la Lune selon \d{T}\=us{\=\i}}
  \end{center}
\end{figure}

La solution de \d{T}\=us{\=\i} est très réussie, en celà qu'elle reproduit de très près la trajectoire prédite par le troisième modèle de Ptolémée, sans plus utiliser ni excentrique, ni point de prosneuse ; voir notre figure \ref{fig023}. La critique que lui adresse {\shatir} peut donc nous surprendre~:
\begin{quote}
  La position mentionnée par Na\d{s}{\=\i}r al-\d{T}\=us{\=\i} dans la
  \emph{Ta\b{d}kira} sur l'élimination des doutes des orbes de
  la Lune est impossible parce qu'on y trouve l'excentrique et que le
  diamètre de l'épicycle suit un point autre que le centre de l'orbe
  portant l'épicycle.\footnote{\textit{Cf. supra} p.~\pageref{doute_tusi_lune}.}
\end{quote}
Nous y reviendrons.

Décrivons à présent le modèle de la Lune proposé par Qu.tb al-D{\=\i}n al-Sh{\=\i}r\=az{\=\i} dans sa \textit{Tu\d{h}fat al-Sh\=ah{\=\i}yat}\footnote{Nous avons consulté le manuscrit de la \emph{Tu\d{h}fat} à la Bibliothèque Nationale (BN~Arabe MS~2516, livre 2, chapitres 10 f.~33r-44v, et 12 f.~50r-56r) en nous aidant de la description faite par Saliba dans \cite{rashed1997}.}
Voici ce que dit {\shatir}~:
\begin{quote}
  La position mentionnée par Qu\d{t}b al-D{\=\i}n
  al-\v{S}{\=\i}r\=az{\=\i} dans la réforme de la configuration des
  orbes de la Lune est impossible car y demeurent l'excentrique et le
  point de prosneuse. L'argument qu'il a avancé (Dieu ait pitié de
  lui) pour permettre le point de prosneuse de la Lune est la plus
  impossible des impossibilités et c'est une image fausse. Il en est
  revenu dans la \emph{Tu\d{h}fa} ; il y indique une autre manière,
  mais elle est aussi impossible.\footnote{\textit{Cf. supra}
    p.~\pageref{shirazi_lune}.}
\end{quote}
Saliba \cite{rashed1997}, suivant peut-être l'avis d'{\shatir}, pense que le modèle décrit dans la \emph{Tu\d{h}fat} est différent des modèles décrits dans deux autres textes de Sh{\=\i}r\=az{\=\i}, les \textit{Ikht{\=\i}y\=ar\=at-i Muzaffar{\=\i}} et la \textit{Nih\=ayat al-Idr\=ak}, qui exposeraient seulement des modèles anciens. Niazi, l'auteur d'une thèse \cite{niazi2011} sur Sh{\=\i}r\=az{\=\i}, a fait une comparaison textuelle des trois ouvrages ; il affirme qu'un modèle semblable au modèle de la \textit{Tu\d{h}fat} est bien présent dans chacun des deux autres ouvrages, mais peut-être pas là où on l'attendrait\footnote{\textit{Cf.} \cite{niazi2011} p.~138-142.}. Ceci nous a aussi été confirmé par Gamini (communication personnelle) qui prépare une édition critique des \textit{Ikht{\=\i}y\=ar\=at}.

\begin{figure}
  \begin{center}
    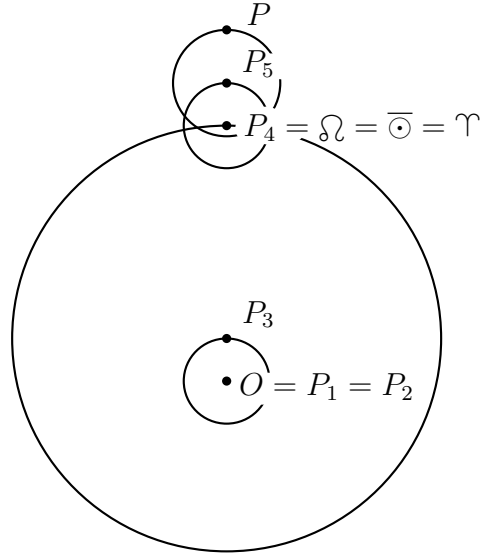
    \caption{\label{fig026}La Lune dans la \textit{Tu\d{h}fat} de Sh{\=\i}r\=az{\=\i}}
  \end{center}
\end{figure}

La figure initiale du modèle de la \textit{Tu\d{h}fat} est notre figure \ref{fig026}. On y voit le centre $P_1$ d'un parécliptique appelé orbe des n{\oe}uds, le centre $P_2$ d'un orbe incliné, le centre $P_3$ d'un orbe excentrique, le centre $P_4$ d'un ``orbe englobant'', le centre $P_5$ d'un épicycle, et la Lune en $P$. Les distances sont~:
$$P_2P_3=\frac{10;19}{2}\qquad P_3P_4=49;41$$
$$P_4P_5=\frac{10;19}{2}\qquad P_5P=5;15$$
La position de la Lune est donnée par la composée de rotations suivante~:
$$P'=R_{P_1,-\lambda_{\mbox{\tiny\ascnode}}}
\,R_{P_1,\mathbf{j},5°}
\,R_{P_2,-(\overline{\eta}-\overline{\lambda}_{\mbox{\tiny\astrosun}}-\lambda_{\mbox{\tiny\ascnode}})}
\,R_{P_3,2\overline{\eta}}
\,R_{P_4,2\overline{\eta}}
\,R_{P_5,-2\overline{\eta}-\overline{\alpha}}(P)$$
Ce modèle est inspiré des modèles de `Ur\d{d}{\=\i} pour les planètes supérieures, mais il diffère substantiellement de son modèle pour la Lune. Au contraire, il est strictement équivalent au modèle de \d{T}\=us{\=\i} sans couple curviligne~: les trajectoires prédites sont égales.

\emph{Démonstration}. Notons $P_i$ ($1\leq i\leq 5$) les centres des orbes de Sh{\=\i}r\=az{\=\i}, et $\tilde{P}_i$ ($1\leq i\leq 7$) les centres des orbes de \d{T}\=us{\=\i}, avec $P_1=P_2=\tilde{P}_1=\tilde{P}_2=\tilde{P}_3$ et $P_4=\tilde{P}_5$. On a alors~:
$$\overrightarrow{P_2P_3}=\overrightarrow{\tilde{P}_4\tilde{P}_5}=\frac{10;19}{2}\mathbf{j}$$
$$\overrightarrow{P_3P_4}=\overrightarrow{\tilde{P}_3\tilde{P}_4}=49;41\mathbf{j}$$
$$P_5=\tilde{P}_6=\tilde{P}_7,\quad\overrightarrow{P_3\tilde{P}_3}=\overrightarrow{P_4\tilde{P}_4}=-\overrightarrow{\tilde{P}_5\tilde{P}_6}$$
Cette dernière égalité permet d'appliquer la proposition 2, et on a alors~:
\begin{align*}
&R_{P_2,-(\overline{\eta}-\overline{\lambda}_{\mbox{\tiny\astrosun}}-\lambda_{\mbox{\tiny\ascnode}})}
\,R_{P_3,2\overline{\eta}}
\,R_{P_4,2\overline{\eta}}
\,R_{P_5,-2\overline{\eta}-\overline{\alpha}}\\
=&R_{\tilde{P}_2,-(\overline{\eta}-\overline{\lambda}_{\mbox{\tiny\astrosun}}-\lambda_{\mbox{\tiny\ascnode}})}
\,\left(
R_{P_3,2\overline{\eta}}
\,R_{P_4,2\overline{\eta}}
\,R_{\tilde{P}_6,-2\overline{\eta}}
\right)
\,R_{\tilde{P}_7,-\overline{\alpha}}\\
=&R_{\tilde{P}_2,-(\overline{\eta}-\overline{\lambda}_{\mbox{\tiny\astrosun}}-\lambda_{\mbox{\tiny\ascnode}})}
\,R_{\tilde{P}_3,2\overline{\eta}}
\,\left(
R_{\tilde{P}_4,-2\overline{\eta}}
\,R_{\tilde{P}_5,4\overline{\eta}}
\,R_{\tilde{P}_6,-2\overline{\eta}}
\right)
\,R_{\tilde{P}_7,-\overline{\alpha}}
\end{align*}
En vertu de la proposition 3, ceci est aussi égal à~:
$$R_{\tilde{P}_2,-(\overline{\eta}-\overline{\lambda}_{\mbox\tiny\astrosun}-\lambda_{\mbox{\tiny\ascnode}})}
\,R_{\tilde{P}_3,2\overline{\eta}}
\,\left(
R_{\tilde{P}_4,2\overline{\eta}}
\,R_{\tilde{P}_5,-4\overline{\eta}}
\,R_{\tilde{P}_6,2\overline{\eta}}
\right)
\,R_{\tilde{P}_7,-\overline{\alpha}},\quad\text{\textit{q. e. d.}}$$

\begin{figure}
  \begin{center}
    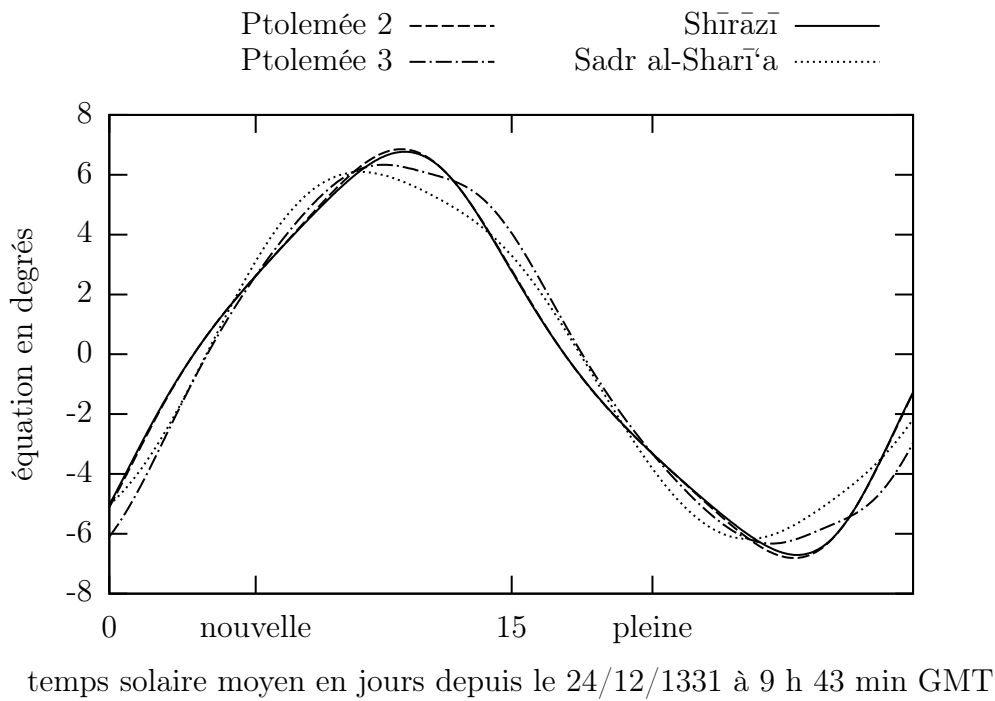
    \caption{\label{fig028}La Lune selon Sh{\=\i}r\=az{\=\i} et Sadr al-Shar{\=\i}`a}
  \end{center}
\end{figure}

Ce modèle produit des trajectoires très proches du deuxième modèle de Ptolémée, et il ne reproduit donc pas les effets du point de prosneuse.\footnote{Ce modèle achevé décrit par Sh{\=\i}r\=az{\=\i} ne tient donc pas compte de la troisième anomalie lunaire. Saliba affirme que, même dans la \textit{Tu\d{h}fat}, ``Qu\d{t}b al-D{\=\i}n ne semble pas avoir réussi à répondre à la seconde objection'' (celle du point de prosneuse), \textit{cf.} \cite{rashed1997},~p.~113. Il resterait cependant à identifier clairement, par une étude complète des trois ouvrages, ``l'argument que [Sh{\=\i}r\=az{\=\i}] a avancé pour permettre le point de prosneuse'' selon {\shatir}.} Dallal et Saliba ont indiqué qu'un contemporain d'{\shatir}, Sadr al-Shar{\=\i}`a (?-1347), a proposé un amendement au modèle de Sh{\=\i}r\=az{\=\i} pour reproduire les effets du point de prosneuse, en ajoutant une petite sphère de rayon $\overrightarrow{P_6P}=0;52\,\mathbf{j}$ et de centre $P_6$ remplaçant le point $P$ dans la figure initiale \ref{fig026}~; il faut alors ajouter une rotation à droite\footnote{\textit{cf.} \cite{dallal1995} p.~373-379 et \cite{rashed1997} p.~112-113. Dans \cite{rashed1997}, le sens de rotation de la petite sphère semble erroné. Dans \cite{dallal1995}, Dallal donne deux interprétations possibles du texte de Sadr al-Shar{\=\i}`a : avec nos notations, soit $\overrightarrow{P_6P}=+0;52\,\mathbf{j}$ et l'angle de rotation est $-2\overline{\eta}$, soit $\overrightarrow{P_6P}=-0;52\,\mathbf{j}$ et l'angle de rotation est $+2\overline{\eta}$. La deuxième interprétation, bien qu'elle ressemble qualitativement au modèle d'{\shatir}, produit un résultat numériquement peu satisfaisant. C'est la première interprétation que nous avons donc retenue dans notre comparaison.}~:
$$P'=R_{P_1,-\lambda_{\mbox{\tiny\ascnode}}}
\,R_{P_1,\mathbf{j},5°}
\,R_{P_2,-(\overline{\eta}-\overline{\lambda}_{\mbox{\tiny\astrosun}}-\lambda_{\mbox{\tiny\ascnode}})}
\,R_{P_3,2\overline{\eta}}
\,R_{P_4,2\overline{\eta}}
\,R_{P_5,-2\overline{\eta}-\overline{\alpha}}
\,R_{P_6,-2\overline{\eta}}(P).$$
Le modèle de Sadr al-Shar{\=\i}`a a bien l'effet escompté, comme on peut le voir dans les octants sur notre figure \ref{fig028}.

\paragraph{Comment expliquer le reproche qu'{\shatir} adresse à \d{T}\=us{\=\i} concernant la Lune~?} Tout d'abord, \d{T}\=us{\=\i} rejette lui-même son propre modèle, à la fin du chapitre II.11 de la \textit{Ta\b{d}kirat}, par une intéressante analyse du taux de variation de l'anomalie due à la prosneuse chez Ptolémée et au couple curviligne chez lui : il montre que son modèle n'a pas exactement les mêmes effets que le troisième modèle de Ptolémée. Est-il donc possible que la critique d'{\shatir} ne s'adresse pas au modèle présenté dans II.11, mais plutôt à la description du troisième modèle de Ptolémée faite par \d{T}\=us{\=\i} dans II.7~? Mais c'est à \d{T}\=us{\=\i} qu'{\shatir} adresse un reproche ; n'aurait-il pas compris qu'il s'agissait là d'un modèle dû à Ptolémée\footnote{ceci pourrait aussi expliquer qu'il ignorait les observations dues à Ptolémée dans les octants, comme on l'a vu plus haut.}~?

Il nous semble plus probable qu'{\shatir} ajoute à l'auto-critique de \d{T}\=us{\=\i} sur la prosneuse un second reproche plus général concernant le modèle présent dans II.11. Mais quel reproche~?

Bien que \d{T}\=us{\=\i} n'utilise pas d'orbe excentrique, il conçoit une combinaison de mouvements épicycliques afin de donner au point $P_6$ une trajectoire \emph{proche du cercle excentrique de Ptolémée dans le plan incliné} ; il prend d'ailleurs soin de le démontrer dans II.11. \d{T}\=us{\=\i} substitue donc à l'excentrique de Ptolémée un dispositif reproduisant la même trajectoire, par le truchement d'un artifice, le couple de \d{T}\=us{\=\i}, produisant un mouvement rectiligne oscillatoire. En tant qu'artifice produisant un mouvement \textit{par essence} interdit dans les cieux, ce dispositif ne devrait-il pas être rejeté, au même titre qu'un point glissant à la circonférence d'un excentrique~? Hélas, {\shatir} lui-même utilise le couple de \d{T}\=us{\=\i} dans certains modèles planétaires... Donc l'argument ne tient pas.

\begin{figure}
  \begin{center}
    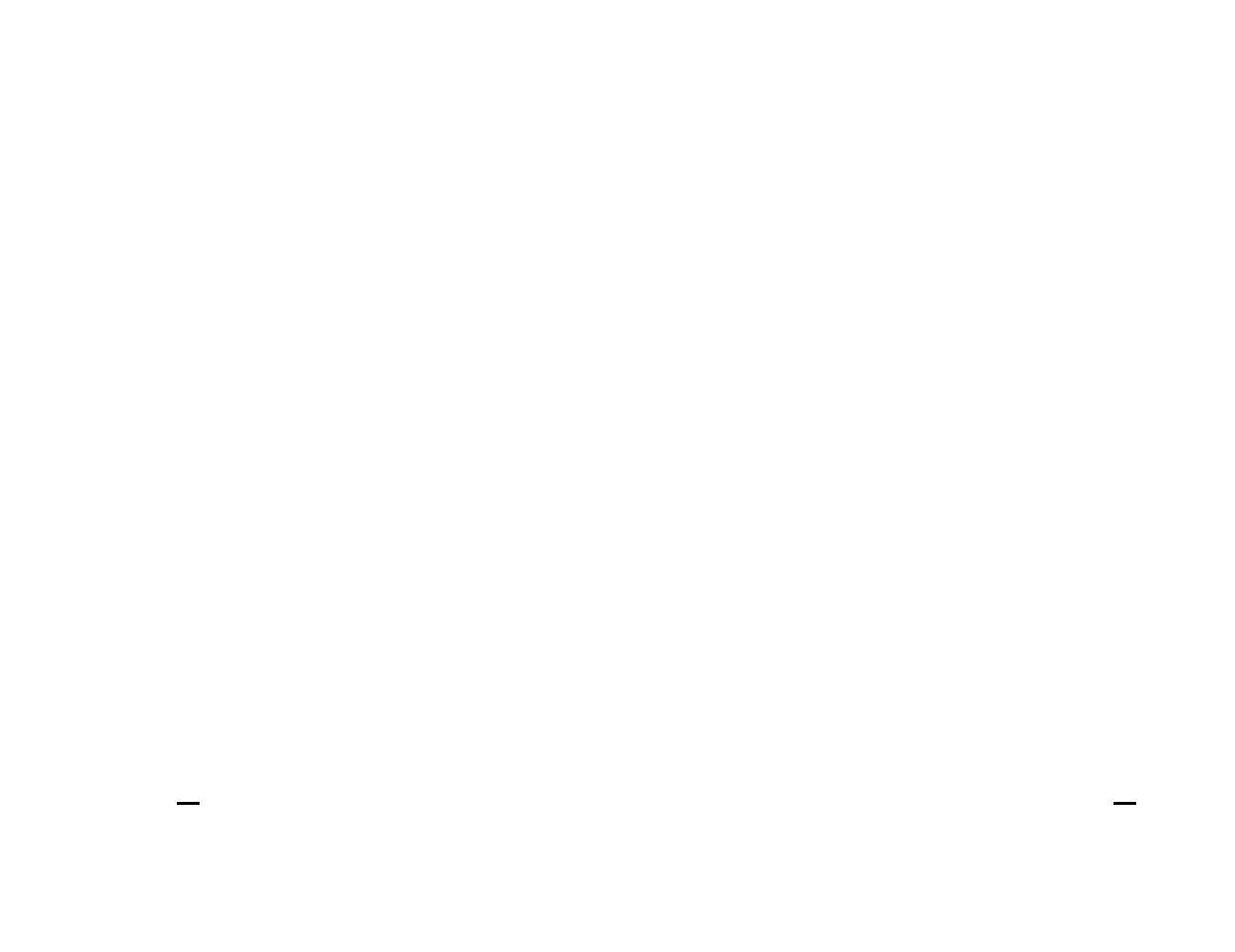
    \caption{\label{fig027}L'équation de la Lune~: meilleurs modèles jusqu'à {\shatir}}
  \end{center}
\end{figure}

\begin{figure}
  \begin{center}
    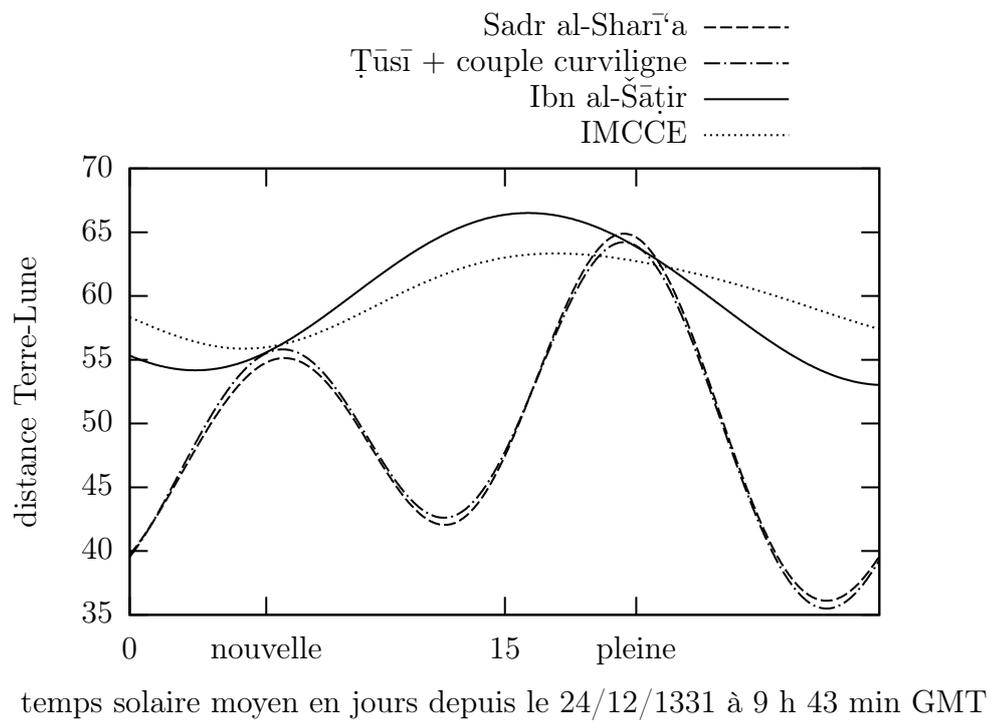
    \caption{\label{fig029}La distance Terre-Lune dans les modèles les meilleurs pour les longitudes}
  \end{center}
\end{figure}

Voici donc la seule explication probable. Le refus des excentriques par {\shatir} n'est pas une objection de principe concernant seulement l'essence physique de l'excentrique ou les mouvements permis par l'astronomie physique. C'est un choix méthodologique dans l'analyse des données de l'observation. Ptolémée aurait commis une erreur dès son \textit{deuxième} modèle lunaire, en changeant la configuration globale des orbes du premier et en intercalant un \textit{grand} excentrique pour rendre compte d'une variation de l'amplitude maximale d'une \textit{petite} anomalie lunaire aux quadratures. Et en effet, l'utilisation d'un excentrique a deux effets joints :
\begin{itemize}
\item produire cette variation d'amplitude de la petite anomalie lunaire
\item induire un changement global dans la configuration des orbes, et en particulier dans la distance Terre-Lune.
\end{itemize}
Seul le premier effet était désirable. {\shatir} concevait certainement ses modèles en suivant une saine méthodologie : partir d'une estimation préalable des mouvements moyens et des distances, pour ensuite ajouter de \textit{petites} corrections en adjoignant des épicycles \textit{de plus en plus petits} perturbant peu les autres caractéristiques des modèles préalables. Il doit avoir expliqué cette méthode dans son ouvrage perdu sur les observations, \textit{Ta`l{\=\i}q al-'ar\d{s}\=ad}. Pour la Lune, la faible variation de la distance Terre-Lune doit compter parmi les faits primitifs que les corrections successives ne peuvent remettre en cause. Nos figures \ref{fig027} et \ref{fig029} comparent entre eux les <<~meilleurs~>> modèles lunaires à partir du troisième modèle de Ptolémée~: celui de Qu\d{t}b al-D{\=\i}n al-Sh\={\i}r\=az{\=\i}, celui d'al-\d{T}\=us{\=\i} avec le couple curviligne, et celui d'{\shatir}. Pour ne pas alourdir ces figures, nous n'avons pas représenté le troisième modèle de Ptolémée~: on a vu que celui d'al-\d{T}\=us{\=\i} en est très proche. \`A partir du deuxième modèle de Ptolémée, tous les modèles lunaires des prédécesseurs d'{\shatir} produisaient des valeurs aberrantes de la distance Terre-Lune, y compris celui de \d{T}\=us{\=\i}, alors même que celui-ci s'était débarrassé de l'excentrique de Ptolémée vis-à-vis de l'astronomie physique. C'est qu'il fallait en fait repartir du \textit{premier} modèle de Ptolémée, et tout reconstruire à neuf.

\paragraph{Calibrage du modèle lunaire d'{\shatir}} Le chemin qu'il a pu suivre pour déterminer les paramètres $P_2P_3$, $P_3P_4$ et $P_4P$ est facile à reconstituer~:
\begin{enumerate}[1)]
\item On utilise trois éclipses pour calculer l'équation maximale $c$ aux syzygies comme l'avait fait Ptolémée. Cela donne
  $$\frac{P_3P_4-P_4P}{P_2P_3}.$$
  {\shatir} utilise $c_2=4;56$ comme équation maximale, valeur connue depuis longtemps.
\item On utilise des observations précises donnant une équation maximale maximale dans les quadratures égale à $c_2=7;40$, déjà connue de Ptolémée. On en déduit
  $$\frac{P_3P_4+P_4P}{P_2P_3}$$
\item Puis on vérifie que le modèle est satisfaisant dans les octants. Pour cela, {\shatir} a dû faire ses propres observations puisqu'il n'avait pas conscience de celles émises par Ptolémée. Enfin, on utilise le modèle pour calculer la distance Terre-Lune et on compare avec des observations anciennes ou récentes qu'{\shatir} semblait bien connaître\footnote{\textit{cf. supra} p.~\pageref{obs_diam_lune}.}.
\end{enumerate}

\begin{figure}
  \begin{center}
    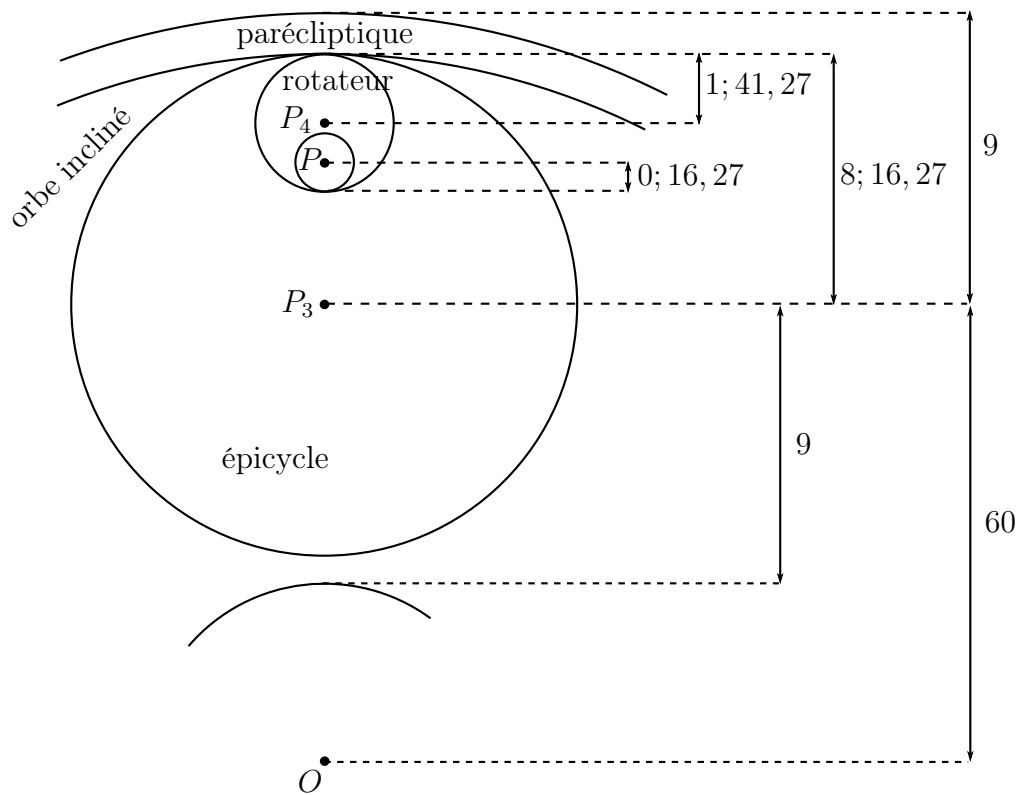
    \caption{\label{fig014}Les orbes solides de la Lune}
  \end{center}
\end{figure}

\begin{figure}
  \begin{center}
    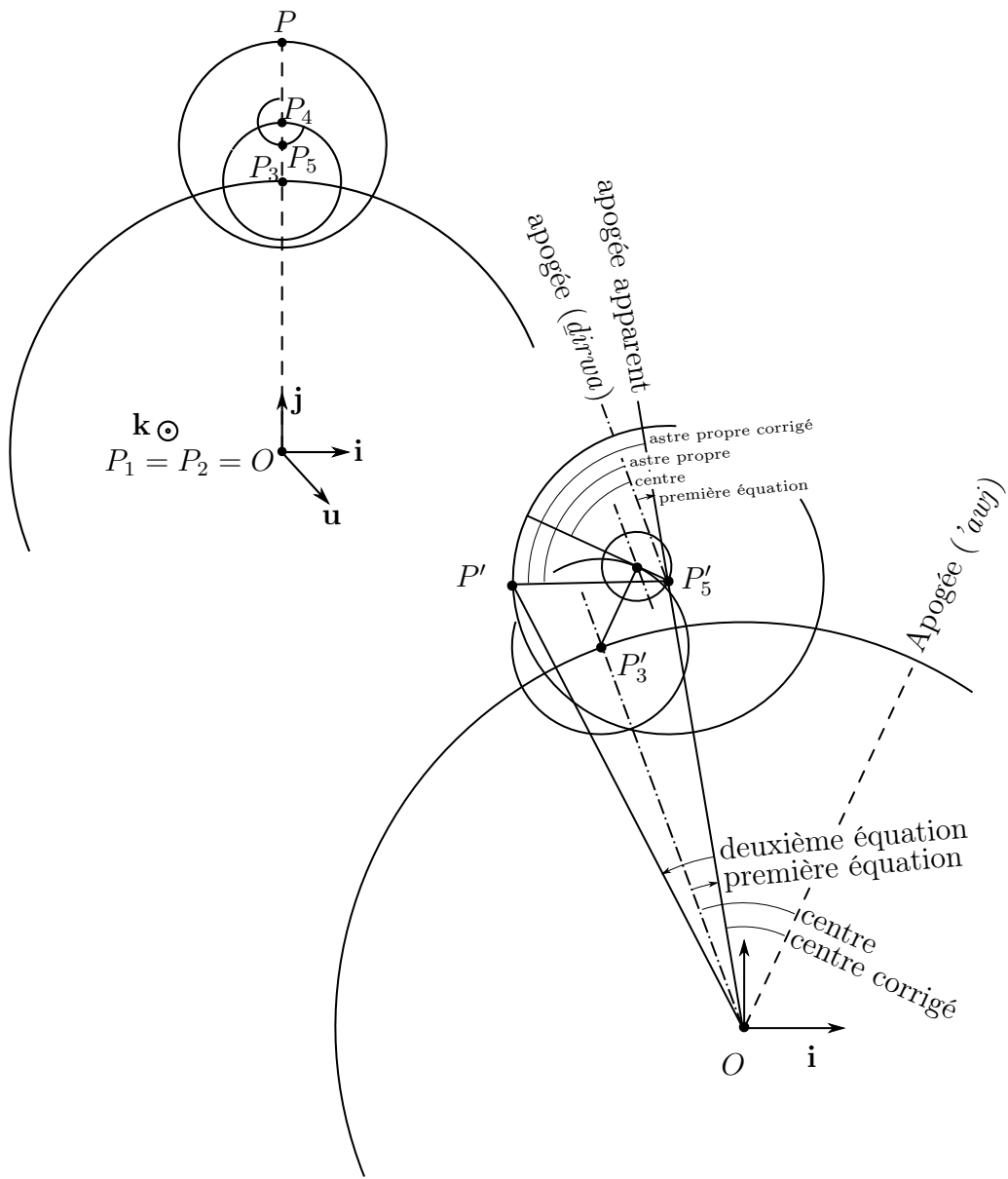
  \end{center}
  \caption{\label{fig040}\label{fig_term}Les orbes de Saturne : figure initiale et terminologie}
\end{figure}

\paragraph{Saturne : figure initiale}\label{begin_saturne}
Dans la figure initiale, les plans des orbes de Saturne sont tous
rabattus dans le plan de l'écliptique $(\mathbf{i},\mathbf{j})$ par
des rotations, les centres des orbes sont tous alignés dans la
direction du vecteur $\mathbf{j}$ elle-même confondue avec la
direction du point vernal~; enfin, la ligne des n{\oe}uds, à l'intersection
des plans de l'orbe incliné et du parécliptique, est orientée dans la
direction du vecteur $\mathbf{u}=\cos(50°)\mathbf{i}-\sin(50°)\mathbf{j}$ qui
pointe vers le n{\oe}ud ascendant. Le point $O$ est le centre du Monde, $P_1$
le centre du parécliptique de Saturne, et $P_2$ le centre de son orbe
incliné. Ces trois points sont confondus $O=P_1=P_2$. Le point $P_3$ est
le centre d'un orbe \emph{déférent} (\textit{falak \d{h}\=amil}), $P_4$ le
centre de l'orbe rotateur, $P_5$ le centre de l'orbe de l'épicycle, et $P$
est le centre du corps sphérique de la planète.
La figure \ref{fig040} précise la position de ces points, mais elle n'est
pas à l'échelle. Posons~: $\overrightarrow{OP_3}=60\ \mathbf{j}$, alors
$$\begin{array}{l}
\overrightarrow{P_3P_4}=5;7,30\ \mathbf{j}\\
\overrightarrow{P_4P_5}=-1;42,30\ \mathbf{j}\\
\overrightarrow{P_5P}=6;30\ \mathbf{j}
\end{array}$$
Il faut appliquer au point $P$ une suite de transformations géométriques
pour obtenir la position de Saturne prédite ou observée à un instant
donné.

\paragraph{Saturne~: transformations géométriques}
Les transformations appliquées à Saturne $P$ sont des rotations paramétrées par trois angles~: l'\emph{Apogée} $\lambda_A$, le \emph{centre moyen} $\overline{\kappa}$, et \emph{l'astre propre} $\overline{\alpha}$. L'\emph{astre moyen} est $\overline{\lambda}=\lambda_A+\overline{\kappa}$, le \emph{n{\oe}ud ascendant} est $\lambda_{\ascnode}=\lambda_A-140°$, et l'\emph{argument de latitude} est $\overline{\lambda}-\lambda_{\ascnode}$. Voici la liste des transformations géométriques utilisées pour produire les mouvements en longitude~:
$$R_{P_5,\overline{\alpha}-\overline{\kappa}},\quad R_{P_4,2\overline{\kappa}},\quad
R_{P_3,-\overline{\kappa}},\quad R_{P_2,\overline{\kappa}}, \quad R_{P_1,\lambda_A}.$$
Comme d'habitude, quand on omet la mention de l'axe de rotation, il s'agit par défaut de la direction vecteur $\mathbf{k}$. Mais il faut se souvenir que la figure \ref{fig2} p.~\pageref{fig2} cas (i) représente pour Saturne la relation entre l'orbe incliné et l'orbe déférent, ainsi que la relation entre l'orbe déférent et l'orbe rotateur. Il faut donc adjoindre des rotations visant à incliner les plans des orbes par rapport aux plans des orbes qui les portent~:
$$R_{P_4,\mathbf{v},-1°},\quad R_{P_3,\mathbf{u},-3°30'},\quad R_{P_2,\mathbf{u},2°30'}.$$
où $\mathbf{v}$ est entraîné par l'orbe rotateur de sorte à ce que sa direction coïncide avec la ligne des n{\oe}uds lorsque $\overline{\kappa}=(\mathbf{j},\mathbf{u})+90°$. Autrement dit~:
$$R_{2((\mathbf{j},\mathbf{u})+90°)}(\mathbf{v})=\mathbf{u}.$$
Pour Saturne, on en déduit que~:
$$\mathbf{v}=R_{2\times(140°-90°)}(\mathbf{u})=\cos(50°)\mathbf{i}+\sin(50°)\mathbf{j}.$$
La théorie des longitudes des chapitres 12 à 14 et la description donnée par {\shatir} dans le chapitre 24 sur les latitudes des planètes supérieures ne laissent aucun doute quant à l'ordre dans lequel appliquer toutes ces transformations.

\begin{figure}
  \begin{center}
    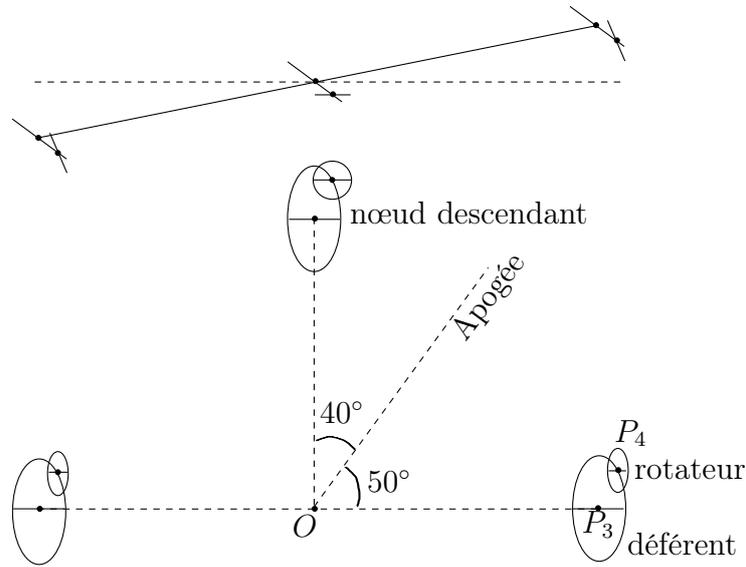
  \end{center}
  \caption{\label{fig041}Saturne, trois positions : projection orthogonale}
\end{figure}

Pour obtenir la configuration des orbes à un instant donné, il faut donc appliquer toutes ces rotations aux points $P_3$, $P_4$, $P_5$, $P$. L'image du point $P_3$ entraîné par les mouvements de l'orbe parécliptique et de l'orbe incliné est~:
$$R_{P_1,\lambda_A}\,R_{P_2,\mathbf{u},2°30'}\,R_{P_2,\overline{\kappa}}(P_3)$$
Quant au point $P_4$, il est aussi entraîné par le mouvement de l'orbe déférent et il devient~:
$$R_{P_1,\lambda_A}\, R_{P_2,\mathbf{u},2°30'}\, R_{P_2,\overline{\kappa}}\, R_{P_3,-\overline{\kappa}}\, R_{P_3,\mathbf{u},-3°30'}(P_4)$$
Le point $P_5$ est aussi entraîné par le mouvement de l'orbe rotateur et il devient~:
$$R_{P_1,\lambda_A}\, R_{P_2,\mathbf{u},2°30'}\, R_{P_2,\overline{\kappa}}\, R_{P_3,-\overline{\kappa}}\, R_{P_3,\mathbf{u},-3°30'}\, R_{P_4,2\overline{\kappa}}\, R_{P_4,\mathbf{v},-1°}(P_5)$$
Enfin le point $P$ est aussi entraîné par l'orbe de l'épicycle dont le plan est confondu avec le plan de l'orbe rotateur, et il devient~:
$$R_{P_1,\lambda_A}\, R_{P_2,\mathbf{u},2°30'}\, R_{P_2,\overline{\kappa}}\, R_{P_3,-\overline{\kappa}}\, R_{P_3,\mathbf{u},-3°30'}\, R_{P_4,2\overline{\kappa}}\, R_{P_4,\mathbf{v},-1°}\, R_{P_5,\overline{\alpha}-\overline{\kappa}}(P)$$
La figure \ref{fig041} montre l'effet des rotations inclinant les plans des orbes, pour trois positions. On a projeté orthogonalement sur le plan du parécliptique (figure du bas) et sur un plan orthogonal à la ligne des n{\oe}uds (figure du haut). On a seulement représenté les ceintures du déférent et du rotateur puisque la ceinture de l'épicycle est dans le plan de la ceinture du rotateur.

\paragraph{Saturne : trajectoire paramétrée}
Le modèle de Saturne est couplé au modèle du Soleil en vertu de la relation suivante~:
$$\overline{\alpha}=\overline{\lambda}_{\astrosun}-\lambda_A-\overline{\kappa}$$
La trajectoire dans l'ensemble des valeurs des paramètres est~:
$$\lambda_A=\dot{\lambda}_At+\lambda_A(0)$$
$$\overline{\kappa}=\dot{\overline{\kappa}}t+\overline{\kappa}(0)$$
$$\overline{\lambda}_{\astrosun}=\dot{\overline{\lambda}}_{\astrosun}t+\overline{\lambda}_{\astrosun}(0)$$
On rappelle que $\overline{\kappa}=\overline{\lambda}-\lambda_A$. On a, d'après {\shatir}~:
$$\overline{\lambda}_{\astrosun}(0)=280;9,0,\quad\dot{\overline{\lambda}}_{\astrosun}=359;45,40\text{ par année persane},$$
$$\overline{\lambda}(0)=157;58,20,\quad\dot{\overline{\lambda}}=12;13,40\text{ par année persane},$$
$$\lambda_A(0)=254;52,\quad\dot{\lambda}_A=0;1\text{ par année persane}.$$

\paragraph{Transformations planes}
Si $R$ et $S$ sont deux rotations, il existe une rotation $T$ telle que $R\circ S=T\circ R$, à savoir, la rotation de même angle que $S$ et dont l'axe est l'image par $R$ de l'axe de $S$. Appliquons cette relation de commutation à toutes les composées de rotations décrites ci-dessus, de sorte à réécrire à droite toutes les rotations dont l'axe est dans la direction du vecteur $\mathbf{k}$. Par exemple~:
$$R_{P_1,\lambda_A}R_{P_2,\mathbf{u},2°30'}=R_{P_2,\mathbf{a},2°30'}R_{P_1,\lambda_A}$$
où le vecteur $\mathbf{a}$ est l'image de $\mathbf{u}$ par $R_{P_1,\lambda_A}$. On fait de même avec les inclinaisons des plans des petits orbes. Il existe donc une composée de rotations $M$ telle que l'image du point $P$ soit~:
$$M\,R_{P_2,\mathbf{a},2°30'}(P')$$
où
$$P'=R_{P_1,\lambda_A}\,R_{P_2,\overline{\kappa}}\,R_{P_3,-\overline{\kappa}}\,
R_{P_4,2\overline{\kappa}}\,R_{P_5,\overline{\alpha}-\overline{\kappa}}(P)$$
On introduit de même les points $P_3'$, $P_4'$, $P_5'$ (voir figure \ref{fig040})~:
$$P_3'=R_{P_1,\lambda_A}\,R_{P_2,\overline{\kappa}}(P_3)$$
$$P_4'=R_{P_1,\lambda_A}\,R_{P_2,\overline{\kappa}}\,R_{P_3,-\overline{\kappa}}(P_4)$$
$$P_5'=R_{P_1,\lambda_A}\,R_{P_2,\overline{\kappa}}\,
R_{P_3,-\overline{\kappa}}\,R_{P_4,2\overline{\kappa}}(P_5)$$
Tous ces points sont dans le plan de la figure initiale~: calculer leurs positions relève entièrement de la géométrie plane.

\paragraph{Les équations de Saturne}\label{equ_sat}
Le calcul de la position du point $P'$ revient à résoudre quatre triangles rectangles. On calcule d'abord la \emph{première équation} dite aussi ``équation du centre''~; c'est l'angle $c_1=(\overrightarrow{OP_3'},\overrightarrow{OP_5'})$. On a~:
$$c_1=-\arcsin\left(\frac{P_3P_4\sin\overline{\kappa}+P_4P_5\sin\overline{\kappa}}{OP_5'}\right)$$
où
$$OP_5'=\sqrt{(P_3P_4\sin\overline{\kappa}+P_4P_5\sin\overline{\kappa})^2+(OP_3+P_3P_4\cos\overline{\kappa}-P_4P_5\cos\overline{\kappa})^2}.$$
On calcule ensuite la \emph{deuxième équation} $c_2=(\overrightarrow{OP_5'},\overrightarrow{OP'})$~:
$$c_2=\arcsin\left(\frac{P_5P\sin(\overline{\alpha}-c_1)}{OP'}\right)$$
où
$$OP'=\sqrt{(P_5P\sin(\overline{\alpha}-c_1))^2+(OP_5'+P_5P\cos(\overline{\alpha}-c_1))^2}$$
La position du point $P'$ est alors donnée par $OP'$ et par l'angle suivant~:
$$(\mathbf{j},\overrightarrow{OP'})=\lambda_A+\overline{\kappa}+c_1+c_2$$
Dans ce qui suit, $c_1$ sera conçue comme une fonction de $\overline{\kappa}$ notée $c_1(\overline{\kappa})$, et $c_2$ comme une fonction de deux variables $\overline{\kappa}$ et $\alpha=\overline{\alpha}-c_1$, notée $c_2(\overline{\kappa},\alpha)$. Plus précisément~:
$$c_2(\overline{\kappa},\alpha)=\arcsin\left(\frac{P_5P\sin\alpha}{\sqrt{(P_5P\sin\alpha)^2+(OP_5'+P_5P\cos\alpha)^2}}\right)$$
où $OP_5'$ est la fonction de $\overline{\kappa}$ vue ci-dessus. On remarque que~:
$$c_1(360°-\overline{\kappa})=-c_1(\overline{\kappa}),$$
$$c_2(360°-\overline{\kappa},360°-\alpha)=-c_2(\overline{\kappa},\alpha).$$

\paragraph{Fonction de deux variables et interpolation}
Les valeurs de $c_1$ et $c_2$ devront être reportées dans des tables. Comme pour la Lune, {\shatir} utilise à cet effet\footnote{\textit{Cf.} chap. 23 de la première partie de la \textit{Nih\=aya}.} la méthode ptoléméenne d'interpolation visant à calculer une valeur approchée de $c_2$ au moyen d'un produit d'une fonction de $\overline{\kappa}$ par une fonction de $\alpha$. \`A $\alpha$ donné, il interpole les valeurs de $c_2$ entre $c_2(0,\alpha)$ et $c_2(180°,\alpha)$ au moyen de la formule suivante~:
$$c_2(\overline{\kappa},\alpha)\simeq c_2(0,\alpha)+\chi(\overline{\kappa})(c_2(180°,\alpha)-c_2(0,\alpha)).$$
Le coefficient d'interpolation $\chi$ est défini par~:
$$\chi(\overline{\kappa})=\frac{\max\vert c_2(\overline{\kappa},\cdot)\vert-\max\vert c_2(0,\cdot)\vert}{\max\vert c_2(180°,\cdot)\vert-\max\vert c_2(0,\cdot)\vert}.$$
Ce coefficient est fonction d'une seule variable~; on peut donc le calculer pour une série de valeurs de $\overline{\kappa}$ et en faire une table. Pour ce faire, on remarque que le maximum $\max\vert c_2(\overline{\kappa},\cdot)\vert$ est atteint quand la droite $(O,P')$ est tangente au cercle de l'orbe de l'épicycle. On a alors\footnote{On le montre aisément par le calcul en posant $z^{-1}=OP_5'+P_5P\cos\alpha$ et en remarquant que la quantité $\tan^2c_2$ à maximiser est alors un polynôme de degré 2 en $z$.}~:
$$\max\vert c_2(\overline{\kappa},\cdot)\vert=\arcsin\frac{P_5P}{OP'_5}.$$

\begin{figure}
  \begin{center}
    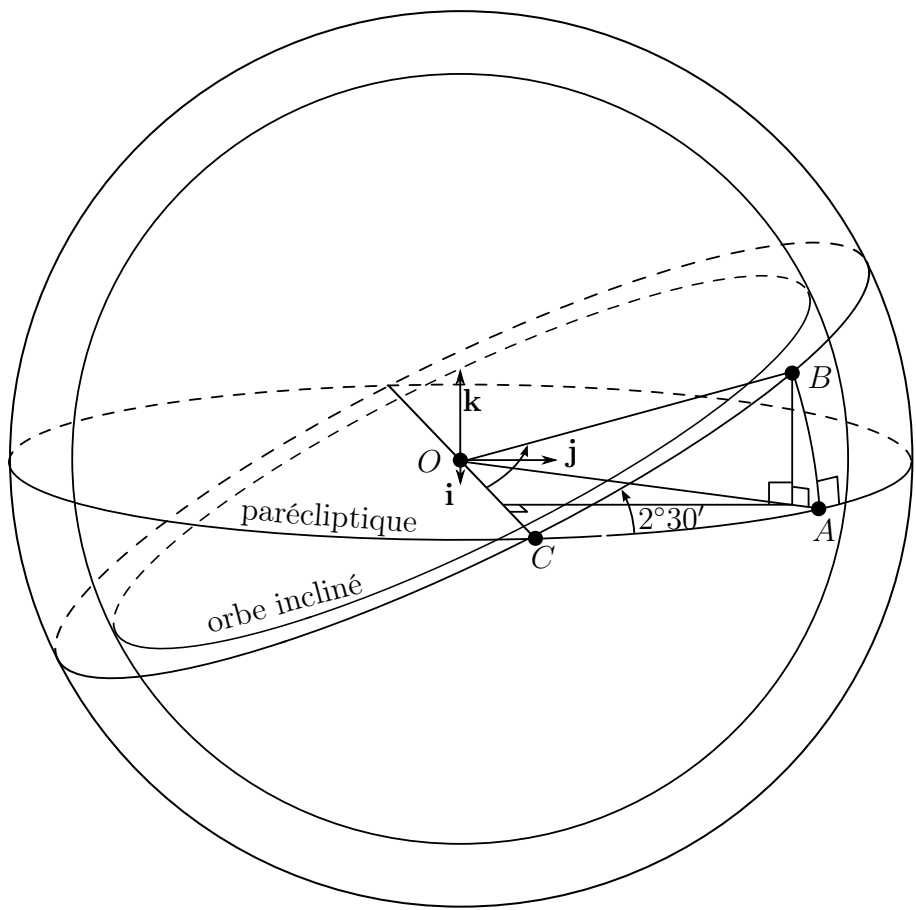
  \end{center}
  \caption{\label{fig042}Saturne~: inclinaison de l'orbe incliné}
\end{figure}

\paragraph{Trigonométrie sphérique}
On va à présent calculer les coordonnées sphériques du point $R_{P_2,\mathbf{a},2°30'}(P')$ par rapport à l'écliptique. L'angle formé entre le vecteur $\mathbf{a}$ et la direction du point $P'$ vaut (modulo $360°$)~:
$$(\mathbf{a},\overrightarrow{OP'})=\overline{\kappa}+c_1+c_2+140°=\overline{\lambda}-\lambda_{\ascnode}+c_1+c_2=\lambda-\lambda_{\ascnode}$$
où $\lambda=\overline{\lambda}+c_1+c_2$. Les égalités suivantes seront prises modulo $360°$. Sur la figure \ref{fig042}, on représente sur la sphère de l'écliptique le point $B$ dans la direction du point $R_{P_2,\mathbf{a},2°30'}(P')$ et le point $C$ dans la direction du vecteur $\mathbf{a}$. En particulier,
$$(\overrightarrow{OC},\overrightarrow{OB})=\lambda-\lambda_{\ascnode}.$$
L'étude du triangle sphérique $ABC$ donne~:
$$\tan(\overrightarrow{OC},\overrightarrow{OA})=\cos(2°30')\times\tan(\lambda-\lambda_{\ascnode}).$$
La longitude du point $R_{P_2,\mathbf{a},2°30'}(P')$ par rapport à l'écliptique, en prenant la direction du point vernal $\mathbf{j}$ comme origine, est donc, quand $\lambda-\lambda_{\ascnode}\in\rbrack -90°,90°\lbrack$~:
$$(\mathbf{j},\overrightarrow{OA})=(\overrightarrow{OC},\overrightarrow{OA})-(\overrightarrow{OC},\mathbf{j})=\arctan(\cos(2°30')\times\tan(\lambda-\lambda_{\ascnode}))+\lambda_{\ascnode}$$
Quand au contraire $\lambda-\lambda_{\ascnode}\in\rbrack 90°,270°\lbrack$, on a~:
$$(\mathbf{j},\overrightarrow{OA})=180°+\arctan(\cos(2°30')\times\tan(\lambda-\lambda_{\ascnode}))+\lambda_{\ascnode}.$$
On rassemble ces deux cas dans la formule suivante, valable pour toutes les valeurs des paramètres~:
$$(\mathbf{j},\overrightarrow{OA})=\lambda+e_n(\lambda-\lambda_{\ascnode})$$
où l'``équation du déplacement'' $e_n(x)$ est définie comme suit sur l'intervalle $\lbrack -90°,270°\rbrack$, et ailleurs par périodicité~:
$$e_n(x)=\left\lbrace\begin{array}{l}
\arctan(\cos(2°30')\times\tan x)-x,\text{ si }x\in\rbrack -90°,90°\lbrack\\
180°+\arctan(\cos(2°30')\times\tan x)-x,\text{ si }x\in\rbrack 90°,270°\lbrack\\
0°,\text{ si }x=\pm 90°
\end{array}\right.$$
Enfin, la latitude du point $R_{P_2,\mathbf{a},2°30'}(P')$ est~:
$$(\overrightarrow{OA},\overrightarrow{OB})=\arcsin(\sin(2°30')\times\sin(\lambda-\lambda_{\ascnode})).$$

\paragraph{Les inclinaisons des petits orbes}
Les rotations qui composent $M$ ont leurs axes contenus dans le plan de l'orbe incliné $OCB$ et sont d'angles petits ($3°30'$ au plus). Elles auront peu d'effet sur la longitude de $R_{P_2,\mathbf{a},2°30'}(P')$. En revanche, l'effet de $M$ et les inclinaisons des plans des petits orbes ont justement été introduits dans ce modèle afin de rendre compte des variations en latitude de Saturne. Hélas, les axes des rotations qui composent $M$ ne passent pas par le point $O$, et de telles rotations sont difficiles à étudier en coordonnées sphériques\footnote{Ibn al-Hayth\=am lui-même, dans un ouvrage inégalé sur l'astronomie mathématique, avait renoncé à une tâche semblable~; \textit{cf.} \cite{rashed2006} p.~444-447.}. Ibn al-\v{S}\=a\d{t}ir se contente ici d'une description qualitative de l'effet en latitude de $M$ (voir chapitre 24).

\begin{figure}
  \begin{center}
    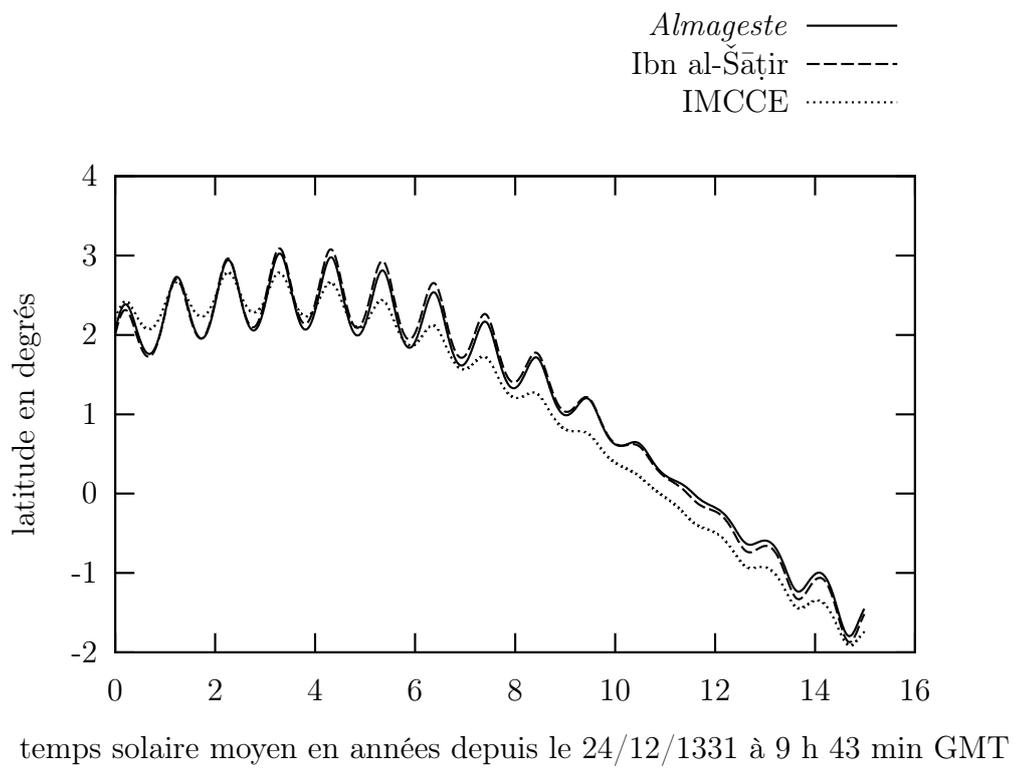
  \end{center}
  \caption{\label{fig044}Latitudes de Saturne sur quinze ans, à partir de l'\'Epoque}
\end{figure}

\paragraph{Comparaisons}
\`A la figure \ref{fig044} nous donnons les latitudes de Saturne pendant quinze ans, à partir de l'\'Epoque choisie par {\shatir}, calculées de trois manières différentes~: (1) par l'IMCCE de l'Observatoire de Paris, (2) en suivant la méthode strictement ptoléméenne de l'\textit{Almageste} mais avec les paramètres d'{\shatir} à l'\'Epoque et ses mouvements moyens\footnote{On a suivi l'interprétation des méthodes de Ptolémée donnée dans \cite{pedersen1974} et \cite{swerdlow2005}.}, (3) au moyen du modèle complet d'{\shatir} décrit par la composée de huit rotations appliquée au point $P$ comme ci-dessus.

Quand {\shatir} indique que le n{\oe}ud ascendant de Saturne est ``à cinq parts du Lion'', il faut peut-être comprendre qu'il est à $5°$ \emph{avant} le début du Lion, à $115°$ du point vernal, en 1339 ap. J.-C. \textit{i. e.} à $t=8$ ans de l'\'Epoque, quand l'Apogée de Saturne est à $254;52+8\times 0;1=255°$ du point vernal\footnote{\textit{Cf.} p.~\pageref{noeud_sup1} \textit{supra}.}~; on aurait alors $\lambda_{\ascnode}=\lambda_A-140$ comme l'affirmait déjà Ptolémée $1200$ ans auparavant. C'est la valeur que nous avons adoptée dans le commentaire ci-dessus. Hélas, l'élongation n{\oe}ud-Apogée n'est pas constante, et ceci explique le décalage important entre le tracé des latitudes selon {\shatir} et la courbe de l'IMCCE sur la figure \ref{fig044}, près du n{\oe}ud ascendant. Nous en concluons qu'{\shatir} n'a pas cherché à corriger la valeur donnée par Ptolémée en confrontant son modèle en latitude avec l'observation directe pour Saturne\footnote{Lire $5°$ \emph{après} le début du Lion, \textit{i. e.} à $125°$ du point vernal, conduirait à un désaccord encore plus grand avec d'hypothétiques observations du temps d'{\shatir}.}. Il semble néanmoins l'avoir fait pour Jupiter~: au même endroit, il affirme que, pour Jupiter, $\lambda_{\ascnode}=\lambda_A-70=181-70=111$, valeur coïncidant avec celle donnée par Ptolémée, mais il ajoute plus loin\footnote{\textit{Cf.} p.~\pageref{noeud_sup2} \textit{supra}.} que ``d'après ses propres observations'' $\lambda_{\ascnode}=\lambda_A-62$. Enfin, pour Mars, il ne mentionne pas la valeur donnée par Ptolémée\footnote{Selon Ptolémée, pour Mars, $\lambda_{\ascnode}=\lambda_A-85°30'$.} et il affirme simplement que $\lambda_{\ascnode}=\lambda_A-90=138-90=48$.

Comme l'a montré Swerdlow\footnote{\textit{Cf.} \cite{swerdlow2005} p.~47.}, la théorie des latitudes des planètes supérieures dans l'\textit{Almageste} n'est guère satisfaisante pour l'observateur moderne~; elle commet souvent des erreurs de l'ordre de $20'$. Ceci tient à trois causes~:

-- Dans le référentiel héliocentrique, le diamètre de l'épicycle indique la direction Soleil-Terre. Le plan de l'épicycle devrait donc toujours rester parallèle à l'écliptique\footnote{à peu près, car Soleil moyen $\simeq$ Soleil vrai.}, contrairement à ce que pense Ptolémée.

-- Les observations transmises par Ptolémée sont peu précises~; elles sont arrondies au degré près~! De plus ces observations des latitudes maximales, lors des oppositions, ou près des conjonctions, devaient être rares pour Saturne dont la période moyenne est de 30 ans environ.

-- Les calculs pour déduire les inclinaisons de l'excentrique et de l'épicycle sont sensibles à la précision de ces observations.

Comme {\shatir} s'appuie sur les observations de Ptolémée et sur ses inférences concernant les inclinaisons maximales des plans des orbes, son modèle présente alors à peu près les mêmes défauts en latitude que celui de Ptolémée.

\begin{table}
  \begin{center}\begin{tabular}{c|c|c|c}
    & Saturne & Jupiter & Mars\\
    $\overrightarrow{P_3P_4}$ & $5;7,30\,\mathbf{j}$ & $4;7,30\,\mathbf{j}$ & $9\,\mathbf{j}$\\
    $\overrightarrow{P_4P_5}$ & $-1;42,30\,\mathbf{j}$ & $-1;22,30\,\mathbf{j}$ & $-3\,\mathbf{j}$\\
    $\overrightarrow{P_5P}$ & $6;30\,\mathbf{j}$ & $11;30\,\mathbf{j}$ & $39;30\,\mathbf{j}$\\
    $(\mathbf{j},\mathbf{u})$ & $-140°$ & $-62°$ & $-90°$\\
    $\overline{\lambda}(0)$ & $157;58,20$ & $272;6,10$ & $292;0,0$\\
    $\dot{\overline{\lambda}}$, par année persane & $12;13,40$ & $30;20,33$ & $191;17,11$\\
    $\lambda_A(0)$ & $254;52$ & $180;52$ & $137;52$ \\
    $\dot{\lambda}_A$, par année persane & $0;1$ & $0;1$ & $0;1$ \\
    inclinaison de l'orbe incliné & $R_{P_2,\mathbf{u},2°30'}$ & $R_{P_2,\mathbf{u},1°30'}$ & $R_{P_2,\mathbf{u},1°}$ \\
    inclinaison du déférent & $R_{P_3,\mathbf{u},-3°30'}$ & $R_{P_3,\mathbf{u},-2°}$ & $R_{P_3,\mathbf{u},-1°37'30''}$ \\
    inclinaison du rotateur & $R_{P_4,\mathbf{v},-1°}$ & $R_{P_4,\mathbf{v},-0°30'}$ & $R_{P_4,\mathbf{v},-0°37'30''}$
  \end{tabular}\end{center}
  \caption{Synopsis des paramètres utilisés par {\shatir} pour les planètes supérieures}
\end{table}

\begin{figure}
  \begin{center}
    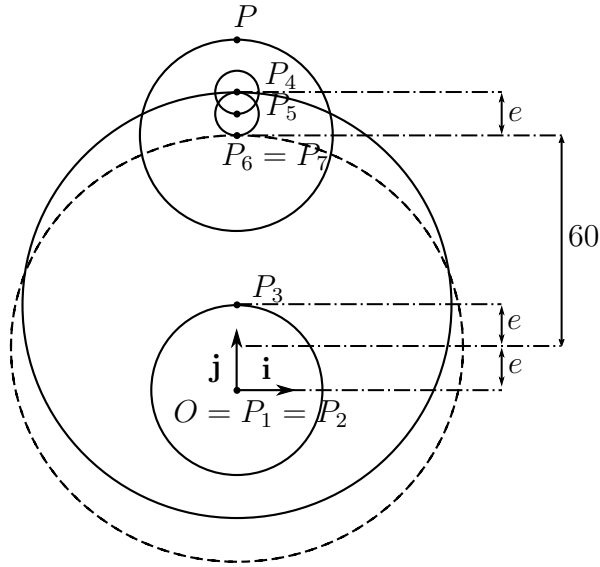
  \end{center}
  \caption{\label{fig046}Saturne, al-\d{T}\=us{\=\i}, figure initiale}
\end{figure}

\paragraph{Les longitudes des planètes supérieures chez \d{T}\=us{\=\i}, `Ur\d{d}{\=\i} et Sh{\=\i}r\=az{\=\i}} Le modèle de \d{T}\=us{\=\i} pour les planètes supérieures comprend sept orbes\footnote{\textit{Cf.} \cite{altusi1993} II.11[10], p.~208.}. Sur la figure initiale (fig.~\ref{fig046}), $P_1$ est le centre du parécliptique, $P_2$ le centre de l'orbe incliné, $P_3$ le centre d'un orbe déférent excentrique, $P_4$ le centre d'une ``grande sphère'', $P_5$ le centre d'une ``petite sphère'', $P_6$ le centre d'un orbe englobant, $P_7$ le centre de l'épicycle, et $P$ le centre du globe planétaire. On a $O=P_1=P_2$, et $P_6=P_7$. Pour Saturne, on note $e=3;25$ et on a~:\footnote{Pour Mars et Jupiter, on a le même modèle que Saturne, seuls $e$ et $P_7P$ diffèrent.}
\begin{align*}
  &\overrightarrow{P_2P_3}=2e\,\mathbf{j}\\
  &\overrightarrow{P_3P_4}=60\,\mathbf{j}\\
  &\overrightarrow{P_4P_5}=\overrightarrow{P_5P_6}=-\dfrac{e}{2}\,\mathbf{j}\\
  &\overrightarrow{P_7P}=6;30\,\mathbf{j}
\end{align*}
\d{T}\=us{\=\i} adjoint trois orbes supplémentaires pour rendre compte des mouvements en latitude au moyen d'un ``couple de \d{T}\=us{\=\i} curviligne''\footnote{\textit{Cf.} \cite{altusi1993} II.11[19], p.~218-220.}~; on les négligera ici. Les rotations décrivant le mouvement en longitude de chaque planète supérieure selon \d{T}\=us{\=\i} sont alors~:
$$R_{P_1,\lambda_A}
\,R_{P_3,\overline{\kappa}}
\,R_{P_4,\overline{\kappa}}
\,R_{P_5,-2\overline{\kappa}}
\,R_{P_6,\overline{\kappa}}
\,R_{P_7,\overline{\alpha}}.$$
On aura remarqué que la ``petite'' et la ``grande'' sphère constituent un couple de \d{T}\=us{\=\i} dont l'effet est décrit par $R_{P_4,\overline{\kappa}}\,R_{P_5,-2\overline{\kappa}}\,R_{P_6,\overline{\kappa}}$.

Le modèle de `Ur\d{d}{\=\i} pour les planètes supérieures comprend cinq orbes dont les centres sont~:\footnote{\textit{Cf.} \cite{saliba1990}, et \cite{rashed1997} p.~119-122.} $O=P_1=P_2$ le centre du Monde, $P_3$ le centre du déférent que Sh{\=\i}r\=az{\=\i} appellera ``déférent corporel'' (car il est bien distinct du déférent de Ptolémée), $P_4$ le centre de l'orbe rotateur, $P_5$ le centre de l'épicycle, et $P$ le centre du globe planétaire. Notant toujours, pour Saturne, $e=3;25$, on pose~:
\begin{align*}
  &\overrightarrow{P_2P_3}=\dfrac{3e}{2}\,\mathbf{j}=5;7,30\,\mathbf{j}\\
  &\overrightarrow{P_3P_4}=60\,\mathbf{j}\\
  &\overrightarrow{P_4P_5}=-\dfrac{e}{2}\,\mathbf{j}\\
  &\overrightarrow{P_5P}=6;30\,\mathbf{j}
\end{align*}
Le mouvement en longitude est alors décrit par les rotations suivantes~:
$$R_{P_1,\lambda_A}
\,R_{P_3,\overline{\kappa}}
\,R_{P_4,\overline{\kappa}}
\,R_{P_5,\overline{\alpha}-\overline{\kappa}}.$$

\begin{figure}
  \begin{center}
    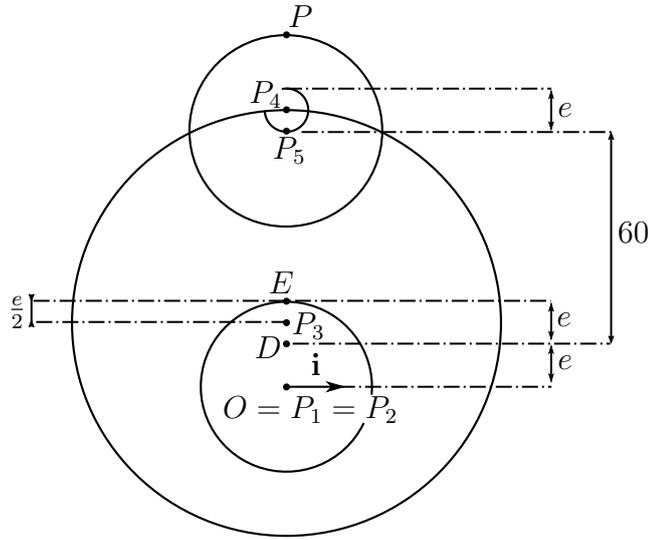
  \end{center}
  \caption{\label{fig047}Saturne, al-`Ur\d{d}{\=\i}, figure initiale}
\end{figure}

Pour la Lune, on a vu Sh{\=\i}r\=az{\=\i} proposer un modèle strictement équivalent au modèle de \d{T}\=us{\=\i} sans couple curviligne, assez différent du modèle erroné qu'offrait `Ur\d{d}{\=\i} dans son \textit{Kit\=ab al-hay'a}. Pour les planètes supérieures, Sh{\=\i}r\=az{\=\i} semble avoir hésité, dans les années 1281-1285, entre un modèle fondé sur un ``principe conjectural'' et le modèle de `Ur\d{d}{\=\i} (ou un modèle équivalent) fondé sur un ``principe déductif''. Ma connaissance des trois ouvrages écrits par Sh{\=\i}r\=az{\=\i} dans ces années-là\footnote{Les \textit{Ikht{\=\i}y\=ar\=at-i muzaffar{\=\i}}, la \textit{Nih\=aya al-idr\=ak} et la \textit{Tu\d{h}fa al-sh\=ah{\=\i}ya}. \textit{Cf.} \cite{rashed1997}, \cite{masoumi2013}, \cite{niazi2011}.} s'appuie sur les travaux de Saliba, Gamini-Masoumi et Niazi~; il me faudra peut-être réviser mes jugements quand les textes originaux seront plus accessibles. Niazi voit un peu la \textit{Nih\=aya al-idr\=ak} comme un commentaire de la \textit{Tadhkira} de \d{T}\=us{\=\i} et les \textit{Ikht{\=\i}y\=ar\=at-i muzaffar{\=\i}} comme un auto-commentaire de la \textit{Nih\=aya al-idr\=ak}, bien que les deux ouvrages aient été rédigés presque simultanément\footnote{\textit{C.} \cite{niazi2011} p. 213.}. Il n'est donc pas surprenant que Sh{\=\i}r\=az{\=\i} reprenne la structure du traité de \d{T}\=us{\=\i} dans sa \textit{Nih\=aya al-idr\=ak} et y décrive d'abord le modèle de Ptolémée pour les planètes supérieures\footnote{\textit{Nih\=aya al-idr\=ak}, livre 2, chapitre 8. \textit{Cf.} \cite{niazi2011} p. 150 et appendices où Niazi donne un résumé analytique du chapitre 8 et une édition du début de ce chapitre. Une édition critique de la fin du chapitre 8 figure dans \cite{masoumi2013}.}. Puis, dans son chapitre sur les latitudes, comme l'avait fait \d{T}\=us{\=\i} avec son couple de \d{T}\=us{\=\i}, Sh{\=\i}r\=az{\=\i} propose un nouveau modèle\footnote{\textit{Cf.} \cite{niazi2011} p. 151. Il serait intéressant d'étudier le comportement du modèle de Sh{\=\i}r\=az{\=\i} pour les latitudes.}. Ce modèle est aussi décrit dans les \textit{Ikht{\=\i}y\=ar\=at} où Sh{\=\i}r\=az{\=\i} explique pourquoi il le préfère au modèle de `Ur\d{d}{\=\i}\footnote{\textit{Cf.} \cite{niazi2011} p.~203-204.}~: dans le modèle de `Ur\d{d}{\=\i}, la trajectoire du centre de l'épicycle dans le plan de l'orbe incliné n'est pas exactement un cercle\footnote{On retrouve évidemment le même défaut dans le modèle de \d{T}\=us{\=\i}, comme nous l'avions déjà remarqué pour la Lune (\textit{cf.} la courbe en pointillés, fig. \ref{fig021} p.~\pageref{fig021} \textit{supra}).}, mais c'en est bien un dans le nouveau modèle de Sh{\=\i}r\=az{\=\i}. Dans les \textit{Ikht{\=\i}y\=ar\=at}, Sh{\=\i}r\=az{\=\i} utilise les termes ``déductif'' et ``conjectural'' pour marquer l'opposition entre les deux modèles\footnote{Le ``principe déductif'', \textit{istinb\=a\d{t}{\=\i}}, désigne le dispositif à l'{\oe}uvre dans le modèle de `Ur\d{d}{\=\i}, et le ``principe conjectural'', \textit{\d{h}ads{\=\i}}, désigne celui à l'{\oe}uvre dans le nouveau modèle de Sh{\=\i}r\=z{\=\i}. \`A cette occasion, il utilise aussi un troisième terme, ``innovant'', \textit{ibd\=a`{\=\i}}, pour désigner un principe que Niazi n'a pas su expliquer (\textit{cf.} \cite{niazi2011} note 52 p.~138). Ce terme se rapporte peut-être aux dispositifs utilisés pour décrire les mouvements en latitude, puisque {\shatir} en parle ainsi dans la \textit{Nih\=aya al-S\=ul} (p.~\pageref{ibdaai} \textit{supra})~:
  
  ``Le principe qu'il a nommé l'invention (\textit{ibd\=a`{\=\i}}) concernant le désordre des latitudes des astres est impossible.''

  Selon Niazi, les termes \textit{istinb\=a\d{t}{\=\i}} et \textit{\d{h}ads{\=\i}} sont aussi utilisé dans la \textit{Nih\=aya al-idr\=ak}, mais ils seront absents de la \textit{Tu\d{h}fa} (\cite{niazi2011} p.~159).}. Pourtant, Niazi montre qu'après de nombreuses révisions, l'état final du chapitre sur les latitudes de la \textit{Nih\=aya al-idr\=ak} semble préférer une solution fondée sur le ``principe déductif'', c'est-à-dire à la manière de `Ur\d{d}{\=\i}~! La \textit{Tu\d{h}fa}, au dire de Niazi, reflète sûrement son choix définitif, et on y trouve le modèle suivant. C'est le seul que nous retiendrons ici~:

\begin{quote}
  ``Le \textit{premier orbe} est le parécliptique~; la partie convexe [du parécliptique] de Saturne touche la partie concave du huitième orbe, sa partie concave touche la partie convexe du parécliptique de Jupiter, la partie concave du parécliptique de Jupiter touche la partie convexe du parécliptique de Mars, et la partie concave du parécliptique de Mars touche la partie convexe du parécliptique de Vénus. Le \textit{deuxième orbe} est le déférent excentrique portant le centre de [l'orbe] englobant~; il est dans la partie supérieure du parécliptique~; la distance de son centre au centre du déférent imaginaire est la moitié de [la distance] entre les centres du Monde et du déférent imaginaire de l'astre~; et sa ceinture est inclinée par rapport à la trajectoire du Soleil, de la grandeur de l'inclinaison de [l'orbe] incliné de l'astre, inclinaison constante. Le \textit{troisième orbe} est [l'orbe] englobant~; il est dans la partie supérieure de l'excentrique~; son axe est perpendiculaire au plan de la ceinture de l'excentrique, et sa ceinture est dans son plan c'est-à-dire dans le plan de [l'orbe] incliné. Le \textit{quatrième orbe} est le fléchisseur\footnote{Nous traduisons ainsi \textit{mumayyil}~; cet orbe est responsable de l'inclinaison variable de la ceinture de l'épicycle.}~; il est dans la partie inférieure de l'orbe englobant~; son axe est parallèle à l'axe de l'orbe englobant et perpendiculaire au plan de l'orbe incliné, et sa ceinture est aussi dans son plan~; et la distance de son centre au centre de l'orbe englobant est égale à celle entre les centres de l'excentrique et du déférent imaginaire de l'astre [...]. Le \textit{cinquième orbe}~: l'astre tourne autour du centre du fléchisseur et d'un axe coupant l'axe du fléchisseur en le centre commun. La ceinture [du cinquième orbe] est inclinée par rapport à la ceinture [du fléchisseur] vers le Nord et vers le Sud, de la grandeur de l'inclinaison, pour l'astre en question, par rapport au plan de l'orbe incliné, et c'est une inclinaison constante.''\footnote{Notre traduction. \textit{Cf.} manuscrit Arabe 2516 à la Bibliothèque Nationale de France~; chapitre 11, f. 45v.}
\end{quote}

\begin{figure}
  \begin{center}
    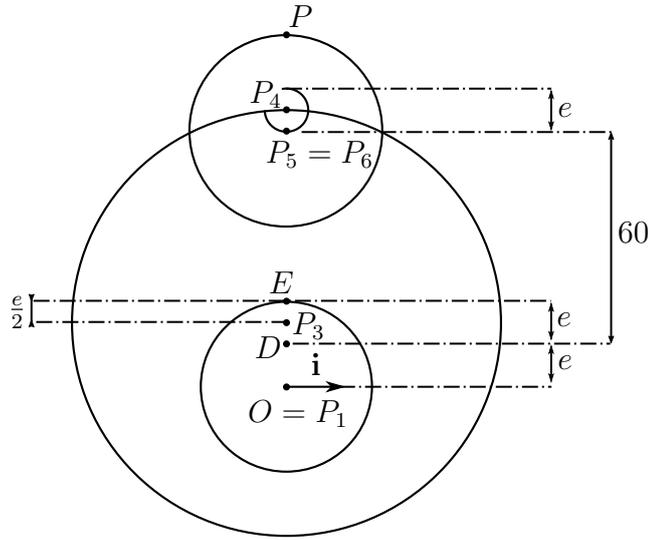
  \end{center}
  \caption{\label{fig048}Saturne, Sh{\=\i}r\=az{\=\i}, figure initiale}
\end{figure}

Les centres des orbes dans le modèle de la \textit{Tu\d{h}fa} sont~: $O=P_1$ centre du parécliptique, $P_3$ centre d'un ``déférent corporel'', $P_4$ centre de l'orbe ``englobant'', $P_5$ centre d'un orbe ``fléchisseur'', $P_6=P_5$ centre de l'épicycle, et $P$ centre du globe planétaire. On notera $D$ le centre du ``déférent imaginaire'' qui n'est autre que le centre de la trajectoire approximativement circulaire du point $P_5$ au sein du parécliptique. On a~:
\begin{align*}
  &\overrightarrow{P_1P_3}=\dfrac{3e}{2}\,\mathbf{j}=5;7,30\,\mathbf{j}\\
  &\overrightarrow{P_3P_4}=60\,\mathbf{j}\\
  &\overrightarrow{P_4P_5}=-\dfrac{e}{2}\,\mathbf{j}=-1;42,30\,\mathbf{j}\\
  &\overrightarrow{P_5P}=6;30\,\mathbf{j}
\end{align*}
Si on néglige les inclinaisons, la composée de rotations entraînant le point $P$ dans son mouvement en longitude est alors~:
$$R_{P_1,\lambda_A}
\,R_{P_3,\overline{\kappa}}
\,R_{P_4,\overline{\kappa}}
\,R_{P_5,-\overline{\kappa}}
\,R_{P_6,\overline{\alpha}}.$$
Pour les longitudes, on voit que ce modèle est équivalent à celui de `Ur\d{d}{\=\i}. Si Sh{\=\i}r\=az{\=\i} semble avoir finalement renoncé à proposer un modèle différent, il a cependant offert une contribution originale à la fin du chapitre 8 du livre 2 de sa \textit{Nih\=aya}. Pour répondre à d'éventuels contradicteurs qui reprocheraient aux astronomes de Maragha d'avoir déplacé le centre du déférent concentrique que Ptolémée avait posé au milieu entre le centre du Monde et le point équant (``bissection de l'excentricité''), Sh{\=\i}r\=az{\=\i} tente d'abord de justifier par des arguments issus de l'observation le choix fait par Ptolémée. Gamini et Masoumi ont récemment analysé cette intéressante réflexion sur les ``arcs rétrogrades''\footnote{Le passage correspondant de la \textit{Nih\=aya al-idr\=ak} et des \textit{Ikht{\=\i}y\=ar\=at} est aussi étudié dans \cite{niazi2011} p.~195-200, mais Niazi a une interprétation légèrement différente de celle de \cite{masoumi2013}. \`A ses yeux il faut voir dans tout ce passage une tentative de justifier la préférence pour \textit{les deux} modèles de `Ur\d{d}{\=\i} et de Sh{\=\i}r\=az{\=\i} vis-à-vis d'éventuels contradicteurs qui souhaiteraient revenir à Ptolémée~; aux yeux de Gamini-Masoumi (\cite{masoumi2013} p.~51) il s'agit plutôt de justifier le refus du modèle de `Ur\d{d}{\=\i}.}. Ce faisant, Sh{\=\i}r\=az{\=\i} caractérise le centre du déférent comme étant le centre de la trajectoire, circulaire chez Ptolémée, du centre de l'épicycle dans le plan de l'orbe incliné, et il caractérise le point équant comme étant le point autour duquel le mouvement du centre de l'épicycle a une vitesse angulaire constante. Dans tous les modèles des astronomes de Maragha ainsi qu'{\shatir} pour les planètes supérieures, il existe bien un tel ``point équant'' et un point proche du centre de la trajectoire circulaire ou presque circulaire du centre de l'épicycle -- même si ces points ne sont pas toujours mentionnés dans les descriptions des modèles. Pour les modèles de `Ur\d{d}{\=\i} et de Sh{\=\i}r\=az{\=\i}, ce sont les points $E$ et $D$ sur nos figures \ref{fig047} et \ref{fig048}. 

En fait le modèle de `Ur\d{d}{\=\i} (et donc aussi celui de Sh{\=\i}r\=az{\=\i}) est strictement équivalent au modèle de \d{T}\=us{\=\i}. Partons en effet du modèle de `Ur\d{d}{\=\i} décrit comme ci-dessus par la composée de rotations $R_{P_1,\lambda_A}
\,R_{P_3,\overline{\kappa}}
\,R_{P_4,\overline{\kappa}}
\,R_{P_5,\overline{\alpha}-\overline{\kappa}}$ appliquée au point $P$ de la fig. \ref{fig047}. $P_1$, $P_2$, $P_3$, $P_4$, $P_5$ désigneront les centres des orbes dans la figure initiale de `Ur\d{d}{\=\i}, et $\tilde{P}_1$, $\tilde{P}_2$, $\tilde{P}_3$, $\tilde{P}_4$, $\tilde{P}_5$, $\tilde{P}_6$, $\tilde{P}_7$ désigneront les centres des orbes dans la figure initiale d'al-\d{T}\=us{\=\i}. On a~:
$$O=P_1=P_2=\tilde{P}_1=\tilde{P}_2,\quad\overrightarrow{P_3\tilde{P}_3}=\frac{e}{2}\,\mathbf{j},\quad P_4=\tilde{P}_5,$$
$$P_5=\tilde{P}_6=\tilde{P}_7,\quad \overrightarrow{P_4\tilde{P}_4}=\frac{e}{2}\,\mathbf{j}.$$
Alors~:
\begin{align*}
R_{P_3,\overline{\kappa}}
\,R_{P_4,\overline{\kappa}}
\,R_{P_5,\overline{\alpha}-\overline{\kappa}}
&= \left(R_{P_3,\overline{\kappa}}\,R_{P_4,-\overline{\kappa}}\right)
R_{P_4,2\overline{\kappa}}
\,R_{P_5,\overline{\alpha}-\overline{\kappa}}\\
&= R_{\tilde{P}_3,\overline{\kappa}}
\left(R_{\tilde{P}_4,-\overline{\kappa}}
\,R_{P_4,2\overline{\kappa}}
\,R_{P_5,\overline{\alpha}-\overline{\kappa}}\right)\quad\text{(proposition 1)}\\
&=R_{\tilde{P}_3,\overline{\kappa}}
\,R_{\tilde{P}_4,\overline{\kappa}}
\,R_{P_4,-2\overline{\kappa}}
\,R_{P_5,\overline{\alpha}+\overline{\kappa}}\quad\text{(proposition 3)}\\
&=R_{\tilde{P}_3,\overline{\kappa}}
\,R_{\tilde{P}_4,\overline{\kappa}}
\,R_{\tilde{P}_5,-2\overline{\kappa}}
\,R_{\tilde{P}_6,\overline{\kappa}}
\,R_{\tilde{P}_7,\overline{\alpha}}.
\end{align*}
\label{equiv_sup}

\begin{figure}
  \begin{center}
    \includegraphics{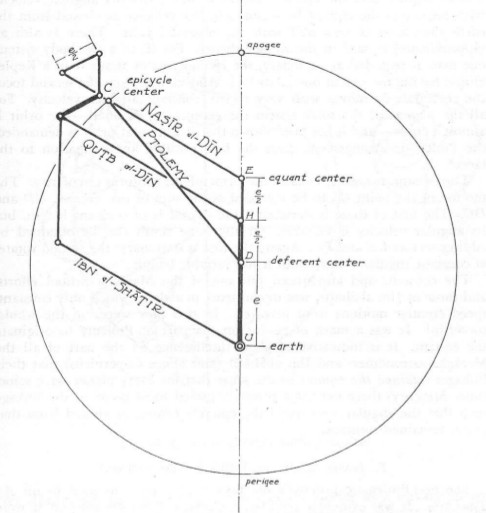}
  \end{center}
  \caption{\label{fig051}Démonstration d'équivalence entre les modèles en longitude de divers savants, selon Kennedy \cite{kennedy1966} p.~367}
\end{figure}

Enfin, pour les longitudes des planètes supérieures, le modèle d'{\shatir} est équivalent aux modèles des trois savants de Maragha. En effet, notons toujours $P_i$ les centres des orbes dans la figure initiale de `Ur\d{d}{\=\i}, mais notons à présent $\tilde{P}_i$ les centres des orbes d'{\shatir}. On a~:
$$O=P_1=P_2=\tilde{P}_1=\tilde{P}_2,\quad P_4=\tilde{P}_4,\quad P_5=\tilde{P}_5,\quad P=\tilde{P}.$$
Comme $\overrightarrow{\tilde{P}_2\tilde{P}_3}=\overrightarrow{P_3P_4}$, en vertu de la proposition 1~:
\begin{align*}
  R_{P_3,\overline{\kappa}}
  \,R_{P_4,\overline{\kappa}}
  \,R_{P_5,\overline{\alpha}-\overline{\kappa}}
  &=\left(R_{P_3,\overline{\kappa}}\,R_{P_4,-\overline{\kappa}}\right)
  \,R_{P_4,2\overline{\kappa}}
  \,R_{P_5,\overline{\alpha}-\overline{\kappa}}\\
  &=R_{\tilde{P}_2,\overline{\kappa}}
  \,R_{\tilde{P}_3,-\overline{\kappa}}
  \,R_{\tilde{P}_4,2\overline{\kappa}}
  \,R_{\tilde{P}_5,\overline{\alpha}-\overline{\kappa}}
\end{align*}
Une démonstration plus intuitive, au moyen d'une figure géométrique, circule beaucoup dans la littérature secondaire depuis Kennedy\footnote{\textit{Cf.} \cite{kennedy1966} p.~367.}. Nous l'avons reproduite, fig.~\ref{fig051} p.~\pageref{fig051}. 

\begin{figure}
  \begin{center}
    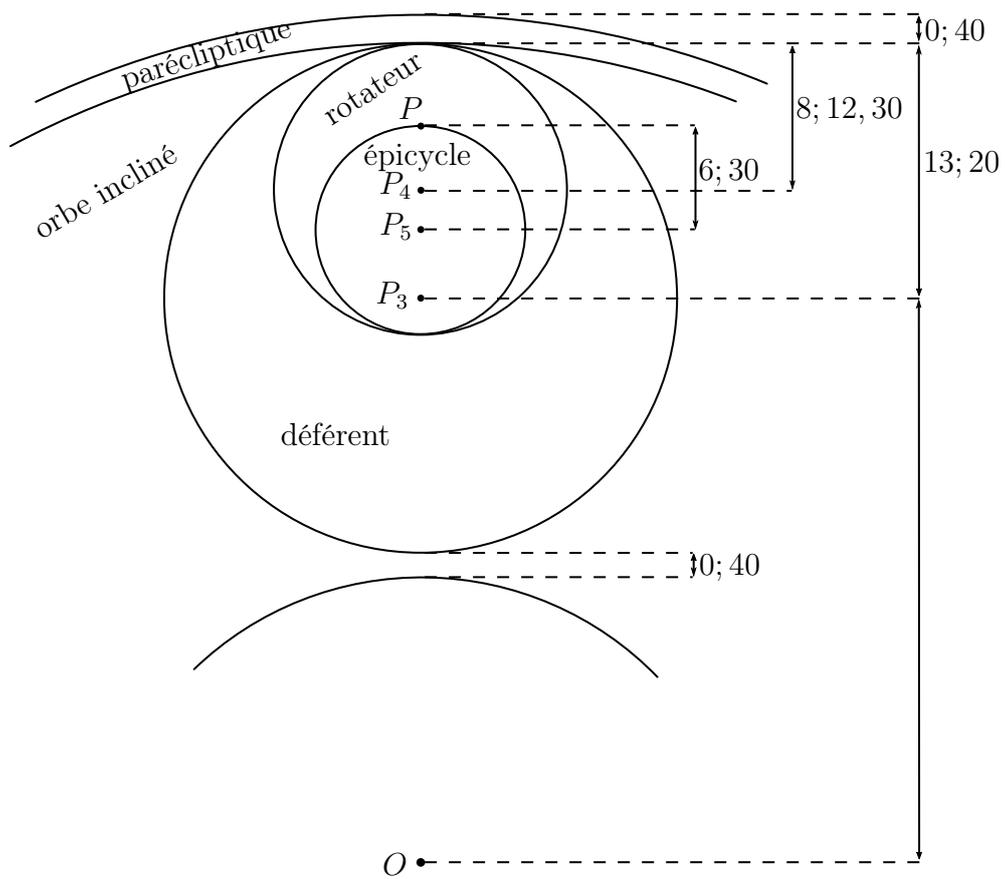
  \end{center}
  \caption{\label{fig049}Saturne, orbes solides}
\end{figure}

\paragraph{Sections planes des orbes solides de Saturne}
Les sections planes des orbes solides de Saturne sont représentées à la figure \ref{fig049} qui respecte à peu près les ordres de grandeurs sauf pour $OP_3$. Les rayons sont~:
\begin{align*}
  \text{rayon du globe planétaire} &= r\\
  \text{rayon de l'orbe de l'épicycle} &=6;30+r\\
  \text{rayon du rotateur} &=1;42,30+6;30+r=8;12,30+r\\
  \text{rayon du déférent} &=5;7,30+1;42,30+6;30+r=13;20+r\\
  \text{rayon extérieur du l'orbe incliné}&=60+5;7,30+1;42,30+6;30+r=73;20+r\\
  \text{rayon intérieur de l'orbe incliné}&=60-(5;7,30+1;42,30+6;30+r)=46;40-r\\
  \text{épaisseur du parécliptique}&=0;40
\end{align*}
Ici comme pour l'orbe total du Soleil, le choix de l'épaisseur $0;40$ du parécliptique est arbitraire, et {\shatir} propose encore d'ajouter un complément à l'épaisseur de l'orbe incliné. Mais pour Saturne, il néglige $r$ quand il calcule les rayons des orbes solides. Il n'en tient compte qu'à la fin de son raisonnement, dans le ``complément'' qu'il ajoute à l'épaisseur de l'orbe incliné, quand il dit~:
\begin{quote}
  ``La distance minimale est à la surface intérieure de l'orbe incliné de Saturne, en parts, quarante-six et deux tiers~; en plus de la réunion des orbes, il faut aussi compter le rayon de l'astre donc cela fait quarante-six parts.''\footnote{\textit{Cf.} p.~\pageref{contiguite_sat} \textit{supra}. Bien que le rayon intérieur de l'orbe incliné passe à $46$, il semble qu'{\shatir} maintienne son rayon extérieur à $73;20$~: c'est l'épaisseur \textit{du parécliptique} qui fait passer le rayon extérieur \textit{du système d'orbes} à $74$. Dans l'orbe incliné, si l'on tient compte du rayon $r$ de la planète, il faut alors supposer que le cercle médian n'est plus un cercle de rayon $60$ mais un cercle de rayon légèrement inférieur (entre $59;40$ et $60$ si $2r<0;40$). En toute rigueur, il faudrait donc réduire proportionnellement le système entier des orbes de Saturne par rapport aux valeurs rapportées ci-dessus.}
\end{quote}
Il vise large en arrondissant le rayon intérieur de l'orbe incliné à $46$. Le rayon intérieur du système d'orbes de Saturne est donc $46$, et son rayon extérieur est, en tenant compte de l'épaisseur du parécliptique, $73;20+0;40=74$. Il ne remettra pas ces valeurs en question dans la conclusion de la \textit{Nihaya}, et il suivra le même procédé pour les quatre autres planètes.

\paragraph{Les orbes solides de Jupiter}
\begin{align*}
  \text{rayon du globe planétaire}&=r\\
  \text{rayon de l'orbe de l'épicycle}&=11;30+r\\
  \text{rayon du rotateur}&=1;22,30+11;30+r=12;52,30+r\\
  \text{rayon de déférent}&=4;7,30+1;22,30+11;30+r=17+r\\
  \text{rayon extérieur de l'orbe incliné}&=60+4;7,30+1;22,30+11;30+r=77+r\\
  \text{rayon intérieur de l'orbe incliné}&=60-(4;7,30+1;22,30+11;30+r)=43-r\\
  \text{épaisseur du parécliptique}&=1
\end{align*}

\paragraph{Les orbes solides de Mars}
\begin{align*}
  \text{rayon du globe planétaire}&=r\\
  \text{rayon de l'orbe de l'épicycle}&=39;30+r\\
  \text{rayon du rotateur}&=3+39;30+r=42;30+r\\
  \text{rayon du déférent}&=9+3+39;30+r=51;30+r\\
  \text{rayon extérieur de l'orbe incliné}&=60+9+3+39;30+r=111;30+r\\
  \text{rayon intérieur de l'orbe incliné}&=60-(9+3+39;30+r)=8;30-r\\
  \text{épaisseur du parécliptique}&=0;30
\end{align*}

\paragraph{Genèse du modèle en latitude pour les planètes supérieures}
Concernant les latitudes des planètes supérieures, {\shatir} s'estime satisfait du progrès qu'il a accompli. Voici la critique qu'il adresse à ses prédécesseurs~:
\begin{quote}
  ``Le rapprochement de la ceinture de l'épicycle et de la ceinture de
  l'orbe incliné, sans un mobile qui ne perturbe les mouvements en
  longitude, dans l'astronomie classique de l'\textit{Almageste} et
  des \emph{Hypothèses} concernant les corps mûs en latitude, est
  impossible. Ce dont se sont efforcés Ptolémée et ses successeurs
  comme Ibn al-Haytham dans son \textit{\'Epître}, Na\d{s}{\=\i}r
  al-\d{T}\=us{\=\i} dans sa \textit{Tadhkira}, Mu'ayyad
  al-`Ur\d{d}{\=\i} dans son \textit{Astronomie}, Qu\d{t}b
  al-Sh{\=\i}r\=az{\=\i} dans sa \uwave{\textit{Muntaha 'adw\=ar}}, et
  ce que quiconque a repris de leurs propos dans son livre [...], rien
  de tout cela ne suffit au but visé pour la longitude et la latitude,
  et quand cela convient pour l'une, cela pêche pour l'autre. Que Dieu
  pardonne à ceux qui ont ici admis leur
  faiblesse.''\footnote{\textit{Cf.} p.~\pageref{critique_lat1}
    \textit{supra}.}
\end{quote}
\`A cette critique s'en ajoutera une seconde~; mais il faut déjà comprendre la première.

Pour expliquer l'oscillation du plan de l'épicycle\footnote{Pour Saturne, selon Ptolémée, l'oscillation du plan de l'épicycle autour du plan de l'excentrique est d'amplitude $i_1+i_2=4°30'$. Il oscille autour d'un plan parallèle à l'écliptique avec une amplitude de $i_2=2°$ puisque l'excentrique est lui-même incliné de $i_1=2°30'$. \textit{Cf.} \cite{swerdlow2005} p.~43.}, Ptolémée avait décrit une sorte de modèle mécanique dans l'\textit{Almageste} en termes équivoques~:
\begin{quote}
  ``En résumé général, suivant les hypothèses, les cercles excentriques des cinq planètes se trouvent inclinés sur le plan du cercle milieu du zodiaque autour du centre du zodiaque, et cette inclinaison est constante dans Saturne, Jupiter et Mars [...]''

  ``Les diamètres apogées des épicycles qui, à \emph{certain point de départ}, étoient dans le plan de l'excentrique, sont transportés par de petits cercles fixés, pour ainsi dire, à leurs extrémités périgées, et d'un rayon propre à représenter les inégalités observées dans les latitudes. Mais ces petits cercles sont \emph{perpendiculaires aux plans des excentriques}~; ils ont leurs centres dans ces plans [...]''

  ``Les diamètres qui coupent à angles droits les diamètres apogées [...] demeurent constamment parallèles au plan de l'écliptique [...]''

  ``Au sujet de ces petits cercles qui font ainsi varier la position des épicycles, il faut d'abord remarquer ce qui suit. Ils sont partagés en deux également \emph{par les plans sur lesquels nous disions que se fait la variation de l'inclinaison} [...]''

  ``Qu'on objecte pas à ces hypothèses, qu'elles sont trop difficiles à saisir [...]''\footnote{\textit{Cf.} \cite{ptolemy1952}, chap. XIII.2. Nos italiques~; nous adoptons ici la traduction française de Halma.}
\end{quote}

\begin{figure}
  \begin{center}
    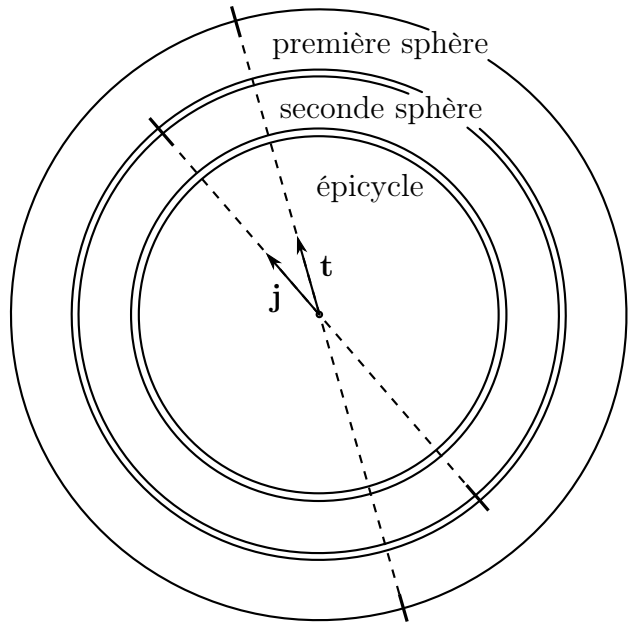
  \end{center}
  \caption{\label{fig052}Ibn al-Haytham's \emph{Iltif\=af}, according to al-\d{T}\=us{\=\i}}
\end{figure}

Peut-être à la recherche d'une \emph{description exacte} du mouvement du plan de l'épicycle, alors même qu'un \emph{calcul exact} des latitudes résultant de l'inclinaison des épicycles semblait si difficile\footnote{Ibn al-Haytham renonce à ce calcul dans sa \textit{Configuration des mouvements}, \textit{cf.} \cite{rashed2006} p.~444-447.}, al-\d{H}asan Ibn al-\d{H}asan Ibn al-Haytham a précisé le modèle de Ptolémée dans un texte aujourd'hui perdu sur \textit{Le mouvement d'enroulement}, probablement antérieur à son traité sur \textit{La configuration des mouvements}. C'est sûrement à ce texte qu'{\shatir} fait allusion. Heureusement pour nous, al-\d{T}\=us{\=\i} a décrit dans la \textit{Ta\b{d}kira} le dispositif conçu par Ibn al-Haytham\footnote{\textit{Cf.} \cite{altusi1993} II.11 [16].}~: il s'agit d'adjoindre à l'épicycle deux orbes ou sphères homocentriques comme sur la figure \ref{fig052}, où l'axe de l'épicycle est perpendiculaire au plan de la section représentée sur la figure. Le mouvement propre de l'épicycle devant être rapporté au référentiel solide constitué par la seconde sphère, c'est un diamètre de la seconde sphère qui jouera le rôle de l'``apogée moyen'' de l'épicycle. Que ce diamètre soit l'axe $\mathbf{j}$ de la seconde sphère. Le mouvement de rotation de la première sphère autour de son axe $\mathbf{t}$ aura pour effet de faire ``osciller'' la direction du vecteur $\mathbf{j}$ autour du vecteur $\mathbf{t}$. Le rôle du mouvement de rotation de la seconde sphère autour de $\mathbf{j}$ est, comme pour la troisième sphère d'un couple de \d{T}\=us{\=\i}, d'annuler le moment angulaire imprimé à l'épicycle par le mouvement de la première sphère. Le plan de l'épicycle oscillera alors autour d'un plan contenant le vecteur $\mathbf{t}$. Reste à comprendre où fixer les pôles de la première sphère dans le système d'orbes assurant les mouvements en longitude. \d{T}\=us{\=\i} reste vague~; il nous faut donc deviner.

\begin{figure}
  \begin{center}
    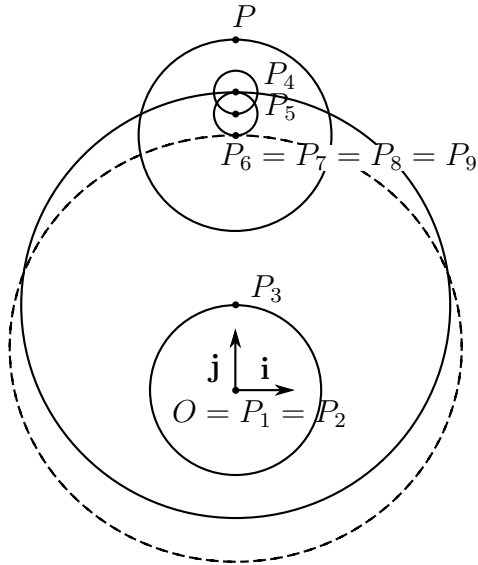
  \end{center}
  \caption{\label{fig053}Saturne d'après la \textit{Ta\b{d}kira} II.11[16]}
\end{figure}

Naïvement, plaçons l'axe de la première 1'axe de la première sphère, c'est-à-dire le vecteur $\mathbf{t}$, dans le plan de la figure initiale (notre fig. \ref{fig046} p.~\pageref{fig046}), et fixons ses pôles dans la ``sphère englobante''. Le modèle de \d{T}\=us{\=\i} auquel on adjoint les deux orbes d'Ibn al-Haytham contient maintenant neuf orbes\footnote{Voir notre fig.~\ref{fig053} p.~\pageref{fig053}, figure plane qui, bien sûr, n'est plus suffisante pour décrire la configuration initiale des orbes.}~:

-- le parécliptique de centre $P_1$,

-- l'orbe incliné en $P_2$,

-- l'orbe déférent en $P_3$,

-- une ``grande sphère'' en $P_4$, une ``petite sphère'' en $P_5$ et une ``sphère englobante'' en $P_6$, formant le couple de \d{T}\=us{\=\i},

-- les deux sphères d'Ibn al-Haytham en $P_7$ et $P_8$, où $P_8=P_7=P_6$,

-- l'épicycle en $P_9=P_8$,

-- Saturne en $P$.

Les rotations décrivant le mouvement de l'astre sont~:
$$R_{P_1,\lambda_A}
\,R_{P_2,\mathbf{u},2°30'}
\,R_{P_3,\overline{\kappa}}
\,R_{P_4,\overline{\kappa}}
\,R_{P_5,-2\overline{\kappa}}
\,R_{P_6,\overline{\kappa}}
\,R_{P_7,\mathbf{t},-(\overline{\kappa}+140°)}
\,R_{P_8,\mathbf{j},(\overline{\kappa}+140°)}
\,R_{P_9,\overline{\alpha}}$$
où $\mathbf{u}=\cos(50°)\mathbf{i}-\sin(50°)\mathbf{j}$, et
$\mathbf{t}=\sin(4°30')\mathbf{i}+\cos(4°30')\mathbf{j}$. Al-\d{T}\=us{\=\i} juge cette solution peu satisfaisante, car~:
\begin{quote}
  ``Les petits cercles mentionnés, en entraînant des inclinaisons latitudinales, causent aussi des inclinaisons en longitude, et les positions des apogées et des périgées sont donc différents de ce qu'elles devraient être quant aux points avec lesquels ils devraient être alignés.''\footnote{\textit{Cf.} \cite{altusi1993} II.11 [15].}
\end{quote}

\begin{figure}
  \begin{center}
    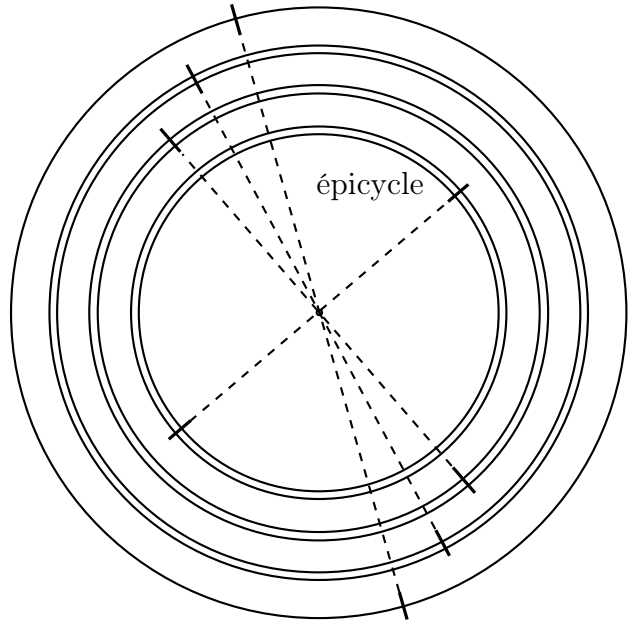
  \end{center}
  \caption{\label{fig054}Le couple curviligne de \d{T}\=us{\=\i} pour les latitudes}
\end{figure}

Dans sa première critique, {\shatir} ne fait que répéter celle qu'al-\d{T}\=us{\=\i} adresse au modèle de Ptolémée précisé par Ibn al-Haytham, et il ajoute que ni al-\d{T}\=us{\=\i}, ni `Ur\d{d}{\=\i}, ni Sh{\=\i}r\=az{\=\i} n'ont su résoudre ce problème~; pourtant al-\d{T}\=us{\=\i} propose une autre solution dans la \textit{Ta\b{d}kira}~: il propose de remplacer la première sphère d'Ibn al-Haytham par \emph{deux} sphères homocentriques (\textit{cf.} figure \ref{fig054}). L'épicycle est toujours contenu dans la seconde sphère d'Ibn al-Haytham, mais elle-même est contenue dans un \emph{petite sphère}, elle-même contenue dans une \emph{grande sphère}, elle-même contenue dans la \emph{sphère englobante} centrée en $P_6$. Leurs centres $P_6=P_7=P_8=P_9=P_{10}$ sont tous confondus. La vitesse de rotation de la petite sphère est double et en sens inverse de celle de la grande sphère~: les axes des trois sphères additionnelles étant presque parallèles, elles agiront presque comme un couple de \d{T}\=us{\=\i}.\label{notations_couple_curviligne}

Quel est le mouvement de l'apogée moyen produit par ce dispositif~? \`A la surface de la sphère englobante, la trace de l'axe de la seconde sphère d'Ibn al-Haytham est une courbe fermée qu'al-\d{T}\=us{\=\i} pensait pouvoir assimiler à un petit arc d'un grand cercle de la sphère englobante, parcouru alternativement dans un sens puis dans l'autre~: l'apogée moyen décrit bien un mouvement oscillatoire\footnote{Ragep l'explique très bien, \textit{cf.} \cite{altusi1993} fig. C26 et note 54 p.~455.}.

Comment insérer ce dispositif dans le modèle d'al-\d{T}\=us{\=\i} pour Saturne~? Comme on l'a fait pour le dispositif d'Ibn al-Haytham, on choisira naïvement d'insérer l'axe de la grande sphère dans le plan de la figure initiale, donc dans le plan de la sphère englobante, dans la direction du vecteur $\mathbf{j}$. On notera~:
\begin{align*}
  \mathbf{u}&=\cos(50°)\mathbf{i}-\sin(50°)\mathbf{j}\\
  \mathbf{v}&=\cos(4°30')\mathbf{j}-\sin(4°30')\mathbf{k}\\
  \mathbf{w}&=\cos(2°15')\mathbf{j}-\sin(2°15')\mathbf{k}
\end{align*}
où $\mathbf{v}$ est l'axe de la seconde sphère d'Ibn al-Haytham et $\mathbf{w}$ est l'axe de la ``petite sphère''. La composée de rotations décrivant le mouvement de Saturne est alors~:
$$R_{P_1,\lambda_A}
\,R_{P_2,\mathbf{u},2°30'}
\,R_{P_3,\overline{\kappa}}
\,R_{P_4,\overline{\kappa}}
\,R_{P_5,-2\overline{\kappa}}
\,R_{P_6,\overline{\kappa}}
\,T\,R_{P_{10},\mathbf{i},-4°30'}\,R_{P_{10},\overline{\alpha}}$$
où $T=R_{P_7,\mathbf{j},\overline{\kappa}+50°}\,R_{P_8,\mathbf{w},-2(\overline{\kappa}+50°)}\,R_{P_9,\mathbf{v},\overline{\kappa}+50°}$ est le ``couple curviligne'' des trois sphères additionnelles.

\begin{figure}
  \begin{center}
    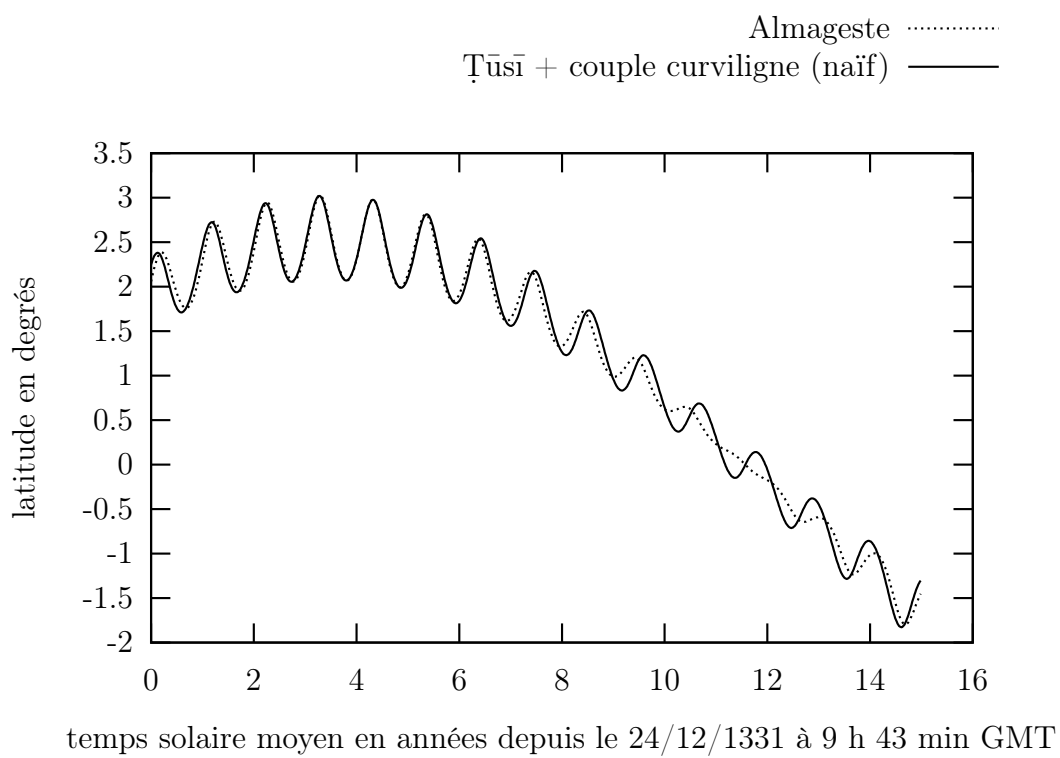
  \end{center}
  \caption{\label{fig055}Interprétation naïve de la solution d'al-\d{T}\=us{\=\i} utilisant son ``couple curviligne'' pour les latitudes de Saturne}
\end{figure}

Pour Saturne dont l'épicycle est petit, la solution d'Ibn al-Haytham et celle du couple curviligne d'al-\d{T}\=us{\=\i} se valent y compris pour les longitudes. Pourtant elles ont toutes deux un défaut notable~: l'épicycle garde une inclinaison par rapport à l'écliptique près des n{\oe}uds, causant les larges oscillations en latitude entre $t=10$ et $t=12$ sur le graphe de la figure \ref{fig055}. Et pour cause~: nous avons, ``naïvement'', fixé les pôles de la plus grande des deux ou trois sphères additionnelles dans le plan de l'orbe englobant, lui-même coïncidant avec le plan de l'excentrique de Ptolémée. L'épicycle va donc osciller autour de l'excentrique~: aux n{\oe}uds, il sera confondu avec le plan de l'excentrique, dont l'inclinaison est de $2°30'$ par rapport à l'écliptique. Il y a donc erreur.

Mais où devrait-on attacher les pôles de la grande sphère dans la sphère englobante, si ce n'est dans le plan de cet orbe~? Al-\d{T}\=us{\=\i} est plutôt vague. Selon lui, les pôles de la première sphère d'Ibn al-Haytham sont à une distance des pôles de sa seconde sphère ``égale à la déviation maximale du diamètre [de l'apogée moyen] de la planète par rapport \emph{au plan dans lequel elle n'a pas de latitude}''\footnote{\textit{Ta\b{d}kira} \cite{altusi1993} II.11[16]. On remarque qu'en décrivant le modèle mécanique de Ptolémée, al-\d{T}\=us{\=\i} est, sur ce point, à peu près aussi vague que Ptolémée~: le modèle résulte selon al-\d{T}\=us{\=\i} en un ``déplacement des extrémités des diamètres de l'épicycle \emph{hors des plans dans lesquels ils n'ont pas d'inclinaison}'' (\textit{Ta\b{d}kira} II.11[14]).}. Il faut donc penser que les axes des deux sphères d'Ibn al-Haytham devraient être inclinés de $2°$ l'un par rapport à l'autre, puisque l'épicycle est incliné de $2°$ au plus par rapport à l'écliptique. Si l'on pouvait alors fixer les pôles de la première sphère dans un plan parallèle à l'écliptique, le modèle produirait une oscillation du plan de l'épicycle de $\pm 2°$ autour d'un plan parallèle à l'écliptique, et il serait confondu avec le plan de l'écliptique aux n{\oe}uds. Hélas, dans le solide de la sphère englobante, il n'existe aucun plan qui reste constamment parallèle à l'écliptique quand les orbes se meuvent~!

Une rotation spatiale n'agit pas seulement sur les directions, elle agit aussi sur les plans affines. {\shatir} l'avait bien compris, car c'est là l'objet de sa seconde critique. Il s'est posé la question de savoir autour de quel plan l'épicycle doit osciller. Selon lui~:
\begin{quote}
  ``L'observation confirme que c'est le plan de l'écliptique et non le
  plan de l'orbe incliné. La plupart des Modernes -- Na\d{s}{\=\i}r
  al-\d{T}\=us{\=\i}, Mu'ayyad al-`Ur\d{d}{\=\i}, Qu\d{t}b
  al-Sh{\=\i}r\=az{\=\i} -- ont pensé que la latitude s'annulait quand
  c'est le plan de l'orbe incliné: c'est impossible
  [...]''\footnote{\textit{Cf.} p.~\pageref{critique_lat2}
    \textit{supra}.}
\end{quote}
En effet, la solution dont nous avons tracé le graphe fig. \ref{fig055}, bien que naïve, est la seule possible, à moins de modifier radicalement le modèle d'al-\d{T}\=us{\=\i}. Pour rendre possible l'oscillation autour d'un plan parallèle à l'écliptique, il faudrait modifier les autres orbes pour garantir l'existence d'un tel plan dans le solide de l'orbe englobant. Il faudrait par exemple incliner, une fois pour toutes, le plan et l'axe de l'orbe englobant par rapport à ceux du déférent. On aurait alors la composée de rotations suivante, toujours en accord avec les notations de la p.~\pageref{notations_couple_curviligne} ci-dessus~:
$$R_{P_1,\lambda_A}
\,R_{P_2,\mathbf{u},2°30'}
\,R_{P_3,\overline{\kappa}}
\,R_{P_4,\overline{\kappa}}
\,R_{P_5,-2\overline{\kappa}}
\,R_{P_6,\mathbf{u},-2°30'}
\,R_{P_6,\overline{\kappa}}
\,T\,R_{P_{10},\mathbf{i},-4°30'}\,R_{P_{10},\overline{\alpha}}$$
où, comme ci-dessus,
$$T=R_{P_7,\mathbf{j},\overline{\kappa}+50°}
\,R_{P_8,\mathbf{w},-2(\overline{\kappa}+50°)}
\,R_{P_9,\mathbf{v},\overline{\kappa}+50°},$$
mais ici $\mathbf{v}=\cos(2°)\mathbf{j}-\sin(2°)\mathbf{k}$ et $\mathbf{w}=\cos(1°)\mathbf{j}-\sin(1°)\mathbf{k}$. Une telle modification rend l'adéquation avec les prédictions ptoléméennes presque parfaite\footnote{Une correction semblable s'appliquerait aussi aisément au dispositif d'Ibn al-Haytham de la fig.~\ref{fig052}.}.

\begin{figure}
  \begin{center}
    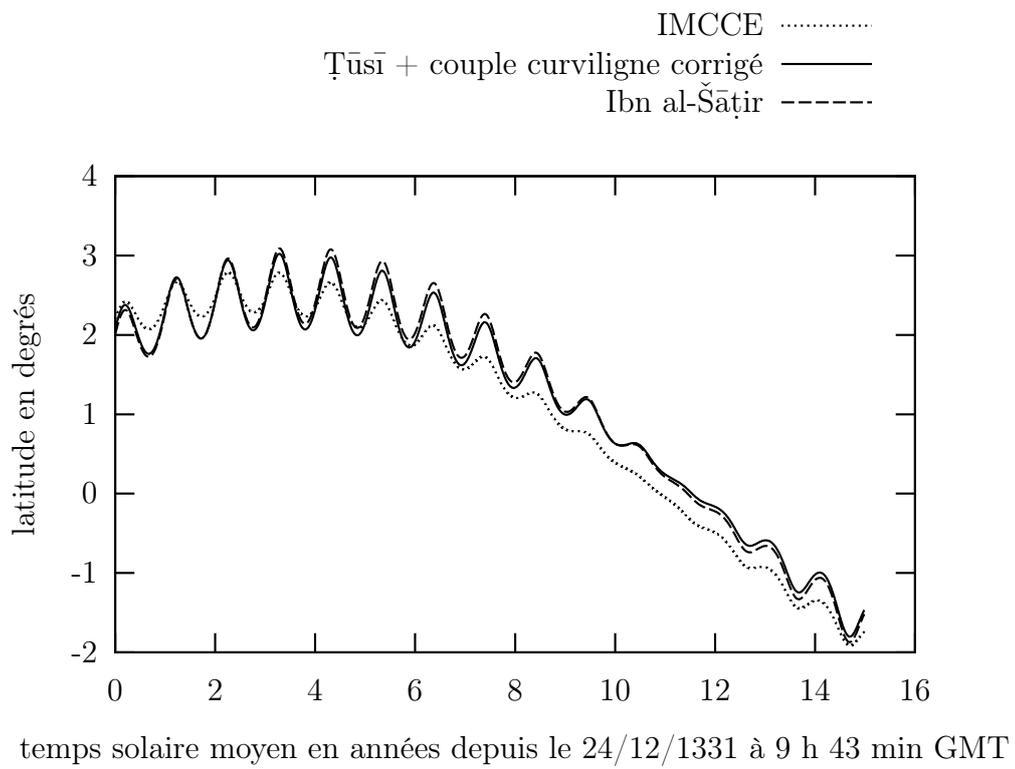
  \end{center}
  \caption{\label{fig056}Comparaison des modèles en latitude pour Saturne}
\end{figure}

Comme on l'a vu plus haut, {\shatir} propose une autre solution, adaptée à son propre modèle en longitude, et consistant simplement à incliner de manière judicieuse les plans et les axes des orbes déjà utilisés. Il n'est besoin d'adjoindre aucun orbe additionnel pour le mouvement en latitude~; cinq orbes suffisent là où al-\d{T}\=us{\=\i} en utilise dix. Son modèle n'est pas tout à fait équivalent à celui d'al-\d{T}\=us{\=\i}, mais il est aussi très proche des prédictions ptoléméennes~; il est donc plus intéressant de comparer les deux modèles aux éphémérides modernes comme nous l'avons fait à la figure \ref{fig056}.

Qu'en conclure~? Au moins dès le XIème siècle avec Ibn al-Haytham, les astronomes avaient engagé l'étude des composées de rotations spatiales. {\shatir} a assimilé les travaux concernant les composées de rotations homocentriques à axes non parallèles (les sphères homocentriques utilisées pour décrire les latitudes), aussi bien que les travaux de l'école de Maragha décrivant les mouvements en longitude par des composées de rotations à axes parallèles. C'est peut-être une meilleure appréhension du concept d'orbe solide, et une compréhension plus fine de l'action des rotations spatiales sur les plans affines, qui lui permettent d'élucider les défauts des solutions d'Ibn al-Haytham et d'al-\d{T}\=us{\=\i}, et d'y remédier.

Un indice montre que la recherche d'{\shatir} s'appuie bien sur les travaux de ses prédécesseurs. Après avoir décrit son modèle pour les latitudes des planètes supérieures, {\shatir} discute qualitativement les effets produits par les inclinaisons des axes des orbes dans la figure initiale, sur les inclinaisons aux n{\oe}uds et aux limites des latitudes. Pour ce faire, il est conduit à identifier $P_3$ et $P_4$ dans ses explications, et sur une de ses figures. On pourra comparer la figure originale p.~\pageref{lat_inf_profil} à notre fig.~\ref{fig041} p.~\pageref{fig041} où nous distinguons $P_3$ et $P_4$. Eu égard au fait que $P_3P_4$ est petit devant les rayons de l'épicycle et de l'excentrique, cette approximation peut se justifier. Mais il s'agit sans doute d'un choix didactique ou d'un résidu d'une réflexion antérieure sur un modèle plus simple. L'usage de l'approximation $P_3\simeq P_4$ montre que les composées de rotations homocentriques avaient déjà valeur de paradigme dans les théories planétaires. C'est en inclinant les axes de plusieurs orbes non homocentriques qu'{\shatir} semble accomplir un geste nouveau.
\label{end_saturne}

\paragraph{Vénus~: figure initiale} La figure initiale pour Vénus est analogue à celle de Saturne (\textit{cf.} fig.~\ref{fig040} p.~\pageref{fig040}), sauf quant aux rayons des orbes et à la position de la ligne des n{\oe}uds. La ligne des n{\oe}uds, à l'intersection des plans de l'orbe incliné et du parécliptique, est ici orientée dans la direction du vecteur $\mathbf{i}$, le n{\oe}ud ascendant étant du même côté que $\mathbf{i}$. Posons $\overrightarrow{OP_3}=60\,\mathbf{j}$, alors
$$\overrightarrow{P_3P_4}=1;41\,\mathbf{j},\qquad
\overrightarrow{P_4P_5}=-0;26\,\mathbf{j},\qquad
\overrightarrow{P_5P}=43;33\,\mathbf{j}.$$

\begin{figure}
  \begin{center}
    \footnotesize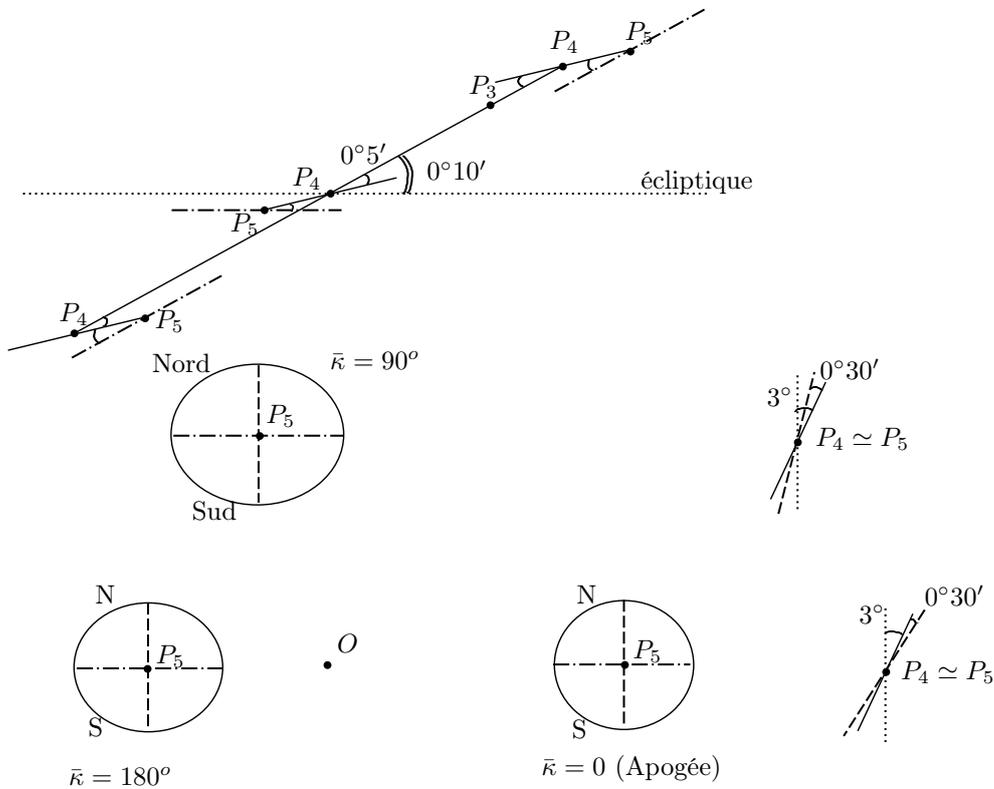
  \end{center}
  \caption{\label{fig057}Vénus, trois positions~: projections orthogonales}
\end{figure}

\paragraph{Vénus~: transformations géométriques}
Les rotations utilisées pour modéliser le mouvement de Vénus en longitude sont analogues à celles des planètes supérieures. En revanche la théorie des latitudes exposée dans le chapitre 25 de la \textit{Nih\=aya} est assez différente. Voici la liste complète des rotations utilisées dans les chapitres 19 (pour les longitudes) et 25 (pour les latitudes)~:
$$R_{P_5,\overline{\alpha}-\overline{\kappa}},\quad R_{P_5,\mathbf{j},0°30'},
\quad R_{P_5,\mathbf{i},0°5'},$$
$$R_{P_4,2\overline{\kappa}},\quad R_{P_4,\mathbf{j},3°},
\quad R_{P_4,\mathbf{i},-0°5'},$$
$$R_{P_3,-\overline{\kappa}},\quad R_{P_2,\overline{\kappa}},\quad
R_{P_2,\mathbf{i},0°10'},\quad R_{P_1,\lambda_A}.$$
L'image du point $P_3$ entraîné par les les mouvements de l'orbe parécliptique et de l'orbe incliné est~:
$$R_{P_1,\lambda_A}\,R_{P_2,\mathbf{i},0°10'}
\,R_{P_2,\overline{\kappa}}(P_3).$$
Quant au point $P_4$, il est aussi entraîné par le mouvement de l'orbe déférent et devient~:
$$R_{P_1,\lambda_A}\,R_{P_2,\mathbf{i},0°10'}
\,R_{P_2,\overline{\kappa}}\,R_{P_3,-\overline{\kappa}}(P_4).$$
Le point $P_5$ est aussi entraîné par le mouvement de l'orbe rotateur et devient~:
$$R_{P_1,\lambda_A} 
\,R_{P_2,\mathbf{i},0°10'} 
\,R_{P_2,\overline{\kappa}}
\,R_{P_3,-\overline{\kappa}}
\,R_{P_4,\mathbf{i},-0°5'}
\,R_{P_4,\mathbf{j},3°}
\,R_{P_4,2\overline{\kappa}}(P_5).$$
Enfin le point $P$, aussi entraîné par l'orbe de l'épicycle, devient~:
$$R_{P_1,\lambda_A} 
\,R_{P_2,\mathbf{i},0°10'} 
\,R_{P_2,\overline{\kappa}}
\,R_{P_3,-\overline{\kappa}}
\,R_{P_4,\mathbf{i},-0°5'}
\,R_{P_4,\mathbf{j},3°}
\,R_{P_4,2\overline{\kappa}}
\,R_{P_5,\mathbf{i},0°5'}
\,R_{P_5,\mathbf{j},0°30'}
\,R_{P_5,\overline{\alpha}-\overline{\kappa}}(P).$$
La figure \ref{fig057} montre l'effet des rotations inclinant les plans des orbes, pour trois positions. On a projeté orthogonalement sur le plan du parécliptique (figure du bas), sur un plan orthogonal à la ligne des n{\oe}uds (figure du haut), et sur un plan orthogonal au plan du parécliptique mais parallèle à la ligne des n{\oe}uds (figure de droite). Sur la figure du bas, on a seulement représenté la ceinture de l'épicycle~; dans les deux autres figures, le petit segment en trait plein de centre $P_4$ est la ceinture du rotateur. Sur la figure de droite, on a confondu $P_4$ et $P_5$, par approximation. La ceinture du déférent reste toujours dans le plan de l'orbe incliné.

\paragraph{Vénus~: trajectoire paramétrée} Le modèle de Vénus est couplé au modèle du Soleil en vertu de la relation suivante~:
$$\overline{\kappa}=\overline{\lambda}_{\astrosun}-\lambda_A.$$
La trajectoire dans l'ensemble des valeurs des paramètres est~:
$$\lambda_A=\dot{\lambda}_At+\lambda_A(0),$$
$$\overline{\lambda}_{\astrosun}=\dot{\overline{\lambda}}_{\astrosun}t+\overline{\lambda}_{\astrosun}(0),$$
$$\overline{\alpha}=\dot{\overline{\alpha}}t+\overline{\alpha}(0).$$
On a, d'après {\shatir}~:
$$\overline{\lambda}_{\astrosun}(0)=280;9,0,\quad\dot{\overline{\lambda}}_{\astrosun}=359;45,40\text{ par année persane},$$
$$\lambda_A(0)=77;52,\quad\dot{\lambda}_A=0;1\text{ par année persane},$$
$$\overline{\alpha}(0)=320;50,19,\quad\dot{\overline{\alpha}}=225;1,48,41\text{ par année persane}.$$

\begin{figure}
  \begin{center}
    \footnotesize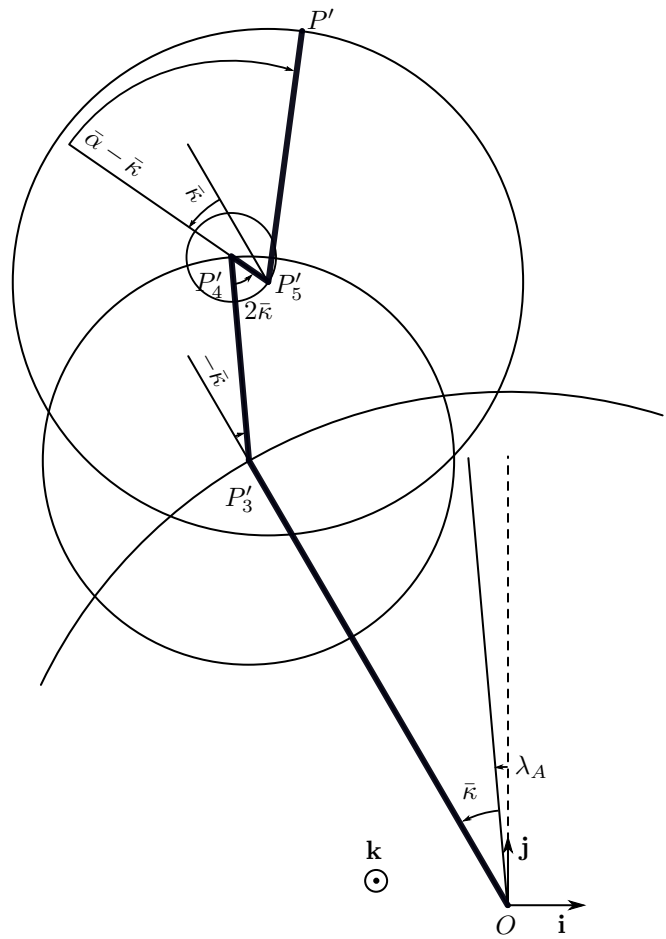
  \end{center}
  \caption{\label{fig058}Les orbes de Vénus à un instant $t$}
\end{figure}

\paragraph{Transformations planes et équations de Vénus}
Comme nous l'avons expliqué pour Saturne, on peut réécrire à droite toutes les rotations dont l'axe est dans la direction du vecteur $\mathbf{k}$, en appliquant des relations de commutation. Il existe donc une composée de rotations $M$ telle que l'image du point $P$ soit~:
$$M\,R_{P_2,\mathbf{a},0°10'}(P')$$
où $\mathbf{a}$ est l'image de $\mathbf{i}$ par la rotation $R_{P_1,\lambda_A}$, et où
$$P'=R_{P_1,\lambda_A}
\,R_{P_2,\overline{\kappa}}
\,R_{P_3,-\overline{\kappa}}
\,R_{P_4,2\overline{\kappa}}
\,R_{P_5,\overline{\alpha}-\overline{\kappa}}(P).$$
On introduit de même les points $P'_3$, $P'_4$, $P'_5$ suivants (\textit{cf.} fig. \ref{fig058}) exactement comme pour Saturne~:
$$P_3'=R_{P_1,\lambda_A}\,R_{P_2,\overline{\kappa}}(P_3),$$
$$P_4'=R_{P_1,\lambda_A}\,R_{P_2,\overline{\kappa}}\,R_{P_3,-\overline{\kappa}}(P_4),$$
$$P_5'=R_{P_1,\lambda_A}\,R_{P_2,\overline{\kappa}}\,
R_{P_3,-\overline{\kappa}}\,R_{P_4,2\overline{\kappa}}(P_5).$$
Les équations pour le mouvement en longitude de Vénus sont donc identiques à celles des planètes supérieures (\textit{cf.} p.~\pageref{equ_sat} ci-dessus).

\paragraph{Vénus~: trigonométrie sphérique} On va à présent calculer les coordonnées du point $R_{P_2,\mathbf{a},0°10'}(P')$ par rapport à l'écliptique. L'angle formé entre le vecteur $\mathbf{a}$ et la direction du point $P'$ vaut (modulo $360°$)~:
$$(\mathbf{a},\overrightarrow{OP'})=\overline{\kappa}+c_1+c_2+90°=\overline{\lambda}-\lambda_{\ascnode}+c_1+c_2=\lambda-\lambda_{\ascnode}$$
où $\lambda=\overline{\lambda}+c_1+c_2$ et $\overline{\lambda}=\overline{\kappa}+\lambda_A$, et où $\lambda_{\ascnode}=\lambda_A-90°$ est la longitude du n{\oe}ud ascendant sur l'écliptique, toujours situé $90°$ avant l'Apogée. Comme on l'a montré pour Saturne, on trouve que la longitude du point $R_{P_2,\mathbf{a},0°10'}(P')$ par rapport à l'écliptique, en prenant la direction du point vernal, c'est-à-dire $\mathbf{j}$, comme origine, vaut~:
$$\lambda+e_n(\lambda-\lambda_{\ascnode})$$
où l'``équation du déplacement'' $e_n(x)$ est définie comme suit sur l'intervalle $\lbrack -90°,270°\rbrack$, et ailleurs par périodicité~:
$$e_n(x)=\left\lbrace\begin{array}{l}
\arctan(\cos(0°10')\times\tan x)-x,\text{ si }x\in\rbrack -90°,90°\lbrack\\
180°+\arctan(\cos(0°10')\times\tan x)-x,\text{ si }x\in\rbrack 90°,270°\lbrack\\
0°,\text{ si }x=\pm 90°.
\end{array}\right.$$
Sa latitude est~:
$$\arcsin(\sin(0°10')\times\sin(\lambda-\lambda_{\ascnode})).$$

\paragraph{Vénus~: les inclinaisons des petits orbes}
Les rotations qui composent $M$ ont leurs axes contenus dans le plan de l'orbe incliné et sont d'angles petits ($3°$ au plus). Elles auront peu d'effet\footnote{\textit{Cf.} \cite{penchevre2016} où nous montrons que l'erreur commise en les négligeant est inférieure à $0°5'$.}. Dans le chapitre 25, {\shatir} se contente d'une description qualitative de l'effet en latitude de $M$.

\paragraph{Vénus~: second modèle} Le modèle que nous avons restitué ci-dessus est le premier modèle décrit par {\shatir} dans le chapitre 25. On voit que les inclinaisons des plans des orbes sont conçues dans l'idée de reproduire les effets décrits dans l'\textit{Almageste}. L'inclinaison de l'``orbe incliné'', variable dans l'\textit{Almageste} mais constante chez {\shatir}, minime, n'a que peu d'influence sur le résultat final. {\shatir} envisage ensuite un second modèle plus proche de la théorie exposée par Ptolémée dans ses \textit{Hypothèses planétaires}. La composée de rotations utilisée est alors la suivante~:
$$R_{P_1,\lambda_A} 
\,R_{P_2,\mathbf{i},0°10'} 
\,R_{P_2,\overline{\kappa}}
\,R_{P_3,-\overline{\kappa}}
\,R_{P_4,\mathbf{i},-0°5'}
\,R_{P_4,\mathbf{j},3°30'}
\,R_{P_4,2\overline{\kappa}}
\,R_{P_5,\mathbf{i},0°5'}
\,R_{P_5,\overline{\alpha}-\overline{\kappa}}.$$

\begin{figure}
  \begin{center}
    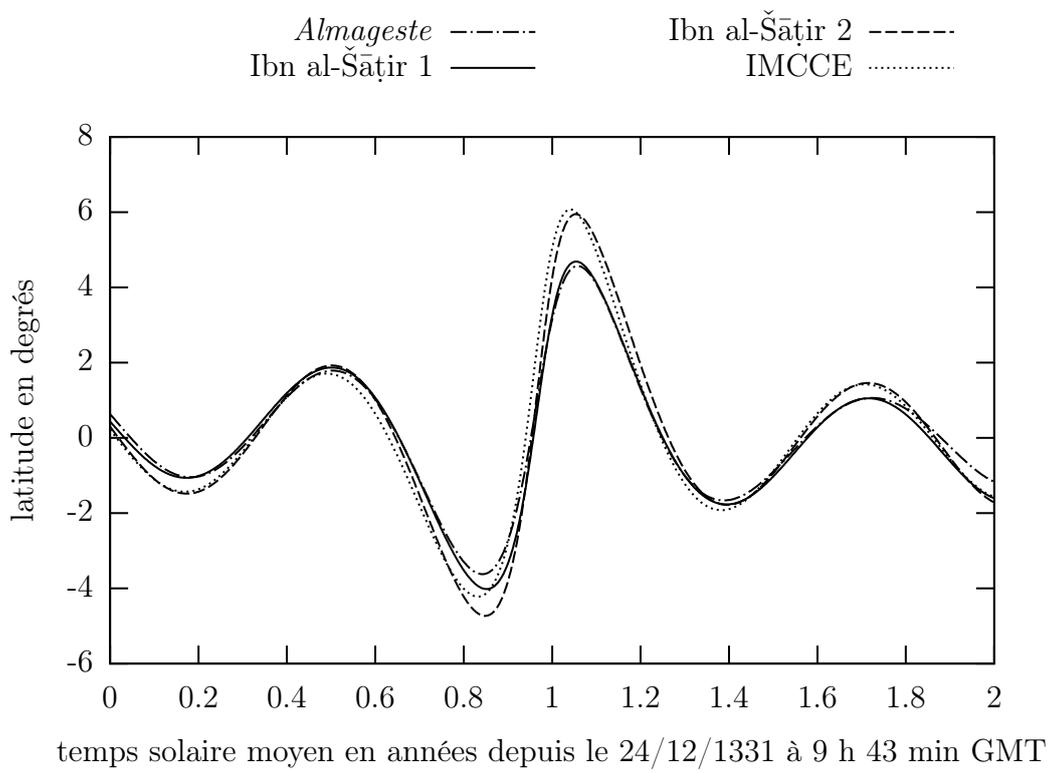
  \end{center}
  \caption{\label{fig059}Vénus~: comparaison des modèles pour les latitudes}
\end{figure}

\paragraph{Comparaison} \`A la figure \ref{fig059}, nous donnons les latitudes de Vénus pendant deux ans, à partir de l'\'Epoque choisie par {\shatir}, calculées de quatre manières différentes~: (1) par l'IMCCE de l'Observatoire de Paris, (2) en suivant la méthode strictement ptoléméenne de l'\textit{Almageste} mais avec les paramètres d'{\shatir} à l'\'Epoque et ses mouvements moyens, (3) au moyen du premier modèle d'{\shatir} censé reproduire les valeurs de l'\textit{Almageste}, et (4) au moyen du second modèle d'{\shatir} imitant les \textit{Hypothèses planétaires}. 

\paragraph{Les orbes solides de Vénus}
La figure est la même que pour les planètes supérieures (figure \ref{fig049}).
\begin{align*}
  \text{rayon du globe planétaire}&=r\\
  \text{rayon de l'orbe de l'épicycle}&=43;33+r\\
  \text{rayon du rotateur}&=0;26+43;33+r=43;59+r\\
  \text{rayon du déférent}&=1;41+0;26+43;33+r=45;40+r\\
  \text{rayon extérieur de l'orbe incliné}&=60+1;41+0;26+43;33+r=105;40+r\\
  \text{rayon intérieur de l'orbe incliné}&=60-(1;41+0;26+43;33+r)=14;20-r\\
  \text{épaisseur du parécliptique}&=0;20
\end{align*}

\paragraph{Le mouvement en longitude de Vénus selon al-\d{T}\=us{\=\i}} Le modèle de Na\d{s}ir al-D{\=\i}n al-\d{T}\=us{\=\i} pour les longitudes de Vénus est aussi analogue à son modèle pour les planètes supérieures\footnote{\textit{Cf.} ci-dessus fig.~\ref{fig046}. Quitte à bien choisir les paramètres, l'équivalence cinématique entre les modèles d'al-\d{T}\=us{\=\i}, al-`Ur\d{d}{\=\i}, al-Sh{\=\i}r\=az{\=\i} et {\shatir}, démontrée p.~\pageref{equiv_sup} \textit{supra}, vaut aussi pour Vénus.}~: sept orbes, dont les centres sont déterminés par les relations suivantes\footnote{\textit{Cf.} \cite{altusi1993} p.~457.},
$$\overrightarrow{P_7P}=43;10\,\mathbf{j},\quad\overrightarrow{P_5P_6}=\overrightarrow{P_4P_5}=-0;37,30\,\mathbf{j},\text{ d'où }\overrightarrow{P_4P_6}=-1;15\,\mathbf{j},$$
$$\overrightarrow{P_3P_4}=60\,\mathbf{j},\qquad\overrightarrow{P_2P_3}=2;30\,\mathbf{j}.$$
Al-\d{T}\=us{\=\i} choisit ces paramètres de façon à reproduire la trajectoire prédite par le modèle de l'\textit{Almageste}. Ainsi $\vert P_2P_3\vert=2e$ et $\vert P_4P_6\vert=e$ où $e=1;15$ est l'excentricité chez Ptolémée, et al-\d{T}\=us{\=\i} garde le rayon de l'épicycle $43;10$ utilisé par Ptolémée.

Au contraire, le modèle d'al-\d{T}\=us{\=\i} serait équivalent au modèle d'{\shatir} si l'on avait~:
$$\overrightarrow{P_7P}=43;33\,\mathbf{j},\quad\overrightarrow{P_5P_6}=\overrightarrow{P_4P_5}=-0;26\,\mathbf{j},\text{ d'où }\overrightarrow{P_4P_6}=-0;52\,\mathbf{j},$$
$$\overrightarrow{P_3P_4}=60\,\mathbf{j},\qquad\overrightarrow{P_2P_3}=2;7\,\mathbf{j}.$$
Sous cette hypothèse, $\vert P_2P_3\vert-\vert P_4P_6\vert=1;15$ serait bien égal à l'excentricité de Vénus selon Ptolémée, mais $\dfrac{\vert P_2P_3\vert}{2}=1;3,30$ serait égal à la demi-excentricité du Soleil \emph{selon {\shatir}}\footnote{Rappelons en effet que dans le modèle du Soleil d'{\shatir}, si l'on suit les notations de la fig.~\ref{fig003} p.~\pageref{fig003} \textit{supra}, le rayon de l'épicycle apparent minimal est $\vert P_3P_4\vert-\vert P_4P\vert=2;7$.}. Al-\d{T}\=us{\=\i} lui-même indiquait que l'excentricité de Vénus est la moitié de celle du Soleil, et que l'excentricité du Soleil valait $2;30$ pour Ptolémée mais environ $2;5$ pour des observateurs plus récents\footnote{\textit{Cf.} \cite{altusi1993} p.~146 et 182. Dans \cite{ghanem1976} p.~64, Kennedy donne les valeurs adoptées par trois autres savants pour l'excentrité de Vénus~: al-B{\=\i}r\=un{\=\i} $1;2,30$, al-Zarqulla $1;3,22$, al-K\=ash{\=\i} $1;3$. Elles sont toutes proches des $1;3,30$ d'{\shatir}.}.

\begin{figure}
  \begin{center}
    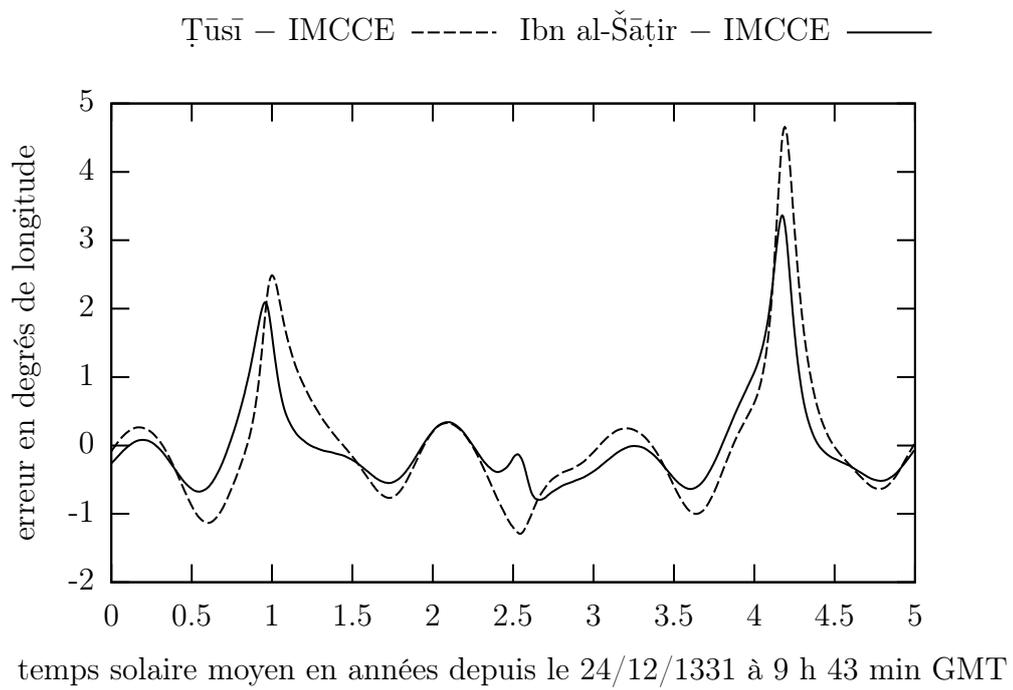
  \end{center}
  \caption{\label{fig060}Erreur en longitude des modèles d'{\shatir} et d'al-\d{T}\=us{\=\i} pour Vénus}
\end{figure}

\paragraph{Comparaison} Pour mieux juger des paramètres adoptés par {\shatir} et al-\d{T}\=us{\=\i}, on donne fig.~\ref{fig060} les longitudes de Vénus pendant cinq ans, à partir de l'\'Epoque choisie par {\shatir}~; comme l'équation de Vénus a une amplitude de l'ordre de $50°$ en coordonnées géocentriques, on a préféré tracer l'\emph{erreur} en degrés de longitude par rapport à la longitude prédite par le serveur d'éphémérides de l'IMCCE~: (1) en adoptant le modèle d'al-\d{T}\=us{\=\i} avec $\vert P_2P_3\vert=2;30$ et $\vert P_7P\vert=43;10$ et les mouvements moyens d'{\shatir}, (2) en adoptant le modèle d'{\shatir}. Dans l'absolu, le choix d'{\shatir} est un peu meilleur. Celui d'al-\d{T}\=us{\=\i} reproduit de très près les prédictions ptoléméennes de l'\textit{Almageste}\footnote{Nous avons omis la courbe représentant la théorie de l'\textit{Almageste}, car à cette échelle elle est indiscernable de la courbe représentant le modèle d'al-\d{T}\=us{\=\i}.}.

\paragraph{Mercure selon {\shatir}}\label{begin_mercure} La théorie d'{\shatir} pour Mercure ressemble à celle de Vénus, à trois différences près~:

-- Dans la figure initiale, la position relative de l'épicycle au sein du rotateur n'est pas la même~: pour Mercure, $\overrightarrow{P_4P_5}$ est orienté dans le même sens que $\mathbf{j}$. Cette particularité n'est d'ailleurs pas indiquée explicitement dans le texte du chapitre 21~; il faut regarder les figures pour le voir, et bien des lecteurs ont dû s'y méprendre.

-- Il y a deux orbes additionnels, l'\emph{orbe englobant} et l'\emph{orbe protecteur}, portés par l'orbe de l'épicycle, agissant comme un couple de \d{T}\=us{\=\i}, et dont le rôle est de faire varier la distance entre l'astre et le centre de l'épicycle (on notera $P_6$ et $P_7$ les centres de ces deux orbes additionnels).

-- Dans le chapitre 25 sur les latitudes, les inclinaisons des plans des orbes sont de signes contraires à ceux de Vénus.

\begin{figure}
  \begin{center}
    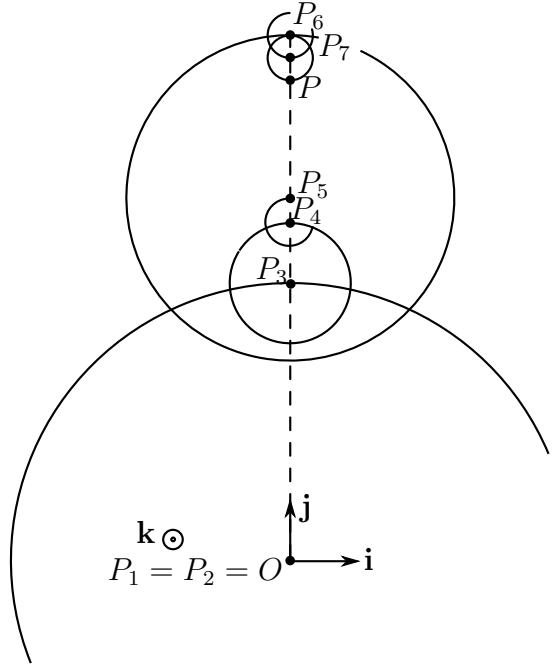
  \end{center}
  \caption{\label{fig061}Mercure~: figure initiale}
\end{figure}

Sur la figure initiale fig.~\ref{fig061}, les centres des orbes sont
définis par~:
$$O=P_1=P_2,\quad \overrightarrow{P_2P_3}=60\,\mathbf{j},\quad
\overrightarrow{P_3P_4}=4;5\,\mathbf{j},$$
$$\overrightarrow{P_4P_5}=0;55\,\mathbf{j},\quad
\overrightarrow{P_5P_6}=22;46\,\mathbf{j},\quad
\overrightarrow{P_6P_7}=\overrightarrow{P_7P}=-0;33\,\mathbf{j}.$$
Comme pour Vénus, {\shatir} présente deux modèles différents pour Mercure dans le chapitre 25. Le premier, censé reproduire les mouvements en latitude décrits dans l'\textit{Almageste}, consiste en la composée de rotations suivante~:
$$\hspace{-2cm}R_{P_1,\lambda_A}
\,R_{P_2,\mathbf{i},-0°10'}
\,R_{P_2,\overline{\kappa}}
\,R_{P_3,-\overline{\kappa}}
\,R_{P_4,\mathbf{i},0°5'}
\,R_{P_4,\mathbf{j},-6°36'30''}
\,R_{P_4,2\overline{\kappa}}
\,R_{P_5,\mathbf{i},-0°5'}
\,R_{P_5,\mathbf{j},-0°22'30''}
\,R_{P_5,\overline{\alpha}-\overline{\kappa}}
\,R_{P_6,2\overline{\kappa}}
\,R_{P_7,-4\overline{\kappa}}.$$
On a, d'après {\shatir}~:
$$\overline{\lambda}_{\astrosun}(0)=280;9,0,\quad\dot{\overline{\lambda}}_{\astrosun}=359;45,40\text{ par année persane},$$
$$\lambda_A(0)=212;52,\quad\dot{\lambda}_A=0;1\text{ par année persane},$$
$$\overline{\alpha}(0)=154;2,\quad\dot{\overline{\alpha}}=1133;57,1\text{ par année persane}.$$

Kennedy\footnote{\textit{Cf.} \cite{ghanem1976} p.~65.} a indiqué deux incohérences, concernant les rayons et les distances, que nous croyons pouvoir résoudre au vu de l'ensemble~:

-- Le ``rayon de l'épicycle apparent'' minimal est $23;52-4\times 0;33=21;40$, alors que les manuscrits indiquent $21;20$ dans le chapitre 21. La ressemblance entre \textit{thulth} ($0;20$) et \textit{thulthay} ($0;40$) a probablement induit un copiste en erreur. De plus, au chapitre 23, {\shatir} écrit bien $21;40$. Nous avons donc opté pour $21;40$, et nous avons laissé $21;20$ dans l'apparat critique.

-- Le texte donne $\vert P_4P_5\vert=0;50$ et non $0;55$. Les valeurs annoncées par {\shatir} pour les distances minimale et maximale de Mercure au centre du Monde quand $\overline{\kappa}=0°$ ou $180°$ présentent une erreur de $0;5$. Kennedy indique que cette erreur disparaît si l'on suppose $\vert P_4P_5\vert=0;55$ ou bien $\vert P_3P_4\vert=4;10$. En fait, les calculs faits par {\shatir} concernant les rayons des orbes solides révèlent qu'il avait lui-même adopté $\vert P_4P_5\vert=0;55$.

\paragraph{\'Equations de Mercure et formule d'interpolation}\label{equ_mercure}
On déduit les équations de Mercure comme on l'a fait pour les planètes supérieures et Vénus~; il faut seulement penser au sens de $\overrightarrow{P_4P_5}$, et surtout remplacer $P_5P$ par le ``rayon de l'épicycle apparent'' $P'_5P'$ que nous calculerons ensuite~:
$$c_1=-\arcsin\left(\frac{P_3P_4\sin\overline{\kappa}-P_4P_5\sin\overline{\kappa}}{OP_5'}\right)$$
où $OP_5'=\sqrt{(P_3P_4\sin\overline{\kappa}-P_4P_5\sin\overline{\kappa})^2+(OP_3+P_3P_4\cos\overline{\kappa}+P_4P_5\cos\overline{\kappa})^2}$.
$$c_2=\arcsin\left(\frac{P_5'P'\sin\alpha}{OP'}\right)$$
où $OP'=\sqrt{(P_5'P'\sin\alpha)^2+(OP_5'+P_5'P'\cos\alpha)^2}$. La formule d'interpolation est~:
$$c_2(\overline{\kappa},\alpha)\simeq c_2(0,\alpha)+\chi(\overline{\kappa})(c_2(180°,\alpha)-c_2(0,\alpha)),$$
$$\chi(\overline{\kappa})=\frac{\max\vert c_2(\overline{\kappa},\cdot)\vert-\max\vert c_2(0,\cdot)\vert}{\max\vert c_2(180°,\cdot)\vert-\max\vert c_2(0,\cdot)\vert},$$
$$\max\vert c_2(\overline{\kappa},\cdot)\vert=\arcsin\frac{P_5'P'}{OP_5'}.$$
Pour calculer $P_5'P'$, on rappelle que les points $P_5'$ et $P'$ sont définis par les transformations planes à appliquer aux points $P_5$ et $P$~:
$$P_5'=R_{P_1,\lambda_A}\,R_{P_2,\overline{\kappa}}\,
R_{P_3,-\overline{\kappa}}\,R_{P_4,2\overline{\kappa}}(P_5),$$
$$P'=R_{P_1,\lambda_A}
\,R_{P_2,\overline{\kappa}}
\,R_{P_3,-\overline{\kappa}}
\,R_{P_4,2\overline{\kappa}}
\,R_{P_5,\overline{\alpha}-\overline{\kappa}}
\,R_{P_6,2\overline{\kappa}}
\,R_{P_7,-4\overline{\kappa}}(P).$$
On a donc
\begin{align*}
  \overrightarrow{P_5'P'}&=\overrightarrow{P_5R_{P_5,\overline{\alpha}-\overline{\kappa}}R_{P_6,2\overline{\kappa}}R_{P_7,-4\overline{\kappa}}(P)}\\
  &=R_{\overline{\alpha}-\overline{\kappa}}(\overrightarrow{P_5R_{P_6,2\overline{\kappa}}R_{P_7,-4\overline{\kappa}}(P)})
\end{align*}
La remarque p.~\pageref{couple_transl} montre que le vecteur $\overrightarrow{P_5R_{P_6,2\overline{\kappa}}R_{P_7,-4\overline{\kappa}}(P)}$ est colinéaire au vecteur $\mathbf{j}$ et que
$$R_{P_6,2\overline{\kappa}}\,R_{P_7,-4\overline{\kappa}}(P)=t_{\overrightarrow{P_7R_{P_6,2\overline{\kappa}}(P_7)}}\,t_{\overrightarrow{PR_{P_7,-2\overline{\kappa}}(P)}}(P).$$
On a donc~:
\begin{align*}
  P_5'P'&=\left<\overrightarrow{P_5R_{P_6,2\overline{\kappa}}R_{P_7,-4\overline{\kappa}}(P)},\mathbf{j}\right>\\
  &=\left<\overrightarrow{P_5P}+\overrightarrow{PR_{P_7,-2\overline{\kappa}}(P)}+\overrightarrow{P_7R_{P_6,2\overline{\kappa}}(P_7)},\mathbf{j}\right>\\
  &=\left<\overrightarrow{P_5P_6}+R_{-2\overline{\kappa}}(\overrightarrow{P_7P})+R_{2\overline{\kappa}}(\overrightarrow{P_6P_7}),\mathbf{j}\right>\\
  &=P_5P_6-(PP_7+P_7P_6)\cos 2\overline{\kappa}
\end{align*}
Finalement~:
$$P_5'P'=P_5P_6-PP_6\cos 2\overline{\kappa}.$$

\begin{figure}
  \begin{center}
    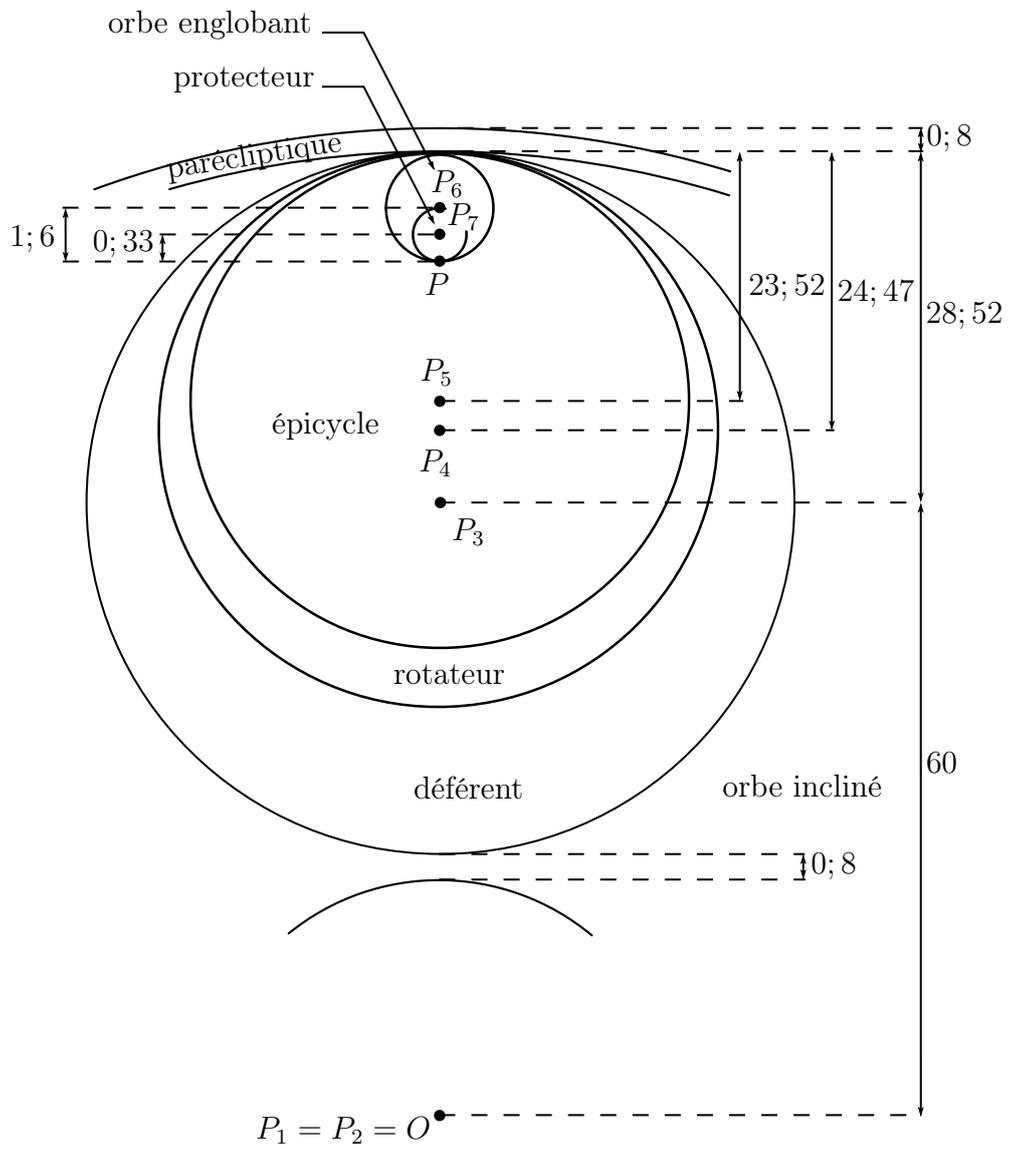
  \end{center}
  \caption{\label{fig050}Mercure, orbes solides}
\end{figure}

\paragraph{Les orbes solides de Mercure}
\textit{Cf.} figure \ref{fig050}.
\begin{align*}
  \text{rayon du globe planétaire}&=r\\
  \text{rayon du protecteur}&=0;33+r\\
  \text{rayon de l'orbe englobant}&=0;33+0;33+r=1;6+r\\
  \text{rayon de l'orbe de l'épicycle}&=22;46+1;6+r=23;52+r\\
  \text{rayon du rotateur}&=0;55+23;52+r=24;47+r\\
  \text{rayon du déférent}&=4;5+24;47+r=28;52+r\\
  \text{rayon extérieur de l'orbe incliné}&=60+28;52+r=88;52+r\\
  \text{rayon intérieur de l'orbe incliné}&=60-(28;52+r)=31;8-r\\
  \text{épaisseur du parécliptique}&=0;8
\end{align*}

\begin{figure}
  \begin{center}
    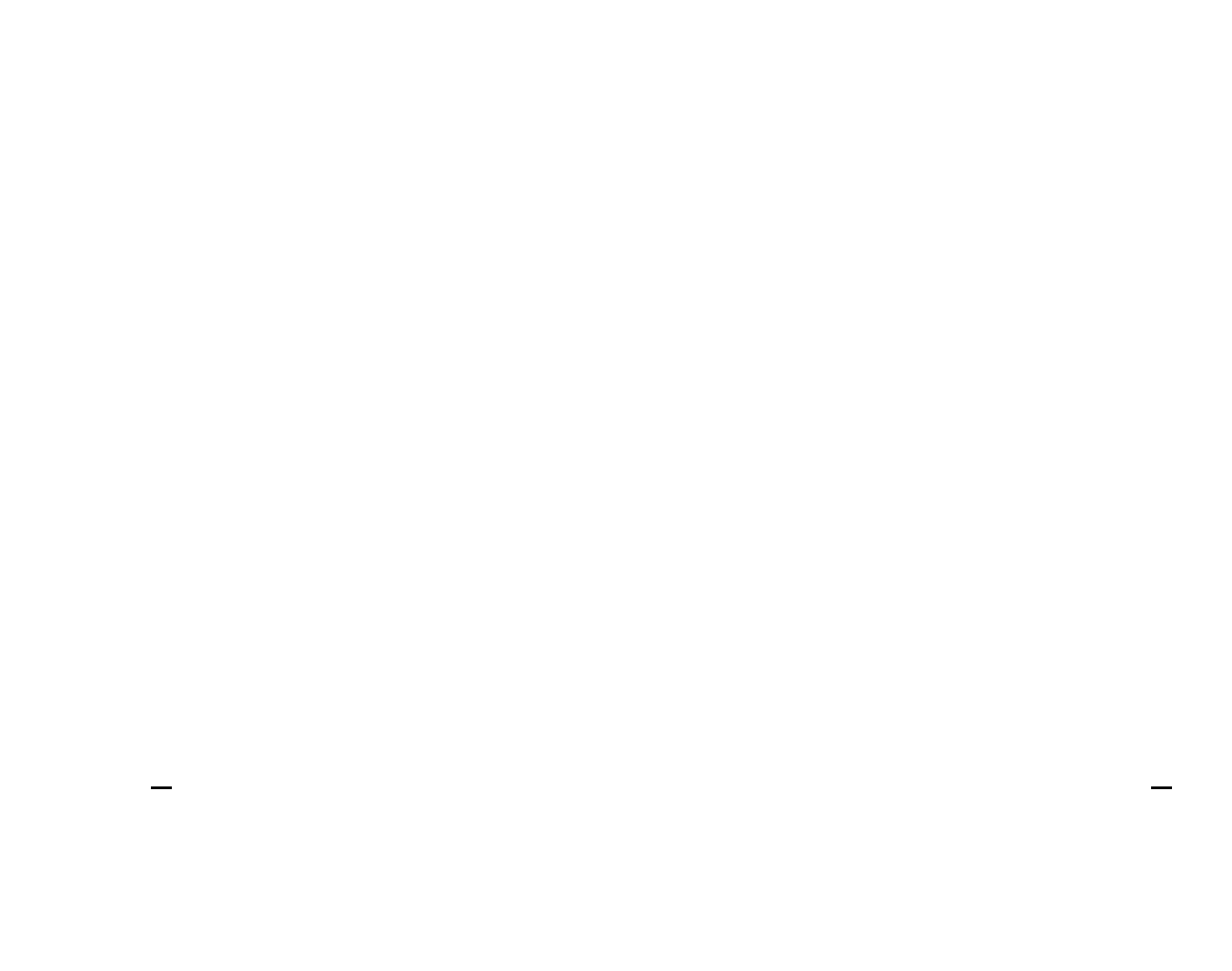
  \end{center}
  \caption{\label{fig062} Les longitudes de Mercure}
\end{figure}

\begin{figure}
  \begin{center}
    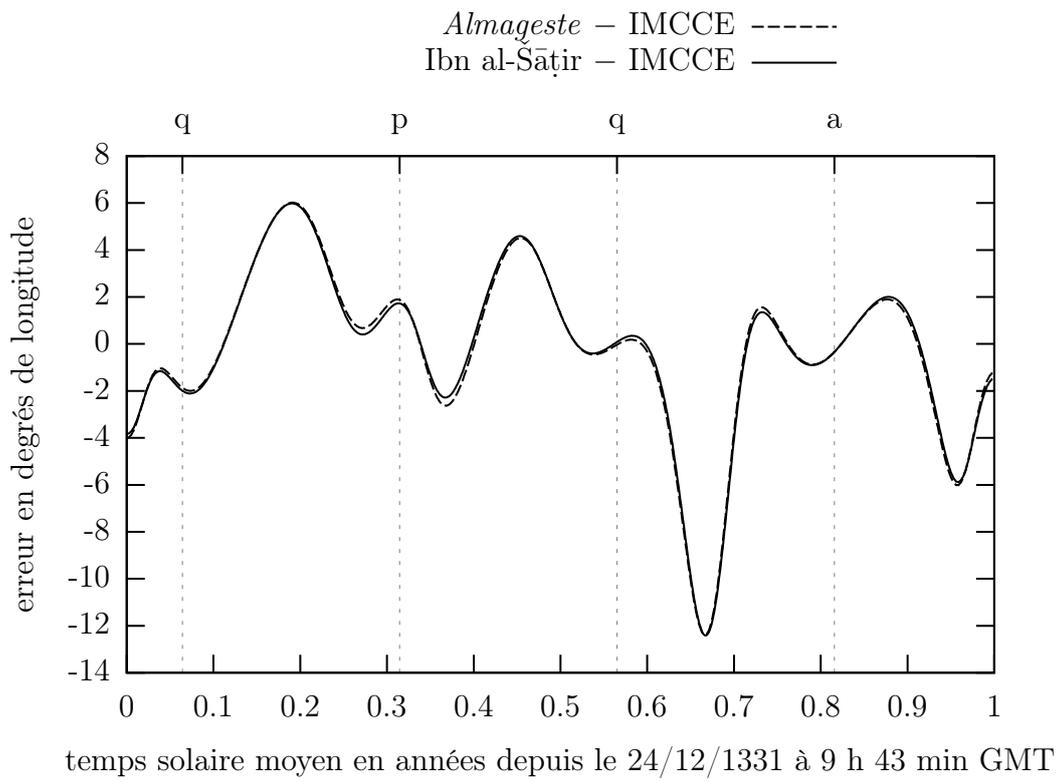
  \end{center}
  \caption{\label{fig063} Erreur en longitude des modèles de Ptolémée et d'{\shatir} pour Mercure}
\end{figure}

\begin{figure}
  \begin{center}
    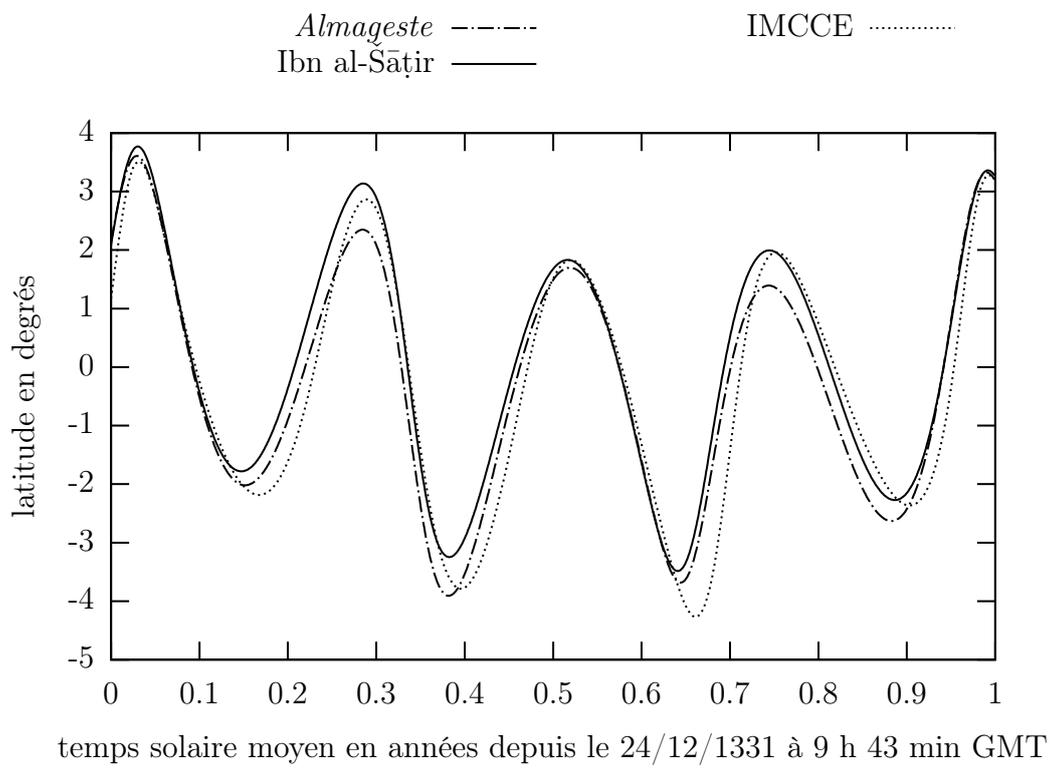
  \end{center}
  \caption{\label{fig064} Les latitudes de Mercure}
\end{figure}

\paragraph{Comparaisons}
\`A la figure \ref{fig062} nous avons tracé l'équation de Mercure sur un an à partir de l'\'Epoque de référence choisie par {\shatir}, (1) calculée au moyen des formules ci-dessus, et (2) calculée au moyen des éphémérides de l'IMCCE. Nous n'avons pas tracé la courbe representant le modèle de l'\textit{Almageste}, car à cette échelle cette courbe est presque confondue avec celle d'{\shatir}. Nous avons repéré sur l'axe supérieur les syzygies et les quadratures au moyen des lettres <<~a~>>, <<~p~>>, <<~q~>>, resp. apogée, périgée, quadrature. Bien que le résultat soit assez bon aux syzygies et aux quadratures, on observe une erreur considérable près des octants~; la figure \ref{fig063} le montre plus clairement, et elle permet de vérifier que le modèle d'{\shatir} n'est pas l'équivalent exact du modèle de l'\textit{Almageste}. Enfin, les latitudes de Mercure sont sur la figure \ref{fig064}.

\paragraph{Le modèle de Qu\d{t}b al-D{\=\i}n al-Sh{\=\i}r\=az{\=\i} pour Mercure} Le modèle de Sh{\=\i}r\=az{\=\i} pour Mercure n'est pas équivalent à celui d'{\shatir}. Il reproduit de manière assez précise les prédictions de l'\textit{Almageste}, mieux encore qu'{\shatir}. Posons $c=3$~; les centres des orbes dans la figure initiale sont définis par\footnote{Le modèle est décrit par Kennedy dans \cite{ghanem1976} p.~101-103.}~:
$$\overrightarrow{P_1P_2}=6\,\mathbf{j},\quad
\overrightarrow{P_2P_3}=60\,\mathbf{j},\quad
\overrightarrow{P_3P_4}=\overrightarrow{P_4P_5}=-\frac{c}{2}\,\mathbf{j}=-1;30\,\mathbf{j},$$
$$\overrightarrow{P_5P_6}=\overrightarrow{P_6P_7}=\frac{c}{2}\,\mathbf{j}=1;30\,\mathbf{j},\quad
\overrightarrow{P_7P_8}=c\,\mathbf{j}=3\,\mathbf{j},\quad
\overrightarrow{P_8P}=22;30\,\mathbf{j}.$$
On a un déférent excentrique de centre $P_2$, un épicycle de centre $P_8$, et les deux triplets d'orbes centrés en $(P_3,P_4,P_5)$ et $(P_5,P_6,P_7)$ constituent deux couples de \d{T}\=us{\=\i}. Les rotations engendrant les mouvements en longitude sont~:
$$R_{P_1,\lambda_A}\,
R_{P_2,\overline{\kappa}}\,
\overbrace{R_{P_3,-\overline{\kappa}}\,
R_{P_4,2\overline{\kappa}}\,
R_{P_5,-\overline{\kappa}}}^{\text{premier couple}}\,
\overbrace{R_{P_5,2\overline{\kappa}}\,
R_{P_6,-4\overline{\kappa}}\,
R_{P_7,2\overline{\kappa}}}^{\text{deuxième couple}}\,
R_{P_7,\overline{\kappa}}\,
R_{P_8,\overline{\alpha}-\overline{\kappa}}$$
$$=R_{P_1,\lambda_A}\,
R_{P_2,\overline{\kappa}}\,
R_{P_3,-\overline{\kappa}}\,
R_{P_4,2\overline{\kappa}}\,
R_{P_5,\overline{\kappa}}\,
R_{P_6,-4\overline{\kappa}}\,
R_{P_7,3\overline{\kappa}}\,
R_{P_8,\overline{\alpha}-\overline{\kappa}}.$$

\paragraph{Calibrage du modèle d'{\shatir} pour Mercure}
En 1974, Willy Hartner \cite{hartner1974} a étudié de près les écarts entre les prédictions ptoléméennes et celles du modèle d'{\shatir} dans les syzygies et les quadratures, et il en a tiré quelques indices concernant la genèse de ce modèle. Nous reprenons ici ses analyses, en formulant d'abord les conditions nécessaires et suffisantes pour qu'un modèle de même forme que celui d'{\shatir} reproduise exactement les élongations maximales observées par Ptolémée dans les syzygies et les quadratures.

Dans le modèle de Ptolémée, notons $R$ le rayon du déférent excentrique portant l'épicycle~: $R$ est la distance entre le centre de l'épicycle et le centre du déférent. Notons $e$ l'excentricité. Le rapport entre les distances Terre--centre de l'épicycle quand $\bar{\kappa}=0°$ et $\bar{\kappa}=180°$ est alors\footnote{\textit{Cf.} \cite{pedersen1974} pour une description de la théorie de Ptolémée.}
$$\frac{R+3e}{R-e}.$$
Dans le modèle d'{\shatir}, ce rapport est
$$\frac{OP_3+(P_3P_4+P_4P_5)}{OP_3-(P_3P_4+P_4P_5)}.$$
Le rayon de l'épicycle apparent d'{\shatir} étant le même en $\bar{\kappa}=0°$ et en $\bar{\kappa}=180°$, le modèle d'{\shatir} reproduira les élongations maximales ($\max\vert c_2(\bar{\kappa},.)\vert$) du modèle de Ptolémée à l'Apogée et au périgée si et seulement si ces deux rapports sont égaux.

Dans les quadratures, Ptolémée avait observé que l'<<~équation du centre~>> est $c_1(90°)=3°$. Dans le modèle de Ptolémée, l'équation du centre se calcule ainsi~:
$$\tan c_1(90°)=\frac{-e}{\sqrt{R^2-e^2}-e}.$$
Dans le modèle d'{\shatir}, on a~:
$$\tan c_1(90°)=\frac{-(P_3P_4-P_4P_5)}{OP_3}.$$
Ces deux grandeurs seront égales si et seulement si
$$P_3P_4-P_4P_5=\frac{e\times OP_3}{\sqrt{R^2-e^2}-e}.$$

Avec $OP_3=R=60$, on obtient, comme Hartner\footnote{\textit{Cf.} \cite{hartner1974} p.~15. Nous reprenons ci-dessous les calculs menés par Hartner, car il les a menés sous l'hypothèse que $P_4P_5=0;50$, au lieu de $0;55$. Ses conclusions restent valables.}, les relations suivantes~:
\begin{align}
  P_3P_4-P_4P_5=\frac{R(C-1)}{C+1}\\
  e=\frac{R(C-1)}{C+3}
\end{align}
où $C=\dfrac{R+3e}{R-e}=\dfrac{OP_3+(P_3P_4+P_4P_5)}{OP_3-(P_3P_4+P_4P_5)}$, et
\begin{align}
  P_3P_4-P_4P_5=\frac{Re}{\sqrt{R^2-e^2}-e}
\end{align}
Cette dernière relation peut être inversé pour exprimer $e$ en fonction de $P_3P_4-P_4P_5$~:
\begin{align}
  e=\frac{R(P_3P_4-P_4P_5)}{\sqrt{R^2+2R(P_3P_4-P_4P_5)+2(P_3P_4-P_4P_5)^2}}
\end{align}
Comme les rayons des orbes choisis par {\shatir} impliquent que\linebreak ${\dfrac{OP_3+(P_3P_4+P_4P_5)}{OP_3-(P_3P_4+P_4P_5)}=\dfrac{65}{55}}$, on aurait, d'après (2), $e=2;36,31$. D'après (3), on aurait alors~:
$$P_3P_4-P_4P_5=2;43,38,$$
d'où
$$c_1(90°)=\arctan\frac{-(P_3P_4-P_4P_5)}{OP_3}=-2;36,18.$$
Au contraire, selon {\shatir}, on a~:
$$P_3P_4-P_4P_5=3;10,$$
d'où
$$c_1(90°)=3;1,16.$$
Hartner en tire une première conclusion~: {\shatir} n'a pas calibré son modèle pour reproduire exactement les élongations maximales observées par Ptolémée aux syzygies et dans les quadratures.

\'Etrangement, Hartner, peut-être trop prudent, juge <<~peu probable~>>\footnote{\textit{Cf.} \cite{hartner1974} p.~23.} qu'{\shatir} n'ait pas eu l'intention de concevoir un véritable <<~substitut géométrique~>> du mécanisme de l'\textit{Almageste}. La suite de son article tend pourtant à convaincre qu'{\shatir} aurait au moins utilisé une observation nouvelle à l'apogée~; et si l'on veut préserver l'adéquation du modèle dans les quadratures, cette observation nouvelle conduit à abandonner les relations (1) à (4). Ptolémée avait observé une élongation maximale à l'apogée, $\max\vert c_2(0,\cdot)\vert=19;3$. Le modèle d'{\shatir} donne\footnote{Hartner trouve $19;30$ sous l'hypothèse que $P_4P_5=0;50$. On trouve $19;28$ en prenant $P_4P_5=0;55$.} $\max\vert c_2(0,\cdot)\vert=19;28$. Cette supposée observation nouvelle sera notre point de départ pour tenter de reconstituer le calibrage du modèle d'{\shatir}.

Quelles sont les données de l'observation~? On a d'abord cette nouvelle observation, qu'on arrondira à $19;30$~:
$$\max\vert c_2(0,\cdot)\vert=\arcsin\frac{P_5P_6-2\times P_6P_7}{OP_3+(P_3P_4+P_4P_5)}=19;30.$$
L'\textit{Almageste} donne aussi une observation de l'élongation maximale au périgée~:
$$\max\vert c_2(180°,\cdot)\vert=\arcsin\frac{P_5P_6-2\times P_6P_7}{OP_3-(P_3P_4+P_4P_5)}=23;15.$$
Enfin, des observations de l'élongation maximale aux quadratures, Ptolémée déduit que $c_1(90°)=3°$ et que $\max\vert c_2(90°,\cdot)\vert=23;15$.
On déduit des deux premières données~:
$$\frac{OP_3-(P_3P_4+P_4P_5)}{OP_3+(P_3P_4+P_4P_5)}=\frac{\sin(19°30')}{\sin(23°15')},$$
d'où, si $OP_3=60$,
$$P_3P_4+P_4P_5=5;1,7.$$
En reportant dans $\max\vert c_2(0,\cdot)\vert$, on en déduit~:
$$P_5P_6-2\times P_6P_7=(60+P_3P_4+P_4P_5)\sin(19°30')=21;42,13.$$
De $c_1(90°)$, on tire~:
$$P_3P_4-P_4P_5=60\tan(3°)=3;8,40,$$
d'où enfin, comme il se doit~:
$$P_3P_4=4;4,53\simeq 4;5,$$
$$P_4P_5=0;56,13\simeq 0;55.$$
Reste à calculer $P_5P_6+2\times P_6P_7$. Or~:
$$\max\vert c_2(90°,\cdot)\vert=23;15=\arcsin\frac{P_5P_6+2\times P_6P_7}{\sqrt{(P_3P_4-P_4P_5)^2+60^2}},$$
d'où
$$P_5P_6+2\times P_6P_7=\sin(23°15')\times\sqrt{(P_3P_4-P_4P_5)^2+60^2}=23;43;2.$$
On devrait alors avoir $2\times P_6P_7=1;1$, d'où $P_6P_7=P_7P=0;30,30$ au lieu de $0;33$ d'{\shatir}~; mais il est bien possible qu'{\shatir} ait observé une équation maximale légèrement supérieure à celle de Ptolémée dans les quadratures. Si $\max\vert c_2(90°,\cdot)\vert=23;30$, sans rien changer par ailleurs, on trouverait en effet $P_6P_7\simeq 0;34$, arrondi qu'on pourrait aussi tronquer à $0;33$.\label{end_mercure}

\paragraph{Tables d'équations en longitude} Les formules p.~\pageref{equ_lune}, \pageref{equ_sat} et \pageref{equ_mercure}, ainsi que les formules d'interpolation afférentes, permettent de restituer les tables d'équations mentionnées aux chapitres 11, 14 et 23 de la première partie de la \textit{Nihaya al-s\=ul} comme nous l'avons fait dans la table \ref{jadwal}. Pour faciliter d'éventuelles comparaisons, nous avons adopté le même format que Fuad Abbud dans \cite{ghanem1976} p.~80. 

\begin{table}
  \begin{center}\begin{tabular}{cc|cccc}
    Lune & $2\bar{\eta}$ ou $\alpha$ & $c_1(2\bar{\eta})$ & $c_2(\alpha,0)$ & $c_2(\alpha,180°)-c_2(\alpha,0)$ & $\chi(2\bar{\eta})$ \\
    & $30$ & $7;32$ & $-2;18$ & $-1;8$ & $0;5$ \\
    & $60$ & $11;48$ & $-4;5$ & $-2;5$ & $0;18$ \\
    & $90$ & $12;9$ & $-4;55$ & $-2;40$ & $0;33$ \\
    & $120$ & $9;33$ & $-4;27$ & $-2;36$ & $0;47$ \\
    & $150$ & $5;11$ & $-2;40$ & $-1;39$ & $0;57$ \\
    \hline
    Saturne & $\bar{\kappa}$ ou $\alpha$ & $c_1(\bar{\kappa})$ & $c_2(0,\alpha)$ & $c_2(180,\alpha)-c_2(0,\alpha)$ & $\chi(\bar\kappa)$ \\
    & $30$ & $-3;6$ & $2;42$ & $0;18$ & $0;3$ \\
    & $60$ & $-5;29$ & $4;50$ & $0;33$ & $0;11$ \\
    & $90$ & $-6;30$ & $5;51$ & $0;42$ & $0;25$ \\
    & $120$ & $-5;48$ & $5;21$ & $0;41$ & $0;41$ \\
    & $150$ & $-3;26$ & $3;13$ & $0;26$ & $0;55$ \\
    \hline
    Jupiter & $\bar{\kappa}$ ou $\alpha$ & $c_1(\bar{\kappa})$ & $c_2(0,\alpha)$ & $c_2(180,\alpha)-c_2(0,\alpha)$ & $\chi(\bar\kappa)$ \\
    & $30$ & $-2;31$ & $4;31$ & $0;22$ & $0;3$ \\
    & $60$ & $-4;26$ & $8;16$ & $0;43$ & $0;12$ \\
    & $90$ & $-5;14$ & $10;23$ & $0;58$ & $0;26$ \\
    & $120$ & $-4;39$ & $9;55$ & $1;2$ & $0;42$ \\
    & $150$ & $-2;44$ & $6;13$ & $0;43$ & $0;55$ \\
    \hline
    Mars & $\bar{\kappa}$ ou $\alpha$ & $c_1(\bar{\kappa})$ & $c_2(0,\alpha)$ & $c_2(180,\alpha)-c_2(0,\alpha)$ & $\chi(\bar\kappa)$ \\
    & $30$ & $-5;15$ & $11;9$ & $1;28$ & $0;2$ \\
    & $60$ & $-9;22$ & $21;45$ & $3;8$ & $0;9$ \\
    & $90$ & $-11;19$ & $30;54$ & $5;17$ & $0;20$ \\
    & $120$ & $-10;20$ & $36;29$ & $8;29$ & $0;36$ \\
    & $150$ & $-6;15$ & $31;51$ & $13;5$ & $0;53$ \\
    \hline
    Vénus & $\bar{\kappa}$ ou $\alpha$ & $c_1(\bar{\kappa})$ & $c_2(0,\alpha)$ & $c_2(180,\alpha)-c_2(0,\alpha)$ & $\chi(\bar\kappa)$ \\
    & $30$ & $-1;0$ & $12;25$ & $0;19$ & $0;4$ \\
    & $60$ & $-1;44$ & $24;26$ & $0;40$ & $0;14$ \\
    & $90$ & $-2;1$ & $35;25$ & $1;8$ & $0;28$ \\
    & $120$ & $-1;46$ & $43;42$ & $1;52$ & $0;44$ \\
    & $150$ & $-1;2$ & $42;47$ & $3;13$ & $0;55$ \\
    \hline
    Mercure & $\bar{\kappa}$ ou $\alpha$ & $c_1(\bar{\kappa})$ & $c_2(0,\alpha)$ & $c_2(180,\alpha)-c_2(0,\alpha)$ & $\chi(\bar\kappa)$ \\
    & $30$ & $-1;25$ & $7;22$ & $0;59$ & $0;12$ \\
    & $60$ & $-2;31$ & $13;54$ & $2;1$ & $0;39$ \\
    & $90$ & $-3;1$ & $18;26$ & $3;4$ & $1;3$ \\
    & $120$ & $-2;44$ & $19;6$ & $3;55$ & $1;11$ \\
    & $150$ & $-1;38$ & $13;11$ & $3;27$ & $1;5$ \\
    & $180$ & $0$ & $0$ & $0$ & $1$ 
  \end{tabular}\end{center}
  \caption{\label{jadwal}Restitution d'un extrait des tables d'équations en longitude d'{\shatir}}
\end{table}

\label{comm_fin}

\pagebreak
\begin{center}
  \Large Numéraux
\end{center}

$$\RL{'a}=1,\,\RL{b}=2,\,\RL{j}=3,\,\RL{d}=4,\,\RL{h h-}=5,$$

$$\RL{w}=6,\,\RL{z}=7,\,\RL{.h}=8,\,\RL{.t}=9,$$

$$\RL{_A y}=10,\,\RL{k}=20,\,\RL{l}=30,\,\RL{m}=40,\,\RL{n}=50,$$

$$\RL{s .s}=60,\,\RL{`}=70,\,\RL{f}=80,\,\RL{.s .z}=90,$$

$$\RL{q}=100,\,\RL{r}=200,\,\RL{^s s}=300,\,\RL{t _t}=400,\,\RL{_t}=500,$$

$$\RL{_h}=600,\,\RL{_d d}=700,\,\RL{.d}=800,\,\RL{.z .g}=900$$

$$\RL{.g ^s}=1000$$

\pagebreak

\newpage
\phantomsection\addcontentsline{toc}{part}{Index}
\printindex

\end{document}